\documentclass[a4paper,11pt]{amsart}

\pdfoutput=1

\usepackage[text={440pt,660pt},centering]{geometry}
\usepackage{amsthm, amssymb, amsmath, amsfonts, mathrsfs}
\usepackage{mathtools}
\usepackage{scalerel} 
\usepackage{stackrel}
\usepackage[colorlinks=true, pdfstartview=FitV, linkcolor=blue, citecolor=blue, urlcolor=blue,pagebackref=false]{hyperref}
\numberwithin{equation}{section}
\usepackage[export]{adjustbox} 
\usepackage{caption}
\usepackage{subcaption}
\usepackage{enumitem}
\usepackage{multicol} 
\usepackage{booktabs} 
\usepackage{subcaption} 

\usepackage{tikz} 
\usetikzlibrary{shapes,arrows,positioning,calc}

\tikzstyle{result} = [rectangle, 
minimum width=4cm, 
minimum height=1cm,  
text width=4cm, 
align=flush center,
font = \footnotesize, 
draw=black]

\tikzstyle{input} = [rectangle, 
minimum width=3cm, 
minimum height=1cm, 
text width=3cm, 
align=flush center,
font = \footnotesize,
draw=black]

\tikzstyle{arrow} = [thick, ->,>=stealth]

\newtheorem{theorem}{Theorem}[section]
\newtheorem{remark}[theorem]{Remark}
\newtheorem{lemma}[theorem]{Lemma}

\newtheorem{prop}[theorem]{Proposition}
\newtheorem{definition}[theorem]{Definition}
\newtheorem{corollary}[theorem]{Corollary}

\definecolor{labelkey}{gray}{.8}
\definecolor{refkey}{gray}{.8}

\definecolor{darkgreen}{rgb}{0,0.5,0}
\definecolor{darkblue}{rgb}{0,0,0.7}
\definecolor{darkred}{rgb}{0.9,0.1,0.1}

\newcommand{\gucomment}[1]{\marginpar{\raggedright\scriptsize{\textcolor{darkblue}{#1}}}}

\newcommand{\cCf}{\cC_{\infty}} 

\newcommand{\EE}{\mathcal{E}} 
\renewcommand{\le}{\leqslant}
\renewcommand{\ge}{\geqslant}
\renewcommand{\leq}{\leqslant}
\renewcommand{\geq}{\geqslant}

\renewcommand{\subset}{\subseteq}
\renewcommand{\supset}{\supseteq}

\newcommand{\G}{\mathcal{G}}
\newcommand{\D}{\mathcal{D}}

\newcommand{\E}{\mathbb{E}}

\newcommand{\bfv}{\mathbf{v}}
\renewcommand{\a}{\mathbf{a}}
\newcommand{\A}{\mathbf{A}} 
\newcommand{\Reff}{\mathsf{R}_{\operatorname{eff}}} 
\newcommand{\Ceff}{\mathsf{C}_{\operatorname{eff}}} 

\newcommand{\Ab}{{\overbracket[1pt][-1pt]{\A}}}
\renewcommand{\b}{\mathbf{b}}

\newcommand{\ab}{{\overbracket[1pt][-1pt]{\a}}}
\newcommand{\bb}{{\overbracket[1pt][-1pt]{\b}}}
\newcommand{\cc}{\mathsf{c}}

\newcommand{\Ll}{\left}
\newcommand{\Rr}{\right}
\newcommand{\lhs}{left-hand side}
\newcommand{\rhs}{right-hand side}

\newcommand{\R}{\mathbb{R}}

\newcommand{\CC}{\mathcal{C}} 
\newcommand{\cC}{\mathscr{C}}
\newcommand{\N}{\mathbb{N}}

\newcommand{\ZZ}{\mathbb{Z}}
\renewcommand{\H}{\mathbb{H}}
\newcommand{\Zd}{{\mathbb{Z}^d}}
\newcommand{\Ed}{E_d}

\renewcommand{\P}{\mathbb{P}}
\newcommand{\PP}{\mathbf{P}} 
\newcommand{\Eb}{\mathbf{E}} 
\newcommand{\Pp}{\mathbb{P}_\p}

\newcommand{\kk}{\mathsf{k}}

\newcommand{\tT}{\mathcal{T}}

\newcommand{\fF}{\mathcal{F}}

\newcommand{\calC}{\mathcal{C}}

\renewcommand{\bar}{\overline}

\renewcommand{\tilde}{\widetilde}
\newcommand{\eps}{\varepsilon}

\renewcommand{\epsilon}{\varepsilon}
\renewcommand{\d}{{\mathrm{d}}}
\newcommand{\var}{\mathbb{V}\!\mathrm{ar}}
\newcommand{\cov}{\mathbb{C}\mathrm{ov}}

\newcommand{\p}{\mathfrak{p}}
\newcommand{\pc}{\p_c}

\newcommand{\Ind}[1]{\mathbf{1}_{\left\{#1\right\}}}
\newcommand{\id}{\mathsf{Id}}

\newcommand{\bracket}[1]{\left\langle{#1}\right\rangle}

\newcommand{\fC}{\mathfrak{C}}
\newcommand{\cK}{\mathcal{K}}

\newcommand{\RD}{\mathcal{RD}}

\newcommand{\piv}{\mathsf{Piv}}
\newcommand{\pivo}{\mathsf{Piv}_{o}}
\newcommand{\spn}{\mathsf{SPN}}
\newcommand{\skl}{\mathsf{SKL}}
\newcommand{\cP}{\mathcal{P}}
\newcommand{\expl}{\mathsf{Explo}}
\newcommand{\tarm}{2\text{-}\mathsf{arm}}
\newcommand{\V}{\mathbf{V}}
\newcommand{\EDGE}{\mathcal{E}}
\newcommand{\DIEDGE}{\mathbf{E}^\rightarrow}
\newcommand{\sE}{\mathscr{E}}
\newcommand{\sV}{\mathscr{V}}
\newcommand{\disE}{\sE^\rightarrow}
\newcommand{\ERE}{\mathcal{E}^{\mathrm{RE}}}
\newcommand{\PTRE}[1]{\partial^{#1\text{-}\mathrm{RE}}}
\newcommand{\PTR}{\partial^{\mathrm{RE}}}

\DeclareMathOperator{\supp}{supp}

\DeclareMathOperator{\dex}{d_{ex}}
\DeclareMathOperator{\din}{d_{in}}

\DeclareMathOperator{\rin}{r_{in}}

\title[
Diffusivity and metric  for slightly supercritical percolation
]{
	Diffusivity and metric exponents for slightly supercritical percolation in high dimensions
}
\author{Chenlin Gu, Tom Hutchcroft, Xiaoran Li}

\address[Chenlin Gu]{Yau Mathematical Sciences Center, Tsinghua University}
\email{gclmath@tsinghua.edu.cn}

\address[Tom Hutchcroft]{Department of Mathematics, Princeton University and The Division of Physics, Mathematics and Astronomy, California Institute of Technology}
\email{t.hutchcroft@princeton.edu}

\address[Xiaoran Li]{The Division of Physics, Mathematics and Astronomy, California Institute of Technology}
\email{xiaoran@caltech.edu}

\begin{document}
\begin{abstract}
	It is a classical fact that random walk on the infinite supercritical Bernoulli $\p$-percolation cluster on $\Zd$ converges to Brownian motion with some diffusivity $\sigma(\p)$, and that the graph metric on the infinite cluster has order $\mu(\p) \|\cdot\|$ at large scales for some constant $\mu(\p)$ known as the time constant.
	We compute the near-critical asymptotics of the two quantities in high dimensions, showing that
	\begin{align*}
		\sigma^2(\p) \asymp  (\p-\p_c)^2 \qquad \text{ and } \qquad \mu(\p) \asymp \frac{1}{\sqrt{\p-\p_c}}
	\end{align*}
	for all $\p>\p_c$. Our results concerning the diffusivity confirm a conjecture of Biskup (2011). We also compute the asymptotic critical probabilities of slabs in the same regime, showing that $\p_c( \mathbb{Z}^{d-1} \times \{0,\ldots,n\}) - \p_c(\Zd) \asymp n^{-2}$ as $n\to \infty$.

	\bigskip
	
	
	\medskip
	

\end{abstract}

\maketitle

\setcounter{tocdepth}{1}
\tableofcontents


\section{Introduction}

\subsection{Diffusion in near-critical high-dimensional percolation} 
Bernoulli bond percolation was introduced by Broadbent and Hammersley \cite{broadbent1957percolation} to study flow in porous materials and now occupies a role of fundamental importance in mathematical statistical mechanics. Within percolation theory, a particular fascination has been reserved for \emph{random walks on percolation clusters}, a topic memorably named ``the ant in the labyrinth'' by de Gennes \cite{de1976percolation}.

Let us now briefly introduce the model.
Let $\Zd$ be the $d$-dimensional integer lattice and $\Ed$ be the associated edges (a.k.a.\ bonds) which connect vertices with $\|x-y\|_1=1$.
We sample bonds independently at random, with probability $\p \in [0,1]$ to be \emph{open} (retained) and probability $(1-\p)$ to be \emph{closed} (deleted), writing $\a=\{\a(e)\}_{e\in \Ed}$ for the associated collection of independent random variables $\a(e)=\mathbf{1}(e$ open$)$.
We are interested in the random graph ${(\Zd, \{e\in\Ed:\a(e)=1\})}$ spanned by the open edges of $\Zd$, whose connected components are referred to as \emph{clusters}. For $d \geq 2$, there exists a \emph{critical point} $\p_c(d) \in (0,1)$ such that there is no infinite cluster almost surely when $\p<\p_c(d)$ and a unique infinite cluster $\cCf$ almost surely when $\p>\p_c(d)$
%
\cite{aizenman1987uniqueness, burton1989density}. These two phases are known as the \emph{subcritical} and \emph{supercritical} phases respectively.
We refer to e.g.\ Grimmett's monograph~\cite{grimmett1999} and the proceedings of Duminil-Copin \cite{duminil2017sixty} for a detailed introduction of the topic, and to e.g.\ \cite{heydenreich2017progress, hutchcroft2026dimension} for introductions to \emph{high-dimensional} percolation specifically.

An important basic feature of the percolation phase transition 
is that the infinite supercritical cluster $\cCf$ has geometry close to that of the underlying lattice $\Zd$ at large scales~\cite{grimmett1990supercritical,antal1996chemical,pete2008note}. 
One particularly interesting aspect of this is that the random walk on $\cCf$ behaves similarly at large scales to the random walk on $\Zd$ and hence converges to Brownian motion under scaling.
One precise formulation concerns the \emph{variable-speed random walk} $\left( X_t \right)_{t \geq 0}$, that is, the continuous-time Markov jump process associated to the generator 
$\nabla \cdot \a \nabla$ defined by 
\begin{equation} \label{eq.VSRW}
	\forall x \in \cCf, \qquad 
	\nabla \cdot \a \nabla u(x) := \sum_{y: \{x,y\} \in \Ed} \a(\{x,y\}) \left( u(y) - u(x) \right).
\end{equation}
It is a now-classical result that this process satisfies a \emph{quenched invariance principle} for all $d\geq 2$ and $\p>\p_c(d)$, meaning that
there exists a deterministic constant $\sigma(\p) > 0$ such that, for $\P_\p$-almost every percolation configuration $\{\a(e)\}_{e \in \Ed}$, the convergence in law 
\begin{equation} \label{eq.IP}
	\Ll(\frac{1}{\sqrt{N}} X_{N t}\Rr)_{t \geq 0} \Longrightarrow ( \sigma(\p) B_{t})_{t \geq 0} \qquad \text{ as $N\to\infty$}
\end{equation}
holds under the Skorokhod topology,
where $(B_{t})_{t \geq 0}$ is the $d$-dimensional Brownian motion.  This was first proven for dimensions $d\geq 4$ by
Sidoravicius and Sznitman \cite{sidoravicius2004quenched} then for all dimensions $d\geq 2$ by Berger and Biskup \cite{berger2007quenched} and Mathieu and Piatnitski \cite{mathieu2007quenched}. (See also the work of  Barlow \cite{barlow2004random} establishing Gaussian bounds on the transition probabilities of $X$.)
This invariance principle has since been extended to other settings; see the survey of Biskup \cite{biskup2011recent} and the very recent survey of Andres \cite{andres2025homogenization}.

The diffusive, Brownian scaling of the random walk on the infinite supercritical cluster is in stark contrast to the \emph{subdiffusive} behaviour of the random walk on \emph{large critical clusters}.
That is, while the $n$-step random walk on the infinite supercritical cluster has typical displacement of order $\sigma(\p) n^{1/2}$ for large $n$, the random walk on a large critical cluster is conjectured to have typical displacement of order $n^{\varkappa+o(1)}$ for some $\varkappa=\varkappa_d<1/2$ when $n$ is not too large relative to the size of the cluster.\footnote{These results are usually formulated in terms of the \emph{incipient infinite cluster} (IIC), i.e. the limit of the critical cluster conditioned to be very large \cite{kesten1986incipient,van2004incipient,chatterjee2025robust}.} This subdiffusivity is best-understood in \emph{high dimensions}, where Kozma and Nachmias \cite{KN09} (see also \cite{barlow2008random}) proved that the typical \emph{intrinsic displacement} of the random walk on a large critical cluster is of order $n^{1/3}$ (meaning displacement with respect to the graph distance on the cluster), which implies that the typical Euclidean displacement is of order $n^{1/6}$ since geodesics in large critical clusters themselves have Brownian scaling in high dimensions \cite{blanc2026super,chatterjee2025limiting,KN11}. The problem is also reasonably well-understood in two dimensions, where Kesten \cite{kesten86} proved that the typical displacement of the random walk on a large critical cluster is $O(n^{1/2-\eps})$ with high probability (see also \cite{ganguly2022chemical}), and 
{\DJ}ankovi{\'c},  Markering, Miller, and Yuan  \cite{djankovic2026scaling} recently proved (among many other things) that the exponent of subdiffusivity is well-defined for critical site percolation on the triangular lattice.
The high dimensional results of \cite{KN09} are expected to hold in all dimensions strictly above the \emph{upper critical dimension} $d_c=6$, but current proof methods \cite{hara1990mean,DCP25a} involve perturbative hypotheses that hold for the standard hypercubic lattice only in very high dimensions, the current record being $d\geq 11$~\cite{fitzner2017mean}.
Analogous phenomena remain poorly understood in dimensions $d=3,4,5,6$. 

In this paper we are interested in the crossover between diffusive scaling of the random walk in the supercritical phase and subdiffusive scaling at criticality, which is captured in part by the behaviour of the diffusivity $\sigma(\p)$ as $\p$ approaches $\p_c$ from above.
As with many other important quantities in statistical physics, the near-critical diffusivity is expected to have an associated critical exponent $t-\beta=t_d-\beta_d$ such that \footnote{This exponent is written as the difference of two exponents as the exponents are normally formulated for the quantities $\ab(\p) \approx (\p-\p_c)^t$ and $\theta(\p)\approx (\p-\p_c)^\beta$, with the diffusivity related to these quantities via the Einstein relations as explained below.}
\begin{align}\label{eq.sigma.conjecture}
	\sigma^2(\p) \approx (\p - \p_c)^{t-\beta}, \qquad \p \searrow \p_c,
\end{align}
where we keep the precise meaning of $\approx$ intentionally vague (in particular, this notation hides conjectured logarithmic corrections to power-law scaling in six dimensions). 
This critical exponent has attracted sustained interest across several disciplines: See \cite[Eq.\ (44)]{Hughes2021} for a survey of the applied mathematics literature and \cite[(9.16)]{torquato2002} for a viewpoint from engineering and material science. Within pure  mathematics, perhaps the best-known conjecture about \eqref{eq.sigma.conjecture} is that $t-\beta = 2$ for $d > 6$, which was proposed by Biskup in \cite[Problem 4.21]{biskup2011recent}  and reiterated by Heydenreich and van der Hofstad in \cite[Open Problem 14.1]{heydenreich2017progress}; see also the related question of Kesten regarding the case $d=2$ in \cite{kestenICM84}.

\smallskip

In this paper we prove this conjecture in very high dimensions. More concretely, we work under the assumption that $d\geq 9$ and that the critical two-point function satisfies the Green's function asymptotics
\begin{equation}
	\label{eq:two_point_assumption}
	\tag{A}
	\mathbb{P}_{\p_c}(x\leftrightarrow y) \sim A\|x-y\|^{-d+2}_2
\end{equation}
for some constant $A>0$, as is conjectured to hold in all dimensions $d>6$ and proven to hold in dimensions $d\geq 11$ \cite{fitzner2017mean,hara1990mean,10.1214/009117907000000231}.

\begin{theorem}
	\label{thm.main}
	Consider nearest-neighbor Bernoulli bond percolation on $\Zd$. If  $d\geq 9$ and \eqref{eq:two_point_assumption} holds then there exist positive constants $C, c$ such that 
	\begin{align}\label{eq.main}
		c(\p-\p_c)^2 \leq \sigma^2(\p) \leq C(\p-\p_c)^2.
	\end{align} 
	for every $\p_c<\p \leq 1$.
\end{theorem}

\begin{remark}
	The restriction of our results to the case $d\geq 9$ is artificial and we plan to extend them to the case $d>6$ in a future version of the paper. This restriction to dimensions $d\geq 9$ enters only through Proposition~\ref{prop.sprinkling.AB.disjoint}, and is not used in the upper bound.
	We expect that our proof methods also apply with minor modifications to \emph{spread out} percolation in all dimensions $d>6$, i.e., percolation on the graph with vertex set $\Zd$ and edge set $\{\{x,y\}:\|x-y\|_1 \leq L\}$ for some constant $L\geq 1$. In this context, the estimate \eqref{eq:two_point_assumption} is known to hold in all dimensions $d>6$ once the parameter $L$ is sufficiently large \cite{hara2003critical}. We restrict attention to nearest-neighbour percolation on $\Zd$ for clarity of exposition.
\end{remark}

\begin{figure}
	\centering
	\includegraphics[width=0.49\textwidth]{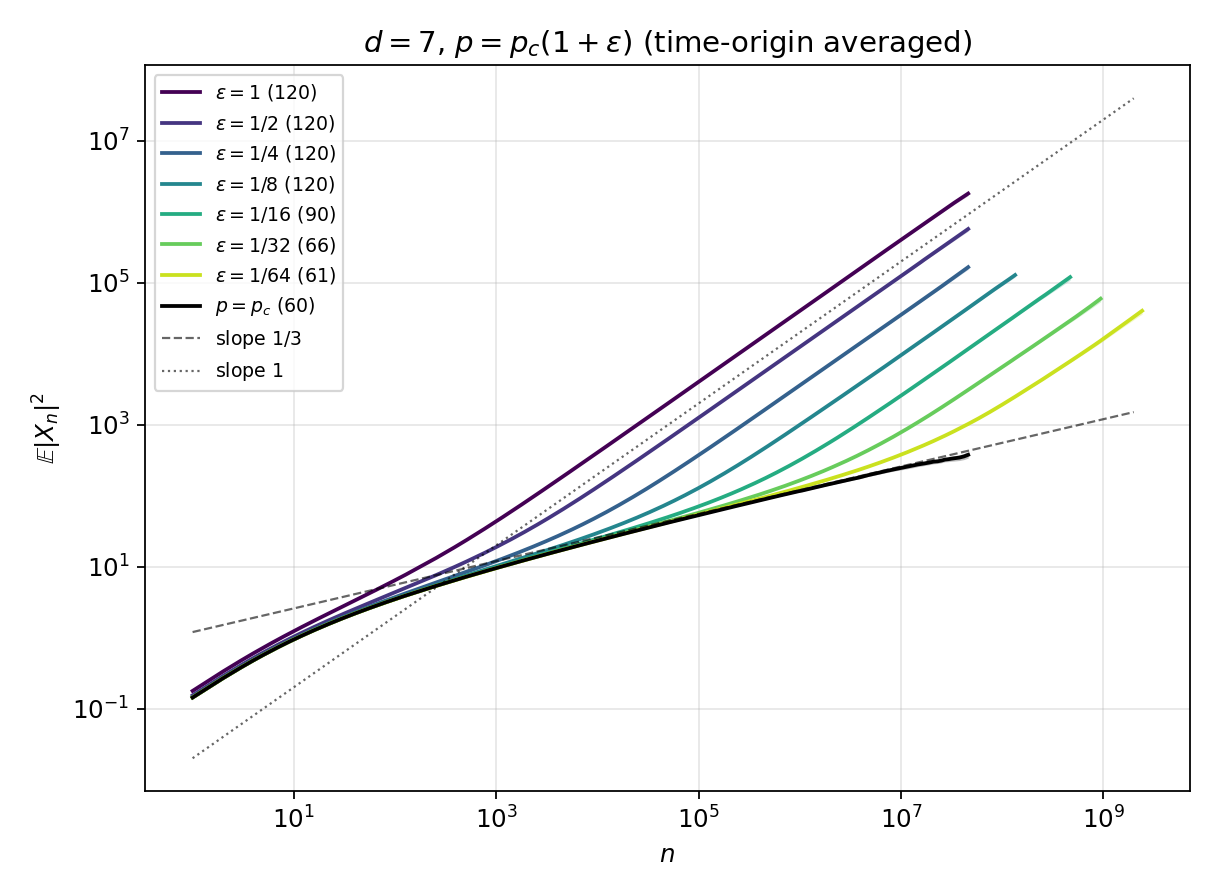}
	\includegraphics[width=0.49\textwidth]{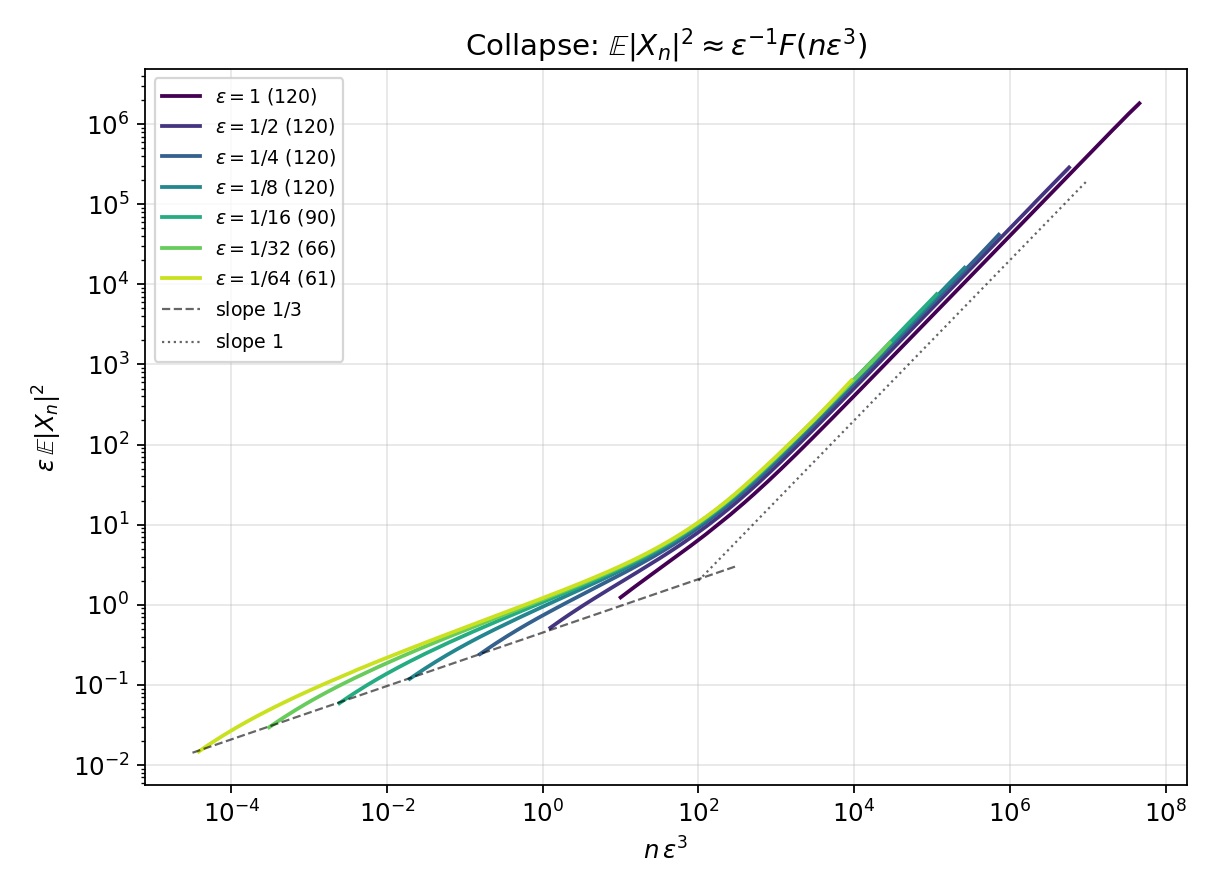}
	\caption{Numerical estimates of $\E_{(1+\eps)\p_c} [\|X_n\|_2^2 \mid X_0 \leftrightarrow \infty]$ on $\mathbb{Z}^7$. Left: log-log plots of the variance  $\E_{(1+\eps)\p_c} [\|X_n\|_2^2 \mid X_0 \leftrightarrow \infty]$ against $n$ for $\eps=1,1/2,1/4,1/16,1/32,0$. Right: Plotting the same data as a log-log plot of the normalized quantities $\eps \E_{(1+\eps)\p_c} [\|X_n\|_2^2 \mid X_0 \leftrightarrow \infty]$ against $\eps^3 n$ highlights the transition between subdiffusive and diffusive scaling around $n\approx 100 \eps^{-3}$. 
		The asymptotic estimate $\sigma^2((1+\eps)\p_c) \asymp \eps^2$ of Theorem~\ref{thm.main} is related to the apparent emergence of a non-degenerate limiting curve asymptotic to a line of slope $1$ at $+\infty$ under this scaling as $\eps\downarrow 0$.
		These simulations use the numerical estimate $\p_c(\mathbb{Z}^7) \approx 0.078 675 230(2)$
		due to Mertens and Moore \cite{mertens2018percolation}  and work by conditioning on the cluster of the origin having an appropriately chosen large finite diameter.}
	\label{fig:numerics}
\end{figure}

The heuristic behind Theorem~\ref{thm.main} is simple and well-known (see \cite[Page~230]{heydenreich2017progress} for example): off-critical percolation in high dimensions is expected to have a single \emph{correlation length} $\xi(\p)$ such that clusters ``look critical'' at scales much smaller than $\xi(\p)$ and are ``obviously off-critical'' at scales much larger than $\xi(\p)$. We warn the reader that there are many reasonable precise definitions of the correlation length that are not obviously (up-to-constants) equivalent. Since paths in large high-dimensional critical clusters look like Brownian motion and hence take order $r^2$ steps to reach distance $r$ \cite{blanc2026super,chatterjee2025limiting,KN11}, logarithmic derivatives of scale-$r$ connection events at the critical point should be of order $r^2$ and off-critical behaviour should become apparent at scale $r$ when $\p$ deviates from $\p_c$ by something of order $r^{-2}$. Taking inverses, this means that the correlation length should be of order $|\p-\p_c|^{-1/2}$ in high dimensions, a fact that is well-understood in the subcritical phase \cite{hutchcroft2023high,hara1990mean,CHS23,blanc2026markovian,liu2026crossover} but much less so in the supercritical phase.
Since the infinite percolation cluster should look like a large critical cluster up to scale $\xi(\p)$, it is reasonable to expect that the typical displacement of the $n$-step random walk obeys subdiffusive $n^{\varkappa_d}$ scaling until $n^{\varkappa_d}$ is of order $\xi(\p)$. On the other hand, at scales much larger than $\xi(\p)$ the random walk should behave like Brownian motion, with different increments approximately independent. These two heuristics together suggest that the variance of the random walk on the infinite cluster satisfies
\[
\operatorname{Var}(X_{k \xi(p)^{1/\varkappa_d}}) \approx k \operatorname{Var}(X_{\xi(p)^{1/\varkappa_d}}) \approx k \xi(\p)^2
\]
when $\p$ is slightly larger than $\p_c$ and $k$ is large (see Figure~\ref{fig:numerics}). If this formula is correct, it yields that
\[
\sigma^2(\p) = \lim_{k\to \infty} \frac{\operatorname{Var}(X_{k \xi(p)^{1/\varkappa_d}})}{k \xi(\p)^{1/\varkappa_d}} \approx \xi(\p)^{2-1/\varkappa_d},
\]
so that $\sigma^2(\p)$ should be of order $\xi(\p)^{-4} \approx |\p-\p_c|^{2}$ in the high-dimensional case, where $\varkappa_d = \frac{1}{6}$ and $\xi(\p) = \vert \p-\p_c \vert^{-1/2}$. The part of this heuristic concerning the random walk on scales smaller than $\xi(\p)$ is not too difficult to justify given known results and techniques about critical percolation in high dimensions; the harder part is showing that Brownian behavior already begins to kick in at this same scale. In fact this is not quite the perspective we take in this paper, and we instead apply a more potential-theoretic approach based on tools from quantitative homogenization theory as we discuss in more detail in Section~\ref{subsec.homogenization_intro}.

\subsection{Quantitative renormalization, the time constant and slab thresholds}
At the heart of the proof of Theorem~\ref{thm.main} is a sharp quantitative understanding of \emph{block renormalization} for slightly supercritical percolation in high dimensions. 
To explain what we mean by this, let us first recall the classical \emph{qualitative} theory of supercritical percolation in $d\geq 3$ dimensions, which is based on \emph{Grimmett--Marstrand theory} \cite{grimmett1990supercritical} (see also \cite{diskin2026supercritical} for a recent alternative approach to this theory). 
It is relatively easy to show that the geometry of the infinite cluster is close to that of the underlying lattice when $\p$ is very close to $1$ via direct counting arguments, and the goal is to prove that the infinite cluster at any supercritical $\p$ behaves like the infinite cluster with $\p$ close to $1$ at large scales.
Roughly speaking, Grimmett--Marstrand theory shows that for each $\p>\p_c$ there exists a length scale $U=U(\p)$ such that the box $[-U,U]^d$ has a high probability to be \emph{good}, meaning that there is a unique cluster crossing the annulus $[-U,U]^d\setminus [-U/2,U/2]^d$ and this cluster touches the box $[-U/8,U/8]^d$ and every face of the box $[-U,U]^d$. We call $U(\p)$ the \emph{local uniqueness scale}.
Since the unique crossing clusters in two good boxes offset by $U/4$ must merge within the union of the two boxes, this allows one to compare the original percolation model to a coarse-grained ``good box'' percolation model on a re-scaled lattice that has an effective value of $\p$ very close to $1$. (Since  the boxes used to define this coarse-grained model overlap, the coarse-grained model is not Bernoulli percolation but instead a \emph{finitely dependent} percolation model; this does not cause any issues thanks in part to a general result of Liggett, Stacey, and Schonmann relating finitely dependent and independent percolation models at high density \cite{liggett1997domination}.) 
Once this is done, one can (with work) show that the infinite cluster is similar to the ambient lattice in various senses on length scales much larger than the local uniqueness scale $U(\p)$  \cite{antal1996chemical,pete2008note,barlow2004random}.

Given this discussion, it is natural to ask how large we must take this length scale $U$ to be. While the proof of Grimmett and Marstrand can be made quantitative, the bounds it yields are extremely poor for $\p\approx \p_c$, with the best known bound in general dimension \cite{duminil2020upper} being of order $\exp( C|\p-\p_c|^{-2})$. As a first guess, it is natural to postulate that the local uniqueness scale $U(\p)$ is of the same order as the correlation length $\xi(\p)$.
In low dimensions it is expected that the correlation length does coincide (up-to-constants) with $U(\p)$, this being related to the conjecture that there are $O(1)$ critical clusters crossing the annulus $[-2r,2r]^d \setminus [-r,r]^d$ with high probability when $d<d_c$: this was understood rigorously in two dimensions by Kesten \cite{kesten1987scaling}, and conditional theorems in this direction are available in other dimensions \cite{borgs1999uniform}. On the other hand, it should \emph{not} be the case that $U(\p)$ and $\xi(\p)$ are of the same order as $\p\downarrow \p_c$ in high dimensions. Indeed, in high dimensions there are roughly $R^{d-6}$ critical clusters crossing the annulus $[-2R,2R]^d\setminus [-R,R]^d$ with high probability \cite{aizenman1997,CHS23}. On scales of order $\xi(\p)$ this large number of crossing clusters should only reduce by a constant factor, with heuristics\footnote{
	More precisely, it is expected that the merging of large clusters 
	from $\p_c$ to $\p_c+\eps$ is described approximately in high dimensions like an Erdos-R\'enyi random graph with average degree of order $\eps R^2$. If this heuristic is correct then at a large multiple of the correlation length $R=\eps^{-1/2}$ we should expect a constant fraction of crossing clusters to merge into a giant, but that there should not be a unique crossing cluster until the larger value $\sqrt{\eps^{-1}\log 1/\eps}$.
} suggesting that $U(\p)$ is of order $\xi(\p)\sqrt{\log \xi(\p)}$.


The core technical contribution behind Theorem~\ref{thm.main} is a sharp coarse-grained renormalization scheme that operates at the correlation-length scale \emph{without} requiring local uniqueness. Rather than merging all the large critical clusters in each box into a single crossing cluster, we retain many of these clusters and use sprinkling to connect them to clusters in neighbouring boxes, showing that critical crossing clusters in adjacent boxes merge in a ``highly supercritical'' way at large constant multiples of the correlation length when we increase the edge probability from $\p_c$ to $\p$. 
For the lower bound on the diffusivity, we must show that this highly supercritical structure above the correlation length supports flows of small energy. Roughly speaking, we do this by constructing order $(n/\xi(\p))^{d-1}\xi(\p)^{d-6}$ paths of length of order $\xi(\p)^2(n/\xi(\p))$ crossing a box of side-length $n\gg\xi(\p)$, and controlling the intersections between these paths so that distributing current equally among them gives the required energy bound.

We remark for experts that our proof does \emph{not} develop any sharp  quantitative theory of \emph{large finite clusters} in slightly supercritical percolation in high dimensions. We point this out since many natural approaches to the problems we study involve first understanding these large finite clusters (which remains open) before using this understanding to deduce results about the infinite cluster, an approach that was taken in e.g.\ the second author's work on nonamenable graphs \cite{hutchcroft2022slightly,hutchcroft2024slightly}. We do however make extensive use of the detailed understanding that has been developed of large clusters in \emph{critical} high-dimensional percolation; a detailed overview of relevant background results is given in Section~\ref{sec.pre}.

\medskip

\noindent \textbf{Intrinsic distances and the time constant.}
The sharp renormalization method we develop also allows us to resolve several other natural problems concerning the geometry of the slightly supercritical infinite cluster (and in fact these other problems are significantly less technical than the computation of the diffusivity). The first such result we highlight concerns \emph{intrinsic distances} (i.e., the graph distance on the cluster).
It follows from classical results of Antal and Pisztora \cite{antal1996chemical} and Garet and Marchand \cite{garet2007large} (relying on Kingman's subadditive ergodic theorem) that for each $\p>\p_c$ there exists a positive constant $\mu(\p)$ and a norm $\|\cdot\|_\p$ with $\|e_i\|_\p=1$ for every standard basis vector $e_i$ such that if $[x]$ denotes the closest point of the infinite cluster to each point $x\in \Zd$ then intrinsic distance $d_\mathrm{in}$ on the infinite cluster satisfies
\[
d_\mathrm{in}([0],[x]) = (1+o(1)) \mu(\p) \|x\|_\p
\]
almost surely as $x\to \infty$, with distances exponentially unlikely to deviate from this asymptotic estimate by a multiplicative factor (see also \cite{dembin2022upper}). The constant $\mu(\p)$ is known as the \emph{time constant}, and is known to be a Lipschitz function of $\p$ throughout the supercritical phase~\cite{cerf2022time}. We remark that the norm $\|\cdot\|_\p$ is conjectured to converge to the Euclidean norm $\|\cdot\|_2$ as $\p\downarrow \p_c$, meaning that the chemical distance becomes approximately invariant under rotations as $\p$ approaches the critical point (instead approaching the $L^1$ norm as $\p\uparrow 1$).

\medskip

Our next theorem computes the order of the time constant as $\p \downarrow \p_c$ in high-dimensional percolation.

\begin{theorem}
	\label{thm.time_constant}
	Consider nearest-neighbor Bernoulli bond percolation on $\Zd$. If  $d\geq 9$ and \eqref{eq:two_point_assumption} holds then there exist positive constants $C, c$ such that 
	\begin{align}\label{eq.main_time_constant}
		c(\p-\p_c)^{-1/2} \leq \mu(\p) \leq C(\p-\p_c)^{-1/2}.
	\end{align} 
	for every $\p_c<\p \leq 1$.
\end{theorem}

A corollary of Theorems~\ref{thm.main} and \ref{thm.time_constant} is that, in high dimensions, the $n$-step random walk on the infinite slightly supercritical cluster has typical intrinsic displacement of order $\sigma(\p) \mu(\p) n^{1/2} \asymp |\p-\p_c|^{1/2} n^{1/2}$ when $n$ is large.

Heuristically, Theorem~\ref{thm.time_constant} corresponds to the assertion that intrinsic geodesics look roughly Brownian on scales much smaller than $\xi(\p) \asymp |\p-\p_c|^{-1/2}$, so that intrinsic distances between the origin and the boundary of the box of radius $\xi(\p)$ are of order $\xi(\p)^2$, while intrinsic geodesics look roughly like straight lines on scales much larger than $\xi(\p)$, so that the intrinsic distance between the origin and the boundary of the box of radius $k\xi(\p)$ is of order $k\xi(\p)^2$.
(Note however that due to the failure of local uniqueness at scale $\xi(\p)$, we expect there to be many points of the infinite cluster with Euclidean distance of order $\xi(\p)$ but intrinsic distance much larger than $\xi(\p)^2$.) Note that the lower bound of Theorem~\ref{thm.time_constant} follows easily from the known exponential decay of the two-point function of slightly \emph{subcritical} high-dimensional percolation at scale $(\pc-\p)^{-1/2}$ \cite{hutchcroft2023high}, while the upper bound is more interesting as it requires us to understand the crossover to Euclidean behaviour above $\xi(\p)$.

\medskip

\noindent \textbf{Slab thresholds.} One of the most important consequences of the Grimmett--Marstrand renormalization theory is that, when $d\geq 3$, the critical probability of the \emph{slab} $\mathbb{Z}^{d-1} \times \{0,\ldots,n\}$ converges as $n\to \infty$ to $\p_c(\Zd)$ \cite{grimmett1990supercritical}. This fact underlies many of the important results concerning the qualitative theory of supercritical percolation, including the exponential decay of the finite-cluster two-point function throughout the supercritical phase \cite{chayes1987bernoulli}. Our methods yield a sharp analysis of the slab critical probability in high dimensions, showing that infinite clusters exist on slabs of thickness significantly larger than the correlation length and do not exist on slabs significantly thinner than the correlation length.

\begin{theorem}
	\label{thm.slab}
	Consider nearest-neighbor Bernoulli bond percolation on $\Zd$. If  $d\geq 9$ and \eqref{eq:two_point_assumption} holds then there exist positive constants $C, c$ such that 
	\begin{align}\label{eq.main_slab}
		c n^{-2} \leq \p_c(\mathbb{Z}^{d-1}\times \{0,\ldots,n\})-\p_c(\Zd) \leq \p_c(\mathbb{Z}^{2}\times \{0,\ldots,n\}^{d-2}) -\p_c(\Zd) \leq C n^{-2}
	\end{align} 
	for every $n\geq 1$.
\end{theorem}

As with our other results, The lower bound follows from one-arm estimates in slabs (Section~\ref{sec.log}), while the machinery developed in this paper yields the upper bound.

\subsection{The homogenization perspective}
\label{subsec.homogenization_intro}

We now discuss how our results and methods relate to \emph{homogenization theory}, where percolation has also been an important example since the work of Zhikov in \cite{zhikov1989}, and in particular to  \emph{quantitative homogenization theory} as developed in the last decade \cite{AKMbook, gloria2015quantification, armstrong2022elliptic}. Significant previous results regarding random walk on supercritical percolation clusters developed using these methods include the 
complete Liouville theorem 
\cite{armstrong2018elliptic} (extending  the earlier results of \cite{benjamini2015disorder}), the computation of first-order asymptotics for correctors \cite{dario2021corrector} and Green's functions \cite{dario2021quantitative}, 
a strong coupling version of the invariance principle \cite{gu2024coupling}, and the smoothness of the mapping $\p \mapsto \sigma(\p)$ \cite{gu2025diffusivity}. There is also very recent development in high-contrast homogenization \cite{armstrong2024renormalization, armstrong2026superdiffusive, armstrong2026homogenization}; see also the ICM proceedings \cite{armstrong2026coarse}.

Homogenization theory provides another approach to understanding invariance principles and the scaling of the diffusivity without explicit consideration of random walks. The key quantity is \emph{the homogenized conductivity}
\begin{align}\label{eq.abm_intro}
	\ab(U, \p) :=  \E_{\p} \Ll[\inf_{v \in \ell_{e_1} + C_0(U)}\frac{1}{\vert U \vert} \sum_{e \in \Ed(U)}  \a(e) \vert \nabla v(e)\vert^2\Rr],
\end{align}
defined for each subset $U \subseteq \Zd$, where $e_1$ is the unit vector of the first coordinate, $\ell_{e_1}(x) : = e_1 \cdot x$, and $C_0(U)$ denotes the set of functions from $U$ to $\R$ that vanish on the boundary of $U$. When $U \nearrow \Zd$, a simple subadditivity argument (Lemma~\ref{lem.subadditive}) shows that 
the quantity $\ab(U, \p)$ converges to the conductivity $\ab(\p)$. The conductivity is related to the diffusivity $\sigma(\p)$ by \emph{the Einstein relation} (Lemma~\ref{lem.basic})
\begin{align*}
	\frac{\sigma^2(\p)}{2} = \frac{\ab(\p)}{\theta(\p)},
\end{align*}
where $\theta(\p)=\mathbb{P}_\p(0\leftrightarrow\infty)$ denotes the percolation probability.
Since 
$\theta(\p)$ is known to obey the asymptotics
	$\theta(\p) \simeq (\p - \p_c)$ as $\p \searrow \p_c$
when $d>6$ and \eqref{eq:two_point_assumption} holds \cite{barsky1991percolation}, our mission is  to prove that 
\begin{align}\label{eq.a.conjecture}
	\ab(\p) \simeq (\p - \p_c)^3 \qquad \p \searrow \p_c.
\end{align}

To prove the upper bound $\sigma^2(\p) \lesssim (\p-\p_c)^2$, it suffices to prove that $\ab(U,\p) \lesssim (\p-\p_c)^3$ for a box $U$ with side-length of order the correlation length $\xi(\p)$, which we can attempt to do by inputting concrete test functions into the variational formula \eqref{eq.abm_intro}. 
As a naive first guess we might try taking the linear function
$\ell_{e_1}(x)=e_1 \cdot x$, but this  yields only the very weak upper bound $\ab(U, \p) \leq \p$. 
Taking the restriction of this linear function to points connected to the boundary of the box yields the still very weak bound 
	$\ab(U, \p) \lesssim \theta(\p)$,
which is not even sufficient to show that $\sigma(\p)\to 0$ as $\p \downarrow \p_c$. To get a better bound we must take the geometry of clusters into account.
When $d>6$ and \eqref{eq:two_point_assumption} holds, it is known that critical percolation in the box of side-length $r$ contains order $r^{d-6}$ \emph{spanning clusters}, which are tree-like and have intrinsic diameter of order $r^2$ \cite{aizenman1997,CHS23}.
%
Inspired by this picture, which we expect to still describe the geometry of the supercritical configuration within the correlation length, we can instead take $v \in \ell_{e_1} + C_0(U)$ to be a linear interpolation of the boundary conditions along the \emph{backbones} of these spanning clusters, i.e., the open paths that cross the box in the $e_1$ direction, taking the function to be piecewise constant off the backbone. Since these paths should have length of order $\xi(\p)^2$ and there should be order $\xi(\p)^{d-6}$ of them, 
we should have that
\begin{align*}
	&\frac{1}{|U|} \sum_{e\in \Ed(U)} \a(e)|\nabla v(e)|^2  \\
	&\hspace{1.7cm}\lesssim  \frac{1}{\vert U \vert}  \cdot \#\text{spanning clusters} \cdot \text{typical length of path} \cdot \left(\frac{\text{correlation length}}{\text{typical length of path}}\right)^2 \\
	&\hspace{1.7cm}\simeq  \xi(\p)^{-d} \cdot  \xi(\p)^{d-6} \cdot \xi(\p)^2 \cdot  (\xi(\p)/\xi(\p)^2)^{2} \simeq \xi(\p)^{-6} \simeq (\p - \p_c)^3
\end{align*}
as desired. In fact this calculation does not really work owing to a large contribution from points near the sides of the box parallel to the $e_1$ direction, which we ignored here, but this can be fixed by taking $U$ to be a cuboid $[0,\xi(\p)] \times [0,L]^{d-1}$ with $L \gg \xi(\p)$.
Our rigorous proof of the upper bound $\sigma^2(\p)$ does not really define an explicit test function but instead estimates the \emph{derivative} of $\ab(U,\p)$ in $\p$ as explained below.


\subsection{Outline of proof} 
We now give a brief overview of the proof of Theorems~\ref{thm.main}, ~\ref{thm.time_constant}, and ~\ref{thm.slab}. We follow the homogenization approach, so that our goal is to prove the estimate \eqref{eq.a.conjecture}. The correlation length is set as $\xi(\p)= |\p-\p_c|^{-1/2} $ in the calculation below. 

\medskip
\textbf{A. Differential inequality for upper bound.}
The proof of upper bound is analytic. Using the subadditivity of conductivity (see \cite[Chapter~1]{AKMbook} and Lemma~\ref{lem.subadditive}), we always have the upper bound
\begin{align}\label{eq.intro.step1}
	\ab(\p) \leq \ab(U,\p)
\end{align}
for every finite $U\subseteq \Zd$.
To avoid boundary layer effects, take $U$ to be a slab with one side-length of order $\xi(\p)$ and all other side-lengths $L\gg \xi(\p)$ (see \eqref{eq.def_cube_new}); the large $L$ limit allows us to ignore the boundary conditions on all faces not orthogonal to $e_1$, and we ignore these boundary conditions for the rest of this discussion.

Since the function $v$ attaining the minimum in \eqref{eq.abm_intro}	is harmonic on the restriction of the percolation configuration to $U$, edges that do not belong to the backbone in the sense that their endpoints are not connected by disjoint paths to the boundary of $U$ do not contribute to $\sum_{e\in \Ed(U)} \a(e)|\nabla v(e)|^2$. 
As such, $\ab(U,\tilde{\p})$ can be expressed equivalently as
\begin{align}\label{eq.intro.step2}
	\ab(U,\tilde{\p}) = \E_{\tilde{\p}}\Ll[\frac{1}{|U|}\sum_{e \in \Ed(U)} \a(e) \vert \nabla v_{U,e_1}(e) \vert^2 \mathbf{1}_{\{e_\bullet \longleftrightarrow \partial U\} \circ \{e_\circ \longleftrightarrow \partial U\}} \Rr],
\end{align}
where $v_{U,e_1}$ is the minimiser in the functional in \eqref{eq.abm_intro} and we write $e_\bullet$ and $e_\circ$ for the two endpoints of the edge $e$. 
A naive application of the maximum principle suggests that $|\nabla v_{U,e_1}(e)|\leq \xi(\p)$, and we can use the BK inequality and known bounds on one-arm probabilities (both in full spaces and half-spaces) \cite{KN11,CH20,van2025one} to easily obtain that
%
\begin{align}\label{eq.intro.aa_bound}
	\ab(U, \tilde{\p}) \lesssim \Ll\{\begin{array}{ll}
		\xi^{-2} (\log{\xi}), & \tilde{\p} \in (\p_c, \p];\\
		\xi^{-2} (\log{\xi}) \exp\Ll(- C^{-1} \xi \sqrt{\p_c - \tilde{\p}}\Rr), & \tilde{\p} \in (0, \p_c].
	\end{array}\Rr.
\end{align}
This is a bad bound for $\p$ close to $\p_c$, where we expect that $|\nabla v_{U,e_1}|$ is very small (and more precisely of order $\xi(\p)^{-1}$) on most backbone edges.

A sharper understanding comes from studying the derivative. Applying Russo's formula and using Cauchy-Schwarz inequality to bound the effect of flipping a pivotal edge on the energy of the harmonic function $v$,  we get the key differential inequality
\begin{align}\label{eq.intro.differential.ab}
	\ab^{(1)}(U,\p) \leq  \frac{1}{\p^{\frac{1}{2}}(1-\p)} \ab^{\frac{1}{2}}(U,\p) \bb^{\frac{1}{2}}(U,\p),
\end{align}
where the quantity
\begin{align}\label{eq.intro.step2b}
	\bb(U,\tilde{\p}) = \E_{\tilde{\p}}\Ll[\frac{1}{|U|}\sum_{e \in \Ed(U)} (1-\a(e)) \vert \nabla v_{U,e_1}(e) \vert^2 \mathbf{1}_{\{e_\bullet \longleftrightarrow \partial U\} \circ \{e_\circ \longleftrightarrow \partial U\}} \Rr]
\end{align}
also satisfies the bound \eqref{eq.intro.aa_bound}. This bound is much sharper for $\bb$ than for $\ab$ since potential differences across closed edges between distinct spanning clusters are typically large. Rewriting the differential  inequality \eqref{eq.intro.differential.ab} as a bound on the derivative of $\ab^{\frac{1}{2}}$, and integrating yields that
\begin{align}\label{eq.integral.ab}
	\ab^{\frac{1}{2}}(U, \p) \lesssim \ab^{\frac{1}{2}}(U, \p') + \int_{\p'}^{\p}  \bb^{\frac{1}{2}}(U,s) \, \d s.
\end{align}
Substituting the bound \eqref{eq.aa_bound} for both $\ab$ and $\bb$ into the right hand side and taking $\p'=\p_c-C \xi^{-2} (\log \xi)^2$ with a large constant $C$ then yields that
\begin{align}\label{eq.intro.step3}
	\ab(U,\p) \lesssim (\p - \p_c)^3  \vert \log (\p - \p_c) \vert,
\end{align}
which is sharp up to a logarithmic factor. 

Finally, to remove the logarithm correction in \eqref{eq.intro.step3}, one  needs a modified version of \eqref{eq.integral.ab}, so that $\ab$ and $\bb$ are more balanced. Precisely, $\bb$ is the term of perturbation, and we should restrict its source in the bulk of $U$. The manipulation then involves the analysis on inhomogeneous percolation.




\begin{remark}
	Experts in the theory of random walk on random fractals may have in mind a different approach to the upper bound on $\sigma^2(\p)$ based on the Markov type $2$ inequality for random walks in reversible random environments in $\Zd$ \cite{Ball1992, NPSS2006,lee2017conformal}. This inequality implies in particular that there exists a universal constant $C$ such that the random walk on the infinite cluster satisfies
	\[
	\mathbb{E}_\p\left[\|X_{k n}\|_2^2 \mid X_0 \leftrightarrow \infty\right] \leq C k \mathbb{E}_\p\left[\|X_{n}\|_2^2 \mid X_0 \leftrightarrow \infty\right]
	\]
	for every $k,n\geq 1$.
	As a consequence of this inequality, to prove $\sigma^2(\p) \lesssim (\p-\p_c)^2$ it suffices to prove that 
	\begin{align*}
		\mathbb{E}_\p\left[\|X_{n}\|_2^2 \mid X_0 \leftrightarrow \infty\right] \lesssim \xi^2, \text{ for } n=\xi^6,
	\end{align*}
	a statement that looks quite close to known estimates for random walks on large critical clusters. Unfortunately, the proof that the random walk on the high-dimensional IIC has typical $n$-step intrinsic displacement of order $n^{1/3}$ \cite{KN09} does \emph{not} show that this holds at the level of moments, and indeed it is known that the geometric properties of the cluster established there do not suffice to prove such moment estimates \cite[Example 2.6]{barlow2008random}. The alternative proof of Lee \cite{lee2017conformal} does establish moment estimates for the random walk on the high-dimensional IIC but also loses logarithmic factors in these estimates.
\end{remark}

\medskip
\textbf{B. Sharp renormalization at correlation length.}
The proof of lower bound is probabilistic. Recall that the effective resistance between two sets in the percolation configuration is defined by
\begin{align}\label{eq.intro.resistance}
	\Reff(A,Z) := \inf \Ll\{\sum_{e} \a^{-1}(e) \theta^2(e): \text{unit flow } \theta \text{ from } A \text{ to } Z\Rr\}.
\end{align}
Let $\Lambda_L$ be the cube of side-length $L$ and let $[\partial \Lambda_{L}]^{-}, [\partial \Lambda_{L}]^{+}$ be two opposite faces of this cube.
The conductivity admits the dual formula  (see \cite[Exercise~2.13]{LP17})
\begin{align}\label{eq.intro.dual}
	\ab(\Lambda_{L}, \p) \geq \frac{L^{2-d} \P_{\p}[E]^2}{\E_{\p}\Ll[\Reff([\partial \Lambda_{L}]^{-}, [\partial \Lambda_{L}]^{+}) \mathbf{1}_E\Rr]},
\end{align}
where $E$ is any event; we will take $E$ to be an appropriate ``good event'' that is not too unlikely and helps us bound the expectation in the denominator.

When $L=\xi(\p)$, then the $\p$-percolation in $\Lambda_{\xi}$ has a geometry close to that of the $\p_c$-percolation. For the latter, recall that in critical high-dimensional percolation a box of side-length $r$ contains roughly $r^{d-6}$ spanning clusters, which look like trees with an intrinsic diameter of order $r^2$ \cite{aizenman1997,CHS23}. Therefore, distributing the unit flow equally over the $\xi^{d-6}$ spanning clusters of $\Lambda_{\xi}$ as in Section~\ref{subsec.homogenization_intro} 
leads easily to the bound
\[
\E_{\p}\Ll[\Reff([\partial \Lambda_{\xi}]^{-}, [\partial \Lambda_{\xi}]^{+}) \mathbf{1}_E\Rr] \lesssim 
\xi^{d-6} \cdot \xi^{-2(d-6)} \cdot \xi^2 = \xi^{-d+8} 
\]
for an event of good probability $E$, and hence \eqref{eq.intro.dual} yields that
\begin{align}\label{eq.intro.lower.correlation}
	\ab(\Lambda_{\xi}, \p) \gtrsim 
	\xi(\p)^{-6}
	\simeq (\p-\p_c)^3.
\end{align}
However, this is not enough for Theorem~\ref{thm.main}, because $\ab(\Lambda_{L}, \p) \searrow \ab(\p)$ as $L \to \infty$ (see \eqref{eq.intro.step1}). One really needs to extend the result to the large scale cases. This is where our sharp quantitative renormalization enters the argument.


\begin{figure}[h!]
	\centering
	\includegraphics[width=\textwidth]{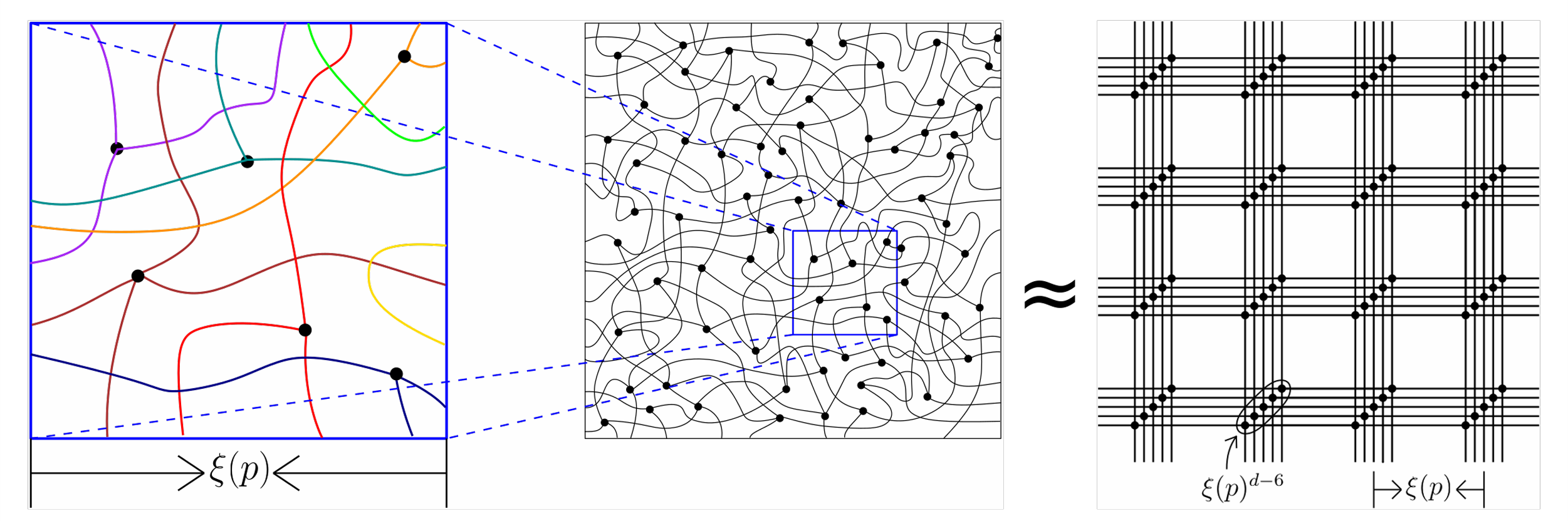}
	\caption{An illustration of the folded coarse lattice. Left: a box of side-length $\xi(\p)$ in critical percolation, with its spanning clusters drawn in different colours and one dot marked on each of them. Middle: a larger region tiled by such boxes, with the dots of all their spanning clusters; the blue square is the box shown on the left. Right: the folded coarse lattice. Every coarse site is replaced by $\xi(\p)^{d-6}$ dots (one such group is circled), and the dots of adjacent coarse sites are joined by parallel channels, each of length $\xi(\p)^2$.}
	\label{fig.folded}
\end{figure}

\smallskip

Now, we present the picture given by the renormalization argument in this paper. We zoom out but just treat the percolation at the length $\xi(\p)$, then the infinite cluster $\cCf$ should look like the following \emph{folded coarse lattice} (see Figure~\ref{fig.folded}): $\Zd$ is divided into boxes of side-length of order $\xi(\p)$, and each box is replaced by $\xi(\p)^{d-6}$ ``dots''. Each dot stands for a spanning cluster in the box, and dots in adjacent boxes have high probability, say $1-\epsilon_0$ with $\epsilon_0\ll1$, to be joined by paths of length $\xi(\p)^2$. The resulting coarse structure is highly supercritical, and hence looks like a folded lattice.


As mentioned in the introduction, this renormalization differs from the ones in previous work \cite{pisztora1996surface, armstrong2018elliptic}, since it does not require the local uniqueness. That is also why the outcome has $\xi(\p)^{d-6}$ folds. Meanwhile, there exist connections between different folds, but we will not use this property and do not quantify in which scale they typically occur. In a word, a renormalization without local uniqueness can prove Theorems~\ref{thm.main}, \ref{thm.time_constant} and~\ref{thm.slab}.

\begin{figure}[t!]
	\centering
	\includegraphics[width=\textwidth]{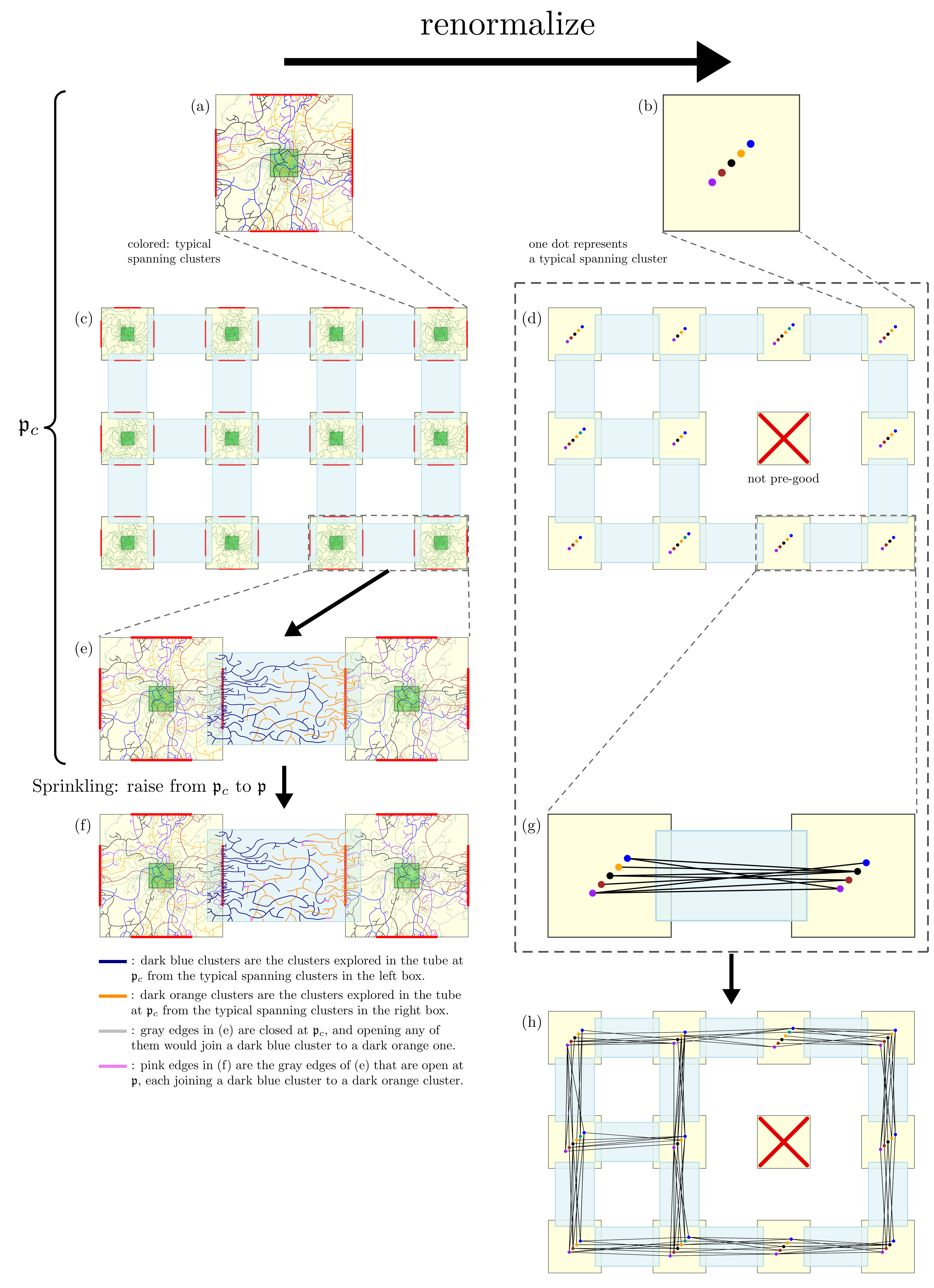}
	\caption{An illustration of the renormalization at correlation length.}\label{fig.ce1}
\end{figure}

We add more details about this renormalization argument. The folded lattice inherits the spirit of coarse-graining, thus it contains three elements: the coarse site, the coarse bond, and the percolation over the coarse structure. We explain them with the help of Figure~\ref{fig.ce1}, where the left column contains the original percolation configuration, and the right column is the coarse model extracted from it.

\emph{Coarse sites.} We place a lattice of disjoint boxes in $\Zd$, with the side-length $\ell_{0}(\p)$
\begin{align}\label{eq.intro.ell0}
	\ell_{0}(\p) = n_0 \xi(\p),
\end{align}
and $n_0$ is a large constant independent of $\p$. Two neighboring boxes also are separated by a distance of $\ell_{0}$; Figure~\ref{fig.ce1}(c). Every box contains a smaller one called its \emph{heart} (the green square in Figure~\ref{fig.ce1}(a) and~(c)). In every box, we explore at first all clusters of the $\p_c$-percolation issued from its heart, and keep the \emph{typical spanning clusters} among them, which are drawn in color in Figure~\ref{fig.ce1}(a). Roughly speaking, \emph{typical spanning clusters} are those with
\begin{itemize}
	\item typical volume $\xi(\p)^4$,
	\item typical intrinsic diameter $\xi(\p)^2$,
	\item typical order of intersection $\xi(\p)^2$ at each face of the boundary.
\end{itemize}
We then replace every typical spanning cluster by a dot as in Figure~\ref{fig.ce1}(b). A box is then called \emph{pre-good} if it contains $\xi(\p)^{d-6}$ typical spanning clusters. Following the same spirit as the seminal work \cite{aizenman1997} of Aizenman, we prove that
\begin{align*}
	\P_{\p}[\Lambda_{\ell_{0}} \text{ is pre-good}] \geq 1-\epsilon_0.
\end{align*}
Moreover, the pre-goodness is independent in different boxes, since this event only depends on the bonds inside. We refer to those pre-good boxes with their dots as \emph{the coarse sites}, while the others are removed as in Figure~\ref{fig.ce1}(d).

\emph{Coarse bonds.} Between every two neighboring boxes we draw a \emph{tube}, which is a rectangular region centered between the two boxes that bridges the gap and reaches a little into both boxes; see the light blue rectangles in Figure~\ref{fig.ce1}(c). Tubes of different pairs of neighbors are disjoint. We keep the tubes between two pre-good boxes, as in Figure~\ref{fig.ce1}(d), and each of them is a \emph{coarse bond}. A coarse bond records which typical spanning clusters of the two boxes are connected through the tube when the parameter is raised from $\p_c$ to $\p$.

\emph{Percolation of the coarse structure.} For the moment, we do not see how the typical spanning clusters in different boxes are connected, but they rely on the clusters in the tubes. This is thus the percolation of the coarse structure. Figures~\ref{fig.ce1}(e) and~(f) magnify one tube together with the two boxes it bridges. At $\p_c$, we explore the tube starting from the boundary points of the typical spanning clusters of both boxes (the dark blue one grows from the left, and the dark orange one grows from the right in Figure~\ref{fig.ce1}(e)). They are not connected as some bonds between them are closed at $\p_c$, but opening any one will make the two growing clusters connected. That is the key role played by the increment of parameter from $\p_c$ to $\p$. See the pink bonds in Figure~\ref{fig.ce1}(f), and its folded lattice version in Figure~\ref{fig.ce1}(g). Repeating this procedure for every typical spanning clusters and in every tube yields the coarse model in Figure~\ref{fig.ce1}(h). This is the folded coarse lattice described above.

We remark that, because the tubes are disjoint, the connections created in different tubes are independent, conditionally on the $\p_c$-configuration in the boxes. We also claim that, one typical spanning cluster $\cC \subset \Lambda_{\ell_{0}}(z)$ can connect to another one $\cC' \subset  \Lambda_{\ell_{0}}(z')$ in its neighbor boxes with high probability
\begin{align}\label{eq.intro.skl}
	\P_{\p}[\exists \cC' \subset  \Lambda_{\ell_{0}}(z'), \text{ s.t. } \cC \leftrightarrow \cC'   \vert \Lambda_{\ell_{0}}(z), \Lambda_{\ell_{0}}(z') \text{ are pre-good}] \geq 1-\epsilon_0.
\end{align}
The version proved in the paper is slightly relaxed, and the detailed justification will be discussed in Part \textbf{C}. Briefly, \eqref{eq.intro.skl} ensures the percolation of the coarse structure very supercritical, so the folded lattice is extended to the whole space. It contains $\xi(\p)^{d-6}$ folds as we can extract a supercritical percolation starting from each typical spanning cluster in the box. 


\medskip

Now, we  do some computation to justify Theorems~\ref{thm.main}, \ref{thm.time_constant} and~\ref{thm.slab} using the folded lattice obtained.

\emph{Intrinsic distances.} For $x,y \in \cCf$ with  an Euclidean distance $\vert x - y\vert = n\gg\xi(\p)$, one needs $ n/\xi(\p)$ coarse steps, and each has a graph distance $\xi(\p)^2$. Hence, the intrinsic distance $\din$ between $x,y$ is of order
\[
\din(x,y) \leq \frac{n}{\xi(\p)}\cdot\xi(\p)^2 \simeq n\,\xi(\p),
\]
that is
\begin{align*}
	\mu(\p) \lesssim (\p-\p_c)^{-1/2},
\end{align*}
in agreement with Theorem~\ref{thm.time_constant}.

\smallskip
 
\begin{figure}[h!]
	\centering
	\includegraphics[width=\textwidth]{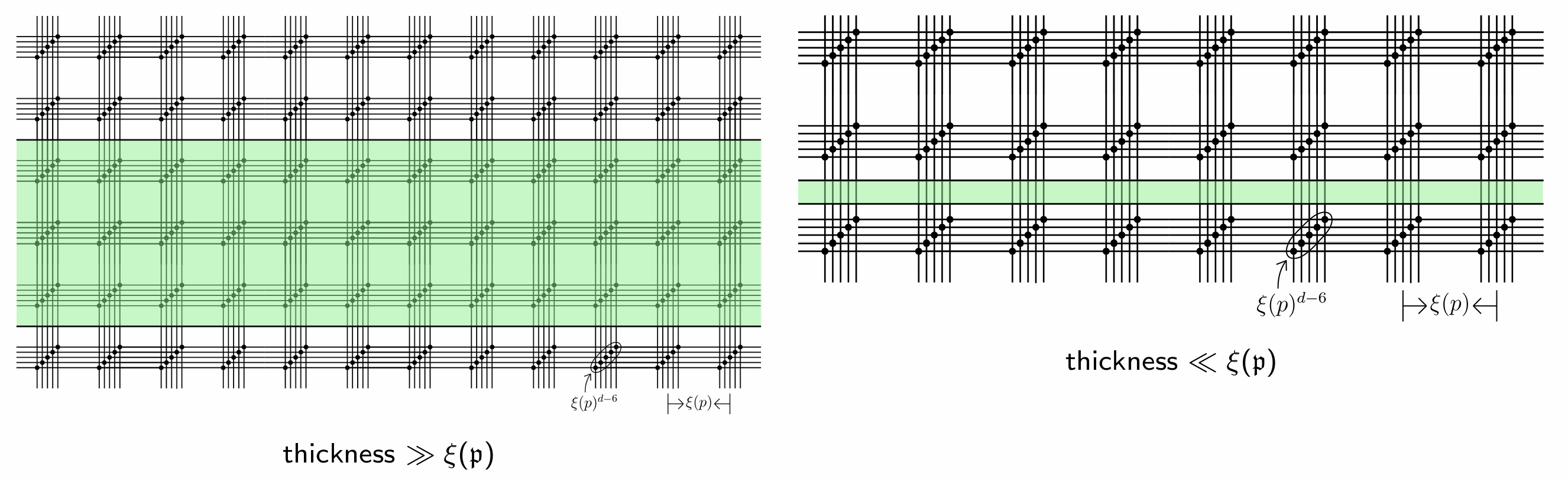}
	\caption{An illustration of the slab threshold.} \label{fig.slabcap}
\end{figure}

\emph{Slab thresholds.} When the thickness $n$ of a slab is a large multiple of $\xi(\p)$, the slab contains several layers of boxes of side-length $\xi(\p)$. Its intersection with the folded lattice contains an infinite component, thus the percolation at $\p$ in it is likely to contain an infinite cluster. This shows $\p$ is supercritical for the slab. A slab much thinner than $\xi(\p)$ does not contain a full layer of boxes of side-length $\xi(\p)$, and thus, its intersection with the folded lattice is likely to have no infinite components. Thus, the percolation at $\p$ is likely to have no infinite cluster. This shows $\p$ is subcritical for such a slab. Hence, $\p$ should be critical for some slab of thickness $ \xi(\p)$, in agreement with Theorem~\ref{thm.slab}. See Figure~\ref{fig.slabcap} for an illustration of this argument.

\smallskip

\emph{Resistance and diffusivity.} As the analysis in \eqref{eq.intro.dual} suggests, we need the following inequality to hold with positive probability,
\begin{equation}\label{eq.heuristic.resistance}
	\Reff\big([\partial\Lambda_L]^-,[\partial\Lambda_L]^+\big)\lesssim L^{2-d}\,\xi(\p)^6 \qquad \text{for } L\gg\xi(\p).
\end{equation}
Consider the folded coarse lattice in $\Lambda_L$, it contains roughly 
\begin{align*}
	N \simeq (L/\xi(\p))^{d-1} \cdot \xi(\p)^{d-6} \simeq L^{d-1}\xi(\p)^{-5}
\end{align*}
crossing paths $\gamma_1,\dots,\gamma_N$, whose corresponding paths in the supercritical infinite cluster are nearly disjoint, and each of them has length 
\begin{align*}
	\forall 1 \leq i \leq N, \quad	\vert \gamma_i \vert \simeq (L/\xi(\p)) \cdot \xi(\p)^2=L\,\xi(\p).
\end{align*}
Thus, we have $L^{d-1}\xi(\p)^{-5}$ nearly disjoint crossing paths of resistance $L \xi(\p)$. The total effective resistance $\Reff\big([\partial\Lambda_L]^-,[\partial\Lambda_L]^+\big)$ is then at most the resistance of such disjoint paths in parallel, which is 
\begin{equation}\label{eq.intro.resistance.upper}
	\Reff\big([\partial\Lambda_L]^-,[\partial\Lambda_L]^+\big) \lesssim \frac{L \xi(\p)}{L^{d-1}\xi(\p)^{-5}} \lesssim L^{2-d} \xi (\p)^6.
\end{equation}
We have to remark that, the heuristic of ``disjoint parallel paths" is oversimplified in some sense. Some more detailed discussion will be addressed in Part \textbf{D}.

We finish this part with several remarks. The computation above gives $\sigma^2(\p)\gtrsim(\p-\p_c)^2$, $\mu(\p)\lesssim(\p-\p_c)^{-1/2}$, and $\p_c(\ZZ^{d-1}\times\{0,\dots,n\})-\p_c\lesssim n^{-2}$. The other three inequalities do not use the renormalization. 
\begin{itemize}
	\item The bound $\sigma^2(\p)\lesssim(\p-\p_c)^2$ is the analytic upper bound as discussed in Part \textbf{A}.
	\item The bound $\mu(\p)\gtrsim(\p-\p_c)^{-1/2}$ follows from the exponential decay of the two-point function of slightly subcritical percolation \cite{hutchcroft2023high}.
	\item The bound $\p_c(\ZZ^{d-1}\times\{0,\dots,n\})-\p_c\gtrsim n^{-2}$ follows from an upper bound on the expected size of clusters in slabs. 
\end{itemize}
The remaining two parts give a panorama of some other key ingredients needed in the renormalization. Actually, the rigorous version will not be defined until Definition~\ref{def.spanning.z}, Definition~\ref{def.goodskl}, and Definition~\ref{def.expl}, which involve other properties. They are proposed after several rounds of forward-backward enhancement of the paper. We believe that, a preview of these technical aspects below helps understand the renormalization better.

\medskip
\textbf{C. Sprinkling property and regular points.} We now turn back to \eqref{eq.intro.skl}, which ensures the folded lattice supercritical. This relies on \emph{the sprinkling property} that the additional randomness from $\p_c$ to $\p$ creates new connections; recall the illustration in Figures~\ref{fig.ce1}(e) and~(f). Following the notations in Part $\mathbf{B}$, given $\cC, \cC'$ two spanning clusters living in two disjoint boxes $\Lambda_{\ell_{0}}(z)$ and $\Lambda_{\ell_{0}}(z')$ of distance $\ell_{0}$, we try to extend them from the boundary of each box and within the tube. Assume that $\cC, \cC'$ are disconnected at $\p_c$. This means that every exploration of open path between $\cC, \cC'$ is blocked by at least one closed bond. We denote by $\piv(\cC, \cC')$ the set of closed pivotal bonds, whose state determines the connection between $\cC, \cC'$. The expectation of the total number of pivotal bonds is  
\begin{align}\label{eq.intro.skl.proba1}
	\E_{\p_c}[\# \piv(\cC, \cC')] \simeq (\ell_{0}^{2})^2 \cdot \ell_{0}^d  \cdot (\ell_{0}^{-d+2})^2 \simeq \ell_{0}^{-d+8}.
\end{align}
Here we use the two-point function estimate and the fact that the spanning cluster typically has a boundary of size $\ell_{0}^2$. Now, let the parameter increase from $\p_c$ to $\p$, which may open some bonds in $\piv(\cC, \cC')$, and thus make $\cC, \cC'$ connected at $\p$. We denote by $\pivo(\cC, \cC')$ the set of open pivotal bonds at $\p$, which admits   
\begin{align}\label{eq.intro.skl.proba11}
	\E_{\p}[\# \pivo(\cC, \cC')] \simeq \ell_{0}^{-d+8} \cdot (\p - \p_c).
\end{align}
In high-dimension percolation, the probability in \eqref{eq.intro.skl.proba11} is very small, thus it is quite rare to find an open pivotal bond between  $\cC$ and $\cC'$. Nevertheless,  we can run $\cC'$ over all the $\ell_{0}^{d-6}$  spanning clusters in $\Lambda_{\ell_{0}}(z')$, so there are roughly
\begin{multline}\label{eq.intro.skl.proba2}
	\sum_{\cC' \text{ typical spanning cluster in } \Lambda_{\ell_{0}}(z')}\E_{\p}[\# \pivo(\cC, \cC')] \\
	\simeq \ell_{0}^{-d+8} \cdot (\p-\p_c) \cdot \ell_{0}^{d-6} \simeq \ell_{0}^{2} \cdot (\p-\p_c) \simeq n_0^2, 
\end{multline}
copies expected to be connected. The last equality utilizes \eqref{eq.intro.ell0} that 
\begin{align*}
	\ell_{0} = n_0 \xi(\p) = n_0 (\p-\p_c)^{-1/2}.
\end{align*}
By choosing the constant $n_0 \gg 1$ but independent of $\p$, thus we obtain a large expectation. The analysis in \eqref{eq.intro.skl.proba2}, together with the additional concentration estimate, then justifies \eqref{eq.intro.skl}.


Another technical input for \eqref{eq.intro.skl} is \emph{the regular point}, because \eqref{eq.intro.skl.proba1} assumes that every boundary point contributes a constant number of pivotal edges. Here, a boundary point is \emph{regular} if, roughly speaking, the explored critical clusters occupy sufficiently little space around it that a new connection can avoid them with good probability.  The consideration of such regular points has been a standard technique in high-dimensional percolation since the work of Kozma and Nachmias~\cite{KN11}. Standard BK arguments show that it is easy for new connections between regular boundary points of two different spanning clusters to avoid these revealed spanning clusters, letting us rigorously estimate the probability of the post-sprinkling connections between spanning clusters discussed above.

\medskip
\textbf{D. Intersection between crossing paths.} As mentioned at the end of Part \textbf{B}, the calculation in \eqref{eq.intro.resistance.upper} is oversimplified in the large scale. Actually, if all the crossing paths were disjoint, then it cannot be like supercritical percolation. Thus, the overlap between the paths needs to be handled. This issue is also well-known from the view point of \emph{random flow} of the resistance, which suggests proving
\begin{align}\label{eq.intro.random_flow}
	\Reff\big([\partial\Lambda_L]^-,[\partial\Lambda_L]^+\big) \leq \Eb\Ll[ \vert \gamma \cap \gamma'\vert \, \vert \a\Rr] \lesssim L^{2-d}\cdot \xi(\p)^{6}.
\end{align}
Here $\gamma,\gamma'$ are two independent unit random flows from $[\partial\Lambda_L]^-$ to $[\partial\Lambda_L]^+$.

We then analyze \eqref{eq.intro.random_flow} using the folded lattice in Part \textbf{B}. The flows $\gamma, \gamma'$ have their coarse version $\Gamma, \Gamma'$ as the projection from the folded lattice to the coarse sites. The intersection between $\Gamma$ and  $\Gamma'$ does not necessarily mean that between $\gamma$ and $\gamma'$, because the latter can live in different folds. The quantity in \eqref{eq.intro.random_flow} then can be decomposed into the product
\begin{align}\label{eq.intro.int}
	L^{2-d}\cdot \xi(\p)^{6} = (L/\xi(\p))^{1-d} \cdot \xi(\p)^{6-d} \cdot L \xi(\p),
\end{align}
and the three terms are interpreted as 
\begin{enumerate}
	\item $(L/\xi(\p))^{1-d}$ for the probability of intersection between $\Gamma,\Gamma'$
	\begin{align*}
		\PP[ \Gamma \cap \Gamma' \neq \emptyset] \simeq (L/\xi(\p))^{1-d},
	\end{align*}
	\item $\xi(\p)^{6-d}$ for the probability of intersection between $\gamma,\gamma'$, conditioned that $\Gamma,\Gamma'$ intersects
	\begin{align*}
		\PP[\gamma \cap \gamma' \neq \emptyset \vert \Gamma \cap \Gamma' \neq \emptyset] \simeq \xi(\p)^{6-d},
	\end{align*}
	\item $L \xi(\p)$ for intrinsic length of crossing path.
\end{enumerate}
Despite this nice interpretation, to make (2) valid is still technically difficult, since $\gamma$ determines the coarse version $\Gamma$ in the couple $(\gamma, \Gamma)$. A more desirable dependence is to preset $\Gamma$ without $\gamma$.

We thus turn to the idea in \cite{BPP98} by Benjamini, Pemantle, and Peres, where they proved the transience of random walk on supercritical percolation of $\Zd, d\geq 3,$ using the random flow. In our setting, their approach allows us to define a measure on the coarse version $\Gamma$ first with \emph{exponential intersection tail}, and then put weights on those crossing paths $\gamma$ of percolation. The presentation of renormalization is also partially modified by this argument. Besides, the estimate \eqref{eq.intro.int} can be seen as a refinement of \cite{BPP98}, in the sense that we obtain the factor $\xi(\p)^{6-d}$ with more effort.

\medskip

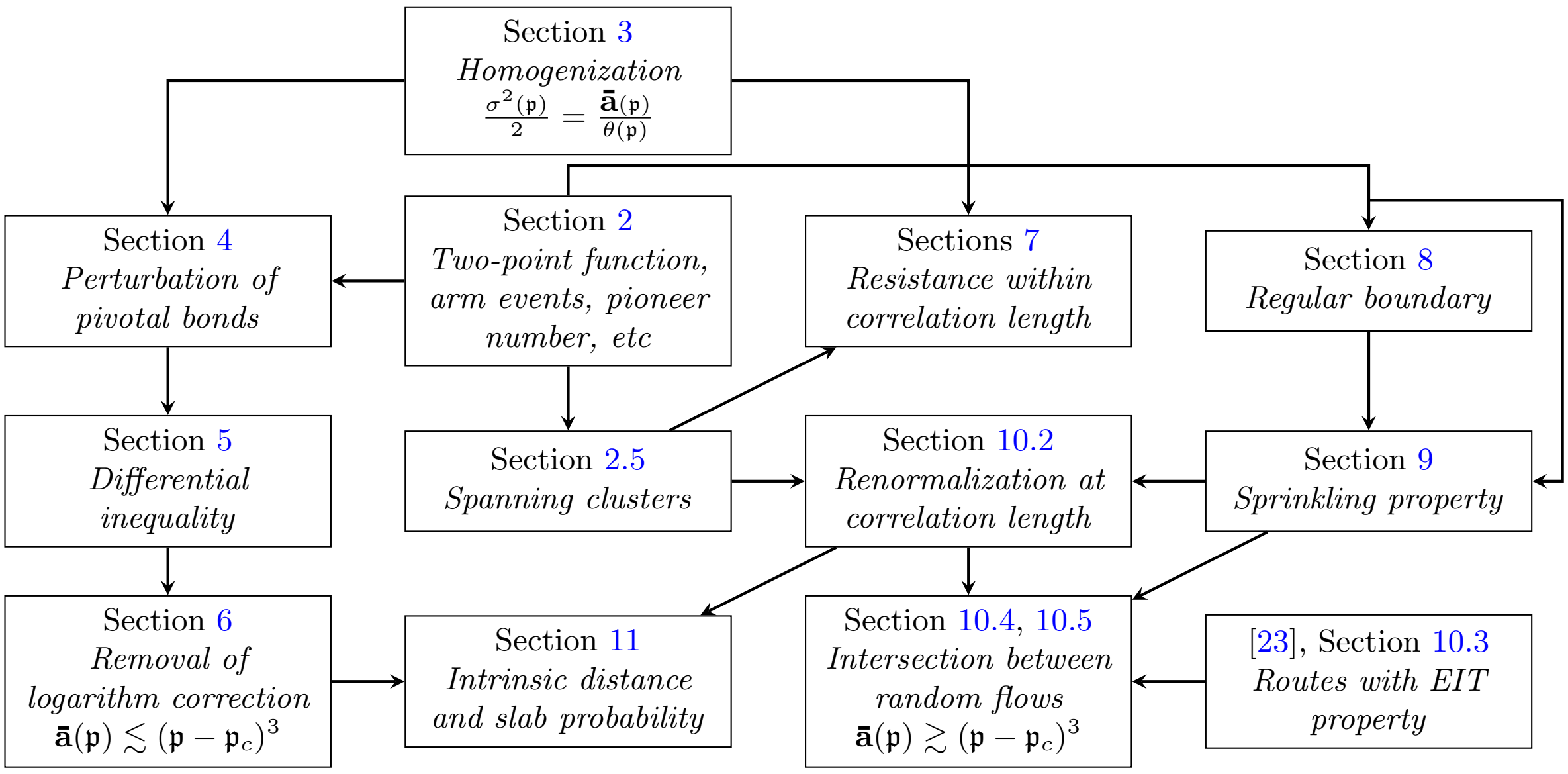
\begin{figure}
	\centering
	\begin{tikzpicture}[node distance=2cm]\label{fig.diagram}
		
		\node (Homo) [input] {Section~\ref{sec.homo} \\ \textit{Homogenization} $\frac{\sigma^2(\p)}{2} = \frac{\ab(\p)}{\theta(\p)}$};
		\node (Perco) [input, below of = Homo] {Section~\ref{sec.pre} \\ \textit{Two-point function, arm events, pioneer number, etc}};
		\node (Span) [input, below of = Perco] {Section~\ref{subsec.spanning_cluster} \\ \textit{Spanning clusters}};
		\node (Metric) [input, below of = Span] {Section~\ref{sec.chem}\\ \textit{Intrinsic distance and slab probability}};

		\node (Pivot) [input, left of = Perco, xshift = -2cm] {Section~\ref{sec.pivotal} \\ \textit{Perturbation of pivotal bonds}};
		\node (Diff) [input, below of = Pivot] {Section~\ref{sec.diff} \\ \textit{Differential inequality}};
		\node (Log) [input, below of = Diff] {Section~\ref{sec.log} \\ \textit{Removal of logarithm correction} $\ab(\p) \lesssim   (\p-\p_c)^3$};
		
		\node (Resistance) [input, right of = Perco, xshift = 2cm] {Sections~\ref{sec.lower.critical}\\ \textit{Resistance within correlation length}};
		\node (RG) [input, below of = Resistance] {Section~\ref{subsec.renormalization.critical} \\ \textit{Renormalization at correlation length}};
		\node (Flow) [input, below of = RG] {Section~\ref{subsec.intersection}, ~\ref{subsec.R.estimate} \\ \textit{Intersection between random flows} $\ab(\p) \gtrsim   (\p-\p_c)^3$};
		
		\node (Regular) [input, right of = Resistance, xshift = 2cm] {Section~\ref{sec.regular}\\ \textit{Regular boundary}};
		\node (Skl) [input, below of = Regular] {Section~\ref{sec.skl}\\ \textit{Sprinkling property}};
		\node (EIT) [input, below of = Skl] {\cite{BPP98}, Section~\ref{subsec.EIT}\\ \textit{Routes with EIT property}};
		
		\draw [arrow] (Perco) -- (Span);	
		\draw [arrow] (Perco) -- (Pivot);	
		\draw [arrow] (Perco.north) -| ($(Perco.north)  + (0cm, 0.3cm)$) -| ($(Regular.north)  + (0cm, 0.3cm)$) -| (Regular.north);	
		\draw [arrow] (Span) -- (RG);	
		\draw [arrow] (Span) -- (Resistance);

		\draw [arrow] (Homo)   -| (Resistance);	
		\draw [arrow] (Homo)   -| (Pivot);	
		
		\draw [arrow] (Pivot) -- (Diff);	
		\draw [arrow] (Diff) -- (Log);	
		\draw [arrow] (Log) -- (Metric);

		\draw [arrow] (RG) -- (Flow);	
		\draw [arrow] (RG) -- (Metric);	
		
		\draw [arrow] (Regular) -- (Skl);	
		\draw [arrow] ($(Regular.north)  + (0cm, 0.3cm)$) -| ($(Regular.north)  + (0.3cm, 0.3cm)$)  -|  ($(Skl.east)  + (0.3cm, 0cm)$) |- (Skl.east);	
		\draw [arrow] (Skl) -- (RG);
		\draw [arrow] (EIT) -- (Flow);		
		\draw [arrow] (Skl) -- (Flow);
	\end{tikzpicture}
	\caption{The outline of proof.}\label{fig.outline}
\end{figure}

The rest of the paper is organized as follows. The basic facts about high-dimensional percolation and homogenization will be recalled in Section~\ref{sec.pre} and Section~\ref{sec.homo}.  The main body of the paper is two-fold. Sections~\ref{sec.pivotal}-\ref{sec.log} are devoted to the upper bound, and Sections~\ref{sec.lower.critical}-\ref{sec.lower.large} treat the lower bound. Section~\ref{sec.chem} presents the proof of the other byproducts on the intrinsic distance and the slab threshold. See Figure~\ref{fig.outline} for the flowchart.


\subsection{Notations}\label{subsec.notation}
We introduce the notations and conventions used in the paper. We use $C(\cdots)$ to indicate the dependence of the constant, and it may vary from line to line in the equations. Without specific explanation, the shorthand notation $A \lesssim B$ means there exists a constant $C(d) < \infty$ such that $A \leq C \cdot B$. We also write $A \asymp B$ if $A \lesssim B$ and $B \lesssim A$.

For every $q \in \R^d$, the affine function is defined as $\ell_{q}(x):= q \cdot x$. 


For nearest neighbors $x,y$, we denote by $x \sim y$. The bonds in $\Ed$ are formed by the nearest neighbors. Given a bond $e \in \Ed$, we usually use $e = \{e_\bullet, e_\circ\}$ or $e = \{e_{-}, e_{+}\}$ as a convention to represent its two endpoints.

For every $U \subset \Zd$, we denote by $\partial U$ the boundary of the set 
\begin{align}\label{eq.Boundary}
	\partial U := \{x \in U \, \vert \, \exists y \in \Zd \setminus U, y \sim x\}.
\end{align}
We denote by $C_0(U)$ the space of functions supported on $U$
\begin{align}\label{eq.C0}
	C_0(U) := \{f:U \to \R \, \vert \, \supp(f) \subseteq U \setminus \partial U\}.
\end{align}

The edge set contained in $U$ is denoted by $\Ed(U)$
\begin{align}\label{eq.edge}
	\Ed(U) := \{\{x,y\} \, \vert \, x,y \in U, x \sim y\}.
\end{align}
We also denote by $\Ed^+(U)$ the edges issued from $U$
\begin{align}\label{eq.edge_plus}
	\Ed^+(U) := \{\{x,x+e_i\} \, \vert \, x\in U, i = 1,2, \cdots,d\}.
\end{align}
The set $\Ed^+(U)$ is sometimes more convenient than $\Ed(U)$ to manipulate. For example, once we partition $U \subset \Zd$ into a disjoint union, we also have the similar disjoint union for the edge sets
\begin{align}\label{eq.union_edge}
	U = \bigsqcup_{i=1}^n U_i \Longrightarrow \Ed^+(U) = \bigsqcup_{i=1}^n \Ed^+(U_i).
\end{align}
Notice that this identity is not necessarily true for the notation $\Ed(U)$. The edge set \eqref{eq.edge_plus} involves more vertices than $U$, so we define
\begin{align}\label{eq.U_plus}
	U^+ := \{x \in \Zd \, \vert \, x \in U \text{ or } x \in U+e_i, i = 1, 2, \cdots, d\}.
\end{align}	

We define the $L^2$-inner product over $U$ for every $f, g: U \to \R$
\begin{align}\label{eq.bracket1}
	\bracket{f,g}_U := \sum_{x \in U} f(x) g(x).
\end{align}	
This inner product also extends naturally to the vector field $F, G: U \to \R^n$
\begin{align}\label{eq.bracket2}
	\bracket{F, G}_U := \sum_{x \in U} F(x) \cdot G(x).
\end{align}
Especially, the finite difference operator in $\Zd$ follows
\begin{equation}\label{eq.defGradient}
	\begin{split}
		\forall x \in \Zd, \qquad (\nabla f)(x) &:= \Ll(\nabla_{e_1} f, \cdots, \nabla_{e_d} f\Rr)(x), \\
		\forall h \in \Zd, x \in \Zd, \qquad (\nabla_h f)(x) &:= f(x+h) - f(x).
	\end{split}
\end{equation}
Then we have 
\begin{align}\label{eq.id_gradient}
	\bracket{\nabla f, \nabla g}_U = \sum_{e \in \Ed^+(U)} \nabla f (e) \nabla g(e).
\end{align}

We denote by $\Lambda_L(z)$ the hypercube of side length $L$ 
\begin{align}\label{eq.cube}
	\Lambda_L(z) := z+\Ll(-\frac{L}{2}, \frac{L}{2}\Rr) \times \Ll(-\frac{L}{2}, \frac{L}{2}\Rr) \cdots \times \Ll(-\frac{L}{2}, \frac{L}{2}\Rr) \cap \Zd,
\end{align}
and set $\Lambda_L \equiv \Lambda_L(0)$ as the default cube centered at $0$. For convenience, we always assume $L$ to be odd number, so $\vert \Lambda_L \vert = L^d$, i.e. the cube contains $L^d$ vertices.

Let $e_1,e_2, \cdots, e_d$ stand for the canonical basis in $\Zd$. We use the notation $[\partial \Lambda_L]^+, [\partial \Lambda_L]^-$ to represent the top and bottom of the boundary
\begin{equation}\label{eq.def.TopBottom}
	\begin{split}
		[\partial \Lambda_L]^+ &:= \Ll\{x \in \partial \Lambda^+_L: x \cdot e_1 = \max_{y \in \Lambda^+_L} y \cdot e_1\Rr\}, \\
		[\partial \Lambda_L]^- &:= \Ll\{x \in \partial \Lambda^+_L: x \cdot e_1 = \min_{y \in \Lambda^+_L} y \cdot e_1\Rr\}, \\
		[\partial \Lambda_L]^{\pm} &:= [\partial \Lambda_L]^+ \cup [\partial \Lambda_L]^-.
	\end{split}
\end{equation}
For technical convenience in the paper, they are defined for parts of the boundary of $\Lambda^+_L$ rather than $\Lambda_L$.

The notation $\H$ stands for the half-space
\begin{align}\label{eq.defH}
	\H := \{x = (x_1, \cdots, x_d) \in \Zd: x_1 \geq 0\},
\end{align}
and $\H_n$ is defined as 
\begin{align}\label{eq.defHn}
	\H_n := \H - n e_1.
\end{align}

We denote by $\dex$ the \emph{extrinsic distance} / $\ell^{\infty}$ distance, while keeping $\din$ for \emph{intrinsic distance/chemical distance} on the percolation. We denote by $\ell_*(\p)$ the correlation length for all $\p > \p_c$ 
\begin{align}\label{eq.critical_length}
	\ell_*(\p) := \lfloor (\p - \p_c)^{-1/2}\rfloor.
\end{align}
We also write $\ell_*$ as a shorthand notation with the default parameter $\p$. 

Given $S \subset \Zd$, we denote by $\{A \stackrel{S}{\longleftrightarrow} B\}$ if there is an open path of $\a$ in $S$ connecting $A$ and $B$. We omit $S$ when $S = \Zd$. However, given cube $\Lambda$ and $A \subset \Lambda$, then we abuse the notation and understand $\{A \longleftrightarrow  [\partial \Lambda]^{\pm}\}$ as a shorthand notation for $\{A \stackrel{\Lambda}{\longleftrightarrow} [\partial \Lambda]^{\pm}\}$.


\section{Percolation background}\label{sec.pre}

We now review relevant results on high-dimensional percolation that we will use from the literature. All the results we quote here are known to hold whenever $d>6$ and \eqref{eq:two_point_assumption} holds, and we will write ``for high dimensional percolation'' as shorthand for this assumption. Many but not all of the results we cite here hold under the weaker assumption that the critical triangle diagram $\sum_{x,y\in \Zd} \P_{\p_c}(0\leftrightarrow x)\P_{\p_c}(x\leftrightarrow y) \P_{\p_c}(y\leftrightarrow 0)$ is finite \cite{aizenman1984tree,hutchcroft2022derivation}.

\subsection{Basic estimates}
Besides the two-point function, some of the most important quantities in percolation include the \emph{susceptibility}, \emph{volume tail}, and \emph{one-arm probability}.
The \emph{susceptibility} $\chi(\p)$ is defined as 
\begin{align}\label{eq.suspectivity}
	\chi(\p) := \E_{\p}[\vert \cC(0)\vert] = \sum_{x \in \Zd} \tau_{\p}(0,x),
\end{align}
where $\tau_\p(x,y)=\P_\p(x\leftrightarrow y)$ denotes the two-point function and $\cC(x)$ denotes the cluster of $x\in \Zd$. For high-dimensional percolation, the susceptibility satisfies the asymptotics
\begin{equation}
	\label{eq.sus}
	\chi(\p) \asymp \frac{1}{\p_c-\p}
\end{equation}
for every $\p<\p_c$ \cite{aizenman1984tree}. The volume tail is simply the probability that the origin belongs to a cluster of large volume, and satisfies
\begin{equation}
	\P_{\p_c}\Ll[\vert \cC(0)\vert \geq n\Rr] \asymp n^{-1/2} \label{eq.volume_decay}
\end{equation}
for every $n \geq 1$ for high-dimensional percolation \cite{barsky1991percolation}.
Finally, the \emph{one-arm probability} is the probability that the origin is connected to Euclidean distance $n$, and for high-dimensional percolation satisfies
\begin{equation}
	\P_{\p_c}\Ll[0 \longleftrightarrow \partial \Lambda_n\Rr] \asymp n^{-2} \label{eq.one_arm}
\end{equation}
for every $n\geq 1$ \cite{KN11} (see also the simpler recent proofs \cite{van2025one,asselah2025capacity,panis2026reversing}). See also \cite{blanc2026markovian,blanc2026super,hara2000scaling,cabezas2026range} for stronger versions of these theorems with precise first-order asymptotic estimates (with the very precise estimates of \cite{hara2000scaling} requiring stronger hypotheses).

\subsection{Crossover estimates}

In this paper we are particularly interested in the crossover between critical and off-critical behaviour. This is already understood for several quantities in the \emph{slightly subcritical} regime.
The most important such estimate for our purposes is the one-arm estimate
\begin{align}\label{eq.onearm_sub}
	\frac{C^{-1}}{n^2} \exp\Ll(- C n \sqrt{\p_c - \p}\Rr) \leq \P_{\p}\Ll[0 \longleftrightarrow \partial \Lambda_n\Rr] \leq \frac{C}{n^2} \exp\Ll(- C^{-1} n \sqrt{\p_c - \p}\Rr).
\end{align}
which holds for all $n\geq 1$ and $\p_c/2\leq \p \leq \p_c$ for high-dimensional percolation as proven in the independent parallel works \cite{hutchcroft2023high,CHS23}. We remark that it is also known similarly that
\begin{align}\label{eq.volume_sub}
	\frac{C^{-1}}{n^{1/2}} \exp\Ll(- C n (\p_c - \p)^2\Rr) \leq \P_{\p}\Ll[|\cC(0)|\geq n\Rr] \leq \frac{C}{n^{1/2}} \exp\Ll(- C^{-1} n (\p_c - \p)^2\Rr).
\end{align}
for every $n\geq 1$ and $\p_c/2\leq \p\leq \p_c$ \cite[Propositions 4.1 and 4.3]{hutchcroft2022slightly} and that
\begin{multline}\label{eq.two_point_sub}
	C^{-1} \|x-y\|^{-d+2} \exp\Ll(- C \|x-y\| \sqrt{\p_c - \p}\Rr) \leq \tau_\p(x,y) 
	\\\leq C \|x-y\|^{-d+2} \exp\Ll(- C^{-1} \|x-y\|\sqrt{\p_c - \p}\Rr).
\end{multline}
for every $\p_c/2\leq \p\leq \p_c$ and $x\neq y$, where the upper bound was proven in \cite{hutchcroft2023high} and the two-sided bound was recently proven by a different method in \cite{blanc2026markovian}. See also \cite{blanc2026markovian,liu2026crossover} for more precise asymptotic estimates and \cite{hutchcroft2022slightly} for similar slightly subcritical estimates on the tail of the intrinsic radius.

There is not yet a good quantitative theory of slightly supercritical high-dimensional percolation of comparable strength to the subcritical theory. What is well understood, however, is that quantities of interest do not change order within the correlation length. This is often best understood via the entropic method of Dewan and Muirhead \cite{DM23} (see also \cite{hutchcroft2022derivation}), as summarised in the following proposition.
\begin{lemma}[Proposition~{1.21} of \cite{DM23}]\label{lem.entropic}
	For every $d \geq 2$, there exists a constant $C_0(d)$ such that for every $\Lambda \subset \Zd$ containing $0$, every event $A$ measurable with respect to 
	\begin{align*}
		{\cC_{\Lambda}(0) := \{x \in \Zd: 0 \stackrel{\Lambda}{\longleftrightarrow} x \}},
	\end{align*}
	and every $\frac{\p_c}{2} \leq \p' < \p \leq \frac{\p_c+1}{2}$, we have that
	\begin{align}\label{eq.entropic}
		\Ll\vert \P_{\p}[A] - \P_{\p'}[A] \Rr\vert \leq C_0 \vert \p - \p'\vert \sqrt{\max \{\P_{\p}[A], \P_{\p'}[A]\} \E_{\p'}[|\cC_{\Lambda}(0)|]}.
	\end{align}
\end{lemma}
This estimate (which holds in every dimension) was used in \cite[Lemma~2.9]{gu2024random} to prove that the within-the-correlation-length slightly supercritical one-arm estimate
\begin{equation}\label{eq.onearm_super}
	\forall \p \in[\p_c, \p_c+n^{-2}], \qquad	\frac{c_0}{n^2} \le \P_{\p}[0 \longleftrightarrow \partial \Lambda_n] \le \frac{C_0}{n^2}
\end{equation}  
holds for high-dimensional percolation.

\subsection{Expected numbers of pioneers}
In this section we discuss estimates for the quantities
$\varphi_\p(S)$ and $\psi_\p(S)$, defined for each set $0\in S\subseteq \Zd$ by
\begin{align}\label{eq.varphi}
	\varphi_\p(S) := \p \sum_{\substack{x \in S \\ y \notin S \\ x \sim y}} \tau^S_{\p}(0,x), \qquad \text{ and } \qquad \psi_\p(S) := \sum_{x \in \partial S}\tau^S_{\p}(0,x).
\end{align}
Here $\tau^S_{\p}(0,x) := \P_{\p}[0 \xleftrightarrow{S} x]$ is the two-point function restricted in $S$.  We call $\psi_\p(S)$ the \emph{expected number of pioneers} with respect to $S$. The quantity $\varphi_\p(S)$ was introduced by Duminil-Copin and Tassion \cite{DCT16}, who noted that $\varphi_{\p_c}(S) \geq 1$ for every finite set $0\in S \subseteq \Zd$ and used the quantity $\varphi_\p(S)$
to prove the sharpness of phase transition. These quantities have since been recognised as useful ways to understand the correlation length in high-dimensional percolation \cite{DCP25a,van2025one,hutchcroft2022derivation}. We define
\emph{the critical length} to be
\begin{align}\label{eq.length_p}
	L(\p) := \min\{n \in \N: \varphi_\p(\Lambda_n) \leq 1/e^2\},
\end{align}
and the \emph{the sharp critical length} to be 
\begin{align}\label{eq.length_p_tilde}
	\tilde{L}(\p) := \min\{n \in \N: \exists S \subset \Lambda_n \text{ containing } 0 \text{ such that } \varphi_\p(S) \leq 1/e^2\}.
\end{align}
It follows from \cite[Proposition 2.3]{van2025one}  that, for high-dimensional percolation, these quantities satisfy
	\begin{align}\label{eq.length_bound}
		c(\p_c - \p)^{-1/2} \leq  \tilde{L}(\p), L(\p) \leq C(\p_c - \p)^{-1/2}.
	\end{align}
	for every $\p_c/2\leq \p< \p_c$.
The pioneer number is a very useful tool in the study of the connectivity probability.
Moreover, for high-dimensional percolation, it was shown in \cite[Theorem 2.3]{hutchcroft2023high} (see also \cite{van2025one,panis2026reversing}) that the expected number of pioneers in boxes and half-spaces satisfy
	\begin{align}
		c \exp\Ll(-C(\p_c - \p)^{\frac{1}{2}}n\Rr) \leq \p \psi_\p(\Lambda_n) \leq \varphi_\p(\Lambda_n) \leq C \exp\Ll(-c(\p_c - \p)^{\frac{1}{2}}n\Rr), \label{eq.expansion_11} \\
		c \exp\Ll(-C(\p_c - \p)^{\frac{1}{2}}n\Rr) \leq \p \psi_\p(\H_n) \leq \varphi_\p(\H_n) \leq C \exp\Ll(-c(\p_c - \p)^{\frac{1}{2}}n\Rr), \label{eq.expansion_12}
	\end{align}
	for every $\p_c/2\leq \p\leq \p_c$
	and in particular that
	\begin{align}\label{eq.expansion_2}
		c \leq \psi_{\p_c}(\H_n) \leq C  
	\end{align}
	for every $n\geq 1$. (The bounds for the box are not proven in \cite{hutchcroft2023high}, but the upper bound for the box follows from the upper bound for the half-space, while the lower bound for the box follows from the bound $\P_\p(0 \leftrightarrow \partial \Lambda_{2n}) \leq \psi_\p(\Lambda_n)\P_\p(0 \leftrightarrow \partial \Lambda_{n})$ and the estimate \eqref{eq.onearm_sub}.)
	The recent work of Panis and Schapira \cite{panis2026reversing} shows more generally that $\phi_{\p_c}(S)$ is bounded over all sets $0\in S \subseteq \Zd$ for high-dimensional percolation.

	\subsection{Estimates in half-space}
	We now state extensions of several of the above estimates to the half-space $\H$. Some of the relevant estimates do not (to our knowledge) appear in previous literature, so that we must supply proofs.
	The first such estimate concerns the half-space susceptibility $\chi^{\H}(\p) := \sum_{x \in \H} \tau^{\H}_{\p}(0,x)$.
	
	\begin{lemma}
		\label{lem:half_space_sus}
		In high-dimensional percolation, there exist finite positive constants $c, C$, such that
		\begin{align}\label{eq.chi_half_space}
			c(\p_c - \p)^{-1/2} \leq  \chi^{\H}(\p) \leq C(\p_c - \p)^{-1/2}
		\end{align}
		for every $\p_c/2 \leq \p < \p_c$.
	\end{lemma}
	
	\begin{proof}[Proof of Lemma~\ref{lem:half_space_sus}]
		It follows by symmetry (or the mass-transport  principle for $\Zd$) that $\chi^{\H}_\p = \sum_{n=0}^\infty \psi_\p(\H_n)$, and the claim follows easily from \eqref{eq.expansion_12}.    \end{proof}

	We next consider half-space arm exponents. It was proven in \cite[Proposition 3.1]{CH20} (see also \cite[Theorem 1.2]{van2025one}) that the probability that the origin is connected to the far side of the slab $\{0,\ldots,n\}\times \mathbb{Z}^{d-1}$ within the same slab in high-dimensional percolation satisfies
	\begin{align}\label{eq.onearm_half}
		\frac{c}{n^3} \leq \P_{\p_c}\Ll[0 \stackrel{\H}{\longleftrightarrow} \H_{-n}\Rr] \leq \frac{C}{n^3}
	\end{align}
	for every $n\geq 1$. 
	We remark that there is now also a good theory of the two-point function in half-spaces as developed in \cite{CH20,panis2026}.
	
	%
	%
	%
	%
	We will make heavy use of the following generalization of this estimate, which gives a sharp estimate of a similar probability when the origin is replaced by a  point at distance $k$ from the left boundary, allows us to take slightly supercritical $\p$ provided that $n$ is within the correlation length, and gives a sharp exponential improvement for slightly subcritical $\p$.
	\begin{prop}\label{prop.onearm_half_plus}
		In high-dimensional percolation, there exist two finite positive constants $C, c$ only depending on the dimension $d$, such that the following estimate holds for all ${0 \leq k \leq \frac{n}{2}}$
		\begin{equation}\label{eq.onearm_half_plus}
			\P_{\p}\Ll[k e_1 \stackrel{\H}{\longleftrightarrow} \H_{-n}\Rr] \leq \Ll\{\begin{array}{ll}
				C(k+1)n^{-3}, & \p \in (\p_c,\p_c + \frac{1}{n^2}], \\
				C(k+1)n^{-3} \exp\Ll(-c(\p_c - \p)^{\frac{1}{2}}n\Rr), & \p \in (0,\p_c].
			\end{array}\Rr.
		\end{equation}
	\end{prop}
	\begin{proof}
		The proof is divided into 3 steps. Throughout the proof, we define the shorthand notation
		\begin{align}\label{eq.theta_H_p}
			\theta^{\H}_n(\p) := \P_{\p}\Ll[0 \stackrel{\H}{\longleftrightarrow} \H_{-n}\Rr].
		\end{align}
		\textit{Step~1: $k=0$ for $\p \in (0,\p_c]$.} As mentioned, the case $\p = \p_c$ and $k = 0$ is already proved in  \cite[(4.15)]{van2025one}, and we treat the subcritical case here. The event $\{0 \stackrel{\H}{\longleftrightarrow} \H_{-n}\}$ implies that there is an open path passing $\partial \H_{\lfloor \frac{-n}{2}\rfloor}$. We denote by $u$ the vertex of the last visit at $\partial \H_{\lfloor \frac{-n}{2}\rfloor}$, and then decompose $\{0 \stackrel{\H}{\longleftrightarrow} \H_{-n}\}$ into the disjoint realization $\{0 \stackrel{\H}{\longleftrightarrow} u\} \circ \{u \stackrel{\H}{\longleftrightarrow} \H_{-n}\}$
		\begin{align*}
			\theta^{\H}_n(\p) &\leq \sum_{u \in \partial \H_{\lfloor \frac{-n}{2}\rfloor}} \P_{\p}\Ll[\{0 \stackrel{\H}{\longleftrightarrow} u\} \circ \{u \xleftrightarrow{\H_{\lfloor \frac{-n}{2}\rfloor}} \H_{-n}\}\Rr] \\
			&\leq \sum_{u \in \partial \H_{\lfloor \frac{-n}{2}\rfloor}} \P_{\p}\Ll[0 \stackrel{\H}{\longleftrightarrow} u\Rr] \P_{\p}\Ll[u \xleftrightarrow{\H_{\lfloor \frac{-n}{2}\rfloor}} \H_{-n}\Rr] \\
			&= \Ll(\sum_{u \in \partial \H_{\lfloor \frac{n}{2}\rfloor}} \P_{\p}\Ll[0 \stackrel{\H_{\lfloor \frac{n}{2}\rfloor}}{\longleftrightarrow} u\Rr]\Rr) \theta^{\H}_{\lfloor \frac{n}{2}\rfloor}(\p). 
		\end{align*}
		The BK inequality is applied in the second line and the stationarity is applied in the third line. Then we notice that 
		\begin{align*}
			\sum_{u \in \partial \H_{\lfloor \frac{n}{2}\rfloor}} \P_{\p}\Ll[0 \stackrel{\H_{\lfloor \frac{n}{2}\rfloor}}{\longleftrightarrow} u\Rr] = \sum_{u \in \partial \H_{\lfloor \frac{n}{2}\rfloor}} \tau_\p^{\H_{\lfloor \frac{n}{2}\rfloor}}(0,u) = \psi_\p(\H_{\lfloor \frac{n}{2}\rfloor}) \lesssim \exp\Ll(-c(\p_c - \p)^{\frac{1}{2}}n\Rr).
		\end{align*}
		The last estimate comes from \eqref{eq.expansion_12}. Therefore, we obtain that 
		\begin{align*}
			\theta^{\H}_n(\p)  \lesssim \exp\Ll(-c(\p_c - \p)^{\frac{1}{2}}n\Rr) \theta^{\H}_{\lfloor \frac{n}{2}\rfloor}(\p) \lesssim n^{-3}\exp\Ll(-c(\p_c - \p)^{\frac{1}{2}}n\Rr).
		\end{align*}
		Here we insert the bound following \eqref{eq.onearm_half} that 
		\begin{align*}
			\theta^{\H}_{\lfloor \frac{n}{2}\rfloor}(\p) \leq \theta^{\H}_{\lfloor \frac{n}{2}\rfloor}(\p_c) \lesssim n^{-3}.
		\end{align*}
		This concludes that 
		\begin{align}\label{eq.onearm_half_zero}
			\forall \p \in (0,\p_c], \qquad \P_{\p}\Ll[0 \stackrel{\H}{\longleftrightarrow} \H_{-n}\Rr] \leq C n^{-3} \exp\Ll(-c(\p_c - \p)^{\frac{1}{2}}n\Rr).
		\end{align}

		\smallskip
		\textit{Step~2: general $k$ for $\p \in (0,\p_c]$.}
		The proof is similar as that of \cite[(4.15)]{van2025one}, and we can prove actually that 
		\begin{align}\label{eq.one_arm_k_times}
			\forall \p \in (0,\p_c], \qquad \P_{\p}\Ll[k e_1 \stackrel{\H}{\longleftrightarrow} \H_{-n}\Rr] \lesssim (k+1) \P_{\p}\Ll[0 \stackrel{\H}{\longleftrightarrow} \H_{ -\lfloor n/2 \rfloor}\Rr].
		\end{align}
		Let $\gamma$ be the path connecting $k e_1$ to $\H_{-n}$, and let $u$ be \emph{the vertex} to attain the minimum in the direction $e_1$, i.e.
		\begin{align*}
			u_1 := u \cdot e_1 = \min\{x \cdot e_1: x \in \gamma\}.
		\end{align*}
		Then we can decompose $\gamma$ into two segments $\gamma_1$ and $\gamma_2$, which realize $\{k e_1 \stackrel{\H_{-u_1}}{\longleftrightarrow} u\}$ and $\{u \stackrel{\H_{- u_1}}{\longleftrightarrow} \H_{-n}\}$ with disjoint bonds; see Figure~\ref{fig.one_arm_half_space} for an illustration. Therefore, we have 
		\begin{align*}
			\P_{\p}\Ll[k e_1 \stackrel{\H}{\longleftrightarrow} \H_{-n}\Rr]
			&\leq \sum_{j=0}^k \sum_{u \in \partial \H_{-j}} \P_{\p}\Ll[\{k e_1 \stackrel{\H_{-j}}{\longleftrightarrow} u\} \circ \{u \stackrel{\H_{-j}}{\longleftrightarrow} \H_{-n}\}\Rr] \\
			&\leq \sum_{j=0}^k \sum_{u \in \partial \H_{-j}} \tau^{\H_{-j}}_{\p}(k e_1, u)\P_{\p}\Ll[u \stackrel{\H_{-j}}{\longleftrightarrow} \H_{-n}\Rr] \\
			&\leq \sum_{j=0}^k \psi_{\p}(\H_{k-j}) \P_{\p}\Ll[0 \stackrel{\H}{\longleftrightarrow} \H_{-(n-j)}\Rr] \\
			&\lesssim (k+1) \P_{\p}\Ll[0 \stackrel{\H}{\longleftrightarrow} \H_{ -\lfloor n/2 \rfloor}\Rr].
		\end{align*}
		Here the BK inequality is applied from the first line to the second line, then we use the stationarity from the second line to the third line. The knowledge about the expected number of pioneer \eqref{eq.expansion_2}
		\begin{align*}
			\forall \p \in (0,\p_c], \qquad \psi_{\p}(\H_{k-j}) \leq \psi_{\p_c}(\H_{k-j}) \leq C,
		\end{align*}
		is utilized in the last line. Thus, the proof of \eqref{eq.one_arm_k_times} is complete. Combining with \eqref{eq.onearm_half_zero}, we finish the proof of \eqref{eq.onearm_half_plus} for $\p \leq \p_c$.
		\begin{figure}[h!]
			\centering
			\includegraphics[width=0.6\textwidth]{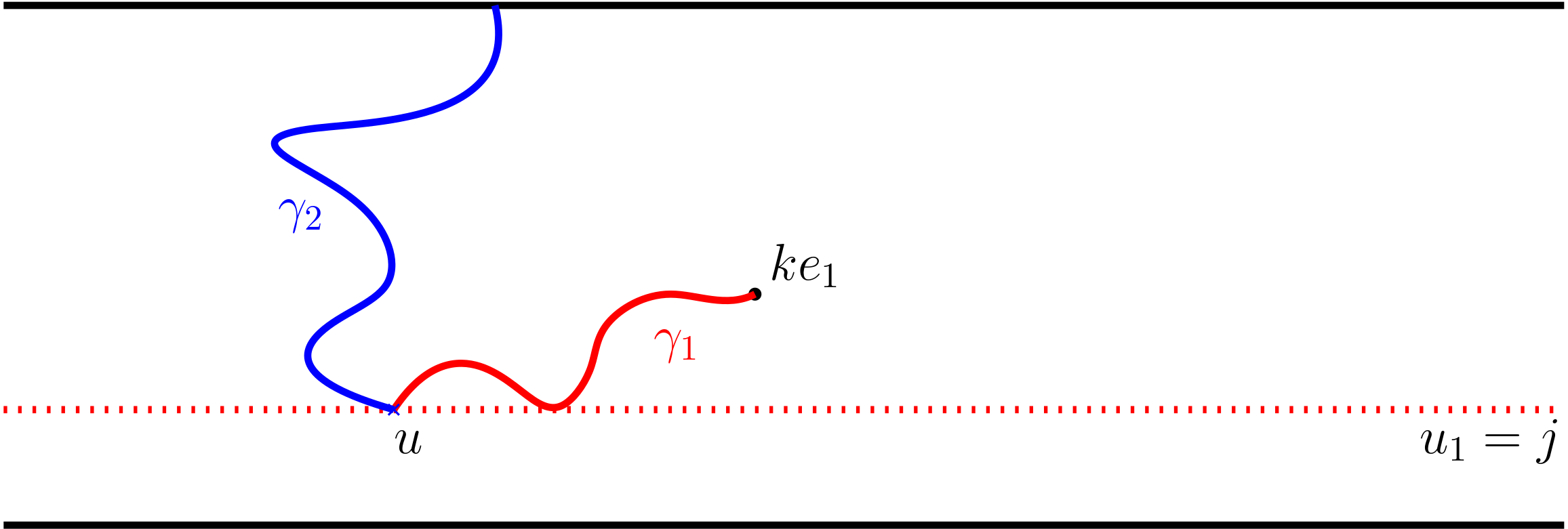}
			\caption{An illustration of the the event $k e_1 \stackrel{\H}{\longleftrightarrow} \H_{-n}$.}
			\label{fig.one_arm_half_space}
		\end{figure}
		
		\smallskip
		\textit{Step~3: $\p \in (\p_c,\p_c + \frac{1}{n^2}]$.}  The argument above does not apply directly to the supercritical regime, as the required estimates on $\psi_\p $ were not yet developed there. Alternatively, we can pass the estimate at $\p_c$ to the slightly supercritical phase via the entropic bound \eqref{eq.entropic}. We denote by $A$ the event 
		\begin{align*}
			A := \{k e_1 \stackrel{\H}{\longleftrightarrow} \H_{-n}\},
		\end{align*}
		and choose $\p'= \p_c - \frac{1}{n^2}$ and $\Lambda = \H$ in \eqref{eq.entropic}
		\begin{align}\label{eq.entropic_one_arm}
			\P_{\p}[A] - \P_{\p'}[A] \leq C_0 (\p-\p')\sqrt{\P_{\p}[A] \E_{\p'}|\cC^\H(k e_1)|},
		\end{align}
		where $\cC^\H(k e_1)$ denotes the cluster of $k e_1$ in the half-space $\mathbb{H}$.
		We refer to \eqref{eq.chi_half_space} for the upper bound of the susceptibility, and the choice of $\p'$ implies that 
		\begin{align*}
			0 \leq \p-\p' \leq 2 n^{-2} \qquad \text{ and } \qquad \E_{\p'}|\cC^\H(k e_1)| \leq C_2 (k+1) n,
		\end{align*}
		where we used the leftmost point decomposition as before to get that 
		\begin{align*}
			\E_{\p'}|\cC^\H(k e_1)| \leq \sum_{i=0}^k \psi_{\p'}(\H_i) \chi^\H(\p') \lesssim (k+1)(\p_c-\p')^{-1/2}.
		\end{align*}  
		Inserting these bounds back into \eqref{eq.entropic_one_arm}, together with $\P_{\p'}[A]$ from \eqref{eq.onearm_half_plus}, we have 
		\begin{align}\label{eq.one_arm_quadratic_formula}
			\P_{\p}[A] \leq C (k+1) n^{-3} + 2 C_0 C_2^{1/2} (k+1)^{1/2} n^{-3/2} \sqrt{\P_{\p}[A]}.
		\end{align}
		
		Notice that a quadratic inequality 
		\begin{align*}
			x^2 \leq C + B x \Longrightarrow x^2 - B x - C \leq 0 ,
		\end{align*}
		with $B,C >0$ will imply that 
		\begin{align}\label{eq.quadratic_bound}
			x \leq \frac{B + \sqrt{B^2 + 4 C}}{2} \leq B + \sqrt{C}.
		\end{align}
		We put this estimate to \eqref{eq.one_arm_quadratic_formula} and get 
		\begin{align*}
			\sqrt{\P_{\p}[A]} \leq (2 C_0 C_2^{1/2} +C^{1/2}) (k+1)^{1/2} n^{-3/2}. 
		\end{align*}
		We take a square and conclude \eqref{eq.onearm_half_plus} in the supercritical regime.
	\end{proof}

	\subsection{Spanning clusters}\label{subsec.spanning_cluster}
	Our analysis of effective resistances relies strongly on understanding the geometry of  \emph{spanning clusters} in critical high-dimensional percolation, a topic first explored by Aizenman  \cite{aizenman1997} and developed further in e.g.\ \cite{CHS23,asselah2024intersection}. 
	Roughly speaking, it is known that, in high-dimensional percolation, a cube of side-length $n$ will contain order $n^{d-6}$ spanning clusters and that these clusters will typically have volume of order $n^4$ and intrinsic diameter of order $n^2$.
	Below we state several estimates in strong quantitative form primarily because we find these forms easier to read than $\eps$-$\delta$ statements, but note that we do not actually need the strong quantitative estimates presented here. The main thing we will need beyond standard results of the literature are lower bounds on the probability that a vertex is connected to \emph{every face} of the box, with order $n^2$ points per face (Lemma~\ref{lem.arm.2dspanning}).
	We write $\cC_n(x)$ to denote the cluster of $x$ in the box $\Lambda_n$.
	

	
	
	In the course of Kozma and Nachmias' proof of the one-arm estimate \eqref{eq.one_arm}, it is shown that it is unlikely for high-dimensional percolation that a critical cluster has Euclidean diameter at least $n$ but volume much smaller than $n^4$. A strong quantitative version of this estimate \cite[Theorem 4]{CHS23} states for high-dimensional percolation that for each $\alpha>3d/2$ there exist positive constants $c,C$ such that
	\begin{align}\label{eq.cluster_size}
		\P_{\p_c}\Ll[\vert \cC_n(0) \vert \leq \lambda n^4 \Big\vert  0 \longleftrightarrow \partial \Lambda_n\Rr] \left\{\begin{array}{ll}
			\leq \exp(-c\lambda^{-\frac{1}{3}}), \qquad & \lambda > (\log n)^\alpha n^{-3};\\
			\geq \exp(-C\lambda^{-\frac{1}{3}}), \qquad & \lambda > C n^{-3}
		\end{array} \right.
	\end{align}
	for every $n\geq 1$. See also \cite{asselah2024intersection} for estimates showing conversely that is is very unlikely for a cluster to have volume at least $n^4$ and Euclidean diameter much smaller than $n$ (a weak form of this following immediately from the estimates $\E_{\p_c} |\cC(0)| \asymp n^2$ and $\P_{\p_c}(0\leftrightarrow \partial \Lambda_n)\asymp n^{-2}$). 
	
	
	It is also known that, in high-dimensional percolation, a cluster is very unlikely to reach Euclidean distance $n$ without the path having length of order $n^2$. A weak version of this follows from the works \cite{KN09,KN11}, while 
		\cite[Theorems~2 and 3]{CHS23} establish the strong quantitative estimates
		%
		%
			\begin{align}\label{eq.chemical_1}
				\P_{\p_c}\Ll[\din(0, \partial \Lambda_n) \leq \lambda n^2 \Big\vert  0 \longleftrightarrow \partial \Lambda_n\Rr] \leq \exp(-c \lambda^{-1}),
			\end{align}
			and
			\begin{align}\label{eq.chemical_2}
				\P_{\p_c}\Ll[\din(0, \partial \Lambda_n) \geq \lambda n^2 \Big\vert  0 \longleftrightarrow \partial \Lambda_n\Rr] \leq \exp(-c \lambda).
			\end{align}
			for every $n\geq 1$ and $\lambda>0$. We also have similarly by Equation (7) of \cite{CHS23} 
			that \begin{align}\label{eq.chemical_point_point}
				\forall x \in \Lambda_n, \qquad	\P_{\p_c}\Ll[\din(0, x) \geq \lambda \vert x \vert^2 \Big\vert  0 \xleftrightarrow{\Lambda_{2n}} x\Rr] \leq \exp(-c \lambda).
			\end{align}
			for every $n\geq 1$, $\lambda>0$, and $x\in \Lambda_n$.

		Finally we come to theorems counting the number of spanning clusters. 
		Sharpening the results of \cite{aizenman1997}, it is shown in \cite{CHS23} that 
		if $\mathscr{S}_n$ denotes the set of clusters in $\Lambda_n$ touching both $[\partial \Lambda_n]^+$ and $[\partial \Lambda_n]^-$ then
		$n^{-d+6}\mathscr{S}_n$ is tight in $(0,\infty)$, meaning that for each $\eps>0$ there exists $C$ such that 
		\begin{align}\label{eq.num_spanning}
			\P_{\p_c}[C^{-1} n^{d-6}\leq \vert \mathscr{S}_n \vert \leq C n^{d-6}] \geq 1-\eps.
		\end{align}
		Moreover, the same statement holds for the set 
		$\mathscr{S}_n'=\{\cC \in \mathscr{S}_n: c' n^4\leq |\cC| \leq C' n^4\}$ for some positive constants $c'$ and $C'$.
		
		
		\medskip
		
		\noindent \textbf{Simultaneous spanning clusters.} In our applications it will be important to use a different notion of spanning cluster than is typically used in the literature. For every $1 \leq i \leq d$, we denote by  $\partial_i \Lambda_n$ and $\partial_{-i} \Lambda_n$ respectively the faces that reaches the maximum and the minimum in the direction $e_i$. That is,
		\begin{align}\label{eq.def.bdei}
			\begin{split}
				\partial_i \Lambda_n &:= \Ll\{x \in \partial \Lambda_n: x \cdot e_i = \max_{y \in \Lambda_n} y \cdot e_i\Rr\}, \\
				\partial_{-i} \Lambda_n &:= \Ll\{x \in \partial \Lambda_n: x \cdot e_i = \min_{y \in \Lambda_n} y \cdot e_i\Rr\}.
			\end{split}
		\end{align}
		
		\begin{definition}\label{def.spanning}
			Given a cube $\Lambda_n$, for every different indexes $i,j \in \{\pm 1, \pm 2, \cdots, \pm d\}$, we call a cluster $\cC$ an $(i,j)$-spanning cluster if $\cC \xleftrightarrow{\Lambda_n} \partial_i \Lambda_n$ and $\cC \xleftrightarrow{\Lambda_n}  \partial_j \Lambda_n$. $\cC$ is called a $2d$-spanning cluster (or simultaneous cluster) if it connects all the $2d$ faces within $\Lambda_n$.
		\end{definition}
		
		In the literature, the term ``spanning cluster'' usually refers to what we call a $(+i,-i)$-spanning cluster. In his work \cite{aizenman1997} estimating the number of (what we call)  $(+i,-i)$-spanning clusters, Aizenman writes  \emph{``It can, however, be easily seen that none of the results will change (except for constants) if that notion is modified by declaring as spanning only those clusters which reach each of the $2d$ boundary faces.''} 
		We give a short proof of this estimate assuming the scaling limit theorem of \cite{blanc2026super}, which is really greatly excess to requirements. Importantly for our applications, we also show that
		the number of spanning clusters does not change order if we require them to have intersections of order $n^2$ with each boundary face of the box.
		
		For each $i\in\{\pm1,\ldots,\pm d\}$ we define $\partial_i^\circ \Lambda_n$ to be the set of points in $\partial_i \Lambda_n$ whose coordinates other than the $i$th coordinate belong to $(-n/4,n/4)$, so that $\partial_i^\circ\Lambda_n$ is at the center of the face  $\partial_i \Lambda_n$ and is well-separated from the boundaries of this face.

		\begin{lemma}\label{lem.arm.2dspanning}
			Consider high-dimensional critical percolation. There exist finite positive constants $c_{\text{\ref{eq.arm.2dspanning}}}, C_{\text{\ref{eq.arm.2dspanning}}}$ only depending on $d$, such that
			\begin{align}\label{eq.arm.2dspanning}
				c_{\text{\ref{eq.arm.2dspanning}}}n^{-2} \leq \P_{\p_c}\Ll[\cC_n(x) \text{ is } 2d\text{-spanning}\Rr] \leq C_{\text{\ref{eq.arm.2dspanning}}}n^{-2}.
			\end{align}
			and moreover
			\begin{multline}\label{eq.num_bd}
				\P_{\p_c}\Big[\forall i\in\{\pm1,\ldots,\pm d\}:\ |\cC_n(x)\cap\partial_i^\circ\Lambda_n|\ge c_{\text{\ref{eq.arm.2dspanning}}}n^2 \\
				\text{ and } |\cC_n(x)\cap\partial_i\Lambda_n|\le C_{\text{\ref{eq.arm.2dspanning}}}n^2\Big]\ \geq c_{\text{\ref{eq.arm.2dspanning}}} n^{-2}.
			\end{multline}
			for every $n\geq 1$ and $x\in \Lambda_{n/16}$.
		\end{lemma}

		\begin{remark}
			This proof will invoke the super-Brownian scaling limit theorem of \cite{blanc2026super}. In that paper, this scaling theorem is proven under a stronger assumption than our usual high-dimensional assumption that $d>6$ and that \eqref{eq:two_point_assumption} holds. More specifically, the results of \cite{blanc2026super} also require a ``two-blob'' condition (B), stating roughly that if we condition two distant points to be connected then, in the large-separation limit, we see two \emph{independent} copies of the IIC locally around the two points. However, it follows from the results of \cite{chatterjee2025robust} that this condition (B) holds whenever $d>6$ and \eqref{eq:two_point_assumption} holds, so that we do not need to invoke (B) as an additional hypothesis.
		\end{remark}
		
		\begin{proof}
			The upper bound follows from \eqref{eq.one_arm}; we focus on the lower bound.
			For each $i\in \{\pm 1,\ldots,\pm d\}$ we write $e_i$ to denote the basis vector $e_i$ if $i>0$ and $e_i=-e_{-i}$ if $i$ is negative and, given a vector $x$, write $x_i=x \cdot e_i$. Fix $x\in \Lambda_{n/16}$. For each $i\in \{\pm 1,\ldots,\pm d\}$ we define $E_{i}(x)$ to be the event that
			$x$ is connected to a point in the box $2ne_i+[-n/8,n/8]^d$ by a path contained in the cuboid 
			$\{z : |z_j| \leq n/8$ for $j\notin \{i,-i\}$, $-3n/8 \leq z_i \leq 3 n\}$.
			We claim that
			\begin{equation}
				\label{eq:2d_spanning_cross}
				\P_{\p_c}\left( \bigcap E_i(x) \right) \gtrsim n^{-2};
			\end{equation}
			this estimate clearly implies the lower bound of \eqref{eq.arm.2dspanning} and can also be used to prove \eqref{eq.num_bd} as we explain below. To keep notation light we prove this estimate only for $x=0$, the proof for general $x\in \Lambda_{n/16}$ being similar.
			
			We give a short proof of the $x=0$ case of \eqref{eq:2d_spanning_cross} from the scaling limit results of \cite{blanc2026super} but note that more direct proofs are also available.
			Let $d_\mathrm{chem}$ denote the chemical distance on $\cC(0)$ (i.e.\ graph distance on the cluster), let $\mu_{\cC(0)}$ denote the counting measure on the cluster, and let $\Phi_{\cC(0)}$ denote the inclusion map from $\cC(0)$ to $\R^d$.
			It is proven in \cite{blanc2026super} (see also \cite{chatterjee2025limiting,cabezas2026range,chatterjee2026convergence}) that, when $d>6$ and \eqref{eq:two_point_assumption} holds, there exists a positive constants $c_1,c_2,c_3$ such that the conditional law of the rescaled cluster
			\[\left(\cC(0), \frac{d_\mathrm{chem}}{c_1 n^2}, \frac{\mu_{\cC(0)}}{c_2 n^4}, \frac{\Phi_{\cC(0)}}{c_3 n}, 0\right)\]
			given $n^4 \leq |\cC(0)| \leq 2n^4$ converges weakly to the law of $(\mathscr{T},\nu,d,\Phi,o)$ where $(\mathscr{T},\nu,d,o)$ is a continuum random tree of random mass between $1$ and $2$ and $\Phi:\mathscr{T}\to \R^d$ is the Brownian embedding of $\mathscr{T}$ into $\R^d$. For our purposes it is more convenient to instead scale the embedding by $n$, and obtain a limit where the Brownian embedding has some constant diffusivity, as we will do from now on. This convergence is proven with respect to a sufficiently strong topology (namely the Gromov--Hausdorff--Prokhorov--function topology) that to prove the lemma it suffices to prove that an appropriate spanning event holds with positive probability for this limit object. (Here we are also using that the event $n^4 \leq |\cC(0)| \leq 2n^4$ has probability of order $n^{-2}$ as follows from \cite[Theorem 1.6]{blanc2026super}.) 
			More specifically, if for each $i\in \{\pm1,\ldots,\pm d\}$ we define $\mathcal{E}_{i}$ to be the event that there exists $y\in \mathscr{T}$ such that $\| \Phi(y) - 2  e_i\|_\infty \leq 1/32$ and 
			$\Phi(w) \in \{(z_1,\ldots,z_d) : |z_j| \leq 1/32$ for $j\notin \{i,-i\}$, $-1/4 \leq z_i \leq 5/2 \}$ for every $w$ on the geodesic between $o$ and $y$ in $\mathscr{T}$, then it suffices to prove that the intersection $\bigcap_{i} \mathcal{E}_i$ occurs with positive probability. (Note that we have reduced all the constants in the definition of the cuboid to make applications of the portmanteau theorem unproblematic.)
			
			We will recall only those properties of this limiting object that are needed for this proof, referring the reader to e.g.\ \cite{le2012spatial} for further background.
			First, it follows from a computation of Aldous \cite{aldous1993continuum} as explained in \cite[Theorem I.5.24]{hutchcroft2025critical} that if, after conditioning on $(\mathscr{T},\nu,d,o)$, we pick $k$ points $X_1,\ldots,X_k$ of $\mathscr{T}$ independently at random from the measure $\nu$ then the law of the smallest subtree containing the points $o,X_1,\ldots,X_k$ is mutually absolutely continuous with respect to the random metric tree constructed by first picking a combinatorial tree with unlabelled internal vertices of degree $3$ and leaves labelled $0,1,\ldots,k$ uniformly at random from the relevant finite set of isomorphism classes then giving each edge an independent Exponential$(1)$ random length. Moreover, once this metric structure is sampled, the Brownian embedding $\Phi$ is constructed by sampling independent Brownian motions indexed by each edge and concatenating these motions according to the tree structure in the natural way. To see that the relevant event has positive probability, take $k=2d$ and note simply that with positive probability:
			\begin{itemize}
				\item Each of the $2d-1$ Brownian motions indexed by an edge of the combinatorial tree not incident to one of the non-root leaves is contained in the box $[-1/64,1/64]^d$.
				\item For each $i\in \{\pm 1,\ldots,\pm d\}$, the Brownian motion indexed by the edge of the tree incident to the leaf labelled $\ell =2|i|-1+\mathbf{1}(i<0)$ starts in $[-1/32,1/32]^d$, ends in $2 e_i+[-1/32,1/32]^d$, and is contained in the cuboid $\{(z_1,\ldots,z_d) : |z_j| \leq 1/32$ for $j\notin \{i,-i\}$, $-1/4 \leq z_i \leq 5/2 \}$.
			\end{itemize}
			This event obviously has positive probability by standard properties of Brownian motion.
			
			It remains to deduce the claim regarding the boundary set from \eqref{eq:2d_spanning_cross}. The fact that we can require the number of boundary points of the cluster in each face to have at most order $n^2$ points follows immediately from \eqref{eq.expansion_11} (or rather its analogue for $0$ not at the center of the box, which follows from \eqref{eq.expansion_2} or the results of \cite{panis2026reversing}) and Markov's inequality; the more interesting part is to show that the number of boundary points is at least of order $n^2$ with good probability. Let us note to this end that if $E_i(x)$ holds, then there must exist a point $y$ in $\cC_n(x) \cap \partial_i^\circ \Lambda_n$ such that $y$ belongs to an open path of $\ell^\infty$ diameter at least $\frac{7n}{8}$ disjoint from $\cC_n(x)$. (Namely, we can take $y$ to be the last visit to $\cC_n(x)$ of an open simple path in the cuboid of the form guaranteed to exist by $E_i(x)$.) Thus, the probability that $E_i(x)$ holds but $|\cC_n(x) \cap \partial_i^\circ \Lambda_n|\leq \delta n^2$ is bounded by
			\[
			\P_{\p_c}(x\leftrightarrow \partial \Lambda_n) \delta n^2 \P_{\p_c}(0 \leftrightarrow \partial \Lambda_{7n/8}) \lesssim \delta n^{-2},
			\]
			which implies the claim in conjunction with \eqref{eq:2d_spanning_cross}.
		\end{proof}

		The next result then generalizes \eqref{eq.num_spanning}. It states that, if a $2d$-spanning cluster has positive probability to have the property $\mathbf{S}$, then typically there are $n^{d-6}$ such spanning clusters in $\Lambda_n$.
		\begin{prop}\label{prop.num_spanning_new}
			Assume there exists a constant $\delta > 0$, and the following lower bound holds uniformly for all $n \in \N_+$ and $x \in \Lambda_n$
			\begin{align}\label{eq.num_spanning_new.condition}
				\P_{\p_c}\Ll[\cC_n(x) \in \mathbf{S} \Big\vert  \cC_n(0) \text{ is } 2d\text{-spanning}\Rr] \geq \delta.
			\end{align}
			Then define the set
			\begin{align}\label{eq.def.spanning.S}
				\mathscr{S} := \{\cC \subset \Lambda_n: \cC \in \mathbf{S}, \text{ is } 2d\text{-spanning}, \,  c_0 n^4 < \vert \cC \vert < C_0 n^4 \}.
			\end{align} 
			For every $\epsilon \in (0,1)$, there exist a positive integer $N_0$ and finite positive constants $c_0, C_0, c_1, C_1$ depending on $\delta, \epsilon, d$ such that
			\begin{align}\label{eq.num_spanning.S}
				\forall n \geq N_0, \qquad \P_{\p_c}\Ll[ c_1 n^{d-6} \leq \vert \mathscr{S} \vert \leq C_1 n^{d-6} \Rr] > 1-\epsilon.
			\end{align}
		\end{prop}
		\begin{proof}
			Because $\mathscr{S}$ in \eqref{eq.def.spanning.S} has more constrains than that in \eqref{eq.num_spanning}, the upper bound is just a corollary. We thus focus on the lower bound. The proof follows the argument in  \cite[Theorem~13.22]{heydenreich2017progress}. Define that
			\begin{align*}
				X_z := \Ind{\cC_n(z) \in \mathbf{S}, \text{ is } 2d\text{-spanning}, \,  c_0 n^4 < \vert \cC_n(z) \vert < C_0 n^4},
			\end{align*}
			and $\vert \mathscr{S} \vert $ follows the study of the sum 
			\begin{align}\label{eq.count.Xz}
				N := \sum_{z \in \Lambda_{n/2}} X_z. 
			\end{align}
			
			The proof can be divided into $3$ steps.
			
			\textit{Step~1: the first moment of $X_z$.} Denote by 
			\begin{align*}
				A(n,c_0, C_0) := \{c_0 n^4 < \vert \cC_n(0) \vert < C_0 n^4 \},
			\end{align*}
			then we have
			\begin{align*}
				&\P_{\p_c}\Ll[ A^c(n,c_0, C_0)  \Big\vert \cC_n(0) \text{ is } 2d\text{-spanning}\Rr]\\
				&= \frac{\P_{\p_c}\Ll[ A^c(n,c_0, C_0) \cap\{ \cC_n(0) \text{ is } 2d\text{-spanning}\}\Rr]}{\P_{\p_c}\Ll[\cC_n(0) \text{ is } 2d\text{-spanning}\Rr]}\\
				&\leq \frac{\P_{\p_c}\Ll[ A^c(n,c_0, C_0) \cap\{ 0 \leftrightarrow \partial \Lambda_n\}\Rr]}{\P_{\p_c}\Ll[0 \leftrightarrow \partial \Lambda_n \Rr]} \frac{\P_{\p_c}\Ll[0 \leftrightarrow \partial \Lambda_n \Rr]}{\P_{\p_c}\Ll[\cC_n(0) \text{ is } 2d\text{-spanning}\Rr]}\\
				&\lesssim \P_{\p_c}\Ll[ A^c(n,c_0, C_0) \, \vert  0 \leftrightarrow \partial \Lambda_n\Rr].
			\end{align*}
			We drop the last factor as the two probabilities are comparable; see \eqref{eq.one_arm} and \eqref{eq.arm.2dspanning}. 	Recall the result in \eqref{eq.cluster_size}, there exists $c_0$ small enough and $C_0$ big enough, such that
			\begin{align*}
				\P_{\p_c}\Ll[ A^c(n,c_0, C_0)  \Big\vert \cC_n(0) \text{ is } 2d\text{-spanning}\Rr] \leq \frac{\delta}{2}.
			\end{align*}
			Now we combine the upper bound above and the condition \eqref{eq.num_spanning_new.condition}, which yields that 
			\begin{align}\label{eq.bound.delta2}
				\P_{\p_c}\Ll[  A(n,c_0, C_0) \cap \{\cC_n(0) \in \mathbf{S}\} \, \vert \cC_n(0) \text{ is } 2d\text{-spanning}\Rr] \geq \frac{\delta}{2}.
			\end{align}
			Then we calculate the expectation of $X_z$ for all $z \in \Lambda_{n/2}$
			\begin{equation}\label{eq.mean.Xz}
				\begin{split}
					&\E_{\p_c}[X_z]\\
					&= \P_{\p_c}\Ll[  A(n,c_0, C_0) \cap \{\cC_n(z) \in \mathbf{S}\} \, \vert \cC_n(z) \text{ is } 2d\text{-spanning}\Rr] \cdot \P_{\p_c}[\cC_n(z) \text{ is } 2d\text{-spanning}]\\
					&\geq \frac{\delta}{2} \cdot c_{\text{\ref{eq.arm.2dspanning}}}n^{-2}.
				\end{split}
			\end{equation}
			The estimate \eqref{eq.bound.delta2} and \eqref{eq.arm.2dspanning} are applied in the last line.

			\textit{Step~2: the first and second moments of $N$.} 
			The mean of $N$ follows \eqref{eq.mean.Xz}
			\begin{align}\label{eq.mean.N}
				\frac{1}{2} c_{\text{\ref{eq.arm.2dspanning}}} \delta n^{d-2} \leq \E_{\p_c}[N] \leq C_{\text{\ref{eq.arm.2dspanning}}}n^{d-2} .
			\end{align}
			Concerning the variance, we follow the trick in \cite[(13.6.13)]{heydenreich2017progress} that
			\begin{equation}\label{eq.var.N}
				\begin{split}
					\var_{\p_c}[N] &\leq \E_{\p_c}[N] + \sum_{\substack{z_1, z_2 \in \Lambda_{n/2} \\ z_1 \neq z_2}} \E_{\p_c}\Ll[X_{z_1} X_{z_2} \Ind{z_1 \stackrel{\Lambda_n}{\leftrightarrow} z_2}\Rr]\\
					&\leq \E_{\p_c}[N] + \sum_{\substack{z_1, z_2 \in \Lambda_{n/2} \\ z_1 \neq z_2}} \tau_{\p_c}(z_1, z_2)\\
					&\lesssim n^{2d-2}.
				\end{split}
			\end{equation}

			Recall that $d > 6$, so \eqref{eq.mean.N} and \eqref{eq.var.N} together imply that
			\begin{align*}
				\lim_{n \to \infty} \frac{\var_{\p_c}[N]}{\E_{\p_c}[N]^2} = \lim_{n \to \infty} \frac{n^{2d-2}}{(n^{d-2})^2} = 0.
			\end{align*}
			Following the Cauchy--Schwarz inequality, for all $\lambda \in (0,1)$, we have 
			\begin{align}\label{eq.PZ}
				\P_{\p_c}\Ll[ \vert N - \E_{\p_c}[N] \vert > \lambda \E_{\p_c}[N]\Rr] \geq \frac{\var_{\p_c}[N]}{\lambda^2 \E_{\p_c}[N^2]} \xrightarrow{n \to \infty} 0.
			\end{align}
			Therefore, the event $\{(1-\lambda) \E_{\p_c}[N] \leq N \leq (1+\lambda)\E_{\p_c}[N]\}$ holds with high probability when $n \to \infty$. This suggests that the value of $N$ is concentrated at $n^{d-2}$.

			\textit{Step~3: counting argument.} 
			Finally, notice that a cluster $\cC \in \mathscr{S}$ will be counted at most $\vert \cC \vert$ times in $N$, which is of order $n^4$ thanks to the constrain $\{ c_0 n^4 < \vert \cC \vert < C_0 n^4 \}$. This suggests that 
			\begin{align*}
				N \leq C_0 n^4 \vert \mathscr{S} \vert,
			\end{align*}
			which implies the following event with high probability
			\begin{align}\label{eq.NS}
				C_0^{-1}(1-\lambda)n^{-4}\E_{\p_c}[N] \leq \vert  \mathscr{S} \vert.
			\end{align} 
			Combining \eqref{eq.PZ}, \eqref{eq.NS}, and \eqref{eq.mean.N}, we choose $\lambda$ small enough and then an associated $c_1, C_1, N_0$, so \eqref{eq.num_spanning} is concluded.
		\end{proof}

		\noindent \textbf{Intrinsic diameter}
		Using the results above, we develop an estimate for the diameter of distance of spanning clusters. Recall that $\din$ stands for the chemical distance. We define the restricted version 
		\begin{align}\label{eq.def.dinS}
			\din(x,y;S) := \text{ chemical distance between } x,y \text{ restricted in }S.
		\end{align}
		Given $\cC$ a $2d$-spanning cluster in $\Lambda_n$, we define its intrinsic diameter 
		\begin{align}\label{eq.def.rin}
			\rin(\cC;\Lambda_n) := \max_{i\neq j \in \{\pm 1, \pm 2, \cdots, \pm d\}} \max_{\substack{x \in \partial_i \Lambda_n, \\ y \in \partial_j \Lambda_n}}\din(x, y;\Lambda_n).
		\end{align} 
		The following lemma will be used several times in the paper. 
		\begin{lemma}\label{lem.spanning.diam}
			For every $\epsilon \in (0,1)$, there exist a positive integer $N_0$ and finite positive constants $c_0, C_0, c_1, C_1, C_2$ depending on $\epsilon$ and $d$, such that the set
			\begin{align}\label{eq.def.spanning.diam}
				\mathscr{S} := \{\cC \subset \Lambda_n: \cC \text{ is } 2d\text{-spanning}, \, \rin(\cC;\Lambda_n) \leq C_2 n^2, \, c_0 n^4 < \vert \cC \vert < C_0 n^4 \}.
			\end{align} 
			satisfies the estimate holds
			\begin{align}\label{eq.num_spanning.diam}
				\forall n \geq N_0, \qquad \P_{\p_c}\Ll[ c_1 n^{d-6} \leq \vert \mathscr{S} \vert \leq C_1 n^{d-6} \Rr] > 1-\epsilon.
			\end{align}
		\end{lemma}

		One key input for Lemma~\ref{lem.spanning.diam} (and also \eqref{eq.chemical_1}, \eqref{eq.chemical_2} in previous work) is the following estimate for the longest self-avoiding open path. 
		\begin{lemma}[(105) of \cite{CHS23}]
			For every $x \in \H$, denote by $\mathfrak{L}(x)$ the length of the longest self-avoiding open path from $x$ to $\partial \H$ if it exists, then there exists a finite positive constant $C(d)$ such that 
			\begin{align}\label{eq.chemical_finite_moment}
				\sup_{x \in \H} \E_{\p_c}[\mathfrak{L}(x)\Ind{x \longleftrightarrow \partial \H}] < C_{\text{\ref{eq.chemical_finite_moment}}}.
			\end{align}
		\end{lemma}
		\begin{proof}[Proof of Lemma~\ref{lem.spanning.diam}]
			The lemma is in the regime of Proposition~\ref{prop.num_spanning_new}, so it suffices to obtain the condition \eqref{eq.num_spanning_new.condition} 
			\begin{align*}
				\cC \in \mathbf{S} \Longleftrightarrow \rin(\cC;\Lambda_n) \leq C_2 n^2.
			\end{align*}
			Therefore, we start by estimating the expectation of the diameter. For every $x \in \Lambda_n$
			let $\gamma(\cC_n(x))$ be the open path to reach the diameter $\rin(\cC_n(x);\Lambda_n)$, then there exist $i,j \in \{\pm 1, \cdots, \pm d\}$ such that $\gamma(\cC_n(x))$ is a simple open path connecting $\partial_i \Lambda_n$ and $\partial_j \Lambda_n$. Let $\gamma_i(x)$ (\emph{resp.} $\gamma_j(x)$) be the simple open path connecting $x$ to $\gamma(\cC_n(x))$, and then follow $\gamma(\cC_n(x))$ to connect  $\partial_i \Lambda_n$ (\emph{resp.} $\partial_j \Lambda_n$). Clearly, we have 
			\begin{align*}
				\rin(\cC_n(x); \Lambda_n) = \vert \gamma(\cC_n(x)) \vert \leq \vert \gamma_i(x) \vert + \vert \gamma_j(x) \vert.
			\end{align*} 
			This entails that 
			\begin{align}\label{eq.chemical.half}
				\begin{split}
					&\E_{\p_c}\Ll[\rin(\cC_n(x); \Lambda_n) \Ind{\cC_n(x) \text{ is } 2d\text{-spanning}}\Rr] \\
					&\leq \E_{\p_c}\Ll[\vert \gamma_i(x) \vert \Ind{\cC_n(x) \text{ is } 2d\text{-spanning}}\Rr] + \E_{\p_c}\Ll[\vert \gamma_j(x) \vert \Ind{\cC_n(x) \text{ is } 2d\text{-spanning}}\Rr]\\
					&\leq \E_{\p_c}\Ll[\vert \gamma_i(x) \vert \Ind{x \longleftrightarrow  \partial_i \Lambda_n}\Rr] + \E_{\p_c}\Ll[\vert \gamma_j(x) \vert \Ind{x \longleftrightarrow  \partial_j \Lambda_n}\Rr]\\
					&\leq 2 \sup_{x \in \H} \E_{\p_c}[\mathfrak{L}(x)\Ind{x \longleftrightarrow \partial \H}] < 2C_{\text{\ref{eq.chemical_finite_moment}}}.
				\end{split}
			\end{align}
			The last two lines make use of the moment estimate of the longest self-avoiding open path in the half-space \eqref{eq.chemical_finite_moment}. 
			
			We then insert the probability of $2d$-spanning cluster \eqref{eq.arm.2dspanning} in \eqref{eq.chemical.half}
			\begin{align*}
				\E_{\p_c}\Ll[\rin(\cC_n(x); \Lambda_n) \, \vert \cC_n(x) \text{ is } 2d\text{-spanning}\Rr] \leq \frac{2C_{\text{\ref{eq.chemical_finite_moment}}}}{\P_{\p_c}[\cC_n(x) \text{ is } 2d\text{-spanning}]} \leq \frac{C n^2}{c_{\text{\ref{eq.arm.2dspanning}}}}.
			\end{align*}
			A Markov inequality yields that 
			\begin{align*}
				\P_{\p_c}\Ll[\rin(\cC_n(x); \Lambda_n) > C_2 n^2 \, \vert \cC_n(x) \text{ is } 2d\text{-spanning}\Rr] \leq  \frac{2C_{\text{\ref{eq.chemical_finite_moment}}}}{C_2 c_{\text{\ref{eq.arm.2dspanning}}}}.
			\end{align*}
			By choosing $C_2$ large enough, the condition $\mathbf{S}$ in \eqref{eq.num_spanning_new.condition} states that 
			\begin{align*}
				\P_{\p_c}\Ll[\rin(\cC_n(x); \Lambda_n) \leq C_2 n^2 \, \vert \cC_n(x) \text{ is } 2d\text{-spanning}\Rr] \geq 1 -  \frac{2C_{\text{\ref{eq.chemical_finite_moment}}}}{C_2 c_{\text{\ref{eq.arm.2dspanning}}}} = \delta > 0.
			\end{align*}
			Then Proposition~\ref{prop.num_spanning_new} applies and we finish the proof of \eqref{eq.num_spanning.diam}.
		\end{proof}

		\section{Recap of homogenization}\label{sec.homo}
		This section summarizes the basic facts from \cite{AKMbook} and \cite{armstrong2018elliptic} about the homogenization on the percolation cluster. Some definitions will be reformulated so that they are easier to analyze in the rest of paper.
		
		\subsection{Definition and subadditivity}
		For every bounded domain $U \subset \Zd$ and $q \in \R^d$, we define the Dirichlet energy associated to the affine boundary condition $\ell_{q}$ 
		\begin{equation}\label{eq.defmu}
			\mathcal{E}_{U,q}(\a) := \inf_{v \in \ell_{q}+C_0(U^+)} \frac{1}{|U|}\Ll \langle\nabla v,\a \nabla v\Rr \rangle_{U},
		\end{equation}
		and denote by $v_{U,q}$ its minimiser. The domain $U^+$ is exactly the support of the functions involved in the functional; see \eqref{eq.bracket2} and \eqref{eq.edge_plus} for the definitions. 	It is clear that $q \mapsto \mathcal{E}_{U,q}(\a)$ is a bilinear application, so we define \emph{the quenched homogenized matrix} $\A(U)$ 
		\begin{align}\label{eq.A}
			q \cdot \A(U) q := \mathcal{E}_{U,q}(\a).
		\end{align}
		Notice that $\A(U)$ implicitly depends on the realization of the percolation $\a$, thus we define the annealed version
		\begin{align}\label{eq.A_anneal}
			\Ab(U,\p) := \E_\p[\A(U)].
		\end{align}
		The specific choice of the domain in \eqref{eq.defmu} allows us to recover the exact subadditivity of $\mathcal{E}_{U,q}(\a)$. This property is the heart of the renormalization approach developed by Armstrong, Kuusi, Mourrat, and Smart \cite{armstrong2016quantitative, AKMbook}, so let us include the complete proof here.
		
		\begin{lemma}\label{lem.subadditive}
			For every bounded domain  $U \subset \Zd$ and its partition $U = \sqcup_{i=1}^n U_i$, the following inequality holds for all the realization $\a$
			\begin{align}\label{eq.subadditive}
				\A(U) \leq \sum_{i=1}^n \frac{\vert U_i \vert}{\vert U \vert} \A(U_i).
			\end{align}
			The same statement also holds for its means $\Ab(U,\p)$.
		\end{lemma}
		\begin{proof}
			We observe that, by gluing the minimiser in every subdomain
			\begin{align}\label{eq.glue}
				\tilde{v} = \ell_{q} + \sum_{i=1}^n (v_{U_i,q} - \ell_{q}),
			\end{align}
			the function $\tilde{v} \in \ell_{q} + C_0(U^+)$, thus provides a sub-minimiser for the functional $\mathcal{E}_{U,q}(\a)$. Therefore, we have 
			\begin{align*}
				\mathcal{E}_{U,q}(\a) &\leq \frac{1}{|U|}\Ll \langle\nabla \tilde{v},\a \nabla \tilde{v}\Rr \rangle_{U}\\
				&= \sum_{i=1}^n \frac{1}{|U|} \Ll \langle\nabla \tilde{v},\a \nabla \tilde{v}\Rr \rangle_{U_i}\\
				&= \sum_{i=1}^n \frac{1}{|U|} \Ll \langle\nabla v_{U_i,q},\a \nabla v_{U_i,q}\Rr \rangle_{U_i}\\
				&= \sum_{i=1}^n \frac{\vert U_i \vert}{\vert U \vert}\mathcal{E}_{U_i,q}(\a).
			\end{align*}
			The second line strongly relies on the nice structure in \eqref{eq.union_edge} and \eqref{eq.id_gradient}
			\begin{align}\label{eq.union_edge_copy}
				\Ed^+(U) = \bigsqcup_{i=1}^n \Ed^+(U_i),
			\end{align}
			and the third line is because of \eqref{eq.glue}, which yields the desired result for $\A(U)$. We take the expectation under $\P_{\p}$ to obtain the result for its mean $\Ab(U,\p)$.
		\end{proof}
		
		One may find the definition of $\mathcal{E}_{U,q}(\a)$ in \eqref{eq.defmu} slightly different from that in \cite{armstrong2018elliptic}, but the difference is a boundary layer and vanishes in the limit. Moreover, the choice of $U^+$ in the variational formula benefits from \eqref{eq.union_edge_copy}, then the subadditivity avoids the boundary layer term in \cite[eq.(4.49)]{armstrong2018elliptic}.
		
		The homogenized conductivity then defines the diffusivity. 
		\begin{lemma}[Section~5 of \cite{armstrong2018elliptic}, eq.(181) of \cite{dario2021quantitative}]\label{lem.basic}
			
			For every $\p \in (\p_c,1]$, the diffusivity for the random walk on the infinite cluster can be expressed as 
			\begin{align}\label{eq.def_sigma}
				\ab(\p) \id := \lim_{m \to \infty} \Ab(\Lambda_{3^m L}, \p), \qquad \frac{\sigma^2(\p)}{2} = \frac{\ab(\p)}{\theta(\p)}.
			\end{align} 
			The first convergence is actually a decreasing to the limit, and the limit coincides for all $L \in \N_+$. 
		\end{lemma}
		\begin{proof}
			The first limit in \eqref{eq.def_sigma} was proved in \cite{armstrong2018elliptic}.	By the stationarity and \eqref{eq.subadditive}, once we choose $L$, the decreasing to the limit holds along the sequence $(\Lambda_{3^m L})_{m \in \N_+}$. To see the limit does not depend on $L$, one needs some nice geometric argument of partition. For example, see \cite[Lemma~A.1]{funaki2024quantitative} for a proof using the Whitney decomposition in \cite[Appendix~J]{Grafakos}.
			
			The second identity in \eqref{eq.def_sigma} is known as Einstein's conductivity-diffusivity relation. We refer to \cite[eq.(181)]{dario2021quantitative} and the paragraph around for the justification and discussion.
		\end{proof}

		\subsection{Boundary layer and water cubes}\label{subsec.watercube}
		Define $\A_{11}(U)$ as its left-top entry of $\A(U)$ and $\Ab_{11}(U,\p)$ as its mean
		\begin{align}\label{eq.def.A11}
			\A_{11}(U) := e_1 \cdot \A(U) e_1, \qquad \Ab_{11}(U,\p):= \E_{\p}[\A_{11}(U)].
		\end{align}
		It is very natural to study $\ab(\p)$ by its approximation $\Ab_{11}(\Lambda_{\ell_*},\p)$ in the correlation length. However, the following observation reveals some issue from the boundary layer.
		
\begin{figure}[h!]
			\centering
			\includegraphics[width=0.3\textwidth]{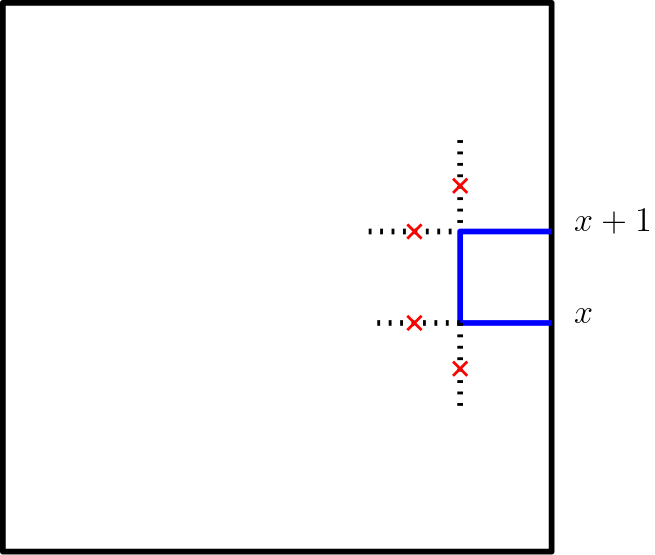}
			\caption{An illustration of the small loop near $\partial\Lambda_{L}$.}\label{fig.loop}
		\end{figure}

		\begin{lemma}\label{lem.boundary}
			For every $\p \in (0,1)$, we have the following estimate that 
			\begin{align}\label{eq.ab_boundary}
				\Ab_{11}(\Lambda_L, \p) \geq 2(d-1) \p^3 (1-\p)^{2d-2}L^{-1}.
			\end{align}
		\end{lemma}
		\begin{proof}
			It suffices to create a small open loop of length $3$ near the boundary  $\partial \Lambda_L$ parallel to $e_1$, and close its interface between other vertices. See Figure~\ref{fig.loop} for an illustration. Then the gradient $\nabla v_{\Lambda_L,e_1} = 1$  on the segment can be obtained by explicit calculation. The total contribution is a boundary layer, which yields the lower bound \eqref{eq.ab_boundary}.  
		\end{proof}
		
		Viewing Lemma~\ref{lem.boundary}, we choose $L = \ell_* = \lfloor(\p - \p_c)^{-1/2}\rfloor$ and obtain that
		\begin{align}\label{eq.boundary_paradox}
			\Ab_{11}(\Lambda_{\ell_*},\p) \geq \p^3 (1-\p)^{2d-2} (\p - \p_c)^{\frac{1}{2}}.
		\end{align}
		This implies that the boundary layer effect is too strong, and we cannot observe the desired rate $(\p-\p_c)^3$ in $\Lambda_{\ell_*}$.

\begin{figure}[h!]
			\centering
			\includegraphics[width=0.6\textwidth]{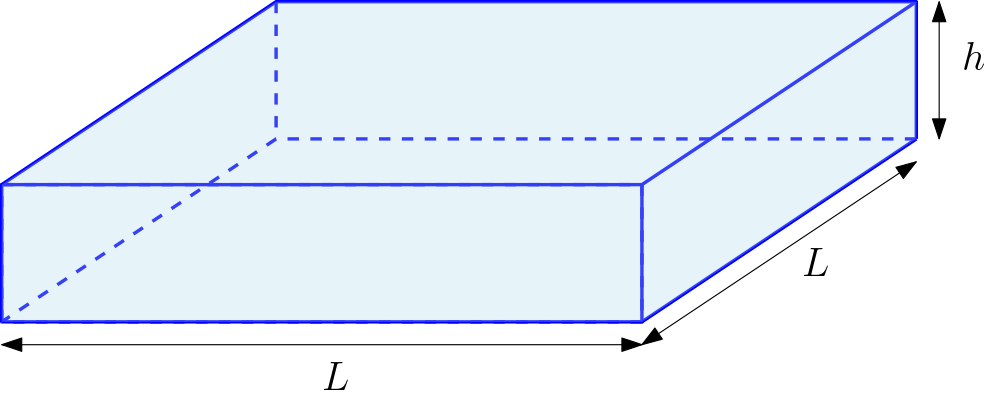}
			\caption{An illustration of the rectangle $\Lambda_{h,L}$.}\label{fig.cube_new}
		\end{figure}

		To handle this problem of the boundary layer, we propose a new family of slab $\Lambda_{h,L}$ we call \emph{the water cube}\footnote{We refer the reader to \url{https://en.wikipedia.org/wiki/Water_Cube} for an introduction of the Water Cube.}
		\begin{align}\label{eq.def_cube_new}
			\Lambda_{h,L} = \Ll(-\frac{h}{2}, \frac{h}{2}\Rr) \times \Ll(-\frac{L}{2}, \frac{L}{2}\Rr) \cdots \times \Ll(-\frac{L}{2}, \frac{L}{2}\Rr) \cap \Zd.
		\end{align}
		See an illustration in Figure~\ref{fig.cube_new}. We call $h$ the height and $L$ the length/width, and usually keep the convention $1 < h \ll L$ for $\Lambda_{h,L}$. 
		Then a similar analysis as \eqref{eq.ab_boundary} will also yield
		\begin{align*}
			\Ab_{11}(\Lambda_{h,L}, \p) \geq 2(d-1) \p^3 (1-\p)^{2d-2}L^{-1}.
		\end{align*}
		This lower bound is independent of the height $h$, so we take $h = \ell_*$ and let $L \to \infty$, then the boundary layer effect vanishes.  Actually, for both two types of cubes,  the creation of small loops around the upper face and the lower face does not contribute the energy, because they are of the same potential. The real contribution comes from the boundary parallel to $e_1$, and $\Lambda_{h, L}$ has more freedom to manipulate.
		
		Viewing the discussions above, the water cube $\Lambda_{h, L}$ will be the basic object where we study the approximation of $\ab$. We give another explanation why it is better than the usual cube $\Lambda_L$. That is, the water cube $\Lambda_{h,L}$ can be seen as a renormalization of $\Lambda_{h}$, so the effect of boundary layer decreases.
		\begin{lemma}\label{lem.chain}
			Given $h \in \N_+$ and $L = N h$ with $N \in \N_+$, we have
			\begin{align}\label{eq.chain}
				\ab(\p) \leq \Ab_{11}(\Lambda_L, \p) \leq \Ab_{11}(\Lambda_{h,L},\p) \leq \Ab_{11}(\Lambda_{h},\p).
			\end{align}
		\end{lemma}
		\begin{proof}
			The first inequality is due to Lemma~\ref{lem.basic}. The decreasing to the limit holds in the partial order of positive symmetric matrix, then we take the vector $e_1$ and obtain that $\Ab_{11}(\Lambda_{3^m L}, \p)$ decreases to the limit $\ab(\p)$. Because $L = N h$, we observe that 
			\begin{align*}
				\Lambda_{h, L} & = \bigsqcup_{z \in h \Zd \cap \Lambda_{h, L}}  \Lambda_{h}(z),\\
				\Lambda_{L} & = \bigsqcup_{k \in h \mathbb{Z} \cap (-\frac{L}{2}, \frac{L}{2})} \Lambda_{h, L}(k e_1).
			\end{align*}
			Then the subadditivity in Lemma~\ref{lem.subadditive} concludes the second and third inequalities in \eqref{eq.chain}.
		\end{proof}
		
		\subsection{An alternative definition}
		We spent efforts discussing the influence of the boundary condition in the last paragraphs. The quantity $\Ab_{11}(\Lambda_{h, L}, \p)$ with $h = \ell_* \ll L$ is one candidate solution. Recall \eqref{eq.def.TopBottom} for the definition of $[\partial \Lambda_{h,L}]^{\pm}$ as the top and bottom faces of the cube. Since the width $L$ is much larger than the height $h$, the influence of the boundary condition from $[\partial \Lambda_{h,L}]^{\pm}$ should be dominant for in the bulk. This is related to another approach to conductivity, which is closer to that in \cite{kestenICM84} by Kesten:
		\begin{align}\label{eq.def.aU}
			\a(\Lambda_{h,L}) \equiv \mathscr{E}_{\Lambda_{h,L}}(\a) := \inf_{v \equiv \ell_{e_1} \text{ on } [\partial \Lambda_{h,L}]^{\pm}} \frac{1}{|\Lambda_{h,L}|}\Ll \langle\nabla v,\a \nabla v\Rr \rangle_{\Lambda_{h,L}}.
		\end{align}
		That is, we relax the boundary condition except that on $[\partial \Lambda_{h,L}]^{\pm}$. We denote by $u_{\Lambda_{h,L}}$ the minimiser in \eqref{eq.def.aU}. The annealed version is defined by
		\begin{align}\label{eq.def.abU}
			\ab(\Lambda_{h,L}) := \E_{\p}[\a(\Lambda_{h,L})].
		\end{align}
		
		The following lemma justifies that, the two quantities can be considered nearly equivalent $\ab(\Lambda_{h,L}) \simeq \Ab_{11}(\Lambda_{h,L})$.
		\begin{lemma}\label{lem.chain.2}
			Given $h, L \in \N_+$ satisfying that $L \geq 2h$, we have
			\begin{align}\label{eq.chain.low}
				\ab(\Lambda_{h,L}, \p) \leq  \Ab_{11}(\Lambda_{h,L}, \p),
			\end{align}
			and
			\begin{align}\label{eq.chain.2}
				\ab(\p) \leq \ab(\Lambda_{h,L}, \p) + h^2 L^{-1}.
			\end{align}
		\end{lemma}
		\begin{proof}
			Because there are more admissible functions in \eqref{eq.def.aU} than \eqref{eq.defmu}, naturally the infimum is smaller and we have \eqref{eq.chain.low}. To prove the inequality in another direction, let us pick $L, \tilde{L} \in \N_+$ such that $\tilde{L} = N L$ with $N \in \N_+$. Then we compare $\A_{11}(\Lambda_{h,\tilde{L}})$ and $\a(\Lambda_{h,L})$, using an argument similar as \eqref{eq.glue}. We propose a sub-minimiser for $\A_{11}(\Lambda_{h,\tilde{L}})$ that, for every $ z \in L \Zd \cap \Lambda_{h, \tilde{L}}$ and $x \in \Lambda_{h, L}(z)$
			\begin{align}\label{eq.glue2}
				\tilde{v}(x) = \Ll\{\begin{array}{ll}
					u_{\Lambda_{h, L}(z)}(x), &	\dex(z, \partial \Lambda_{h, \tilde{L}})  > L, \\
					v_{\Lambda_{h, L}(z),e_1}(x), &	\dex(z, \partial \Lambda_{h, \tilde{L}})  < L.
				\end{array}
				\Rr.
			\end{align}
			\begin{figure}[h!]
				\centering
				\includegraphics[width=0.8\textwidth]{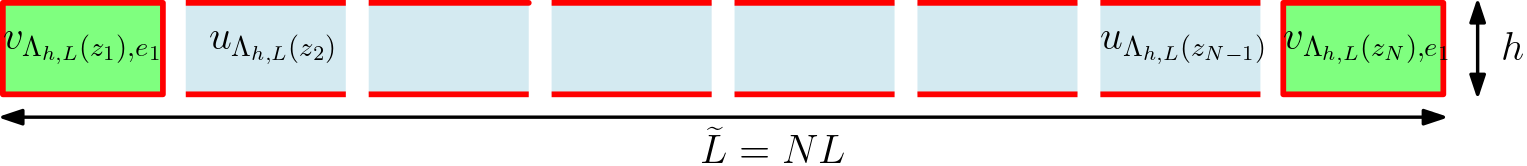}
				\caption{An illustration of \eqref{eq.glue2}. The red segment represents the boundary condition $\ell_{e_1}$.}\label{fig.glue}
			\end{figure}
			That is, we put the minimisers of $\A_{11}(\Lambda_{h,L}(z))$ near $\partial \Lambda_{h, \tilde{L}}$ to meet the boundary conditions of $\A_{11}(\Lambda_{h,\tilde{L}})$, and glue the minimisers of $\a(\Lambda_{h,L}(z))$ in the interior to satisfy the boundary conditions on $[\partial \Lambda_{h,\tilde{L}}]^{\pm}$. See Figure~\ref{fig.glue} for an illustration. Therefore, we have $\tilde{v} \in \ell_{e_1}+C_0(\Lambda_{h,\tilde{L}}^+)$, which implies that
			\begin{align*}
				\A_{11}(\Lambda_{h,\tilde{L}}) & \leq \frac{1}{|\Lambda_{h,\tilde{L}}|}\Ll \langle\nabla \tilde{v},\a \nabla \tilde{v}\Rr \rangle_{\Lambda_{h,\tilde{L}}}\\
				&\leq \frac{\vert \Lambda_{h, L}\vert}{\vert \Lambda_{h, \tilde{L}}\vert} \Ll(\sum_{\substack{z \in L \Zd \cap \Lambda_{h, \tilde{L}}: \\
						\dex(z, \partial \Lambda_{h, \tilde{L}})  > L}} \a(\Lambda_{h,L}(z)) + \sum_{\substack{z \in L \Zd \cap \Lambda_{h, \tilde{L}}: \\
						\dex(z, \partial \Lambda_{h, \tilde{L}})  < L}} \A_{11}(\Lambda_{h,L}(z))\Rr) + h^2 L^{-1}.
			\end{align*}
			In the analysis, the term $h^2 L^{-1}$ is the error generated in the interfaces between different cubes $\Lambda_{h,L}(z)$. Since there is not much information for the functions there, and we just utilize the bound $\vert \nabla \tilde{v} \vert \leq h$ from the maximum principle. We take the expectation under $\P_{\p}$ and obtain that 
			\begin{align}\label{eq.sub.Aa}
				\Ab_{11}(\Lambda_{h,\tilde{L}}, \p) \leq \ab(\Lambda_{h,L}, \p) + L \tilde{L}^{-1} \Ab_{11}(\Lambda_{h,L}, \p) + h^2 L^{-1}.
			\end{align}
			Combining the equations \eqref{eq.sub.Aa} and \eqref{eq.chain}, then we take $\tilde{L} \to \infty$
			\begin{align*}
				\ab(\p) &\leq \lim_{\tilde{L} \to \infty} \Ab_{11}(\Lambda_{h,\tilde{L}}, \p) \\
				&\leq \lim_{\tilde{L} \to \infty}\ab(\Lambda_{h,L}, \p) + L \tilde{L}^{-1} \Ab_{11}(\Lambda_{h,L}, \p) + h^2 L^{-1} \\
				&= \ab(\Lambda_{h,L}, \p) + h^2 L^{-1}.
			\end{align*}
			This is the desired result in \eqref{eq.chain.2}.
		\end{proof}
		Thanks to the upper bound \eqref{eq.chain.2}, the quantity $\ab(\Lambda_{h,L}, \p)$ will be the key object discussed in the following sections.
		
		\section{Upper bound I: Analysis of pivotal edges}\label{sec.pivotal}
		Similar as the study of the connectivity $\theta(\p)$, a useful technique is to analyze the pivotal edges in the conductivity $\ab(\p)$. We present this idea and give a weak \emph{a priori} bound for the conductivity $\ab(\Lambda_{h,L}, \p)$.
		

		\subsection{Pivotal event in Dirichlet energy}\label{subsec.pivotal}
		Recall that a bond $e'$ is pivotal for the event $E$ if $E$ happens when $\a(e') = 1$ and does not happen when $\a(e') = 0$. The following observation relates the harmonic function and the two-arm event. Since the lemma does not involve any numerical analysis, we let $\Lambda$ stand for a general water cube in \eqref{eq.def_cube_new} without specific height and width.
		
		\begin{lemma}\label{lem.two_arm}
			Let $\Lambda$ be a cube and  $u_{\Lambda}$ be defined in \eqref{eq.def.aU}. Given $e = \{e_\bullet, e_\circ\} \subset \Ed^+(\Lambda)$ satisfying
			\begin{align}\label{eq.boundary_connection}
				\{e_\bullet \longleftrightarrow [\partial \Lambda]^{\pm}\} \cap \{e_\circ \longleftrightarrow  [\partial \Lambda]^{\pm}\},
			\end{align}
			then $u_{\Lambda}(e_\bullet) \neq  u_{\Lambda}(e_\circ)$ implies the two-arm event $\tarm(e, \Lambda)$ defined as  
			\begin{multline}\label{eq.two_arm}
				\tarm(e, \Lambda) := \Ll(\{e_\bullet \longleftrightarrow [\partial \Lambda]^-\} \circ \{e_\circ \longleftrightarrow [\partial \Lambda]^+\}\Rr)\\
				\cup  \Ll(\{e_\bullet \longleftrightarrow [\partial \Lambda]^+\} \circ \{e_\circ \longleftrightarrow [\partial \Lambda]^-\}\Rr).
			\end{multline}
		\end{lemma}
		\begin{proof}
			We divided the proof into two steps.
			
			\emph{Step~1: disjoint realization.} Firstly, we strengthen \eqref{eq.boundary_connection} as
			\begin{align}\label{eq.two_arm_baby}
				\{e_\bullet \longleftrightarrow [\partial \Lambda]^{\pm}\} \circ \{e_\circ \longleftrightarrow \partial [\partial \Lambda]^{\pm}\}.
			\end{align}
			The argument is a  proof by contradiction. Concerning the sets of the pivotal edges
			\begin{align*}
				S(e_\bullet, \Lambda) &:= \Ll\{ e' \in \Ed^+(\Lambda) : e' \text{ is pivotal for } \{e_\bullet \longleftrightarrow [\partial \Lambda]^{\pm}\} \Rr\},\\
				S(e_\circ, \Lambda) &:= \Ll\{ e' \in \Ed^+(\Lambda) : e' \text{ is pivotal for } \{e_\circ \longleftrightarrow [\partial \Lambda]^{\pm}\} \Rr\},
			\end{align*}
			and assume that 
			\begin{align}\label{eq.S.assumption}
				S(e_\bullet, \Lambda) \cap S(e_\circ, \Lambda) \neq \emptyset,
			\end{align}
			then there exists an edge $e'$ such that all the open paths connecting $e_\bullet, e_\circ$ to $[\partial \Lambda]^{\pm}$ must pass. We denote by $e'_*$ one endpoint of this edge.  Let $(X_n)_{n \in \N}$ be the simple random walk on the percolation,  and $P^{\a}_x$ stand for the probability space of random walk starting from $x$ under the realization $\a$. We also define $\tau_{0}$ and $\tau_{e'_*}$ respectively as the hitting time on the boundary $[\partial \Lambda]^{\pm}$ and $e'_*$
			\begin{align*}
				\tau_{0} := \inf\{n \in \N: X_n \in [\partial \Lambda]^{\pm}\}, \qquad \tau_{e'_*} := \inf\{n \in \N: X_n = e'_*\}.
			\end{align*}
			Then the probabilistic representation of harmonic function yields
			\begin{align*}
				u_{\Lambda}(e_\bullet) = E^{\a}_{e_\bullet}[X_{\tau_{0}}], \qquad u_{\Lambda}(e_\circ) = E^{\a}_{e_\circ}[X_{\tau_{0}}], \qquad u_{\Lambda}(e'_*) = E^{\a}_{e'_*}[X_{\tau_{0}}].
			\end{align*}
			By the definition of pivotal edge, we have
			\begin{align*}
				\tau_{e'_*} < \tau_{0}.
			\end{align*}   
			Using the strong Markov property, the process after $X_{\tau_{e'_*}}$ can be seen as a Markov chain starting from $e'_*$, so we have
			\begin{align*}
				u_{\Lambda}(e_\bullet) = E^{\a}_{e_\bullet}[X_{\tau_{0}}] = E^{\a}_{e_\bullet}[E^{\a}_{e_\bullet}[X_{\tau_{0}} \vert X_{\tau_{e'_*}} ]] = E^{\a}_{e'_*}[X_{\tau_{0}}] = u_{\Lambda}(e'_*).
			\end{align*}
			This justifies that
			\begin{align*}
				u_{\Lambda}(e_\bullet) = u_{\Lambda}(e_\circ) = u_{\Lambda}(e'_*),
			\end{align*}
			which is a contradiction.
			
			In conclusion, the assumption \eqref{eq.S.assumption} does not hold and we have $S(e_\bullet, \Lambda) \cap S(e_\circ, \Lambda) = \emptyset$. Let $\gamma(e_\bullet)$ be an open path connecting $e_\bullet$ to $[\partial \Lambda]^{\pm}$, then all the edges along $\gamma(e_\bullet)$ are not pivotal for $\{e_\circ  \longleftrightarrow [\partial \Lambda]^{\pm}\}$. Therefore, even when closing $\gamma(e_\bullet)$, we can still find another open path $\gamma(e_\circ)$ connecting $e_\circ$ and  $[\partial \Lambda]^{\pm}$. This concludes \eqref{eq.two_arm_baby}.
			
			\smallskip
			
			\emph{Step~2: classification.} The event \eqref{eq.two_arm_baby} can be further classified into $4$ cases. See Table~\ref{t.classification} and Figure~\ref{fig.classification} for an illustration.
			\begin{equation}\label{eq.pivot.4cases}
				\begin{split}
					\text{case a):} \qquad \{e_\bullet \longleftrightarrow [\partial \Lambda]^-\} \circ \{e_\circ \longleftrightarrow [\partial \Lambda]^+\},\\
					\text{case b):} \qquad \{e_\bullet \longleftrightarrow [\partial \Lambda]^+\} \circ \{e_\circ \longleftrightarrow [\partial \Lambda]^-\},\\
					\text{case c):} \qquad \{e_\bullet \longleftrightarrow [\partial \Lambda]^-\} \circ \{e_\circ \longleftrightarrow [\partial \Lambda]^-\},\\
					\text{case d):} \qquad  \{e_\bullet \longleftrightarrow [\partial \Lambda]^+\} \circ \{e_\circ \longleftrightarrow [\partial \Lambda]^+\}.\\
				\end{split}
			\end{equation}
			The two-arm event $\tarm(e, \Lambda)$ in \eqref{eq.two_arm} counts a) and b), which is a subset of \eqref{eq.two_arm_baby}. We thus study the cases c) and d). Actually, if $[\partial \Lambda]^-$ or $[\partial \Lambda]^+$ is the only face that the edge $e$ connects to, then $u_{\Lambda}(e_\bullet) = u_{\Lambda}(e_\circ)$ by harmonic extension from this unique face. Therefore, they are reduced to 
			\begin{align*}
				\text{case c):} \qquad \Ll(\{e_\bullet \longleftrightarrow [\partial \Lambda]^-\} \circ \{e_\circ \longleftrightarrow [\partial \Lambda]^-\}\Rr) \cap \{\{e_\bullet, e_\circ\} \longleftrightarrow [\partial \Lambda]^{+}\},
			\end{align*}
			or
			\begin{align*}
				\text{case d):} \qquad \Ll(\{e_\bullet \longleftrightarrow [\partial \Lambda]^+\} \circ \{e_\circ \longleftrightarrow [\partial \Lambda]^+\}\Rr) \cap \{\{e_\bullet, e_\circ\} \longleftrightarrow [\partial \Lambda]^{-}\}.
			\end{align*}
			By symmetry, we just study the case c). We claim that, this case can be reduced to a) and b) in \eqref{eq.two_arm}, i.e.
			\begin{multline}\label{eq.two_arm_reduction}
				\Ll(\{e_\bullet \longleftrightarrow [\partial \Lambda]^-\} \circ \{e_\circ \longleftrightarrow [\partial \Lambda]^-\}\Rr) \cap \{\{e_\bullet, e_\circ\} \longleftrightarrow [\partial \Lambda]^{+}\} \\
				\subset \Ll(\{e_\bullet \longleftrightarrow [\partial \Lambda]^+\} \circ \{e_\circ \longleftrightarrow [\partial \Lambda]^-\}\Rr) 
				\cup  \Ll(\{e_\bullet \longleftrightarrow [\partial \Lambda]^-\} \circ \{e_\circ \longleftrightarrow [\partial \Lambda]^+\}\Rr).
			\end{multline}
			Let $\gamma_1, \gamma_2$ be two disjoint open paths such that $e_\bullet \stackrel{\gamma_1}{\longleftrightarrow} [\partial \Lambda]^-$ and $e_\circ \stackrel{\gamma_2}{\longleftrightarrow} [\partial \Lambda]^-$, and $\gamma_3$ be the third open paths connecting $e_\bullet$ or $e_\circ$ to $[\partial \Lambda]^{+}$. Walking along the path $\gamma_3$ from $e_\bullet$ or $e_\circ$ to $[\partial \Lambda]^{+}$, and let $z$ be the last vertex of $\gamma_3$ that intersects $\gamma_1 \cup \gamma_2$, and $\gamma'_3$ be the segment starting from $z$. Clearly, we have 
			\begin{align*}
				\gamma'_3 \cap (\gamma_1 \cup \gamma_2) = z.
			\end{align*}
			\begin{enumerate}
				\item If $z \in \gamma_1$, then the paths $e_\bullet \stackrel{\gamma_1}{\longleftrightarrow} z \stackrel{\gamma'_3}{\longleftrightarrow}[\partial \Lambda]^{+}$ and $e_\circ \stackrel{\gamma_2}{\longleftrightarrow} [\partial \Lambda]^-$ are disjoint.
				\item If $z \in \gamma_2$, then the paths $e_\circ \stackrel{\gamma_2}{\longleftrightarrow} z \stackrel{\gamma'_3}{\longleftrightarrow}[\partial \Lambda]^{+}$ and $e_\bullet \stackrel{\gamma_1}{\longleftrightarrow} [\partial \Lambda]^-$ are disjoint.
			\end{enumerate}
			Therefore, we find another open path connecting to $[\partial \Lambda]^{+}$, which justifies \eqref{eq.two_arm_reduction}. Moreover, it also concludes \eqref{eq.two_arm}.
			
			\begin{table}[ht]
				\begin{center}
					\begin{tabular}{ c|cccccc } 
						& a) & b) & c) & d)  \\
						\toprule
						$e_\bullet$ & $-$ & $+$ & $-$ & $+$ \\
						\midrule 
						$e_\circ$ & $+$ & $-$ & $-$ & $+$ \\ 
						\bottomrule
					\end{tabular}
					\caption{Classification of the connecting faces.}
					\label{t.classification}
				\end{center}
			\end{table}
			
			\begin{figure}[ht]
				\centering
				\includegraphics[width=0.45\textwidth]{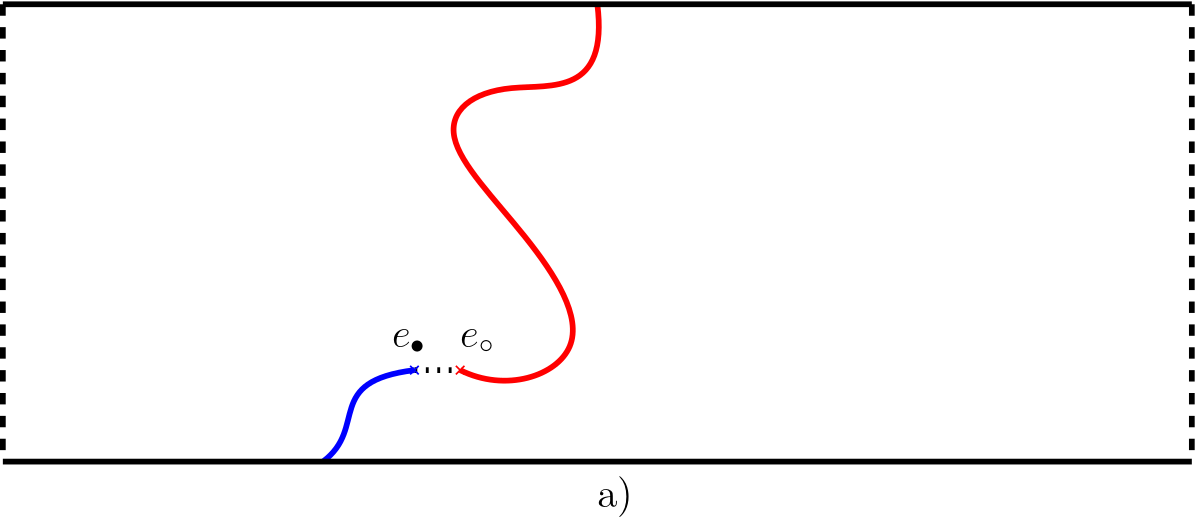} \quad
				\includegraphics[width=0.45\textwidth]{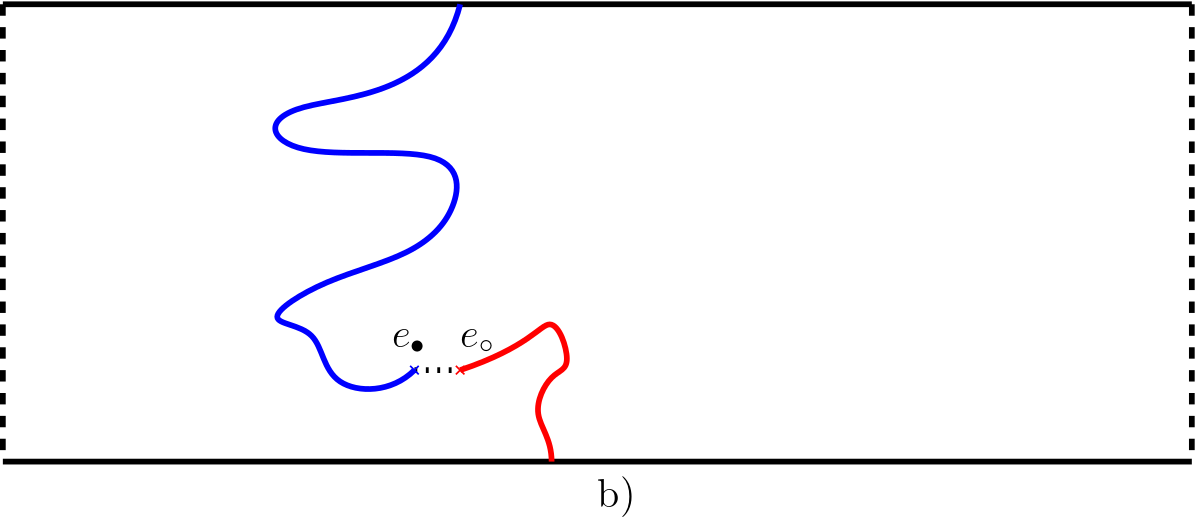}
				\\
				\vspace{0.5cm}
				\includegraphics[width=0.45\textwidth]{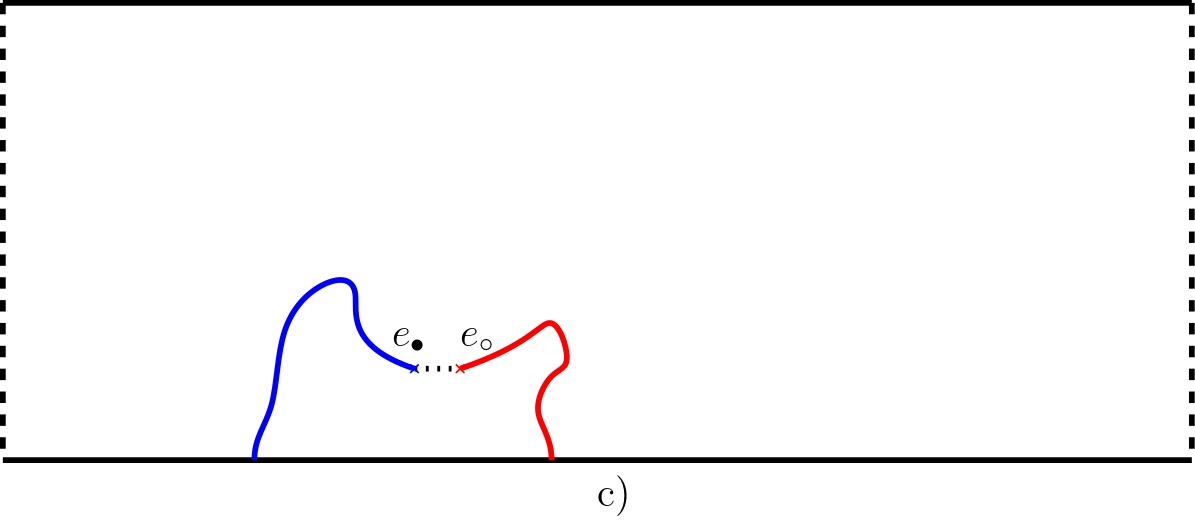} \quad
				\includegraphics[width=0.45\textwidth]{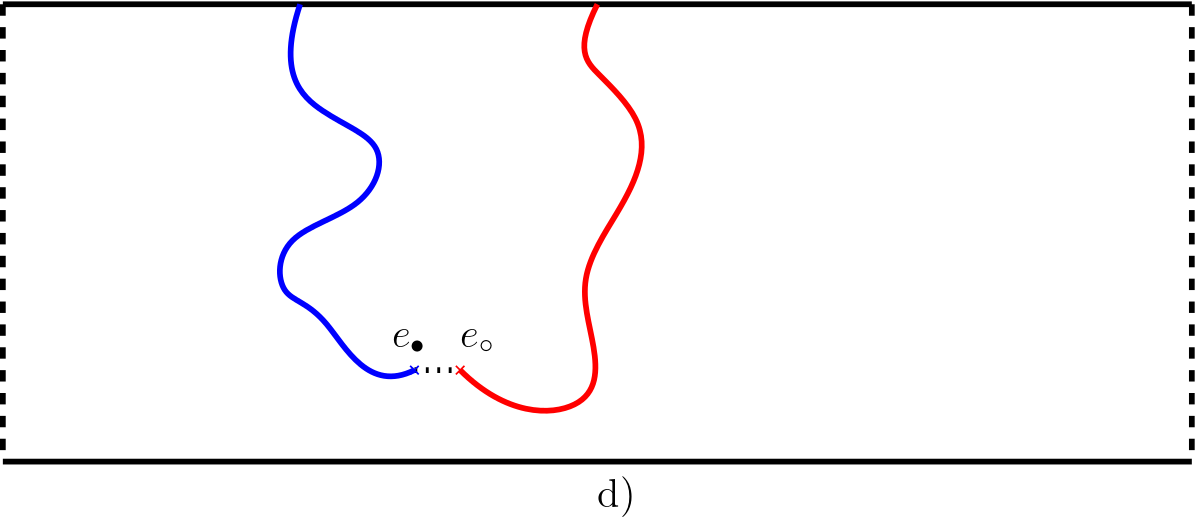}
				\\
				\caption{The illustrations for the cases a)-d).}
				\label{fig.classification}
			\end{figure}
		\end{proof}
		
		As a direct consequence, the two-arm event localizes the contribution of homogenized conductivity in $\a(\Lambda)$.
		\begin{corollary}
			The homogenized conductivity $\a(\Lambda)$ in \eqref{eq.def.aU} can be reformulated as
			\begin{align}\label{eq.defaa}
				\a(\Lambda) :=  \frac{1}{|\Lambda|}\sum_{e \in \Ed^+(\Lambda)} \a(e) \vert \nabla u_{\Lambda}(e) \vert^2 \Ind{\tarm(e, \Lambda)} .
			\end{align}
		\end{corollary}
		\begin{proof}
			The original expression of $\a(\Lambda)$ is 
			\begin{align*}
				\a(\Lambda) = \frac{1}{|\Lambda|}\sum_{e \in \Ed^+(\Lambda)} \a(e) \vert \nabla u_{\Lambda}(e) \vert^2.
			\end{align*}
			The percolation $\a$ restricts the sum on the open bonds $e = \{e_\bullet, e_\circ\}$. If the bond does not connect to the boundary $[\partial \Lambda]^{\pm}$, then $u_{\Lambda}$ is constant on the associated component. Therefore, the contribution comes from the bonds satisfying the condition \eqref{eq.boundary_connection}
			\begin{align*}
				\a(\Lambda) := \frac{1}{|\Lambda|}\sum_{e \in \Ed^+(\Lambda)} \a(e) \vert \nabla u_{\Lambda}(e) \vert^2 \Ind{\{e_\bullet \longleftrightarrow [\partial \Lambda]^{\pm}\} \cap \{e_\circ \longleftrightarrow \partial [\partial \Lambda]^{\pm}\}}.
			\end{align*}
			Furthermore, we notice only the bonds $\nabla u_{\Lambda}(e) \neq 0$ contributes, thus Lemma~\ref{lem.two_arm} applies and yields \eqref{eq.defaa}.
		\end{proof}
		
		\subsection{A pirori bound of conductivity}
		Using the observation of pivotal edges, we give an \emph{a priori} bound for the conductivity.

		\begin{lemma}\label{lem.aa_bound}
			In high-dimensional percolation, there exists a finite positive constant $C(d)$ such that
			\begin{align}\label{eq.aa_bound}
				\ab(\Lambda_{h, L}, \tilde{\p}) \leq \Ll\{\begin{array}{ll}
					C h^{-2} (\log{h}), & \tilde{\p} \in (\p_c, \p_c + h^{-2}];\\
					C h^{-2} (\log{h}) \exp\Ll(- C^{-1} h \sqrt{\p_c - \tilde{\p}}\Rr), & \tilde{\p} \in (0, \p_c].
				\end{array}\Rr.
			\end{align}
		\end{lemma}
		\begin{proof}
			We start from the reformulation in \eqref{eq.defaa} and keep using the notation $e = \{e_\bullet, e_\circ\}$
			\begin{equation}\label{eq.aa_bound_step1}
				\begin{split}
					\ab(\Lambda_{h, L}, \tilde{\p}) &=  \E_{\tilde{\p}}\Ll[ \frac{1}{|\Lambda_{h, L}|}\sum_{e \in \Ed^+(\Lambda_{h, L})} \a(e) \vert \nabla u_{\Lambda_{h, L}}(e) \vert^2 \Ind{\tarm(\Lambda_{h, L},e)}\Rr] \\
					&\leq \frac{h^2}{|\Lambda_{h, L}|}\sum_{e \in \Ed^+(\Lambda_{h, L})} \P_{\tilde{\p}}\Ll[ \tarm(\Lambda_{h, L},e) \Rr]\\
					&\leq \frac{2 h^2}{|\Lambda_{h, L}|}\sum_{e \in \Ed^+(\Lambda_{h, L})} \P_{\tilde{\p}}\Ll[ \{e_\bullet \longleftrightarrow [\partial \Lambda_{h, L}]^-\} \circ \{e_\circ \longleftrightarrow [\partial \Lambda_{h, L}]^+\} \Rr]\\
					&\leq \frac{2 h^2}{|\Lambda_{h, L}|}\sum_{e \in \Ed^+(\Lambda_{h, L})} \P_{\tilde{\p}}\Ll[e_\bullet \longleftrightarrow [\partial \Lambda_{h, L}]^-\Rr] \P_{\tilde{\p}}\Ll[e_\circ \longleftrightarrow [\partial \Lambda_{h, L}]^+\Rr].
				\end{split}
			\end{equation}
			From the first line to the second line, we make use of the maximum principle $\vert \nabla u_{\Lambda_{h, L}}(e) \vert \leq h$. From the third line to the fourth line, the BK inequality is applied.

			We then need the estimate of the one-arm event. Without loss of generality, we assume that $\dex(e, [\partial \Lambda_{h, L}]^{-}) = k \leq \frac{h}{2}$. Then we have an upper bound for $\P_{\tilde{\p}}\Ll[e_\bullet \longleftrightarrow [\partial \Lambda_{h, L}]^-\Rr]$ using the one-arm estimate in slightly supercritical regime in \eqref{eq.onearm_super}
			\begin{align}\label{eq.onearm_low}
				\P_{\tilde{\p}}\Ll[e_\bullet \longleftrightarrow [\partial \Lambda_{h, L}]^-\Rr] \leq \P_{\p_c+h^{-2}}\Ll[e_\bullet \longleftrightarrow \partial \Lambda_k(e_\bullet)\Rr] \lesssim \frac{1}{(k+1)^2}.
			\end{align}
			Concerning the probability $\P_{\tilde{\p}}\Ll[e_\circ \longleftrightarrow [\partial \Lambda_{h, L}]^+\Rr]$, as $h \ll L$ and $k$ can be very small, this event is similar to that in the half-space. We need a generalized version of \eqref{eq.onearm_half}, which is proved in Proposition~\ref{prop.onearm_half_plus}
			\begin{equation}\label{eq.onearm_up}
				\begin{split}
					\P_{\tilde{\p}}\Ll[e_\circ \longleftrightarrow [\partial \Lambda_{h, L}]^+\Rr] &\leq \P_{\tilde{\p}}\Ll[k e_1 \stackrel{\H}{\longleftrightarrow} \H_{-h}\Rr] \\
					&\lesssim \Ll\{\begin{array}{ll}
						(k+1) h^{-3} , & \tilde{\p} \in (\p_c, \p_c + h^{-2}];\\
						(k+1) h^{-3}  \exp\Ll(- C^{-1} h \sqrt{\p_c - \tilde{\p}}\Rr), & \tilde{\p} \in (0, \p_c].
					\end{array}\Rr.
				\end{split}
			\end{equation}
			
			Inserting estimates \eqref{eq.onearm_low} and \eqref{eq.onearm_up} back to \eqref{eq.aa_bound_step1}, we have the following calculation. Here we just take $\tilde{\p} \in (\p_c, \p_c + h^{-2}]$ for example, while the case $\tilde{\p} \in (0, \p_c]$ has an additional exponential factor
			\begin{equation}\label{eq.aa_bound_step2}
				\begin{split}
					\ab(\Lambda_{h, L}, \tilde{\p}) &\leq \frac{4h^2}{|\Lambda_{h,L}|}\sum_{\substack{e \in \Ed^+(\Lambda_{h,L}): \\
							\dex(e, [\partial \Lambda]^{-}) \leq \frac{h}{2}}}\P_{\tilde{\p}}\Ll[e_\bullet \longleftrightarrow [\partial \Lambda_{h, L}]^-\Rr] \P_{\tilde{\p}}\Ll[e_\circ \longleftrightarrow [\partial \Lambda_{h, L}]^+\Rr]\\
					&\lesssim \frac{4h^2}{|\Lambda_{h,L}|}\sum_{k=0}^{\lfloor \frac{h}{2} \rfloor} L^{d-1} \cdot (k+1)^{-2} \cdot (k+1) h^{-3} \\
					&\lesssim \frac{1}{|\Lambda_{h,L}|} L^{d-1} \cdot h^{-1} \cdot \log h\\
					&\lesssim h^{-2} \log h. 
				\end{split}
			\end{equation}
			This concludes the desired result \eqref{eq.aa_bound}.
		\end{proof}
		
		\begin{remark}
			Let us comment the bound in Lemma~\ref{eq.aa_bound}. We pick $\p \in (\p_c, 1]$ and correlation length $h = \ell_*(\p) := \lfloor(\p - \p_c)^{-1/2}\rfloor$ with a very large $L$. As we know from the the definition \eqref{eq.def.aU}, a very rough upper bound is obtained by testing $\ell_{e_1}$ in \eqref{eq.defmu}
			\begin{align*}
				\ab(\Lambda_{\ell_*,L},\tilde{\p}) &\leq \E_{\tilde{\p}}\Ll[\frac{1}{|\Lambda_{\ell_*,L}|}\sum_{e \in \Ed^+(\Lambda_{\ell_*,L})} \a(e) \vert \nabla \ell_{e_1}(e) \vert^2 \Rr]\\
				&\leq \frac{1}{|\Lambda_{\ell_*,L}|} \E_{\tilde{\p}}\Ll[\{x \in \Lambda_{\ell_*,L}: x \xleftrightarrow{\a} \partial \Lambda^+_{\ell_*,L}\}\Rr] \simeq \ell_*^{-1}.
			\end{align*}
			Thus, even when taking $\tilde{\p} < \p_c$, we do not see the decay in subcritical regime. This is because the calculation above counts the loops like in Lemma~\ref{lem.boundary} but lying near $[\partial \Lambda^+_{\ell_*,L}]^{\pm}$, which do not really contribute. The analysis in Lemma~\ref{lem.aa_bound} explains that the contribution comes from the clusters connecting two different faces, and \eqref{eq.aa_bound} yields 
			\begin{align}\label{eq.aa_bound_2}
				\ab(\Lambda_{\ell_*,L},\tilde{\p}) \lesssim   (\p-\p_c) \vert \log (\p - \p_c) \vert \exp\Ll(- C^{-1} \ell_* \sqrt{\p_c - \tilde{\p}}\Rr).
			\end{align}
			Then it captures the decay in the subcritical regime. Meanwhile, taking $\tilde{\p} = \p$, the bound in \eqref{eq.aa_bound} states that 
			\begin{align*}
				\ab(\Lambda_{\ell_*,L}, \p) \lesssim   (\p-\p_c) \vert \log (\p - \p_c) \vert.
			\end{align*}
			This is still relatively weak, and we will improve it in next section.
		\end{remark}

		\section{Upper bound II: The differential inequality}\label{sec.diff}
		The main result in this section is an upper bound for the conductivity.
		\begin{prop}\label{prop.upper_log}
			In high-dimensional percolation, for every $\p \in (\p_c, 1]$, we have
			\begin{align}\label{eq.upper_log}
				\ab(\p) \lesssim   (\p-\p_c)^3 \vert \log (\p - \p_c) \vert.  
			\end{align}
		\end{prop}
		We will study it via a differential inequality. The tools in Sections~\ref{subsec.Glauber},~\ref{subsec.diff} hold for the Bernoulli percolation in all dimensions $d \geq 2$, and the results of high-dimensional percolation only enter in Section~\ref{subsec.upper_log}. 
		
		\subsection{Glauber derivative}\label{subsec.Glauber}
		For every $e \in \Ed$, we define $\a^e$ to be the configuration opening the edge $e$ 
		\begin{align}\label{eq.defae}
			\a^e(e') = \Ll\{ \begin{array}{ll}
				\a(e') & e' \neq e;\\
				1 & e' = e.
			\end{array}\Rr.
		\end{align}
		For a measurable function $f: \{0,1\}^{\Ed} \to \R $, we define the shorthand notation
		\begin{align}\label{eq.deffe}
			f^e(\a) := f(\a^e).
		\end{align}
		
		The derivative of the $\ab$ can be expressed via Russo's formula, and the {\rhs} of \eqref{e:derivative_formula} can be seen as the Glauber derivative. We refer to \cite[Proposition~3.1]{gu2025diffusivity} for the proof of the following lemma, which actually applies to a more general case involving the higher-order Glauber derivatives.
		\begin{lemma}[Russo's formula]\label{pr:expansion_init}
			For every bounded cube $\Lambda \subset \Zd$, the mapping $\p \mapsto \ab(\Lambda,\p)$ is polynomial, and its first derivative $\ab^{(1)}(\Lambda,\p) := \frac{\d}{\d \p}\ab(\Lambda,\p)$ has an explicit expression
			\begin{equation}
				\label{e:derivative_formula}
				\forall \p \in (0,1), \qquad \ab^{(1)}(\Lambda,\p) =  \frac{1}{(1-\p)} \sum_{e \in \Ed^+(\Lambda)} \E_{\p}\Ll[\mathscr{E}_{\Lambda}(\a^e) - \mathscr{E}_{\Lambda}(\a) \Rr].
			\end{equation}
		\end{lemma}

		Opening the edge increases the conductivity. The following comparison argument will be quite useful. Especially, the identity in the second line of \eqref{eq.difference_comparison} allows to transform the difference of the Dirichlet energy to an inner product. This identity can be found in \cite[Lemma~3.1]{dg1} and \cite[Lemma~4.2]{GGMN}. We also attribute Wenhao Zhao for the version on percolation, which appeared in the early manuscript of \cite{gu2025diffusivity}. 
		\begin{lemma}\label{lem.comparison}
			For two configurations $\a \leq  \tilde{\a}$, we have $\mathscr{E}_{\Lambda}(\a) \leq  \mathscr{E}_{\Lambda}(\tilde{\a})$. Moreover, we have the estimate for their difference
			\begin{multline}\label{eq.difference_comparison}
				\frac{1}{|\Lambda|}\Ll \langle\nabla u_{\Lambda}(\tilde{\a}), (\tilde{\a} - \a) \nabla u_{\Lambda}(\tilde{\a})\Rr \rangle_{\Lambda} \\
				\leq \mathscr{E}_{\Lambda}(\tilde{\a}) - \mathscr{E}_{\Lambda}(\a) = \frac{1}{|\Lambda|}\Ll \langle\nabla u_{\Lambda}(\tilde{\a}), (\tilde{\a} - \a) \nabla u_{\Lambda}(\a)\Rr \rangle_{\Lambda} \\
				\leq \frac{1}{|\Lambda|}\Ll \langle\nabla u_{\Lambda}(\a), (\tilde{\a} - \a) \nabla u_{\Lambda}(\a)\Rr \rangle_{\Lambda}. 
			\end{multline}
		\end{lemma}
		\begin{proof}
			The comparison $\mathscr{E}_{\Lambda}(\a) \leq  \mathscr{E}_{\Lambda}(\tilde{\a})$ can be deduced from the first line of \eqref{eq.difference_comparison}. Using the definition, we have 
			\begin{align*}
				\mathscr{E}_{\Lambda}(\a) \leq \frac{1}{|\Lambda|}\Ll \langle\nabla u_{\Lambda}(\tilde{\a}), \a \nabla u_{\Lambda}(\tilde{\a})\Rr \rangle_{\Lambda}, \qquad \mathscr{E}_{\Lambda}(\tilde{\a}) = \frac{1}{|\Lambda|}\Ll \langle\nabla u_{\Lambda}(\tilde{\a}), \tilde{\a} \nabla u_{\Lambda}(\tilde{\a})\Rr \rangle_{\Lambda}.
			\end{align*}
			This is because $u_{\Lambda}(\tilde{\a})$ is not the minimizer in the variational problem of $\mathscr{E}_{\Lambda}(\a)$.
			Combining the two formulas will yield the first inequality in \eqref{eq.difference_comparison}. Similarly, we get the second inequality in \eqref{eq.difference_comparison} by using the following two formulas
			\begin{align*}
				\mathscr{E}_{\Lambda}(\tilde{\a}) \leq \frac{1}{|\Lambda|}\Ll \langle\nabla u_{\Lambda}(\a), \tilde{\a} \nabla u_{\Lambda}(\a)\Rr \rangle_{\Lambda}, \qquad \mathscr{E}_{\Lambda}(\a) = \frac{1}{|\Lambda|}\Ll \langle\nabla u_{\Lambda}(\a), \a \nabla u_{\Lambda}(\a)\Rr \rangle_{\Lambda}.
			\end{align*}  
			
			The identity in \eqref{eq.difference_comparison} comes from the harmonic equations. Viewing the definition \eqref{eq.def.aU}, for every $\varphi \equiv 0$ on $[\partial \Lambda]^{\pm}$, we have 
			\begin{align*}
				\Ll \langle\nabla \varphi, \a \nabla u_{\Lambda}(\a)\Rr \rangle_{\Lambda} = 0.
			\end{align*}
			Therefore, the quantity $\mathscr{E}_{\Lambda}(\tilde{\a})$ can be expressed with an arbitrary $\varphi \equiv 0$ on $[\partial \Lambda]^{\pm}$
			\begin{align*}
				\mathscr{E}_{\Lambda}(\a) = \frac{1}{|\Lambda|}\Ll \langle\nabla (\ell_{e_1} + \varphi), \a \nabla u_{\Lambda}(\a)\Rr \rangle_{\Lambda}.
			\end{align*}
			We take $\varphi = u_{\Lambda}(\tilde{\a}) - \ell_{e_1}$, which yields that 
			\begin{align}\label{eq.harmonic_1}
				\mathscr{E}_{\Lambda}(\a) = \frac{1}{|\Lambda|}\Ll \langle\nabla u_{\Lambda}(\tilde{\a}), \a \nabla u_{\Lambda}(\a)\Rr \rangle_{\Lambda}.
			\end{align}
			Similarly, we get another expression for $\mathscr{E}_{\Lambda}(\tilde{\a})$ that 
			\begin{align}\label{eq.harmonic_2}
				\mathscr{E}_{\Lambda}(\tilde{\a}) = \frac{1}{|\Lambda|}\Ll \langle\nabla u_{\Lambda}(\tilde{\a}), \tilde{\a} \nabla u_{\Lambda}(\a)\Rr \rangle_{\Lambda}.
			\end{align}
			The difference between  \eqref{eq.harmonic_1} and \eqref{eq.harmonic_2} is the identity in \eqref{eq.difference_comparison}. 
		\end{proof}

		The following result is a counterpart of Lemma~\ref{lem.two_arm}, which relates the Glauber derivative to the two-arm event.
		\begin{lemma}\label{lem.pivot_trick}
			Given an edge $e = \{e_\bullet, e_\circ\} \subset \Ed^+(\Lambda)$, then $\mathscr{E}_{\Lambda}(\a^e) \neq  \mathscr{E}_{\Lambda}(\a)$ implies that ${\a(e) = 0}$ and the two-arm event $\tarm(e, \Lambda)$ defined in \eqref{eq.two_arm}.
		\end{lemma}
		\begin{proof}
			
			The necessary condition $\a(e) = 0$ is obvious. Otherwise, $\a^e \equiv \a$ and $\mathscr{E}_{\Lambda}(\a^e) =  \mathscr{E}_{\Lambda}(\a)$. We then make the discussion case by case.
			
			\textit{Case~1: $e_\bullet \nleftrightarrow [\partial \Lambda]^{\pm}$ and $e_\circ \nleftrightarrow [\partial \Lambda]^{\pm}$.} Notice that the Dirichlet energy in $\mathscr{E}_{\Lambda}(\a)$ is determined by the clusters connecting to the boundary. For this case, then opening $e$ will not change the clusters connecting to the boundary, and we have $\mathscr{E}_{\Lambda}(\a^e) = \mathscr{E}_{\Lambda}(\a)$.
			
			\smallskip
			
			\textit{Case~2:  $e_\bullet \leftrightarrow [\partial \Lambda]^{\pm}$ and $e_\circ  \nleftrightarrow [\partial \Lambda]^{\pm}$.}
			If one endpoint connects to the boundary $[\partial \Lambda]^{\pm}$ and another does not, then opening $e$ will enlarge the clusters connecting to the boundary $[\partial \Lambda]^{\pm}$. Let $\cC$ be this enlargement. Then we observe that $u^e_{\Lambda}$ is obtained from $u_{\Lambda}$ via a constant extension on $\cC$ with the value $u_{\Lambda}(e_\bullet)$. This justifies again $\mathscr{E}_{\Lambda}(\a^e) = \mathscr{E}_{\Lambda}(\a)$ for such case.
			
			\smallskip
			
			\textit{Case~3: $e_\bullet \leftrightarrow [\partial \Lambda]^{\pm}$ and $e_\circ  \leftrightarrow [\partial \Lambda]^{\pm}$.} Applying Lemma~\ref{lem.comparison} to $\a^e$, we obtain that 
			\begin{align*}
				\mathscr{E}_{\Lambda}(\a^e) - \mathscr{E}_{\Lambda}(\a) = \frac{1}{|\Lambda|}\Ll \langle\nabla u^e_{\Lambda}, (\a^e - \a) \nabla u_{\Lambda}\Rr \rangle_{\Lambda}.
			\end{align*}
			The condition $\mathscr{E}_{\Lambda}(\a^e) \neq  \mathscr{E}_{\Lambda}(\a)$ then implies that $u_{\Lambda}(e_\bullet) \neq u_{\Lambda}(e_\circ)$. Thus the setting of Lemma~\ref{lem.two_arm} is satisfied, and we can further conclude $\tarm(e, \Lambda)$.
		\end{proof}

		\begin{figure}[h!]
			\centering
			\includegraphics[width=0.3\textwidth]{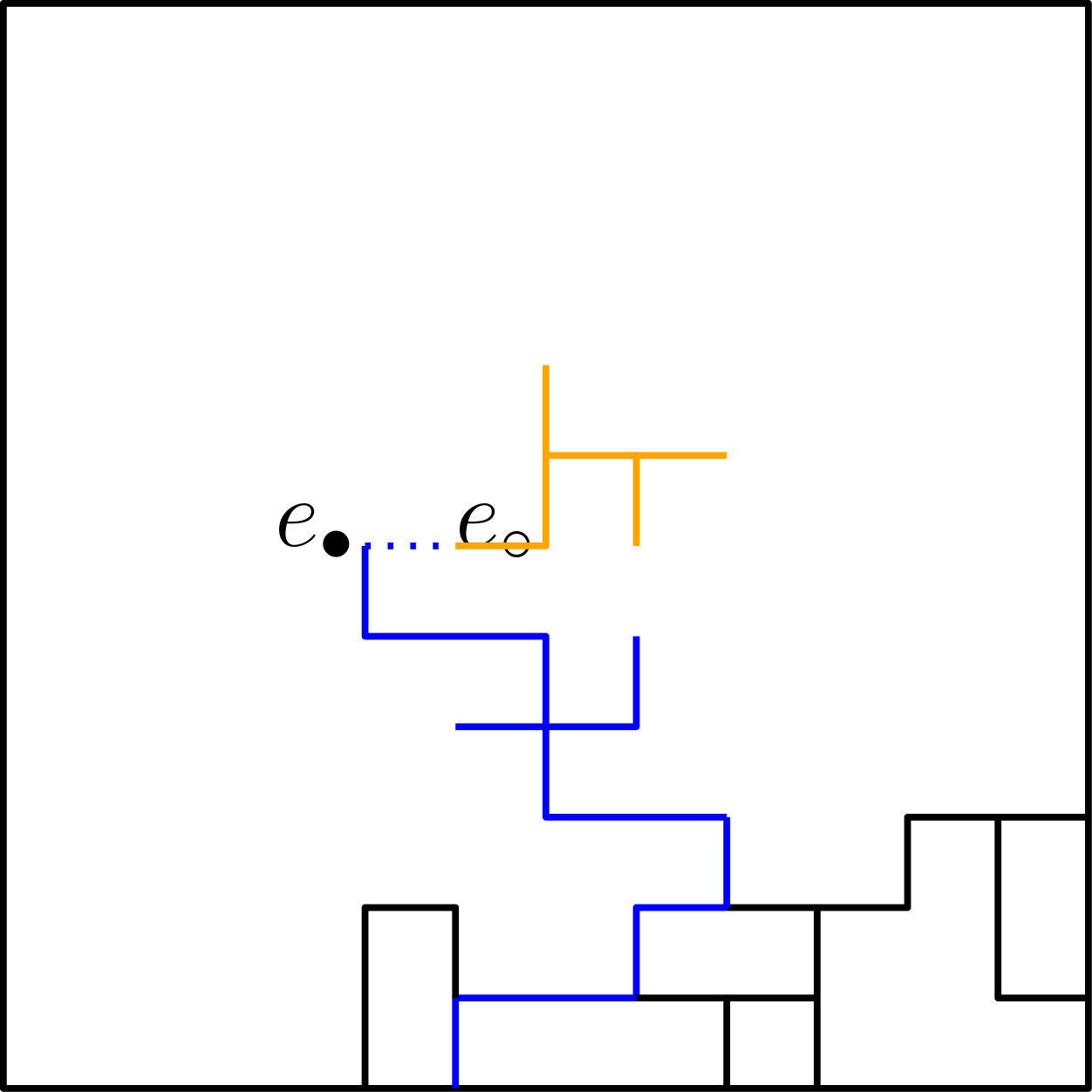} \hspace{0.5cm}
			\includegraphics[width=0.3\textwidth]{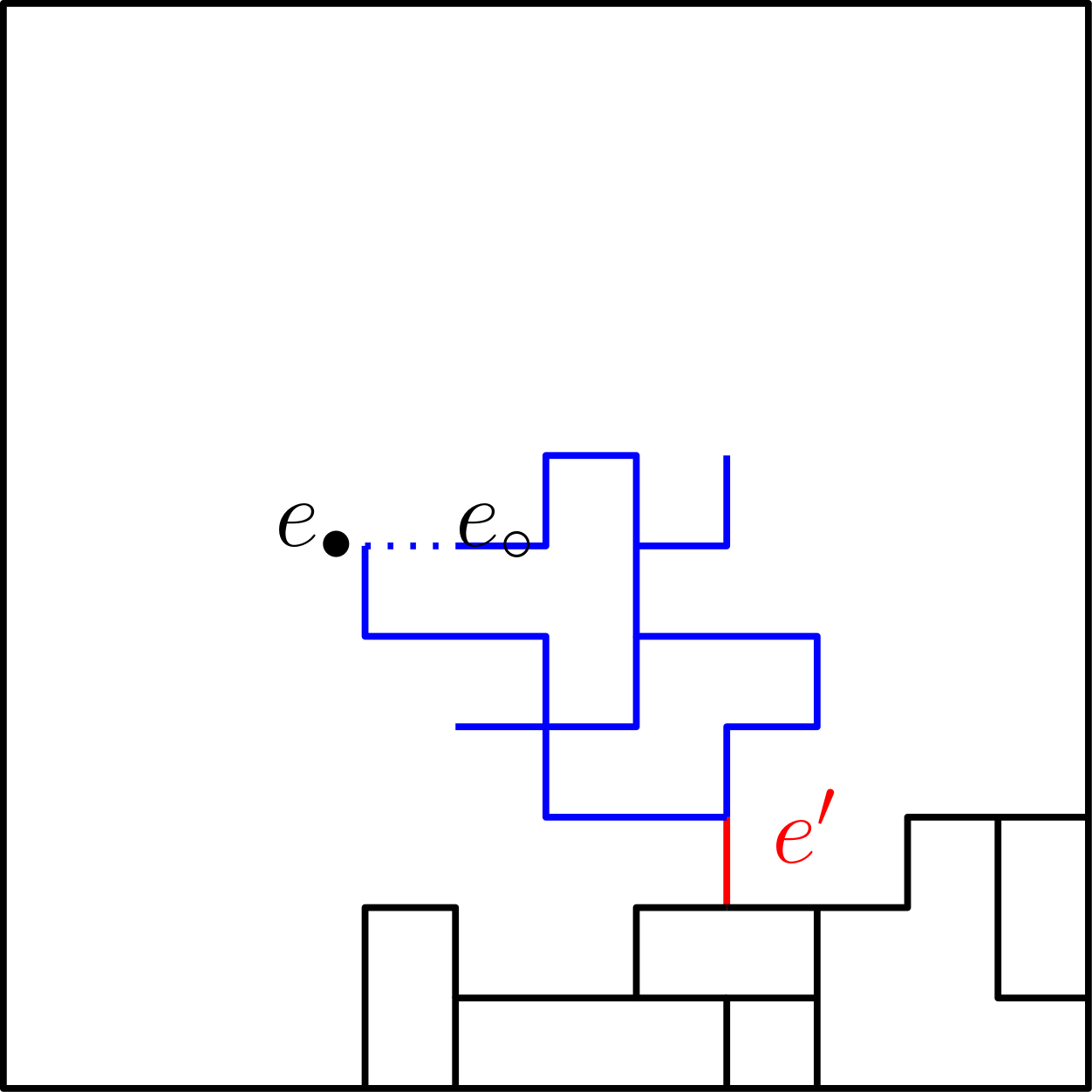} \hspace{0.5cm}
			\includegraphics[width=0.3\textwidth]{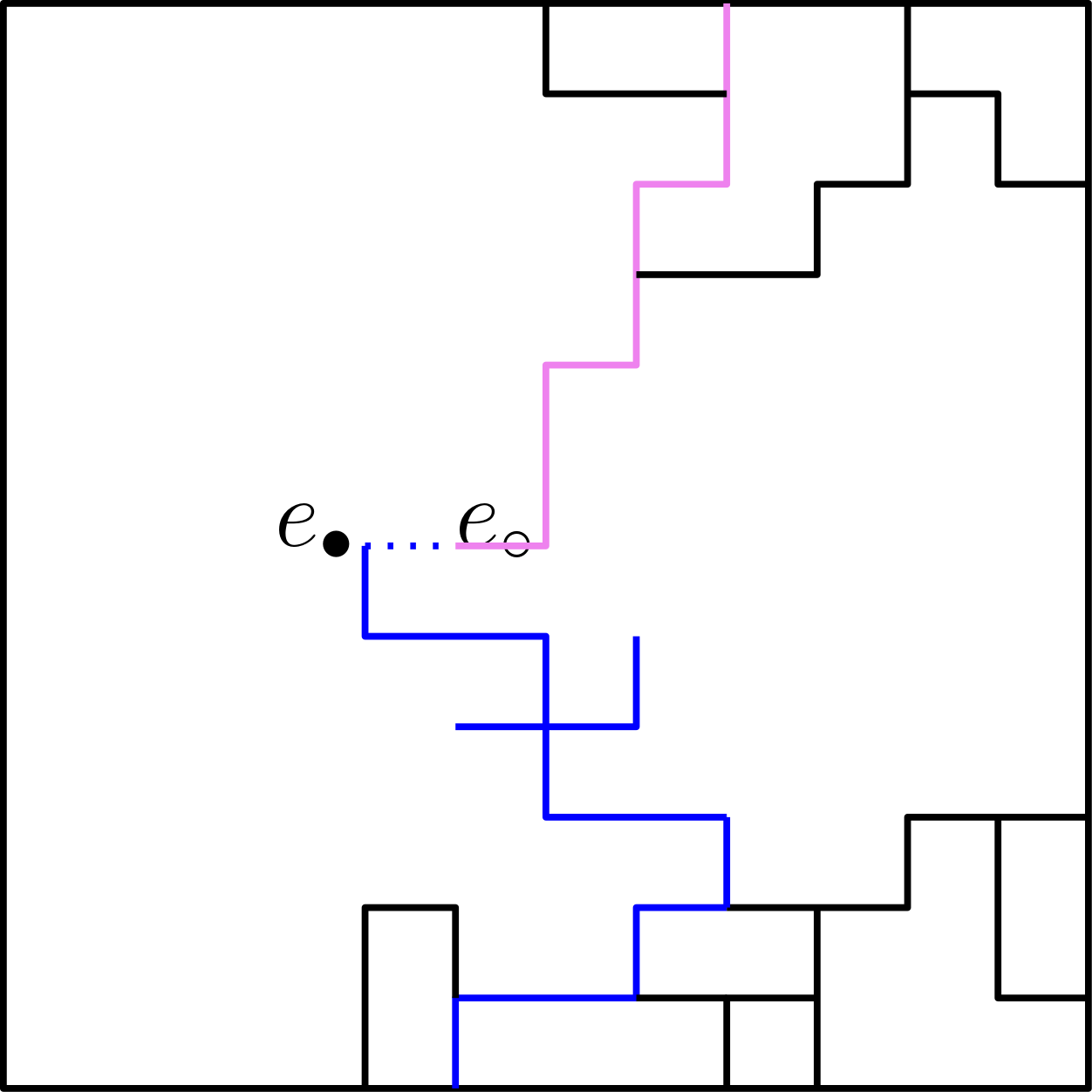}
			\caption{There are  $3$ cases in the perturbation, and only the third case with two disjoint paths contributes to the increment of Dirichlet energy.}
		\end{figure}

		We finish this part by introducing Mecke's identity. This identity is well-known in the analysis of point process, and we are inspired by \cite[Theorem~4.1]{bookpoisson}.
		\begin{lemma}[Mecke's identity]
			For every measurable function $F: \{0,1\}^{\Zd} \times E(\Zd) \to \R$, the following identity holds that
			\begin{align}\label{eq.Mecke}
				\E_{\p}\Ll[F(\a, e) \a(e)\Rr] = \p \E_{\p}\Ll[F(\a^e, e)\Rr].
			\end{align}
		\end{lemma}
		\begin{proof}
			We make the calculation that
			\begin{align*}
				\E_{\p}\Ll[F(\a, e) \a(e)\Rr] &= \E_{\p}\Ll[F(\a, e) \a(e) \Ind{\a(e)=1}\Rr] + \underbrace{\E_{\p}\Ll[F(\a, e) \a(e)  \Ind{\a(e)=0}\Rr]}_{=0} \\
				&= \E_{\p}\Ll[F(\a^e, e)  \Ind{\a(e)=1}\Rr]\\
				&= \E_{\p}\Ll[F(\a^e, e) \Rr] \P_{\p}[\a(e)=1]\\
				&= \p \E_{\p}\Ll[F(\a^e, e) \Rr].
			\end{align*}
			The first line is just the decomposition. From the second line to the third line, we make use of the fact that the configuration $\a^e$ is independent of $\a(e)$.
		\end{proof}

		\subsection{Differential inequality}\label{subsec.diff}
		We combine all the results in previous subsection and derive the following differential inequality.
		\begin{prop}\label{prop.inequality_diff}
			For every cube $\Lambda$ and all $\p \in (0,1)$, we have
			\begin{align}\label{eq.inequality_diff}
				\ab^{(1)}(\Lambda,\p) \leq  \frac{1}{\p^{\frac{1}{2}}(1-\p)} \ab^{\frac{1}{2}}(\Lambda,\p) \bb^{\frac{1}{2}}(\Lambda,\p).
			\end{align}
			Here $\bb(\Lambda,\p)$ is a positive quantity defined as 
			\begin{align}\label{eq.defbb}
				\bb(\Lambda,\p) := \E_{\p}\Ll[\frac{1}{|\Lambda|}\sum_{e \in \Ed^+(\Lambda)} (1-\a(e)) \vert \nabla u_{\Lambda}(e) \vert^2 \Ind{\tarm(e, \Lambda)} \Rr].
			\end{align}
		\end{prop}

		\begin{proof}
			We combine Russo's formula \eqref{e:derivative_formula} with the identity in \eqref{eq.difference_comparison}
			\begin{align*}
				\ab^{(1)}(\Lambda,\p) &\stackrel{\text{\eqref{e:derivative_formula}}}{=}  \frac{1}{(1-\p)} \sum_{e \in \Ed^+(\Lambda)} \E_{\p}\Ll[\mathscr{E}_{\Lambda}(\a^e) - \mathscr{E}_{\Lambda}(\a) \Rr]\\
				&\stackrel{\text{\eqref{eq.difference_comparison}}}{=} \frac{1}{(1-\p)} \sum_{e \in \Ed^+(\Lambda)} \E_{\p}\Ll[\frac{1}{|\Lambda|}\Ll \langle\nabla u_{\Lambda}(\a^e), (\a^e - \a) \nabla u_{\Lambda}(\a)\Rr \rangle_{\Lambda} \Rr]\\
				&\leq  (1- \p)^{-1} (\mathbf{I} \cdot \mathbf{II})^{\frac{1}{2}}.
			\end{align*}
			In the third line, we make use of the observation in Lemma~\ref{lem.pivot_trick} and Cauchy--Schwarz inequality, and obtain the two terms defined respectively as
			\begin{equation}\label{eq.def_12_aCS}
				\begin{split}
					\mathbf{I} &:=   \E_{\p}\Ll[\frac{1}{|\Lambda|}\sum_{e \in \Ed^+(\Lambda)}\Ll \langle\nabla u^e_{\Lambda}, (\a^e - \a) \nabla u^e_{\Lambda}\Rr \rangle_{\Lambda} \Ind{\tarm(e, \Lambda)} \Rr], \\
					\mathbf{II} &:= \E_{\p}\Ll[\frac{1}{|\Lambda|}\sum_{e \in \Ed^+(\Lambda)} \Ll \langle\nabla u_{\Lambda}, (\a^e - \a) \nabla u_{\Lambda}\Rr \rangle_{\Lambda} \Ind{\tarm(e, \Lambda)} \Rr]. 
				\end{split}
			\end{equation} 
			
			We notice that 
			\begin{align*}
				\forall e' \neq e, \qquad \a^e(e') = \a(e'). 
			\end{align*}
			Then the term $\mathbf{II}$ is exactly $\bb(\Lambda,\p)$ defined in \eqref{eq.defbb}, and we only need to treat the term $\mathbf{I}$, which can be simplified as 
			\begin{align*}
				\mathbf{I} &=   \E_{\p}\Ll[\frac{1}{|\Lambda|}\sum_{e \in \Ed^+(\Lambda)} (1 - \a(e)) \vert \nabla u^e_{\Lambda}(e)\vert^2 \Ind{\tarm(e, \Lambda)} \Rr]\\
				&\leq \E_{\p}\Ll[\frac{1}{|\Lambda|}\sum_{e \in \Ed^+(\Lambda)}  \vert \nabla u^e_{\Lambda}(e)\vert^2 \Rr]\\
				&\stackrel{\text{\eqref{eq.Mecke}}}{=} \p^{-1}\E_{\p}\Ll[\frac{1}{|\Lambda|}\sum_{e \in \Ed^+(\Lambda)} \a(e) \vert \nabla u_{\Lambda}(e)\vert^2 \Rr]\\
				&= \p^{-1} \ab(\Lambda,\p).
			\end{align*}
			From the first line to the second line, we use the fact $1-\a(e) \leq 1$ and drop the indicator. The second line to the third is an application of Mecke's identity \eqref{eq.Mecke} to $F(\a, e) = \vert \nabla u_{\Lambda}(\a, e)\vert^2$, and the last line is the definition. This confirms the desired result \eqref{eq.inequality_diff}.
		\end{proof}

		The study of the term $\bb(\Lambda_{h, L}, \tilde{\p})$ is also an important step to understand $\ab(\Lambda_{h, L}, \tilde{\p})$. It was known in \cite{gu2025diffusivity} this quantity is uniformly bounded with respect to the domain $U$, but we are more interested in its dependence on the second parameter $\p$. 
		\begin{lemma}\label{lem.bb_bound}
			In high-dimensional percolation, for every $\p \in (\p_c, 1]$, then there exists a finite positive constant $C(d)$ such that
			\begin{align}\label{eq.bb_bound}
				\bb(\Lambda_{h, L}, \tilde{\p}) \leq \Ll\{\begin{array}{ll}
					C h^{-2} (\log{h}) , & \tilde{\p} \in (\p_c, \p_c + h^{-2}];\\
					C h^{-2} (\log{h}) \exp\Ll(- C^{-1} h \sqrt{\p_c - \tilde{\p}}\Rr)  , & \tilde{\p} \in (0, \p_c].
				\end{array}\Rr.
			\end{align}
		\end{lemma}
		\begin{proof}
			Viewing the definition \eqref{eq.defbb} and \eqref{eq.defaa}, the proof follows that of Lemma~\ref{lem.aa_bound}, mutatis mutandis.
		\end{proof}

		\subsection{Proof of upper bound with logarithm correction}\label{subsec.upper_log}
		
		\begin{proof}[Proof of Proposition~\ref{prop.upper_log}]
			Given $\p \in (\p_c, 1]$, we fix the correlation length $\ell_* = \lfloor(\p - \p_c)^{-1/2}\rfloor$ and $L = N \ell_*$, with $N$ to be chosen. Using the subadditivity, we obtain a comparison in \eqref{eq.chain.2}
			\begin{align}\label{eq.upper_step_1}
				\ab(\tilde{\p}) \leq \ab(\Lambda_{\ell_*,L},\tilde{\p}) + \ell_*^2 L^{-1},
			\end{align}
			and then we have a differential inequality from \eqref{eq.inequality_diff} 
			\begin{align}\label{eq.upper_step_2}
				\forall s \in (0, \p), \qquad \ab^{(1)}(\Lambda_{\ell_*,L},s) \leq  \frac{1}{s^{\frac{1}{2}}(1-s)} \ab^{\frac{1}{2}}(\Lambda_{\ell_*,L},s) \bb^{\frac{1}{2}}(\Lambda_{\ell_*,L},s).
			\end{align}
			Rewriting \eqref{eq.upper_step_2} as a bound on the derivative of $\ab(\Lambda_{\ell_*,L},s)^{\frac{1}{2}}$ and integrating over the neighborhood $[\p',\p]$ of $\p_c$, we  obtain that
			\begin{align}\label{eq.upper_step_2plus}
				\ab^{\frac{1}{2}}(\Lambda_{\ell_*,L}, \p) \lesssim \ab^{\frac{1}{2}}(\Lambda_{\ell_*,L}, \p') + \int_{\p'}^{\p}  \bb^{\frac{1}{2}}(\Lambda_{\ell_*,L},s) \, \d s.
			\end{align}
			Here we should choose $\p'$ close to $\p_c$, so that the factor $\frac{1}{s^{\frac{1}{2}}(1-s)}$ in \eqref{eq.upper_step_2} can be bounded uniformly by a constant only depending on $d$. We then need to clarify the contribution of two parts in \eqref{eq.upper_step_2plus}.
			
			\smallskip
			\emph{Step~1: contribution of $\int_{\p'}^{\p}  \bb^{\frac{1}{2}}(\Lambda_{\ell_*,L},s) \, \d s$.} We divide the integral into supercritical regime and subcritical regime.

			In the supercritical regime $s \in (\p_c, \p]$, because $s \leq \p_c + \ell^{-2}_*$, we insert the estimate of $\bb$ in \eqref{eq.bb_bound}
			\begin{align*}
				\forall s \in (\p_c, \p], \qquad \bb^{\frac{1}{2}}(\Lambda_{\ell_*,L},s) \lesssim  \ell^{-1}_* \log^{\frac{1}{2}} \ell_*.
			\end{align*}
			An integration from $\p_c$ to $\p$ yields that
			\begin{align}\label{eq.upper_step_3b_sup}
				\int_{\p_c}^{\p}  \bb^{\frac{1}{2}}(\Lambda_{\ell_*,L},s) \, \d s \lesssim  (\p-\p_c) \ell^{-1}_* \log^{\frac{1}{2}} \ell_*.
			\end{align}
			
			In the subcritical regime, since the bound $\bb$ in \eqref{eq.bb_bound} changes, we have a different differential inequality 
			\begin{align*}
				\forall s \in (0, \p_c], \qquad \bb^{\frac{1}{2}}(\Lambda_{\ell_*,L},s) \lesssim \ell^{-1}_* \log^{\frac{1}{2}} \ell_* \exp\Ll(- C^{-1} \ell_* \sqrt{\p_c - s}\Rr),
			\end{align*}
			which yields that 
			\begin{align*}
				\int_{\p'}^{\p_c}  \bb^{\frac{1}{2}}(\Lambda_{\ell_*,L},s) \, \d s \lesssim \int_{\p'}^{\p_c} \ell^{-1}_* \log^{\frac{1}{2}} \ell_*  \exp\Ll(- C^{-1} \ell_* \sqrt{\p_c - s}\Rr) \, \d s.
			\end{align*}
			The extra exponential factor in this estimate allows us to make the integral beyond the window of the correlation length. To see this, let us make a change of variable $s = \p_c - \ell_*^{-2} t^2$ with $t \in [0, \ell_* (\p_c - \p')^{\frac{1}{2}}]$, and we obtain that
			\begin{equation}\label{eq.upper_step_3b_sub}
				\begin{split}
					\int_{\p'}^{\p_c}  \bb^{\frac{1}{2}}(\Lambda_{\ell_*,L},s) \, \d s  &\lesssim  \int_{0}^{\ell_*(\p_c - \p')^{\frac{1}{2}}} \ell^{-3}_* \log^{\frac{1}{2}} \ell_*  t \exp\Ll(- C^{-1} t\Rr)  \, \d t\\
					&\lesssim \ell^{-3}_* \log^{\frac{1}{2}} \ell_*.
				\end{split}
			\end{equation}
			Combining \eqref{eq.upper_step_3b_sup} and \eqref{eq.upper_step_3b_sub}, we obtain that 
			\begin{align}\label{eq.upper_step_3b}
				\int_{\p'}^{\p}  \bb^{\frac{1}{2}}(\Lambda_{\ell_*,L},s) \, \d s  \lesssim  \ell^{-3}_* \log^{\frac{1}{2}} \ell_*.
			\end{align}

			\emph{Step~2: contribution of $\ab^{\frac{1}{2}}(\Lambda_{\ell_*,L}, \p')$.}
			It suffices to pick $\p' = \p_c - \tilde{C} \ell_*^{-2}(\log \ell_*)^2$ with a very large $\tilde{C}$, then Lemma~\ref{lem.aa_bound} entails that
			\begin{equation}\label{eq.upper_step_3a}
				\begin{split}
					\ab(\Lambda_{\ell_*,L},\p')  	&\lesssim \ell_*^{-2} (\log{\ell_*}) \exp\Ll(- C^{-1} \tilde{C}^{\frac{1}{2}} \log \ell_*\Rr)\\
					&\lesssim o(\ell_*^{-6}).
				\end{split}
			\end{equation}
			
			\smallskip
			Finally, we combine \eqref{eq.upper_step_1}, \eqref{eq.upper_step_2plus}, \eqref{eq.upper_step_3b}  with \eqref{eq.upper_step_3a}, with the choice of  ${\p' = \p_c - \tilde{C} \ell_*^{-2}(\log \ell_*)^2}$, and obtain that 
			\begin{align*}
				\ab(\p) \lesssim  \ell^{-6}_* \log \ell_* +  \ell^2_*  L^{-1}
			\end{align*}
			Recall $\ell_* = \lfloor(\p - \p_c)^{-1/2}\rfloor$, and take $L \to \infty$. Then $\ell^{-6}_* \log \ell_*$ is the main contribution, while the other terms are negligible. Therefore, we conclude the  result \eqref{eq.upper_log}.
		\end{proof}
		
		\section{Upper bound III: Removal of the logarithm}\label{sec.log}
		In this section, we refine the estimate and remove the logarithm correction in Proposition~\ref{prop.upper_log}.
		\begin{prop}\label{prop.upper}
			In high-dimensional percolation, for every $\p \in (\p_c, 1]$, we have
			\begin{align}\label{eq.upper}
				\ab(\p) \lesssim   (\p-\p_c)^3.  
			\end{align}
		\end{prop}
		Let us revisit the proof of Proposition~\ref{prop.upper_log}. The key formulas are \eqref{eq.upper_step_1} and \eqref{eq.upper_step_2plus}
		\begin{align*}
			\ab^{\frac{1}{2}}(\p) \leq \ab^{\frac{1}{2}}(\Lambda_{\ell_*,L}, \p) + \ell_* L^{-\frac{1}{2}} \lesssim \ab^{\frac{1}{2}}(\Lambda_{\ell_*,L}, \p') + \int_{\p'}^{\p}  \bb^{\frac{1}{2}}(\Lambda_{\ell_*,L},s) \, \d s + \ell_* L^{-\frac{1}{2}}.
		\end{align*}
		The philosophy is to shift $\ab^{\frac{1}{2}}(\Lambda_{\ell_*,L}, \p)$ to the subcritical regime, and then to clarify that 
		\begin{enumerate}
			\item the contribution $\ab^{\frac{1}{2}}(\Lambda_{\ell_*,L}, \p')$ at subcritical regime is very small; 
			\item the perturbation $\int_{\p'}^{\p}  \bb^{\frac{1}{2}}(\Lambda_{\ell_*,L},s) \, \d s$ is of correct order.
		\end{enumerate}
		We also use the estimates in Lemmas~\ref{lem.aa_bound}, ~\ref{lem.bb_bound} as inputs. Although these two lemmas were proved exactly by the same strategy, their final objects are different as mentioned above. Therefore, we can make a better balance between them. That is, we slightly reduce the order of the perturbation involving $\bb$, and still keep the term $\ab$ negligible at some point $\p'<\p_c$.
		
		A detailed verification of the proof of Lemmas~\ref{lem.aa_bound}, ~\ref{lem.bb_bound} confirmed the idea above. Actually, the logarithm correction only occurs in \eqref{eq.aa_bound_step2} there, which comes from the boundary layer near the top or bottom. We can tolerate such correction for $\ab$ at $\p'$, since it is always a small quantity; see \eqref{eq.upper_step_3a}. Nevertheless, we let the term $\bb$ only count the perturbation in the bulk sense, thus there \eqref{eq.aa_bound_step2} will no longer have the logarithm correction.
		
		For the completeness, we formulate the discussion above as a mathematical proof in Section~\ref{subsec.bulk}. A technical one-arm estimate in inhomogeneous percolation will be addressed in Section~\ref{subsec.inhomo}.
		
		\subsection{Perturbation in bulk}\label{subsec.bulk}
		Let us define the bulk of $\Lambda_{h,L}$ with a parameter $\delta \in (0,1)$
		\begin{align}\label{eq.cube.delta}
			\Lambda_{h,L}^\delta := \Ll\{x \in \Lambda_{h,L}: \dex(x, [\Lambda_{h,L}]^{\pm}) > \frac{\delta h}{2}\Rr\}.
		\end{align}
		We construct another probability space $\P_{\p, \tilde{\p}}^{\delta}$ that samples independent Bernoulli random variables of parameter $\tilde{\p}$ for $\Ed(\Lambda_{h,L}^\delta)$, and independent Bernoulli random variable of parameters $\p$ for the other bonds in $\Ed$, i.e. 
		\begin{align}\label{eq.def.Pdelta}
			\P_{\p, \tilde{\p}}^{\delta} := \mathrm{Bernoulli}(\p)^{\otimes \Ed(\Lambda_{h,L}^\delta)} \otimes \mathrm{Bernoulli}(\tilde{\p})^{\otimes (\Ed \setminus \Ed(\Lambda_{h,L}^\delta))}.
		\end{align}
		Let $\E_{\p, \tilde{\p}}^{\delta}$ be the associated expectation. We then follow \eqref{eq.def.aU} and \eqref{eq.def.abU}, and define the homogenized conductivity 
		\begin{align}\label{eq.def.aU.delta}
			\ab^\delta(\Lambda_{h,L},\p; \tilde{\p}) = \E_{\p, \tilde{\p}}^{\delta}[\a(\Lambda_{h,L})].
		\end{align}
		This quantity is related to the original one \eqref{eq.A_anneal} as 
		\begin{align}\label{eq.proba.identity}
			\P_{\p, \p}^{\delta} = \P_{\p} \qquad \Longrightarrow \qquad \ab^\delta(\Lambda_{h,L},\p; \p) = \ab(\Lambda_{h,L},\p).
		\end{align}
		The derivatives over the parameter $\tilde{\p}$ is then transformed to the perturbation over $\Lambda_{h,L}^\delta$, which yields a counterpart of Proposition~\ref{prop.inequality_diff}. 
		
		\begin{prop}\label{prop.inequality_diff_delta}
			For every $\p \in (0,1)$, we have
			\begin{align}\label{eq.inequality_diff_delta}
				\frac{\d}{\d \tilde{\p}} \ab^\delta(\Lambda_{h,L},\p; \tilde{\p}) \leq  \frac{1}{\tilde{\p}^{\frac{1}{2}}(1-\tilde{\p})} \ab^\delta(\Lambda_{h,L},\p; \tilde{\p})^{\frac{1}{2}} \bb^\delta(\Lambda_{h,L},\p; \tilde{\p})^{\frac{1}{2}}.
			\end{align}
			Here $ \bb^\delta(\Lambda_{h,L},\p; \tilde{\p})$ is a positive quantity defined as 
			\begin{align}\label{eq.defbb.delta}
				\bb^\delta(\Lambda_{h,L},\p; \tilde{\p})  := \E_{\p, \tilde{\p}}^{\delta}\Ll[\frac{1}{|\Lambda_{h,L}|}\sum_{e \in \Ed(\Lambda_{h,L}^{\delta})} (1-\a(e)) \vert \nabla u_{\Lambda_{h,L}}(e) \vert^2  \Ind{\tarm(e, \Lambda_{h,L})}\Rr].
			\end{align}
		\end{prop}
		\begin{proof}
			The main argument is the following calculation
			\begin{align*}
				&\frac{\d}{\d \tilde{\p}} \ab^\delta(\Lambda_{h,L},\p; \tilde{\p}) \\
				&=  \frac{1}{(1-\p)} \sum_{e \in \Ed(\Lambda_{h,L}^{\delta})} \E_{\p, \tilde{\p}}^{\delta}\Ll[\mathscr{E}_{\Lambda_{h,L}}(\a^e) - \mathscr{E}_{\Lambda_{h,L}}(\a) \Rr]\\
				&= \frac{1}{(1-\p)} \sum_{e \in \Ed(\Lambda_{h,L}^{\delta})} \E_{\p, \tilde{\p}}^{\delta}\Ll[\frac{1}{|\Lambda_{h,L}|}\Ll \langle\nabla u_{\Lambda_{h,L}}(\a^e), (\a^e - \a) \nabla u_{\Lambda_{h,L}}(\a)\Rr \rangle_{\Lambda_{h,L}} \Rr]\\
				&= \frac{1}{(1-\p)} \sum_{e \in \Ed(\Lambda_{h,L}^{\delta})} \E_{\p, \tilde{\p}}^{\delta}\Ll[\frac{1}{|\Lambda_{h,L}|}\Ll \langle\nabla u_{\Lambda_{h,L}}(\a^e), (\a^e - \a) \nabla u_{\Lambda_{h,L}}(\a)\Rr \rangle_{\Lambda_{h,L}} \Ind{\tarm(e, \Lambda_{h,L})} \Rr]
			\end{align*}
			The second line makes use of Russo's formula in Lemma~\ref{pr:expansion_init}. Because the derivative is acted on $\tilde{\p}$, only the bonds in $\Ed(\Lambda_{h,L}^{\delta})$ are involved; see the definition of $\P_{\p, \tilde{\p}}$ in \eqref{eq.def.Pdelta}. The third line follows Lemma~\ref{lem.comparison}, which is deterministic. The analysis of pivotal edges is then applied. Lemma~\ref{lem.pivot_trick}, which is also deterministic, suggests an indicator of two-arm event in the forth line.  The rest of analysis is similar as \eqref{eq.def_12_aCS} and we obtain the two terms in \eqref{eq.inequality_diff_delta}.
		\end{proof}

		We then need an upper bound for $\ab^\delta$ and $\bb^\delta$. One key estimate needed is a version of Proposition~\ref{prop.onearm_half_plus} under $\P_{\p, \tilde{\p}}^{\delta}$.
		\begin{prop}\label{prop.onearm_half_plus_delta}
			In high-dimensional percolation, there exists two finite positive constants $c_*(d,\delta), C(d,\delta)$, such that the following estimate holds 
			for every $0 \leq k \leq \frac{h}{2}$ and ${\p \in (\p_c + c_* h^{-2}]}$
			\begin{multline}\label{eq.onearm_half_plus_delta}
				\P_{\p, \tilde{\p}}^{\delta}\Ll[k e_1 \longleftrightarrow  [ \partial \Lambda_{h,L}]^{+}\Rr] \\ 
				\leq \Ll\{\begin{array}{ll}
					C(k+1)h^{-3}, & \tilde{\p} \in (\p_c + c_* h^{-2}], \\
					C(k+1)h^{-3} \exp\Ll(-C^{-1} h \sqrt{\p_c - \tilde{\p}}\Rr), & \tilde{\p} \in (0,\p_c].
				\end{array}\Rr.
			\end{multline}
		\end{prop}
		The case $\tilde{\p} \in (\p_c + c_* h^{-2}]$ is trivial, as it just needs a relaxation of $\p$ and $\tilde{\p}$ to $\in (\p_c + c_* h^{-2}]$ and the monotonicity of event. Concerning the case $\tilde{\p} \in  (0,\p_c]$, the estimate above is heuristically valid, because the open path has to cross a band of width $(1-\delta)h$ of subcritical phase. Meanwhile, we should also notice the potential difficulty, as the boundary is supercritical and brings additional chance. The rigorous proof of Proposition~\ref{prop.onearm_half_plus_delta} requires more efforts and is technical, thus we postpone it to Section~\ref{subsec.inhomo}.

		Using this key one-arm probability, we develop the counterpart of Lemmas~\ref{lem.aa_bound} and ~\ref{lem.bb_bound}.
		\begin{lemma}\label{lem.aabb_bound_new}
			In high-dimensional percolation, for the constant $c_*$ in Proposition~\ref{prop.onearm_half_plus_delta} and  every ${\p \in (\p_c, \p_c + c_* h^{-2}]}$, then there exists a finite positive constant $C(d, \delta)$ such that
			\begin{equation}\label{eq.bb_bound_new}
				\bb^\delta(\Lambda_{h,L},\p; \tilde{\p})  \leq \Ll\{\begin{array}{ll}
					C h^{-2} , & \tilde{\p} \in (\p_c, \p_c + c_* h^{-2}];\\
					C h^{-2}  \exp\Ll(- C^{-1} h \sqrt{\p_c - \tilde{\p}}\Rr)  , & \tilde{\p} \in (0, \p_c].
				\end{array}\Rr.
			\end{equation}
			and
			\begin{equation}\label{eq.aa_bound_new}
				\ab^\delta(\Lambda_{h,L},\p; \tilde{\p})  \\
				\leq \Ll\{\begin{array}{ll}
					C h^{-2} (\log{h}) , & \tilde{\p} \in (\p_c, \p_c + c_* h^{-2}];\\
					C h^{-2} (\log{h}) \exp\Ll(- C^{-1}  h \sqrt{\p_c - \tilde{\p}}\Rr)  , & \tilde{\p} \in (0, \p_c].
				\end{array}\Rr.
			\end{equation}
		\end{lemma}

		\begin{proof}
			For the term $\ab^\delta(\Lambda_{h,L},\p; \tilde{\p})$, it suffices to replace \eqref{eq.onearm_half_plus} by \eqref{eq.onearm_half_plus_delta}. Then we obtain exactly the same estimate.
			
			Concerning the term $\bb^\delta(\Lambda_{h,L},\p; \tilde{\p})$, since the estimate for $\tilde{\p} \in (\p_c, \p_c + c_* h^{-2}]$ in \eqref{eq.onearm_half_plus_delta} is the same as \eqref{eq.onearm_half_plus}, the corresponding upper bound does not change. It only remains to clarify the part $\tilde{\p} \in  (0,\p_c]$ in \eqref{eq.bb_bound_new}. We follow the proof of Proposition~\ref{prop.onearm_half_plus} starting from \eqref{eq.aa_bound_step2}
			\begin{align*}
				\bb^\delta(\Lambda_{h,L},\p; \tilde{\p}) \leq \frac{2h^2}{|\Lambda_{h,L}|}\sum_{\substack{e \in \Ed(\Lambda_{h,L}^{\delta})}}\P_{\p, \tilde{\p}}^{\delta}\Ll[e_\bullet \longleftrightarrow [\partial \Lambda_{h, L}]^-\Rr] \P_{\p, \tilde{\p}}^{\delta}\Ll[e_\circ \longleftrightarrow [\partial \Lambda_{h, L}]^+\Rr].
			\end{align*}
			For $\dex(e, [\partial \Lambda_{h, L}]^{-}) = k \leq \frac{h}{2}$, the one-arm probability of $\P_{\p, \tilde{\p}}^{\delta}\Ll[e_\bullet \longleftrightarrow [\partial \Lambda_{h, L}]^-\Rr]$ follows the the monotonicity and \eqref{eq.onearm_super}
			\begin{align}\label{eq.onearm_low_delta}
				\P_{\p, \tilde{\p}}^{\delta}\Ll[e_\bullet \longleftrightarrow [\partial \Lambda_{h, L}]^-\Rr] \leq \P_{\p_c + c_* h^{-2}}\Ll[e_\bullet \longleftrightarrow \Lambda_k(e_\bullet)\Rr] \lesssim \frac{1}{(k+1)^2}.
			\end{align}
			The probability of $\P_{\tilde{\p}}\Ll[e_\circ \longleftrightarrow [\partial \Lambda_{h, L}]^+\Rr]$ relies on the part $\tilde{\p} \in (0, \p_c]$ in \eqref{eq.onearm_half_plus_delta}
			\begin{equation}\label{eq.onearm_up_delta}
				\P_{\p, \tilde{\p}}^{\delta}\Ll[e_\circ \longleftrightarrow [\partial \Lambda_{h, L}]^+\Rr] \leq \P_{\p, \tilde{\p}}^{\delta} \Ll[k e_1 \stackrel{\H}{\longleftrightarrow} \H_{-h}\Rr] \leq C(k+1)h^{-3} \exp\Ll(-C^{-1} h \sqrt{\p_c - \tilde{\p}}\Rr).
			\end{equation}
			These estimates entail that 
			\begin{equation}\label{eq.bb_bound_step2}
				\begin{split}
					\bb^\delta(\Lambda_{h,L},\p; \tilde{\p}) &\leq \frac{4h^2}{|\Lambda_{h,L}|}\sum_{\substack{e \in \Ed(\Lambda_{h,L}^{\delta}): \\
							\dex(e, [\partial \Lambda]^{-}) \leq \frac{h}{2}}}\P_{\p, \tilde{\p}}^{\delta}\Ll[e_\bullet \longleftrightarrow [\partial \Lambda_{h, L}]^-\Rr] \P_{\p, \tilde{\p}}^{\delta}\Ll[e_\circ \longleftrightarrow [\partial \Lambda_{h, L}]^+\Rr]\\
					&\lesssim \frac{4h^2}{|\Lambda_{h,L}|}\sum_{k=\lfloor \frac{\delta h}{2} \rfloor}^{\lfloor \frac{h}{2} \rfloor} L^{d-1} \cdot (k+1)^{-2} \cdot (k+1) h^{-3}  \exp\Ll(- C^{-1}  h \sqrt{\p_c - \tilde{\p}}\Rr)\\
					&\lesssim \frac{1}{|\Lambda_{h,L}|} L^{d-1} \cdot h^{-1} \cdot \Ll( \log(h/2) - \log(\delta h /2) \Rr) \exp\Ll(- C^{-1}  h \sqrt{\p_c - \tilde{\p}}\Rr)\\
					&\lesssim h^{-2} \exp\Ll(- C^{-1}  h \sqrt{\p_c - \tilde{\p}}\Rr). 
				\end{split}
			\end{equation}
			The notation $\lesssim$ in this calculation is up to a constant depending on $d$ and $\delta$. The major difference is that the contribution all comes from the bonds in $\Ed(\Lambda_{h,L}^{\delta})$, which have a distance at least $\frac{\delta h }{2}$ from the boundary. As a consequence, the logarithm correction disappears here. This concludes that of \eqref{eq.bb_bound_new}.
		\end{proof}
		
		With the help of the new bounds $\ab^\delta, \bb^\delta$, the differential inequality gives the refined result.
		\begin{proof}[Proof of Proposition~\ref{prop.upper}]
			In order to make use of the one-arm probability in Proposition~\ref{prop.onearm_half_plus_delta}, we choose the correlation length
			\begin{align*}
				\ell_*(\p) := \lfloor c_*^{-1/2} \cdot (\p - \p_c)^{-1/2}\rfloor.
			\end{align*}
			Recall \eqref{eq.chain.2} and \eqref{eq.proba.identity}, then we obtain that 
			\begin{equation}\label{eq.upper_step1_delta}
				\begin{split}
					\ab(\p)^{\frac{1}{2}} &\leq \ab(\Lambda_{\ell_*,L}, \p)^{\frac{1}{2}} + \ell_* L^{-\frac{1}{2}} \\
					&= \ab^\delta(\Lambda_{\ell_*,L},\p; \p)^{\frac{1}{2}} +  \ell_* L^{-\frac{1}{2}} \\
					&\lesssim \ab(\Lambda_{\ell_*,L},\p; \p')^{\frac{1}{2}} + \int_{\p'}^{\p}  \bb^\delta(\Lambda_{\ell_*,L},\p; s)^{\frac{1}{2}} \, \d s +  \ell_* L^{-\frac{1}{2}}.
				\end{split}
			\end{equation}
			We then insert the bound of $\ab^\delta$ and $\bb^\delta$ in Lemma~\ref{lem.aabb_bound_new}, and the rest of calculation is similar as that of Proposition~\ref{prop.upper_log}: in Step~1, there is no longer the logarithm correction in $\bb^\delta$, so we get the correct order of perturbation that
			\begin{align}\label{eq.upper_step_3b_delta}
				\int_{\p'}^{\p}  \bb^\delta(\Lambda_{\ell_*,L},\p; s)^{\frac{1}{2}} \, \d s  \lesssim  \ell^{-3}_*.
			\end{align}
			In Step~2, the quantity $\ab^\delta$ is exactly same as $\ab$ there, which yields that 
			\begin{align}\label{eq.upper_step_3a_delta}
				\ab(\Lambda_{\ell_*,L},\p; \p')  \lesssim \ell_*^{-2} (\log{\ell_*}) \exp\Ll(- C^{-1} \tilde{C}^{\frac{1}{2}} \log \ell_*\Rr) \lesssim o(\ell_*^{-6}),
			\end{align}
			by taking  ${\p' = \p_c - \tilde{C} \ell_*^{-2}(\log \ell_*)^2}$ with a very large $ \tilde{C} $. Finally, combining \eqref{eq.upper_step1_delta}, \eqref{eq.upper_step_3b_delta}, \eqref{eq.upper_step_3a_delta}, and taking $L \to \infty$, we conclude the desired result in \eqref{eq.upper}.
		\end{proof}

		\subsection{One-arm probability in inhomogeneous percolation}\label{subsec.inhomo}
		This subsection is devoted to the proof of Proposition~\ref{prop.onearm_half_plus_delta}. The main argument is still the differential inequality, and the one-arm event will be analyzed in various contexts. 
		
		For every $m,n \in \mathbb{Z} \cup \{-\infty, + \infty\}$ satisfying $m < n$, we denote by $S_{m,n}$ the slab and $L_j$ the hyperplane
		\begin{align}\label{eq.def.Smn_Lj}
			S_{m,n} := ([m,n] \cap \ZZ ) \times \mathbb{Z}^{d-1} , \qquad L_j := \{j\} \times \mathbb{Z}^{d-1}.
		\end{align}
		We define $A_p(m,n;k,j)$ as the expected number of pioneers in $L_j$ issued from $k e_1$ within $S_{m,n}$, i.e. for every $m \leq k,j \leq n$
		\begin{align}\label{eq.def.Amnkj}
			A_p(m,n;k,j) := \sum_{u \in L_j} \tau_p^{S_{m,n}}(ke_1, u) = \sum_{u \in L_j} \P_{p}\Ll[ k e_1 \xleftrightarrow{S_{m,n}} u \Rr].
		\end{align}
		This quantity extends  $\psi_p, \varphi_p$ defined in \eqref{eq.varphi}, since the spatial information is richer.

		Let us state at first several basic properties about the quantity $A_p$.
		\begin{prop}\label{prop.A.basic}
			\begin{enumerate}
				\item The following identities always hold
				\begin{align}
					A_p(m,n;k,j) &= A_p(m,n;j,k), \label{eq.A.identity1}\\
					A_p(0,n;k,j) &= A_p(0,n;n-k,n-j). \label{eq.A.identity2}
				\end{align}
				\item $A_p$ satisfies the convolution-type inequality
				\begin{align}\label{eq.A.convolution}
					A_p(m,n;k,j) \leq \sum_{i=k \vee j}^n A_p(m,i;k,i) A_p(m,i;j,i).
				\end{align}
				\item $A_p$ satisfies the differential inequality
				\begin{align}\label{eq.A.differential}
					\frac{\d }{\d p} A_p(m,n;k,j) \lesssim p \sum_{i=m}^n A_p(m,n;k,i) A_p(m,n;i,j).
				\end{align}
				\item In high-dimensional percolation, there exists a finite positive constant $C(d)$ such that the following estimate holds for all $p \in (0,\p_c]$
				\begin{align}\label{eq.A.basis}
					A_p(-\infty,n;k,n) \leq C\exp\Ll(-C^{-1} (n-k) \sqrt{\p_c - p}\Rr).
				\end{align}
			\end{enumerate}
		\end{prop}
		\begin{proof}
			In (1), the identity \eqref{eq.A.identity2} is a consequence of the symmetry of slab. For \eqref{eq.A.identity1}, we have 
			\begin{align*}
				A_p(m,n;k,j) &= \sum_{u \in L_j} \tau_p^{S_{m,n}}(ke_1, u) \\
				&= \sum_{u \in L_j} \tau_p^{S_{m,n}}(u, ke_1) \\
				&= \sum_{u' \in L_k} \tau_p^{S_{m,n}}(j e_1,  u').
			\end{align*}
			The second line makes use of the reversibility of the open path. Then the third line is just a change of variable and relies on the stationarity of the percolation.
			
			Concerning (2), for every $u \in L_j$, we consider the highest hyperplane that the open path $\{k e_1 \xleftrightarrow{S_{m,n}} u\}$ can reach. It can be $\{L_i\}_{ (k \vee j) \leq i \leq n}$ and the open path can be decomposed by its first hit of $L_i$ at $v$ 
			\begin{align*}
				\{k e_1 \xleftrightarrow{S_{m,n}} u\} \subset \bigcup_{i = k \vee j }^n \{k e_1 \xleftrightarrow{S_{m,i}} v\} \circ \{v \xleftrightarrow{S_{m,i}} u\}.
			\end{align*} 
			Here we can restrict the slab from $S_{m,n}$ to $S_{m,i}$ because $L_i$ is the highest hyperplane. Applying the BK inequality to the decomposition above, we obtain the convolution-type inequality
			\begin{align*}
				A_p(m,n;k,j) &\leq \sum_{i = k \vee j }^n \sum_{u \in L_j} \sum_{v  \in L_i}\P_p\Ll[ \{k e_1 \xleftrightarrow{S_{m,i}} v\} \circ \{v \xleftrightarrow{S_{m,i}} u\}\Rr]\\
				&\leq \sum_{i = k \vee j }^n \Ll( \sum_{v  \in L_i}\tau_p^{S_{m,i}}(ke_1, v)\Rr) \Ll(\sum_{u \in L_j} \tau_p^{S_{m,i}}(v, u)\Rr)\\
				&\leq \sum_{i=k \vee j}^n A_p(m,i;k,i) A_p(m,i;j,i).
			\end{align*}
			
			Concerning (3), we apply Russo's formula to the event $\{k e_1 \xleftrightarrow{S_{m,n}} u\}$
			\begin{align*}
				\frac{\d}{\d p} \tau_p^{S_{m,n}}(ke_1, u) = \sum_{e \in \Ed(S_{m,n})} \P_p\Ll[\{k e_1 \xleftrightarrow{S_{m,n}} u\} \cap \{e \text{ is pivotal}\}  \Rr].
			\end{align*}
			An edge $e = \{e_\circ, e_\bullet\}$ is pivotal for $\{k e_1 \xleftrightarrow{S_{m,n}} u\}$ implies that, the open path can be decomposed into three disjoint segments 
			\begin{align*}
				\{k e_1 \xleftrightarrow{S_{m,n}} e_\circ\} \circ \{e \text{ is open }\}  \circ \{e_\bullet \xleftrightarrow{S_{m,n}} u\}.
			\end{align*}
			Therefore, we use the BK inequality to obtain the differential inequality
			\begin{align*}
				\frac{\d }{\d p} A_p(m,n;k,j) &\leq \sum_{u \in L_j} \sum_{i=m}^n \sum_{e_\circ \in L_i, e_\bullet \sim e_\circ}\P_p\Ll[ \{k e_1 \xleftrightarrow{S_{m,n}} e_\circ\} \circ \{e \text{ is open }\}  \circ \{e_\bullet \xleftrightarrow{S_{m,n}} u\}\Rr]\\
				&\leq p \sum_{i=m}^n \Ll(\sum_{e_\circ \in L_i, e_\bullet \sim e_\circ} \tau_p^{S_{m,n}}(ke_1, e_\circ)\Rr) \Ll(\sum_{u \in L_j}  \tau_p^{S_{m,n}}(e_\bullet, u)\Rr)\\
				&\lesssim p \sum_{i=m}^n A_p(m,n;k,i) A_p(m,n;i,j).
			\end{align*}

			(4) is a consequence of \eqref{eq.expansion_12}
			\begin{align*}
				A_p(-\infty,n;k,n) \leq \psi_p(\H_{n-k}) \leq C\exp\Ll(-C^{-1} (n-k) \sqrt{\p_c - p}\Rr).
			\end{align*}
		\end{proof}
		
		We aim to deduce more general estimates starting from Proposition~\ref{prop.A.basic}. The first observation is an upper bound for general $A_p(m,n;k,j)$ in subcritical regime. 
		\begin{lemma}\label{lem.A.upperbound}
			In high-dimensional percolation, there exists a finite positive constant $C(d)$ such that the following estimate holds for all $p \in (0,\p_c]$
			\begin{multline}\label{eq.A.upperbound}
				A_p(m,n;k,j) \leq C\exp\Ll(-C^{-1} \vert j - k \vert \sqrt{\p_c - p}\Rr) \\
				\cdot \min\Ll\{(n+1-j), (n+1-k), (k+1-m), (j+1-m), (\p_c - p)^{-1/2}\Rr\}.
			\end{multline}
		\end{lemma}
		\begin{proof}
			Without loss of generality, we assume that $m \leq k \leq j \leq n$. The upper bound comes from \eqref{eq.A.convolution}. It transforms the general $A_p(m,n;k,j)$ to the convolution of diagonal version, and the latter is known from \eqref{eq.A.basis}
			\begin{align*}
				A_p(m,n;k,j) &\leq \sum_{i=j}^n A_p(m,i;k,i) A_p(m,i;j,i)\\
				&\lesssim \sum_{i=j}^n  \exp\Ll(-C^{-1} (i-k) \sqrt{\p_c - p}\Rr) \exp\Ll(-C^{-1} (i-j) \sqrt{\p_c - p}\Rr)\\
				&\lesssim \exp\Ll(-C^{-1} (j-k) \sqrt{\p_c - p}\Rr)  \sum_{i=j}^n  \exp\Ll(-C^{-1} (i-j) \sqrt{\p_c - p}\Rr) \\
				&\lesssim \exp\Ll(-C^{-1} (j-k) \sqrt{\p_c - p}\Rr) \cdot \min\{(n+1-j), (\p_c - p)^{-1/2} \}.
			\end{align*}
			In the third line, we use the fact $i \geq j$, so the first exponential factor is bounded by $\exp\Ll(-C^{-1} (j-k) \sqrt{\p_c - p}\Rr)$. We then make a sum of the second exponential factor. It can be bounded either by its integral, or trivially by the number of the terms $(n+1-j)$.
			
			Notice that $(n+1-j)$ is the distance between $L_k$ or $L_j$ to one boundary. Using the symmetry in \eqref{eq.A.identity2}, we can choose a minimum among all the possible distances. This concludes \eqref{eq.A.upperbound}.
		\end{proof}
		
		We then aim to extend the result to supercritical regime.
		\begin{lemma}\label{lem.A.upperbound.super}
			In high-dimensional percolation, there exists a finite positive constant $c_{**}(d)$, such that for all $p \leq \p_c + c_{**} n^{-2}$, we have
			\begin{align}\label{eq.A.upperbound.super}
				A_p(0,n;k,j) \lesssim \min\{(n+1-j), (n+1-k), j+1, k+1\}.
			\end{align} 
		\end{lemma}
		\begin{proof}
			The estimate \eqref{eq.A.upperbound.super} is valid for $p = \p_c$. In order to study the value $p > \p_c$, we make use of the differential inequality and maximal function. An auxiliary diagonal maximal function $M_p (n)$ is defined as 
			\begin{align}\label{eq.M.def}
				M_p(n)  := \max_{0 \leq k \leq j \leq n} A_p(0,j;k,j).
			\end{align}
			The maximal function also gives a bound for the general $A_p(0,n;k,j)$ via the convolution-type inequality \eqref{eq.A.convolution}
			\begin{align}\label{eq.A.maximal_1}
				A_p(0,n;k,j) \leq \sum_{i=k \vee j}^n A_p(0,i;k,i) A_p(0,i;j,i) \leq \Ll(n+1- (k \vee j)\Rr) M_p(n)^2.
			\end{align}
			In the beginning, we only have the information of $M_p(n)$ for $p = \p_c$ thanks to \eqref{eq.A.basis} 
			\begin{align}\label{eq.M.pc}
				M_{\p_c}(n) \leq C_{\text{\ref{eq.M.pc}}}.
			\end{align}
			Therefore, we study its value at $p$ via a differential inequality. Viewing the definition \eqref{eq.M.def}, let us first apply \eqref{eq.A.differential} to the derivative to $A_p(0,j;k,j)$, $0 \leq j \leq n$
			\begin{align*}
				\frac{\d}{\d p} A_p(0,j;k,j)  &\lesssim p \sum_{i=m}^n A_p(0,j;k,i) A_p(0,j;i,j) \\
				&\lesssim \sum_{i=m}^n \Ll(n+1- (k \vee i)\Rr) M_p(n)^2 \cdot M_{p}(n)\\
				&\lesssim n^2 M_p(n)^3.
			\end{align*}
			In the second line, \eqref{eq.A.maximal_1} is applied to $A_p(0,j;k,i)$ and \eqref{eq.M.def} is applied to $A_p(0,j;i,j)$. The uniform estimate above then implies a differential inequality of $M_p(n)$ that
			\begin{align}\label{eq.M.diff}
				\frac{\d}{\d p}  M_p(n) \leq c_{\text{\ref{eq.M.diff}}} n^2 M_p(n)^3,
			\end{align}
			which is equivalent to 
			\begin{align*}
				\frac{\d}{\d p}  M_p(n)^{-2} \geq -2 c_{\text{\ref{eq.M.diff}}} n^2.
			\end{align*}
			An integral from $\p_c$ to $p$ yields that 
			\begin{align*}
				M_p(n)^{-2} \geq M_{\p_c}(n)^{-2}  -2 c_{\text{\ref{eq.M.diff}}} n^2 (p - \p_c) \geq C^{-2}_{\text{\ref{eq.M.pc}}} -2 c_{\text{\ref{eq.M.diff}}} n^2 (p - \p_c).
			\end{align*} 
			Taking $c_{**} = \frac{1}{4} c^{-1}_{\text{\ref{eq.M.diff}}} C^{-2}_{\text{\ref{eq.M.pc}}}$, then we have the estimate 
			\begin{align}\label{eq.M.upper}
				\forall p \leq \p_c + c_{**} n^{-2}, \qquad M_p(n) \leq 2 C_{\text{\ref{eq.M.pc}}}.
			\end{align}
			Finally, we put the bound back to \eqref{eq.A.maximal_1} and use the symmetry in \eqref{eq.A.identity2}, then \eqref{eq.A.upperbound.super} is proved.
		\end{proof}
		
		The estimate \eqref{eq.A.upperbound.super} is slightly too rough as it fails to capture the effect of the half-plane. The following bound, sharing the same philosophy as Proposition~\ref{prop.onearm_half_plus}, is more precise when $k$ is very close $0$.
		Moreover, $A_p (0,n;k,j)$ is finite implies that $p$ is a subcritical parameter in $S_{0,n}$, which gives the lower bound of the slab threshold (Theorem~\ref{thm.slab}), see Section~\ref{sec.chem.slab}.
		\begin{lemma}\label{lem.A.upperbound.half}
			In high-dimensional percolation, there exists a finite positive constant $C(d)$ such that the following estimate holds
			\begin{equation}\label{eq.A.upperbound.half}
				A_p(0,n;k,n) \leq \Ll\{\begin{array}{ll}
					C \Ll(\frac{k+1}{n+1}\Rr) , & p \in (\p_c, \p_c + c_{**} h^{-2}];\\
					C \Ll(\frac{k+1}{n+1}\Rr)  \exp\Ll(- C^{-1} (n-k) \sqrt{\p_c - p}\Rr)  , & p \in (0, \p_c].
				\end{array}\Rr.
			\end{equation}
			Here $c_{**}$ is the same constant in Lemma~\ref{lem.A.upperbound.super}.
		\end{lemma}
		\begin{proof}
			The bound is only interesting when $0 \leq k \leq n$, otherwise we can just cite \eqref{eq.A.upperbound.super} and \eqref{eq.A.upperbound}. 
			
			We denote by $A'_p(0,n;k,n) = \frac{\d}{\d p} A_p(0,n;k,n)$ the derivative. Under the condition ${0 \leq k \leq n}$, we calculate $A_p(0,n;k,n)$ as an integral
			\begin{align}\label{eq.A.integral}
				A_p(0,n;k,n) = \int_0^p A'_s(0,n;k,n) \, \d s \lesssim  \int_0^p \sum_{i=0}^n A_s(0,n;k,i) A_s(0,n;i,n) \, \d s.
			\end{align}
			The inequality comes from \eqref{eq.A.differential}. We then divide it into two cases.
			
			\emph{Case~1: $p \in (0, \p_c]$.} We insert the bound \eqref{eq.A.upperbound}
			\begin{align*}
				&A_p(0,n;k,n) \\
				&\lesssim  \int_0^p \sum_{i=0}^n  \exp\Ll(-C^{-1} \Ll(\vert i - k \vert + \vert i - n \vert\Rr) \sqrt{\p_c - s}\Rr)  \min \{(k+1), (\p_c -s)^{-1/2}\}\, \d s\\
				&\lesssim (k+1)(n+1)  \int_0^p   \exp\Ll(-C^{-1} (n-k) \sqrt{\p_c - s}\Rr)  \, \d s \\
				&\lesssim (k+1)(n+1)  \int_{\sqrt{\p_c -p}}^{\sqrt{\p_c}}   t \exp\Ll(-C^{-1} (n-k) t\Rr)  \, \d t \\
				&\lesssim \frac{(n+1)(k+1)}{(n-k)^2}\exp\Ll(-C^{-1} (n-k) \sqrt{\p_c -p}\Rr) .
			\end{align*}
			The passage to the third line relies on the fact that ${\vert i - k \vert + \vert i - n \vert \geq (n-k)}$. Then we make a change of variable $s = \p_c - t^2$. Finally, recall the condition ${0 \leq k \leq n}$, then we have ${\frac{(n+1)(k+1)}{(n-k)^2} \leq \frac{(k+1)}{(n+1)}}$, which yields the desired bound as \eqref{eq.A.upperbound.half} for $p \in (0, \p_c]$.
			
			\emph{Case~2: $p \in (\p_c, \p_c + c_{**} n^{-2}]$.} For such regime, \eqref{eq.A.upperbound.super} can be applied to \eqref{eq.A.integral}
			\begin{align*}
				A_s(0,n;k,i) &\lesssim \min\{(n+1-j), (n+1-k), i+1, k+1\} \lesssim (k+1), \\
				A_s(0,n;i,n) &\lesssim 1.
			\end{align*}
			Then an integral along $(\p_c, p]$ yields that 
			\begin{align*}
				A_p(0,n;k,n) &\lesssim A_{\p_c}(0,n;k,n) + \int_{\p_c}^p \sum_{i=0}^n  A_s(0,n;k,i) A_s(0,n;i,n) \, \d s \\
				&\lesssim A_{\p_c}(0,n;k,n) +  \int_{\p_c}^p   n(k+1) \, \d s  \\
				&\lesssim \frac{(k+1)}{(n+1)} +  c_{**} n^{-2} \cdot n(k+1) \\
				&\lesssim \frac{(k+1)}{(n+1)}.
			\end{align*}
			In the third line, we insert the condition $p \in (\p_c, \p_c + c_{**} n^{-2}]$ and the estimate at $\p_c$ from Case~1. This finishes the proof.
		\end{proof}
		
		With the help of the extensions above, we turn to the inhomogeneous slab. We denote by $\P_{p,q}$ the probability space that samples Bernoulli percolation of parameter $p$ in $S_{0,m}$ and of parameter $q$ in $S_{m,n}$. Then, for $0 \leq k,j \leq n$, we define $A_{p,q}(0,m,n;k,j)$ as a version of $A_{p}(m,n;k,j)$ under $\P_{p,q}$
		\begin{align}\label{eq.def.A.inhomo}
			A_{p,q}(0,m,n;k,j) := \sum_{u \in L_j} \P_{p,q}\Ll[ k e_1 \xleftrightarrow{S_{0,n}} u \Rr].
		\end{align} 
		If $n \leq m$, then we just treat it as no $q$-percolation sampled, thus $A_{p,q}$ is reduced to homogeneous percolation and we have $A_{p,q}(0,m,n;k,j) = A_{p}(0,n;k,j)$.

		We have an estimate for the quantity $A_{p,q}(0,m,n;k,j)$ defined in \eqref{eq.def.A.inhomo}. In the formula below, $d_m$ is a distance truncated above $L_m$  
		\begin{align}\label{eq.def.dm}
			d_m(k,j) := \vert (j-m)_+ - (k-m)_+ \vert.
		\end{align}
		It captures the minimum distance that an open path from $k e_1$ to $L_j$ has to cross in the subcritical percolation $S_{m,n}$.

		\begin{lemma}\label{lem.A.upperbound.inhomo}
			In high-dimensional percolation, there exist finite positive constants $C(d), \tilde{c}_{**}(d)$, such that for every $q \in (0,\p_c]$ and $p \leq \p_c + \tilde{c}_{**} n^{-1} m^{-1}$, we have 
			\begin{align}\label{eq.A.upperbound.inhomo}
				A_{p,q}(0,m,n;k,j) \leq C (n+1-j) \exp\Ll(-C^{-1} d_m(k,j) \sqrt{\p_c - q}\Rr).
			\end{align}
		\end{lemma}
		\begin{proof}
			The quantity $A_{p,q}(0,m,n;k,j)$ satisfies the similar properties as \eqref{eq.A.identity1}
			\begin{align}\label{eq.A.identity1_new}
				A_{p,q}(0,m,n;k,j) = A_{p,q}(0,m,n;j,k),
			\end{align}
			and as \eqref{eq.A.convolution}
			\begin{align}\label{eq.A.convolution_new}
				A_{p,q}(0,m,n;k,j) \leq \sum_{i=k \vee j}^n A_{p,q}(0,m,i;k,i)  A_{p,q}(0,m,i;j,i).
			\end{align}
			Therefore, it suffices to study the diagonal terms like $A_{p,q}(0,m,i;k,i)$. We divide the proof into two steps.
			
			\emph{Step~1: differential inequality.}	The proof follows that of Lemma~\ref{lem.A.upperbound.super}. We also define a maximal function 
			\begin{align}\label{eq.M.def.mn}
				M_{p,q}(m,n)  := \max_{0 \leq k , j \leq n} A_{p,q}(0,m,j;k,j) \exp\Ll(c_{\text{\ref{eq.M.def.mn}}} d_m(k,j) \sqrt{\p_c - q}\Rr).
			\end{align}
			Here the constant $c_{\text{\ref{eq.M.def.mn}}}$ is small and will be determined later.

			Using \eqref{eq.A.convolution_new}, the counterpart of \eqref{eq.A.maximal_1} is 		
			\begin{equation}\label{eq.A.maximal_2}
				\begin{split}
					&A_{p,q}(0,m,n;k,j) \\
					&\leq \sum_{i=k \vee j}^n A_{p,q}(0,m,i;k,i)  A_{p,q}(0,m,i;j,i) \\
					&\leq \sum_{i=k \vee j}^n M_{p,q}(m,n)^2 \exp\Ll(-c_{\text{\ref{eq.M.def.mn}}} \Ll( d_m(k,i) + d_m(i,j)\Rr)\sqrt{\p_c - q}\Rr)\\
					&\leq \Ll(n+1- (k \vee j)\Rr) M_{p,q}(m,n)^2 \exp\Ll(-c_{\text{\ref{eq.M.def.mn}}} d_m(k,j) \sqrt{\p_c - q}\Rr).
				\end{split}
			\end{equation}
			In the third line, the triangle inequality $d_m(k,j) \leq d_m(k,i) + d_m(i,j)$ is applied. The differential inequality \eqref{eq.A.differential} also yields that 
			\begin{align*}
				&\frac{\d}{\d p} A_{p,q}(0,m,j;k,j)  \\
				&\lesssim p \sum_{i=0}^m  A_{p,q}(0,m,j;k,i) A_{p,q}(0,m,j;i,j) \\
				&\lesssim \sum_{i=0}^m \Ll(n+1- (k \vee i)\Rr) M_{p,q}(m,n)^3 \exp\Ll(-c_{\text{\ref{eq.M.def.mn}}} \Ll(d_m(k,i) + d_m(i,j)\Rr) \sqrt{\p_c - q}\Rr) \\
				&\lesssim m n M_{p,q}(m,n)^3 \exp\Ll(-c_{\text{\ref{eq.M.def.mn}}} d_m(k,j) \sqrt{\p_c - q}\Rr).
			\end{align*}
			We insert the bounds \eqref{eq.A.maximal_2} and \eqref{eq.M.def.mn} in the third line above. Here the major difference is the percolation of parameter $p$ is only sampled in $S_{0,m}$. The estimate above implies that 
			\begin{align}\label{eq.M.diff_2}
				\frac{\d}{\d p} M_{p,q}(m,n) \leq c_{\text{\ref{eq.M.diff_2}}} m n M_{p,q}(m,n)^3.
			\end{align}
			
			\emph{Step~2: starting point $p=\p_c$.} It is natural to follow the steps after \eqref{eq.M.diff}, then the choice of starting point is necessary. A naive choice is $p = q$, then we have $A_{q,q}(0,m,j;k,j) = A_{q}(0,j;k,j)$ and \eqref{eq.A.basis} suggests that $M_{q,q}(m,n) \leq C$. However, following the steps \eqref{eq.M.diff} and \eqref{eq.M.upper}, this only allows to obtain a result for $p \leq q + c  n^{-1} m^{-1}$, which is weaker than the statement in Lemma~\ref{lem.A.upperbound.inhomo}. We thus aim to prove the basis for $p = \p_c$
			\begin{align}\label{eq.M.mnbound}
				M_{\p_c,q}(m,n) \leq C_{\text{\ref{eq.M.mnbound}}}.
			\end{align}
			We turn to analyze the upper bound of $A_{\p_c,q}(0,m,j;k,j)$ and it can be divided into $3$ cases concerning the position of $k,j$. Without loss of generality, we just assume $0 \leq k \leq j \leq n$.
			
			\emph{Case~a): $0 \leq k,j \leq m$.} For this case, we use the convention around \eqref{eq.def.A.inhomo} 
			\begin{align*}
				A_{\p_c,q}(0,m,j;k,j) = A_{\p_c}(0,j;k,j) \leq C.
			\end{align*}
			The last bound comes from \eqref{eq.A.basis}. Since $d_m(k,j) = 0$ for such case, we then obtain that
			\begin{align}\label{eq.M.mnbound.a}
				\max_{0 \leq k \leq j \leq n} A_{\p_c,q}(0,m,j;k,j) \exp\Ll(c_{\text{\ref{eq.M.def.mn}}} d_m(k,j) \sqrt{\p_c - q}\Rr) \leq C.
			\end{align}
			See Figure~\ref{fig.inhomo_a} for an illustration.
			\begin{figure}[h!]
				\centering
				\includegraphics[width=0.6\textwidth]{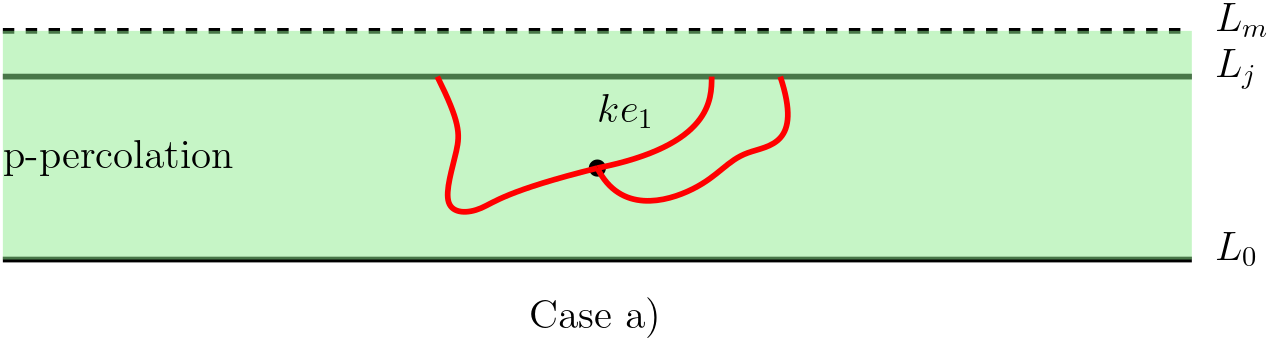}
				\caption{An illustration of $A_{p,q}(0,m,j;k,j)$, $0 \leq k,j \leq m$.}\label{fig.inhomo_a}
			\end{figure}
			
			\emph{Case~b): $0 \leq k \leq m \leq j$.} For every open path $\{k e_1 \xleftrightarrow{S_{0,j}} u\}, u \in L_j$, let $v$ be its last hit at $L_m$. We decompose the path into disjoint realization
			\begin{align*}
				\{k e_1 \xleftrightarrow{S_{0,j}} u\} \subset \bigcup_{v \in L_m} \{k e_1 \xleftrightarrow{S_{0,j}} v\} \circ \{v \xleftrightarrow{S_{m,j}} u\}.
			\end{align*} 
			Then we apply the BK inequality and obtain that 
			\begin{equation}\label{eq.M.mnbound.b.step}
				\begin{split}
					A_{\p_c,q}(0,m,j;k,j) &\leq \sum_{v \in L_m} \P_{\p_c,q}\Ll[ k e_1 \xleftrightarrow{S_{0,j}} v \Rr] \sum_{u \in L_j} \P_{\p_c,q}\Ll[ v \xleftrightarrow{S_{m,j}} u \Rr]\\
					&= A_{\p_c,q}(0,m,j;k,m) A_{q}(m,j;m,j)\\
					&\leq A_{\p_c}(0,j;k,m) A_{q}(m,j;m,j). 
				\end{split}
			\end{equation}
			Here we also rely on the fact that $S_{m,j}$ samples a percolation of parameter $q$, and relax $\P_{\p_c, q}$ to $\P_{\p_c}$ for the first term. Using the results in homogeneous slab, we obtain that 
			\begin{align*}
				A_{\p_c}(0,j;k,m) &\stackrel{\text{\eqref{eq.A.upperbound}}}{\lesssim} (j-m+1), \\
				A_{q}(m,j;m,j) &\stackrel{\text{\eqref{eq.A.upperbound.half}}}{\lesssim} \frac{1}{(j-m+1)}\exp\Ll(-C^{-1} (j-m) \sqrt{\p_c - q}\Rr).
			\end{align*}
			Notice that $(j-m) = d_m(j,k)$. Putting these bounds back to \eqref{eq.M.mnbound.b.step}, we conclude that 
			\begin{multline}\label{eq.M.mnbound.b}
				\max_{0 \leq k \leq m \leq j} A_{\p_c,q}(0,m,j;k,j) \exp\Ll(c_{\text{\ref{eq.M.def.mn}}} d_m(k,j) \sqrt{\p_c - q}\Rr)\\
				\leq \exp\Ll(\Ll(c_{\text{\ref{eq.M.def.mn}}} -C^{-1}\Rr) d_m(k,j) \sqrt{\p_c - q}\Rr) \leq C.
			\end{multline}
			Here we need to take $c_{\text{\ref{eq.M.def.mn}}} \leq C^{-1}$. See Figure~\ref{fig.inhomo_a} for an illustration.
			\begin{figure}[h!]
				\centering
				\includegraphics[width=0.6\textwidth]{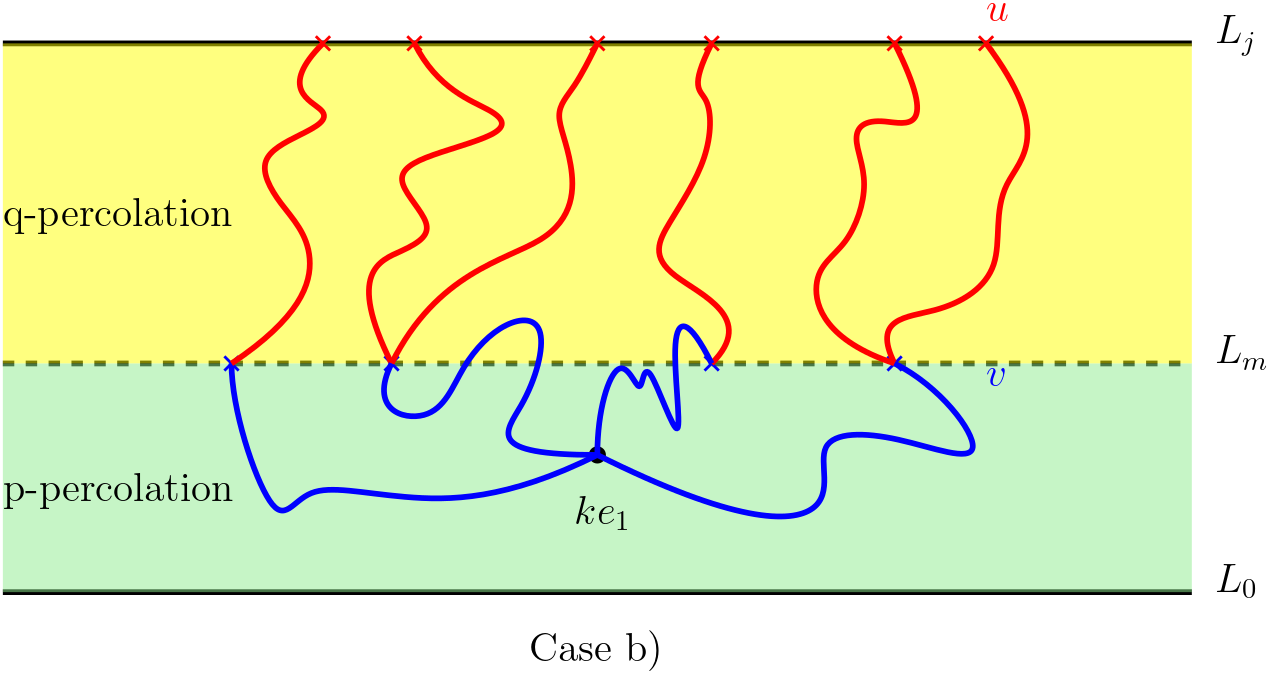}
				\caption{An illustration of $A_{p,q}(0,m,j;k,j)$, $0 \leq k \leq m \leq j$.}\label{fig.inhomo_b}
			\end{figure}
			
			\emph{Case~c): $0 \leq m \leq k \leq j$.} Recall the symmetry \eqref{eq.A.identity1_new}, then we have $A_{\p_c,q}(0,m,j;k,j) = A_{\p_c,q}(0,m,j;j,k)$, so we analyze the latter. Similar as the proof of Case~b), for every open path $\{j e_1 \xleftrightarrow{S_{0,j}} u\}, u \in L_k$, let $v$ be its first hit at $L_{\lfloor \frac{k+j}{2} \rfloor}$. Then we decompose the path into disjoint realization
			\begin{align*}
				\{j e_1 \xleftrightarrow{S_{0,n}} u\} \subset \bigcup_{v \in L_{\lfloor \frac{k+j}{2} \rfloor}} \{j e_1 \xleftrightarrow{S_{\lfloor \frac{k+j}{2} \rfloor,j}} v\} \circ \{v \xleftrightarrow{S_{0,j}} u\}.
			\end{align*} 
			The BK inequality then applies and yields that 
			\begin{equation}\label{eq.M.mnbound.c.step}
				\begin{split}
					A_{\p_c,q}(0,m,j;j,k) &\leq \sum_{v \in L_{\lfloor \frac{k+j}{2} \rfloor}} \P_{\p_c,q}\Ll[ j e_1 \xleftrightarrow{S_{\lfloor \frac{k+j}{2} \rfloor,j}} v \Rr] \sum_{u \in L_k} \P_{\p_c,q}\Ll[ v \xleftrightarrow{S_{0,j}} u \Rr]\\
					&=A_{q}\Ll(\lfloor \frac{k+j}{2} \rfloor,j;j, \lfloor \frac{k+j}{2} \rfloor\Rr) A_{\p_c,q}\Ll(0,m,j;\lfloor \frac{k+j}{2} \rfloor,k\Rr) \\
					&=A_{q}\Ll(\lfloor \frac{k+j}{2} \rfloor,j;\lfloor \frac{k+j}{2} \rfloor,j \Rr) A_{\p_c,q}\Ll(0,m,j;k, \lfloor \frac{k+j}{2} \rfloor\Rr) \\
					&\leq A_{q}\Ll(\lfloor \frac{k+j}{2} \rfloor,j;\lfloor \frac{k+j}{2} \rfloor,j \Rr) A_{\p_c}\Ll(0,j;k, \lfloor \frac{k+j}{2} \rfloor\Rr). 
				\end{split}
			\end{equation}
			In the second line, we use the fact that $S_{\lfloor \frac{k+j}{2} \rfloor,j}$ samples a percolation of parameter $q$. In the third line, \eqref{eq.A.identity1_new} and \eqref{eq.A.identity1} are applied. We then relax $\P_{\p_c, q}$ to $\P_{\p_c}$ for the second term in the forth line. Afterwards, we apply the results in  homogeneous slab once again
			\begin{align*}
				A_{q}\Ll(\lfloor \frac{k+j}{2} \rfloor,j;\lfloor \frac{k+j}{2} \rfloor,j \Rr) &\stackrel{\text{\eqref{eq.A.upperbound.half}}}{\lesssim} \frac{2}{(j-k+1)}\exp\Ll(-(2C)^{-1} (j-k) \sqrt{\p_c - q}\Rr)\\
				A_{\p_c}\Ll(0,j;k, \lfloor \frac{k+j}{2} \rfloor\Rr) &\stackrel{\text{\eqref{eq.A.upperbound}}}{\lesssim} \frac{j-k+1}{2}. \\
			\end{align*}
			These estimates, together with \eqref{eq.M.mnbound.c.step}, conclude that 
			\begin{multline}\label{eq.M.mnbound.c}
				\max_{0 \leq m \leq k \leq j} A_{\p_c,q}(0,m,j;k,j) \exp\Ll(c_{\text{\ref{eq.M.def.mn}}} d_m(k,j) \sqrt{\p_c - q}\Rr)\\
				\leq \exp\Ll(\Ll(c_{\text{\ref{eq.M.def.mn}}} -(2C)^{-1}\Rr) d_m(k,j) \sqrt{\p_c - q}\Rr) \leq C.
			\end{multline}
			The choice of parameter is  $c_{\text{\ref{eq.M.def.mn}}} \leq (2C)^{-1}$. See Figure~\ref{fig.inhomo_c} for an illustration.
			\begin{figure}[h!]
				\centering
				\includegraphics[width=0.6\textwidth]{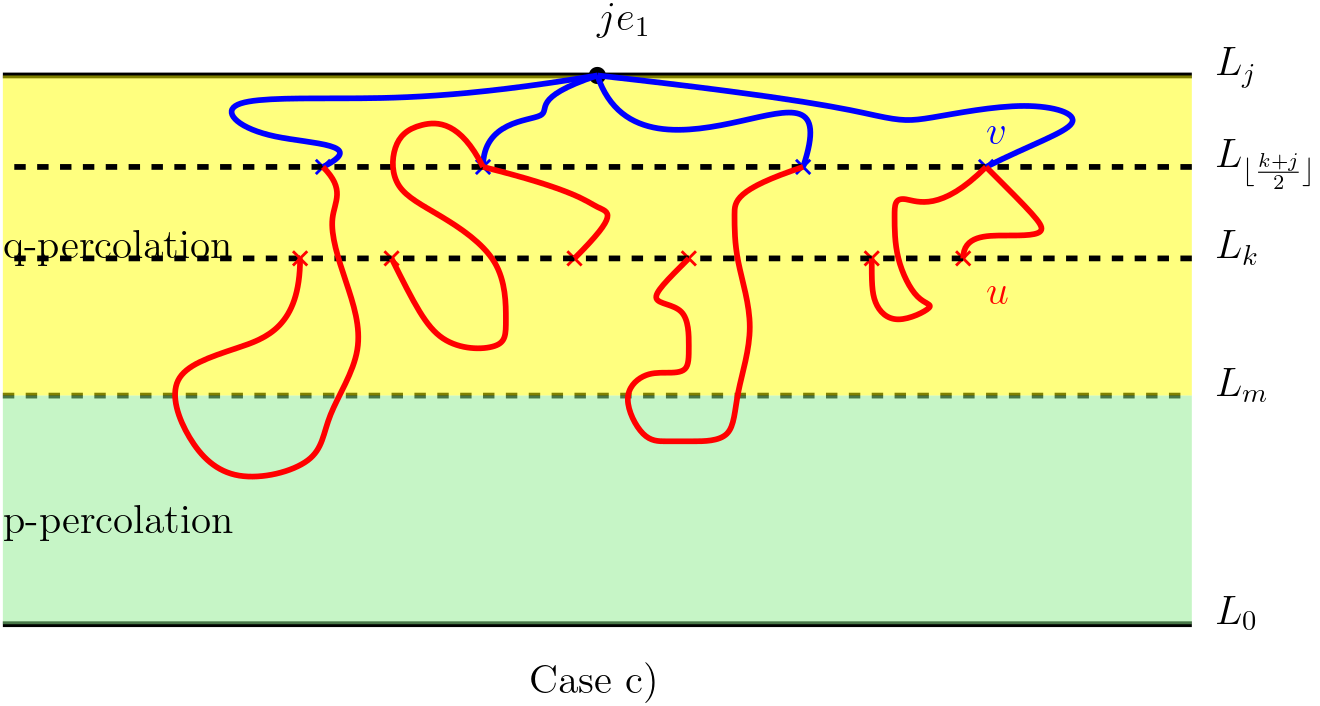}
				\caption{An illustration of $A_{p,q}(0,m,j;k,j)$, $0 \leq m \leq k \leq j$.}\label{fig.inhomo_c}
			\end{figure}
			
			In conclusion, \eqref{eq.M.mnbound.a}, \eqref{eq.M.mnbound.b}, \eqref{eq.M.mnbound.c} together justify \eqref{eq.M.mnbound}. Therefore, \eqref{eq.M.mnbound} and \eqref{eq.M.diff_2} suggest the existence of $\tilde{c}_{**}(d)$, such that for every $q \in (0,\p_c]$ and $p \leq \p_c + \tilde{c}_{**} n^{-1} m^{-1}$.
			\begin{align*}
				M_{\p_c,q}(m,n) \leq C_{\text{\ref{eq.M.mnbound}}}/2.
			\end{align*}
			Putting this bound back to \eqref{eq.A.maximal_2}, then \eqref{eq.A.upperbound.inhomo} is proved.
		\end{proof}

		We are now ready to prove Proposition~\ref{prop.onearm_half_plus_delta}.
		\begin{proof}[Proof of Proposition~\ref{prop.onearm_half_plus_delta}]
			For convenience, we use the notations in \eqref{eq.def.Smn_Lj} and it suffices to focus on the situation $\tilde{\p} \in (\p_c + c_* h^{-2}]$. The proof can be divided into two cases.
			
			\emph{Case~1: $0 \leq k \leq \delta h /2$.} For every open path $\{k e_1 \longleftrightarrow  [ \partial \Lambda_{h,L}]^{+}\}$, we denote by $u_1$ its first hit at $L_{\lfloor \frac{\delta h}{2}\rfloor}$, and $u_2$ its first hit at $L_{\lfloor (1-\frac{\delta}{2})h \rfloor}$. These two vertices break $\{k e_1 \longleftrightarrow  [ \partial \Lambda_{h,L}]^{+}\}$ into disjoint realization 
			\begin{multline*}
				\{k e_1 \longleftrightarrow  [ \partial \Lambda_{h,L}]^{+}\} \\
				\subset \bigcup_{u_1 \in L_{\lfloor \frac{\delta h}{2}\rfloor}} \bigcup_{u_2 \in L_{\lfloor (1-\frac{\delta}{2})h \rfloor}} \{k e_1 \xleftrightarrow{S_{0, \lfloor \frac{\delta h}{2}\rfloor}}  u_1\} \circ \{u_1 \xleftrightarrow{S_{0, \lfloor (1-\frac{\delta}{2})h \rfloor}}  u_2\} \circ \{u_2 \xleftrightarrow{S_{0, h}}  [ \partial \Lambda_{h,L}]^{+}\}.
			\end{multline*}
			Then the BK inequality applies 
			\begin{align*}
				&\P_{\p, \tilde{\p}}^{\delta}\Ll[k e_1 \longleftrightarrow  [ \partial \Lambda_{h,L}]^{+}\Rr] \\
				&\leq \Ll(\sum_{u_1 \in L_{\lfloor \frac{\delta h}{2}\rfloor}} \P_{\p}\Ll[ k e_1 \xleftrightarrow{S_{0, \lfloor \frac{\delta h}{2}\rfloor}}  u_1\Rr]  \Rr) \Ll(\sum_{u_2 \in L_{\lfloor (1-\frac{\delta}{2})h \rfloor}}  \P_{\p, \tilde{\p}}^{\delta}\Ll[u_1 \xleftrightarrow{S_{0, \lfloor (1-\frac{\delta}{2})h \rfloor}}  u_2\Rr] \Rr)\\
				& \qquad \qquad \cdot  \P_{\p}\Ll[u_2 \xleftrightarrow{S_{0, h}}  [ \partial \Lambda_{h,L}]^{+}\Rr]\\
				&\leq A_{\p}\Ll(0, \lfloor \delta h / 2 \rfloor; k, \lfloor \delta h / 2 \rfloor\Rr) A_{\p, \tilde{\p}}\Ll(0, \lfloor \delta h / 2 \rfloor, \lfloor (1-\delta/2)h \rfloor; \lfloor \delta h / 2 \rfloor, \lfloor (1-\delta/2)h \rfloor\Rr)\\
				& \qquad \qquad \cdot\P_{\p}\Ll[\lfloor (1-\delta/2)h \rfloor e_1 \xleftrightarrow{S_{0, h}}  [ \partial \Lambda_{h,L}]^{+}\Rr].
			\end{align*}
			Here we also relax the parameter of probability space in function of the domain. Afterwards, Lemma~\ref{lem.A.upperbound.inhomo} is applied to the inhomogeneous slab and every single term can be bounded respectively by
			\begin{align*}
				A_{\p}\Ll(0, \lfloor \delta h / 2 \rfloor; k, \lfloor \delta h / 2 \rfloor\Rr) &\stackrel{\text{\eqref{eq.A.upperbound.half}}}{\lesssim} \frac{(k+1)}{\delta h/ 2}, \\
				A_{\p, \tilde{\p}}\Ll(0, \lfloor \delta h / 2 \rfloor, \lfloor (1-\delta/2)h \rfloor; \lfloor \delta h / 2 \rfloor, \lfloor (1-\delta/2)h \rfloor\Rr) &\stackrel{\text{\eqref{eq.A.upperbound.inhomo}}}{\lesssim} \exp\Ll(-C^{-1} (1-\delta) h \sqrt{\p_c - \tilde{\p} }\Rr), \\
				\P_{\p}\Ll[\lfloor (1-\delta/2)h \rfloor e_1 \xleftrightarrow{S_{0, h}}  [ \partial \Lambda_{h,L}]^{+}\Rr] &\stackrel{\text{\eqref{eq.onearm_super}}}{\lesssim} \frac{4}{(\delta h)^2}.
			\end{align*}
			It concludes that 
			\begin{align}\label{eq.onearm_half_plus_delta_case1}
				\P_{\p, \tilde{\p}}^{\delta}\Ll[k e_1 \longleftrightarrow  [ \partial \Lambda_{h,L}]^{+}\Rr] \lesssim (k+1)(\delta h )^{-3}\exp\Ll(-C^{-1} (1-\delta) h \sqrt{\p_c - \tilde{\p} }\Rr).
			\end{align}
			
			\emph{Case~2: $\delta h /2 \leq k \leq h /2$.} The argument is similar as above, while we only need to use $u_2$ as the first hit at $L_{\lfloor (1-\frac{\delta}{2})h \rfloor}$. Then we break the path and  apply the BK inequality 
			\begin{align*}
				&\P_{\p, \tilde{\p}}^{\delta}\Ll[k e_1 \longleftrightarrow  [ \partial \Lambda_{h,L}]^{+}\Rr] \\
				&\leq \Ll(\sum_{u_2 \in L_{\lfloor (1-\frac{\delta}{2})h \rfloor}}  \P_{\p, \tilde{\p}}^{\delta}\Ll[ \{k e_1  \xleftrightarrow{S_{0, \lfloor (1-\frac{\delta}{2})h \rfloor}}  u_2\}\Rr]  \Rr) \P_{\p}\Ll[u_2 \xleftrightarrow{S_{0, h}}  [ \partial \Lambda_{h,L}]^{+}\Rr]\\
				&\leq A_{\p, \tilde{\p}}\Ll(0, \lfloor \delta h / 2 \rfloor, \lfloor (1-\delta/2)h \rfloor; k, \lfloor (1-\delta/2)h \rfloor\Rr) \P_{\p}\Ll[\lfloor (1-\delta/2)h \rfloor e_1 \xleftrightarrow{S_{0, h}}  [ \partial \Lambda_{h,L}]^{+}\Rr].
			\end{align*}
			The estimate of inhomogeneous slab is then applied
			\begin{align*}
				A_{\p, \tilde{\p}}\Ll(0, \lfloor \delta h / 2 \rfloor, \lfloor (1-\delta/2)h \rfloor; k, \lfloor (1-\delta/2)h \rfloor\Rr) &\stackrel{\text{\eqref{eq.A.upperbound.inhomo}}}{\lesssim} \exp\Ll(-(2C)^{-1} (1-\delta) h \sqrt{\p_c - \tilde{\p} }\Rr), \\
				\P_{\p}\Ll[\lfloor (1-\delta/2)h \rfloor e_1 \xleftrightarrow{S_{0, h}}  [ \partial \Lambda_{h,L}]^{+}\Rr] &\stackrel{\text{\eqref{eq.onearm_super}}}{\lesssim} \frac{4}{(\delta h)^2}.
			\end{align*}
			Her the range $\delta h /2 \leq k \leq h /2$ is also utilized. As a consequence, we obtain that
			\begin{equation}\label{eq.onearm_half_plus_delta_case2}
				\begin{split}
					\P_{\p, \tilde{\p}}^{\delta}\Ll[k e_1 \longleftrightarrow  [ \partial \Lambda_{h,L}]^{+}\Rr] &\lesssim (\delta h )^{-2}\exp\Ll(-(2C)^{-1} (1-\delta) h \sqrt{\p_c - \tilde{\p} }\Rr) \\
					&\lesssim (k+1) (\delta h )^{-3}\exp\Ll(-(2C)^{-1} (1-\delta) h \sqrt{\p_c - \tilde{\p} }\Rr).
				\end{split}
			\end{equation}
			The last inequality also relies on the range $\delta h /2 \leq k \leq h /2$. Therefore, \eqref{eq.onearm_half_plus_delta_case1} and \eqref{eq.onearm_half_plus_delta_case2} are unified, and they justify \eqref{eq.onearm_half_plus_delta}
		\end{proof}


		\section{Lower bound I: Within the correlation length}\label{sec.lower.critical}
		In the rest of the paper, our main focus turns to the lower bound for the conductivity $\ab(\p)$. As a warm-up, we show in this section a first lower bound at the critical point or in the correlation length. 
		\begin{prop}\label{prop.lower_water}
			In high-dimensional percolation, the following estimate holds 
			\begin{align}\label{eq.lower_water_general}
				\ab(\Lambda_{n}, \p_c) \gtrsim  n^{-6}.  
			\end{align}
			Especially, for every $\p \in (\p_c, 1]$ and $\ell_* = \lfloor (\p -\p_c)^{-1/2} \rfloor$, we have 
			\begin{align}\label{eq.lower_water}
				\ab(\Lambda_{\ell_*}, \p) \gtrsim  (\p-\p_c)^3.  
			\end{align}
		\end{prop}
		
		The argument will make rigorous the heuristic in Section~\ref{subsec.homogenization_intro}. We transform the study using the language of resistance in Section~\ref{subsec.resistance}, and prove the lower bound in Section~\ref{subsec.low_correlation}. Some heuristics and challenges about the generalization in large scale are also discussed in Section~\ref{subsec.discussion}, in order to motivate the analysis in the next sections.
		
		\subsection{Effective resistance}\label{subsec.resistance}
		The effective resistance is another notation in the literature, which is defined as the minimization of energy
		\begin{align}\label{eq.resistance}
			\Reff(A,Z) := \inf \Ll\{\mathscr{E}_*(\a,\theta): \text{unit flow } \theta \text{ from } A \text{ to } Z\Rr\},
		\end{align}
		where the energy of flow is defined as 
		\begin{align}\label{eq.def.energy.flow}
			\mathscr{E}_*(\a,\theta) := \sum_{e} \a^{-1}(e) \theta^2(e).
		\end{align}
		The flow is a function $\theta$ defined on the oriented edge such that 
		\begin{align*}
			\theta(\overrightarrow{xy}) = \theta(\overrightarrow{yx}),
		\end{align*}
		and a unit flow from $A$ to $Z$ satisfies that 
		\begin{align}\label{eq.def.flow}
			\begin{split}
				\sum_{y \sim x} \theta(\overrightarrow{xy}) = 0, \qquad \forall x \notin A \cup Z,\\
				\sum_{x \in A} \sum_{y \sim x}\theta(\overrightarrow{xy}) = \sum_{y \in Z} \sum_{x \sim y}\theta(\overrightarrow{xy}) = 1.
			\end{split}
		\end{align}
		The definition of \eqref{eq.resistance} also depends on the graph where it is defined, and we will explain it in the concrete context. We refer to \cite[Chapter~3.1]{LP17} for more properties about the effective resistance. 
		
		The key relation between the effective resistance in \eqref{eq.resistance} and the conductivity \eqref{eq.A} is stated as follows. The main interest of Lemma~\ref{lem.resistance_key} is to transfer the lower bound of the conductivity to the upper bound of the resistance.
		\begin{lemma}\label{lem.resistance_key}
			For every $n \in \N_+$, we have 
			\begin{align}\label{eq.resistance_key}
				\a(\Lambda_{n}) = \frac{n^{2-d}}{\Reff([\partial \Lambda_{n}]^{-}, [\partial \Lambda_{n}]^{+})}.
			\end{align}
			Moreover, for every $\a$-measurable event $A$, we have 
			\begin{align}\label{eq.resistance_key_mean}
				\ab(\Lambda_{n}, \p) \geq \frac{n^{2-d} \P_{\p}[A]^2}{\E_{\p}\Ll[\Reff([\partial \Lambda_{n}]^{-}, [\partial \Lambda_{n}]^{+}) \mathbf{1}_A\Rr]}.
			\end{align}
		\end{lemma}
		\begin{remark}
			When the graph is not connected, the resistance can be infinite. Then the restriction on the good event $A$ in \eqref{eq.resistance_key_mean} can help exclude such case and yield non-trivial lower bound. Besides, viewing \eqref{eq.resistance_key_mean}, it suffices to study the upper bound of the expectation of resistance.
		\end{remark}
		\begin{proof}[Proof of Lemma~\ref{lem.resistance_key}]
			In the literature of electric network, there exists another notation of effective conductivity defined as 
			\begin{align}\label{eq.conductivity_NEW}
				\Ceff(A,Z) := \inf \Ll\{\sum_{e} \a(e) v^2(e): v_{|A} = 0, v_{|Z}=1\Rr\},
			\end{align}
			and it satisfies (see \cite[Exercise~2.13]{LP17})
			\begin{align}\label{eq.inverse}
				\Reff(A,Z) &= \Ceff(A,Z)^{-1}.
			\end{align}
			We clarify the difference between \eqref{eq.conductivity_NEW} and \eqref{eq.def.aU}
			\begin{equation}\label{eq.resistance_long_calculation}
				\begin{split}
					\a(\Lambda_{n}) &= \inf_{v:v \equiv \ell_{e_1} \text{ on } [\partial \Lambda_{n}]^{\pm}} \frac{1}{|\Lambda_{n}|}\Ll \langle\nabla v,\a \nabla v\Rr \rangle_{\Lambda_{n}} \\
					& = n^2 \cdot (n^d)^{-1}\inf_{\substack{v:v \equiv 0 \text{ on } [\partial \Lambda_{n}]^{-}, \\
							v \equiv 1 \text{ on } [\partial \Lambda_{n}]^{+}}} \Ll \langle\nabla v,\a \nabla v\Rr \rangle_{\Lambda_{n}}\\
					& = n^{2-d} \Ceff([\partial \Lambda_{n}]^{-}, [\partial \Lambda_{n}]^{+})\\
					& = n^{2-d} \Reff([\partial \Lambda_{n}]^{-}, [\partial \Lambda_{n}]^{+})^{-1}.
				\end{split}
			\end{equation}
			The first line is just the definition \eqref{eq.def.aU}. The second line to the third line makes use of the scaling, and then it is reduced to the definitions \eqref{eq.conductivity_NEW} and \eqref{eq.resistance}. This finishes the proof of \eqref{eq.resistance_key}.
			
			To derive \eqref{eq.resistance_key_mean}, we take the expectation under $\P_{\p}$
			\begin{equation}\label{eq.lower_water_step1}
				\begin{split}
					\ab(\Lambda_{n}, \p) &\geq n^{2-d} \E_{\p}\Ll[\Reff([\partial \Lambda_{n}]^{-}, [\partial \Lambda_{n}]^{+})^{-1}\Rr] \\
					&\geq n^{2-d} \E_{\p}\Ll[\Reff([\partial \Lambda_{n}]^{-}, [\partial \Lambda_{n}]^{+})^{-1} \mathbf{1}_A\Rr].
				\end{split}
			\end{equation}
			From the first line to the second line, we just restrict the expectation on the event $A$. 	The Cauchy--Schwarz inequality implies that
			\begin{align*}
				\E_{\p}\Ll[\Reff([\partial \Lambda_{n}]^{-}, [\partial \Lambda_{n}]^{+})^{-1} \mathbf{1}_A\Rr] \E_{\p}\Ll[\Reff([\partial \Lambda_{n}]^{-}, [\partial \Lambda_{n}]^{+}) \mathbf{1}_A\Rr] \geq \P_{\p}[A]^2. 
			\end{align*}
			Therefore, \eqref{eq.lower_water_step1} can finally be reformulated as the desired result \eqref{eq.resistance_key_mean}.
		\end{proof}

		\subsection{Lower bound in correlation length}\label{subsec.low_correlation}
		\begin{proof}[Proof of Proposition~\ref{prop.lower_water}]
			The result \eqref{eq.lower_water} is the special case of \eqref{eq.lower_water_general} by using the monotonicity of $\p \mapsto \ab(\Lambda_n,\p)$ and by setting $n = \ell_* = \lfloor (\p -\p_c)^{-1/2} \rfloor$
			\begin{align*}
				\ab(\Lambda_{\ell_*}, \p) \geq \ab(\Lambda_{\ell_*}, \p_c) \gtrsim \ell_*^{-6} \geq (\p-\p_c)^3.
			\end{align*}
			We thus aim to prove \eqref{eq.lower_water_general}, which relies on \eqref{eq.resistance_key_mean}
			\begin{align}\label{eq.resistance_key_mean.2}
				\ab(\Lambda_{n}, \p_c) \geq \frac{n^{2-d} \P_{\p_c}[A]^2}{\E_{\p_c}\Ll[\Reff([\partial \Lambda_{n}]^{-}, [\partial \Lambda_{n}]^{+}) \mathbf{1}_A\Rr]}.
			\end{align}
			The problem is thus reduced to the upper bound of 
			\begin{align*}
				\E_{\p_c}\Ll[\Reff([\partial \Lambda_{n}]^{-}, [\partial \Lambda_{n}]^{+}) \mathbf{1}_A\Rr].
			\end{align*}
			Recall the spanning cluster introduced in Definition~\ref{def.spanning}. As we know, if there is no spanning cluster, then the resistance will be infinite, which is not expected. The restriction on event $A$ aims to exclude such case. We use the set $\mathscr{S}$ in \eqref{eq.def.spanning.diam}
			\begin{align}\label{eq.def.spanning.2}
				\mathscr{S} = \{\cC \subset \Lambda_n: \cC \text{ is } 2d\text{-spanning}, \, \rin(\cC;\Lambda_n) \leq C_2 n^2, \, c_0 n^4 < \vert \cC \vert < C_0 n^4 \},
			\end{align}
			and define the good event $A$
			\begin{align}\label{eq.many_spanning_cluster}
				A := \{ c_1 n^{d-6} \leq \vert \mathscr{S} \vert \leq C_1 n^{d-6}  \}.
			\end{align}
			By Lemma~\ref{lem.spanning.diam}, there exist finite positive numbers $c_0, C_0, c_1, C_1, C_2, N_0$ such that 
			\begin{align}\label{eq.lower_water_step2}
				\forall n \geq N_0, \qquad \P_{\p_c}[A] \geq \frac{1}{2}.
			\end{align}

			Viewing \eqref{eq.resistance}, we look for a flow to give an upper bound of the resistance. Here is our plan: we distribute the flow uniformly on all the $2d$-spanning clusters in $\mathscr{S}$, and let it run along the geodesic path connecting $[\partial \Lambda_{n}]^{-}$ and $[\partial \Lambda_{n}]^{+}$. We denote by $\gamma(\cC)$ the geodesic path for every $\cC \in \mathscr{S}$, that is 
			\begin{align}\label{eq.choice_theta}
				\theta(e) := \Ll\{\begin{array}{ll}
					\vert \mathscr{S} \vert^{-1}, & e \in \gamma(\cC), \cC \in \mathscr{S},\\
					0, & \text{otherwise.}
				\end{array}\Rr.
			\end{align}
			We then have
			\begin{align}\label{eq.resistance_representation}
				\Reff([\partial \Lambda_{n}]^{-}, [\partial \Lambda_{n}]^{+}) \leq \sum_{\cC \in \mathscr{S}} \sum_{e \in \gamma(\cC)} \theta^2(e) = \vert \mathscr{S} \vert^{-2} \sum_{\cC \in \mathscr{S}} \vert \gamma(\cC) \vert.
			\end{align}

			Notice that $\vert \gamma(\cC) \vert \leq \rin(\cC; \Lambda_n)$, which is the diameter of intrinsic distance in \eqref{eq.def.rin}. Thus, we take the expectation of \eqref{eq.resistance_representation} over the event $A$ and obtained
			\begin{equation}\label{eq.lower_water_step3}
				\begin{split}
					\E_{\p_c}\Ll[\Reff([\partial \Lambda_{n}]^{-}, [\partial \Lambda_{n}]^{+})\mathbf{1}_A\Rr] 
					&\leq \E_{\p_c}\Ll[\vert \mathscr{S} \vert^{-2}\sum_{\cC \in \mathscr{S}} \rin(\cC; \Lambda_n) \mathbf{1}_A\Rr]\\
					&\leq (c_1 n^{d-6})^{-2} \E_{\p_c}\Ll[\sum_{\cC \in \mathscr{S}} \rin(\cC; \Lambda_n)\Rr]\\
					&\leq (c_1 n^{d-6})^{-2} \cdot  C_1 n^{d-6} \cdot C_2 n^2 \\
					& = c_1^{-2} C_1 C_2 n^{8-d}.
				\end{split}
			\end{equation}
			Here the factors $c_1 n^{d-6}$ and $C_1 n^{d-6}$ come from the restriction over $A$ and \eqref{eq.many_spanning_cluster}, and $C_2 n^2$ is the upper bound of the intrinsic diameter $\rin(\cC; \Lambda_n)$ in \eqref{eq.def.spanning.2}. 
			
			We combine the results in \eqref{eq.lower_water_step3} and \eqref{eq.lower_water_step2}, and put them back to \eqref{eq.resistance_key_mean.2}, 
			\begin{align*}
				\ab(\Lambda_{n}, \p_c) \geq \frac{n^{2-d} \P_{\p_c}[A]^2}{\E_{\p_c}\Ll[\Reff([\partial \Lambda_{n}]^{-}, [\partial \Lambda_{n}]^{+}) \mathbf{1}_A\Rr]} \geq \frac{n^{2-d}}{4 c_1^{-2} C_1 C_2 n^{8-d}} \gtrsim n^{-6},
			\end{align*}
			then the proof of \eqref{eq.lower_water_general} is complete.
		\end{proof}

		\subsection{Some discussion for the lower bound in large scale}\label{subsec.discussion}
		Firstly,  Proposition~\ref{prop.lower_water} is not enough to conclude the lower bound of $\ab(\p)$. We recall that $\ab(\Lambda_n, \p)$ converges to $\ab(\p)$ from above, so the lower bound cannot be deduced just in the correlation length $\ell_{*}$ like the upper bound. The main task is then to go further to large scale, which is a study $\ab(\Lambda_n, \p)$ for $n \gg \ell_{*}$.
		
		Secondly, as Lemma~\ref{lem.resistance_key} transforms the study lower bound of $\ab(\Lambda_n, \p)$ to the upper bound of the resistance $\Reff$, one may expect a subadditivity for th latter. This is partially true because the resistance defines a metric, i.e. for all $x,y,z \in \Zd$
		\begin{align*}
			\Reff(x,z) \leq \Reff(x,y) + 	\Reff(y,z). 
		\end{align*}
		See \cite[Exercise~2.67]{LP17}. Here we also abuse the notation by identifying $\Reff(x,z) \equiv \Reff(\{x\}, \{z\})$. However, the triangle inequality does \emph{not} generalize to a set-valued version like 
		\begin{align*}
			\Reff(A,C) \leq \Reff(A,B) + 	\Reff(B,C).
		\end{align*}
		The main challenge is, the flow from $A$ to $B$ does not necessarily match that from $B$ to $C$. The extreme case can be $\Reff(A,B), \Reff(B,C) < \infty$, but $\Reff(A,C) = + \infty$.
		
		Another observation also hints the potential difficulty. If the proof in this section could conclude the lower bound, then what is the role of $\p$ in the proof ? Notice that the the property $\p > \p_c$ is nearly missing in the proof of Proposition~\ref{prop.lower_water}. 
		
\begin{figure}[h!]
			\centering
			\includegraphics[width=0.4\textwidth]{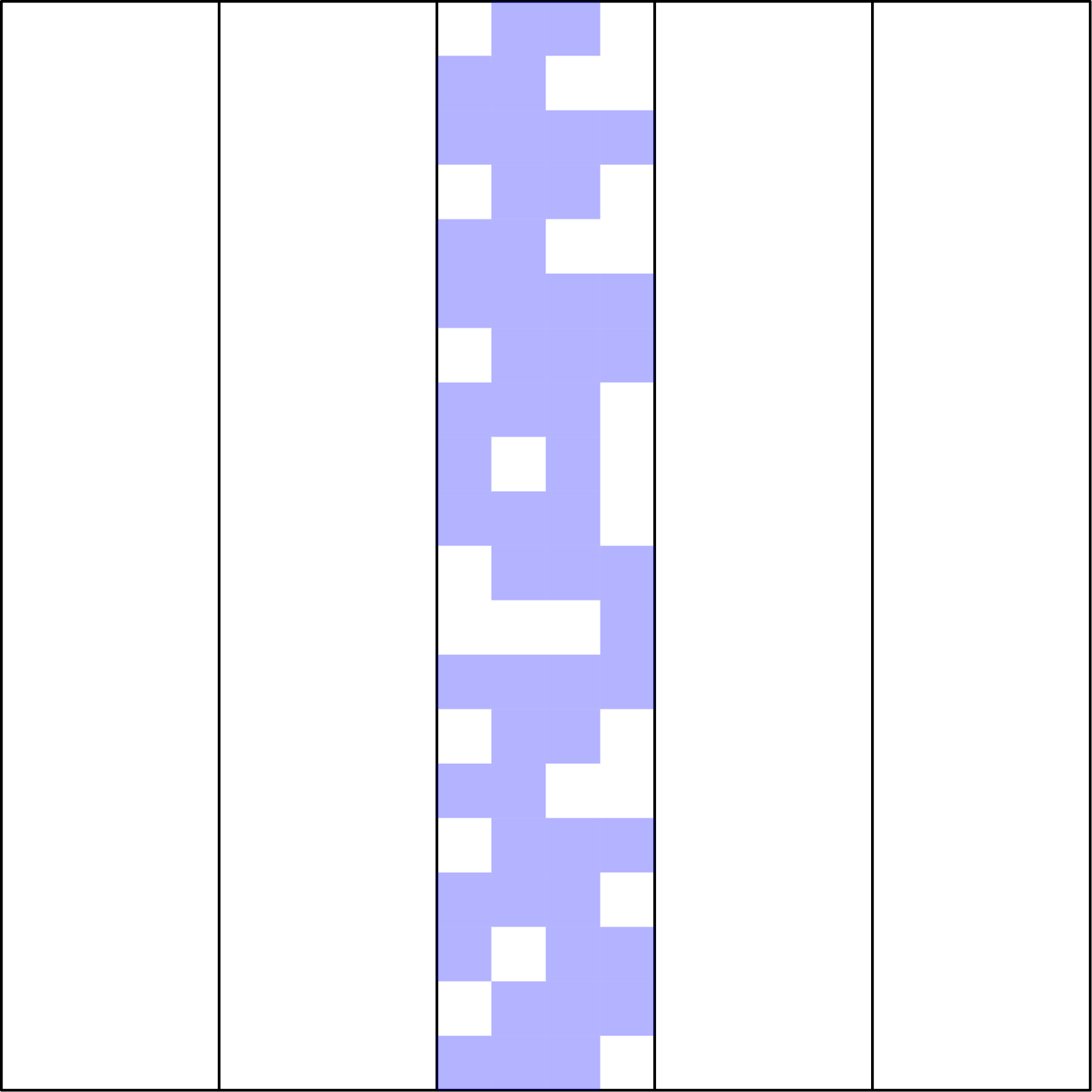}\hspace{1cm}
			\includegraphics[width=0.4\textwidth]{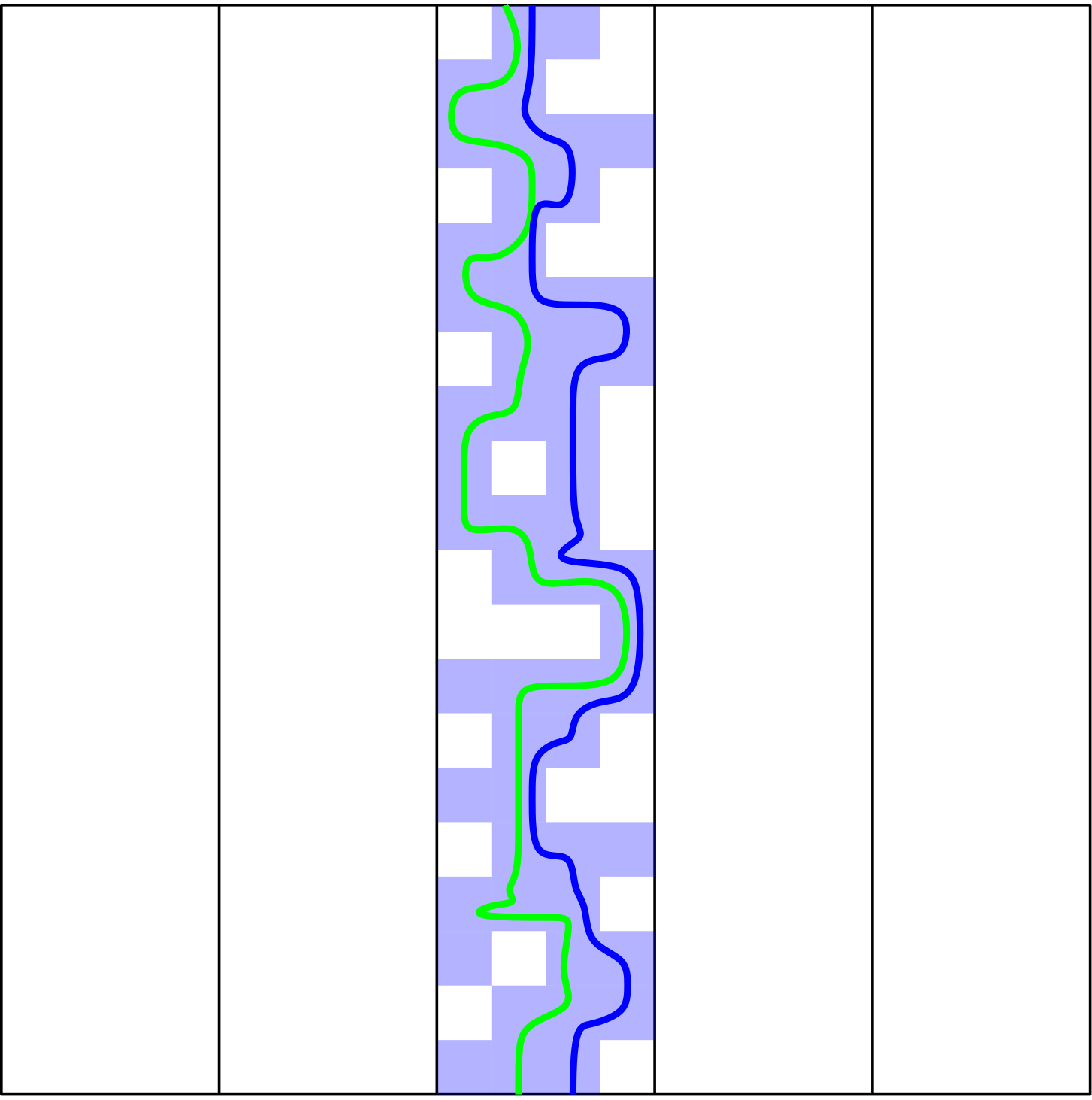}
			\caption{An illustration of sprinkling and coarse-graining. The pre-good cubes are already supercritical. Then the spanning clusters merge via the sprinkling property, and choose their own paths in it.}\label{fig.CoarseGrainNew}
		\end{figure}

		Finally, let us mention one key heuristic for the lower bound. It is \emph{the sprinkling property}, which tells us how the critical percolation evolves to supercritical percolation: let $\P_{\p}$ be the space of monotone coupling space of percolation. We sample $\p_c$-percolation first, then as the proof of Proposition~\ref{prop.lower_water}, we can typically, say with probability $(1-\epsilon)$, find $\ell_{*}^{d-6}$ spanning clusters in a cube $\Lambda_{\ell_{*}}$ satisfying correct resistance estimate. The $\p_c$-percolation is critical, so we do not see giant cluster. Afterward, an extra probability $(\p - \p_c)$ is sampled, which open more bonds. As a consequence, the spanning clusters in different cubes $\{\Lambda_{\ell_{*}}(z)\}_{z \in \ell_{*} \Zd}$ merge, and form the giant cluster.  See Figure~\ref{fig.CoarseGrainNew} for an illustration.
		
		The full and rigorous description about the passage of percolation from $\p_c$ to $\p > \p_c$ will be addressed in Section~\ref{sec.lower.large}, and the lower bound in large scale will be proved as a consequence. Before that, Sections~\ref{sec.regular}, \ref{sec.skl} are devoted to detailed geometric sprinkling properties between two cubes of correlation length.

		\section{Lower bound II: Regular boundary points}\label{sec.regular}
		
		This section is devoted to the study of the regular points, which are those well-separatedly  lying on the boundary of the spanning clusters. As illustrated in Figure~\ref{fig.span_skl} in the introduction, we just focus on one cube and its neighborhood in the section. For convenience, the following notations will be used in this section. 
		\begin{align}\label{eq.def.SQH_regular}
			Q := \Lambda_n, \qquad S:= \Lambda_{3n/4}, \qquad H:= \Lambda_{\cc n}.
		\end{align}
		$Q$ will be the basic cube in the renormalization in the future section. $S$ means a \emph{small} cube, and $H$ stands for the \emph{heart} of cube. We then denote by $\fC_{\cc}$
		\begin{align}\label{eq.def.fC_regular}
			\fC_{\cc} := \{y \in Q: y \xleftrightarrow{Q} H\},
		\end{align}
		the clusters explored from $H$ and restricted in $Q$. One major mission in this section is to determine the reasonable value of $\cc  \in (0, 1/5)$, so that various properties of regular points hold. 
		
\begin{figure}[h!]
			\centering
			\includegraphics[width=0.6\textwidth]{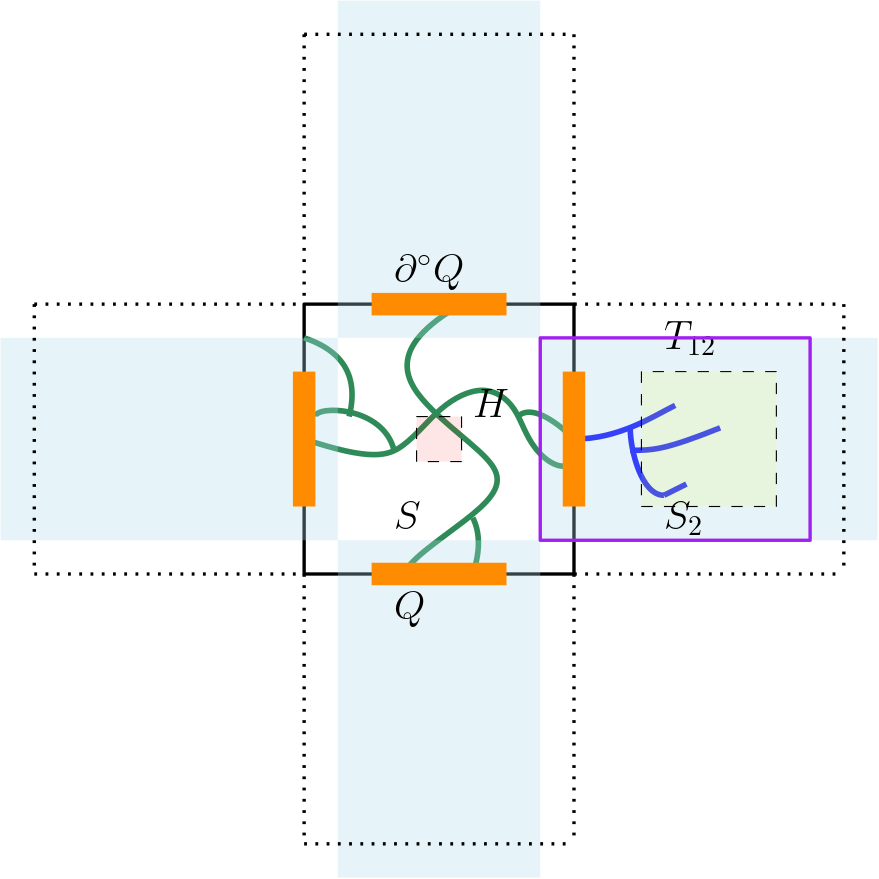}
			\caption{An illustration of the $(Q,S,H,T)$-structure. The heart $H$ of the cube is in light red, and the tube $T$ is in light blue. The boundary $\partial^\circ Q$ is in orange.}\label{fig.regular}
		\end{figure}

		The clusters in $\fC_{\cc}$ can be further extended outside $Q$. We restrict such exploration from the interior of the face in every direction $i \in \{ \pm 1, \pm 2, \cdots, \pm d\}$
		\begin{align}\label{eq.def.boundary.o}
			\begin{split}
				\partial_{i}^\circ Q &:= \Ll\{ y \in \partial_i Q:  y_j \in \Ll(-\frac{n}{4},  \frac{n}{4}\Rr), \forall j \neq \vert i \vert \Rr\}, \\
				\partial^\circ Q &:= \bigcup_{i \in \{ \pm 1, \pm 2, \cdots, \pm d\}}\partial_{i}^\circ Q.
			\end{split}
		\end{align}
		The exploration is also restricted in a tube. In the present section, we just focus on one direction $e_1$, thus the tube $T$ and its first segment $T_{12}$ are respectively defined by 
		\begin{align}\label{eq.def.tube_regular}
			\begin{split}
				T &:= n e_1 + \Big( -\frac{5n}{8}, \frac{5n}{8} \Big) \times \Big( -\frac{3n}{8}, \frac{3n}{8} \Big)^{d-1}, \\
				T_{12} &:= n e_1 + \Big( -\frac{5n}{8}, \frac{3n}{8} \Big) \times \Big( -\frac{3n}{8}, \frac{3n}{8} \Big)^{d-1}.
			\end{split}
		\end{align}
		We then set $S_2$ as another small cube in the interior of $T_{12}$
		\begin{align}\label{eq.def.S2}
			S_2 := n e_1 + \Lambda_{n/2}.
		\end{align}
		We will see, in the next section, that $S_2$ is the regime of pivotal bonds. The structure above is called $(Q,S,H,T)$-structure; see Figure~\ref{fig.regular} for an illustration. We also strongly suggest the readers to take a look at Figures~\ref{fig.ce1} and ~\ref{fig.span_skl} for motivations.

		

		In high dimension, a spanning cluster should contain enough points on $\partial Q$ such that the cluster has a good probability to be extended from each of them. We will call such points \emph{regular boundary points}. The version in this paper is characterized by two analytic inequalities.
		\begin{definition}[Regular points]\label{def.regular}
			Given the configuration $\fC_{\cc}$ and for every $r \in \N$ and $\delta >0$, a point $x \in \fC_{\cc} \cap \partial Q$ is called $(r,\delta)$-regular if the following two inequalities hold:
			\begin{equation}\label{def:ERE1}
				\ERE_1 (x, \fC_{\cc},r):= \sum_{y \in \fC_{\cc} \setminus  S} \max \{r, |x-y|\}^{-d+2} \leq \delta,
			\end{equation}
			\begin{equation}\label{def:ERE2}
				\ERE_2 (x, \fC_{\cc},r):= \sum_{y \in \fC_{\cc} \cap \partial Q} \max \{r, |x-y|\}^{-d+4} \leq \delta.
			\end{equation}
		\end{definition}
		The definition is a little technical, so let us preview this section at first. The parameters $\cc, r, \delta$ will also be determined along the way. 
		\begin{itemize}
			\item Section~\ref{subsec.regular.ext} introduces the key property about the regular points. Roughly speaking, as the neighborhood of the regular points is quite empty, it has a good probability to be the unique one at the boundary for the two-point estimate. This will also be one important input for Section~\ref{sec.skl}.
			\item Sections~\ref{subsec.regular.prob}-\ref{subsec.regular.var} contain the calculation of expectations and variance of $\ERE_1, \ERE_2$ above. The parameter $\delta$ is not used here.
			\item Section~\ref{subsec.regular.clt} aims to prove that the number of regular points is rich: typically $Q$ contains $n^{d-6}$ spanning clusters with $n^2$ regular boundary points. It relies not only the moment estimate, but also a careful choice of the parameters $\cc, r, \delta$.
		\end{itemize}

		\subsection{Properties of regular boundary points}\label{subsec.regular.ext}
		

		The most important property about the regular point is \emph{the unique extension property} as follows. It states that the regular point has a good probability to be the unique one to make the exploration in the tube.  
		\begin{prop}[Unique extension property]\label{prop.extension}
			There is $\delta_{\mathrm{ext}}(d) >0$ such that, for every $r \in \N_+$, there exists $c_{\ref{prop.extension}} (r)>0$ satisfying the following property: for sufficiently large $n \in \N$, every $ \cc  >0$ and $\delta \in (0, \delta_{\mathrm{ext}})$, every $z \in S_2$, all the $(r,\delta)$-regular boundary point $x \in \fC_{\cc} \cap \partial_{1}^\circ Q$ has the property
			\begin{equation}
				\P_{\p_c} \Ll[ \Big\{ x \stackrel{T_{12}\setminus \EDGE (\fC_{\cc})}{\longleftrightarrow}z \Big\} \setminus \bigcup_{y \in (\fC_{\cc} \cap \partial Q) \setminus \{x \}} \Big\{ x \stackrel{T_{12} \setminus \EDGE (\fC_{\cc})}{\longleftrightarrow} y \Big\} \bigg\vert \fC_{\cc} \Rr] \geq c_{\ref{prop.extension}} n^{-d+2}.
			\end{equation}
		\end{prop}
		\begin{proof}
			We shift $x$ by a $r e_1$ and denote by 
			\begin{align*}
				x^{+r}:= x+r  e_1,
			\end{align*}
			and consider the following two events:
			\begin{align*}
				A(x) &:= \Big\{ x \stackrel{T_{12}\setminus \EDGE (\fC_{\cc})}{\longleftrightarrow}z \Big\} \setminus \bigcup_{y \in \fC_{\cc} \cap \partial Q \setminus \{x\}} \Big\{ x \stackrel{T_{12} \setminus \EDGE (\fC_{\cc})}{\longleftrightarrow} y \Big\}, \\
				A(x^{+r}) &:= \Big\{ x^{+r} \stackrel{T_{12}\setminus \EDGE (\fC_{\cc})}{\longleftrightarrow}z \Big\} \setminus \bigcup_{y \in \fC_{\cc} \cap \partial Q} \Big\{ x^{+r} \stackrel{T_{12} \setminus \EDGE (\fC_{\cc})}{\longleftrightarrow} y \Big\}.
			\end{align*}
			We will study $A(x)$ via the  event $A(x^{+r})$, which benefits from the shift.
			
			\medskip
			\textit{Step 1: A surgery lemma.} We first show that there is $c(r)>0$, depending only on $d$ and $r$, such that
			\begin{equation}\label{eq:surgery}
				\P_{\p_c} \big[ A(x) \mid \fC_{\cc} \big] \geq c(r)\, \P_{\p_c} \big[ A(x^{+r}) \mid \fC_{\cc} \big].
			\end{equation}
			
			Let $I:=I(x,r) := \{ x, x+ e_1, x+2 e_1, \dots, x+r e_1 \}$ denote the $r+1$ vertices on the segment from $x$ to $x^{+r}$, and let
			\begin{equation}\label{eq:Udef}
				U:=U(x,r) := \big\{ e: e \text{ is an edge in } T_{12} \setminus \EDGE (\fC_{\cc}) \text{ with at least one endpoint in } I \big\}
			\end{equation}
			In particular $U$ contains the $r$ edges of the segment $x \to x^{+r}$ itself, and $|U|$ is finite and depends only on $d$ and $r$ (indeed $|U| \leq (r+1) \cdot 2d$). See Figure~\ref{fig.surgery}(a) for an illustration of $x$, $x^{+r}$, $I$, and $U$.
			
			Given $\fC_{\cc}$, write $\omega_U$ and $\omega_{U^c}$ for the restrictions of the percolation configuration to $U$ and to its complement $U^c := T_{12} \setminus (\fC_{\cc} \cup U)$; these two restrictions are independent, and $\omega_U$ is a family of i.i.d. $\mathrm{Bernoulli}(\p_c)$ variables. Define the event
			\begin{equation}
				B(x) := \Big\{ \omega_{U^c} \ : \ \exists\, \eta \in \{0,1\}^{U} \text{ such that } A(x) \text{ occurs under } (\omega_{U^c}, \eta) \Big\}.
			\end{equation}
			By construction, $B(x)$ is measurable with respect to $\omega_{U^c}$ alone.
			
			\smallskip
			\noindent \emph{Claim: $A(x^{+r}) \subset B(x)$.}
			
			Let $\omega = (\omega_{U^c}, \omega_U)$ be such that $A(x^{+r})$ occurs and let $\gamma$ be an open path in $T_{12} \setminus \fC_{\cc}$ from $x^{+r}$ to $z$ under the configuration $\omega$ (see Figure~\ref{fig.surgery}(b)). Define $\eta \in \{0,1\}^U$ as follows:
			\begin{itemize}
				\item every edge on the segment from $x$ to $x^{+r}$ is declared \emph{open};
				\item the last edge $e_\gamma$ in $\gamma$ which is also in $U$ remains \emph{open};
				\item every other edge of $U$ is declared \emph{closed}.
			\end{itemize}
			Under the configuration $(\omega_{U^c}, \eta)$, all the vertices of $I$ are connected to $z$ in $T_{12} \setminus \fC_{\cc}$ through the edge $e_\gamma$ and the part of $\gamma$ after $e_\gamma$.
			Moreover, if we let $e_\gamma^+$ be the endpoint of $e_\gamma$ not in $I$, then the cluster of $e_\gamma^+$ inside $T_{12} \setminus (\fC_{\cc} \cup U)$ under the configuration $\omega$ is the same as the one under the configuration $(\omega_{U^c},\eta)$, which is smaller than the cluster of $e_\gamma^+$ inside $T_{12} \setminus \fC_{\cc}$ under the configuration $\omega$. Since $e_\gamma$ is the only open edge in $(\omega_{U^c},\eta)$ to go out from $I$, we have
			\begin{equation}\label{eihnbu24}
				K_x (\omega_{U^c},\eta) =K_{e_\gamma^+} (\omega_{U^c},\eta) \subseteq K_{e_\gamma^+} (\omega) \cup I =K_{x^{+r}} (\omega) \cup I
			\end{equation}
			where $K_x (\omega_{U^c},\eta)$ is the cluster of $x$ in $T_{12} \setminus \fC_{\cc}$ under the configuration $(\omega_{U^c},\eta)$, and another two notations are similar.
			Since $A(x^{+r})$ occurs for $\omega$, we have $y \notin K_{x^{+r}} (\omega)$ for every $y \in \fC_{\cc} \cap \partial Q$. Then \eqref{eihnbu24} implies $y \notin K_{x} (\omega_{U^c}, \eta)$ for every $y \in \fC_{\cc} \cap \partial Q \setminus \{x\}$. Hence $A(x)$ occurs for $(\omega_{U^c}, \eta)$, proving the claim (see Figure~\ref{fig.surgery}(c) for an illustration of the resulting configuration).
			
			\begin{figure}[h!]
				\centering
				\includegraphics[width=0.8\textwidth]{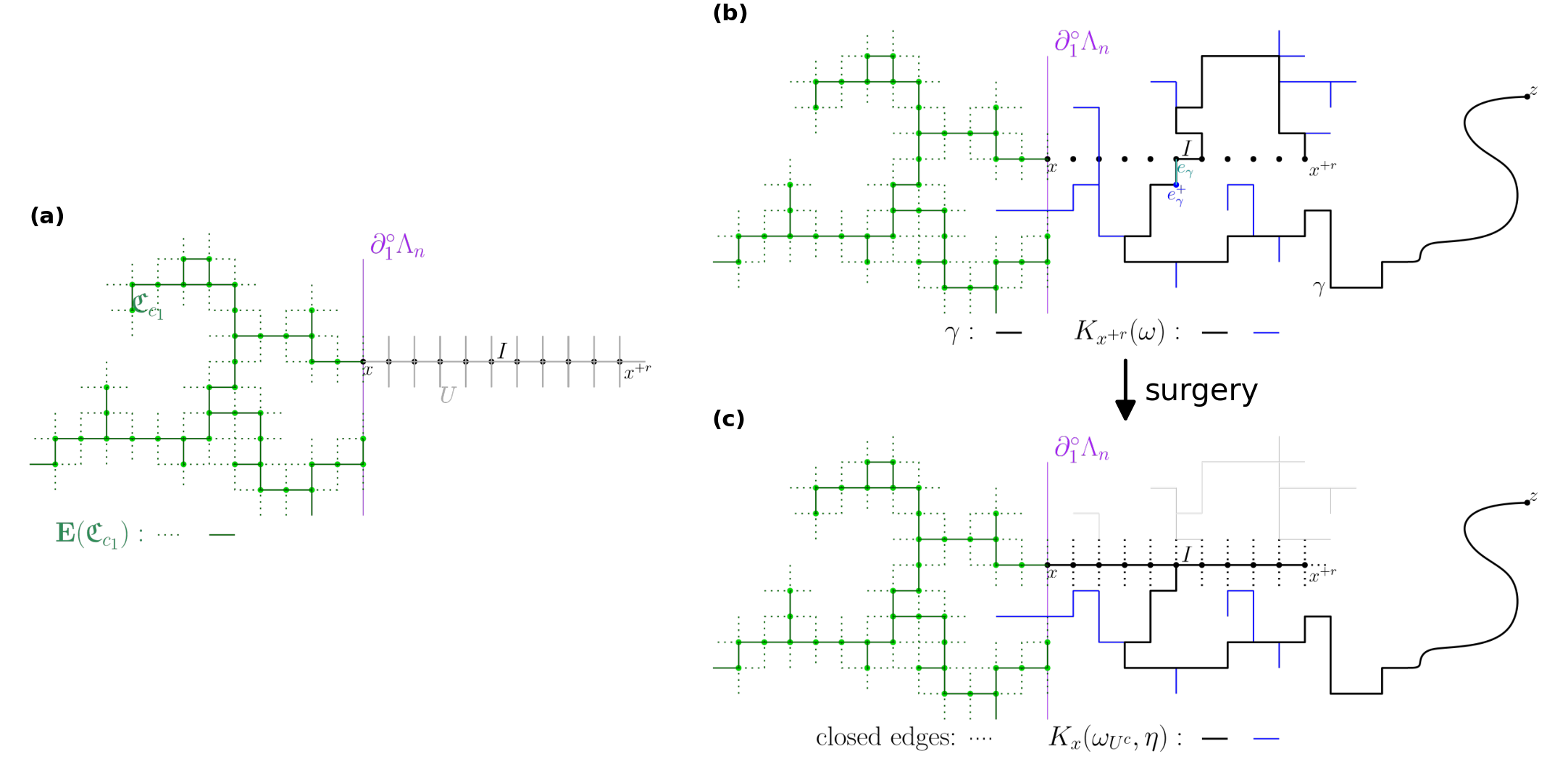}
				\caption{Illustration of the surgery construction in Step 1. (a) The cluster $\fC_{\cc}$, the boundary segment $\partial_{1}^\circ \Lambda_n$, the point $x$, the segment $I$ from $x$ to $x^{+r}$, and the edge set $U$ (the tube around $I$). (b) A sample configuration $\omega$ realizing $A(x^{+r})$: the path $\gamma$ from $x^{+r}$ to $z$, its last edge $e_\gamma$ lying in $U$, the far endpoint $e_\gamma^+$, and the resulting cluster $K_{x^{+r}}(\omega)$. (c) The surgered configuration $(\omega_{U^c},\eta)$: closing every edge of $U$ except $e_\gamma$ and opening the segment $I$ connects all of $I$ (in particular $x$) to $z$, giving the cluster $K_x(\omega_{U^c},\eta) \subseteq K_{x^{+r}}(\omega)\cup I$.}\label{fig.surgery}
			\end{figure}
			
			\smallskip
			Since $U$ is finite, for each $\omega_{U^c} \in B(x)$ we may fix a distinguished witnessing configuration $\eta(\omega_{U^c}) \in \{0,1\}^U$ (e.g. the lexicographically smallest one for which $A(x)$ occurs), and partition
			\begin{equation}
				B(x) = \bigsqcup_{\eta \in \{0,1\}^U} B_\eta, \qquad B_\eta := \big\{ \omega_{U^c} \in B(x) : \eta(\omega_{U^c}) = \eta \big\}.
			\end{equation}
			By definition, $A(x)$ occurs on $B_\eta \times \{\omega_U = \eta\}$. Using the independence of $\omega_U$ and $\omega_{U^c}$ given $\fC_{\cc}$,
			\begin{align*}
				\P_{\p_c} \big[ A(x) \mid \fC_{\cc} \big] 
				&\geq \sum_{\eta \in \{0,1\}^U} \P_{\p_c} \big[ \omega_{U^c} \in B_\eta \mid \fC_{\cc} \big] \, \P_{\p_c} \big[ \omega_U = \eta \mid \fC_{\cc} \big] \\
				&\geq \Big( \min_{\eta \in \{0,1\}^U} \P_{\p_c} \big[ \omega_U = \eta \mid \fC_{\cc} \big] \Big) \sum_{\eta \in \{0,1\}^U} \P_{\p_c} \big[ \omega_{U^c} \in B_\eta \mid \fC_{\cc} \big] \\
				&= \min\{ \p_c, 1-\p_c \}^{|U|} \, \P_{\p_c} \big[ B(x) \mid \fC_{\cc} \big] \\
				&\geq \min\{ \p_c, 1-\p_c \}^{|U|} \, \P_{\p_c} \big[ A(x^{+r}) \mid \fC_{\cc} \big],
			\end{align*}
			where the last inequality follows from the Claim. Since $|U|$ depends only on $d$ and $r$, setting $c(r) := \min\{ \p_c, 1-\p_c \}^{|U|}$ establishes \eqref{eq:surgery}.
			
			\medskip
			\textit{Step 2: Estimate of $\P_{\p_c} [ A(x^{+r}) \mid \fC_{\cc} ]$.} By BK inequality we have
			\begin{align}\label{nqo803bj}
				\begin{split}
					\P_{\p_c} \Ll[ x^{+r} \stackrel{T_{12}\setminus \EDGE (\fC_{\cc})}{\longleftrightarrow}z \Big\vert \fC_{\cc} \Rr] &\geq \P_{\p_c} \Ll[ x^{+r} \stackrel{T_{12}}{\longleftrightarrow}z  \Rr] -\sum_{y \in \fC_{\cc} \cap T_{12}} \P_{\p_c} \Ll[ \{x^{+r} \stackrel{T_{12}}{\longleftrightarrow}y\} \circ \{y \stackrel{T_{12}}{\longleftrightarrow}z\}  \Rr]\\
					&\geq \P_{\p_c} \Ll[ x^{+r} \stackrel{T_{12}}{\longleftrightarrow}z  \Rr] -\sum_{y \in \fC_{\cc} \cap T_{12}} \P_{\p_c} \Ll[ x^{+r} \longleftrightarrow y  \Rr] \P_{\p_c} \Ll[ y \longleftrightarrow z  \Rr].
				\end{split}
			\end{align}
			Here, in the second term, we decompose $x^{+r} \stackrel{T_{12}}{\longleftrightarrow}z$ by its last intersection at some $y \in \fC_{\cc} \cap T_{12}$. We then need to estimate its outcomes.
			
			Concerning the first term in \eqref{nqo803bj},  we notice that,  for sufficiently large $n$, $x^{+r}$ and $z$ are in $T_{12}$ and their distance to $\partial T_{12}$ are both at least $\frac{1}{8}n-1$. Moreover,
			\begin{equation}
				\vert x^{+r} -z \vert \geq \vert x-z\vert -r \geq \frac{n}{4}-r.
			\end{equation}
			Thus, Lemma~\ref{lem.box_twopoint}, applied after translating $x^{+r}$ to the origin (with $M,(\alpha_i),(\beta_i)$ determined by the fixed shape of $T_{12}$, so that the resulting constant depends only on $d$): there is $c_{\ref{lem.box_twopoint}} >0$ such that for any sufficiently large $n$,
			\begin{equation}\label{vlsohe}
				\P_{\p_c} \Ll[ x^{+r} \stackrel{T_{12}}{\longleftrightarrow}z  \Rr] \geq c_{\ref{lem.box_twopoint}} n^{-d+2}.
			\end{equation}
			
			Concerning the second term in \eqref{nqo803bj}, because 
			\begin{align*}
				\vert x^{+r}-y\vert \geq \max \{r, \vert x-y \vert\}, \qquad \vert y-z\vert \geq z_1 -y_1 \geq \frac{n}{4},
			\end{align*}
			there is $C_1  >0$ such that
			\begin{equation}\label{nrbow}
				\P_{\p_c} \Ll[ y \longleftrightarrow z  \Rr] \leq C_1 n^{-d+2},\qquad \P_{\p_c} \Ll[ x^{+r} \longleftrightarrow y  \Rr] \leq C_1 \max \{r, \vert x-y \vert\}^{-d+2}.
			\end{equation}
			Applying \eqref{nrbow} and \eqref{vlsohe} into the inequality \eqref{nqo803bj}, we get
			\begin{equation}\label{vbnoeonn}
				\P_{\p_c} \Ll[ x^{+r} \stackrel{T_{12}\setminus \EDGE (\fC_{\cc})}{\longleftrightarrow}z \Big\vert \fC_{\cc} \Rr] \geq (c_{\ref{lem.box_twopoint}} -C_1 \ERE_1 (x,\fC_{\cc},r) \big) n^{-d+2}.
			\end{equation}
			
			Using \eqref{vbnoeonn} and Lemma~\ref{lem.threepoint_upper}, we have
			\begin{align*}
				&\P_{\p_c} \Ll[ A(x^{+r}) \Big\vert \fC_{\cc} \Rr] \\
				&\geq \P_{\p_c} \Ll[ x^{+r} \stackrel{T_{12}\setminus \EDGE (\fC_{\cc})}{\longleftrightarrow}z \Big\vert \fC_{\cc} \Rr] -\sum_{y \in \fC_{\cc} \cap \partial Q} \P_{\p_c} \Ll[ x^{+r} \stackrel{T_{12}\setminus \EDGE (\fC_{\cc})}{\longleftrightarrow}z \stackrel{T_{12}\setminus \EDGE (\fC_{\cc})}{\longleftrightarrow}y \Big\vert \fC_{\cc} \Rr]\\
				& \geq (c_{\ref{lem.box_twopoint}} -C_1 \ERE_1 (x,\fC_{\cc},r) \big) n^{-d+2} -C_{\ref{lem.threepoint_upper}} \sum _{y \in \fC_{\cc} \cap \partial Q} n^{-d+2} \vert x^{+r} -y\vert^{-d+4}\\
				& \geq (c_{\ref{lem.box_twopoint}} -C_1 \ERE_1 (x,\fC_{\cc},r) \big) n^{-d+2} -C_{\ref{lem.threepoint_upper}} \sum _{y \in \fC_{\cc} \cap \partial Q} n^{-d+2} \max \{r, \vert x -y\vert\}^{-d+4}\\
				& \geq (c_{\ref{lem.box_twopoint}} -C_1 \ERE_1 (x,\fC_{\cc},r) -C_{\ref{lem.threepoint_upper}} \ERE_2 (x,\fC_{\cc},r) \big) n^{-d+2}.
			\end{align*}
			Let $\delta_{\mathrm{ext}} := \dfrac{c_{\ref{lem.box_twopoint}}}{4\max\{C_1, C_{\ref{lem.threepoint_upper}}\}}$. Then, for all $\delta \in (0, \delta_{\mathrm{ext}})$ and  $x$ as $(r,\delta)$-regular point, we have
			\begin{equation}\label{eq:step2conclusion}
				\P_{\p_c} \big[ A(x^{+r}) \mid \fC_{\cc} \big] \geq \frac{c_{\ref{lem.box_twopoint}}}{2} n^{-d+2}.
			\end{equation}
			
			Finally, combining \eqref{eq:surgery} and \eqref{eq:step2conclusion},
			\begin{equation*}
				\P_{\p_c} \big[ A(x) \mid \fC_{\cc} \big] \geq c(r) \cdot \frac{c_{\ref{lem.box_twopoint}}}{2} n^{-d+2}.
			\end{equation*}
			Since $A(x)$ is exactly the event in the statement of the proposition, setting $c_{\ref{prop.extension}} (r) := \frac{c_{\ref{lem.box_twopoint}}}{2}\, c(r)$ completes the proof.
		\end{proof}	
		
		To end this subsection, we prove the following basic estimate, which will be used later in this paper.
		\begin{lemma}\label{lem.regular.baby}
			For every $(r,\delta)$-regular point $x \in \fC_{\cc} \cap \partial Q$, it satisfies that
			\begin{equation}\label{eq.regular.baby}
				\sum_{y \in (\fC_{\cc} \cap \partial Q) \setminus \{ x \}} |x-y|^{-d+4} \leq \delta r^{d-4}.
			\end{equation}
		\end{lemma}
		
		\begin{proof}
			For $x,y \in \Zd$ with $x \neq y$, we observe that
			\begin{align*}
				| x-y | \geq \frac{1}{r} \max \{ r, |x-y| \},
			\end{align*}
			which implies that
			\begin{equation*}
				|x-y|^{-d+4} \leq r^{d-4} \max \{ r, |x-y| \}^{-d+4}.
			\end{equation*}
			Then \eqref{def:ERE2} implies \eqref{eq.regular.baby}.
		\end{proof}

		\subsection{The upper and lower bounds for $\P_{\p_c} [x \in \fC_{\cc}]$}\label{subsec.regular.prob}
		In this subsection, we prove the following proposition.
		
		\begin{prop}\label{PxinfC}
			There exist $c_{\ref{PxinfC}}, C_{\ref{PxinfC}} >0$ such that, for every $ \cc  \in (0, 1/5)$ and every $n\in\N_+$,
			\begin{equation}\label{eq.PxinfC_upper}
				\P_{\p_c} [u \in \fC_{\cc}] \leq C_{\ref{PxinfC}}\,(1+l_u)\, \cc ^{d-4} n^{-3}, \qquad \forall\,u\in Q\setminus\Lambda_{n/2},
			\end{equation}
			and
			\begin{equation}\label{veiouvoh}
				\P_{\p_c} [u \in \fC_{\cc}] \geq c_{\ref{PxinfC}}\, \cc ^{d-4} n^{-3}, \qquad \forall\,u\in\partial_{1}^\circ  Q.
			\end{equation}
		\end{prop}
		
		This proposition directly implies
		\begin{corollary}\label{cor.PxinfC}
			For every $ \cc  \in (0, 1/5)$ and every $x \in \partial_{1}^\circ  Q$,
			\begin{equation*}
				c_{\ref{PxinfC}}\,  \cc ^{d-4} n^{-3} \leq \P_{\p_c} [x \in \fC_{\cc}] \leq C_{\ref{PxinfC}}\,  \cc ^{d-4} n^{-3}.
			\end{equation*}
		\end{corollary}
		
		\begin{proof}
			We first prove the upper bound \eqref{eq.PxinfC_upper}, for $u\in Q\setminus\Lambda_{n/2}$; the lower bound \eqref{veiouvoh} occupies the remainder of this subsection. In the proof, we will use 
			\begin{align*}
				l_x:=n-\vert x\vert,
			\end{align*}
			for the distance from $x$ to $\partial Q$.

			Note that once $u \in \fC_{\cc}$ happens, there is a path from $u$ to $\partial  H$, and, since $\vert u\vert>n/2>2 \cc n$, this path must have a last intersection with $\partial \Lambda_{2 \cc n}$. Splitting the path at this point and applying the BK inequality, we obtain
			\begin{equation}\label{meoucweu}
				\P_{\p_c} [u \in \fC_{\cc}] \leq \sum_{y \in \partial \Lambda_{2 \cc n}} \P_{\p_c} \big[ u \stackrel{ Q}{\longleftrightarrow} y \big] \P_{\p_c} \big[ y \stackrel{\Lambda_{2 \cc n}}{\longleftrightarrow} \partial  H \big].
			\end{equation}
			
			Since $u \in  Q\setminus\Lambda_{n/2}$ and $y \in \Lambda_{2 \cc n} \subseteq \Lambda_{2n/5}$, we have $|u-y| \geq n/10$. Lemma~\ref{lem.box_twopoint_onesided} gives, with $C_{\ref{lem.box_twopoint_onesided}}>0$ as in that lemma,
			\begin{equation}\label{nevo02c}
				\P_{\p_c} \big[ u \stackrel{ Q}{\longleftrightarrow} y \big] \leq C_{\ref{lem.box_twopoint_onesided}}\,(1+l_u)\,\vert u-y\vert^{1-d} \lesssim  (1+l_u)\,n^{1-d}.
			\end{equation}
			For $y \in \partial \Lambda_{2 \cc n}$ with $y \cdot e_1 =-2 \cc n$, let $\H_y: =\{ z \in \ZZ^d: z \cdot  e_1 \geq -2 \cc n \} \supseteq \Lambda_{2 \cc n}$. If $y \stackrel{\Lambda_{2 \cc n}}{\longleftrightarrow} \partial  H$ happens, then $y \stackrel{\H_y}{\longleftrightarrow} \partial  H$ also happens. Since $y \in \partial \Lambda_{2 \cc n}$, the distance between $y$ and $\partial  H$ is at least $ \cc n$, so
			\begin{equation*}
				\P_{\p_c} \big[ y \stackrel{\Lambda_{2 \cc n}}{\longleftrightarrow} \partial  H \big] \leq \P_{\p_c} \big[ y \stackrel{\H_y}{\longleftrightarrow} \partial  H \big] \leq \P_{\p_c} \big[ y \stackrel{\H_y}{\longleftrightarrow} \H_y \setminus (y + H) \big].
			\end{equation*}	
			By Lemma~\ref{lem.onearm_half}, with $C_{\ref{lem.onearm_half}}>0$ as in that lemma,
			\begin{equation}\label{ooppjoi}
				\P_{\p_c} \big[ y \stackrel{\Lambda_{2 \cc n}}{\longleftrightarrow} \partial  H \big] \leq \P_{\p_c} \big[ y \stackrel{\H_y}{\longleftrightarrow} \H_y \setminus (y + H) \big] \leq C_{\ref{lem.onearm_half}}\,( \cc n)^{-3},
			\end{equation}
			which, by symmetry, holds for all $y \in \partial \Lambda_{2 \cc n}$.
			Applying \eqref{nevo02c} and \eqref{ooppjoi} into \eqref{meoucweu} and noting that there are $\sim_d ( \cc n)^{d-1}$ points in $\partial \Lambda_{2 \cc n}$, we get
			\begin{equation*}
				\P_{\p_c} [u \in \fC_{\cc}] \lesssim ( \cc n)^{d-1}\, (1+l_u)\,n^{1-d}\, ( \cc n)^{-3} =(1+l_u)\, \cc ^{d-4} n^{-3},
			\end{equation*}
			which proves the upper bound \eqref{eq.PxinfC_upper}.

			For the lower bound, we fix $x \in \partial_{1}^\circ  Q$ and define in this subsection
			\begin{equation}
				Y:= \sum_{y \in \partial  H} \mathbf{1} (x \stackrel{ Q}{\longleftrightarrow} y).
			\end{equation}
			Then $\{x\in\fC_{\cc}\} =\{Y>0\}$ and the Cauchy--Schwarz inequality gives
			\begin{equation}\label{eq.PxinfC_CS}
				\P_{\p_c}[x\in\fC_{\cc}]=\P_{\p_c}[Y>0]\geq\frac{\E_{\p_c}[Y]^2}{\E_{\p_c}[Y^2]}.
			\end{equation}
			Then the lower bound in \eqref{veiouvoh} immediately follows from a lower bound on $\E_{\p_c}[Y]$ (Lemma~\ref{lem.Y_firstmoment}) and an upper bound on $\E_{\p_c}[Y^2]$ (Lemma~\ref{lem.Y_secondmoment}).
		\end{proof}

		\begin{lemma}\label{lem.Y_firstmoment}
			There exists $c_{\ref{lem.Y_firstmoment}}>0$ such that, for every $ \cc \in(0,1/5)$, every $n\in\N_+$, and every $x\in\partial_{1}^\circ  Q$,
			\begin{equation*}
				\E_{\p_c}[Y]\geq c_{\ref{lem.Y_firstmoment}}\, \cc ^{d-1}.
			\end{equation*}
		\end{lemma}
		
		\begin{proof}
			For $y\in\partial H$ we have $\vert y\vert= \cc n$, so $l_y=n- \cc n\geq(1-1/5)n\geq n/2$. Moreover, since $x\cdot e_1=n$ while $x\cdot e_j\in[n/4,n/2]$ for $j\neq1$, the distance from $x$ to $\partial\partial Q$ (the set of points in $\partial Q$ with at least two coordinates equal to $\pm n$) is at least $n/2$. Proposition~\ref{lem.box_twopoint_lower}, applied with $\alpha=1/2$, then gives, with $c_{\ref{lem.box_twopoint_lower}}(1/2)>0$ as in that proposition,
			\begin{equation}\label{bnrpoq}
				\P_{\p_c}\Ll[x\stackrel{ Q}{\longleftrightarrow}y\Rr]\geq c_{\ref{lem.box_twopoint_lower}}\,n^{1-d},\qquad\forall\,y\in\partial H.
			\end{equation}
			Since $\vert\partial H\vert\geq c_2\,( \cc n)^{d-1}$ for some $c_2>0$, summing \eqref{bnrpoq} over $y\in\partial H$ gives
			\begin{equation*}
				\E_{\p_c}[Y]\geq c_{\ref{lem.box_twopoint_lower}}\,n^{1-d}\,\vert\partial H\vert\geq c_{\ref{lem.box_twopoint_lower}}\,c_2\, \cc ^{d-1},
			\end{equation*}
			which proves the lemma with $c_{\ref{lem.Y_firstmoment}}:=c_{\ref{lem.box_twopoint_lower}}(1/2)\,c_2$.
		\end{proof}
		
		\begin{lemma}\label{lem.Y_secondmoment}
			There exists $C_{\ref{lem.Y_secondmoment}}>0$ such that, for every $ \cc \in(0,1/5)$, every $n\in\N_+$, and every $x\in\partial_{1}^\circ  Q$,
			\begin{equation*}
				\E_{\p_c}[Y^2]\leq C_{\ref{lem.Y_secondmoment}}\, \cc ^{d+2}\,n^3.
			\end{equation*}
		\end{lemma}

		\begin{proof}
			Fix distinct $y,y'\in\partial H$ and write $k:=1+ \vert y-y'\vert$. By the BK inequality, decomposing $\{x\stackrel{ Q}{\longleftrightarrow}y,\ x\stackrel{ Q}{\longleftrightarrow}y'\}$ at the branch point $u\in Q$ of the two connections,
			\begin{equation}\label{ekbnbd}
				\P_{\p_c}\Ll[x\stackrel{ Q}{\longleftrightarrow}y,\ x\stackrel{ Q}{\longleftrightarrow}y'\Rr]\leq\sum_{u\in Q}\P_{\p_c}\Ll[x\stackrel{ Q}{\longleftrightarrow}u\Rr]\,\P_{\p_c}\Ll[u \longleftrightarrow y\Rr]\,\P_{\p_c}\Ll[u \longleftrightarrow y'\Rr],
			\end{equation}
			We split this sum according to $\vert x-u\vert\gtrless n/4$.
			
			\emph{Case $\vert x-u\vert\geq n/4$.}  By Lemma~\ref{lem.box_twopoint_onesided},
			\begin{equation*}
				\P_{\p_c}\Ll[x\stackrel{ Q}{\longleftrightarrow}u\Rr]\leq C_{\ref{lem.box_twopoint_onesided}}\,\vert x-u\vert^{1-d}\leq C_{\ref{lem.box_twopoint_onesided}}\,(n/4)^{1-d}\lesssim n^{1-d},
			\end{equation*}
			uniformly over such $u$. Bounding the remaining sum over $u$ by the sum over all of $\Zd$,
			\begin{equation*}
				\sum_{\substack{u\in Q\\ \vert x-u\vert\geq n/4}}\P_{\p_c}\Ll[x\stackrel{ Q}{\longleftrightarrow}u\Rr]\P_{\p_c}\Ll[u \longleftrightarrow y\Rr]\P_{\p_c}\Ll[u \longleftrightarrow y'\Rr]\lesssim n^{1-d}\sum_{u\in\Zd}\P_{\p_c}\Ll[u \longleftrightarrow y\Rr]\P_{\p_c}\Ll[u \longleftrightarrow y'\Rr].
			\end{equation*}
			There is $C_1>0$ such that $\P_{\p_c}[v \longleftrightarrow w]\leq C_1\,\vert v-w\vert^{2-d}$ for every $v\neq w\in\Zd$. Splitting the last sum according to whether $\vert u-y\vert\leq k/2$ or $\vert u-y'\vert\leq k/2$ (in either case the other factor is $\leq C_1(k/2)^{2-d}\lesssim k^{2-d}$ by the triangle inequality, while $\#\{u:\vert u-y\vert=m\}\lesssim m^{d-1}$ gives $\sum_{m\leq k/2}m^{d-1}\cdot C_1m^{2-d}\sim_dk^2$ for the remaining factor), and bounding the contribution from $u$ far from both $y,y'$ by a further dyadic decomposition (which converges since $d>4$), one gets $\sum_{u\in\Zd}\P_{\p_c}[u \longleftrightarrow y]\,\P_{\p_c}[u \longleftrightarrow y']\lesssim k^{4-d}$. Hence
			\begin{equation}\label{eq.PxinfC_region1}
				\sum_{\substack{u\in Q\\ \vert x-u\vert\geq n/4}}\P_{\p_c}\Ll[x\stackrel{ Q}{\longleftrightarrow}u\Rr]\,\P_{\p_c}\Ll[u \longleftrightarrow y\Rr]\,\P_{\p_c}\Ll[u \longleftrightarrow y'\Rr]\lesssim  n^{1-d}k^{4-d}.
			\end{equation}
			
			\emph{Case $\vert x-u\vert< n/4$.} Since $\vert x\vert=n$, the triangle inequality gives $\vert u\vert>3n/4$, so, using $\vert y\vert= \cc n\leq n/5$,
			\begin{equation*}
				\vert u-y\vert\geq\vert u\vert-\vert y\vert>3n/4-n/5>n/2,
			\end{equation*}
			and likewise $\vert u-y'\vert>n/2$; by the basic two-point bound above, 
			\begin{align*}
				\P_{\p_c}[u \longleftrightarrow y],\P_{\p_c}[u \longleftrightarrow y']\leq C_1(n/2)^{2-d}\lesssim n^{2-d}
			\end{align*}
			uniformly over such $u$. On the other hand, Lemma~\ref{lem.box_twopoint_onesided} gives $\P_{\p_c}[x\stackrel{ Q}{\longleftrightarrow}u]\leq C_{\ref{lem.box_twopoint_onesided}}\vert x-u\vert^{1-d}$, so
			\begin{equation*}
				\sum_{\substack{u\in Q\\ \vert x-u\vert<n/4}}\P_{\p_c}\Ll[x\stackrel{ Q}{\longleftrightarrow}u\Rr]\lesssim \sum_{m=1}^{n/4}m^{d-1}\cdot m^{1-d}\sim_dn.
			\end{equation*}
			Hence
			\begin{equation}\label{eq.PxinfC_region2}
				\sum_{\substack{u\in Q\\ \vert x-u\vert<n/4}}\P_{\p_c}\Ll[x\stackrel{ Q}{\longleftrightarrow}u\Rr]\,\P_{\p_c}\Ll[u \longleftrightarrow y\Rr]\,\P_{\p_c}\Ll[u \longleftrightarrow y'\Rr]\lesssim n\cdot(n^{2-d})^2=n^{5-2d}.
			\end{equation}
			
			Combining \eqref{eq.PxinfC_region1} and \eqref{eq.PxinfC_region2},
			\begin{equation*}
				\P_{\p_c}\Ll[x\stackrel{ Q}{\longleftrightarrow}y,\ x\stackrel{ Q}{\longleftrightarrow}y'\Rr]\lesssim  n^{1-d}k^{4-d}+n^{5-2d}.
			\end{equation*}
			Since $y,y'\in\partial H$ and $ \cc <1/5$, we have $k\leq2 \cc n<n$; as $4-d<0$, this gives $k^{4-d}\geq n^{4-d}$, so $n^{5-2d}\leq n^{1-d}k^{4-d}$, and the bound above simplifies to
			\begin{equation}\label{eq.PxinfC_pair}
				\P_{\p_c}\Ll[x\stackrel{ Q}{\longleftrightarrow}y,\ x\stackrel{ Q}{\longleftrightarrow}y'\Rr]\lesssim n^{1-d}k^{4-d}.
			\end{equation}
			
			Summing \eqref{eq.PxinfC_pair} over distinct $y,y'\in\partial H$: fixing $y$, the number of $y'\in\partial H$ at distance exactly $k-1$ from $y$ is $\lesssim k^{d-2}$, and $k$ ranges over $[1,2 \cc n]$, so
			\begin{equation*}
				\sum_{y'\in\partial H}k^{4-d}\lesssim \sum_{k=1}^{2 \cc n}k^{d-2}\cdot k^{4-d}=\sum_{k=1}^{2 \cc n}k^2\lesssim ( \cc n)^3.
			\end{equation*}
			Summing further over $y\in\partial H$ (using again $\vert\partial H\vert\sim_d( \cc n)^{d-1}$),
			\begin{equation}\label{eq.PxinfC_offdiag}
				\sum_{y, y'\in\partial H}\P_{\p_c}\Ll[x\stackrel{ Q}{\longleftrightarrow}y,\ x\stackrel{ Q}{\longleftrightarrow}y'\Rr]\lesssim n^{1-d}\cdot( \cc n)^{d-1}\cdot( \cc n)^3= \cc ^{d+2}n^3.
			\end{equation}
			This proves the lemma.
		\end{proof}

		\subsection{Expectation of $\ERE_1 (x,\fC_{\cc},r)$ and $\ERE_2 (x,\fC_{\cc},r)$}\label{subsec.regular.E}
		
		\begin{prop}\label{prop:Ef1}
			For any $ \cc  \in (0,1/5)$, there are two constants $C_{\ref{prop:Ef1}}>0$ such that for any $r \in \N$, there is $n_{\ref{prop:Ef1}} ( \cc ,r) \in \N$ such that for any $n \geq n_{\ref{prop:Ef1}}$ and any point $x \in \partial_{1}^\circ  Q$,
			\begin{equation}\label{eq:prop:Ef1}
				\E_{\p_c} \Ll[ \ERE_1 (x, \fC_{\cc},r) \Big\vert x \in \fC_{\cc} \Rr] \leq C_{\ref{prop:Ef1}} ( \cc ^{d-4} + r^{-d+6}).
			\end{equation}
		\end{prop}

		\begin{prop}\label{prop:Ef2}
			For any $ \cc  \in (0,1/5)$, there are two constants $C_{\ref{prop:Ef2}}>0$ such that for any $r \in \N$, there is $n_{\ref{prop:Ef2}} ( \cc ,r)  \in \N$ such that for any $n \geq n_{\ref{prop:Ef2}}$ and any point $x \in \partial_{1}^\circ  Q$,
			\begin{equation}\label{eq:prop:Ef2}
				\E_{\p_c} \Ll[ \ERE_2 (x, \fC_{\cc},r) \Big\vert x \in \fC_{\cc} \Rr] \leq C_{\ref{prop:Ef2}} ( \cc ^{d-4}+ r^{-d+6}).
			\end{equation}
		\end{prop}

		Recall from the Appendix~\ref{subsct:Upper bound for the 2-point function in cubes and annulus} that, for $ \cc\in(0,1/5)$, $\operatorname{int} H:= H\setminus\partial H$ and $\mathrm{Ann}_{ \cc n}:= Q\setminus\operatorname{int} H$; for $w\in\mathrm{Ann}_{ \cc n}$, $r_w$ denotes the minimal $\ell^\infty$-distance between $w$ and $\partial H$.
		
		\begin{proof}[Proof of Proposition~\ref{prop:Ef1}]
			We first compute the conditional expectation:
			\begin{equation}\label{eq.Ef1_ratio}
				\E_{\p_c}\big[\ERE_1(x,\fC_{\cc},r)\mid x\in\fC_{\cc} \big] = \frac{1}{\P_{\p_c}[x\in\fC_{\cc}]}\sum_{y\in Q\setminus S}\max\{r,\vert x-y\vert\}^{2-d}\,\P_{\p_c}\big[x,y\in\fC_{\cc}\big].
			\end{equation}
			Fix $y\in Q\setminus S$.
			Note that $\{ x \in \fC_{\cc} \} = \{ x \stackrel{\mathrm{Ann}_{ \cc n}}{\longleftrightarrow} \partial  H \}$ and $\{ y \in \fC_{\cc} \} = \{ y \stackrel{\mathrm{Ann}_{ \cc n}}{\longleftrightarrow} \partial  H \}$.
			Splitting according to whether $x$ and $y$ lie in the same cluster inside $\mathrm{Ann}_{ \cc n}$,
			\begin{equation}\label{eq.Ef1_split}
				\P_{\p_c}\big[x,y\in\fC_{\cc}\big] \leq  \P_{\p_c}\big[x\stackrel{\mathrm{Ann}_{ \cc n}}\longleftrightarrow y, x\stackrel{\mathrm{Ann}_{ \cc n}}\longleftrightarrow  H\big] + \P_{\p_c}\big[ \{x\in\fC_{\cc}\} \circ \{y\in\fC_{\cc}\}\big].
			\end{equation}
			
			For the second term in the right-side of \eqref{eq.Ef1_split}, BK inequality implies that its contribution to~\eqref{eq.Ef1_ratio} is
			\begin{align}
				&\qquad \frac{1}{\P_{\p_c}[x\in\fC_{\cc}]}\sum_{y\in Q\setminus S}\max\{r,\vert x-y\vert\}^{2-d}\, \P_{\p_c}[x\in\fC_{\cc}]\cdot\P_{\p_c}[y\in\fC_{\cc}]\\
				&=\sum_{y\in Q\setminus S}\max\{r,\vert x-y\vert\}^{2-d}\,\P_{\p_c}[y\in\fC_{\cc}].\label{niebbnr}
			\end{align}
			Since $y\in Q\setminus S$, Proposition~\ref{PxinfC} gives
			\begin{equation*}
				\P_{\p_c}[y\in\fC_{\cc}] \leq C_{\ref{PxinfC}} \,(1+l_y)\, \cc^{d-4}\,n^{-3}.
			\end{equation*}
			Substituting into~\eqref{niebbnr}, the contribution of the second term in the right-side of \eqref{eq.Ef1_split} to~\eqref{eq.Ef1_ratio} is at most
			\begin{equation*}
				C_{\ref{PxinfC}}\, \cc^{d-4}\,n^{-3}\sum_{y\in Q\setminus S}\max\{r,\vert x-y\vert\}^{2-d}\,(1+l_y).
			\end{equation*}
			Writing $m:=\vert x-y\vert$ and using $l_y \leq m$ and that the number of points of $ Q$ at distance exactly $m$ from $x$ is $\lesssim m^{d-1}$,
			\begin{equation*}
				\sum_{y\in Q\setminus S}\max\{r,\vert x-y\vert\}^{2-d}\,(1+l_y) \lesssim  \sum_{m=1}^{2n}m^d\max\{r,m\}^{2-d}.
			\end{equation*}
			Splitting at $m\lessgtr r$: for $m\leq r$, $\sum_{m\leq r}m^d\,r^{2-d}\lesssim r^3$; for $m>r$, $\sum_{r<m\leq2n}m^2\lesssim n^3$. Since $r\leq2n$, both terms are $\lesssim n^3$, so
			\begin{equation*}
				\sum_{y\in Q\setminus S}\max\{r,\vert x-y\vert\}^{2-d}\,(1+l_y) \lesssim n^3.
			\end{equation*}
			Thus, the second term in the right-side of \eqref{eq.Ef1_split} contributes, for some $C_1 >0$, at most $C_1\, \cc^{d-4}$ to~\eqref{eq.Ef1_ratio}, matching the first term on the right-side of \eqref{eq:prop:Ef1}.

			It remains to bound the contribution of the first term on the right of~\eqref{eq.Ef1_split} to~\eqref{eq.Ef1_ratio}. Fix $y\in Q\setminus S$. 
			By the BK inequality, decomposing inside $\mathrm{Ann}_{ \cc n}$ the event $\{x\leftrightarrow y\}\cap\{x \leftrightarrow\partial H\}$ at the branch point $u\in\mathrm{Ann}_{ \cc n}$, we obtain
			\begin{equation}\label{eq.Ef1_branchpoint}
				\P_{\p_c}\big[x\stackrel{\mathrm{Ann}_{ \cc n}}\longleftrightarrow y, x\stackrel{\mathrm{Ann}_{ \cc n}}\longleftrightarrow  H\big] \leq \sum_{u\in\mathrm{Ann}_{ \cc n}} \P_{\p_c}\big[x\stackrel{ Q}\longleftrightarrow u\big]\,\P_{\p_c}\big[u\stackrel{\Zd}\longleftrightarrow y\big]\,\P_{\p_c}\big[u\stackrel{\mathrm{Ann}_{ \cc n}}\leftrightarrow\partial H\big].
			\end{equation}

			Since $ \cc<1/5<1/2$, we have $\mathrm{Ann}_{ \cc n}=\big(\mathrm{Ann}_{ \cc n}\cap\Lambda_{n/2}\big)\cup\big( Q\setminus\Lambda_{n/2}\big)$, and we bound the sum in~\eqref{eq.Ef1_branchpoint} separately on each part.
			
			\medskip
			\noindent\emph{Case $u\in Q\setminus\Lambda_{n/2}$.}
			Since $\mathrm{Ann}_{ \cc n}\subset Q$, Proposition~\ref{PxinfC} gives
			\begin{equation}\label{eq.Ef1_farleg3}
				\P_{\p_c}\big[u\stackrel{\mathrm{Ann}_{ \cc n}}\longleftrightarrow\partial H\big] \leq \P_{\p_c}\big[u\stackrel{ Q}\longleftrightarrow H\big] = \P_{\p_c}[u\in\fC_{\cc}] \leq C_{\ref{PxinfC}}\,(1+l_u)\, \cc^{d-4}\,n^{-3}.
			\end{equation}
			For the first factor on the right-side of \eqref{eq.Ef1_branchpoint}, since $x\in\partial_{1}^\circ  Q$ satisfies $l_x=0$, Lemma~\ref{lem.box_twopoint_onesided} gives, with $C_{\ref{lem.box_twopoint_onesided}}>0$ as in that lemma,
			\begin{equation*}
				\P_{\p_c}\big[x\stackrel{ Q}\longleftrightarrow u\big] \leq C_{\ref{lem.box_twopoint_onesided}}\,\vert x-u\vert^{1-d}.
			\end{equation*}
			Since $1+l_u\leq1+\vert x-u\vert\lesssim \vert x-u\vert$, we have
			\begin{equation*}
				\P_{\p_c}\big[x\stackrel{ Q}\longleftrightarrow u\big]\,\P_{\p_c}\big[u\stackrel{\mathrm{Ann}_{ \cc n}}\longleftrightarrow\partial H\big] \lesssim   \cc^{d-4}\,n^{-3}\,\vert x-u\vert^{2-d}.
			\end{equation*}
			Multiplying by the middle factor on the right-side of \eqref{eq.Ef1_branchpoint} and summing over $u\in\Zd$, we recover the standard bubble bound: writing $k:=\vert x-y\vert$, splitting according to whether $\vert u-x\vert\leq k/2$ or $\vert u-y\vert\leq k/2$ (using $\P_{\p_c}[u\stackrel{\Zd}\longleftrightarrow y]\lesssim \vert u-y\vert^{2-d}$ and the triangle inequality in either case, together with $\#\{u:\vert u-x\vert=m\}\lesssim m^{d-1}$), gives $\sum_{u\in\Zd}\vert x-u\vert^{2-d}\P_{\p_c}[u\leftrightarrow y]\lesssim k^{4-d}$ (using $d>4$ for convergence of the far tail). Hence
			\begin{equation}\label{eq.Ef1_farsum}
				\sum_{u\in Q\setminus\Lambda_{n/2}} \P_{\p_c}\big[x\stackrel{ Q}\longleftrightarrow u\big]\,\P_{\p_c}\big[u\stackrel{\Zd}\longleftrightarrow y\big]\,\P_{\p_c}\big[u\stackrel{\mathrm{Ann}_{ \cc n}}\longleftrightarrow\partial H\big] \lesssim   \cc^{d-4}\,n^{-3}\,k^{4-d}.
			\end{equation}
			Multiplying by $\max\{r,k\}^{2-d}$ and summing over $y\in Q\setminus S$, using $\#\{y:\vert x-y\vert=k\}\lesssim k^{d-1}$ and splitting at $k\lessgtr r$ (using $d>6$ for convergence of the tail),
			\begin{equation*}
				\sum_y\max\{r,k\}^{2-d}\cdot  \cc^{d-4}n^{-3}k^{4-d} \lesssim   \cc^{d-4}n^{-3}\sum_kk^3\max\{r,k\}^{2-d}\lesssim   \cc^{d-4}\,n^{-3}\,r^{6-d}.
			\end{equation*}
			Dividing by the lower bound~\eqref{veiouvoh} of Proposition~\ref{PxinfC}, $\P_{\p_c}[x\in\fC_{\cc}]\geq c_{\ref{PxinfC}} \cc^{d-4}n^{-3}$, this case's contribution to~\eqref{eq.Ef1_ratio} is $\lesssim  r^{6-d}$, matching the second term on the right-side of \eqref{eq:prop:Ef1}.
			
			\medskip
			\noindent\emph{Case $u\in\mathrm{Ann}_{ \cc n}\cap\Lambda_{n/2}$.}
			We first bound the third factor in~\eqref{eq.Ef1_branchpoint}. For $u\in\mathrm{Ann}_{ \cc n}$, since $u\notin\operatorname{int} H$, the $\ell^\infty$-distance from $u$ to $ H$ equals $r_u$, so
			\begin{equation}\label{eq.Ef1_onearm}
				\P_{\p_c}\big[u\stackrel{\mathrm{Ann}_{ \cc n}}\longleftrightarrow\partial H\big] \leq \P_{\p_c}\Ll[u\stackrel{\Zd}\longleftrightarrow\Zd\setminus\big(u+\Lambda_{r_u}\big)\Rr] \lesssim (1+r_u)^{-2},
			\end{equation}
			where the last inequality is by the full-space one-arm exponent \eqref{eq.one_arm}

			Here $\vert x-u\vert\geq\vert x\vert-\vert u\vert>n-n/2=n/2$, so Lemma~\ref{lem.box_twopoint_onesided} gives $\P_{\p_c}[x\stackrel{ Q}\longleftrightarrow u]\lesssim n^{1-d}$, uniformly over such $u$. Likewise $\vert u-y\vert\geq\vert y\vert-\vert u\vert>3n/4-n/2=n/4$, so $\P_{\p_c}[u\stackrel{\Zd}\longleftrightarrow y]\lesssim n^{2-d}$, uniformly over such $u,y$. Combined with~\eqref{eq.Ef1_onearm},
			\begin{equation*}
				\text{the right-side of~\eqref{eq.Ef1_branchpoint}} \lesssim  n^{1-d}\cdot n^{2-d}\cdot(1+r_u)^{-2} = n^{3-2d}(1+r_u)^{-2},
			\end{equation*}
			uniformly in $y$. In particular this leaves the sum over $u$ and the sum over $y$ separable.
			
			For the sum over $u$, write $\rho:=\vert u\vert\in[ \cc n,n/2)$ and $\tau:=r_u=\rho- \cc n$; since $\#\{u:\vert u\vert=\rho\}\lesssim \rho^{d-1}$,
			\begin{equation*}
				\sum_{u\in\mathrm{Ann}_{ \cc n}\cap\Lambda_{n/2}}(1+r_u)^{-2} \lesssim  \sum_{\tau=0}^{n/2- \cc n}( \cc n+\tau)^{d-1}(1+\tau)^{-2}.
			\end{equation*}
			Splitting at $\tau\lessgtr  \cc n$: for $\tau\leq  \cc n$, $( \cc n+\tau)^{d-1}\lesssim ( \cc n)^{d-1}$, and $\sum_{\tau\geq0}(1+\tau)^{-2}$ converges to a universal constant, so this range contributes $\lesssim ( \cc n)^{d-1}$; for $\tau> \cc n$, $( \cc n+\tau)^{d-1}\lesssim \tau^{d-1}$, so this range contributes $\lesssim \sum_{\tau> \cc n}\tau^{d-3}\lesssim  n^{d-2}$. Hence
			\begin{equation}\label{eq.Ef1_usum}
				\sum_{u\in\mathrm{Ann}_{ \cc n}\cap\Lambda_{n/2}}(1+r_u)^{-2} \lesssim   \cc^{d-1}n^{d-1} + n^{d-2} \lesssim   \cc^{d-1}n^{d-1}
			\end{equation}
			for sufficiently large $n$.
			
			For the sum over $y$, writing $m:=\vert x-y\vert$ and using $\#\{y\in Q:\vert x-y\vert=m\}\lesssim m^{d-1}$ and splitting at $m\lessgtr r$,
			\begin{equation}\label{eq.Ef1_ysum}
				\sum_{y\in Q\setminus S}\max\{r,\vert x-y\vert\}^{2-d} \lesssim  \sum_{m=1}^{2n}m^{d-1}\max\{r,m\}^{2-d} \lesssim  r^2+n^2 \lesssim  n^2
			\end{equation}
			for sufficiently large $n$.
			
			Combining~\eqref{eq.Ef1_usum} and~\eqref{eq.Ef1_ysum}, this case contributes, to the numerator of~\eqref{eq.Ef1_ratio},
			\begin{equation*}
				\lesssim  n^{3-2d} \cdot  \cc^{d-1}n^{d-1}\cdot n^2 =  \cc^{d-1}n^{4-d}.
			\end{equation*}
			Dividing by $\P_{\p_c}[x\in\fC_{\cc}]\geq c_{\ref{PxinfC}} \cc^{d-4}n^{-3}$, this case's contribution to~\eqref{eq.Ef1_ratio} is, for some $C_3>0$,
			\begin{equation}\label{eq.Ef1_caseB}
				\leq C_3\, \cc^{3}\,n^{7-d}.
			\end{equation}
			If $d=7$, this matches the first term on the right-side of \eqref{eq:prop:Ef1}. If $d \geq 8$, this term is smaller than $ \cc^{d-4}$ for sufficiently large $n$. Therefore, the proposition holds.
		\end{proof}

		\begin{proof}[Proof of Proposition~\ref{prop:Ef2}]
			Recall that $\ERE_2(x,\fC_{\cc},r) = \sum_{y\in\fC_{\cc}\cap\partial Q}\max\{r,\vert x-y\vert\}^{4-d}$, and that $x\in\partial_{1}^\circ  Q\subset\partial Q$. We first compute the conditional expectation:
			\begin{equation}\label{eq.Ef2_ratio}
				\E_{\p_c}\big[\ERE_2(x,\fC_{\cc},r)\mid x\in\fC_{\cc}\big] = \frac{1}{\P_{\p_c}[x\in\fC_{\cc}]}\sum_{y\in\partial Q}\max\{r,\vert x-y\vert\}^{4-d}\,\P_{\p_c}\big[x,y\in\fC_{\cc}\big].
			\end{equation}
			
			As in the proof of Proposition~\ref{prop:Ef1}, $\{x\in\fC_{\cc}\}=\{x\stackrel{\mathrm{Ann}_{ \cc n}}{\longleftrightarrow}\partial H\}$ and $\{y\in\fC_{\cc}\}=\{y\stackrel{\mathrm{Ann}_{ \cc n}}{\longleftrightarrow}\partial H\}$. Splitting according to whether $x$ and $y$ lie in the same cluster inside $\mathrm{Ann}_{ \cc n}$:
			\begin{equation}\label{eq.Ef2_split}
				\P_{\p_c}\big[x,y\in\fC_{\cc}\big] \leq \P_{\p_c}\big[x\stackrel{\mathrm{Ann}_{ \cc n}}\longleftrightarrow y,\ x\stackrel{\mathrm{Ann}_{ \cc n}}\longleftrightarrow H\big] + \P_{\p_c}\big[\{x\in\fC_{\cc}\}\circ\{y\in\fC_{\cc}\}\big].
			\end{equation}

			By the BK inequality, the second term on the right of~\eqref{eq.Ef2_split} contributes to~\eqref{eq.Ef2_ratio}
			\begin{equation*}
				\sum_{y\in\partial Q\setminus\{x\}}\max\{r,\vert x-y\vert\}^{4-d}\,\P_{\p_c}[y\in\fC_{\cc}].
			\end{equation*}
			Since $y\in\partial Q$, we have $l_y=0$, so the upper bound~\eqref{eq.PxinfC_upper} of Proposition~\ref{PxinfC} gives $\P_{\p_c}[y\in\fC_{\cc}]\leq C_{\ref{PxinfC}}\, \cc^{d-4}\,n^{-3}$. Writing $k:=\vert x-y\vert$, using  $\#\{y\in\partial Q:\vert x-y\vert=k\}\lesssim k^{d-2}$, and splitting at $k\lessgtr r$: for $k\leq r$, $\sum_{k\leq r}k^{d-2}r^{4-d}\lesssim r^3$; for $k>r$, $\sum_{r<k\leq2n}k^{d-2}k^{4-d}=\sum k^2\lesssim n^3$. Thus,
			\begin{equation*}
				\sum_{y\in\partial Q\setminus\{x\}}\max\{r,\vert x-y\vert\}^{4-d} \lesssim n^3
			\end{equation*}
			for any sufficiently large $n$.
			Thus, the different-branch term contributes, for some $C_1>0$, at most $C_1 \cc^{d-4}$ to~\eqref{eq.Ef2_ratio}, matching the first term on the right-side of \eqref{eq:prop:Ef2}.
			
			It remains to bound the contribution of the first term on the right of~\eqref{eq.Ef2_split} to~\eqref{eq.Ef2_ratio}. Fix $y\in Q\setminus S$. 
			By the BK inequality, decomposing inside $\mathrm{Ann}_{ \cc n}$ the event $\{x\leftrightarrow y\}\cap\{x\leftrightarrow\partial H\}$ at the branch point $u\in\mathrm{Ann}_{ \cc n}$, and relaxing the first two factors from $\mathrm{Ann}_{ \cc n}$ to $ Q$,
			\begin{equation}\label{eq.Ef2_branchpoint}
				\P_{\p_c}\big[x\stackrel{\mathrm{Ann}_{ \cc n}}\longleftrightarrow y,\ x\stackrel{\mathrm{Ann}_{ \cc n}}\longleftrightarrow H\big] \leq \sum_{u\in\mathrm{Ann}_{ \cc n}} \P_{\p_c}\big[x\stackrel{ Q}\longleftrightarrow u\big]\,\P_{\p_c}\big[u\stackrel{ Q}\longleftrightarrow y\big]\,\P_{\p_c}\big[u\stackrel{\mathrm{Ann}_{ \cc n}}\longleftrightarrow\partial H\big].
			\end{equation}
			
			Since $ \cc<1/5<1/2$, $\mathrm{Ann}_{ \cc n}=\big(\mathrm{Ann}_{ \cc n}\cap\Lambda_{n/2}\big)\cup\big( Q\setminus\Lambda_{n/2}\big)$, and we bound the sum in~\eqref{eq.Ef2_branchpoint} separately on each part.
			
			\medskip
			\noindent\emph{Case $u\in Q\setminus\Lambda_{n/2}$.}
			As in the proof of Proposition~\ref{prop:Ef1}, $\P_{\p_c}[u\stackrel{\mathrm{Ann}_{ \cc n}}\longleftrightarrow\partial H]\leq\P_{\p_c}[u\in\fC_{\cc}]\leq C_{\ref{PxinfC}}(1+l_u) \cc^{d-4}n^{-3}$. For the first two factors, since $l_x=0$ and $l_y=0$, applying Lemma~\ref{lem.box_twopoint_twosided}, with $C_{\ref{lem.box_twopoint_twosided}}>0$ as in that lemma,
			\begin{equation*}
				\P_{\p_c}\big[x\stackrel{\mathrm{Ann}_{ \cc n}}\longleftrightarrow u\big] \leq C_{\ref{lem.box_twopoint_twosided}}(1+l_u)\vert x-u\vert^{-d}, \qquad \P_{\p_c}\big[u\stackrel{ Q}\longleftrightarrow y\big] \leq C_{\ref{lem.box_twopoint_twosided}}(1+l_u)\vert u-y\vert^{-d}.
			\end{equation*}
			Multiplying all three factors, the summand in~\eqref{eq.Ef2_branchpoint} is at most a constant multiple of
			\begin{equation*}
				(1+l_u)^3\,\vert x-u\vert^{-d}\vert u-y\vert^{-d}\, \cc^{d-4}n^{-3}.
			\end{equation*}
			Since $l_u\leq\vert x-u\vert$ and $l_u\leq\vert u-y\vert$, writing $k:=\vert x-y\vert$, we bound
			\begin{equation*}
				\Sigma:=\sum_{u\in Q\setminus\Lambda_{n/2}}(1+l_u)^3\vert x-u\vert^{-d}\vert u-y\vert^{-d}
			\end{equation*}
			by splitting according to $\vert x-u\vert\lessgtr k/2$.
			
			If $\vert x-u\vert\leq k/2$: then $\vert u-y\vert\geq k/2$ (triangle inequality), so $\vert u-y\vert^{-d}\lesssim k^{-d}$; using $l_u\leq\vert x-u\vert$, the summand is $\lesssim (1+\vert x-u\vert)^3\vert x-u\vert^{-d}k^{-d}\lesssim \vert x-u\vert^{3-d}k^{-d}$. Summing over such $u$ using $\#\{u:\vert x-u\vert=m\}\lesssim m^{d-1}$,
			\begin{equation*}
				\sum_{\vert x-u\vert\leq k/2}\vert x-u\vert^{3-d}k^{-d} \lesssim k^{-d}\sum_{m=1}^{k/2}m^{d-1}m^{3-d}=k^{-d}\sum_mm^2\lesssim k^{-d}\cdot k^3=k^{3-d}.
			\end{equation*}
			The case $\vert u-y\vert\leq k/2$ is symmetric (using $l_u\leq\vert u-y\vert$, valid since $l_y=0$) and gives the same bound $k^{3-d}$.
			
			If $\vert x-u\vert>k/2$ and $\vert u-y\vert>k/2$: bounding $(1+l_u)^3\lesssim(1+\vert x-u\vert)^3$, the summand is $\lesssim \vert x-u\vert^{3-d}\vert u-y\vert^{-d}$.
			Note that by the triangle inequality, we have $|u-y| \geq \max \{ k/2, |x-u|-k \} \gtrsim |x-u|$. Thus, the summand is $\lesssim  \rho^{3-2d}$, where $\rho:= |x-u| >k/2$.
			Using $\#\{u:\vert x-u\vert=\rho\}\lesssim \rho^{d-1}$, this case's contribution to the sum $\Sigma$ is
			\begin{equation*}
				\lesssim  \sum_{\rho=k/2}^{+\infty} \rho^{d-1} \cdot \rho^{3-2d} =\sum_{\rho=k/2}^{+\infty} \rho^{2-d} \lesssim  k^{3-d}.
			\end{equation*}
			
			Therefore,
			\begin{equation}\label{eq.Ef2_farsum}
				\Sigma \lesssim k^{3-d}.
			\end{equation}
			Multiplying by $ \cc^{d-4}n^{-3}$ and by $\max\{r,k\}^{4-d}$, and summing over $y\in\partial Q$ using $\#\{y:\vert x-y\vert=k\}\lesssim k^{d-2}$,
			\begin{equation*}
				\sum_yk^{d-2}\cdot  \cc^{d-4}n^{-3}k^{3-d}\max\{r,k\}^{4-d}= \cc^{d-4}n^{-3}\sum_kk\max\{r,k\}^{4-d}.
			\end{equation*}
			Splitting at $k\lessgtr r$: for $k\leq r$, $\sum_{k\leq r}k\cdot r^{4-d}\sim r^{6-d}$; for $k>r$, $\sum_{k>r}k^{5-d}\lesssim r^{6-d}$ (converges since $d>6$). Hence this case contributes $\lesssim  \cc^{d-4}n^{-3}r^{6-d}$ to the numerator of~\eqref{eq.Ef2_ratio}. Dividing by the lower bound~\eqref{veiouvoh} of Proposition~\ref{PxinfC}, this case's contribution to~\eqref{eq.Ef2_ratio} is
			\begin{equation}\label{eq.Ef2_caseA}
				\lesssim r^{6-d},
			\end{equation}
			matching the second term on the right-side of \eqref{eq:prop:Ef2}.
			
			\medskip
			\noindent\emph{Case $u\in\mathrm{Ann}_{ \cc n}\cap\Lambda_{n/2}$.}
			As in the proof of Proposition~\ref{prop:Ef1}, the third factor in~\eqref{eq.Ef2_branchpoint} satisfies $\P_{\p_c}[u\stackrel{\mathrm{Ann}_{ \cc n}}\longleftrightarrow\partial H]\lesssim (1+r_u)^{-2}$ by \eqref{eq.one_arm}. Here $\vert x-u\vert\geq\vert x\vert-\vert u\vert>n-n/2=n/2$ and $\vert u-y\vert\geq\vert y\vert-\vert u\vert>n-n/2=n/2$ (using $\vert x\vert=\vert y\vert=n$); so Lemma~\ref{lem.box_twopoint_twosided} gives $\P_{\p_c}[x\stackrel{ Q}\longleftrightarrow u]\lesssim  n^{1-d}$ and $\P_{\p_c}[u\stackrel{ Q}\longleftrightarrow y]\lesssim n^{1-d}$, uniformly over such $u,y$. Combined with the one-arm bound,
			\begin{equation*}
				\text{the right-side of \eqref{eq.Ef2_branchpoint}} \lesssim n^{1-d}\cdot n^{1-d}\cdot(1+r_u)^{-2}=n^{2-2d}(1+r_u)^{-2},
			\end{equation*}
			uniformly in $y$, leaving the sum over $u$ and the sum over $y$ separable.
			
			Recall from the proof of Proposition~\ref{prop:Ef1} that $\sum_{u\in\mathrm{Ann}_{ \cc n}\cap\Lambda_{n/2}}(1+r_u)^{-2}\lesssim  \cc^{d-1}n^{d-1}$ for sufficiently large $n$. For the sum over $y$, writing $k=\vert x-y\vert$ and using $\#\{y\in\partial Q:\vert x-y\vert=k\}\lesssim k^{d-2}$ and splitting at $k\lessgtr r$,
			\begin{equation*}
				\sum_{y\in\partial Q}\max\{r,\vert x-y\vert\}^{4-d} \lesssim \sum_{k=1}^{2n}k^{d-2}\max\{r,k\}^{4-d} \lesssim r^3+n^3\lesssim n^3
			\end{equation*}
			for sufficiently large $n$.
			
			Combining these two sums, this case contributes, to the numerator of~\eqref{eq.Ef2_ratio},
			\begin{equation*}
				\lesssim n^{2-2d}\cdot  \cc^{d-1}n^{d-1}\cdot n^3 =  \cc^{d-1}n^{4-d}.
			\end{equation*}
			Dividing by $\P_{\p_c}[x\in\fC_{\cc}]\geq c_{\ref{PxinfC}} \cc^{d-4}n^{-3}$, this case's contribution to~\eqref{eq.Ef2_ratio} is, for some $C_2>0$,
			\begin{equation}\label{eq.Ef2_caseB}
				\leq C_2\, \cc^3\,n^{7-d}.
			\end{equation}
			If $d=7$, this matches the first term on the right-side of \eqref{eq:prop:Ef2}. If $d\geq8$, this term is smaller than $ \cc^{d-4}$ for sufficiently large $n$. Therefore, the proposition holds.
		\end{proof}

		\subsection{Variance of $\ERE_1 (x,\fC_{\cc},r)$ and $\ERE_2 (x,\fC_{\cc},r)$}\label{subsec.regular.var}
		\begin{prop}\label{prop:VARf1}
			For any $ \cc \in (0,1/5)$, there is $C_{\ref{prop:VARf1}} ( \cc) >0$ such that for any $r \in \N$, there is $n_{\ref{prop:VARf1}} ( \cc,r) \in \N$ such that for any $n \geq n_{\ref{prop:VARf1}}$ and any point $x \in \partial_{1}^\circ  Q$,
			\begin{equation}
				\var_{\p_c} \Ll[ \ERE_1 (x, \fC_{\cc},r) \Big\vert x \in \fC_{\cc} \Rr] \leq C_{\ref{prop:VARf1}} r^{-d+6}.
			\end{equation}
		\end{prop}

		\begin{prop}\label{prop:VARf2}
			For any $ \cc \in (0,1/5)$, there is $C_{\ref{prop:VARf2}} ( \cc) >0$ such that for any $r \in \N$, there is $n_{\ref{prop:VARf2}} ( \cc,r) \in \N$ such that for any $n \geq n_{\ref{prop:VARf2}}$ and any point $x \in \partial_{1}^\circ  Q$,
			\begin{equation}
				\var_{\p_c} \Ll[ \ERE_2 (x, \fC_{\cc},r) \Big\vert x \in \fC_{\cc} \Rr] \leq C_{\ref{prop:VARf2}} r^{-d+6}.
			\end{equation}
		\end{prop}

		We first establish two key lemmas to prove the above two propositions.
		In this subsection, we fix $x \in \partial_{1}^\circ  Q$ and define the event
		\begin{equation}\label{biejeobsr}
			A(y_1,y_2) := \{ x,y_1, y_2 \in \fC_{\cc} \} \setminus \big\{ \{x,y_1 \in \fC_{\cc}\} \circ \{y_2 \in \fC_{\cc}\} \big\}
		\end{equation}
		for $y_1, y_2 \in  Q$.
		
		\begin{lemma}\label{nisnmfs}
			For any function $f:  Q \rightarrow [0, +\infty)$ and any $x \in  Q$,
			\begin{equation}
				\var\Ll[ \sum_{y \in  Q} f(y) \bigg\vert x \in \fC_{\cc} \Rr] \leq \frac{\sum_{y_1,y_2 \in  Q} f(y_1) f(y_2) \P_{\p_c} [A(y_1,y_2)]}{\P_{\p_c} [x \in \fC_{\cc}]}
			\end{equation}
		\end{lemma}
		
		\begin{proof}
			Let $Z:= \sum_{y \in  Q} f(y)$.
			Let us calculate the conditional variance:
			\begin{align}
				&\qquad \var\Ll[ Z | x \in \fC_{\cc} \Rr]\label{dkjsghkt}\\
				&= \E_{\p_c} \Ll[ Z^2 | x \in \fC_{\cc} \Rr] -\E_{\p_c} \Ll[ Z | x \in \fC_{\cc} \Rr]^2\\
				&=\frac{\P_{\p_c} [x \in \fC_{\cc}] \E_{\p_c} \Ll[ Z^2 \cdot \mathbf{1} (x \in \fC_{\cc}) \Rr]- \E_{\p_c} \Ll[ Z \cdot \mathbf{1} (x \in \fC_{\cc}) \Rr]^2}{\P_{\p_c} [x \in \fC_{\cc}]^2}.\label{suheboi}
			\end{align}
			Further calculate:
			\begin{equation}
				\qquad \E_{\p_c} \Ll[ Z^2 \cdot \mathbf{1} (x \in \fC_{\cc}) \Rr] =\sum_{y_1,y_2 \in  Q} f(y_1) f(y_2) \P_{\p_c} [x,y_1,y_2 \in \fC_{\cc}];
			\end{equation}
			\begin{equation}
				\begin{aligned}
					&\qquad \E_{\p_c} \Ll[ Z \cdot \mathbf{1} (x \in \fC_{\cc}) \Rr]^2 = \Big( \sum_{y \in  Q} f(y) \P [x,y \in \fC_{\cc}] \Big)^2\\
					&= \sum_{y_1,y_2 \in  Q} f(y_1) f(y_2) \P_{\p_c} [x,y_1 \in \fC_{\cc}] \P_{\p_c} [x,y_2 \in \fC_{\cc}].
				\end{aligned}
			\end{equation}
			Thus, the numerator of \eqref{suheboi} equals
			\begin{align}
				\sum_{y_1,y_2 \in  Q} &f(y_1) f(y_2) \cdot\\
				&\big( \P_{\p_c} [x \in \fC_{\cc}] \P_{\p_c} [x,y_1,y_2 \in \fC_{\cc}] -\P_{\p_c} [x,y_1 \in \fC_{\cc}] \P_{\p_c} [x,y_2 \in \fC_{\cc}] \big).\label{eothopb}
			\end{align}
			
			Recalling the definition of $A(y_1,y_2)$ from \eqref{biejeobsr} and applying BK inequality, we obtain
			\begin{align}
				\P_{\p_c} [x,y_1,y_2 \in \fC_{\cc}] &= \P_{\p_c} \big[ \{x,y_1 \in \fC_{\cc}\} \circ \{y_2 \in \fC_{\cc}\} \big] +\P_{\p_c} [A(y_1,y_2)]\label{pnoowp}\\
				&\leq \P_{\p_c} [x,y_1 \in \fC_{\cc}] \P_{\p_c} [y_2 \in \fC_{\cc}]+ \P_{\p_c} [A(y_1,y_2)].\label{pnoowp2}
			\end{align}
			Using \eqref{pnoowp}--\eqref{pnoowp2} and the FKG inequality, we obtain
			\begin{align*}
				\eqref{eothopb} &\leq \P_{\p_c} [x \in \fC_{\cc}] \P_{\p_c} [x,y_1 \in \fC_{\cc}] \P_{\p_c} [y_2 \in \fC_{\cc}]+ \P_{\p_c} [x \in \fC_{\cc}] \P_{\p_c} [A(y_1,y_2)]\\
				&\qquad \qquad \qquad -\P_{\p_c} [x,y_1 \in \fC_{\cc}] \P_{\p_c} [x,y_2 \in \fC_{\cc}]\\
				&\leq \P_{\p_c} [x \in \fC_{\cc}] \P_{\p_c} [A(y_1,y_2)].
			\end{align*}
			
			Hence, the numerator of \eqref{suheboi} equals
			\begin{equation}\label{vekboqqeo}
				\P_{\p_c} [x \in \fC_{\cc}] \sum_{y_1,y_2 \in  Q} f(y_1) f(y_2) \cdot \P_{\p_c} [A(y_1,y_2)].
			\end{equation}
			Substituting \eqref{vekboqqeo} into \eqref{dkjsghkt}--\eqref{suheboi} proves the lemma.
		\end{proof}

		
		For $D\subset\Zd$ and $u,v\in D$, write $\tau_D(u,v):=\P_{\p_c}\big[u\stackrel{D}{\longleftrightarrow}v\big]$, and for $u\in\mathrm{Ann}_{ \cc n}$ let
		\begin{align*}
			g(u) := \P_{\p_c}\Ll[u\stackrel{\mathrm{Ann}_{ \cc n}}{\longleftrightarrow}\partial H\Rr].
		\end{align*}
		
		For $y_1,y_2\in\mathrm{Ann}_{ \cc n}$, define
		\begin{align*}
			\EE_{\ref{lem.Ay1y2_diagrams}.1}(y_1,y_2) &:= \sum_{u,w\in\mathrm{Ann}_{ \cc n}} g(w)\,\tau_{\mathrm{Ann}_{ \cc n}}(w,u)\,\tau_{\mathrm{Ann}_{ \cc n}}(w,y_1)\,\tau_{\mathrm{Ann}_{ \cc n}}(u,x)\,\tau_{\mathrm{Ann}_{ \cc n}}(u,y_2),\\
			\EE_{\ref{lem.Ay1y2_diagrams}.2}(y_1,y_2) &:= \sum_{u,w\in\mathrm{Ann}_{ \cc n}} g(u)\,\tau_{\mathrm{Ann}_{ \cc n}}(u,w)\,\tau_{\mathrm{Ann}_{ \cc n}}(u,y_2)\,\tau_{\mathrm{Ann}_{ \cc n}}(w,x)\,\tau_{\mathrm{Ann}_{ \cc n}}(w,y_1),\\
			\EE_{\ref{lem.Ay1y2_diagrams}.3}(y_1,y_2) &:= \sum_{u,w\in\mathrm{Ann}_{ \cc n}} g(u)\,\tau_{\mathrm{Ann}_{ \cc n}}(u,w)\,\tau_{\mathrm{Ann}_{ \cc n}}(u,x)\,\tau_{\mathrm{Ann}_{ \cc n}}(w,y_2)\,\tau_{\mathrm{Ann}_{ \cc n}}(w,y_1),\\
			\EE_{\ref{lem.Ay1y2_diagrams}.4}(y_1,y_2) &:= g(y_1)\sum_{u\in\mathrm{Ann}_{ \cc n}} g(u)\,\tau_{\mathrm{Ann}_{ \cc n}}(u,x)\,\tau_{\mathrm{Ann}_{ \cc n}}(u,y_2),\\
			\EE_{\ref{lem.Ay1y2_diagrams}.5}(y_1,y_2) &:= g(x)\sum_{w\in\mathrm{Ann}_{ \cc n}} g(w)\,\tau_{\mathrm{Ann}_{ \cc n}}(w,y_1)\,\tau_{\mathrm{Ann}_{ \cc n}}(w,y_2).
		\end{align*}
		
		\begin{figure}[h!]
			\centering
			\includegraphics[width=\textwidth]{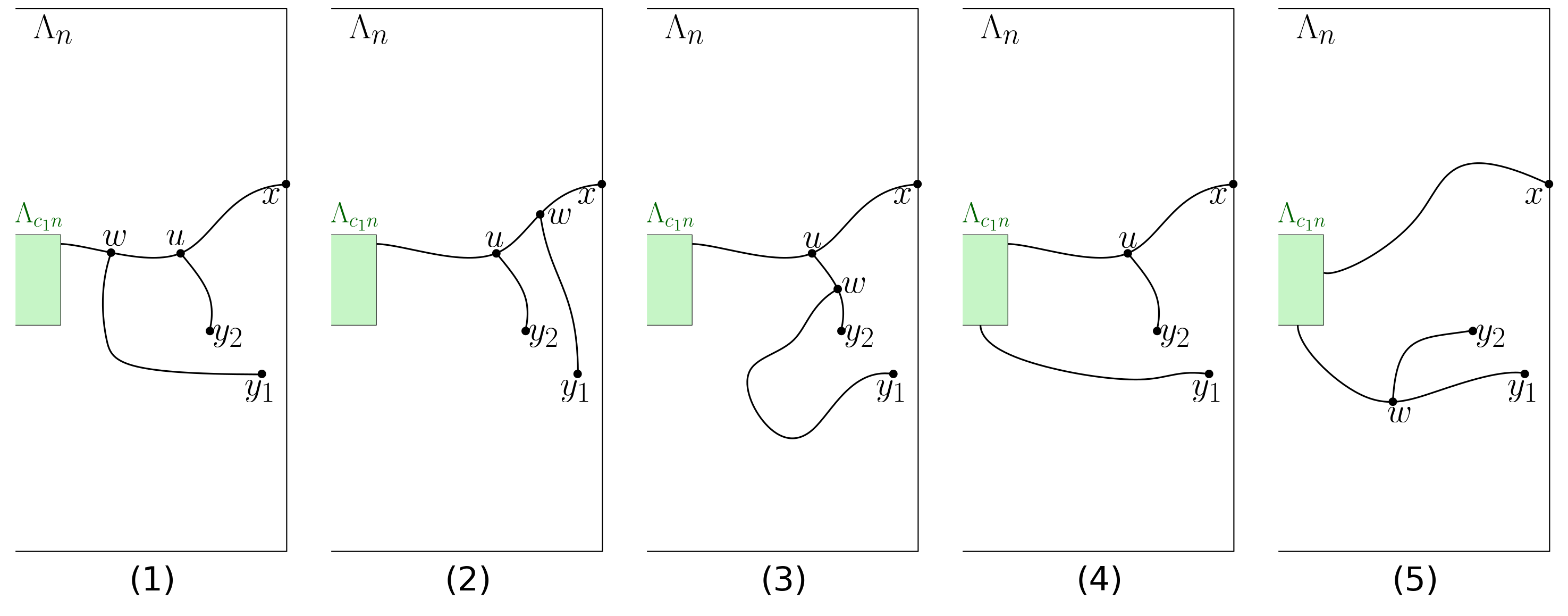}
			\caption{The diagrams (1)--(5), corresponding to $\EE_{\ref{lem.Ay1y2_diagrams}.1},\ldots,\EE_{\ref{lem.Ay1y2_diagrams}.5}$, respectively. The green box is $ H$, and each curve is an open path in $\mathrm{Ann}_{ \cc n}$ (a curve ending on the green box ends on $\partial H$).}\label{fig.Ay1y2_diagrams}
		\end{figure}
		
		\begin{lemma}\label{lem.Ay1y2_diagrams}
			For every $ \cc\in(0,1/5)$, every $n\in\N_+$, and every $y_1,y_2\in\mathrm{Ann}_{ \cc n}$,
			\begin{align*}
				\P_{\p_c}\big[A(y_1,y_2)\big] \leq \sum_{i=1}^{5}\EE_{\ref{lem.Ay1y2_diagrams}.i}(y_1,y_2),
			\end{align*}
			where the five terms correspond, in this order, to the five diagrams in Figure~\ref{fig.Ay1y2_diagrams}.
		\end{lemma}
		
		\begin{proof}
			Throughout this proof, write $E:=\{x,y_1\in\fC_{\cc}\}$ and $R:=\partial H$. Recall $\{v\in\fC_{\cc}\}=\{v\stackrel{\mathrm{Ann}_{ \cc n}}{\longleftrightarrow}R\}$ for $v\in\mathrm{Ann}_{ \cc n}$. Fix a configuration $\omega\in A(y_1,y_2)$, and open paths $\pi_x,\pi_{y_1},\pi_{y_2}$ in $\mathrm{Ann}_{ \cc n}$ from $x,y_1,y_2$ to $R$, respectively. If $\pi_{y_2}$ intersected neither $\pi_x$ nor $\pi_{y_1}$, then $\pi_x\cup\pi_{y_1}$ would witness $E$ and be edge-disjoint from $\pi_{y_2}$, which witnesses $\{y_2\in\fC_{\cc}\}$; thus $E\circ\{y_2\in\fC_{\cc}\}$ would occur, contradicting $\omega\in A(y_1,y_2)$. Hence $\pi_{y_2}$ intersects $\pi_x$ or $\pi_{y_1}$.
			
			\emph{Case 1: $\pi_x$ and $\pi_{y_2}$ intersect.} Let $u$ be the first vertex of $\pi_{y_2}$ (starting from $y_2$) lying on $\pi_x$, and let $T:=\pi_x\cup\pi_{y_2}[y_2,u]$. Then $T$ is the union of three edge-disjoint open paths, from $u$ to $x$, to $y_2$, and to $R$, which meet only at $u$. If $\pi_{y_1}$ does not intersect $T$, then $\omega$ contains these three paths together with $\pi_{y_1}$, which are pairwise edge-disjoint; this is Figure~\ref{fig.Ay1y2_diagrams}(4). Otherwise, let $w$ be the first vertex of $\pi_{y_1}$ (starting from $y_1$) lying on $T$; then $\pi_{y_1}[y_1,w]$ is edge-disjoint from $T$, and $w$ lies on one of the three legs of $T$ (possibly at $u$). If $w$ is on the leg from $u$ to $R$, we get Figure~\ref{fig.Ay1y2_diagrams}(1); if $w$ is on the leg from $u$ to $x$, we get Figure~\ref{fig.Ay1y2_diagrams}(2); if $w$ is on the leg from $u$ to $y_2$, we get Figure~\ref{fig.Ay1y2_diagrams}(3).
			
			\emph{Case 2: $\pi_x$ and $\pi_{y_2}$ do not intersect.} Then $\pi_{y_1}$ intersects $\pi_{y_2}$. Let $w$ be the first vertex of $\pi_{y_1}$ (starting from $y_1$) lying on $\pi_{y_2}$, and let $T':=\pi_{y_2}\cup\pi_{y_1}[y_1,w]$, the union of three edge-disjoint open paths from $w$ to $y_1$, to $y_2$, and to $R$, which meet only at $w$. We claim that $\pi_{y_1}[y_1,w]$ does not intersect $\pi_x$. Otherwise, let $v$ be the first vertex of $\pi_{y_1}$ lying on $\pi_x$; then $v$ comes before $w$ on $\pi_{y_1}$ (as $w\in\pi_{y_2}$ and $\pi_x\cap\pi_{y_2}=\emptyset$), so $\pi_{y_1}[y_1,v]$ does not intersect $\pi_{y_2}$. Therefore $\pi':=\pi_{y_1}[y_1,v]\cup\pi_x[v,R]$ is an open path from $y_1$ to $R$ that does not intersect $\pi_{y_2}$, and $\pi_x\cup\pi'$ witnesses $E$ and is edge-disjoint from $\pi_{y_2}$. This contradicts $\omega\in A(y_1,y_2)$. Hence $\pi_x$ does not intersect $T'$, and $\omega$ contains $\pi_x$ together with the three legs of $T'$, which are pairwise edge-disjoint. This is Figure~\ref{fig.Ay1y2_diagrams}(5).
			
			In each case, $\omega$ contains pairwise edge-disjoint open paths in $\mathrm{Ann}_{ \cc n}$ forming one of the five diagrams. Applying the BK inequality, the probability that Figure~\ref{fig.Ay1y2_diagrams}($i$) occurs is at most $\EE_{\ref{lem.Ay1y2_diagrams}.i}(y_1,y_2)$. Summing over $i=1,\ldots,5$ proves the lemma.
		\end{proof}

		\begin{proof}[Proof of Proposition~\ref{prop:VARf1}]
			Fix $x \in \partial_{1}^\circ  Q$. Applying Lemma~\ref{nisnmfs} with
			\begin{equation*}
				f(y):= \max\{ r, |x-y| \}^{-d+2} \cdot \mathbf{1} (y \notin  S),
			\end{equation*}
			we get
			\begin{equation}
				\begin{aligned}
					\var_{\p_c} &\Ll[ \ERE_1 (x, \fC_{\cc},r) \Big\vert x \in \fC_{\cc} \Rr]\\
					&\leq \frac{\sum_{y_1,y_2 \in  Q \setminus  S} \max \{r, |x-y_1|\}^{-d+2} \cdot \max \{r, |x-y_2|\}^{-d+2} \cdot \P_{\p_c} [A(y_1,y_2)]}{\P_{\p_c} [x \in \fC_{\cc}]}.
				\end{aligned}
			\end{equation}
			
			Then applying Lemma~\ref{lem.Ay1y2_diagrams} and noting that $\EE_{\ref{lem.Ay1y2_diagrams}.2} (y_1,y_2) =\EE_{\ref{lem.Ay1y2_diagrams}.1} (y_2,y_1)$, we obtain
			\begin{equation}\label{brouhl}
				\begin{aligned}
					\var_{\p_c} &\Ll[ \ERE_1 (x, \fC_{\cc},r) \Big\vert x \in \fC_{\cc} \Rr]\\
					&\leq \frac{1}{\P_{\p_c} [x \in \fC_{\cc}]} \sum_{y_1,y_2 \in  Q \setminus  S} \max \{r, |x-y_1|\}^{-d+2} \cdot \max \{r, |x-y_2|\}^{-d+2} \cdot \\
					&\qquad \qquad \qquad \big( 2\EE_{\ref{lem.Ay1y2_diagrams}.1}(y_1,y_2) +\EE_{\ref{lem.Ay1y2_diagrams}.3}(y_1,y_2) +\EE_{\ref{lem.Ay1y2_diagrams}.4}(y_1,y_2) +\EE_{\ref{lem.Ay1y2_diagrams}.5}(y_1,y_2) \big).
				\end{aligned}
			\end{equation}
			
			Recall $g(u) := \P_{\p_c} \big[ u\stackrel{\mathrm{Ann}_{ \cc n}}{\longleftrightarrow}\partial H \big]$.
			Define
			\begin{equation*}
				\rho_1 := \sum_{y \in  Q \setminus  S} \max \{ r,|x-y| \}^{-d+2} g(y)
			\end{equation*}
			and define
			\begin{equation*}
				F_1 (u) := \sum_{y \in  Q \setminus  S} \max \{ r,|x-y| \}^{-d+2} \tau_{\mathrm{Ann}_{ \cc n}} (u,y)
			\end{equation*}
			for $u \in \mathrm{Ann}_{ \cc n}$.
			Expanding $\EE_{\ref{lem.Ay1y2_diagrams}.i} (y_1,y_2)$ in \eqref{brouhl} and summing over $y_1,y_2$, we get
			\begin{align*}
				&\sum_{y_1,y_2 \in  Q \setminus  S} \max \{r, |x-y_1|\}^{-d+2} \cdot \max \{r, |x-y_2|\}^{-d+2} \cdot \EE_{\ref{lem.Ay1y2_diagrams}.1}(y_1,y_2) \\
				&\qquad \leq \sum_{u,w \in \mathrm{Ann}_{ \cc n}} g(w) \tau_{\mathrm{Ann}_{ \cc n}} (w,u) \tau_{\mathrm{Ann}_{ \cc n}} (u,x) F_1 (w) F_1 (u),\\
				&\sum_{y_1,y_2 \in  Q \setminus  S} \max \{r, |x-y_1|\}^{-d+2} \cdot \max \{r, |x-y_2|\}^{-d+2} \cdot \EE_{\ref{lem.Ay1y2_diagrams}.3}(y_1,y_2) \\
				&\qquad \leq \sum_{u,w \in \mathrm{Ann}_{ \cc n}} g(u) \tau_{\mathrm{Ann}_{ \cc n}} (u,w) \tau_{\mathrm{Ann}_{ \cc n}} (u,x) F_1 (w)^2,\\
				&\sum_{y_1,y_2 \in  Q \setminus  S} \max \{r, |x-y_1|\}^{-d+2} \cdot \max \{r, |x-y_2|\}^{-d+2} \cdot \EE_{\ref{lem.Ay1y2_diagrams}.4}(y_1,y_2) \\
				&\qquad \leq \rho_1 \sum_{u \in \mathrm{Ann}_{ \cc n}} g(u) \tau_{\mathrm{Ann}_{ \cc n}} (u,x) F_1 (u),\\
				&\sum_{y_1,y_2 \in  Q \setminus  S} \max \{r, |x-y_1|\}^{-d+2} \cdot \max \{r, |x-y_2|\}^{-d+2} \cdot \EE_{\ref{lem.Ay1y2_diagrams}.5}(y_1,y_2) \\
				&\qquad \leq \P_{\p_c} [x \in \fC_{\cc}] \sum_{w \in \mathrm{Ann}_{ \cc n}} g(w) F_1 (w)^2,
			\end{align*}
			and, thus, by Corollary~\ref{cor.PxinfC},
			\begin{align}
				\var_{\p_c} &\Ll[ \ERE_1 (x, \fC_{\cc},r) \Big\vert x \in \fC_{\cc} \Rr] \lesssim_{ \cc} \\
				&n^3 \sum_{u,w \in \mathrm{Ann}_{ \cc n}} g(w) \tau_{\mathrm{Ann}_{ \cc n}} (w,u) \tau_{\mathrm{Ann}_{ \cc n}} (u,x) F_1 (w) F_1 (u)\label{sunzi1} \\
				+&n^3 \sum_{u,w \in \mathrm{Ann}_{ \cc n}} g(u) \tau_{\mathrm{Ann}_{ \cc n}} (u,w) \tau_{\mathrm{Ann}_{ \cc n}} (u,x) F_1 (w)^2\\
				+&n^3 \rho_1 \sum_{u \in \mathrm{Ann}_{ \cc n}} g(u) \tau_{\mathrm{Ann}_{ \cc n}} (u,x) F_1 (u)\label{vskfnob}\\
				+& \sum_{w \in \mathrm{Ann}_{ \cc n}} g(w) F_1 (w)^2.\label{sunzi4}
			\end{align}
			
			In the remainder of this proof, we only bound the part of \eqref{vskfnob} with $u\in Q\setminus\Lambda_{n/2}$, which gives the dominant contribution among the terms above. The remaining terms are bounded in the Appendix.
			
			
			\medskip
			\noindent\emph{Estimate of $\rho_1$.} Since $\{y\in\fC_{\cc}\}=\{y\stackrel{\mathrm{Ann}_{ \cc n}}{\longleftrightarrow}\partial H\}$, we have $g(y)=\P_{\p_c}[y\in\fC_{\cc}]$, so Proposition~\ref{PxinfC} gives, for every $y\in Q\setminus S$,
			\begin{equation*}
				g(y)\leq C_{\ref{PxinfC}}\,(1+l_y)\, \cc^{d-4}\,n^{-3}.
			\end{equation*}
			Since $\vert x\vert=n$, we have $l_y \leq\vert x-y\vert$. Using $\#\{y\in Q:\vert x-y\vert=m\}\lesssim m^{d-1}$, and splitting at $m\lessgtr r$,
			\begin{equation*}
				\rho_1\lesssim_{ \cc} n^{-3} \sum_{m=1}^{2n}m \cdot m^{d-1}\max\{r,m\}^{2-d} \leq n^{-3}\big(1+r^3+n^3\big)\lesssim_{ \cc} 1.
			\end{equation*}
			
			\medskip
			\noindent\emph{Estimate of $F_1$.} We claim that, for every $u\in \mathrm{Ann}_{ \cc n}$,
			\begin{equation}\label{eq.F1_bound}
				F_1(u)\lesssim\max\{r,\vert x-u\vert\}^{4-d}.
			\end{equation}
			Write $k:=\vert x-u\vert$ and $M:=\max\{r,k\}$. Since $\tau_{\mathrm{Ann}_{ \cc n}} (u,y)\leq \tau_{\Zd}(u,y)\lesssim(1+\vert u-y\vert)^{2-d}$,
			\begin{equation*}
				F_1(u)\lesssim \sum_{y\in\Zd}\max\{r,\vert x-y\vert\}^{2-d}\,(1+\vert u-y\vert)^{2-d}.
			\end{equation*}
			We split the sum into three parts.
			\begin{itemize}
				\item $\vert y-u\vert\leq M/2$. Here $\max\{r,\vert x-y\vert\}^{2-d}\lesssim M^{2-d}$: if $k\leq r$ this holds because $\max\{r,\vert x-y\vert\}\geq r=M$, and if $k>r$ because $\vert x-y\vert\geq k-M/2=M/2$. Since $\sum_{\vert y-u\vert\leq M/2}(1+\vert u-y\vert)^{2-d}\lesssim M^2$, this part contributes $\lesssim M^{4-d}$.
				\item $\vert y-u\vert>M/2$ and $\vert x-y\vert\geq\vert y-u\vert/2$. The summand is $\lesssim \vert y-u\vert^{4-2d}$, so this part contributes $\lesssim \sum_{m>M/2}m^{d-1}m^{4-2d}\lesssim M^{4-d}$ (using $d>4$).
				\item $\vert y-u\vert>M/2$ and $\vert x-y\vert<\vert y-u\vert/2$. This part is empty if $k=0$. Otherwise, the triangle inequality gives $2k/3<\vert y-u\vert<2k$ and $\vert x-y\vert<k$, so $(1+\vert u-y\vert)^{2-d}\lesssim k^{2-d}$, and
				\begin{equation*}
					\sum_{\vert x-y\vert<k}\max\{r,\vert x-y\vert\}^{2-d}\lesssim \begin{cases}k^dr^{2-d}&\text{if }k\leq r,\\k^2&\text{if }k>r.\end{cases}
				\end{equation*}
				Hence this part contributes $\lesssim k^2r^{2-d}\leq r^{4-d}$ if $k\leq r$, and $\lesssim k^{4-d}$ if $k>r$; in both cases $\lesssim M^{4-d}$.
			\end{itemize}
			Summing the three parts proves \eqref{eq.F1_bound}.
			
			\medskip
			\noindent\emph{Estimate of \eqref{vskfnob} for $u\in Q\setminus\Lambda_{n/2}$.} For such $u$, Proposition~\ref{PxinfC} gives $g(u)=\P_{\p_c}[u\in\fC_{\cc}]\lesssim_{ \cc} (1+l_u) n^{-3}$.  Lemma~\ref{lem.box_twopoint_onesided} gives $\tau_{\mathrm{Ann}_{ \cc n}}(u,x)\lesssim (1+\vert x-u\vert)^{1-d}$. Together with $l_u\leq\vert x-u\vert$, this yields
			\begin{equation*}
				g(u)\,\tau_{\mathrm{Ann}_{ \cc n}}(u,x)\lesssim_{ \cc} n^{-3}\,(1+\vert x-u\vert)^{2-d}.
			\end{equation*}
			Combining this with the bound $\rho_1\lesssim_{ \cc} 1$ and with \eqref{eq.F1_bound}, and using $\#\{u\in Q:\vert x-u\vert=k\}\lesssim (1+k)^{d-1}$,
			\begin{align*}
				n^3\rho_1\sum_{u\in Q\setminus\Lambda_{n/2}}g(u)\,\tau_{\mathrm{Ann}_{ \cc n}}(u,x)\,F_1(u)&\lesssim_{ \cc} \sum_{u\in Q}(1+\vert x-u\vert)^{2-d}\max\{r,\vert x-u\vert\}^{4-d}\\
				&\lesssim \sum_{k=0}^{2n}(1+k)\max\{r,k\}^{4-d}.
			\end{align*}
			For $k\leq r$ the summand is $\lesssim(1+k)\,r^{4-d}$, so these $k$ contribute $\lesssim r^{6-d}$; for $k>r$ the summand is $\lesssim k^{5-d}$, and $\sum_{k>r}k^{5-d}\lesssim r^{6-d}$ since $d>6$. Therefore,
			\begin{equation*}
				n^3\rho_1\sum_{u\in Q\setminus\Lambda_{n/2}}g(u)\,\tau_{\mathrm{Ann}_{ \cc n}}(u,x)\,F_1(u)\lesssim_{ \cc} r^{6-d},
			\end{equation*}
			which is as claimed in the proposition.
		\end{proof}

		\begin{proof}[Proof of Proposition~\ref{prop:VARf2}]
			Fix $x \in \partial_{1}^\circ  Q$.
			
			Applying Lemma~\ref{nisnmfs} with
			\begin{equation*}
				f(y):= \max\{ r, |x-y| \}^{-d+4} \cdot \mathbf{1} (y \in \partial  Q),
			\end{equation*}
			we get
			\begin{equation}
				\begin{aligned}
					\var_{\p_c} &\Ll[ \ERE_2 (x, \fC_{\cc},r) \Big\vert x \in \fC_{\cc} \Rr]\\
					&\leq \frac{\sum_{y_1,y_2 \in \partial  Q} \max \{r, |x-y_1|\}^{-d+4} \cdot \max \{r, |x-y_2|\}^{-d+4} \cdot \P_{\p_c} [A(y_1,y_2)]}{\P_{\p_c} [x \in \fC_{\cc}]}.
				\end{aligned}
			\end{equation}
			
			Then applying Lemma~\ref{lem.Ay1y2_diagrams} and noting that $\EE_{\ref{lem.Ay1y2_diagrams}.2} (y_1,y_2) =\EE_{\ref{lem.Ay1y2_diagrams}.1} (y_2,y_1)$, we obtain
			\begin{equation}
				\begin{aligned}
					\var_{\p_c} &\Ll[ \ERE_2 (x, \fC_{\cc},r) \Big\vert x \in \fC_{\cc} \Rr]\\
					&\leq \frac{1}{\P_{\p_c} [x \in \fC_{\cc}]} \sum_{y_1,y_2 \in \partial  Q} \max \{r, |x-y_1|\}^{-d+4} \cdot \max \{r, |x-y_2|\}^{-d+4} \cdot \\
					&\qquad \qquad \qquad \big( 2\EE_{\ref{lem.Ay1y2_diagrams}.1}(y_1,y_2) +\EE_{\ref{lem.Ay1y2_diagrams}.3}(y_1,y_2) +\EE_{\ref{lem.Ay1y2_diagrams}.4}(y_1,y_2) +\EE_{\ref{lem.Ay1y2_diagrams}.5}(y_1,y_2) \big).
				\end{aligned}
			\end{equation}
			
			Recall $g(u) := \P_{\p_c} \big[ u\stackrel{\mathrm{Ann}_{ \cc n}}{\longleftrightarrow}\partial H \big]$.
			Define
			\begin{equation*}
				\rho_2 := \sum_{y \in \partial  Q} \max \{ r,|x-y| \}^{-d+4} g(y)
			\end{equation*}
			and define
			\begin{equation*}
				F_2 (u) := \sum_{y \in \partial  Q} \max \{ r,|x-y| \}^{-d+4} \tau_{\mathrm{Ann}_{ \cc n}} (u,y)
			\end{equation*}
			for $u \in \mathrm{Ann}_{ \cc n}$.
			Expanding $\EE_{\ref{lem.Ay1y2_diagrams}.i} (y_1,y_2)$ and summing over $y_1,y_2$ exactly as in the proof of Proposition~\ref{prop:VARf1} (with $\rho_1$ and $F_1$ replaced by $\rho_2$ and $F_2$), and then using Corollary~\ref{cor.PxinfC}, we get
			\begin{align}
				\var_{\p_c} &\Ll[ \ERE_2 (x, \fC_{\cc},r) \Big\vert x \in \fC_{\cc} \Rr] \lesssim_{ \cc} \\
				&n^3 \sum_{u,w \in \mathrm{Ann}_{ \cc n}} g(w) \tau_{\mathrm{Ann}_{ \cc n}} (w,u) \tau_{\mathrm{Ann}_{ \cc n}} (u,x) F_2 (w) F_2 (u)\label{sunzi5} \\
				+&n^3 \sum_{u,w \in \mathrm{Ann}_{ \cc n}} g(u) \tau_{\mathrm{Ann}_{ \cc n}} (u,w) \tau_{\mathrm{Ann}_{ \cc n}} (u,x) F_2 (w)^2\\
				+&n^3 \rho_2 \sum_{u \in \mathrm{Ann}_{ \cc n}} g(u) \tau_{\mathrm{Ann}_{ \cc n}} (u,x) F_2 (u)\label{vskfnob2}\\
				+& \sum_{w \in \mathrm{Ann}_{ \cc n}} g(w) F_2 (w)^2.\label{sunzi8}
			\end{align}
			
			In the remainder of this proof, we only bound the part of \eqref{vskfnob2} with $u\in Q\setminus\Lambda_{n/2}$, which gives the dominant contribution among the terms above. The remaining terms are bounded in the Appendix.
			
			\medskip
			\noindent\emph{Estimate of $\rho_2$.} For $y\in\partial Q$ we have $l_y=0$, so Proposition~\ref{PxinfC} gives $g(y)=\P_{\p_c}[y\in\fC_{\cc}]\lesssim_{ \cc}n^{-3}$. Since $\#\{y\in\partial Q:\vert x-y\vert = m\}\lesssim m^{d-2}$, we get
			\begin{equation*}
				\rho_2\lesssim_{ \cc}n^{-3}\sum_{y\in\partial Q}\max\{r,\vert x-y\vert\}^{4-d}\lesssim n^{-3}\big(1+r^3+n^3\big)\lesssim_{ \cc}1.
			\end{equation*}
			
			\medskip
			\noindent\emph{Estimate of $F_2$.} We claim that, for every $u\in Q$,
			\begin{equation}\label{eq.F2_bound}
				F_2(u)\lesssim\max\{r,\vert x-u\vert\}^{4-d}.
			\end{equation}
			Write $k:=\vert x-u\vert$ and $M:=\max\{r,k\}$. Since $x\in\partial Q$, we have $l_u \leq k$, so $1+l_u\lesssim M$. For $y\in\partial Q$, Lemma~\ref{lem.box_twopoint_twosided} gives $\tau_{\mathrm{Ann}_{ \cc n}}(u,y)\leq\tau_{ Q}(u,y)\lesssim(1+l_u)(1+\vert u-y\vert)^{-d}$, and therefore
			\begin{equation*}
				F_2(u)\lesssim(1+l_u)\sum_{y\in\partial Q}\max\{r,\vert x-y\vert\}^{4-d}\,(1+\vert u-y\vert)^{-d}.
			\end{equation*}
			We will use that $\vert u-y\vert\geq l_u$ for every $y\in\partial Q$, and that $\#\{y\in\partial Q:\vert u-y\vert=m\}\lesssim(1+m)^{d-2}$.
			We split the sum into three parts.
			\begin{itemize}
				\item $\vert y-u\vert\leq M/2$. Here $\max\{r,\vert x-y\vert\}^{4-d}\lesssim M^{4-d}$: if $k\leq r$ this holds because $\max\{r,\vert x-y\vert\}\geq r=M$, and if $k>r$ because $\vert x-y\vert\geq k-M/2=M/2$. Moreover,
				\begin{equation*}
					\sum_{y\in\partial Q}(1+\vert u-y\vert)^{-d}\lesssim\sum_{m\geq l_u}(1+m)^{d-2}(1+m)^{-d}\lesssim(1+l_u)^{-1}.
				\end{equation*}
				Hence this part contributes $\lesssim M^{4-d}$ to $F_2(u)$.
				\item $\vert y-u\vert>M/2$ and $\vert x-y\vert\geq\vert y-u\vert/2$. Here $\max\{r,\vert x-y\vert\}^{4-d}\leq\vert x-y\vert^{4-d}\lesssim\vert y-u\vert^{4-d}$, so the summand is $\lesssim(1+\vert u-y\vert)^{4-2d}$, and
				\begin{equation*}
					\sum_{\substack{y\in\partial Q\\ \vert y-u\vert>M/2}}(1+\vert u-y\vert)^{4-2d}\lesssim\sum_{m>M/2}(1+m)^{d-2}(1+m)^{4-2d}\lesssim M^{3-d}.
				\end{equation*}
				Together with the factor $1+l_u\lesssim M$, this part contributes $\lesssim M^{4-d}$ to $F_2(u)$.
				\item $\vert y-u\vert>M/2$ and $\vert x-y\vert<\vert y-u\vert/2$. The triangle inequality gives $2k/3<\vert y-u\vert<2k$ and $\vert x-y\vert<k$, so $(1+\vert u-y\vert)^{-d}\lesssim k^{-d}$, and $1+l_u\lesssim k$. If $k\leq r$, then $\max\{r,\vert x-y\vert\}^{4-d}\leq r^{4-d}$ and $\#\{y\in\partial Q:\vert x-y\vert<k\}\lesssim k^{d-1}$, so this part contributes $\lesssim(1+l_u)\,k^{-d}\,k^{d-1}\,r^{4-d}\lesssim r^{4-d}$. If $k>r$, using $\#\{y\in\partial Q:\vert x-y\vert=m\}\lesssim m^{d-2}$ and splitting at $m\lessgtr r$,
				\begin{equation*}
					\sum_{\substack{y\in\partial Q\\ \vert x-y\vert<k}}\max\{r,\vert x-y\vert\}^{4-d}\lesssim \sum_{m=1}^{r}m^{d-2}r^{4-d}+\sum_{m=r+1}^{k}m^{d-2}m^{4-d}\lesssim k^3.
				\end{equation*}
				Hence this part contributes $\lesssim(1+l_u)\,k^{-d}\,k^3\lesssim k^{4-d}$. In both cases this part contributes $\lesssim M^{4-d}$ to $F_2(u)$.
			\end{itemize}
			Summing the three parts proves \eqref{eq.F2_bound}.
			
			\medskip
			\noindent\emph{Estimate of \eqref{vskfnob2} for $u\in Q\setminus\Lambda_{n/2}$.} By the same argument as in the proof of Proposition~\ref{prop:VARf1}, using Proposition~\ref{PxinfC}, Lemma~\ref{lem.box_twopoint_onesided} and $l_u\leq\vert x-u\vert$, we have
			\begin{equation*}
				g(u)\,\tau_{\mathrm{Ann}_{ \cc n}}(u,x)\lesssim_{ \cc}n^{-3}\,(1+\vert x-u\vert)^{2-d}.
			\end{equation*}
			Combining this with $\rho_2\lesssim_{ \cc}1$ and with \eqref{eq.F2_bound}, and using $\#\{u\in Q:\vert x-u\vert=k\}\lesssim(1+k)^{d-1}$,
			\begin{align*}
				n^3\rho_2\sum_{u\in Q\setminus\Lambda_{n/2}}g(u)\,\tau_{\mathrm{Ann}_{ \cc n}}(u,x)\,F_2(u)&\lesssim_{ \cc}\sum_{u\in Q}(1+\vert x-u\vert)^{2-d}\max\{r,\vert x-u\vert\}^{4-d}\\
				&\lesssim\sum_{k=0}^{2n}(1+k)\max\{r,k\}^{4-d}\lesssim r^{6-d}.
			\end{align*}
			This is as claimed in the proposition.
		\end{proof}

		\subsection{Regular clusters and regular cubes}\label{subsec.regular.clt}
		The mission in this part is to prove that $Q$ has high probability to have many regular boundary points. We start from the expectation of irregular boundary points.
		\begin{prop}\label{prop.choice_c}
			Given $\delta, \epsilon >0$, there is $\cc(d, \delta) \in (0,1/5)$ and $r_{\ref{prop.choice_c}} (d, \delta, \epsilon) >0$ satisfying
			\begin{equation}\label{voebliu2e}
				\E_{\p_c} [\#\{x \in \fC_{\cc} \cap \partial_{1}^\circ  Q: x \text{ is not } (r_{\ref{prop.choice_c}}, \delta)\text{-regular}\}] \leq \epsilon n^{d-4}.
			\end{equation}
		\end{prop}
		
		\begin{proof}
			The notation $\lesssim$ in the proof exceptionally depends on $\delta$. We first choose a small number $\cc >0$ and a large integer $r_* \in \N$ such that
			\begin{equation}\label{eq.constrain.cr}
				\begin{split}
					\max \{ C_{\ref{prop:Ef1}}, C_{\ref{prop:Ef2}} \} \cdot \cc^{d-4} &= \delta /4, \\
					\max \{ C_{\ref{prop:Ef1}}, C_{\ref{prop:Ef2}} \} \cdot r_*^{6-d} &= \delta /4.
				\end{split}
			\end{equation}
			Then by Propositions~\ref{prop:Ef1} and~\ref{prop:Ef2}, for any $x \in \partial_{1}^\circ  Q$ and any $r \geq r_*$,
			\begin{align*}
				\E_{\p_c} \Ll[ \ERE_1 (x, \fC_{\cc}, r) \Big\vert x \in \fC_{\cc} \Rr] &\leq C_{\ref{prop:Ef1}} (\cc^{d-4} +r^{6-d}) <\frac{\delta}{2},\\
				\E_{\p_c} \Ll[ \ERE_2 (x, \fC_{\cc}, r) \Big\vert x \in \fC_{\cc} \Rr] &\leq C_{\ref{prop:Ef2}} (\cc^{d-4} +r^{6-d}) <\frac{\delta}{2}.
			\end{align*}
			Then applying Propositions~\ref{prop:VARf1} gives
			\begin{align*}
				&\P_{\p_c} \Ll[ \ERE_1 (x, \fC_{\cc}, r) >\delta \Big\vert x \in \fC_{\cc} \Rr]\\
				\leq\ &\P_{\p_c} \Ll[ \ERE_1 (x, \fC_{\cc}, r) -\E_{\p_c} \Ll[ \ERE_1 (x, \fC_{\cc}, r) \big\vert x \in \fC_{\cc} \Rr] >\frac{\delta}{2} \Big\vert x \in \fC_{\cc} \Rr]\\
				\leq\ &\var_{\p_c} \Ll[ \ERE_1 (x, \fC_{\cc}, r) \Big\vert x \in \fC_{\cc} \Rr] \cdot \Big( \frac{\delta}{2} \Big)^{-2}\\
				\lesssim\ & \delta^{-2} r^{-d+6}.
			\end{align*}
			Similarly, Proposition~\ref{prop:VARf2} implies
			\begin{equation*}
				\P_{\p_c} \Ll[ \ERE_1 (x, \fC_{\cc}, r) >\delta \Big\vert x \in \fC_{\cc} \Rr] \lesssim \delta^{-2} r^{-d+6}.
			\end{equation*}
			Then
			\begin{equation*}
				\begin{aligned}
					&\P_{\p_c} \Ll[ x \text{ is not } (r,\delta) \text{-regular} \Big\vert x \in \fC_{\cc} \Rr]\\
					\leq\ &\P_{\p_c} \Ll[ \ERE_1 (x, \fC_{\cc}, r) >\delta \Big\vert x \in \fC_{\cc} \Rr] +\P_{\p_c} \Ll[ \ERE_2 (x, \fC_{\cc}, r) >\delta \Big\vert x \in \fC_{\cc} \Rr]\\
					\lesssim\ & \delta^{-2} r^{-d+6}.
				\end{aligned}
			\end{equation*}
			The estimate above and Corollary~\ref{cor.PxinfC} together yield that
			\begin{equation}\label{eq.irr.one}
				\begin{aligned}
					&\P_{\p_c} [x \in \fC_{\cc} \text{ and } x \text{ is not } (r,\delta) \text{-regular}]\\
					=\ &\P_{\p_c} [x \in \fC_{\cc}] \P_{\p_c} \Ll[ x \text{ is not } (r,\delta) \text{-regular} \Big\vert x \in \fC_{\cc} \Rr]\\
					\lesssim\ & \delta^{-2} \cc^{d-4} r^{-d+6} n^{-3},
				\end{aligned}
			\end{equation}
			We set $r = \theta r_*$ with $\theta \geq 1$, and sum \eqref{eq.irr.one} over $x \in \partial_{1}^\circ  Q$ and using $|\partial_{1}^\circ  Q| \lesssim n^{d-1}$
			\begin{align*}
				&\E_{\p_c} [\#\{x \in \fC_{\cc} \cap \partial_{1}^\circ  Q: x \text{ is not } (r_{\ref{prop.choice_c}}, \delta)\text{-regular}\}] \\
				&\lesssim \delta^{-2} \cc^{d-4} r^{-d+6} n^{d-4} \\
				&\lesssim \theta^{-d+6} n^{d-4}.
			\end{align*}
			In the last line, we use \eqref{eq.constrain.cr} for the choice of $\cc$ and $r_*$. Thus, for each $\epsilon >0$, we can find sufficiently large $\theta$ such that \eqref{voebliu2e} holds.
		\end{proof}
		\begin{remark}
			From here to the end of the main text, we always fix $\cc(d,\delta)$ associated to $d$ and $\delta$ in the $(r,\delta)$-regular point, so that Proposition~\ref{prop.choice_c} holds.
		\end{remark}

		Now we define the regular cluster and regular cube. For every  $K \subseteq \fC_{\cc}$, we write
		\begin{equation}\label{eq.def.regular_bd}
			\PTRE{(r,\delta)}_{i} K := \{ y \in K \cap \partial_{i}^\circ Q: y \text{ is } (r,\delta) \text{-regular} \}.
		\end{equation}
		The set $\PTRE{(r,\delta)}_{i} K$ depends not only on $K$, but also on the full configuration of $\fC_{\cc}$.
		
		\begin{definition}[Regular cluster and regular cube]
			A cluster $K \subseteq \fC_{\cc}$ is called \emph{$(r,\delta)$-regular} if for every $i \in\{\pm 1, \pm 2, \cdots, \pm d\}$, we have
			\begin{align}\label{eq.def.reg.K}
				|\PTRE{(r,\delta)}_{i} K| \geq \delta n^{2}.
			\end{align}
			An admissible $\fC_{\cc}$ is called \emph{$(r,\delta)$-regular} if 
			\begin{align}\label{eq.def.reg.C}
				\#\{K \subset \fC_{\cc}: K \text{ is } (r,\delta)\text{-regular}\} \geq \delta n^{d-6}.
			\end{align}
		\end{definition}

		The following proposition states that a cube has high probability to be regular. Recall $\delta_{\mathrm{ext}}(d)$ defined in Proposition~\ref{prop.extension}.
		
		\begin{prop}\label{prop.regular_cube}
			For any $\epsilon >0$, there is $r_{\ref{prop.regular_cube}} (\epsilon) \in \N$ and $\delta_{\ref{prop.regular_cube}} (\epsilon) \in (0, \delta_{\mathrm{ext}}(d))$ such that for every sufficiently large $n \in \N$,
			\begin{equation}
				\P_{\p_c} [\fC_{\cc} \text{ is } (r_{\ref{prop.regular_cube}},\delta_{\ref{prop.regular_cube}}) \text{-regular}] \geq 1-\epsilon.
			\end{equation}
		\end{prop}
		\begin{proof}
			We define two events 
			\begin{align*}
				\mathsf{Reg} (r,\delta) := \{\fC_{\cc} \text{ is } (r,\delta) \text{-regular}\}, 
			\end{align*}
			and 
			\begin{multline*}
				A(\delta):= \{\text{there are } \delta n^{d-6} \text{ spanning clusters } K \subset \fC_{\cc} \\
				\text{ s.t. } |K \cap \partial_{i}^\circ   Q| \geq \delta_2 n^{2}, \forall i = \pm 1, \pm 2, \cdots, \pm d \}.
			\end{multline*}
			By Lemma~\ref{lem.arm.2dspanning} and Proposition~\ref{prop.num_spanning_new}, we can find $\delta_2 (\epsilon) >0$ such that for sufficiently large $n$,
			\begin{equation}\label{iboipjk2}
				\P_{\p_c} [A (\delta_2)] >1- (\epsilon /2).
			\end{equation}

			We set the parameter
			\begin{align*}
				\delta_{\ref{prop.regular_cube}} (\epsilon):= \min\{\delta_2 /2, \delta_{\mathrm{ext}}(d)\},
			\end{align*}
			and are going to estimate
			\begin{equation*}
				\P_{\p_c} [A (\delta_2) \setminus \mathsf{Reg} (r,\delta_{\ref{prop.regular_cube}})].
			\end{equation*}
			
			If $\mathsf{Reg} (r,\delta_{\ref{prop.regular_cube}})$ does not hold, there will be at most $\delta_{\ref{prop.regular_cube}} n^{d-6}$ many $(r, \delta_{\ref{prop.regular_cube}})$-regular clusters inside $ Q$. If, moreover, $A (\delta_2)$ holds, there will be at least $\delta_{\ref{prop.regular_cube}} n^{d-6}$ clusters $K$ inside $ Q$ with the following properties:
			\begin{itemize}
				\item $|K \cap \partial_{i}^\circ   Q| \geq \delta_2 n^{2} \geq 2 \delta_{\ref{prop.regular_cube}} n^{2}$ for all $i \in\{\pm 1, \pm 2, \cdots, \pm d\}$;
				\item not $(r, \delta_{\ref{prop.regular_cube}})$-regular, i.e., there is some $i \in\{\pm 1, \pm 2, \cdots, \pm d\}$ such that
				\begin{align*}
						|\PTRE{(r, \delta_{\ref{prop.regular_cube}})}_{i} K| < \delta_{\ref{prop.regular_cube}} n^{2}.
				\end{align*}
			\end{itemize}
			Therefore, we find $\delta_{\ref{prop.regular_cube}}^2 n^{d-4}$ boundary points in $\fC_{\cc} \cap \partial^\circ  Q$ that are not $(r, \delta_{\ref{prop.regular_cube}})$-regular, which gives that
			\begin{equation}\label{iewbufv}
				\begin{aligned}
					&\P_{\p_c} [A (\delta_2) \setminus \mathsf{Reg} (r,\delta_2/2)]\\
					\leq\ & \P_{\p_c} \Ll[\#\{x \in \fC_{\cc} \cap \partial^\circ  Q: x \text{ is not } (r, \delta_{\ref{prop.regular_cube}})\text{-regular}\} \geq \delta_{\ref{prop.regular_cube}}^2 n^{d-4}\Rr]\\
					\leq\ & \sum_{1 \leq |i| \leq d} \P_{\p_c} \Ll[\#\{x \in \fC_{\cc} \cap \partial^\circ_i  Q: x \text{ is not } (r, \delta_{\ref{prop.regular_cube}})\text{-regular}\} \geq \frac{\delta_{\ref{prop.regular_cube}}^2 n^{d-4}}{2 d}\Rr]\\
					=\ &2d\ \P_{\p_c} \Ll[\#\{x \in \fC_{\cc} \cap \partial^\circ_1  Q: x \text{ is not } (r, \delta_{\ref{prop.regular_cube}})\text{-regular}\} \geq \frac{\delta_{\ref{prop.regular_cube}}^2 n^{d-4}}{2 d}\Rr],
				\end{aligned}
			\end{equation}
			where the last equality is by symmetry.
			
			The parameter $r_{\ref{prop.regular_cube}} (\epsilon)$ is set as
			\begin{align*}
				r_{\ref{prop.regular_cube}} (\epsilon) := r_{\ref{prop.choice_c}} \Ll(d, \delta_{\ref{prop.regular_cube}}, \frac{\delta_{\ref{prop.regular_cube}}^2 \cdot \epsilon}{8 d^2}\Rr).
			\end{align*}
			Then Proposition~\ref{prop.choice_c} and Markov inequality imply that 
			\begin{equation*}
				\P_{\p_c} \Ll[\#\{x \in \fC_{\cc} \cap \partial^\circ_1  Q: x \text{ is not } (r_{\ref{prop.regular_cube}}, \delta_{\ref{prop.regular_cube}})\text{-regular}\} \geq \frac{\delta_{\ref{prop.regular_cube}}^2 n^{d-4}}{2 d}\Rr] \leq \frac{\epsilon}{4d}.
			\end{equation*}
			Applying this into \eqref{iewbufv} and recalling $\delta_{\ref{prop.regular_cube}} =\delta_2 /2$, we get
			\begin{equation}
				\P_{\p_c} [A (\delta_2) \setminus \mathsf{Reg} (r_{\ref{prop.regular_cube}},\delta_{\ref{prop.regular_cube}})] \leq \frac{\epsilon}{2}.
			\end{equation}
			Recalling \eqref{iboipjk2}, we obtain
			\begin{equation}
				\begin{aligned}
					\P_{\p_c} [\mathsf{Reg} (r_{\ref{prop.regular_cube}}, \delta_{\ref{prop.regular_cube}})] &\geq \P_{\p_c} [A (\delta_2) \cap \mathsf{Reg} (r_{\ref{prop.regular_cube}}, \delta_{\ref{prop.regular_cube}})]\\
					& =\P_{\p_c} [A (\delta_2)] -\P_{\p_c} [A (\delta_2) \setminus \mathsf{Reg} (r_{\ref{prop.regular_cube}}, \delta_{\ref{prop.regular_cube}})]\\
					& \geq 1- \frac{\epsilon}{2} -\frac{\epsilon}{2} =1-\epsilon.
				\end{aligned}
			\end{equation}
			This proves the proposition.
		\end{proof}

		\section{Lower bound III: The sprinkling property}\label{sec.skl}
		This section is devoted to the sprinkling phenomena: under the monotone coupling of the percolation,  when the parameter increases from $\p_c$ to $\p$, a spanning cluster of $\p_c$-percolation in correlation length $\ell_{*}$ will have large probability to explore another distance $\ell_{*}$. Then this explains how the percolation evolves from critical phase to supercritical phase.
		
\begin{figure}[h!]
			\centering
			\includegraphics[width=0.8\textwidth]{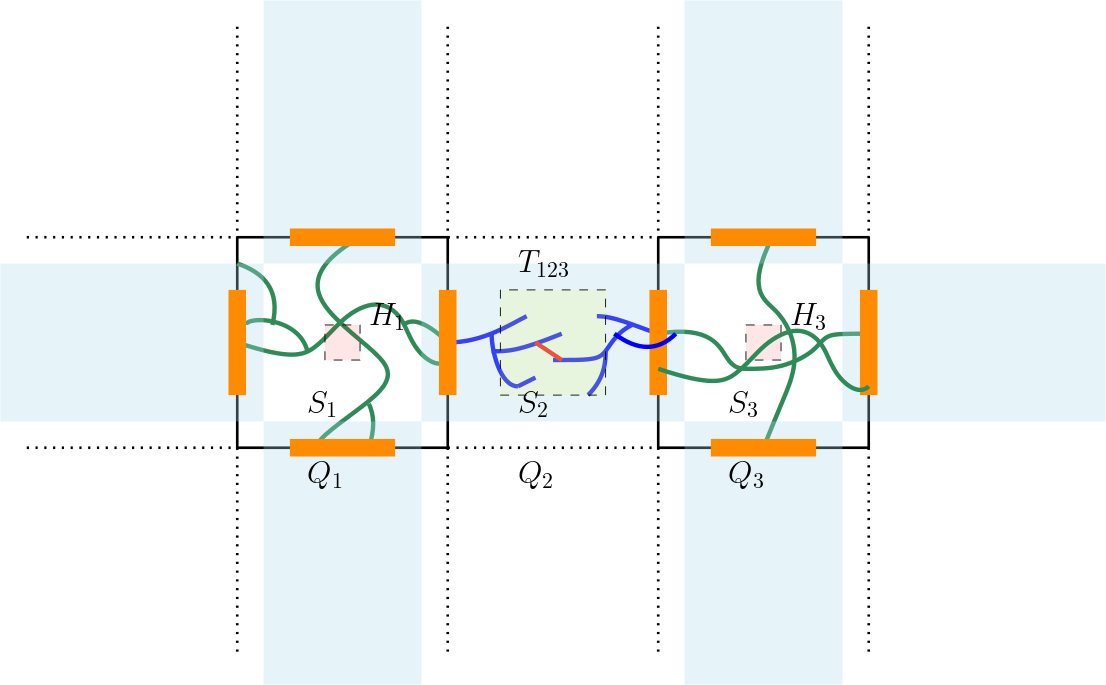}
			\caption{An illustration of sprinkling clusters between the neighbor cubes.}\label{fig.skl}
		\end{figure}

		To make the heuristic above more rigorous, we focus on a small neighborhood consisting of $3$ cubes. For convenience, we fix $n \in \N_+$ and follow the $(Q,S,H,T)$-structure in \eqref{eq.def.SQH_regular} 
		\begin{align}\label{eq.def.SQH}
			Q_1 := \Lambda_n, \qquad S_1:= \Lambda_{3n/4}, \qquad H_1:= \Lambda_{\cc n}.
		\end{align}
		This structure is shifted to the neighbor cubes centered at $n e_1$, with $S_2$ slightly smaller
		\begin{align}\label{eq.def.SQH2}
			Q_2 := \Lambda_n(n e_1), \qquad S_2:= \Lambda_{n/2}(n e_1),
		\end{align}
		and those centered at $2n e_1$
		\begin{align}\label{eq.def.SQH3}
			Q_3 := \Lambda_n(2n e_1), \qquad S_3 := \Lambda_{3n/4}(2n e_1), \qquad H_3 := \Lambda_{\cc n}(2n e_1).
		\end{align}
		See Figure~\ref{fig.skl} for an illustration.

		Throughout the section, all the clusters refer to those in $\p_c$-percolation without specific explanation. We assume the information of the clusters containing $H_1$ and living in $Q_1$ are already known, and denote their configuration by $\fC_{1}$. The counterpart in $Q_3$ is defined similarly.
		\begin{align}\label{eq.def.fC}
			\fC_{1} := \{y \in Q_1: y \xleftrightarrow{Q_1} H_1\}, \qquad \fC_{3} := \{y \in Q_3: y \xleftrightarrow{Q_3} H_3\}. 
		\end{align}
		The following shorthand notations will also be used 
		\begin{equation}\label{eq.fC_shorthand}
			\fC_{13} \equiv \fC_{1} \cup \fC_{3}. 
		\end{equation}
		We will extend them from the boundary $\partial_{1}^{\circ} Q_1, \partial_{-1}^\circ Q_3$, which follows \eqref{eq.def.boundary.o}. The sprinkling clusters will live in the tube $T_{12d}$ defined as 
		\begin{align}\label{eq.def.tube}
			T_{123} := n e_1 + \Big( -\frac{5n}{8}, \frac{5n}{8} \Big) \times \Big( -\frac{3n}{8}, \frac{3n}{8} \Big)^{d-1},
		\end{align}
		and its restricted segment in different cubes are defined as 
		\begin{align}\label{eq.def.tube.12_23}
			\begin{split}
				T_{12} &:= n e_1 + \Big( -\frac{5n}{8}, \frac{3n}{8} \Big) \times \Big( -\frac{3n}{8}, \frac{3n}{8} \Big)^{d-1}, \\
				T_{23} &:= n e_1 + \Big( -\frac{3n}{8}, \frac{5n}{8} \Big) \times \Big( -\frac{3n}{8}, \frac{3n}{8} \Big)^{d-1}.
			\end{split}
		\end{align}

		One essential object in sprinkling is once again the pivotal bond. For technical reason, its definition slightly differs from that in previous sections, and is defined as follows. We let $\cC_{12}(x), \cC_{23}(x)$ respectively stand for the cluster containing $x$ that lives in $T_{12}$ or $T_{23}$
		\begin{align}\label{eq.def.C12_23}
			\cC_{12}(x) := \{y \in T_{12}: y \xleftrightarrow{T_{12} \setminus \mathcal{E}(\fC_{1})} x\}, \qquad \cC_{23}(x) := \{y \in T_{23}: y \xleftrightarrow{T_{23} \setminus \mathcal{E}(\fC_{3})} x\}.
		\end{align}
		Here $\xleftrightarrow{T_{12} \setminus \mathcal{E}(\fC_{1})}$ means the connection of open path in $\p_c$-percolation without using the bonds in $\fC_{1}$ (but the intersection in vertex is still allowed).	For every $A \subset \partial_{1}^\circ Q_1 \cap \fC_{1}$ and $B \subset \partial_{-1}^\circ Q_3 \cap \fC_{3}$, the pivotal edges between them are then defined as
		\begin{multline}\label{eq.def.piv_AB}
			\piv(A,B) := \Ll\{(e_{-}, e_+) \in \overrightarrow{E}(S_2): e_{-} \xleftrightarrow{T_{12} \setminus \mathcal{E}(\fC_{1})} A, \quad e_+ \xleftrightarrow{T_{23} \setminus \mathcal{E}(\fC_{3})} B, \Rr.\\
			\Ll. \text{ and } \cC_{12}(e_{-}) \cap \cC_{23}(e_{+}) = \emptyset \Rr\}.
		\end{multline}
		Especially, this definition implies that $\a(\{e_{-}, e_+\}) = 0$. $\piv(A,B)$ means that, in $\p_c$-percolation, one explores $A$ and $B$ respectively in $T_{12}$ and $T_{23}$, then $\{e_{-}, e_+\}$ is a closed edge to make $A$ connecting $B$ fail. Afterwards, we also let $\pivo(A,B)$ stand for the open pivotal edge in $\p$-percolation
		\begin{align}\label{eq.def.pivo_AB}
			\pivo(A,B) := \Ll\{(e_{-}, e_+) \in \piv(A,B): \{e_{-}, e_+\} \text{ is open in } \p\text{-percolation}\Rr\}.
		\end{align}
		\begin{definition}[Sprinkling cluster and paths]\label{def.skl}
			The \emph{sprinkling cluster} connecting $A$ and $B$ is defined as 
			\begin{align}\label{eq.def.skl_AB}
				\cC_{s}(A,B) := \bigcup_{(e_{-}, e_+) \in \pivo(A,B)} \cC_{12}(e_{-}) \cup \cC_{23}(e_+). 
			\end{align}
			More precisely, \eqref{eq.def.skl_AB} stands for the set of vertices, and $\cC_{s}(A,B)$ if formed by
			\begin{itemize}
				\item  the two clusters $\cC_{12}(e_{-}), \cC_{23}(e_+)$ of $\p_c$-percolation,
				\item  and one open bond $\{e_{-}, e_{+}\}$ in $\p$-percolation. 
			\end{itemize}
			In $\cC_{s}(A,B)$, every open path in  connecting $A$ and $B$, and containing a unique $(e_{-}, e_+) \in \pivo(A,B)$ as pivotal edge is called \emph{a sprinkling path}.
		\end{definition}
		
		Thus $\cC_{s}(A,B)$ connects $A$ and $B$ in a very specific way. By default, we assume $\cC_{s}(A,B) \neq \emptyset$ when using this notation. $\piv(A,B), \pivo(A,B)$ and  $\cC_{s}(A,B)$ are the objects we will study in this section. 
		
		Let us state one major proposition in this section. All the regular points in this section refers to $(r,\delta)$-regular as in Proposition~\ref{prop.extension}, so there is an implicit dependence in function of $r$ and $\delta$ there. As this dependence will not be used in this section, we just omit it to simplify the notation.

		\begin{prop}\label{prop.sprinkling.AB}
			Given $\epsilon > 0$, there exists a finite positive constants $c_{\text{\ref{eq.sprinkling.AB.condition}}}(\epsilon)$,  such that for all large integer $n$ and all regular sets $A \subset \partial_{1}^\circ Q_1 \cap \fC_{1}$ and $B \subset \partial_{-1}^\circ Q_3 \cap \fC_{3}$ satisfying
			\begin{align}\label{eq.sprinkling.AB.condition}
				\vert A \vert, \vert B \vert \geq c_{\text{\ref{eq.sprinkling.AB.condition}}} n^2, \qquad  \vert A \vert \cdot \vert B \vert \geq c_{\text{\ref{eq.sprinkling.AB.condition}}} n^{d-2} , 
			\end{align}
			we have 
			\begin{equation}\label{eq.sprinkling.AB}
				\begin{split}
					M := \E_{\p_c} \Ll[\# \piv(A,B) \, \Big\vert \fC_{13} \Rr] &\asymp n^{-d+4} |A| \cdot |B|, \\
					\P_{\p_c} \Ll[ \frac{M}{2} \leq \# \piv(A,B) \leq \frac{3 M}{2}  \, \Big\vert \fC_{13} \Rr] &\geq 1-\epsilon.
				\end{split}
			\end{equation}
		\end{prop}
		In another word, the typical value of $ \# \piv(A,B) $ is $n^{-d+4} |A| \cdot |B|$.
		Viewing Proposition~\ref{prop.regular_cube}, the spanning cluster has typical regular boundary of order $n^2$, and there are typically $n^{d-6}$ spanning clusters in a cube. Thus, Proposition~\ref{prop.sprinkling.AB} roughly ensures the spanning cluster has high probability to connect some other spanning cluster in the neighbor cube.
		
		The rest of this section is organized as follows. The proof of Proposition~\ref{prop.sprinkling.AB} relies on the first moment and second moment method, so Sections~\ref{subsec.e.pivAB}, \ref{subsec.v.pivAB} are respectively devoted to the expectation and variance of $\piv(A,B)$. We prove the main proposition in Section~\ref{subsec.skl}, together with more related properties. 
		
		\subsection{Expectation of $\#\piv(A,B)$}\label{subsec.e.pivAB}
		
		We first give the upper bound and lower bound for the expectation of pivotal bonds. In the proof, we will use the notation $\overrightarrow{\partial} \cC$, which denotes oriented edges on the boundary of the cluster $\cC$
		\begin{align}\label{eq.def.bd_clt}
			\overrightarrow{\partial} \cC := \Ll\{(e_{-}, e_+) \in \overrightarrow{E}: e_{-} \in \cC, e_+ \notin \cC \Rr\}.
		\end{align}

		\begin{prop}\label{prop.mean.pivAB}
			There exist finite positive constants $c (d)$ and $C(d)$, such that conditioned on $\fC_{13}$, for all regular sets $A \subset \partial_{1}^\circ Q_1 \cap \fC_{1}$ and $B \subset \partial_{-1}^\circ Q_3 \cap \fC_{3}$, the following estimate holds for $n \gg 1$
			\begin{equation}\label{eq.mean.pivAB}
				c n^{-d+4} |A| \cdot |B| \leq \E_{\p_c} \Ll[ \# \piv(A,B) \big\vert \fC_{13} \Rr] \leq C n^{-d+4} |A| \cdot |B|.
			\end{equation}
		\end{prop}
		\begin{proof}
			The proof of upper bound is direct, because the definition of pivotal edge in \eqref{eq.def.piv_AB} implies the disjoint realization of open path connecting $A, B$ respectively
			\begin{align*}
				&\E_{\p_c} \Ll[ \# \piv(A,B) \big\vert \fC_{13} \Rr]\\
				&\leq \sum_{(e_{-}, e_+) \in \overrightarrow{E}(S_2)} \P_{\p_c} \Ll[ \exists x \in A, y \in B, \{x \xleftrightarrow{T_{12} \setminus \mathcal{E}(\fC_{1})} e_{-} \} \circ \{y \xleftrightarrow{T_{23} \setminus \mathcal{E}(\fC_{3})} e_+ \}  \big\vert \fC_{13} \Rr] \\
				&\leq \sum_{(e_{-}, e_+) \in \overrightarrow{E}(S_2)} \sum_{x \in A} \sum_{y \in B}\P_{\p_c} \Ll[ x \xleftrightarrow{T_{12} \setminus \mathcal{E}(\fC_{1})} e_{-}  \big\vert \fC_{13} \Rr] \P_{\p_c} \Ll[ y \xleftrightarrow{T_{23} \setminus \mathcal{E}(\fC_{3})} e_+   \big\vert \fC_{13} \Rr]\\
				&\leq \sum_{(e_{-}, e_+) \in \overrightarrow{E}(S_2)} \sum_{x \in A} \sum_{y \in B}\P_{\p_c} \Ll[ x \xleftrightarrow{T_{12} } e_{-} \Rr] \P_{\p_c} \Ll[ y \xleftrightarrow{T_{23}} e_+  \Rr].
			\end{align*}
			From the third line to the forth line, we relax the constrain, so there are more choices for the open paths. Then we insert the two-point function estimate 
			\begin{align*}
				\E_{\p_c} \Ll[ \# \piv(A,B) \big\vert \fC_{13} \Rr] \leq \sum_{(e_{-}, e_+) \in \overrightarrow{E}(S_2)} \sum_{x \in A} \sum_{y \in B} n^{-d+2} \cdot n^{-d+2} \lesssim n^{-d+4} |A| \cdot |B|.
			\end{align*}
			This is the desired upper bound in \eqref{eq.mean.pivAB}.
			
			The proof of lower bound is more technical, as the condition $\cC_{12}(e_{-}) \cap \cC_{23}(e_{+}) = \emptyset$ has more constrain than disjoint realization. The proof can be divided into $3$ steps.
			
			\textit{Step~1: reduction.} We abuse the notation and write $\piv(x,B) \equiv \piv(\{x\},B)$ when $A = \{x\}$ only contains a singe vertex. Then it suffices to study the following subset of $\piv(x,B)$
			\begin{align}\label{eq.def.piv_AB_star}
				\piv_*(x,B) := \Ll\{(e_{-}, e_+) \in \piv(x,B): \cC_{12}(e_{-}) \cap \fC_{1} \cap \partial Q_1 = \{x\}\Rr\}.
			\end{align}
			In another word, the vertex $x$ is the unique one on the boundary connecting to the pivotal bonds. This brings a disjoint union structure 
			\begin{align*}
				\bigsqcup_{x \in A} \piv_*(x,B) \subset \piv(A,B), 
			\end{align*}
			
			We aim to prove that
			\begin{align}\label{eq.mean.pivAB.unique}
				\E_{\p_c} \Ll[ \# \piv_*(x,B) \big\vert \fC_{13} \Rr] \gtrsim n^{-d+4}  |B|,
			\end{align}
			which will entail the lower bound in \eqref{eq.mean.pivAB}
			\begin{equation}\label{eq.mean.pivAB.reduction}
				\begin{split}
					\E_{\p_c} \Ll[ \# \piv(A,B) \big\vert \fC_{13} \Rr] &\geq \E_{\p_c} \Ll[ \# \bigsqcup_{x \in A} \piv_*(x,B) \big\vert \fC_{13} \Rr] \\
					&= \sum_{x \in A} \E_{\p_c} \Ll[ \# \piv_*(x,B) \big\vert \fC_{13} \Rr]\\
					&\gtrsim \sum_{x \in A}  n^{-d+4} |B| \\
					&\gtrsim  n^{-d+4}  |A|  \cdot |B|. 
				\end{split}
			\end{equation}

			\textit{Step~2: two-point function without hitting $\cC_{12}(x)$.} We thus study \eqref{eq.mean.pivAB.unique} and first express the expectation conditioned on the cluster $\cC_{12}(x)$
			\begin{multline}\label{eq.id.pivot123}
				\E_{\p_c} \Ll[ \# \piv_*(x,B) \big\vert \cC_{12}(x), \fC_{13} \Rr] = \E_{\p_c} \Big[ \Ind{\cC_{12}(x) \cap \fC_{1} \cap \partial Q_1 = \{x\} }  \\
				\sum_{ \substack{ (e_{-}, e_+) \in \overrightarrow{\partial} \cC_{12}(x) \\ e_{-} \in S_2 }} \P_{\p_c} \Ll[\{e_+ \xleftrightarrow{T_{23} \setminus \cC_{12}(x) \setminus \mathcal{E}(\fC_{3})} B \} \setminus \{e_+ \xleftrightarrow{T_{23} \setminus \mathcal{E}(\fC_{3})}  \partial_{1}S_2 \cap \cC_{12}(x)\} \Big\vert \fC_{13} , \cC_{12}(x) \Rr] \Big\vert \fC_{13} \Big].
			\end{multline}
			Let us justify the identity \eqref{eq.id.pivot123}. Conditioned on $\cC_{12}(x)$, the counting of $\piv_*(x,B)$ is naturally over the boundary $\overrightarrow{\partial} \cC_{12}(x)$ and in $ S_2 $. The indicator $\Ind{\cC_{12}(x) \cap \fC_{1} \cap \partial Q_1 = \{x\} } $ also gives the unique intersection condition. Afterwards, viewing  $e_+ \xleftrightarrow{T_{23} \setminus \mathcal{E}(\fC_{3})} B$ and ${\cC_{12}(e_{-}) \cap \cC_{23}(e_{+}) = \emptyset}$ in \eqref{eq.def.piv_AB}, the event $\{e_+ \xleftrightarrow{T_{23} \setminus \cC_{12}(x) \setminus \mathcal{E}(\fC_{3})} B \}$ is a necessary. We also need to exclude other cases of intersection between $\cC_{12}(e_{-})$ and $\cC_{23}(e_{+})$. The boundary $ \overrightarrow{\partial} \cC_{12}(x)$ has already been revealed conditioned on $\cC_{12}(x)$, thus the only case is  $\{e_+ \xleftrightarrow{T_{23} \setminus \mathcal{E}(\fC_{3})}  \partial_{1}S_2 \cap \cC_{12}(x) \}$, which explains the second part in \eqref{eq.id.pivot123}.
			
			The lower bound of the first conditional expectation is technical, as we need to avoid touching a random cluster. This part is treated in Lemma~\ref{lem.avoid} below separately, and we cite directly the conclusion in \eqref{eq.avoid}, which yields that
			\begin{multline*}
				\E_{\p_c} \Ll[ \# \piv_*(x,B) \big\vert \cC_{12}(x), \fC_{13} \Rr] \gtrsim n^{-d+2} \vert B \vert \times\\
				\E_{\p_c} \Ll[  \frac{ \Ind{\cC_{12}(x) \cap \fC_{1} \cap \partial Q_1 = \{x\} }\Ll(\vert Q_2 \cap \cC_{12}(x)\vert - n^{-d+4} \vert Q_2 \cap \cC_{12}(x)\vert \cdot \vert \partial_{1}S_2 \cap \cC_{12}(x)\vert\Rr)_+^2}{ \sum_{ \substack{ (e_{-}, e_+) \in \overrightarrow{\partial} \cC_{12}(x) \\ e_{-} \in S_2 }} \Ll(\sum_{z \in Q_2 \cap \cC_{12}(x)} \P_{\p_c}\Ll[{z \xleftrightarrow{T_{12}} e_{-}}  \Rr] \Rr)^2}\Big\vert \cC_{12}(x), \fC_{13} \Rr].
			\end{multline*}
			Here the notation $(\cdot)_+$ is defined as $(\cdot)_+ := \max\{\cdot, 0\}$. For two random variables $X,Y \geq 0$, the Cauthy--Schwarz inequality gives
			\begin{equation*}
				\E \Ll[ \frac{X^2}{Y} \Rr] \E[Y] \geq \E [X]^2,
			\end{equation*}
			thus
			\begin{equation*}
				\E \Ll[ \frac{X^2}{Y} \Rr] \geq \frac{\E [X]^2}{\E [Y]}.
			\end{equation*}
			We apply this to the estimate above under $\E_{\p_c} \Ll[ \, \cdot \, \Big\vert \fC_{13} \Rr]$, and get
			\begin{align}\label{eq.mean.pivAB.middle}
				\E_{\p_c} \Ll[ \# \piv_*(x,B) \big\vert \fC_{13} \Rr] 	\gtrsim n^{-d+2} \vert B \vert \cdot  \frac{\E [X \vert \fC_{13}]^2}{\E [Y \vert \fC_{13}]},
			\end{align}
			with the terms $X, Y$ defined as
			\begin{align*}
				X &:= \Ind{\cC_{12}(x) \cap \fC_{1} \cap \partial Q_1 = \{x\} } \Ll( \vert Q_2 \cap \cC_{12}(x)\vert - n^{-d+4} \vert Q_2 \cap \cC_{12}(x)\vert \cdot \vert \partial_{1}S_2 \cap \cC_{12}(x)\vert\Rr)_+, \\
				Y &:= \sum_{ \substack{ (e_{-}, e_+) \in \overrightarrow{\partial} \cC_{12}(x) \\ e_{-} \in S_2 }} \Ll(\sum_{z \in Q_2 \cap \cC_{12}(x)} \P_{\p_c}\Ll[ {z \xleftrightarrow{T_{12}} e_{-}}  \Rr] \Rr)^2.
			\end{align*}
			
			\textit{Step~3: lower bound of nominator.}
			It remains to estimate the nominator and denominator in \eqref{eq.mean.pivAB.middle}. Concerning the nominator, we can write it as 
			\begin{align*}
				X &= (X_1 - X_2)_+, \\
				X_1 &:=  \Ind{\cC_{12}(x) \cap \fC_{1} \cap \partial Q_1 = \{x\} }\vert Q_2 \cap \cC_{12}(x)\vert, \\
				X_2 &:= \Ind{\cC_{12}(x) \cap \fC_{1} \cap \partial Q_1 = \{x\} } n^{-d+4} \vert Q_2 \cap \cC_{12}(x)\vert \cdot \vert \partial_{1}S_2 \cap \cC_{12}(x)\vert.
			\end{align*}
			
			The lower bound of $X_1$ is a consequence of regular points and Proposition~\ref{prop.extension}
			\begin{equation}\label{eq.mean.pivAB.nominator1}
				\begin{split}
					\E_{\p_c} \Ll[X_1 \Big\vert \fC_{13} \Rr] &=\E_{\p_c} \Ll[ \Ind{\cC_{12}(x) \cap \fC_{1} \cap \partial Q_1 = \{x\} } \vert Q_2 \cap \cC_{12}(x)\vert \Big\vert \fC_{13} \Rr] \\
					& = \sum_{z \in Q_2} \P_{\p_c} \Ll[ z \xleftrightarrow{T_{12} \setminus \mathcal{E}(\fC_{1})} x \setminus \bigcup_{y \in \fC_{1} \cap \partial Q_1 \setminus \{x \}} y \xleftrightarrow{T_{12} \setminus \mathcal{E}(\fC_{1})}  x  \Big\vert \fC_{13} \Rr]\\
					&\gtrsim n^d \cdot n^{-d+2}  \gtrsim n^2.
				\end{split}
			\end{equation}
			The upper bound of $X_2$ can be estimated by the tree-graph inequality
			\begin{equation}\label{eq.mean.pivAB.nominator2}
				\begin{split}
					\E_{\p_c} \Ll[X_2 \Big\vert \fC_{13} \Rr] &\leq \E_{\p_c} \Ll[n^{-d+4} \vert Q_2 \cap \cC_{12}(x)\vert \cdot \vert \partial_{1}S_2 \cap \cC_{12}(x)\vert \Big\vert \fC_{13} \Rr] \\
					& = n^{-d+4}\sum_{z \in Q_2} \sum_{y \in \partial_{1}S_2 }\P_{\p_c} \Ll[ x \xleftrightarrow{T_{12} \setminus \mathcal{E}(\fC_{1})} z \xleftrightarrow{T_{12} \setminus \mathcal{E}(\fC_{1})} y  \Big\vert \fC_{13} \Rr]\\
					&\lesssim n^{-d+4}\sum_{z \in Q_2} \sum_{y \in \partial_{1}S_2 } \tau^{T_{12}}_{\p_c}(x,z,y)\\
					&\lesssim n^{-d+4}\sum_{z \in Q_2} \sum_{y \in \partial_{1}S_2 } \sum_{u \in T_{12}} \tau^{T_{12}}_{\p_c}(x,u)\tau^{T_{12}}_{\p_c}(u,y)\tau^{T_{12}}_{\p_c}(u,z)\\
					&\lesssim n^{-d+4}\sum_{z \in Q_2} \sum_{y \in \partial_{1}S_2 } n^{-2d+5} \lesssim n^{-d+8}.
				\end{split}
			\end{equation}
			The tree-graph inequality is inserted in the forth line. The bound relies on the fact that the distance between every pair in $x, y, z$ is of order $n$, and the half-plane bound
			\begin{align*}
				\tau^{T_{12}}_{\p_c}(u,y) \lesssim \vert u - y \vert^{-d+1}.
			\end{align*}
			We refer to the recent reference \cite{panis2026} for this estimate.
			
			As the high-dimensional percolation lives in $\Zd$ with $d > 6$,  we compare the results \eqref{eq.mean.pivAB.nominator1} and \eqref{eq.mean.pivAB.nominator2}, and observe that for $n \gg 1$,
			\begin{align*}
				\E_{\p_c} \Ll[X_1 \vert \fC_{13} \Rr] \gtrsim n^2 \gg n^{-d+8} \gtrsim \E_{\p_c} \Ll[X_2 \vert \fC_{13} \Rr].
			\end{align*}
			Therefore, we conclude the lower bound of the nominator when $n \gg 1$
			\begin{align*}
				\E_{\p_c} \Ll[X \vert \fC_{13} \Rr] = \E_{\p_c} \Ll[(X_1 - X_2 )_+ \vert \fC_{13} \Rr] \geq \E_{\p_c} \Ll[X_1 - X_2 \vert \fC_{13} \Rr] \gtrsim n^2, 
			\end{align*}
			which implies that 
			\begin{align}\label{eq.mean.pivAB.nominator}
				\E_{\p_c} \Ll[X \vert \fC_{13} \Rr]^2 \gtrsim n^4.
			\end{align}
			
			\textit{Step~4: upper bound of denominator.} The upper bound of the denominator will be calculated as follows
			\begin{align*}
				&\E_{\p_c} \Ll[\sum_{ \substack{ (e_{-}, e_+) \in \overrightarrow{\partial} \cC_{12}(x) \\ e_{-} \in S_2 }} \Ll(\sum_{z \in Q_2 \cap \cC_{12}(x)} \P_{\p_c}\Ll[ {z \xleftrightarrow{T_{12}} e_{-}}  \Rr] \Rr)^2\Big\vert \fC_{13} \Rr]\\
				&= \E_{\p_c} \Ll[\sum_{ \substack{ (e_{-}, e_+) \in \overrightarrow{\partial} \cC_{12}(x) \\ e_{-} \in S_2 }} \sum_{z,z' \in Q_2 \cap \cC_{12}(x)} \P_{\p_c}\Ll[ {z \xleftrightarrow{T_{12}} e_{-}} \Rr] \P_{\p_c}\Ll[ {z' \xleftrightarrow{T_{12}} e_{-}}  \Rr] \Big\vert \fC_{13} \Rr]\\
				&\lesssim \sum_{w, z,z' \in Q_2} \P_{\p_c}\Ll[ w, z,z' \in \cC_{12}(x) \Big\vert \fC_{13} \Rr] \P_{\p_c}\Ll[ {z \xleftrightarrow{T_{12}} w} \Rr] \P_{\p_c}\Ll[ {z' \xleftrightarrow{T_{12}} w}  \Rr]. 
			\end{align*}
			In the last line, we make a change of variable $w = e_{-}$. Notice the quantity $ \P_{\p_c}\Ll[ w, z,z' \in \cC_{12}(x) \Big\vert \fC_{13} \Rr] $ can be expresses via the $4$-point function
			\begin{align*}
				&\P_{\p_c}\Ll[ w, z,z' \in \cC_{12}(x) \Big\vert \fC_{13} \Rr]\\
				&=  \P_{\p_c}\Ll[ x \xleftrightarrow{T_{12} \setminus \mathcal{E}(\fC_{1})} w \xleftrightarrow{T_{12} \setminus \mathcal{E}(\fC_{1})} z \xleftrightarrow{T_{12} \setminus \mathcal{E}(\fC_{1})} z' \Big\vert \fC_{13} \Rr]\\
				&\leq  \P_{\p_c}\Ll[ x \xleftrightarrow{T_{12} } w \xleftrightarrow{T_{12}} z \xleftrightarrow{T_{12}} z' \Rr]\\
				&= \tau^{T_{12}}_{\p_c}(x, w, z, z').
			\end{align*}
			In the third line above, we relax the domain of connection. Therefore,  we obtain a preliminary upper bound for the denominator
			\begin{multline}\label{eq.mean.pivAB.denominator.pre}
				\E_{\p_c} \Ll[\sum_{ \substack{ (e_{-}, e_+) \in \overrightarrow{\partial} \cC_{12}(x) \\ e_{-} \in S_2 }} \Ll(\sum_{z \in Q_2 \cap \cC_{12}(x)} \P_{\p_c}\Ll[ {z \xleftrightarrow{T_{12}} e_{-}}  \Rr] \Rr)^2\Big\vert \fC_{13} \Rr] \\
				\lesssim \sum_{w, z,z' \in Q_2} \tau^{T_{12}}_{\p_c}(x, w, z, z') \tau^{T_{12}}_{\p_c}(w, z) \tau^{T_{12}}_{\p_c}(w, z').
			\end{multline}
			
			To further handle the $4$-point function on the {\rhs} of \eqref{eq.mean.pivAB.denominator.pre}, we use the tree-graph inequality that  
			\begin{align*}
				\E_{\p_c} \Ll[\sum_{ \substack{ (e_{-}, e_+) \in \overrightarrow{\partial} \cC_{12}(x) \\ e_{-} \in S_2 }} \Ll(\sum_{z \in Q_2 \cap \cC_{12}(x)} \P_{\p_c}\Ll[ {z \xleftrightarrow{T_{12}} e_{-}}  \Rr] \Rr)^2\Big\vert \fC_{13} \Rr] \leq \mathbf{I} + \mathbf{II},
			\end{align*}
			where the terms $\mathbf{I}$ and $\mathbf{II}$ are 
			\begin{equation}\label{eq.mean.pivAB.denominator.decom}
				\begin{split}
					\mathbf{I} := \sum_{w, z,z' \in Q_2} \sum_{u,v \in T_{12}} \tau_{\p_c}(x, u) \tau_{\p_c}(u,v) \tau_{\p_c}(v, z) \tau_{\p_c}(v, w) \tau_{\p_c}(w, z) \tau_{\p_c}(u,z') \tau_{\p_c}(w, z'), \\
					\mathbf{II} := \sum_{w, z,z' \in Q_2} \sum_{u,v \in T_{12}} \tau_{\p_c}(x, u) \tau_{\p_c}(u,v) \tau_{\p_c}(v, z) \tau_{\p_c}(v, z') \tau_{\p_c}(u, w) \tau_{\p_c}(w, z)  \tau_{\p_c}(w, z').
				\end{split}
			\end{equation}
			When developing the inequality above, we keep the domain $T_{12}$ at first for the adding vertices, and then use $\tau^{T_{12}}_{\p_c} \leq \tau_{\p_c}$ to remove the superscript.  As these two terms were analyzed in \cite[Proposition~1.8, Figure~1]{hutchcroft2022derivation}, we cite the result there
			\begin{align*}
				\mathbf{I} , \mathbf{II} \leq \sum_{u  \in T_{12}} \tau_{\p_c}(x, u) \nabla_{\p_c}^2 \lesssim n^2.
			\end{align*}
			Therefore, we obtain an upper bound for the denominator 
			\begin{align}\label{eq.mean.pivAB.denominator}
				\E_{\p_c} \Ll[\sum_{ \substack{ (e_{-}, e_+) \in \overrightarrow{\partial} \cC_{12}(x) \\ e_{-} \in S_2 }} \Ll(\sum_{z \in Q_2 \cap \cC_{12}(x)} \P_{\p_c}\Ll[ {z \xleftrightarrow{T_{12}} e_{-}}  \Rr] \Rr)^2\Big\vert \fC_{13} \Rr] \lesssim n^2.
			\end{align}
			
			Finally, we combine \eqref{eq.mean.pivAB.nominator} and \eqref{eq.mean.pivAB.denominator}, and put them back to \eqref{eq.mean.pivAB.middle}
			\begin{align*}
				\E_{\p_c} \Ll[ \# \piv_*(x,B) \big\vert \fC_{13} \Rr]  \gtrsim \frac{n^{-d+2} \vert B \vert \cdot n^4}{n^2} \gtrsim n^{-d+4} \vert B \vert.
			\end{align*}
			As discussed at the very beginning of the proof, this justifies the claim \eqref{eq.mean.pivAB.unique} and thus proves \eqref{eq.mean.pivAB}.		\qedhere

		\end{proof}

		The following technical estimate about lower bound is inspired by \cite[Lemma~2.1]{hutchcroft2022derivation}.
		\begin{lemma}\label{lem.avoid}
			Under the same setting as Proposition~\ref{prop.mean.pivAB}, we have the following estimate
			\begin{align}
				\sum_{ \substack{ (e_{-}, e_+) \in \overrightarrow{\partial} \cC_{12}(x) \\ e_{-} \in S_2 }} \P_{\p_c} \Ll[\{e_+ \xleftrightarrow{T_{23} \setminus \cC_{12}(x) \setminus \mathcal{E}(\fC_{3})} B \} \setminus \{e_+ \xleftrightarrow{T_{23} \setminus \mathcal{E}(\fC_{3})}  \partial_{1}S_2 \cap \cC_{12}(x)\} \Big\vert \fC_{13} , \cC_{12}(x) \Rr] \nonumber  \\
				\gtrsim \frac{n^{-d+2} \vert B \vert \ \Ll(\vert Q_2 \cap \cC_{12}(x)\vert - n^{-d+4} \vert Q_2 \cap \cC_{12}(x)\vert \cdot \vert \partial_{1}S_2 \cap \cC_{12}(x)\vert\Rr)_+^2}{ \sum_{ \substack{ (e_{-}, e_+) \in \overrightarrow{\partial} \cC_{12}(x) \\ e_{-} \in S_2 }} \Ll(\sum_{z \in Q_2 \cap \cC_{12}(x)} \P_{\p_c}\Ll[{z \xleftrightarrow{T_{12}} e_{-}}  \Rr] \Rr)^2}. \label{eq.avoid}
			\end{align}
			The notation is defined as $(\cdot)_+ := \max\{\cdot, 0\}$. 
		\end{lemma}
		\begin{proof}
			For convenience, we define an event  in the proof 
			\begin{align}\label{eq.def.E_avoid_boundary}
				E^{S}(z, B; F; V) := \{z \xleftrightarrow{S \setminus  F} B \} \setminus \{z \xleftrightarrow{S} V\}.
			\end{align}
			Roughly, we are interested in th event $\{z \leftrightarrow B\}$, and require it happens on $S \setminus F$ and its cluster does not touch $V$. Using the following shorthand notation
			\begin{align}\label{eq.avoid.shorthand}
				\mathcal{S} := T_{23} \setminus \mathcal{E}(\fC_{3}),  \qquad \mathcal{V} := \partial_{1}S_2 \cap \cC_{12}(x),
			\end{align}
			the probability on the {\lhs} of \eqref{eq.avoid} is a concrete example $z = e_+$, $F = \cC_{12}(x)$, $S = \mathcal{S}$, and $V = \mathcal{V}$ 
			\begin{align}\label{eq.avoid.LHS}
				E^{\mathcal{S}}(e_+, B; \cC_{12}(x); \mathcal{V}).
			\end{align}
			
			To study the lower bound \eqref{eq.avoid.LHS}, let us equip the probability space with $\{\a'(e)\}_{e \in \Ed}$ as another independent $\p_c$-percolation. The trick in \cite[Lemma~2.1]{hutchcroft2022derivation} suggests to shift the event to another vertex $z \in  Q_2 \cap \cC_{12}(x)$ and consider the probability of
			\begin{align}\label{eq.avoid.object}
				E^{\a', \mathcal{S}}(z, B; \emptyset; \mathcal{V}),
			\end{align}
			which is 
			\begin{align*}
				\{z \xleftrightarrow{\a', \mathcal{S}} B\} \setminus \{z \xleftrightarrow{\a', \mathcal{S}} \mathcal{V}\}.
			\end{align*}
			Here $``\xleftrightarrow{\a', \mathcal{S}} "$ means that the connection by the open path of $\a'$ lives in $\mathcal{S}$. We will derive an upper bound and a lower bound of \eqref{eq.avoid.object}, and the comparison then yields \eqref{eq.avoid}. The proof can be divided naturally in $3$ steps.

			\textit{Step~1: upper bound.}
			Since  $\a'$ is another copy of percolation, its path can still hit the explored cluster in $\a$. We thus decompose the open path of ${\{z \xleftrightarrow{\a', \mathcal{S}} B  \}}$ at its last intersection at some $e_{-} \in \cC_{12}(x) \cap S_2$, i.e.
			\begin{align*}
				\{z \xleftrightarrow{\a', \mathcal{S}} B  \} \subset \bigcup_{ \substack{ (e_{-}, e_+) \in \overrightarrow{\partial} \cC_{12}(x) \\ e_{-} \in S_2 }}	\{z \xleftrightarrow{\a', \mathcal{S}} e_{-}  \} \circ \{\a'(\{e_{-}, e_+\}) = 1\} \circ \{e_+ \xleftrightarrow{\a', \mathcal{S} \setminus \cC_{12}(x)} B \}.
			\end{align*}
			Thanks to the later intersection, the path after $e_+$ does not contain any intersection with $\cC_{12}(x)$. After further excluding the path intersecting $\mathcal{V}$ and using the BK inequality, we obtain that 
			\begin{align*}
				&\P_{\p_c} \Ll[E^{\a', \mathcal{S}}(z, B; \emptyset; \mathcal{V})  \Big\vert \fC_{13} , \cC_{12}(x) \Rr]\\ 
				&\leq \sum_{\substack{ (e_{-}, e_+) \in \overrightarrow{\partial} \cC_{12}(x) \\ e_{-} \in S_2 }} \p_c \P_{\p_c} \Ll[z \xleftrightarrow{\a', \mathcal{S}} e_{-}  \Big\vert \fC_{13} , \cC_{12}(x) \Rr] \P_{\p_c} \Ll[\{e_+ \xleftrightarrow{\a', \mathcal{S} \setminus \cC_{12}(x)} B\} \setminus \{e_+ \xleftrightarrow{\a', \mathcal{S}} \mathcal{V} \}  \Big\vert \fC_{13} , \cC_{12}(x) \Rr] \\
				&\leq \sum_{\substack{ (e_{-}, e_+) \in \overrightarrow{\partial} \cC_{12}(x) \\ e_{-} \in S_2 }} \p_c \P_{\p_c} \Ll[z \xleftrightarrow{T_{23}} e_{-}\Rr]  \P_{\p_c} \Ll[\{e_+ \xleftrightarrow{\mathcal{S} \setminus \cC_{12}(x)} B\} \setminus \{e_+ \xleftrightarrow{\mathcal{S}} \mathcal{V} \}  \Big\vert \fC_{13} , \cC_{12}(x) \Rr] .
			\end{align*}
			From the second line to the third line, we drop the constrain for the first term. For the second term, because the path will no longer touch $\cC_{12}(x)$, the expectation is same under $\a$ or $\a'$. This is the desired term on the {\lhs} of \eqref{eq.avoid}. The upper bound of \eqref{eq.avoid.object} is obtained
			\begin{multline}\label{eq.avoid_upper}
				\P_{\p_c} \Ll[E^{\a', \mathcal{S}}(z, B; \emptyset; \mathcal{V})  \Big\vert \fC_{13} , \cC_{12}(x) \Rr] \\
				\leq \sum_{\substack{ (e_{-}, e_+) \in \overrightarrow{\partial} \cC_{12}(x) \\ e_{-} \in S_2 }} \p_c \P_{\p_c} \Ll[z \xleftrightarrow{T_{23}} e_{-}\Rr]  \P_{\p_c} \Ll[E^{\mathcal{S}}(e_+, B; \cC_{12}(x); \mathcal{V}) \Big\vert \fC_{13} , \cC_{12}(x) \Rr]. 
			\end{multline}

			\textit{Step~2: lower bound.} Now we turn to the lower bound of \eqref{eq.avoid.object}, which can be written as 
			\begin{align}\label{eq.avoid.object.again}
				E^{\a', \mathcal{S}}(z, B; \emptyset; \mathcal{V}) = 	\{z \xleftrightarrow{\a', \mathcal{S}} B \} \setminus \{z \xleftrightarrow{\a', \mathcal{S}} \mathcal{V} \xleftrightarrow{\a', \mathcal{S}} B \}.
			\end{align} 
			
			For its first part of \eqref{eq.avoid.object.again}, we have
			\begin{equation}\label{eq.avoid_lower1}
				\begin{split}
					&\P_{\p_c} \Ll[z \xleftrightarrow{\a', \mathcal{S}} B  \Big\vert \fC_{13} , \cC_{12}(x) \Rr]\\
					&=\P_{\p_c} \Ll[z \xleftrightarrow{\mathcal{S}} B  \Big\vert \fC_{13}\Rr]\\
					&\geq \sum_{y \in B} \P_{\p_c} \Ll[z \xleftrightarrow{\mathcal{S}} y,  \cC_{23}(z) \cap \fC_{3} \cap \partial Q_3 = \{y\} \Big\vert \fC_{13}\Rr]\\
					&\gtrsim n^{-d+2} \vert B\vert.
				\end{split}
			\end{equation}
			Here the second line makes use of the independence. The argument above is quite similar as \eqref{eq.mean.pivAB.reduction}, and Proposition~\ref{prop.extension} is applied in the last line.
			
			For the second part of \eqref{eq.avoid.object.again}, it has an upper bound
			\begin{equation}\label{eq.avoid_lower2}
				\begin{split}
					&\P_{\p_c} \Ll[z \xleftrightarrow{\a', \mathcal{S}} \mathcal{V} \xleftrightarrow{\a', \mathcal{S}} B  \Big\vert \fC_{13} , \cC_{12}(x) \Rr]\\
					&= \P_{\p_c} \Ll[ \exists y \in B, v \in \mathcal{V}, z \xleftrightarrow{\a', \mathcal{S}} v \xleftrightarrow{\a', \mathcal{S}} y  \Big\vert \fC_{13} , \cC_{12}(x) \Rr]\\
					&\leq \sum_{y \in B}\sum_{v \in \mathcal{V}}\P_{\p_c} \Ll[ z \xleftrightarrow{\a', \mathcal{S}} v \xleftrightarrow{\a', \mathcal{S}} y  \Big\vert \fC_{13} , \cC_{12}(x) \Rr]\\
					&\leq \sum_{y \in B}\sum_{v \in \mathcal{V}} \sum_{u \in \mathcal{S}} \tau_{\p_c}(z, u) \tau_{\p_c}(u, y) \tau_{\p_c}(u, v) \\
					&\lesssim \vert B\vert \cdot \vert \mathcal{V}\vert n^{-2d+6}.
				\end{split}
			\end{equation}
			Recall the definition of $\mathcal{V} = \partial_{1}S_2 \cap \cC_{12}(x)$ in \eqref{eq.avoid.shorthand}, which  is $ \cC_{12}(x)$-measurable, so we have the third line. Then the tree-graph inequality is applied in the forth line, and it gives the bound as the distance between every pair in $z, y, v$ is larger than $n$.
			
			Therefore, we get a lower bound of \eqref{eq.avoid.object}
			\begin{align}\label{eq.avoid_lower}
				\P_{\p_c} \Ll[E^{\a', \mathcal{S}}(z, B; \emptyset; \mathcal{V})  \Big\vert \fC_{13} , \cC_{12}(x) \Rr] \gtrsim n^{-d+2} \vert B\vert \Ll(1 - n^{-d+4} \vert \mathcal{V}\vert\Rr).
			\end{align}

			\textit{Step~3: comparison.}
			Combining \eqref{eq.avoid_lower} and \eqref{eq.avoid_upper}, the double percolation configuration argument helps develop an inequality 
			\begin{multline*}
				n^{-d+2} \vert B\vert \Ll(1 - n^{-d+4} \vert \mathcal{V}\vert\Rr) \\
				\lesssim 	\sum_{\substack{ (e_{-}, e_+) \in \overrightarrow{\partial} \cC_{12}(x) \\ e_{-} \in S_2 }} \P_{\p_c} \Ll[z \xleftrightarrow{T_{23}} e_{-}\Rr]  \P_{\p_c} \Ll[E^{\mathcal{S}}(e_+, B; \cC_{12}(x); \mathcal{V}) \Big\vert \fC_{13} , \cC_{12}(x) \Rr].
			\end{multline*}
			We then sum over all $z \in Q_2 \cap \cC_{12}(x)$ and obtain 
			\begin{multline}\label{eq.avoid.middle}
				n^{-d+2} \vert B\vert \Ll(1 - n^{-d+4} \vert \mathcal{V}\vert\Rr) \vert Q_2 \cap \cC_{12}(x) \vert  \\
				\lesssim \sum_{\substack{ (e_{-}, e_+) \in \overrightarrow{\partial} \cC_{12}(x) \\ e_{-} \in S_2 }} \Ll(\sum_{z \in Q_2 \cap \cC_{12}(x)}\P_{\p_c} \Ll[z \xleftrightarrow{T_{23}} e_{-}\Rr]\Rr) \P_{\p_c} \Ll[E^{\mathcal{S}}(e_+, B; \cC_{12}(x); \mathcal{V}) \Big\vert \fC_{13} , \cC_{12}(x) \Rr].
			\end{multline}
			
			It suffices to reformulate \eqref{eq.avoid.middle} as that in \eqref{eq.avoid}. We further apply a Cauchy--Schwarz inequality to its {\rhs}
			\begin{align*}
				&	n^{-d+2} \vert B\vert \Ll(1 - n^{-d+4} \vert \mathcal{V}\vert\Rr) \vert Q_2 \cap \cC_{12}(x) \vert  \\
				&\lesssim \Ll(\sum_{\substack{ (e_{-}, e_+) \in \overrightarrow{\partial} \cC_{12}(x) \\ e_{-} \in S_2 }} \Ll(\sum_{z \in Q_2 \cap \cC_{12}(x)}\P_{\p_c} \Ll[z \xleftrightarrow{T_{23}} e_{-}\Rr]\Rr)^2 \Rr)^{\frac{1}{2}} \\
				& \qquad \times \Ll(\sum_{\substack{ (e_{-}, e_+) \in \overrightarrow{\partial} \cC_{12}(x) \\ e_{-} \in S_2 }}\P_{\p_c} \Ll[E^{\mathcal{S}}(e_+, B; \cC_{12}(x); \mathcal{V}) \Big\vert \fC_{13} , \cC_{12}(x) \Rr]^2\Rr)^{\frac{1}{2}}.
			\end{align*}
			We replace one $\P_{\p_c} \Ll[E^{\mathcal{S}}(e_+, B; \cC_{12}(x); \mathcal{V}) \Big\vert \fC_{13} , \cC_{12}(x) \Rr]$ by its upper bound 
			\begin{align*}
				\P_{\p_c} \Ll[E^{\mathcal{S}}(e_+, B; \cC_{12}(x); \mathcal{V}) \Big\vert \fC_{13} , \cC_{12}(x) \Rr] \leq \sum_{y  \in B}\P_{\p_c} \Ll[e_+ \xleftrightarrow{T_{23}} y  \Rr] \lesssim n^{-d+2} \vert B\vert,
			\end{align*}
			then we get 
			\begin{align*}
				&	n^{-d+2} \vert B\vert \Ll(1 - n^{-d+4} \vert \mathcal{V}\vert\Rr) \vert Q_2 \cap \cC_{12}(x) \vert  \\
				&\lesssim \Ll(\sum_{\substack{ (e_{-}, e_+) \in \overrightarrow{\partial} \cC_{12}(x) \\ e_{-} \in S_2 }} \Ll(\sum_{z \in Q_2 \cap \cC_{12}(x)}\P_{\p_c} \Ll[z \xleftrightarrow{T_{23}} e_{-}\Rr]\Rr)^2 \Rr)^{\frac{1}{2}} \\
				& \qquad \times \Ll(\sum_{\substack{ (e_{-}, e_+) \in \overrightarrow{\partial} \cC_{12}(x) \\ e_{-} \in S_2 }}\P_{\p_c} \Ll[E^{\mathcal{S}}(e_+, B; \cC_{12}(x); \mathcal{V}) \Big\vert \fC_{13} , \cC_{12}(x) \Rr]\Rr)^{\frac{1}{2}} \times \Ll(n^{-d+2} \vert B\vert\Rr)^{\frac{1}{2}}.
			\end{align*}
			In order to make a square of the two sides, the {\lhs} should be positive. We thus add a truncation $(\cdot)_+$ to handle this issue, then \eqref{eq.avoid} is obtained. 
		\end{proof}

		\subsection{Variance of $\#\piv(A,B)$}\label{subsec.v.pivAB}

		\begin{prop}\label{prop.var.pivAB}
			Under the same setting as Proposition~\ref{prop.mean.pivAB}, we have
			\begin{multline}\label{eq.var.pivAB}
				\var_{\p_c} \Ll[ \# \piv(A,B) \big\vert \fC_{13} \Rr]  \lesssim  n^{-2d+10} \Ll(\vert A \vert \cdot \vert B \vert^2 + \vert A \vert^2 \cdot \vert B \vert\Rr) \\
				+ \Ll(n^{-2d+6} + n^{-3d+14}\Rr)  \vert A \vert^2 \cdot \vert B \vert^2  +  \vert A \vert \cdot \vert B \vert \Ll\{\begin{array}{ll}
					n^{-d+5} & d=7,\\
					n^{-d+4} \log n & d=8, \\
					n^{-d+4} & d \geq 9.
				\end{array}\Rr. 
			\end{multline}
		\end{prop}	
		\begin{proof}
			We use the notation $e = (e_{-}, e_{+}), e' = (e'_{-}, e'_{+})$ for oriented edges in the proof.
			\begin{align*}
				&\var_{\p_c} \Ll[ \# \piv(A,B) \big\vert \fC_{13} \Rr] \\
				&= \sum_{e \in \overrightarrow{E}(S_2)} \var_{\p_c} \Ll[  \Ind{e \in \piv(A,B)} \big\vert \fC_{13} \Rr] + \sum_{ \substack{e , e' \in \overrightarrow{E}(S_2), \\ e \neq e' }} \cov_{\p_c} \Ll[ \Ind{e \in \piv(A,B)},  \Ind{e' \in \piv(A,B)}\big\vert \fC_{13} \Rr]\\
				&\leq 	\E_{\p_c} \Ll[ \# \piv(A,B) \big\vert \fC_{13} \Rr] + \sum_{ \substack{e , e' \in \overrightarrow{E}(S_2), \\ e \neq e' }} \cov_{\p_c} \Ll[ \Ind{e \in \piv(A,B)},  \Ind{e' \in \piv(A,B)}\big\vert \fC_{13} \Rr].
			\end{align*}
			We then focus on the covariance term. Viewing the definition of \eqref{eq.def.piv_AB}, the event ${\{e \in \piv(A,B)\}}$ is measurable by 
			\begin{align}\label{eq.def.Cpiv}
				\cC_{\piv}(e) := \cC_{12}(e_{-}) \cup \cC_{23}(e_{+}).
			\end{align}
			Then a classical trick for the variance of observable in percolation is the following decomposition
			\begin{multline*}
				\P_{\p_c} \Ll[ e ,e' \in \piv(A,B)\big\vert \fC_{13} \Rr] 	= \P_{\p_c} \Ll[ e ,e' \in \piv(A,B), \cC_{\piv}(e)  \cap \cC_{\piv}(e') \neq \emptyset \big\vert \fC_{13} \Rr] \\
				+ \P_{\p_c} \Ll[ e ,e' \in \piv(A,B), \cC_{\piv}(e)  \cap \cC_{\piv}(e') = \emptyset \big\vert \fC_{13} \Rr].
			\end{multline*}
			The part $\cC_{\piv}(e)  \cap \cC_{\piv}(e') = \emptyset$ has the inclusion
			\begin{align*}
				\{e ,e' \in \piv(A,B), \cC_{\piv}(e)  \cap \cC_{\piv}(e') = \emptyset\} \subset  \{e \in \piv(A,B)\} \circ \{e' \in \piv(A,B)\},
			\end{align*}
			because the events $\{e \in \piv(A,B)\}$ and $\{e' \in \piv(A,B)\}$ are witnessed by disjoint sets. Then BK inequality applies
			\begin{multline*}
				\P_{\p_c} \Ll[ e ,e' \in \piv(A,B), \cC_{\piv}(e)  \cap \cC_{\piv}(e') = \emptyset \big\vert \fC_{13} \Rr] \\
				\leq \P_{\p_c} \Ll[ e \in \piv(A,B)\big\vert \fC_{13} \Rr] \P_{\p_c} \Ll[ e' \in \piv(A,B)\big\vert \fC_{13} \Rr].
			\end{multline*}
			The observation above yields that 
			\begin{multline}\label{eq.var.pivAB.pre}
				\cov_{\p_c} \Ll[ \Ind{e \in \piv(A,B)},  \Ind{e' \in \piv(A,B)}\big\vert \fC_{13} \Rr] \\
				\leq \P_{\p_c}\Ll[ e ,e' \in \piv(A,B), \cC_{\piv}(e)  \cap \cC_{\piv}(e') \neq \emptyset \big\vert \fC_{13} \Rr].
			\end{multline}
			
			We then further discuss the situation of the intersection $\cC_{\piv}(e)  \cap \cC_{\piv}(e') \neq \emptyset$, because each of them is the union of two clusters. We claim that $\cC_{12}(e_{-})$ can intersect at most one between $\cC_{12}(e'_{-})$ and $\cC_{23}(e'_{+})$. Otherwise, we would have 
			\begin{align*}
				\{\cC_{12}(e_{-}) \cap \cC_{12}(e'_{-}) \neq \emptyset\} \Longrightarrow \{\cC_{12}(e_{-}) = \cC_{12}(e'_{-})\},
			\end{align*}
			as the cluster is a class of equivalence. Afterwards, the equality above implies that
			\begin{align*}
				\{\cC_{12}(e_{-}) \cap \cC_{23}(e'_{+}) \neq \emptyset\} \Longleftrightarrow \{\cC_{12}(e'_{-}) \cap \cC_{23}(e'_{+}) \neq \emptyset\}.
			\end{align*}
			The last statement contradicts the condition in \eqref{eq.def.pivo_AB}, so the claim is valid. The claim also naturally generalizes to all the components in $\cC_{\piv}(e), \cC_{\piv}(e')$. For convenience, we define a matrix of intersection
			\begin{align}\label{eq.def.Matrix}
				M := \begin{pmatrix}
					\Ind{\cC_{12}(e_{-}) \cap \cC_{12}(e'_{-}) \neq \emptyset} & \Ind{\cC_{12}(e_{-}) \cap \cC_{23}(e'_{+})\neq \emptyset} \\
					\Ind{\cC_{23}(e_{+}) \cap \cC_{12}(e'_{-}) \neq \emptyset} & \Ind{\cC_{23}(e_{+}) \cap \cC_{23}(e'_{+}) \neq \emptyset}
				\end{pmatrix}.
			\end{align}
			Then the claim in this paragraph justifies that, each row and column of the matrix contains at most one non-zero entry. We put this observation in \eqref{eq.var.pivAB.pre}, and obtain that 
			\begin{align}\label{eq.var.pivAB.decom}
				\var_{\p_c} \Ll[ \# \piv(A,B) \big\vert \fC_{13} \Rr] \leq \E_{\p_c} \Ll[ \# \piv(A,B) \big\vert \fC_{13} \Rr] + \mathbf{I} + \mathbf{II} + \mathbf{III} + \mathbf{IV} + \mathbf{V} + \mathbf{VI},
			\end{align}
			where the $6$ terms are defined respectively as 
			\begin{equation}\label{eq.var.pivAB.decom}
				\begin{split}
					\mathbf{I} &:= \sum_{ \substack{e , e' \in \overrightarrow{E}(S_2), \\ e \neq e' }} \P_{\p_c}\Ll[ e ,e' \in \piv(A,B), \, M = \begin{pmatrix}
						1 & 0 \\
						0 & 0
					\end{pmatrix} \Big\vert \fC_{13} \Rr],\\
					\mathbf{II} &:= \sum_{ \substack{e , e' \in \overrightarrow{E}(S_2), \\ e \neq e' }} \P_{\p_c}\Ll[ e ,e' \in \piv(A,B), \, M = \begin{pmatrix}
						1 & 0 \\
						0 & 1
					\end{pmatrix} \Big\vert \fC_{13} \Rr],\\
					\mathbf{III} &:= \sum_{ \substack{e , e' \in \overrightarrow{E}(S_2), \\ e \neq e' }} \P_{\p_c}\Ll[ e ,e' \in \piv(A,B), \, M = \begin{pmatrix}
						0 & 1 \\
						0 & 0
					\end{pmatrix} \Big\vert \fC_{13} \Rr],\\
					\mathbf{IV} &:= \sum_{ \substack{e , e' \in \overrightarrow{E}(S_2), \\ e \neq e' }} \P_{\p_c}\Ll[ e ,e' \in \piv(A,B), \, M = \begin{pmatrix}
						0 & 1 \\
						1 & 0
					\end{pmatrix} \Big\vert \fC_{13} \Rr],\\
					\mathbf{V} &:= \sum_{ \substack{e , e' \in \overrightarrow{E}(S_2), \\ e \neq e' }} \P_{\p_c}\Ll[ e ,e' \in \piv(A,B), \, M = \begin{pmatrix}
						0 & 0 \\
						1 & 0
					\end{pmatrix} \Big\vert \fC_{13} \Rr],\\
					\mathbf{VI} &:= \sum_{ \substack{e , e' \in \overrightarrow{E}(S_2), \\ e \neq e' }} \P_{\p_c}\Ll[ e ,e' \in \piv(A,B), \, M = \begin{pmatrix}
						0 & 0 \\
						0 & 1
					\end{pmatrix} \Big\vert \fC_{13} \Rr].			
				\end{split}
			\end{equation}
			In the rest of the proof, the $6$ terms will be treated one by one.
			
			\smallskip

			\textit{Step~1: $M = \begin{pmatrix}
					1 & 0 \\
					0 & 0
				\end{pmatrix}$.}
			The event in $\mathbf{I}$ has an inclusion that 
			\begin{multline*}
				\Ll\{e ,e' \in \piv(A,B), \, M = \begin{pmatrix}
					1 & 0 \\
					0 & 0
				\end{pmatrix}\Rr\} \\
				\subset \bigcup_{x \in A} \bigcup_{\substack{ y,y' \in B \\ y \neq y'}}\{x, e_{-}, e'_{-} \in \cC_{12}(e_{-})\} \circ \{ e_+ \xleftrightarrow{T_{23} \setminus \mathcal{E}(\fC_{3})} y \} \circ  \{ e'_+ \xleftrightarrow{T_{23} \setminus \mathcal{E}(\fC_{3})} y' \}.
			\end{multline*}
			Using the BK inequality,  we have an upper bound for $\mathbf{I}$ 
			\begin{align*}
				\mathbf{I} \leq \sum_{ \substack{e , e' \in \overrightarrow{E}(S_2) \\ e \neq e' }}\sum_{x \in A} \sum_{\substack{ y,y' \in B \\ y \neq y'}} \tau^{T_{12}}_{\p_c}(x, e_{-}, e'_{-}) \tau^{T_{23}}_{\p_c}(e_+,y)\tau^{T_{23}}_{\p_c}(e'_+,y').
			\end{align*}
			For convenience, we also use the notation $\tau^S_{\p_c}(x,e) := \P_{\p_c}[x \xleftrightarrow{S} e]$ for vertex-edge two-point function in the proof. The $3$-point function requires an application of the tree-graph inequality; see Figure~\ref{fig.diagram.I_II} for an illustration. 
			\begin{align*}
				\mathbf{I} &\lesssim \sum_{ \substack{e , e' \in \overrightarrow{E}(S_2) \\ e \neq e' }}\sum_{x \in A} \sum_{\substack{ y,y' \in B \\ y \neq y'}} \sum_{u \in T_{12}} \tau^{T_{12}}_{\p_c}(x, u) \tau^{T_{12}}_{\p_c}(u, e) \tau^{T_{12}}_{\p_c}(u, e')\tau^{T_{23}}_{\p_c}(e,y)\tau^{T_{23}}_{\p_c}(e',y') \\
				&\lesssim n^{-2d+4} \vert B\vert^2\sum_{x \in A}\sum_{u \in T_{12}} \tau_{\p_c}(x, u) \Ll( \sum_{e \in \overrightarrow{E}(S_2)}\tau_{\p_c}(u, e)\Rr) \Ll( \sum_{e' \in \overrightarrow{E}(S_2)}\tau_{\p_c}(u, e')\Rr)\\
				&\lesssim  n^{-2d+8}  \vert B\vert^2 \sum_{x \in A}\Ll(\sum_{u \in T_{12}} \tau_{\p_c}(x, u)\Rr)\\
				&\lesssim  n^{-2d+10}   \vert A\vert \cdot \vert B\vert^2 . 
			\end{align*}
			Because $y,y' \in B \subset \partial Q_3$ and $e , e' \in \overrightarrow{E}(S_2)$,  we have $\vert y - e_+\vert, \vert y' - e'_+\vert \gtrsim n$, which yields the factor $n^{-2d+4}$ in the second line. Afterwards, we do integration of two-point function one by one, and get a bound 
			\begin{align}\label{eq.var.pivAB.I}
				\mathbf{I}  \lesssim  n^{-2d+10} \vert A \vert \cdot \vert B \vert^2. 
			\end{align}
			
			\smallskip
			
			\textit{Step~2: $M = \begin{pmatrix}
					1 & 0 \\
					0 & 1
				\end{pmatrix}$. }
			The event in $\mathbf{II}$ has an inclusion that 
			\begin{multline*}
				\Ll\{e ,e' \in \piv(A,B), \, M = \begin{pmatrix}
					1 & 0 \\
					0 & 1
				\end{pmatrix}\Rr\} \\
				\subset \bigcup_{x \in A} \bigcup_{y \in B}\{x, e_{-}, e'_{-} \in \cC_{12}(e_{-})\} \circ \{ y, e_{+}, e'_{+} \in \cC_{23}(e_{+})\}.
			\end{multline*}
			This gives an upper bound of the term $\mathbf{II}$. See Figure~\ref{fig.diagram.I_II} the diagram using the tree-graph inequality.
			\begin{align*}
				\mathbf{II} & \leq  \sum_{ \substack{e , e' \in \overrightarrow{E}(S_2) \\ e \neq e' }}\sum_{x \in A} \sum_{y \in B} \tau^{T_{12}}_{\p_c}(x, e_{-}, e'_{-}) \tau^{T_{23}}_{\p_c}(y, e_{+}, e'_{+})\\
				&\lesssim  \sum_{ \substack{e , e' \in \overrightarrow{E}(S_2) \\ e \neq e' }}\sum_{x \in A} \sum_{y \in B} \sum_{u \in T_{12}} \sum_{v \in T_{23}} \tau_{\p_c}(x,u) \tau_{\p_c}(u,e) \tau_{\p_c}(u,e') \tau_{\p_c}(e, v) \tau_{\p_c}(e', v) \tau_{\p_c}(v,y)\\
				&\lesssim  \sum_{x \in A} \sum_{y \in B}  \Ll(\tau_{\p_c} \star (\tau_{\p_c} \star \tau_{\p_c})^2 \star \tau_{\p_c} \Rr) (x,y).
			\end{align*}
			We then insert the bound of two-point function 
			\begin{align*}
				\tau_{\p_c}(x, u) \lesssim \vert x - u\vert^{-d+2}, \qquad (\tau_{\p_c} \star \tau_{\p_c})^2(u, y) \lesssim \vert u - y\vert^{-2d+8}.
			\end{align*}
			Their convolution depends on the dimension, and gives the bound of $\mathbf{II}$
			\begin{align}\label{eq.var.pivAB.II}
				\mathbf{II}  \lesssim \vert A \vert \cdot \vert B \vert \Ll\{\begin{array}{ll}
					n^{-d+5} & d=7,\\
					n^{-d+4} \log n & d=8, \\
					n^{-d+4} & d \geq 9.
				\end{array}\Rr. 
			\end{align}

			\begin{figure}[h!]
				\centering
				\includegraphics[width=0.8\textwidth]{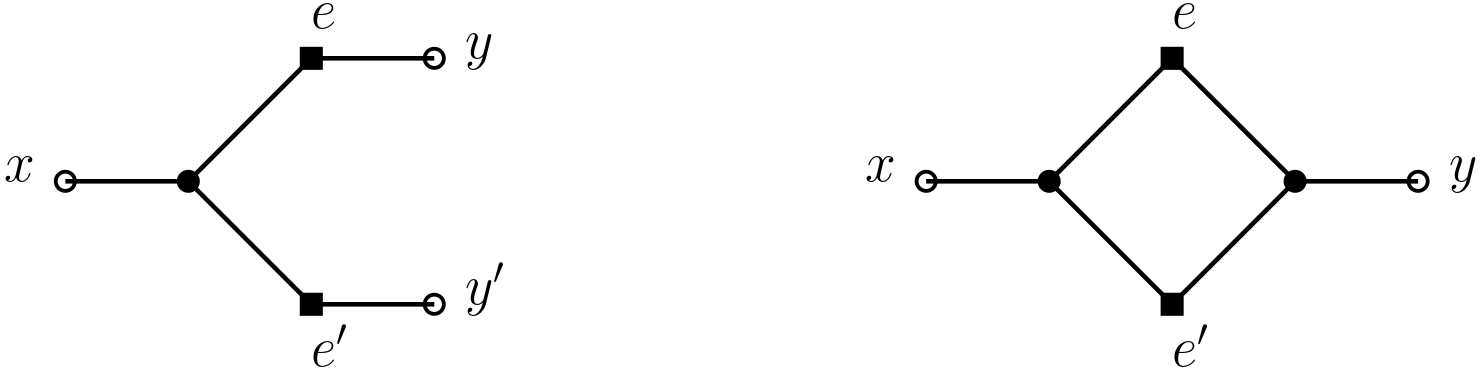}
				\caption{The diagrams for the terms $\mathbf{I}$ (left) and $\mathbf{II}$ (right).}\label{fig.diagram.I_II}
			\end{figure}
			
			\smallskip
			
			\textit{Step~3: $M = \begin{pmatrix}
					0 & 1 \\
					0 & 0
				\end{pmatrix}$. }
			The event in $\mathbf{III}$ has an inclusion that 
			\begin{multline*}
				\Ll\{e ,e' \in \piv(A,B), \, M = \begin{pmatrix}
					0 & 1 \\
					0 & 0
				\end{pmatrix}\Rr\} \\
				\subset \bigcup_{\substack{ x,x' \in A \\ x \neq x'}} \bigcup_{\substack{ y,y' \in B \\ y \neq y'}}\{x, e_{-}, y', e'_+ \text{ are in the same cluster in } T_{123}\} \circ \{ e_+ \xleftrightarrow{T_{23} \setminus \mathcal{E}(\fC_{3})} y \}  \circ \{ e'_{-} \xleftrightarrow{T_{12} \setminus \mathcal{E}(\fC_{1})} x' \}.
			\end{multline*}
			Then the BK inequality gives an upper bound of $\mathbf{III}$
			\begin{align*}
				\mathbf{III} &\leq \sum_{ \substack{e , e' \in \overrightarrow{E}(S_2) \\ e \neq e' }}\sum_{\substack{ x,x' \in A \\ x \neq x'}} \sum_{\substack{ y,y' \in B \\ y \neq y'}} \tau^{T_{123}}_{\p_c}(x, e_{-}, y', e'_+ ) \tau^{T_{23}}_{\p_c}(e_+,y) \tau^{T_{12}}_{\p_c}(x',e'_{-})\\
				&\lesssim \mathbf{III.1} + \mathbf{III.2} + \mathbf{III.3}.
			\end{align*}
			Concerning the $4$-point function, the tree-graph inequality will then yield $3$ diagrams; see Figure~\ref{fig.diagram.III} for an illustration.
			
			The term $\mathbf{III.1}$ is treated as follows.
			\begin{align*}
				\mathbf{III.1} &= \sum_{ \substack{e , e' \in \overrightarrow{E}(S_2) \\ e \neq e' }}\sum_{\substack{ x,x' \in A \\ x \neq x'}} \sum_{\substack{ y,y' \in B \\ y \neq y'}} \sum_{u,v \in T_{123}} \tau_{\p_c}(x,u) \tau_{\p_c}(u,y') \tau_{\p_c}(u,v) \tau_{\p_c}(v, e) \tau_{\p_c}(v, e') \tau_{\p_c}(e,y) \tau_{\p_c}(e',x') \\
				&\lesssim n^{-2d+4} \vert A\vert \cdot \vert B\vert \sum_{\substack{ x \in A}} \sum_{\substack{ y' \in B }} \sum_{u,v \in T_{123}} \tau_{\p_c}(x,u) \tau_{\p_c}(u,y') \tau_{\p_c}(u,v) \\
				& \qquad \times \Ll( \sum_{e \in \overrightarrow{E}(S_2)}\tau_{\p_c}(v, e)\Rr) \Ll( \sum_{e' \in \overrightarrow{E}(S_2)}\tau_{\p_c}(v, e')\Rr) \\
				&\lesssim n^{-2d+8} \vert A\vert \cdot \vert B\vert \sum_{\substack{ x \in A}} \sum_{\substack{ y' \in B }} \sum_{u \in T_{123}} \tau_{\p_c}(x,u) \tau_{\p_c}(u,y') \Ll(\sum_{v \in T_{123}}  \tau_{\p_c}(u,v)\Rr) \\
				&\lesssim   n^{-2d+10} \vert A\vert \cdot \vert B\vert \sum_{\substack{ x \in A}} \sum_{\substack{ y' \in B }} \Ll(\tau_{\p_c} \star \tau_{\p_c}\Rr)(x,y'). 
			\end{align*}
			Because $\vert y - e_+\vert, \vert x' - e'_{-}\vert \gtrsim n$, we have $\tau_{\p_c}(e,y), \tau_{\p_c}(e',x') \lesssim n^{-d+2}$ and a factor $n^{-2d+4}$ is  added  in the second line. Then we do integration of the two-point function one by one. Finally, using the estimate of convolution 
			\begin{align*}
				\Ll(\tau_{\p_c} \star \tau_{\p_c}\Rr)(x,y') \lesssim \vert x - y'\vert^{-d+4},
			\end{align*}
			we obtain a bound for the term $\mathbf{III.1}$
			\begin{align}\label{eq.var.pivAB.III.1}
				\mathbf{III.1} \lesssim n^{-3d+14} \vert A\vert^2 \cdot \vert B\vert^2.
			\end{align}
			
			The same argument works for the term $\mathbf{III.2}$
			\begin{align*}
				\mathbf{III.2} &= \sum_{ \substack{e , e' \in \overrightarrow{E}(S_2) \\ e \neq e' }}\sum_{\substack{ x,x' \in A \\ x \neq x'}} \sum_{\substack{ y,y' \in B \\ y \neq y'}} \sum_{u,v \in T_{123}} \tau_{\p_c}(x,u) \tau_{\p_c}(u,v) \tau_{\p_c}(v,y')  \tau_{\p_c}(u, e) \tau_{\p_c}(v, e') \tau_{\p_c}(e,y) \tau_{\p_c}(e',x') \\
				&\lesssim n^{-2d+4} \vert A\vert \cdot \vert B\vert \sum_{\substack{ x \in A}} \sum_{\substack{ y' \in B }} \sum_{u,v \in T_{123}} \tau_{\p_c}(x,u) \tau_{\p_c}(u,v) \tau_{\p_c}(v,y')  \\
				& \qquad \times \Ll( \sum_{e \in \overrightarrow{E}(S_2)}\tau_{\p_c}(u, e)\Rr) \Ll( \sum_{e' \in \overrightarrow{E}(S_2)}\tau_{\p_c}(v, e')\Rr) \\
				&\lesssim n^{-2d+8} \vert A\vert \cdot \vert B\vert \sum_{\substack{ x \in A}} \sum_{\substack{ y' \in B }}  ( \tau_{\p_c} \star \tau_{\p_c} \star \tau_{\p_c})(x,y').
			\end{align*}
			Using the estimate of convolution
			\begin{align*}
				( \tau_{\p_c} \star \tau_{\p_c} \star \tau_{\p_c})(x,y') \lesssim \vert x - y'\vert^{-d+6},
			\end{align*}
			\begin{align}\label{eq.var.pivAB.III.2}
				\mathbf{III.2} \lesssim n^{-3d+14} \vert A\vert^2 \cdot \vert B\vert^2.
			\end{align}
			
			The term $\mathbf{III.3}$ works exactly like the term $\mathbf{III.2}$. Combining \eqref{eq.var.pivAB.III.1} and \eqref{eq.var.pivAB.III.2}, we conclude that
			\begin{align}\label{eq.var.pivAB.III}
				\mathbf{III} \lesssim n^{-3d+14} \vert A\vert^2 \cdot \vert B\vert^2.
			\end{align}
			
			\begin{figure}[h!]
				\centering
				\includegraphics[width=0.8\textwidth]{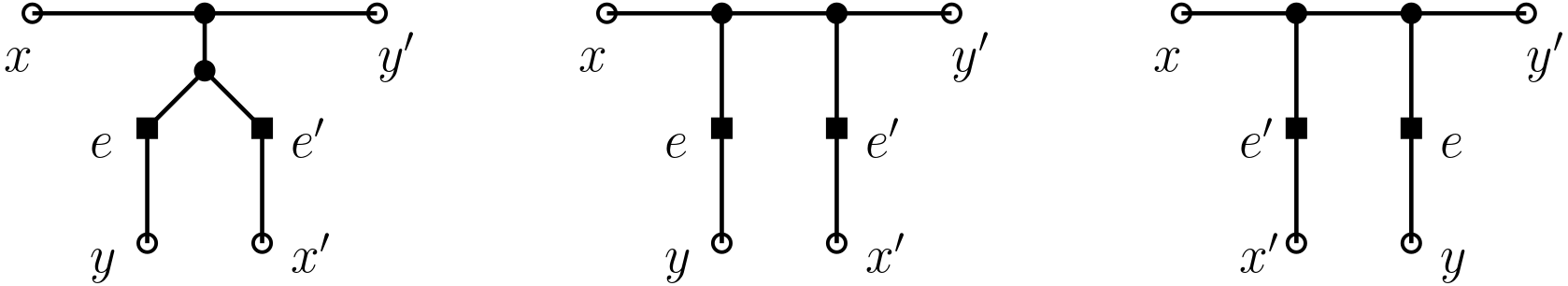}
				\caption{The diagrams for the term $\mathbf{III}$, which are indexed as $\mathbf{III.1}$, $\mathbf{III.2}$, and $\mathbf{III.3}$ from left to right.}\label{fig.diagram.III}
			\end{figure}
			
			\smallskip
			
			\textit{Step~4: $M = \begin{pmatrix}
					0 & 1 \\
					1 & 0
				\end{pmatrix}$.}
			The event in $\mathbf{IV}$ has an inclusion that 
			\begin{multline*}
				\Ll\{e ,e' \in \piv(A,B), \, M = \begin{pmatrix}
					0 & 1 \\
					1 & 0
				\end{pmatrix}\Rr\} \\
				\subset \bigcup_{\substack{ x,x' \in A \\ x \neq x'}} \bigcup_{\substack{ y,y' \in B \\ y \neq y'}}\{x, e_{-}, y', e'_+ \text{ are in the same cluster}\} \circ \{x', e'_{-}, y, e_+ \text{ are in the same cluster}\}.
			\end{multline*}
			Then the BK inequality gives an upper bound of $\mathbf{IV}$
			\begin{align*}
				\mathbf{IV} &\leq \sum_{ \substack{e , e' \in \overrightarrow{E}(S_2) \\ e \neq e' }}\sum_{\substack{ x,x' \in A \\ x \neq x'}} \sum_{\substack{ y,y' \in B \\ y \neq y'}} \tau^{T_{123}}_{\p_c}(x, e_{-}, y', e'_+ ) \tau^{T_{123}}_{\p_c}(x', e'_{-}, y, e_+ ).
			\end{align*}
			To further analyze this inequality, we need to apply the tree-graph inequality to the two terms above. Since they involve the $4$-point function,  it will generate $9$ possible diagrams. See Figure~\ref{fig.diagram.IV} for the whole classification, and we label the diagrams by $\mathbf{III.ij}$ with $i,j \in \{1,2,3\}$ as the index of the row and column in the table. 
			
			\begin{figure}[h!]
				\centering
				\includegraphics[width=0.8\textwidth]{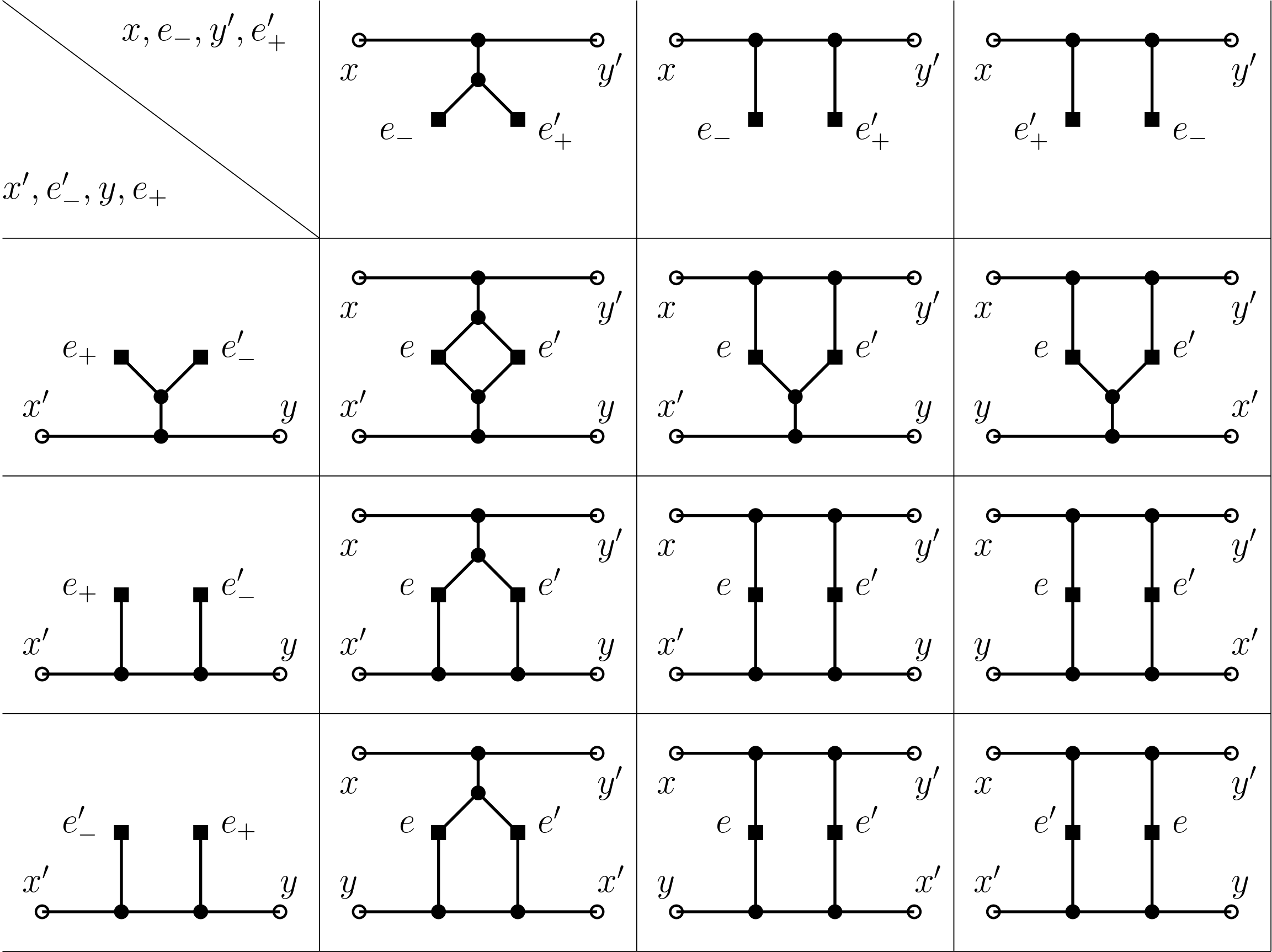}
				\caption{The diagrams for the term $\mathbf{IV}$.}\label{fig.diagram.IV}
			\end{figure}

			For the term $\mathbf{IV.11}$, it contains a subset as the term $\mathbf{II}$. For convenience, we define the kernel
			\begin{align*}
				K(u,v) := \Ll(\tau_{\p_c} \star (\tau_{\p_c} \star \tau_{\p_c})^2 \star \tau_{\p_c} \Rr)(u,v).
			\end{align*}
			Then the term $\mathbf{IV.11}$ can be written as 
			\begin{align*}
				\mathbf{IV.11} = \sum_{\substack{ x,x' \in A \\ x \neq x'}} \sum_{\substack{ y,y' \in B \\ y \neq y'}} \sum_{u,v \in T_{123}} \tau_{\p_c}(x,u) \tau_{\p_c}(u,y')  K(u,v) \tau_{\p_c}(v,x') \tau_{\p_c}(v,y).
			\end{align*}
			Because $\vert x - y'\vert \gtrsim n$, at least one between $\vert x - u\vert, \vert u - y'\vert$ should be of order $n$. Without loss of generality, we assume that $\vert u - y'\vert \gtrsim n$, then use the estimate of two-point function that 
			\begin{align*}
				\mathbf{IV.11} &\leq n^{-d+2} \vert B\vert \sum_{\substack{ x,x' \in A \\ x \neq x'}} \sum_{\substack{ y \in B}} \sum_{u,v \in T_{123}} \tau_{\p_c}(x,u)  K(u,v) \tau_{\p_c}(v,x') \tau_{\p_c}(v,y) \\
				&\lesssim n^{-d+2} \vert B\vert \sum_{\substack{ x,x' \in A \\ x \neq x'}} \sum_{\substack{ y \in B}} \sum_{v \in T_{123}} (\tau_{\p_c} \star K) (x,v) \tau_{\p_c}(v,x') \tau_{\p_c}(v,y)\\
				&\lesssim n^{-d+2} \vert B\vert \sum_{\substack{ x,x' \in A \\ x \neq x'}} \sum_{\substack{ y \in B}} \sum_{v \in T_{123}} \Ll( \frac{(\tau_{\p_c} \star K) (x,v)^p}{p}   +  \frac{\tau_{\p_c}(v,x')^{^{p'}}}{p'}  \Rr) \tau_{\p_c}(v,y). 
			\end{align*}
			The Young inequality is applied in the third line above and $p,p' \geq 1$ satisfy that $\frac{1}{p} + \frac{1}{p'} = 1$. The parameter $p$ is chosen to make the sum more balanced. As we have already estimated the bound of $K$ in \eqref{eq.var.pivAB.II}, then we have 
			\begin{align*}
				(\tau_{\p_c} \star K) (x,v) \lesssim \vert v - x\vert^{-d + \alpha}, \qquad \alpha = \Ll\{\begin{array}{ll}
					7 & d=7,\\
					6+ & d=8,\\
					6 & d\geq 9.
				\end{array}\Rr.
			\end{align*}
			Meanwhile, the decay rate of $\tau_{\p_c}(v,x')$ is of rate $\vert v - x\vert^{-d + 2}$. Thus, we solve the equation 
			\begin{align*}
				(-d + \alpha)p = (-d+2) p', \qquad \frac{1}{p} + \frac{1}{p'} = 1,
			\end{align*}
			which yields that
			\begin{align}\label{eq.choice.p}
				p = 1 + \frac{-d+2}{-d+\alpha}.
			\end{align} 
			The case $d = 7$ implies $p = + \infty$, which should be treated separately. For $d \geq 8$, we get the bound 
			\begin{align*}
				\mathbf{IV.11} &\lesssim  n^{-d+2} \vert B\vert \sum_{\substack{ x,x' \in A \\ x \neq x'}} \sum_{\substack{ y \in B}} \sum_{v \in T_{123}} \frac{(\tau_{\p_c} \star K) (x,v)^p}{p}  \tau_{\p_c}(v,y)\\
				&\lesssim  n^{-d+2} \vert A\vert \cdot \vert B\vert \sum_{\substack{ x \in A }} \sum_{\substack{ y \in B}} ((\tau_{\p_c} \star K)^p \star \tau_{\p_c})(x,y) \\
				&\lesssim \vert A\vert^2 \cdot \vert B\vert^2  \Ll\{\begin{array}{ll}
					n^{-2d + 4} \log n  & d=8,\\
					n^{-2d + 4}   & d\geq 9.
				\end{array}\Rr.
			\end{align*}	
			For the case $d = 7$, we use the fact $	(\tau_{\p_c} \star K) (x,v) \lesssim 1$ and make the calculation directly 
			\begin{align*}
				\mathbf{IV.11} &\lesssim  n^{-d+2} \vert B\vert \sum_{\substack{ x,x' \in A \\ x \neq x'}} \sum_{\substack{ y \in B}} \sum_{v \in T_{123}} \tau_{\p_c}(x', v) \tau_{\p_c}(v,y)\\
				&\lesssim  n^{-d+2} \vert A\vert \cdot \vert B\vert \sum_{\substack{ x \in A }} \sum_{\substack{ y \in B}} \sum_{v \in T_{123}} (\tau_{\p_c} \star \tau_{\p_c})(x', y) \\
				&\lesssim n^{-2d + 6} \vert A\vert^2 \cdot \vert B\vert^2. 
			\end{align*}	
			We summarize the bounds above as 
			\begin{align}\label{eq.var.pivAB.IV11}
				\mathbf{IV.11} \lesssim n^{-2d+6}   \vert A\vert^2 \cdot \vert B\vert^2.
			\end{align}

			The remaining diagrams in $\mathbf{IV}$ can be divided into two classes of equivalence $\{12, 13, 21, 31\}$ and $\{22, 23, 32, 33\}$, so it suffices to study $\mathbf{IV.12}$ and $\mathbf{IV.22}$. We state the main argument for these two terms, and leave the details as an exercise to the reader. The main observation is $x,x',y,y'$ are all of distance order $n$ from $e,e'$, so we can bound at least $4$ terms of two-point function by $n^{-4(d+2)}$. Afterwards, we treat the remaining $6$ correlations as if they were of typical distance $n$, then we would have a bound 
			\begin{align*}
				n^{-4(d+2)} \cdot n^{6 \times 2} =  n^{-4d+20}.
			\end{align*}  
			This heuristic is correct, if the remaining structure is a tree or a forest after removing $x,x',y,y'$. Otherwise, it still involves a diagram of $5$-point cycle in $\mathbf{IV.12}$ or a diagram of $6$-point cycle in $\mathbf{IV.22}$. Like the term $\mathbf{II}$, the outcome depends on dimension
			\begin{align}\label{eq.var.pivAB.IV12}
				\mathbf{IV.12} \lesssim  \vert A\vert^2 \cdot \vert B\vert^2 \Ll\{\begin{array}{ll}
					n^{-4d+20} & 7 \leq d < 10,\\
					n^{-3d+10} \log n & d=10,\\
					n^{-3d+10}  &  d>10,
				\end{array}\Rr.
			\end{align}
			and
			\begin{align}\label{eq.var.pivAB.IV22}
				\mathbf{IV.22} \lesssim  \vert A\vert^2 \cdot \vert B\vert^2 \Ll\{\begin{array}{ll}
					n^{-4d+20} & 7 \leq d < 12,\\
					n^{-3d+8} \log n & d=12,\\
					n^{-3d+8}  &  d>12.
				\end{array}\Rr.
			\end{align}

			We combine \eqref{eq.var.pivAB.IV11}, \eqref{eq.var.pivAB.IV12}, \eqref{eq.var.pivAB.IV22}, and conclude that 
			\begin{align}\label{eq.var.pivAB.IV}
				\mathbf{IV} \lesssim n^{-2d+6} \vert A\vert^2 \cdot \vert B\vert^2.
			\end{align}
			Here we unify the expression using the high-dimensional condition $d > 6$.

			\smallskip
			
			\textit{Step~5: $M = \begin{pmatrix}
					0 & 0 \\
					1 & 0
				\end{pmatrix}$.} This case is symmetry as the case $\mathbf{III}$, and it has an estimate that 
			\begin{align}\label{eq.var.pivAB.V}
				\mathbf{V} \lesssim n^{-3d+14} \vert A\vert^2 \cdot \vert B\vert^2.
			\end{align}
			
			\smallskip
			
			\textit{Step~6: $M = \begin{pmatrix}
					0 & 0 \\
					0 & 1
				\end{pmatrix}$. } This case is symmetry as the case $\mathbf{I}$, and it admits an estimate that 
			\begin{align}\label{eq.var.pivAB.VI}
				\mathbf{VI} \lesssim  n^{-2d+10} \vert A \vert^2 \cdot \vert B \vert. 
			\end{align}
			
			\smallskip
			
			Finally, a combination of \eqref{eq.var.pivAB.decom}, \eqref{eq.var.pivAB.I}, \eqref{eq.var.pivAB.II}, \eqref{eq.var.pivAB.III}, \eqref{eq.var.pivAB.IV}, \eqref{eq.var.pivAB.V}, \eqref{eq.var.pivAB.VI} yields the desired result \eqref{eq.var.pivAB}
		\end{proof}

		We are now ready to prove the main result in this section.
		
		\begin{proof}[Proof of Proposition~\ref{prop.sprinkling.AB}]
			The estimate of $M := \E_{\p_c} \Ll[ \# \piv(A,B) \big\vert \fC_{13} \Rr]$ is given in Proposition~\ref{prop.mean.pivAB}. For the variance, we the Markov inequality and the estimates \eqref{eq.mean.pivAB}, \eqref{eq.var.pivAB} to obtain that 
			\begin{align*}
				&\P_{\p_c} \Ll[ \vert \# \piv(A,B)  - M \vert \geq \frac{1}{2} M \big\vert \fC_{13} \Rr] \\
				&\leq \frac{4\var_{\p_c} \Ll[ \# \piv(A,B) \big\vert \fC_{13} \Rr] }{M^2} \\
				&\lesssim \frac{ n^{-d+5} \vert A \vert \cdot \vert B \vert  + n^{-2d+10} \Ll(\vert A \vert \cdot \vert B \vert^2 + \vert A \vert^2 \cdot \vert B \vert\Rr) 
					+ \Ll(n^{-2d+6} + n^{-3d+14}\Rr) \vert A \vert^2 \cdot \vert B \vert^2}{n^{-2d+8} \vert A\vert^2 \cdot \vert B \vert^2}.
			\end{align*}
			The crucial terms are $ n^{-d+5} \vert A \vert \cdot \vert B \vert$ and $n^{-2d+10} \Ll(\vert A \vert \cdot \vert B \vert^2 + \vert A \vert^2 \cdot \vert B \vert\Rr)$,  and we respectively need the condition $\vert A \vert \cdot \vert B \vert \geq c n^{d-2}$ and $\vert A \vert, \vert B \vert \geq c n^2$ to help get
			\begin{align}\label{eq.sprinkling.AB.Markov.1}
				\P_{\p_c} \Ll[ \vert \# \piv(A,B)  - M \vert \geq \frac{1}{2} M \big\vert \fC_{13} \Rr] \lesssim \frac{1}{c} \leq \epsilon.
			\end{align}
			Here we choose $c(\epsilon)$ large enough in function of $\epsilon$. 
			
		\end{proof}
		
		\subsection{Geometric properties of sprinkling clusters}\label{subsec.skl}
		We develop more properties about the sprinkling cluster. The first result below strengthens Proposition~\ref{prop.sprinkling.AB}. We notice that $\# \piv(A,B)$ grows linearly in function of $\vert B \vert$, so the pivotal edges should be distributed equally on $B$. This can be measured by the following quantity for $B = \bigsqcup_{i \in I} B_i$, 
		\begin{align}\label{eq.def.NAB}
			\mathcal{N}(A, (B_i)_{i \in I}) = \sum_{i \in I} \Ind{\# \pivo(A, B_i) > 0}.
		\end{align}
		The proposition states that the sprinkling clusters are distributed in several subsets among $(B_i)_{i \in I}$.
		
		\begin{prop}\label{prop.sprinkling.AB.disjoint}
			Fix $d \geq 9$. Given $\epsilon > 0$, there exists $\mathbf{c}_{\text{\ref{eq.sprinkling.AB.condition.disjoint}}}(\epsilon)$ to make the following statement hold:  for every $C_{\text{\ref{eq.sprinkling.AB.conditionB}}} \in \R_+$ and $m_{\text{\ref{eq.sprinkling.AB.disjoint}}} \in \N_+$,  there exists a parameter $\lambda(\epsilon, C_{\text{\ref{eq.sprinkling.AB.conditionB}}}, m_{\text{\ref{eq.sprinkling.AB.disjoint}}})$,  such that under $\P_{\p}[\, \cdot \, \vert \fC_{13}]$ with $\p = \p_c + \lambda n^{-2}$, for all large integer $n$ and all regular sets $A \subset \partial_{1}^\circ Q_1 \cap \fC_{1}$ and $B = \bigsqcup_{i \in I} B_i \subset \partial_{-1}^\circ Q_3 \cap \fC_{3}$ satisfying 
			\begin{align}\label{eq.sprinkling.AB.condition.disjoint}
				\vert A \vert, \vert B \vert \geq \mathbf{c}_{\text{\ref{eq.sprinkling.AB.condition.disjoint}}} n^2, \qquad  \vert A \vert \cdot \vert B \vert \geq \mathbf{c}_{\text{\ref{eq.sprinkling.AB.condition.disjoint}}} n^{d-2} , 
			\end{align}
			and 
			\begin{align}\label{eq.sprinkling.AB.conditionB}
				\max_{i \in I}\vert B_i \vert \leq C_{\text{\ref{eq.sprinkling.AB.conditionB}}} n^2,
			\end{align}
			and
			\begin{align}\label{eq.sprinkling.AB.conditionA}
				\vert A\vert \lesssim o(n^{d-6}),
			\end{align}
			we have 
			\begin{align}\label{eq.sprinkling.AB.disjoint}
				\P_{\p} \Ll[ \mathcal{N}(A, (B_i)_{i \in I}) > m_{\text{\ref{eq.sprinkling.AB.disjoint}}} \, \Big\vert \fC_{13} \Rr] \geq 1-\epsilon.
			\end{align}
		\end{prop}
		\begin{proof}
			As \eqref{eq.def.NAB} considers the sprinkling in $\P_{\p}$, we study the version under $\P_{\p_c}$ at first. Let us denote by $M_q$ the moment 
			\begin{align*}
				M_q := \sum_{i \in I} \# \piv(A, B_i)^q.
			\end{align*}
			Using the Cauchy--Schwarz inequality, we have 
			\begin{align}\label{eq.sprinkling.AB.disjoint.CS}
				(M_1)^2 \leq M_0 M_2 \Longrightarrow M_0 \geq \frac{(M_1)^2}{M_2}.
			\end{align}
			This suggests us to study the lower bound of $M_1$ and the upper bound of $M_2$.
			
			Concerning $M_1$, because $B = \bigsqcup_{i \in I} B_i$, we have 
			\begin{align*}
				M_1 = \sum_{i \in I} \# \piv(A, B_i) \geq \# \piv(A, B).
			\end{align*}
			$\mathbf{c}_{\text{\ref{eq.sprinkling.AB.condition.disjoint}}}(\epsilon)$ is the first constant to be determined. We choose $\mathbf{c}_{\text{\ref{eq.sprinkling.AB.condition.disjoint}}}(\epsilon) = c_{\text{\ref{eq.sprinkling.AB.condition}}}(\frac{\epsilon}{4})$ and apply Proposition~\ref{prop.sprinkling.AB}
			\begin{align}\label{eq.moment.M1}
				\P_{\p_c} \Ll[ M_1 > C_1 n^{-4+d} \vert A \vert \cdot \vert B \vert \, \Big\vert \fC_{13} \Rr] \geq 1-\frac{\epsilon}{4}.
			\end{align}
			
			For the moment $M_2$, its upper bound can be calculated by the mean-variance decomposition
			\begin{align*}
				& \E_{\p_c} \Ll[ M_2  \, \Big\vert \fC_{13} \Rr] \\
				&= \sum_{i \in I}\E_{\p_c} \Ll[ \# \piv(A, B_i)^2  \, \Big\vert \fC_{13} \Rr]  \\
				&=  \sum_{i \in I} \Ll(\E_{\p_c} \Ll[ \# \piv(A, B_i) \, \Big\vert \fC_{13} \Rr]^2 + \var_{\p_c} \Ll[ \# \piv(A, B_i) \, \Big\vert \fC_{13} \Rr] \Rr) \\
				&\lesssim \sum_{i \in I} \Ll( n^{-d+4} \vert A \vert \cdot \vert B_i \vert  + n^{-2d+10} \Ll(\vert A \vert \cdot \vert B_i \vert^2 + \vert A \vert^2 \cdot \vert B_i \vert\Rr) + \Ll(n^{-2d+8} + n^{-3d+14}\Rr) \vert A \vert^2 \cdot \vert B_i \vert^2 \Rr).
			\end{align*}
			In the forth line, we apply \eqref{eq.mean.pivAB} and \eqref{eq.var.pivAB} for $d \geq 9$, and group the terms. We then further simplify the last equation using \eqref{eq.sprinkling.AB.conditionB}
			\begin{align*}
				\E_{\p_c} \Ll[ M_2  \, \Big\vert \fC_{13} \Rr] \lesssim \Ll(n^{-d+4}  + n^{-2d+12} \Rr)\vert A \vert \cdot \vert B \vert + \Ll(n^{-2d+10}  + n^{-3d+16}\Rr)\vert A \vert^2 \cdot \vert B \vert.
			\end{align*}
			By Markov inequality and the fact $d \geq 9$, we have 
			\begin{align*}
				\P_{\p_c} \Ll[ M_2 > C_2 n^{-4+d} \vert A \vert \cdot \vert B \vert \, \Big\vert \fC_{13} \Rr] \lesssim \frac{1}{C_2} + \frac{1}{C_2}\Ll(n^{-d+6}  + n^{-2d+12}\Rr)\vert A \vert.
			\end{align*}
			Notice the condition \eqref{eq.sprinkling.AB.conditionA}, we thus choose a large constant $C_2$ to conclude that 
			\begin{align}\label{eq.moment.M2}
				\P_{\p_c} \Ll[ M_2 > C_2 n^{-4+d} \vert A \vert \cdot \vert B \vert \, \Big\vert \fC_{13} \Rr] \leq \frac{\epsilon}{4}.
			\end{align}
			
			In summary,  on the event 
			\begin{align*}
				E := \Ll\{ M_1 > C_1 n^{-4+d} \vert A \vert \cdot \vert B \vert \Rr\} \cap \Ll\{M_2 \leq  C_2 n^{-4+d} \vert A \vert \cdot \vert B \vert \Rr\},
			\end{align*}
			\eqref{eq.sprinkling.AB.disjoint.CS} and \eqref{eq.sprinkling.AB.condition.disjoint} implies  
			\begin{align*}
				M_0 \geq \frac{(M_1)^2}{M_2} \geq (C_1)^2 (C_2)^{-1} n^{-4 + d} \vert A \vert \cdot \vert B \vert \geq \mathbf{c}_{\text{\ref{eq.sprinkling.AB.condition.disjoint}}} (C_1)^2 (C_2)^{-1}  n^2.
			\end{align*}
			Therefore, we have $E \subset E'$
			\begin{align}\label{eq.sprinkling.AB.Eprim}
				E' := \Ll\{M_0  \geq  \mathbf{c}_{\text{\ref{eq.sprinkling.AB.condition.disjoint}}} (C_1)^2 (C_2)^{-1}  n^2 \Rr\}.
			\end{align}
			Furthermore, combining the two results \eqref{eq.moment.M1}  and \eqref{eq.moment.M2} above, we get that 
			\begin{align}\label{eq.moment.M0}
				\P_{\p_c} \Ll[  M_0 \geq  \mathbf{c}_{\text{\ref{eq.sprinkling.AB.condition.disjoint}}} (C_1)^2 (C_2)^{-1}  n^2 \, \Big\vert \fC_{13} \Rr] \geq \P_{\p_c} \Ll[  E \, \Big\vert \fC_{13} \Rr] \geq 1 - \frac{\epsilon}{2}.
			\end{align}
			In summary, the choice of the parameter $\mathbf{c}_{\text{\ref{eq.sprinkling.AB.condition.disjoint}}}(\epsilon)$ ensures $M_0$ of order $n^2$.

			Finally, we need to choose a second parameter $\lambda(\epsilon, C_{\text{\ref{eq.sprinkling.AB.conditionB}}}, m_{\text{\ref{eq.sprinkling.AB.disjoint}}})$ to further establish the statement. Conditioned on the $\p_c$-percolation and on the event $E'$, all the pivotal bonds in $\piv(A,B)$ independently have probability $\frac{\p - \p_c}{\p_c}$ to open. For each $B_i$ such that $\piv(A, B_i) \neq \emptyset$, we pick one pivotal bond. Then by FKG inequality and Poisson approximation, we have
			\begin{align*}
				\mathcal{N}(A, (B_i)_{i \in I})  \stackrel{stoc}{\succcurlyeq} \textrm{Poisson}(\lambda \mathbf{c}_{\text{\ref{eq.sprinkling.AB.condition.disjoint}}} (C_1)^2 (C_2 \p_c)^{-1}).
			\end{align*}
			Here ``$\stackrel{stoc}{\succcurlyeq}$" means stochastic domination, and the parameter is obtained from \eqref{eq.sprinkling.AB.Eprim}  and  $\p - \p_c = \lambda n^{-2}$
			\begin{align*}
				\mathbf{c}_{\text{\ref{eq.sprinkling.AB.condition.disjoint}}} (C_1)^2 (C_2)^{-1}  n^2 \cdot \Ll(\frac{\p - \p_c}{\p_c}\Rr) = \mathbf{c}_{\text{\ref{eq.sprinkling.AB.condition.disjoint}}} (C_1)^2 (C_2 \p_c)^{-1}  n^2 \cdot \lambda n^{-2} \p_c^{-1} = \lambda \mathbf{c}_{\text{\ref{eq.sprinkling.AB.condition.disjoint}}} (C_1)^2 (C_2 \p_c)^{-1}.
			\end{align*}
			Using the concentration of Poisson distribution, it suffices to pick $\lambda$ of value satisfying ${\lambda \mathbf{c}_{\text{\ref{eq.sprinkling.AB.condition.disjoint}}} (C_1)^2 (C_2 \p_c)^{-1} \gg m_{\text{\ref{eq.sprinkling.AB.disjoint}}}}$ to conclude the desired result \eqref{eq.sprinkling.AB.disjoint}.
		\end{proof}
		\begin{remark}\label{rmk.c_epsilon}
			The relation $\mathbf{c}_{\text{\ref{eq.sprinkling.AB.condition.disjoint}}}(\epsilon) \asymp \epsilon^{-1}$ can be tracked from \eqref{eq.sprinkling.AB.Markov.1} and  \eqref{eq.moment.M1}.
		\end{remark}

		We are interested in more geometric properties about the sprinkling cluster defined in \eqref{eq.def.skl_AB}. As a counterpart of \eqref{eq.def.rin}, we define the intrinsic diameter of $\cC_{s}(A,B)$
		\begin{align}\label{eq.def.rin_skl_AB}
			\rin(\cC_{s}(A,B)) := 1+ \max_{e \in \pivo(A,B)} \Ll(\max_{x \in A \cap \cC_{12}(e_{-})} \din(e_{-}, x) +    \max_{y \in B\cap \cC_{23}(e_{+})} \din(e_{+}, y)\Rr).
		\end{align}
		Here $\din$ is the intrinsic distance on $\p_c$-percolation. Therefore, $\rin(\cC_{s}(A,B))$ is actually the longest path to connect $A$ and $B$ via $\cC_{s}(A,B)$. The following lemma states that, the typical intrinsic diameter of $\cC_{s}(A,B)$ is of order $n^2$. For convenience, we define the event in Proposition~\ref{prop.sprinkling.AB.disjoint} as 
		\begin{align}\label{eq.def.Event.m}
			E_m := \{ \mathcal{N}(A, (B_i)_{i \in I}) > m \}.
		\end{align}
		
		\begin{lemma}\label{lem.sprinkling.AB.geo}
			Under the same setting as Proposition~\ref{prop.sprinkling.AB.disjoint} and further assume $\vert A \vert \cdot \vert B \vert \lesssim n^{d-2}$, then there exists a finite positive constant $K(\epsilon)$, such that 
			\begin{align}\label{eq.sprinkling.AB.geo1}
				\P_{\p} \Ll[ \rin(\cC_{s}(A,B)) \leq K n^2 \Big\vert \fC_{13},  E_m \Rr] \geq 1-\epsilon.
			\end{align}
		\end{lemma} 
		\begin{proof}
			We use the first moment argument to prove \eqref{eq.sprinkling.AB.geo1}
			\begin{align*}
				\E_{\p} \Ll[ \rin(\cC_{s}(A,B)) \leq K n^2 \Big\vert \fC_{13},  E_m \Rr] =  \frac{	\E_{\p} \Ll[ \rin(\cC_{s}(A,B))  \Ind{E_m} \Big\vert \fC_{13} \Rr]}{ \P_{\p_c}[E_m \vert \fC_{13} ]}.
			\end{align*}
			The denominator has a lower bound $\P_{\p_c}[E_m \vert \fC_{13} ] \geq 1 - \epsilon$ by \eqref{eq.sprinkling.AB}, so it suffices to get an upper bound of the nominator.
			\begin{align*}
				&\E_{\p} \Ll[ \rin(\cC_{s}(A,B))  \Ind{E_m} \Big\vert \fC_{13} \Rr] \\
				& \leq 1 + \E_{\p} \Ll[  \sum_{e \in \pivo(A,B)} \Ll(\max_{x \in A \cap \cC_{12}(e_{-})} \din(e_{-}, x) +    \max_{y \in B\cap \cC_{23}(e_{+})} \din(e_{+}, y)\Rr)  \Big\vert \fC_{13} \Rr]\\
				& \leq 1 + \sum_{e \in  \overrightarrow{E}(S_2) }\sum_{x \in A} \sum_{y \in B}   (\p - \p_c)\E_{\p_c} \Ll[  \Ll(\din(e_{-}, x) +  \din(e_{+}, y)\Rr) \Ind{\{x \xleftrightarrow{T_{12} \setminus \mathcal{E}(\fC_{1})} e_{-} \} \circ \{y \xleftrightarrow{T_{23} \setminus \mathcal{E}(\fC_{3})} e_+ \} } \Big\vert \fC_{13} \Rr] .
			\end{align*}
			In the third line, we use the fact that the pivotal edge independently has probability $(\p-\p_c)$ to open. The maximum is also relaxed to the sum there. We just treat one term 
			\begin{align*}
				&\E_{\p_c} \Ll[  \din(e_{-}, x)  \Ind{\{x \xleftrightarrow{T_{12} \setminus \mathcal{E}(\fC_{1})} e_{-} \} \circ \{y \xleftrightarrow{T_{23} \setminus \mathcal{E}(\fC_{3})} e_+ \} } \Big\vert \fC_{13} \Rr]  \\
				&= \E_{\p_c} \Ll[   \Ll(\din(e_{-}, x) \Ind{x \xleftrightarrow{T_{12} \setminus \mathcal{E}(\fC_{1})} e_{-} }\Rr) \circ \Ind{y \xleftrightarrow{T_{23} \setminus \mathcal{E}(\fC_{3})} e_+ } \Big\vert \fC_{13} \Rr]\\
				&\leq  \E_{\p_c} \Ll[ \din(e_{-}, x) \Ind{x \xleftrightarrow{T_{12} \setminus \mathcal{E}(\fC_{1})} e_{-} } \Big\vert \fC_{13} \Rr]  \P_{\p} \Ll[ y \xleftrightarrow{T_{23} \setminus \mathcal{E}(\fC_{3})} e_+   \Big\vert \fC_{13} \Rr]\\
				&\leq \E_{\p_c} \Ll[ \din(e_{-}, x) \vert x \leftrightarrow  e_{-}   \Rr] \P_{\p} \Ll[x \leftrightarrow  e_{-} \Rr]  \P_{\p} \Ll[y \leftrightarrow e_+  \Rr]\\
				&\lesssim n^2 \cdot n^{-d+2} \cdot n^{-d+2} \lesssim n^{-2d+6}. 
			\end{align*}
			The passage from the second line to the third line is an application of BK inequality, as $\din(e_{-}, x)$ is measurable by $\cC_{12}(e_{-})$, which is disjoint with $\cC_{23}(e_{+})$. Then we relax the domain. The conditional expectation of chemical distance $\E_{\p_c} \Ll[ \din(e_{-}, x) \vert x \leftrightarrow  e_{-}   \Rr] \lesssim n^2$ comes from \eqref{eq.chemical_point_point}. Therefore, we get an upper bound that 
			\begin{align*}
				\E_{\p} \Ll[ \rin(\cC_{s}(A,B))  \Ind{E_m} \Big\vert \fC_{13} \Rr] &\leq  1 + \sum_{e \in  \overrightarrow{E}(S_2) }\sum_{x \in A} \sum_{y \in B}   (\p - \p_c) n^{-2d+6}\\
				&\lesssim 1 + (\p - \p_c) n^{-d+6} \vert A \vert \cdot \vert B \vert \\
				&\lesssim n^2.
			\end{align*}
			In the last line, we use the condition $\p - \p_c = \lambda n^{-2}$ and $ \vert A \vert \cdot \vert B \vert \lesssim n^{d-2}$ in \eqref{eq.sprinkling.AB.condition}. Afterwards, a Markov inequality with a choice of large $K$ concludes \eqref{eq.sprinkling.AB.geo1}.
			
		\end{proof}

		\subsection{Coalescence between sprinkling clusters}\label{subsec.coal}
		All the results above aim to say that we have a lot of sprinkling clusters. Meanwhile, as the sprinkling clusters are very sparse in high-dimensional space, there should remain enough room for more sprinklings. This subsection is devoted to the study of coalescence between sprinkling clusters, and will prepare several useful estimates for the later paragraph.
		
		\begin{prop}[Coalescence]\label{prop.coal}
			Fix $\p = \p_c + n^{-2}$. For every regular boundary $A, A' \subset \partial_{1}^\circ Q_1 \cap \fC_{1}$ and $B \subset \partial_{-1}^\circ Q_3 \cap \fC_{3}$ satisfying 
			\begin{align*}
				A \cap A' = \emptyset.
			\end{align*}
			Define the number of coalescence $\mathcal{N}(A, A';B)$ as 
			\begin{align}\label{eq.def.AAB}
				\mathcal{N}(A, A';B) := \sum_{(x,x') \in A \times A'} \Ind{\cC_{s}(x,B) \cap \cC_{s}(x',B) \neq \emptyset}.
			\end{align}
			It satisfies an estimate that 
			\begin{align}\label{eq.AAB}
				\E_{\p} \Ll[ \mathcal{N}(A, A';B) \, \Big\vert \fC_{13}\Rr] \lesssim n^{-d+2} \min\{\vert A \vert, \vert A' \vert\} \cdot \vert B \vert +   n^{-2d+6} \vert A \vert \cdot \vert A' \vert \cdot \vert B \vert.
			\end{align}
			Moreover, if $\dex(A,A') \gtrsim n$, the estimate can be improved as 
			\begin{align}\label{eq.AAB_plus}
				\E_{\p} \Ll[ \mathcal{N}(A, A';B) \, \Big\vert \fC_{13}\Rr] \lesssim  n^{-2d+6} \vert A \vert \cdot \vert A' \vert \cdot \vert B \vert.
			\end{align}
		\end{prop}
		
		\begin{proof}
			We study \eqref{eq.AAB} restricted on different cases.
			
			\smallskip
			\textit{Case~1: case $\cC_{12}(x) \cap \cC_{12}(x') \neq \emptyset$, for some $(x,x') \in A \times A'$.} For this case, the contribution to  $\E_{\p} \Ll[ \mathcal{N}(A, A';B) \, \Big\vert \fC_{13}\Rr]$ is 
			\begin{align*}
				\mathbf{Coal.1} =  \sum_{(x,x') \in A \times A'} \P_{\p} \Ll[ \cC_{s}(x,B) \cap \cC_{s}(x',B) \neq \emptyset, \cC_{12}(x) \cap \cC_{12}(x') \neq \emptyset \, \Big\vert \fC_{13}\Rr].
			\end{align*}
			It satisfies an estimate
			\begin{align*}
				\mathbf{Coal.1} &\leq n^{-2} \sum_{x \in A, x' \in A'} \sum_{y \in B} \sum_{e \in  \overrightarrow{E}(S_2)} \P_{\p_c}\Ll[\{x \xleftrightarrow{T_{12}} x' \xleftrightarrow{T_{12}} e_{-}\} \circ \{e_+ \xleftrightarrow{T_{23}} y\} \, \Big\vert \fC_{13}\Rr]  \\
				&\lesssim n^{-2} \sum_{x \in A, x' \in A'} \sum_{y \in B} \sum_{e \in  \overrightarrow{E}(S_2)}\sum_{u \in T_{12}} \tau_{\p_c}(x,u)\tau_{\p_c}(x',u)\tau_{\p_c}(u,e)\tau_{\p_c}(e,y).
			\end{align*}
			Here the factor $n^{-2}$ is because the pivotal edge $e$ is open in $\p$-percolation. In the second line, we use BK inequality and the tree-graph inequality. Notice that $\vert y - e_+\vert \asymp n$, we continue to develop the upper bound using the two-point function
			\begin{align*}
				\mathbf{Coal.1} &\lesssim n^{-d} \vert B \vert \sum_{x \in A, x' \in A'} \sum_{e \in  \overrightarrow{E}(S_2)}\sum_{u \in T_{12}} \tau_{\p_c}(x,u)\tau_{\p_c}(x',u)\tau_{\p_c}(u,e) \\
				&\lesssim n^{-d+2} \vert B \vert \sum_{x \in A, x' \in A'} \sum_{e \in  \overrightarrow{E}(S_2)} \tau_{\p_c}(x,u)\tau_{\p_c}(x',u)\\
				&=n^{-d+2} \vert B \vert \sum_{x \in A, x' \in A'}  (\tau_{\p_c} \star \tau_{\p_c})(x,x').
			\end{align*}
			Notice the bound 
			\begin{align}\label{eq.coal.1_conv}
				(\tau_{\p_c} \star \tau_{\p_c})(x,x') \lesssim \vert x - x'\vert^{-d+4}.		 
			\end{align}
			The basic estimate \eqref{eq.regular.baby} about the regular points implies that 
			\begin{align*}
				\forall x \in A, \qquad \sum_{x' \in A'}(\tau_{\p_c} \star \tau_{\p_c})(x,x') \lesssim 1, \\
				\forall x' \in A', \qquad \sum_{x \in A}(\tau_{\p_c} \star \tau_{\p_c})(x,x') \lesssim 1.
			\end{align*}
			This concludes that 
			\begin{align}\label{eq.coal.1}
				\mathbf{Coal.1} \lesssim n^{-d+2} \min\{\vert A \vert, \vert A' \vert\} \cdot \vert B \vert.
			\end{align}
			
			\smallskip
			\textit{Case~2: $\cC_{12}(x) \cap \cC_{12}(x') = \emptyset$ for all $(x,x') \in A \times A'$, and $\cC_{12}(x') \cap \cC_{23}(e_+) \neq \emptyset$ for some $e \in \pivo(A,B)$.}  In this case, the contribution is bounded by 
			\begin{align*}
				\mathbf{Coal.2} = n^{-2}\sum_{x \in A, x' \in A'} \sum_{y \in B} \sum_{e \in  \overrightarrow{E}(S_2)} \P_{\p_c}\Ll[\{x \xleftrightarrow{T_{12}} e_{-} \} \circ \{x', e_+, y \text{ in the same cluster}\} \, \Big\vert \fC_{13}\Rr].
			\end{align*} 
			Afterwards, we analyze it using the BK inequality, the tree-graph inequality, and two-point function
			\begin{align*}
				\mathbf{Coal.2} &\leq n^{-2}\sum_{x \in A, x' \in A'} \sum_{y \in B} \sum_{e \in  \overrightarrow{E}(S_2)} \tau_{\p_c}(x,e_{-}) \tau_{\p_c}(x',e_{+}, y) \\
				&\lesssim n^{-2}\sum_{x \in A, x' \in A'} \sum_{y \in B} \sum_{e \in  \overrightarrow{E}(S_2)} \sum_{u \in T_{123}}\tau_{\p_c}(x,e) \tau_{\p_c}(e, u) \tau_{\p_c}(u,x') \tau_{\p_c}(u,y) \\
				&\lesssim n^{-d} \vert A \vert\sum_{x' \in A'} \sum_{y \in B} \sum_{e \in  \overrightarrow{E}(S_2)} \sum_{u \in T_{123}}\tau_{\p_c}(e, u) \tau_{\p_c}(u,x') \tau_{\p_c}(u,y) \\
				&\lesssim n^{-d+2} \vert A \vert\sum_{x' \in A'} \sum_{y \in B} (\tau_{\p_c}\star \tau_{\p_c})(x',y).
			\end{align*}
			In the second line, we apply the two-point function estimate to $\tau_{\p_c}(x,e)$ as $\vert x - e_{-}\vert \asymp n$. In the third line, we make an integration with respect to $e$ at first, and then use the convolution of function. Since $x'$ and $y$ have a distance of order $n$, we have
			\begin{align*}
				(\tau_{\p_c} \star \tau_{\p_c})(x',y) \lesssim n^{-d+4}.		 
			\end{align*}
			This yields that
			\begin{align}\label{eq.coal.2}
				\mathbf{Coal.2} \lesssim n^{-2d+6} \vert A \vert \cdot \vert A' \vert \cdot \vert B \vert.
			\end{align}
			We remark there is another symmetric case by exchanging the role of $x$ and $x'$ above, but the contribution is same.

			\textit{Step~3: other cases.} The contribution of $\E_{\p} \Ll[ \mathcal{N}(A, A';B) \, \Big\vert \fC_{13}\Rr]$ on the complement of Cases~1 and ~2 is  
			\begin{align*}
				\mathbf{Coal.3} = 	n^{-4}\sum_{\substack{x \in A \\ x' \in A'}}  \sum_{ \substack{e , e' \in \overrightarrow{E}(S_2), \\ e \neq e' }} \P_{\p_c}\Ll[ e \in \piv(x,B), e' \in \piv(x',B), \,  M  = \begin{pmatrix}
					0 & 0 \\
					0 & 1
				\end{pmatrix}\Big\vert \fC_{13} \Rr].
			\end{align*}
			where the matrix $M$ defined in \eqref{eq.def.Matrix}
			\begin{align*}
				M \equiv \begin{pmatrix}
					M_{11} & M_{12} \\
					M_{21} & M_{22}
				\end{pmatrix} = \begin{pmatrix}
					\Ind{\cC_{12}(e_{-}) \cap \cC_{12}(e'_{-}) \neq \emptyset} & \Ind{\cC_{12}(e_{-}) \cap \cC_{23}(e'_{+})\neq \emptyset} \\
					\Ind{\cC_{23}(e_{+}) \cap \cC_{12}(e'_{-}) \neq \emptyset} & \Ind{\cC_{23}(e_{+}) \cap \cC_{23}(e'_{+}) \neq \emptyset}
				\end{pmatrix}.
			\end{align*}
			Actually, after excluding Case~1 and Case~2, the only possible non-zero entry in $M$ is $M_{22}$, which corresponds to $\mathbf{Coal.3}$. Because $e \neq e'$ and they are sampled in independent additional randomness from $(\p - \p_c)$, there is a factor $n^{-4}$ above. Using the similar estimates as in  \eqref{eq.var.pivAB.VI}, we then conclude that 
			\begin{align}\label{eq.coal.3}
				\mathbf{Coal.3} \lesssim n^{-2d+6} \vert A \vert \cdot \vert A' \vert \cdot \vert B \vert.
			\end{align}
			Combining \eqref{eq.coal.1}, \eqref{eq.coal.2}, \eqref{eq.coal.3}, we prove the desired result \eqref{eq.AAB}.
			
			\smallskip
			The improvement in \eqref{eq.AAB_plus} just comes from  $\mathbf{Coal.1}$. If $\dex(A,A') \gtrsim n$ is satisfied, then the upper bound in \eqref{eq.coal.1_conv} is direct, which yields
			\begin{align*}
				\mathbf{Coal.1} &\lesssim n^{-d+2} \vert B \vert \sum_{x \in A, x' \in A'}  (\tau_{\p_c} \star \tau_{\p_c})(x,x') \\
				&\lesssim n^{-d+2} \vert B \vert \sum_{x \in A, x' \in A'}  n^{-d+4} \\
				&\lesssim n^{-2d+6} \vert A \vert \cdot \vert A' \vert \cdot \vert B \vert.
			\end{align*}
		\end{proof}
		\begin{remark}
			As a consequence, we have
			\begin{multline}\label{eq.proba.coal}
				\P_{\p} \Ll[ \cC_{s}(A,B) \cap \cC_{s}(A',B) = \emptyset \, \Big\vert \fC_{13}\Rr] \\
				\lesssim n^{-d+2} \min\{\vert A \vert, \vert A' \vert\} \cdot \vert B \vert +   n^{-2d+6} \vert A \vert \cdot \vert A' \vert \cdot \vert B \vert.
			\end{multline}
			There will be variants of \eqref{eq.proba.coal} for the probability of coalescence between sprinkling clusters in different cubes. However, if the coalescence occurs in different cubes, it may only cover a subset of $\mathbf{Coal.1}, \mathbf{Coal.2}, \mathbf{Coal.3}$.
		\end{remark}

		\section{Lower bound IV: Large scale resistances}\label{sec.lower.large}
		We aim to prove the general lower bound in this section. 
		\begin{prop}\label{prop.lower}
			In high-dimensional percolation, for every $\p \in (\p_c, 1]$, we have
			\begin{align}\label{eq.lower}
				\ab(\p) \gtrsim   (\p-\p_c)^3.  
			\end{align}
		\end{prop}
		
		The main task is to estimate the effective resistance in large scale, and this section is organized as follows.
		\begin{itemize}
			\item Section~\ref{subsec.flow} recalls the main tool of resistance, and reformulates it as \emph{random flow}. A sufficient condition for Proposition~\ref{prop.lower} is also given, with two main challenges pointed out: 
			\begin{enumerate}
				\item One needs to construct the crossing path in large scale, and to understand its dependence in function of $(\p-\p_c)$.
				\item The flow should be distributed on the crossing path, with very small probability of intersection.
			\end{enumerate}
			\item Section~\ref{subsec.renormalization.critical} treats the first challenge. We combine the spanning clusters and sprinkling clusters to build the giant clusters in large scale. The main argument relies on an exploration process, and the outcome will yield a structure closer to the supercritical percolation.
			\item Section~\ref{subsec.EIT} handles the second challenge in macroscopic scale. We adapt the idea of Benjamini, Pemantle, and Peres in \cite{BPP98}, who defined a measure on oriented path with the \emph{exponential tail intersection} property. The random flow is also constructed in this part.
			\item Section~\ref{subsec.intersection} further analyzes the second challenge in microscopic scale. We will prove a mixing property: even two independent random flows may visit  the same cube of correlation length $\ell_{*}$, the probability of the real intersection is of order $\ell_{*}^{6-d}$.  
			\item Finally, we give the estimate of resistance in Section~\ref{subsec.R.estimate}. 
		\end{itemize}

		\smallskip
		The idea of coarse-graining is the heart of this section, and we will switch the viewpoint below and above the correlation length from time to time. The following notations will be utilized throughout the section. Recall the correlation length that ${\ell_*:= \lfloor (\p - \p_c)^{-1/2} \rfloor}$ and we set
		\begin{align}\label{eq.ell0}
			\ell_{0} := n_0 \ell_* = n_0 \lfloor (\p - \p_c)^{-1/2} \rfloor,
		\end{align} 
		as another basic scale of correlation length. The exact value of $n_0 \in \N_+$ will be determined later, but its role is to make the percolation slightly more supercritical when $n_0$ is large, because 
		\begin{align*}
			\frac{\p - \p_c}{\ell_{0}^{-2}} = \frac{\p - \p_c}{( n_0 \lfloor (\p - \p_c)^{-1/2} \rfloor )^{-2}} \simeq n_0^2.
		\end{align*}
		For every $z \in \ell_{0} \Zd$, we call $\Lambda_{\ell_{0}}(z)$  a \emph{coarsened site}  or \emph{coarsened cube}. Two cubes $\Lambda_{\ell_{0}}(z),\Lambda_{\ell_{0}}(z')$ with $z,z' \in \ell_{0} \Zd$ are neighbors if $\vert z - z'\vert = \ell_{0}$, and we denote it by $z \stackrel{\ell_{0}}{\sim} z'$. We then define the edge set associated to the relation $\stackrel{\ell_{0}}{\sim}$
		\begin{align}\label{eq.def.E0}
			E_0 := \{\{z,z'\}: z,z' \in \ell_{0} \Zd, z \stackrel{\ell_{0}}{\sim} z'\}.
		\end{align}
		As we study the resistance in $\Lambda_L$, we also define $M \in \N_+$ to be the number of coarsened sites in each ordinate 
		\begin{align}\label{eq.def.M}
			M := \frac{L}{\ell_{0}}.
		\end{align}
		The coarsened site inside is then defined by 
		\begin{align}\label{eq.def.seeds}
			Z_{M, \ell_{0}} := \ell_{0} \Zd  \cap \Lambda_L.
		\end{align}
		$Z_{M, \ell_{0}}$ thus is a coarse-grained version of $\Lambda_L$. 	As we switch the viewpoint from time to time between $Z_{M, \ell_{0}}$ and $\Lambda_L$, for a self-avoided path, we call it simply \emph{path} in $\Lambda_L$, and call it \emph{route} when it lives in $Z_{M, \ell_{0}}$.
		
		As there will be more randomness in this section besides the percolation, we write $\PP$ for the probability of the larger space and $\Eb$ the associated expectation.
		
		\begin{figure}[t!]
			\centering
			\includegraphics[width=0.6\textwidth]{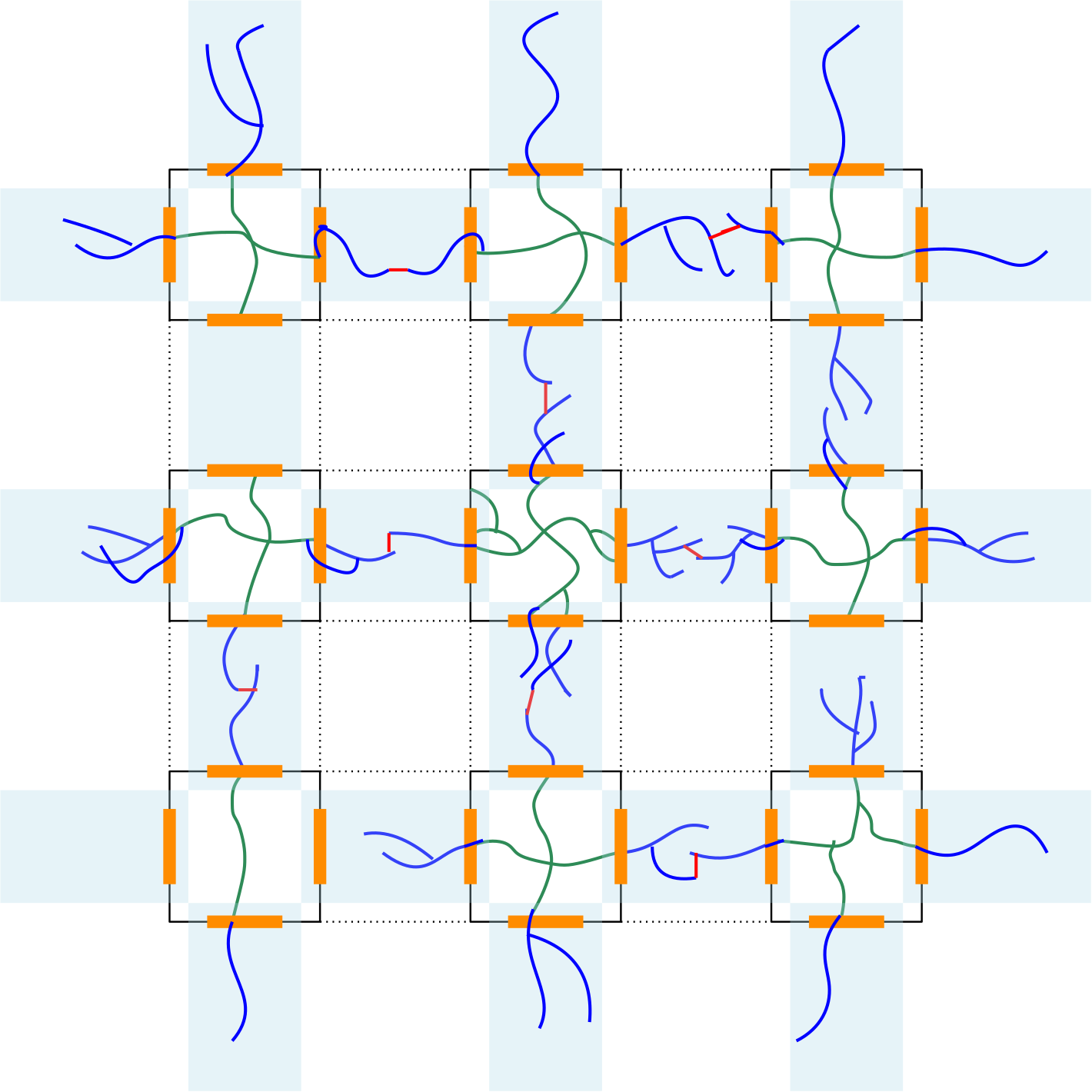}
			\caption{An illustration of the critical renormalization. The color green refers to the spanning clusters, and blue represents the sprinkling clusters at $\p_c$. The red bonds are the open pivotal bonds for sprinkling clusters at $\p$.}\label{fig.span_skl}
		\end{figure}

		\subsection{Random flow}\label{subsec.flow}
		Recall the energy of flow defined in \eqref{eq.def.energy.flow} and write its version in $\Lambda_L$ as 
		\begin{align}\label{eq.def.energy.flow.cube}
			\mathscr{E}_*(\Lambda_L, \a,\theta) := \sum_{e \in \Ed(\Lambda_L)} \a^{-1}(e) \theta^2(e).
		\end{align}
		The main task throughout the section is to verify the two uniform bounds \eqref{eq.lower.Alower} and \eqref{eq.lower.Rupper} below.

		\begin{prop}\label{prop.lower.reduction}
			In high-dimensional percolation, given $\p \in (\p_c, 1]$, for every $L \in \N_+$, there exists a measurable event $A_L$ and a unit flow $\theta$ from $[\partial \Lambda_{L}]^{-}$ to $[\partial \Lambda_{L}]^{+}$, such that the following estimates hold 
			\begin{align}\label{eq.lower.Alower}
				\liminf_{L \to \infty}  \P_{\p}\Ll[A_L\Rr] \geq c_A > 0,
			\end{align}
			and
			\begin{align}\label{eq.lower.Rupper}
				\limsup_{L \to \infty}\E_{\p}\Ll[\mathscr{E}_*(\Lambda_L, \a,\theta) \mathbf{1}_{A_L}\Rr] L^{d-2} \lesssim \ell_*(\p)^6.
			\end{align}
		\end{prop}
		
		This proposition is sufficient for the lower bound estimate of the conductivity.
		\begin{proof}[Proof of Proposition~\ref{prop.lower} assuming Proposition~\ref{prop.lower.reduction}]
			This is just an application of \eqref{eq.resistance_key_mean} 
			\begin{align}\label{eq.resistance_key_mean_mega}
				\ab(\Lambda_{L}, \p) \geq \frac{L^{2-d} \P_{\p}[A_L]^2}{\E_{\p}\Ll[\mathscr{E}_*(\Lambda_L, \a,\theta) \mathbf{1}_{A_L}\Rr]}.
			\end{align}
			We combine it with the uniform bounds \eqref{eq.lower.Alower}, \eqref{eq.lower.Rupper}, and deduce the following estimate for $L$ large enough
			\begin{align}\label{eq.lowerObject}
				\ab(\Lambda_{L}, \p) \geq \frac{L^{2-d} \P_{\p}[A_L]^2}{\E_{\p}\Ll[\mathscr{E}_*(\Lambda_L, \a,\theta) \mathbf{1}_{A_L}\Rr]} \gtrsim \frac{L^{2-d} c_A}{L^{2-d} \ell_*(\p)^6} \gtrsim \ell_*(\p)^{-6} \gtrsim (\p-\p_c)^3.
			\end{align} 
			We then pass $L \to \infty$, then \eqref{eq.def_sigma} and \eqref{eq.chain.low} conclude the desired bound in \eqref{eq.lower}
			\begin{align*}
				\ab(\p) = \lim_{L \to \infty} \Ab_{11}(\Lambda_{L}, \p) \geq \lim_{L \to \infty} \ab(\Lambda_{L}, \p) \gtrsim (\p-\p_c)^3.
			\end{align*}
		\end{proof}

		As we mentioned at the beginning of this section, there are two main challenges. The first one is easy to understand, as we need at first enough crossing paths. This is also the focus of Sections~\ref{sec.regular} and ~\ref{sec.skl}. We now explain the second challenge. This concerns the construction of the unit flow. A useful argument is \emph{random flow} argument: sample a random path $\gamma \in S$ among all the possible self-avoiding crossing paths between $[\partial \Lambda_{L}]^{-}$ to $[\partial \Lambda_{L}]^{+}$, which gives a random flow
		\begin{align*}
			\Ind{e \in \gamma}.
		\end{align*} 
		Then we propose a unit flow using its average
		\begin{align}\label{eq.flow_by_path}
			\theta(e) := \Eb\Ll[\Ind{e \in \gamma} \vert \a\Rr]. 
		\end{align}
		Here $\Eb$ stands for the expectation under a larger space containing both the percolation and the random path $\gamma$. The relation between the intersection and resistance is stated in the following lemma.

		\begin{lemma}\label{lem.flow}
			The flow $\theta$ defined in \eqref{eq.flow_by_path} satisfies the conditions in \eqref{eq.def.flow}. Moreover, let $\gamma, \gamma'$ be sampled as two independent random flows conditioned on $\a$, then we have
			\begin{align}\label{eq.flow_intersection}
				\mathscr{E}_*(\Lambda_L, \a,\theta) = \Eb\Ll[ \vert \gamma \cap \gamma'\vert \, \vert \a\Rr].
			\end{align}
		\end{lemma}	
		\begin{proof}
			Every random flow $\Ind{e \in \gamma}$ satisfies \eqref{eq.def.flow}, because it has one source and one sink, and every other sites are only passed once. This also implies that the average $\theta$ satisfies the same property. Concerning \eqref{eq.flow.intersection},  we have the following trick by independence
			\begin{align}\label{eq.flow.trick}
				\Ll(\Eb\Ll[\Ind{e \in \gamma} \vert \a\Rr]\Rr)^2 = \Eb\Ll[\Ind{e \in \gamma} \vert \a\Rr]\Eb\Ll[\Ind{e \in \gamma'} \vert \a\Rr] = \Eb\Ll[\Ind{e \in \gamma \cap \gamma'} \vert \a\Rr].
			\end{align}
			Then the resistance can be reformulated as
			\begin{align}\label{eq.flow.intersection}
				\begin{split}
					\mathscr{E}_*(\Lambda_L, \a,\theta) &= \sum_{e \in \Ed(\Lambda_L)} \a^{-1}(e) \theta^2(e)\\
					&= \sum_{e \in \Ed(\Lambda_L)} \a^{-1}(e) \Eb\Ll[\Ind{e \in \gamma \cap \gamma'} \vert \a\Rr] \\
					&=  \Eb\Ll[ \vert \gamma \cap \gamma'\vert \, \vert \a\Rr].
				\end{split}
			\end{align}
			In the third line, we utilize Fubini's lemma and the fact $\a(e) = 1$ for all $e \in \gamma \cap \gamma'$.
		\end{proof}
		\begin{remark}\label{rmk.heurictic.mixing}
			Inspired by Proposition~\ref{prop.lower_water}, it is not hard to establish a heuristic about the random flow for \eqref{eq.lower.Rupper}: in $[\partial \Lambda_L]^{-}$, it contains $M^{d-1}$ coarsened cube, and every cube contains $\ell_{0}^{d-6}$ spanning clusters. If a positive proportion of these spanning clusters can finally be extended to  $[\partial \Lambda_L]^{+}$, then we choose uniformly one among them as $\gamma$, and obtain
			\begin{align*}
				\Eb\Ll[ \vert \gamma \cap \gamma'\vert \, \vert \a\Rr] &\simeq \Ll(M^{d-1} \cdot \ell_{0}^{d-6}\Rr)^{-1} \cdot M \ell_{0}^2 \\
				&\simeq M^{2-d} \cdot \ell_{0}^{2-d} \cdot \ell_{0}^{6} \\
				&\simeq L^{2-d}\cdot \ell_{0}^{6}.
			\end{align*}
			Here the factor $\Ll(M^{d-1} \cdot \ell_{0}^{d-6}\Rr)^{-1}$ suggests a mean-field mixing in the intersection, and $M \ell_{0}^2 $ is the typical length of one crossing path. We need to rigorously construct the crossing paths and the property about intersection.
		\end{remark}

		\subsection{Critical renormalization}\label{subsec.renormalization.critical}
		We implement a renormalization in correlation length $\ell_0$, and it will give a clear structure about the geometry of $\p$-percolation in slightly supercritical regime.  It relies on the sprinkling properties in Section~\ref{sec.skl}, and the construction here also take the random flow in Section~\ref{subsec.flow} into account. 
		
		We generalize the $(Q,S,H,T)$-structure in \eqref{eq.def.SQH} to the coarsened sites $z \in \ell_{0} \Zd$, with a choice $n = \frac{\ell_{0}}{2}$ there
		\begin{align}\label{eq.def.SQH.ell0.QSH}
			Q(z) := \Lambda_{\frac{\ell_{0}}{2}}(z), \qquad S(z) := \Lambda_{\frac{3\ell_{0}}{8}}(z), \qquad H(z):= \Lambda_{\frac{\cc \ell_{0}}{2}}(z).
		\end{align}
		The sprinkling cluster is extended from the interior of the boundary of $Q(z)$
		\begin{align}\label{eq.def.SQH.ell0.bd.o}
			\forall i \in \{\pm 1, \cdots, \pm d\}, \qquad  \partial_i^\circ Q(z) := \Ll\{ y \in \partial_i Q(z):  (y-z)_j \in \Ll(-\frac{\ell_{0}}{8},  \frac{\ell_{0}}{8}\Rr), \forall j \neq \vert i \vert \Rr\}.
		\end{align}
		The notation $S(z,z')$ is conserved to the domain of the pivotal bonds between neighbor cubes
		\begin{align}\label{eq.def.SQH.ell0.S}
			\forall z \stackrel{\ell_{0}}{\sim} z', \qquad  S(z,z') := \Lambda_{\frac{\ell_{0}}{4}}\Ll(\frac{z+z'}{2}\Rr),
		\end{align}
		and $T(z,z')$ stands for the tube connecting the neighbors
		\begin{align}\label{eq.def.SQH.ell0.T}
			\forall z \stackrel{\ell_{0}}{\sim} z', \qquad  T(z,z') := 	\frac{z+z'}{2} + \Big( -\frac{5 \ell_0}{16}, \frac{5 \ell_0}{16} \Big) \times \Big( -\frac{3\ell_{0}}{16}, \frac{3\ell_{0}}{16} \Big)^{d-1}.
		\end{align}
		See Figure~\ref{fig.span_skl} for an illustration, and compare it with its local version in Figure~\ref{fig.skl}.
		
		We always work in the natural monotone coupling space, then the clusters in $Q(z)$ issued from $H(z)$ are defined as 
		\begin{align}\label{eq.def.fC.ell0}
			\fC(z) := \{y \in Q(z): y \xleftrightarrow{Q(z)} H(z) \text{ in } \p_c\text{-percolation}\}.
		\end{align}
		We also denote by $\mathcal{G}$ its union 
		\begin{align}\label{eq.sigma_pc}
			\mathcal{G} := \bigcup_{z \in \ell_{0} \Zd} \fC(z).
		\end{align}

		The general idea of the renormalization at correlation length can be stated as follows: we sample at first the $\p_c$-percolation $\mathcal{G}$, then for every neighbor pair $\{z,z'\} \in E_0$, we establish the sprinkling clusters as Definition~\ref{def.skl}, using the pivotal bonds in $S(z,z')$ and contained in the tube
		\begin{align}\label{eq.skl.support}
			T(z,z'). 
		\end{align}
		Notice the sprinkling cluster also depends on the subset of boundary in $\partial^{\circ} Q (z)$, and this will be the $2d$-spanning clusters introduced in Section~\ref{subsec.spanning_cluster} and their regular boundary in Section~\ref{sec.regular}.
		
		We record one important observation about the sprinkling clusters. It also explains the choice of parameters in $(Q,S,H,T)$-structure.
		\begin{lemma}\label{lem.skl.independent}
			Conditioned on $\G$, the sprinkling clusters in different neighbor pairs $(\{z,z'\})_{\{z,z'\} \in E_0}$ are independent.
		\end{lemma}
		\begin{proof}
			As mentioned in \eqref{eq.skl.support}, for every neighbor pair $\{z,z'\} \in E_0$, the sprinkling clusters connecting $Q(z)$ and $Q(z')$ is contained in $T(z,z')$. The tube $T(z,z')$ can be seen as the convex hull generated between  $S(z)$ and $S(z')$. Therefore, for different neighbor pairs $\{z_1,z_1'\}, \{z_2, z_2'\}$, the tubes are disjoint
			\begin{align*}
				T(z_1,z_1') \cap T(z_2,z_2') = \emptyset.
			\end{align*}
			See Figure~\ref{fig.span_skl} for an illustration. 	Since the percolation $\{\a(e)\}_{e \in \Ed}$ is sampled following the product Bernoulli measure, naturally $\{\a(e)\}_{e \in \Ed(T(z_1,z_1'))}$ and $\{\a(e)\}_{e \in \Ed(T(z_2,z_2'))}$ are independent. One can add more pairs and repeat this argument, then the statement is justified.
		\end{proof}

		\subsubsection{Pre-good cubes and typical spanning / sprinkling clusters}\label{subsubsec.good_skl}
		
		This paragraph summarizes several properties about spanning clusters in Section~\ref{subsec.spanning_cluster} and sprinkling clusters in Section~\ref{sec.skl}. 
		
		The \emph{typical spanning} clusters refer to those $2d$-spanning clusters with correct order of physical quantities. Recall the definition of $2d$-spanning in Definition~\ref{def.spanning}, the intrinsic diameter in \eqref{eq.def.rin}, and the regular boundary in Definition~\ref{def.regular}. We follow  \eqref{eq.def.regular_bd} and further define a more precise boundary 
		\begin{align}\label{eq.def.SQH.ell0.bd.2}
			\begin{split}
				\partial^\circ_{i,+} Q &:=\Ll\{ y \in \partial_i Q:  y_j \in \Ll(\frac{\ell_{0}}{16}, \frac{\ell_{0}}{8}\Rr), \forall j \neq \vert i \vert \Rr\}, \\
				\partial^\circ_{i,-} Q &:=\Ll\{ y \in \partial_i Q:  y_j \in \Ll(-\frac{\ell_{0}}{8},  -\frac{\ell_{0}}{16}\Rr), \forall j \neq \vert i \vert \Rr\}.
			\end{split}
		\end{align}
		This segment of boundary also generalizes the notion of regular boundary points
		\begin{align}\label{eq.def.regular_bd_2}
			\begin{split}
				\PTRE{(r,\delta)}_{i, +} K := \Ll\{ y \in \partial^\circ_{i,+} Q: y \text{ is } (r,\delta) \text{-regular} \Rr\},\\ 
				\PTRE{(r,\delta)}_{i, -} K := \Ll\{ y \in \partial^\circ_{i,-} Q: y \text{ is } (r,\delta) \text{-regular} \Rr\}. 
			\end{split}
		\end{align}

		\begin{definition}[Typical spanning clusters]\label{def.spanning.z}
			In $\p_c$-percolation of $Q(z)$, a spanning cluster $\cC \subset \fC(z)$ is called \emph{typical} if it satisfies the following properties:
			\begin{enumerate}
				\item \emph{$2d$-spanning:} $\cC$ is $2d$-spanning in $Q(z)$.
				\item \emph{Size:} $c_{\mathrm{v}} \ell_{0}^4 < \vert \cC \vert < C_{\mathrm{v}} \ell_{0}^4 $.
				\item \emph{Intrinsic diameter:} $\rin(\cC) \leq K_{\mathrm{spn}} \ell_{0}^2$.
				\item \emph{Regular boundary:}  for every $i \in \{\pm 1, \cdots, \pm d\}$ 
				\begin{align*}
					c_{\mathrm{b}} \ell_{0}^2 \leq \vert \PTRE{(r,\delta)}_{i, +} \cC  \vert ,\vert \PTRE{(r,\delta)}_{i, -} \cC  \vert \leq |\cC\cap\partial_i Q(z)|\leq C_{\mathrm{b}}\ell_0^2.
				\end{align*}
			\end{enumerate} 
			We denote by $\mathscr{S}(z)$ the set of typical spanning clusters
			\begin{equation}\label{eq.def.spanning.z}
				\mathscr{S}(z) := \Ll\{\cC: \cC \text{ is a } \text{typical spanning cluster} \text{ in } Q(z)\Rr\}.
			\end{equation}
		\end{definition}
		The item (4) above is a little more precise than \eqref{eq.def.regular_bd} discussed in the previous sections, but it follows the same philosophy. Moreover, (4) also highlights that the regular points are distributed equally on the boundary.

		A cube is pre-good if it contains enough number of typical spanning clusters.
		\begin{definition}[Pre-good cube]\label{def.pregood}
			Given $z \in \ell_{0} \Zd$, a cube $Q(z)$ is pre-good if and only if
			\begin{align}\label{eq.pregood}
				c_{\ref{eq.pregood}} \ell^{d-6}_{0} \leq \# \mathscr{S}(z) \leq C_{\ref{eq.pregood}} \ell^{d-6}_{0}.
			\end{align}
		\end{definition}

		As discussed in Section~\ref{subsec.skl}, there is a high probability to obtain pre-good cube.
		\begin{prop}\label{prop.proba.pregood}
			Fix $\epsilon_0 \in (0,1)$ satisfying $\epsilon_0 \ll 1$, there exist finite positive constants 
			\begin{align}\label{eq.vector.pregood}
				\Ll(r, \cc, \delta, c_{\mathrm{v}}, C_{\mathrm{v}}, c_{\mathrm{b}}, C_{\mathrm{b}}, K_{\mathrm{spn}}, c_{\ref{eq.pregood}}, C_{\ref{eq.pregood}} \Rr),
			\end{align}
			only depending on $\epsilon_0,d$, such that for all $z \in \ell_{0} \Zd$
			\begin{align}\label{eq.proba.pregood}
				\P_{\p_c}\Ll[Q(z) \text{ is pre-good }\Rr] \geq 1 - \frac{\epsilon_0}{4}.
			\end{align}
		\end{prop}
		\begin{proof}
			The main argument follows Proposition~\ref{prop.num_spanning_new}. Once we prove that the $2d$-spanning cluster containing $0$ satisfies the property $\mathbf{S}$ with a positive probability, then we can find parameters $(c_{\mathrm{v}}, C_{\mathrm{v}}, c_{\ref{eq.pregood}}, C_{\ref{eq.pregood}})$ to realize $\ell^{d-6}_{0}$ copies of $2d$-spanning clusters of the property (2) and $\mathbf{S}$ with probability $1 - \frac{\epsilon_0}{4}$; see \eqref{eq.num_spanning.S}. The concrete property $\mathbf{S}$ is (3) and (4) in Definition~\ref{def.pregood}, and they were justified in previous sections:
			\begin{itemize}
				\item (3) was proved in Lemma~\ref{lem.spanning.diam}, where we tune the parameter $K_{\mathrm{spn}}$.
				\item (4) was proved in Proposition~\ref{prop.regular_cube}, where the parameters $r,\delta,c_{\mathrm{b}}, C_{\mathrm{b}}$ were determined there. The parameter $\cc$ was fixed as prerequisite in Proposition~\ref{prop.choice_c}.
			\end{itemize}
			Then we conclude \eqref{eq.proba.pregood}.
		\end{proof}
		\begin{remark}\label{rmk.proba.k-nest}
			In the rest of this section, without further explanation, we always set the constants $\Ll(r, \cc, \delta, c_{\mathrm{v}}, C_{\mathrm{v}}, c_{\mathrm{b}}, C_{\mathrm{b}}, K_{\mathrm{spn}}, c_{\ref{eq.pregood}}, C_{\ref{eq.pregood}} \Rr)$ in function of $\epsilon_0$ and $d$ as in Proposition~\ref{prop.proba.pregood}. As $r,\delta$ are fixed, we just write $\PTR$ as a shorthand notation for $\PTRE{(r,\delta)}$.
		\end{remark}

		We also have a similar definition for the \emph{typical sprinkling clusters}. Following the notation $\cC_{s}(\cdot, \cdot)$ in \eqref{eq.def.skl_AB}, we define the sprinkling cluster between $A \subset \fC(z)$ and $\fC(z+e_i)$, via the $(r, \delta)$-regular boundary of the opposite faces, with the pivotal bonds in $S(z,z+e_i)$, within $T(z,z+e_i)$ 
		\begin{align}\label{eq.def.goodskl}
			\skl(A, e_i) \equiv \cC_{s}\Ll(\PTR_{i} A , \, \PTR_{-i} \bigcup_{\cC' \in \mathscr{S}(z+e_i)} \cC' \Rr).
		\end{align}
		As a convention, when using this notation, it also means the sprinkling cluster exists. We define the set of offspring clusters in the direction $e_i$ of $A$
		\begin{align}\label{eq.def.S_child}
			\mathscr{S}_{\mathrm{off}}( A, e_i) := \{\cC' \in \mathscr{S}(z+e_i) : \PTR_{i} A \xleftrightarrow{\skl(A, e_i)} \cC'\}.
		\end{align}
		We also introduce the set of coalescence in the direction $e_i$
		\begin{align}\label{eq.def.V_coal}
			\mathcal{V}_{\mathrm{coal}}(A, e_i) := \Ll\{x \in \PTR_{i} \bigcup_{\cC' \in \mathscr{S}(z)} \cC' \setminus A: \skl(x, e_i) \cap \skl(A, e_i) \neq \emptyset\Rr\}.
		\end{align}
		Recall the intrinsic diameter $\rin$ of sprinkling cluster in \eqref{eq.def.rin_skl_AB}. The typical sprinkling clusters should branch and coalesce with a reasonable number of clusters, and have correct order of intrinsic diameter.
		
		\begin{definition}[Typical sprinkling cluster]\label{def.goodskl}
			Given a pre-good cubes $Q(z)$  and $A \subset \fC(z)$. For every direction $i \in \{\pm 1, \cdots , \pm d\}$, given $Q(z+e_i)$ pre-good, then $\skl(A, e_i)$ is called typical if all the following properties hold:
			\begin{enumerate}
				\item \emph{Offspring clusters}:  $10 \leq m_{\mathrm{off}} \leq  \vert \mathscr{S}_{\mathrm{off}}(A,e_i) \vert  \leq M_{\mathrm{off}} < \infty$.
				\item \emph{Coalescence on the boundary}: $\vert \mathcal{V}_{\mathrm{coal}}(A,e_i) \vert \leq M_{\mathrm{coal}} < \infty$.
				\item \emph{Intrinsic diameter}: $\vert \rin(\skl(A, e_i)) \vert \leq K_{\mathrm{skl}} \ell_{0}^2 $.
			\end{enumerate}
		\end{definition}
		As a summary of the results from Sections~\ref{subsec.skl} and ~\ref{subsec.coal}, we claim that every sprinkling cluster $\skl(A, e_i)$ has high probability to be typical with a reasonable choice of parameters.
		\begin{corollary}\label{cor.proba.goodskl}
			Given $\epsilon_0 \in (0,1)$ satisfying $\epsilon_0 \ll 1$, then the finite constant $\mathbf{c}_{\text{\ref{eq.sprinkling.AB.condition.disjoint}}}(\epsilon_0/8)$ in Proposition~\ref{prop.sprinkling.AB.disjoint}  makes the following statement hold: for every $m_{\mathrm{off}} > 10$, there exist constants
			\begin{align}\label{eq.vector.goodskl}
				\Ll(n_0, K_{\mathrm{skl}}, M_{\mathrm{off}}, M_{\mathrm{coal}}, C_{\text{\ref{eq.proba.goodskl.condition}}}\Rr),
			\end{align}
			only depending on $\epsilon_0,d$, such that for all $z \in \ell_{0} \Zd$ and $\p \searrow \p_c$, every $i \in \{\pm 1, \cdots , \pm d\}$, and every $A \subset \fC(z)$ satisfying
			\begin{align}\label{eq.proba.goodskl.condition}
				\mathbf{c}_{\text{\ref{eq.sprinkling.AB.condition.disjoint}}} \ell_0^2 \leq \vert \PTR_i A\vert \leq C_{\text{\ref{eq.proba.goodskl.condition}}} \ell_0^2,
			\end{align}
			we have the estimate when $Q(z), Q(z+e_i)$  are pre-good
			\begin{align}\label{eq.proba.goodskl}
				\P_{\p}\Ll[\skl(A,e_i) \text{ is typical } \vert \fC(z), \fC(z+e_i)  \Rr] \geq 1 - \frac{\epsilon_0}{4}.
			\end{align}
		\end{corollary}
		\begin{proof}
			Take the direction as $e_1$ and assume $Q(z+e_1)$ pre-good. The condition (1) is the most important one, and the main argument is Proposition~\ref{prop.sprinkling.AB.disjoint}: we take the partition of $B$ as 
			\begin{align*}
				B_j = \PTR_{-1} \cC', \qquad \cC' \in \mathscr{S}(z+e_1).
			\end{align*}
			The choice of $C_{\text{\ref{eq.proba.goodskl.condition}}}$ can be any finite constant larger than $\mathbf{c}_{\text{\ref{eq.sprinkling.AB.condition.disjoint}}}$. The conditions \eqref{eq.sprinkling.AB.conditionB} and \eqref{eq.sprinkling.AB.conditionA} are naturally satisfied thanks to the property of typical spanning cluster; see (4) of Definition~\ref{def.spanning.z}.
			
			We then determine the value of $n_0$. We set 
			\begin{align*}
				m_{\text{\ref{eq.sprinkling.AB.disjoint}}} := m_{\mathrm{off}}, \qquad C_{\text{\ref{eq.sprinkling.AB.conditionB}}} := C_{\mathrm{b}},
			\end{align*}
			and a natural choice for $\lambda$ in Proposition~\ref{prop.sprinkling.AB.disjoint} should satisfy that
			\begin{align*}
				\p - \p_c = \lambda(\epsilon_0 / 8,  C_{\mathrm{b}} ,  m_{\mathrm{off}}) \ell_0^{-2}.
			\end{align*}
			Using the definition \eqref{eq.ell0}, the value of $n_0$ is determined
			\begin{align}\label{eq.n0_choice}
				n_0 = \lambda(\epsilon_0 / 8,  C_{\mathrm{b}} ,  m_{\mathrm{off}})^{\frac{1}{2}}.
			\end{align}
			Under the choice of parameters above, Proposition~\ref{prop.sprinkling.AB.disjoint} implies that 
			\begin{align*}
				\P_{\p} \Ll[ \vert \mathscr{S}_{\mathrm{off}}(A,e_i) \vert > m_{\mathrm{off}} \, \Big\vert \fC(z), \fC(z+e_1) \Rr] \geq 1-\frac{\epsilon_0}{8}.
			\end{align*}
			
			Conditions (2),(3) can all be verified by the first moment argument with a choice of $K_{\mathrm{skl}}, M_{\mathrm{off}}, M_{\mathrm{coal}}$ large, which was discussed in Lemma~\ref{lem.sprinkling.AB.geo}, and Proposition~\ref{prop.coal}. Here we highlight the part $M_{\mathrm{coal}}$, which is an application of Proposition~\ref{prop.coal} by choosing 
			\begin{align*}
				A' = \PTR_{i} \fC(z) \setminus A, \qquad B = \PTR_{-i} \bigcup_{\cC' \in \mathscr{S}(z+e_i)} \cC'.
			\end{align*}
			Then the estimate \eqref{eq.AAB} suggests that 
			\begin{align}\label{eq.Vcoal}
				\begin{split}
					&\E_{\p} \Ll[ \vert \mathcal{V}_{\mathrm{coal}}(A, e_i) \vert  \, \Big\vert \fC(z), \fC(z+e_1) \Rr] \\
					&\lesssim \ell_0^{-d+2} \cdot C_{\text{\ref{eq.proba.goodskl.condition}}} \ell_0^2 \cdot C_{\mathrm{b}} C_{\ref{eq.pregood}} \ell^{d-4}_{0}  +  \ell_0^{-2d+6}  C_{\text{\ref{eq.proba.goodskl.condition}}} \ell_0^2 \cdot \Ll(C_{\mathrm{b}} C_{\ref{eq.pregood}} \ell^{d-4}_{0}\Rr)^2\\
					&\lesssim 1.
				\end{split}
			\end{align}
			Therefore, we can choose a large $M_{\mathrm{coal}}$ and use Markov inequality to derive the result.
		\end{proof}

		It is quite natural to assign $A$ as a typical spanning cluster in Definition~\ref{def.spanning.z}. However, there is a competition of the minimum regular boundary between \eqref{eq.proba.goodskl.condition} and (4) of Definition~\ref{def.spanning.z}. Moreover, there is no way to cheat, as to tune the parameter $c_{\mathrm{b}}$ will also change the number of typical spanning clusters. For this reason, we have to define \emph{the $\kk$-nest} \footnote{We refer to \url{https://en.wikipedia.org/wiki/Beijing_National_Stadium} for the introduction of the Bird's Nest.} as 
		\begin{align}\label{eq.def.k-nest}
			\CC := \bigcup_{j=1}^{\kk} \cC_{j}, \qquad \cC_j \in \mathscr{S}(z).
		\end{align}
		Then, as the regular boundary is enlarged, the $\kk$-nest has the desired property as in Corollary~\ref{cor.proba.goodskl}. This approach also has to pay some price, as the $\kk$-nest is not a connected component. We then need to make a good choice of parameters, to ensure the consecutive sprinkling can form a real giant cluster.
		
		We also introduce two variants of $\skl(A, e_i)$, which indicate the sprinkling cluster issued from $\partial^\circ_{i, +}$ and $\partial^\circ_{i, -}$
		\begin{align}\label{eq.def.goodskl.pm}
			\forall \sigma \in \{+,-\}, \qquad \skl^\sigma (A, e_i) &\equiv \cC_{s}\Ll(\PTR_{i,\sigma} A , \, \PTR_{-i} \bigcup_{\cC' \in \mathscr{S}(z+e_i)} \cC' \Rr).
		\end{align}
		We also define $\mathscr{S}^{\sigma}_{\mathrm{off}}( A, e_i), \mathcal{V}^{\sigma}_{\mathrm{off}}( A, e_i)$ similarly for $\sigma \in \{+,-\}$.
		
		\begin{prop}\label{prop.proba.k-nest}
			Fix $\epsilon_0 \in (0,1)$ satisfying $\epsilon_0 \ll 1$, there exist finite positive constants 
			\begin{align}\label{eq.vector.k-nest}
				\Ll(n_0, \kk, K_{\mathrm{skl}}, m_{\mathrm{off}}, M_{\mathrm{off}}, M_{\mathrm{coal}}\Rr),
			\end{align}
			only depending on $\epsilon_0,d$ and satisfying
			\begin{align}\label{eq.condition.k-nest}
				\kk \geq 10, \qquad m_{\mathrm{off}} = 2 \kk^2.
			\end{align}
			such that for all $z \in \ell_{0} \Zd$ and $\p \searrow \p_c$, we have the following three statements:
			\begin{itemize}
				\item For every $i \in \{\pm 1, \cdots , \pm d\}$, when $Q(z)$, $Q(z+e_i)$ are pre-good, then for every $\kk$-nest $\CC \subset \fC(z)$, we have 
				\begin{align}\label{eq.proba.k-nest}
					\P_{\p}\Ll[\skl(\CC, e_i) \text{ is typical } \vert \fC(z), \fC(z+e_i) \Rr] \geq 1 - \frac{\epsilon_0}{4}.
				\end{align}
				\item If $\skl(\CC, e_i)$ is typical, then there exits at least one typical spanning cluster $\cC \subset \CC$ satisfying 
				\begin{align}\label{eq.pigeonhole.k-nest}
					\vert \mathscr{S}_{\mathrm{off}}(\cC,e_i) \vert \geq 2 \kk.
				\end{align}
				We call all such spanning cluster(s) \emph{good spanning cluster(s)}.
				\item \eqref{eq.proba.k-nest} and \eqref{eq.pigeonhole.k-nest} also hold for $\skl^+(\CC, e_i)$ and $\skl^-(\CC, e_i)$.
				
			\end{itemize}
		\end{prop}
		\begin{proof}
			Let us start with the constant $\kk$. Following Corollary~\ref{cor.proba.goodskl}, we fix $\mathbf{c}_{\text{\ref{eq.sprinkling.AB.condition.disjoint}}}$ as
			\begin{align*}
				\mathbf{c}_{\text{\ref{eq.sprinkling.AB.condition.disjoint}}} \equiv \mathbf{c}_{\text{\ref{eq.sprinkling.AB.condition.disjoint}}}(\epsilon_0 / 8).
			\end{align*}
			Afterwards, notice that every typical spanning cluster has at least $c_{\mathrm{b}} \ell_{0}^2$ regular points (see (4) of Definition~\ref{def.spanning.z}), and every pre-good cube contains at least $c_{\ref{eq.pregood}} \ell_{0}^{d-6}$ typical spanning clusters (see \eqref{eq.pregood}), we need the following two inequalities to make $\CC$ verify the conditions in \eqref{eq.proba.goodskl.condition}
			\begin{align*}
				\kk c_{\mathrm{b}} \ell_{0}^2 &\geq  \mathbf{c}_{\text{\ref{eq.sprinkling.AB.condition.disjoint}}}(\epsilon_0 / 8) \ell_{0}^2, \\
				\kk c_{\mathrm{b}} \ell_{0}^2 \cdot c_{\ref{eq.pregood}} \ell_{0}^{d-6} \cdot  c_{\mathrm{b}} \ell_{0}^2 &\geq  \mathbf{c}_{\text{\ref{eq.sprinkling.AB.condition.disjoint}}}(\epsilon_0 / 8) \ell_{0}^{d-2}.
			\end{align*}
			This suggests a choice
			\begin{align}\label{eq.k_choice}
				\kk := \max \Ll\{10, \frac{ \mathbf{c}_{\text{\ref{eq.sprinkling.AB.condition.disjoint}}}(\epsilon_0 / 8)}{c_{\mathrm{b}}}, \frac{ \mathbf{c}_{\text{\ref{eq.sprinkling.AB.condition.disjoint}}}(\epsilon_0 / 8)}{c_{\ref{eq.pregood}} c^2_{\mathrm{b}}}\Rr\}.
			\end{align}
			Once the value $\kk$ is determined, we set $m_{\mathrm{off}} = 2 \kk^2$, and apply Corollary~\ref{cor.proba.goodskl} to determine other parameters and conclude \eqref{eq.proba.k-nest}.
			
			The result \eqref{eq.pigeonhole.k-nest} is purely deterministic. For $\CC = \bigcup_{j=1}^{\kk} \cC_{j}$,  notice that 
			\begin{align*}
				\mathscr{S}_{\mathrm{off}}(\CC,e_i) =  \bigcup_{j=1}^{\kk} \mathscr{S}_{\mathrm{off}}( \cC_{j},e_i). 
			\end{align*}
			The subadditivity implies that 
			\begin{align*}
				2 \kk^2 = m_{\mathrm{off}} \leq \vert \mathscr{S}_{\mathrm{off}}(\CC,e_i) \vert \leq \sum_{j=1}^{\kk} \vert \mathscr{S}_{\mathrm{off}}( \cC_{j},e_i) \vert. 
			\end{align*}
			Then the pigeonhole principle implies \eqref{eq.pigeonhole.k-nest}.
			
			Finally, (4) of Definition~\ref{def.spanning.z} actually ensures two sets of regular boundary of size $c_{\mathrm{b}}\ell_{0}^2$, thus both $\PTR_{i, +} \CC$ and $\PTR_{i, -} \CC$ make Corollary~\ref{cor.proba.goodskl} valid. Therefore, they also satisfy \eqref{eq.proba.k-nest} and \eqref{eq.pigeonhole.k-nest}.
		\end{proof}
		
		\begin{remark}\label{rmk.proba.pregood}
			In the rest of this section, without further explanation, we always set the constants $	\Ll(n_0, \kk, K_{\mathrm{skl}}, m_{\mathrm{off}}, M_{\mathrm{off}}, M_{\mathrm{coal}}\Rr)$ in function of $\epsilon_0$ and $d$ as in Proposition~\ref{prop.proba.k-nest}. 
		\end{remark}
		
		\begin{remark}\label{rmk.two_good_cluster}
			Actually, with some more work, one can argue that a $\kk$-nest $\CC$ has probability $1-\frac{\epsilon_0}{2}$ to generate two good sprinkling clusters thanks to $\PTR_{i, +} \CC$ and $\PTR_{i, -} \CC$. This illustrates the power of the definition \eqref{eq.def.regular_bd_2}. We still choose to present the version about $\CC$ for convenience, while the stronger version about $\PTR_{i, +} \CC$ and $\PTR_{i, -} \CC$, and the choice of parameter \eqref{eq.condition.k-nest} will be useful in Section~\ref{subsec.intersection}.
		\end{remark}

		\subsubsection{Exploration along the route}\label{subsubsec.algo}
		We are now ready to implement a renormalization at correlation length $\ell_{0}$. The main idea relies on Proposition~\ref{prop.proba.pregood} and Proposition~\ref{prop.proba.k-nest}, which connects the spanning clusters and the sprinkling clusters. This is better understood by an \emph{exploration process} below. 
		
		In the exploration process, a route $\Gamma \subset Z_{M, \ell_{0}}$ is given and assume pre-good for all the coarsened cube. We start from a $\kk$-nest, and explore its neighbors to connect other typical spanning clusters using the sprinkling property, and then repeat this procedure. If the exploration continues to the end, its output will be a self-avoided path  $\gamma \subset \Lambda_L$ in $\p$-percolation crossing $\Gamma$.

		\begin{definition}[Exploration process]\label{def.expl}
			The exploration process runs as follows. 
			
			\smallskip
			\textbf{Setting.} 
			\begin{itemize}
				\item A self-avoiding route $\Gamma= \{z_0,\dots,z_n\} \subset Z_{M, \ell_{0}}$, with $v_t$ the unit vector satisfying $z_{t+1}=z_t+\ell_0v_t$.
				\item Two sets of random variables independent with the percolation:
				\begin{itemize}
					\item $\sigma = \{\sigma_0, \dots, \sigma_n\}$ i.i.d. taking value in $\{+,-\}$ with probability $\frac{1}{2}$;
					\item $U = \{U_0, \dots, U_n\}$ i.i.d. taking value uniformly in $[0,1]$.
				\end{itemize}
			\end{itemize}
			
			\smallskip
			\textbf{Exploration.} 
			Assume that every cube $Q(z_t)$, $0\le t\le n$, is pre-good.
			\begin{enumerate}
				\item Using $U_0$, let $\CC_0$ be a uniformly chosen $\kk$-nest of $\mathscr S(z_0)$.
				\item Given the $\kk$-nest $\CC_t\subset\mathscr S(z_t)$ with $t<n$:
				\begin{itemize}
					\item if $\skl^{\sigma_t}(\CC_t,v_t)$ is not typical, the exploration stops;
					\item otherwise, implement two actions:
					\begin{itemize}
						\item let $\cC^*_t$ be the first cluster of $\CC_t$, in a fixed deterministic order, that is good in $\skl^{\sigma_t}(\CC_t,v_t)$;
						\item using $U_{t}$, let $\CC_{t+1}$ be a uniformly chosen $\kk$-nest of $\mathscr S^{\sigma_t}_{\mathrm{off}}(\cC^*_t,v_t)$.
					\end{itemize}
				\end{itemize}
				\item Let $\cC^*_n$ be the first cluster of $\CC_n$ if the exploration reaches $n$, and then it is called \emph{survival}.
			\end{enumerate}
			
			We denote by the exploration process by
			\begin{align*}
				\expl(\Gamma, \a; \sigma, U) := \{\CC_0, \skl^{\sigma_0}(\CC_0,v_0), \CC_1, \skl^{\sigma_1}(\CC_1,v_1),\cdots\}.
			\end{align*}
			
		\end{definition}

		We then extract an open path in $\p$-percolation from the exploration process.
		\begin{definition}[Microscopical path]\label{def.path}
			Given that the exploration process $\expl(\Gamma, \a; \sigma, U)$ survives, in a fixed deterministic order,  choose a simple path 
			\begin{align}\label{eq.gamma.skl}
				\gamma^{\skl}_t \subset \skl^{\sigma_t}(\cC^*_t,v_t),
			\end{align}
			connecting $\cC^*_t$ and $\cC^*_{t+1}$ for all $0 \leq t \leq n-1$, and choose a sprinkling path (Definition~\ref{def.skl})
			\begin{align}\label{eq.gamma.spn}
				\gamma^{\spn}_t \subset \cC^*_t,  
			\end{align}
			connecting $\gamma^{\skl}_{t-1}$ and  $\gamma^{\skl}_{t}$. For the endpoint $t = 0$ or $n$, $\gamma^{\spn}_t$ only needs to verify one condition. Set $\gamma(\Gamma, \a; \sigma, U)$ as 
			\begin{align}\label{eq.def.path}
				\gamma(\Gamma, \a; \sigma, U) :=  \bigcup_{t = 0}^{n-1} (\gamma^{\spn}_t  \cup \gamma^{\skl}_t) \cup \gamma^{\spn}_n.
			\end{align}
		\end{definition}

		We state several properties about the exploration process and the its associated path.
		\begin{prop}\label{prop.path}
			The following properties hold for $\gamma \equiv \gamma(\Gamma, \a; \sigma, U)$.
			\begin{enumerate}
				\item $\gamma$ is a self-avoiding path in $\p$-percolation.
				\item If $\Gamma$ is a route connecting $[\partial Z_{M, \ell_{0}} ]^{-}$ and $[\partial  Z_{M, \ell_{0}}]^{+}$, then $\gamma$ connects $[\partial \Lambda_{L-\ell_{0}/2}]^{-}$ and $[\partial  \Lambda_{L-\ell_{0}/2}]^{+}$. 
				\item $\gamma$ satisfies the estimate of length that 
				\begin{align}\label{eq.path.length.spn_skl}
					\vert \gamma^{\spn}_t \vert \leq K_{\mathrm{spn}}\ell_{0}^2, \qquad \vert \gamma^{\skl}_t \vert \leq K_{\mathrm{skl}}\ell_{0}^2,
				\end{align}
				and
				\begin{align}\label{eq.path.length}
					\vert \gamma \vert \leq \Ll(K_{\mathrm{spn}} + K_{\mathrm{skl}}\Rr)\ell_{0}^2 \vert \Gamma \vert.
				\end{align}
			\end{enumerate}
		\end{prop}
		\begin{proof}
			(1) it is clear that the construction in Definition~\ref{def.path} gives a path:
			\begin{itemize}
				\item Because $\skl^{\sigma_t}(\cC^*_t,v_t)$ connect $\cC^*_t$, and $\CC_{t+1}$ and $\cC^*_{t+1} \subset \CC_{t+1}$, there exists a path $\gamma^{\skl}_t$ to connect $\cC^*_t$ and $\cC^*_{t+1}$.
				\item Similarly, by the observation above, because $\cC^*_t$ is connected by $\gamma^{\skl}_t$ and $\gamma^{\skl}_{t-1}$, and $\cC^*_t$ is a cluster, there exists a path $\gamma^{\spn}_t$ to connect two sprinkling clusters.
			\end{itemize}
			After the construction, the spanning paths and sprinkling paths are connected alternatively, so they form a crossing path. It is self-avoiding as each segment like $\gamma^{\spn}_t$ and $\gamma^{\skl}_t$ are simple paths. Moreover, there is no intersection of bonds between the sprinkling clusters and spanning clusters from the definition; see Definition~\ref{def.skl}. Finally, the spanning clusters in different cubes are disjoint. This is also the case for the sprinkling clusters; see Lemma~\ref{lem.skl.independent}.
			
			(2) is a corollary from (1). It only connects the boundary of $\Lambda_{L-\ell_{0}/2}$, because the path near the boundary lives in a spanning cluster of $Q(z)$ defined in \eqref{eq.def.SQH.ell0.QSH}, which has distance $\frac{\ell_{0}}{4}$ at each end.
			
			(3) Recall the definition \eqref{eq.gamma.spn} and \eqref{eq.gamma.skl}, and the fact $\cC^*_t$ is typical spanning cluster and that $\skl^{\sigma_t}(\CC_t,v_t)$ is typical sprinkling cluster, then equation \eqref{eq.path.length.spn_skl} follows the point-wise bound in Definition~\ref{def.spanning.z} and Definition~\ref{eq.def.goodskl}
			\begin{align*}
				\vert \gamma^{\skl}_t \vert \leq \rin(\skl^{\sigma_t}(\cC^*_t,v_t)) &\leq K_{\mathrm{skl}}\ell_{0}^2, \\
				\vert  \gamma^{\spn}_t \vert \leq \rin(\cC^*_t)   &\leq K_{\mathrm{spn}} \ell_{0}^2.
			\end{align*}
			Therefore, we conclude that 
			\begin{align*}
				\vert \gamma \vert \leq \sum_{z_t \in \Gamma} \Ll(K_{\mathrm{spn}} \ell_{0}^2 + K_{\mathrm{skl}}\ell_{0}^2\Rr) \leq (K_{\mathrm{spn}} + K_{\mathrm{skl}}) \vert \Gamma \vert \ell_{0}^2.
			\end{align*}
		\end{proof}

		\begin{remark}\label{rmk.renormalization}
			In the end of this section, let us make some quick discussion. One can, of course, state the exploration process in a more general way. For example, if we do not fix a route $\Gamma$, but let it make exploration like Definition~\ref{def.expl} and visit every site in $\ell_{0} \Zd$ only once, then it will form a supercritical percolation in coarsened sites. Then we can further use the renormalization in supercritical regime, and choose some crossing path. See Figure~\ref{fig.Algorithm} for an illustration.
			
			Nevertheless, the approach above has a hidden technical issue. When analyzing the intersection following Lemma~\ref{lem.flow}, the crossing route chosen by the supercritical renormalization has already revealed the information in microscopical scale. To make use of the properties of sprinkling clusters, especially those in Section~\ref{subsec.coal}, it is better to have a deterministic route.
			
			Then the competition appears. There is no way why the exploration process always survives in a given route, even in very supercritical percolation. One needs a plan to pick out those explorations which can survive.
			
			\begin{figure}
				\centering
				\includegraphics[width=0.25\textwidth]{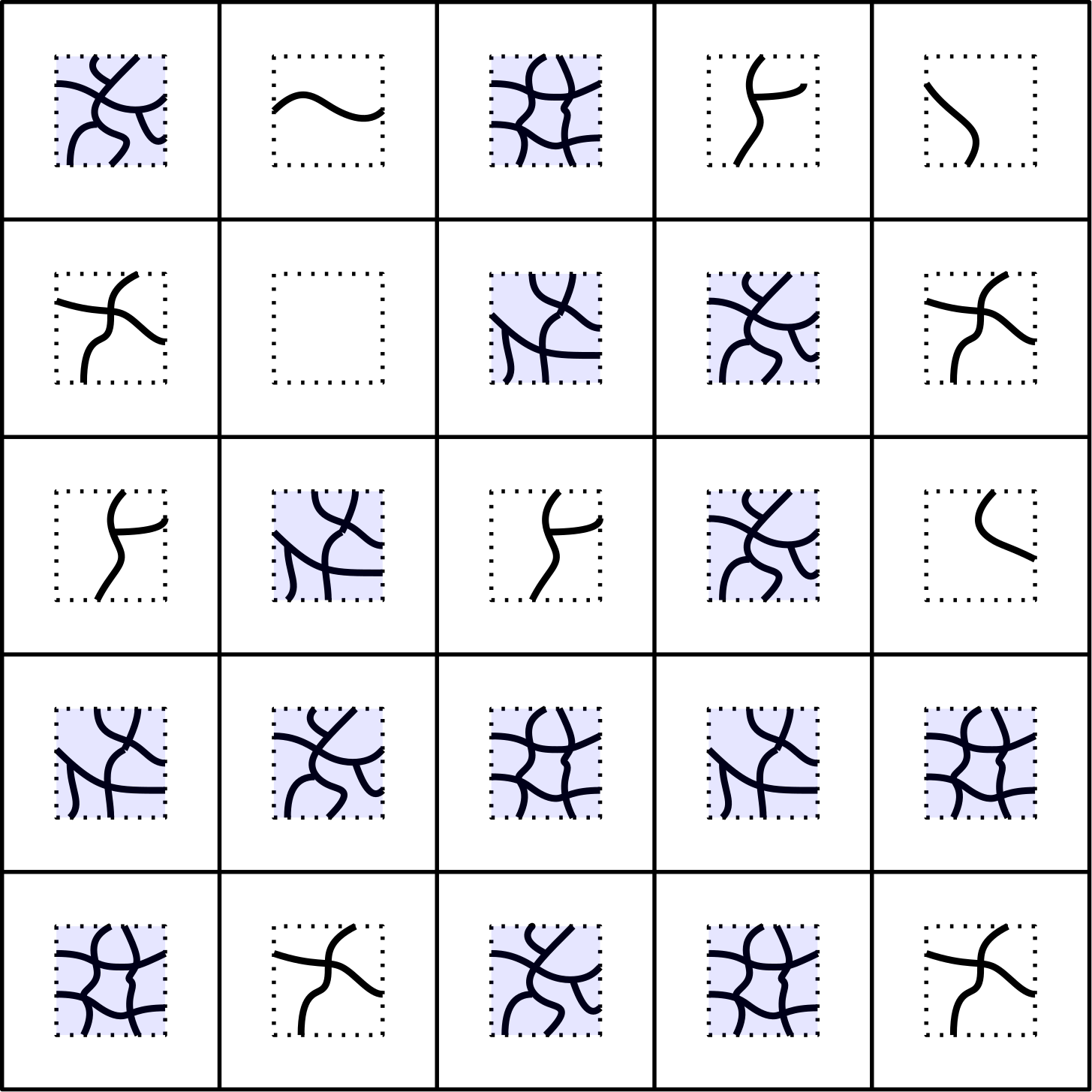}\hspace{1cm}
				\includegraphics[width=0.25\textwidth]{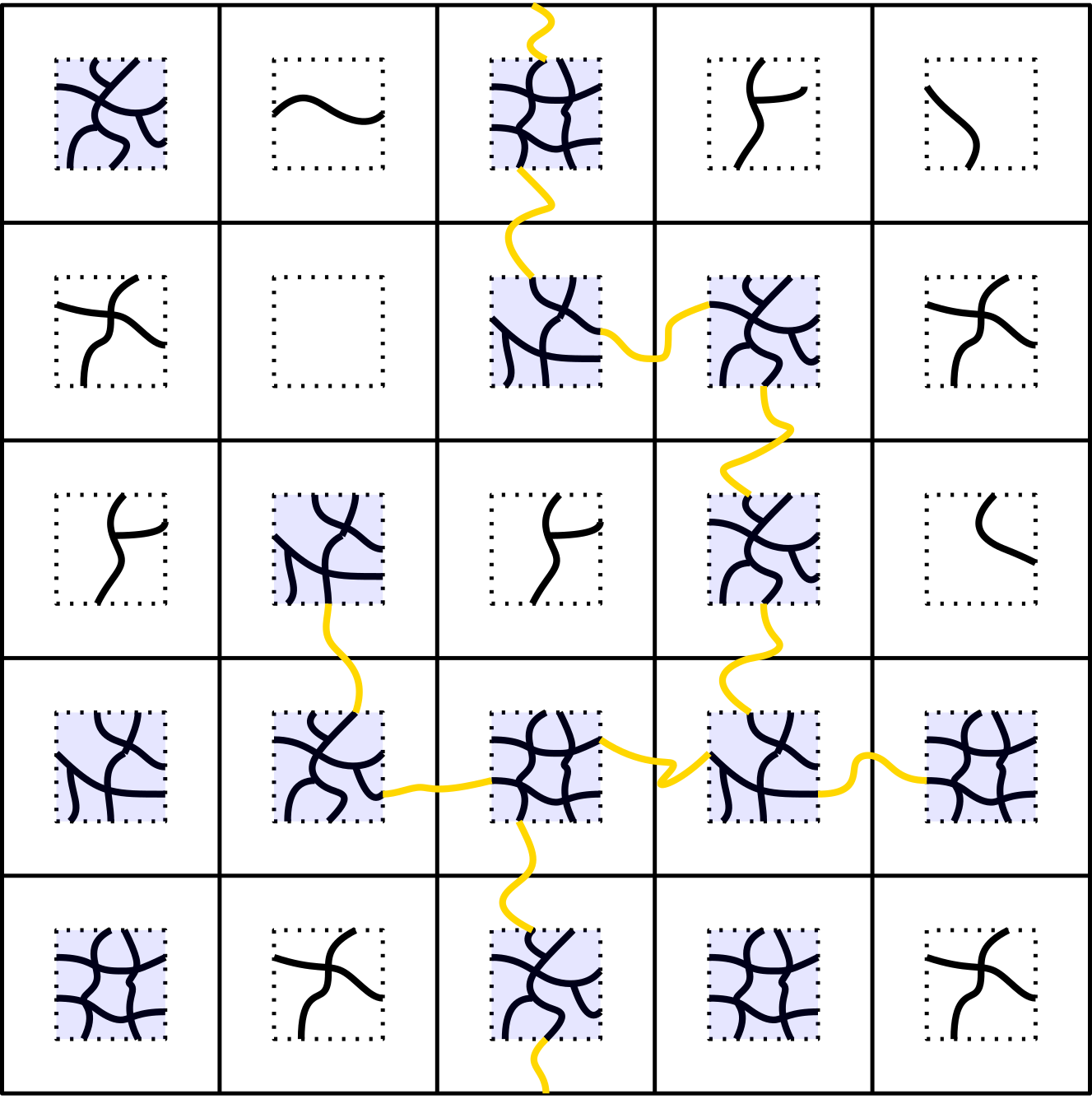}\\
				\vspace{1cm}
				\includegraphics[width=0.25\textwidth]{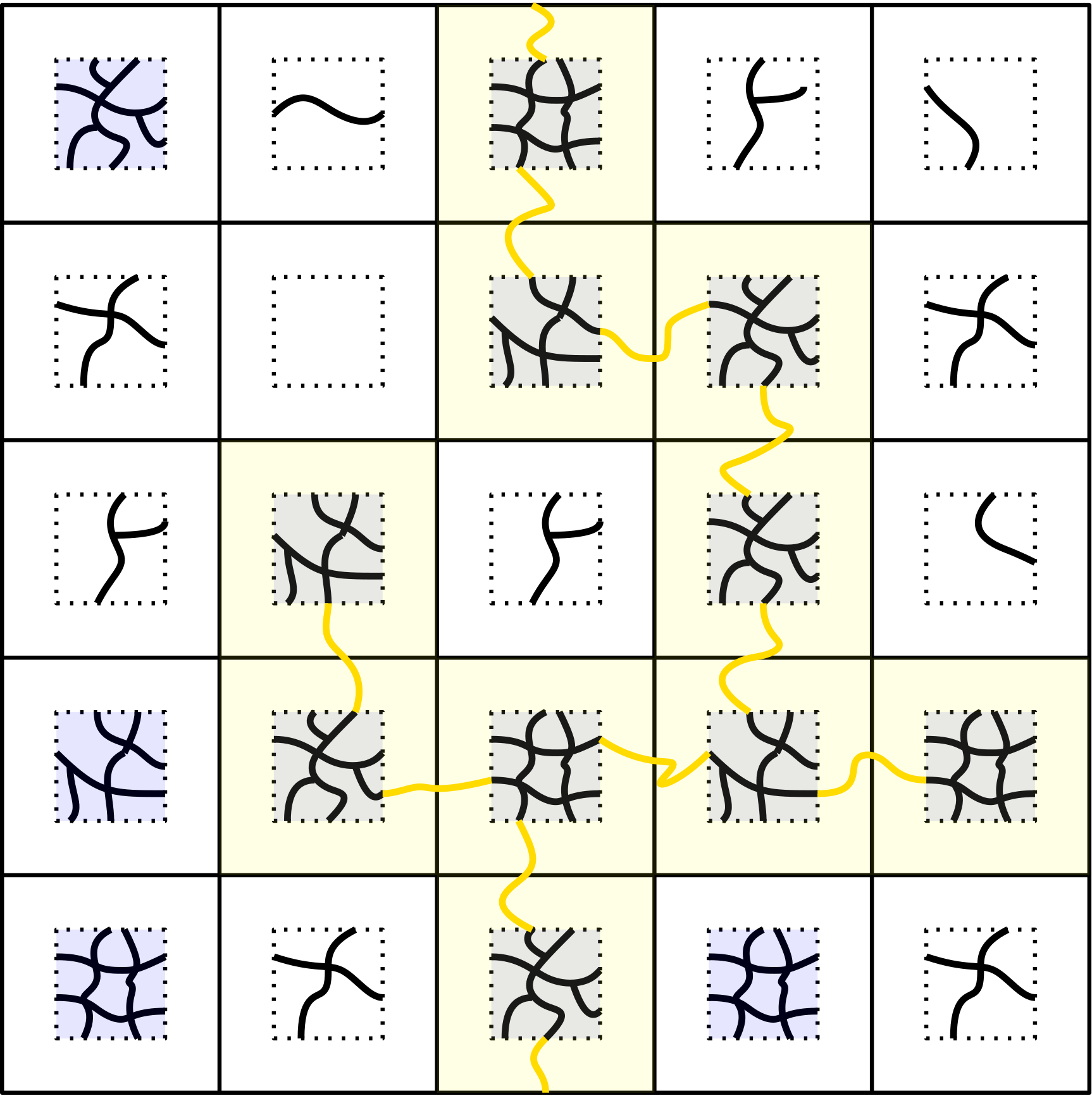}\hspace{1cm}
				\includegraphics[width=0.25\textwidth]{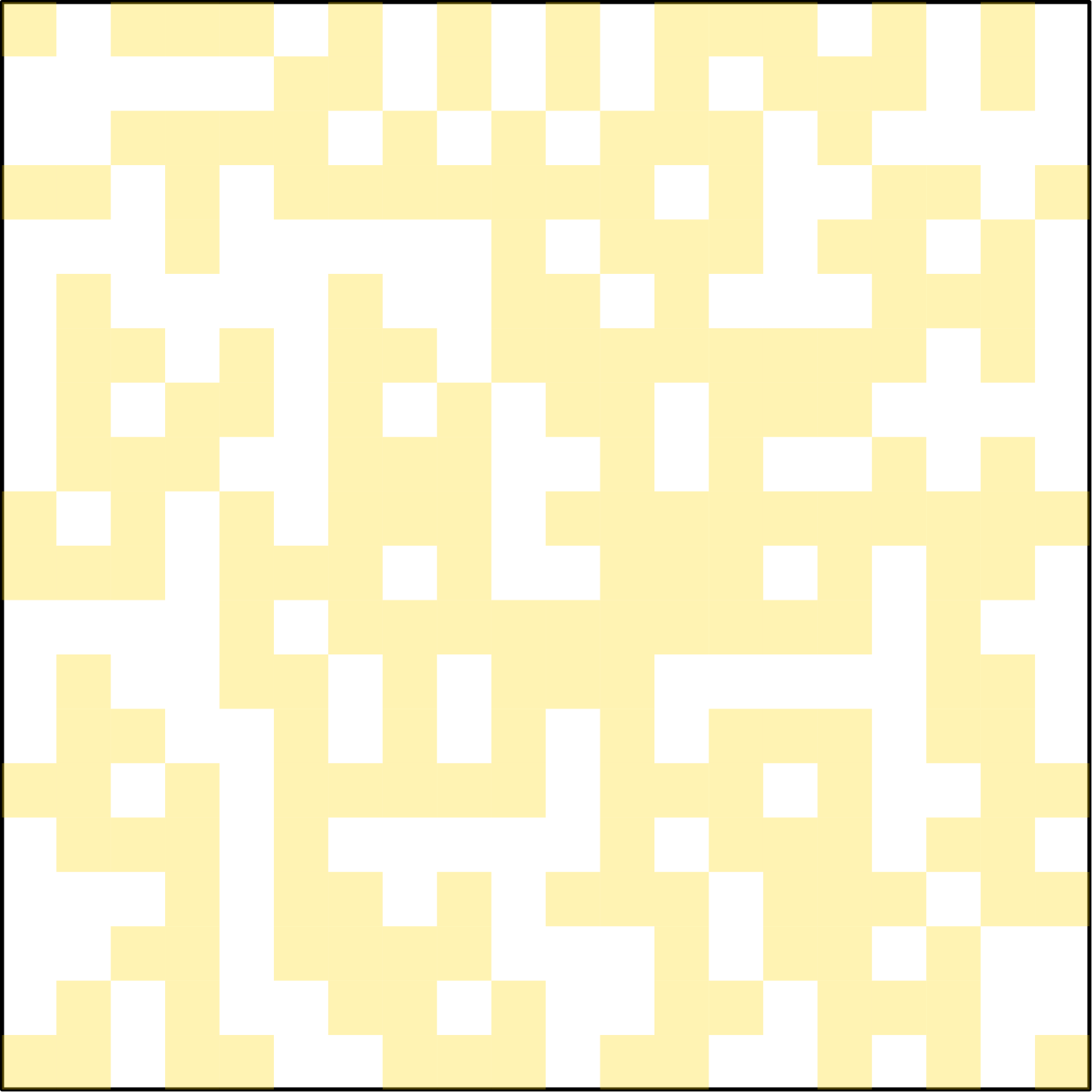}\\
				\vspace{1cm}
				\includegraphics[width=0.25\textwidth]{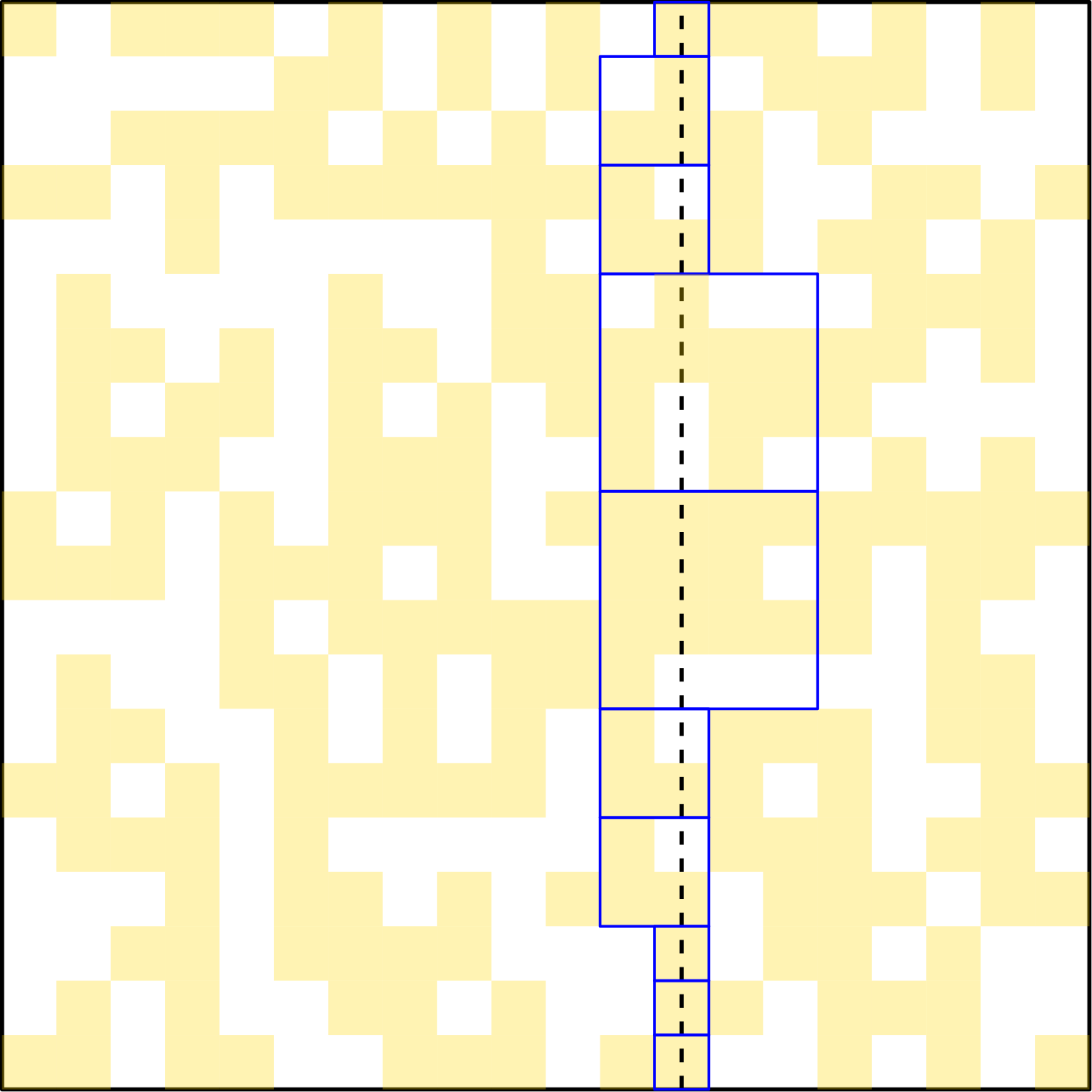}\hspace{1cm}
				\includegraphics[width=0.25\textwidth]{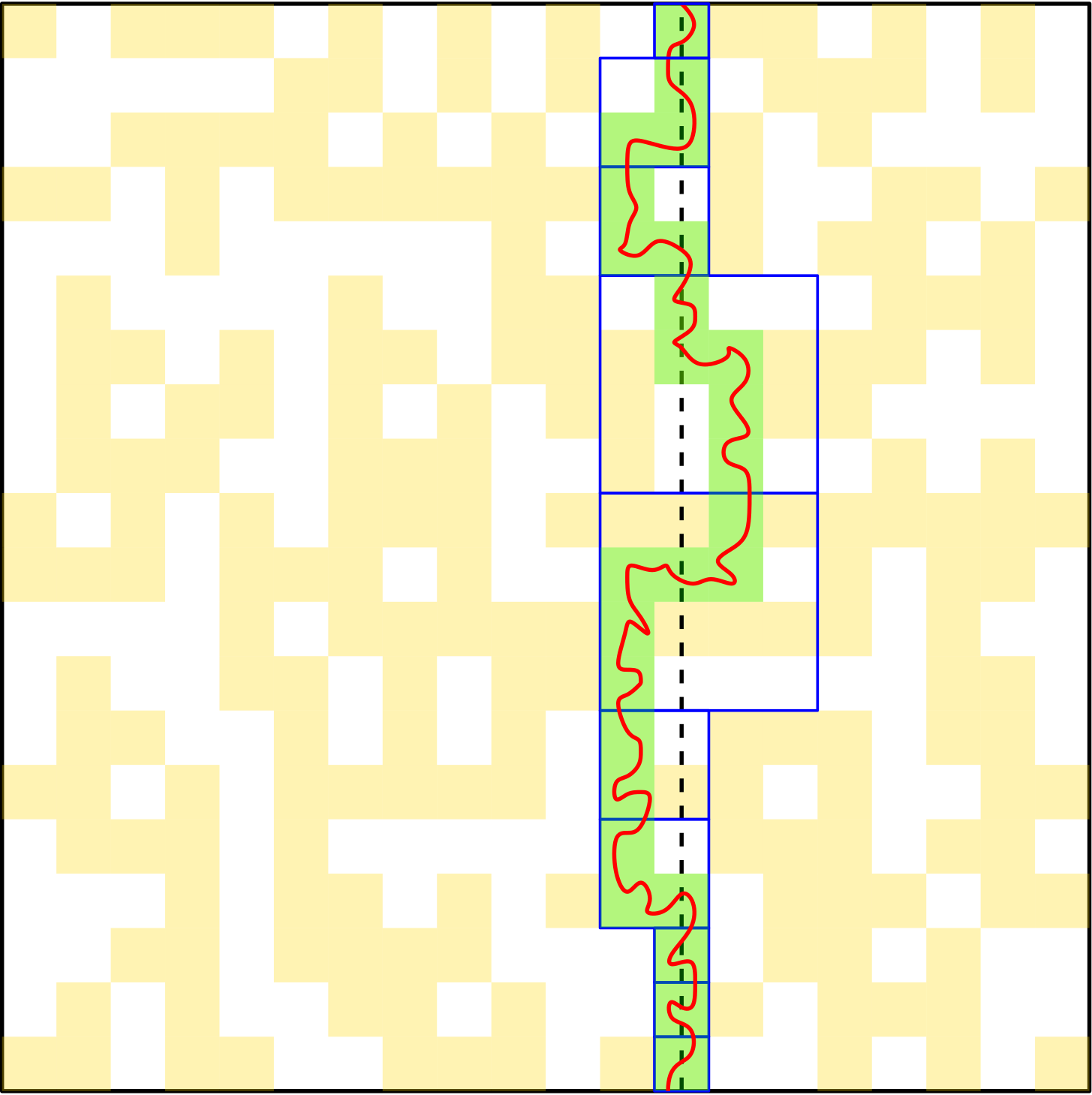}
				\caption{An illustration of the generalized exploration process, which does not have a given route but explores all the possible directions.}\label{fig.Algorithm}
			\end{figure}
		\end{remark}

		\subsection{Allocation of flow with EIT}\label{subsec.EIT}
		This subsection answers the question mentioned in last paragraph, and the strategy follows the work \cite{BPP98} by Benjamini, Pemantle, and Peres. Their work constructed a measure on random oriented path on $\Zd$, so that two independent paths have an \emph{exponential intersection tail} (EIT). Then, this measure is applied to the resistance of supercritical percolation, and proved the transience of the random walk there for $\Zd, d \geq 3$. 
		
		This subsection also contains these two ingredients. We first construct a route measure $\mu$, so that the intersections in macroscopic level is well-controlled. Afterwards, we use the route measure to re-weight the exploration process, and design the random flow.

		\subsubsection{Route measure}
		
		We first introduce the oriented route in $Z_{M, \ell_{0}}$.
		\begin{definition}[Oriented route]\label{def.route.o}
			A route $\Gamma = \{z_0, z_1, \cdots \}$ is called \emph{oriented route} if the unit vector in ${z_{t+1} = z_t + \ell_{0} v_t}$ satisfies
			\begin{align*}
				v_t = \Ll\{\begin{array}{ll}
					e_1 & t \text{ is even};\\
					e_i, \quad i \in \{2, \cdots, d\} & t \text{ is odd}.
				\end{array}\Rr.
			\end{align*} 
			A step is \emph{vertical} if $v_t=e_1$,  and \emph{horizontal} otherwise. 	Let the set $\mathscr{G}_{M, \ell_{0}}$ represents all the crossing oriented routes in $Z_{M, \ell_{0}}$, i.e. 
			\begin{align*}
				\mathscr{G}_{M, \ell_{0}} := \{\Gamma \subset Z_{M, \ell_{0}}: \Gamma \text{ is an oriented route connecting } [\partial  Z_{M, \ell_{0}}]^{-} \text{ and } [\partial  Z_{M, \ell_{0}}]^{+} \}.
			\end{align*}
		\end{definition}
		In another word, the oriented route  alternatively makes the vertical steps and horizontal steps. Therefore, given the value $t$, we know the height $z_t \cdot e_1$ in the orient route. 
		
		We then construct a probability measure $\mu$ on $\mathscr{G}_{M, \ell_{0}}$, and we are especially interested in the intersections between two oriented routes $\Gamma, \Gamma'$ sampled independently under $\mu$. The following proposition states that the intersection is possible to show an exponential tail. For the moment, there is no percolation in the statement and proof.
		\begin{prop}[Exponential intersection tails]\label{prop.EIT}
			There exist $q_0(d)\in(0,1)$ and $C(d)<\infty$ such that, for all $M$ large enough, there is a probability measure $\mu$ over $\mathscr{G}_{M, \ell_{0}}$ satisfying
			\begin{equation}\label{eq.nw.EIT}
				\forall k \in \N_+, \qquad \sum_{\substack{\Gamma, \Gamma' \in \mathscr{G}_{M, \ell_{0}} \\ \vert \Gamma \cap \Gamma' \vert \ge k}}\mu\otimes\mu\big(\Gamma, \Gamma'\big) \le C M^{2-d} q_0^{k}.
			\end{equation}
		\end{prop}
		
		\begin{proof}
			The proof is divided into $2$ parts.
			
			\textit{Step~1: oriented route by random walk.} We start the construction of oriented route by random walk and study the intersection. For simplicity, we work on a cube $B_M$ of side length $M$, which is assumed odd
			\begin{align*}
				B_M = \{0,1, \cdots, M-1\} \times \Ll(\Ll\lfloor \frac{-M}{2} \Rr\rfloor, \cdots, \Ll\lfloor \frac{M}{2} \Rr\rfloor\Rr)^{d-1}.
			\end{align*}
			In probability space $P$, we sample a lazy random walk $X = (X_j)_{j \geq 0}$ in $\mathbb{Z}^{d-1}$ with the following transition probability. By default, we regard $(e_2, \cdots, e_d)$ the canonical basis of $\mathbb{Z}^{d-1}$.
			\begin{align*}
				P[X_{j+1} = y \vert X_{j}=x] = \Ll\{\begin{array}{ll}
					\frac{1}{2} & y = x;\\
					\frac{1}{2(d-1)} & y = x + e_i, \quad i \in \{2, \cdots, d\} .
				\end{array}\Rr.
			\end{align*}
			Moreover, we sample the initial condition $X_0$ uniformly in $\{x \in \mathbb{Z}^{d-1}: \vert x \vert \leq \frac{M}{10}\}$. The random walk then defines an oriented route in the following way
			\begin{align*}
				L(X) := \big((0,X_0),\ (1,X_0),\ (1,X_1),\ (2,X_1),\ (2,X_2),\ \dots,\ (M-2,X_{M-2}),\ (M-1,X_{M-2})\big).
			\end{align*}
			It is clear the route $L(X)$ is self-avoiding, as it makes alternatively the vertical and horizontal moves. 
			
			We then sample two independent copies $X, X'$, and study the intersection of $L(X) \cap L(X')$. It is reduced to the quantity
			\begin{align*}
				N(X,X') := \#\{0 \leq j \leq M-1: X_j = X'_j, \text{ or } X_j = X'_{j-1}, \text{ or } X_{j-1} = X'_j\}.
			\end{align*}
			That is, $N(X,X')$ counts the number of height where $L(X)$ intersects $L(X')$. By the definition of oriented route, there are at most $2$ intersections at every height. Therefore, we have clearly
			\begin{align}\label{eq.NXX.1}
				N(X,X') \leq \vert L(X) \cap L(X') \vert \leq 2 N(X,X').
			\end{align}
			Then quantity $N(X,X')$ can be further simplified. We define that
			\begin{align*}
				D_j:=X_j-X'_j.
			\end{align*}
			Then the intersection of two routes is dominated by the event $\vert D_j \vert \leq 2$, because every step the random walk moves at most $1$. Therefore, we have
			\begin{align}\label{eq.NXX.2}
				N(X,X') \leq \#\{0 \leq j \leq M-1: D_j \in {\mathsf B}\},
			\end{align}
			with $\mathsf B:=\{y \in \ZZ^{d-1}: |y| \le 2\}$.

			The problem of intersection is then reduced to $(D_j)_{j \geq 0}$, which is a symmetric random walk of in $\ZZ^{d-1}$. Clearly, it is transient, because $d > 6$. Its number of visit is closely related to the hitting probability and Green function. We define the hitting time $\tau_{{\mathsf B}}$
			\begin{align*}
				\tau_{{\mathsf B}}:=\inf\{j \geq 1: D_j \in {\mathsf B}\}.
			\end{align*}
			For $x \in {\mathsf B}^c$, it is equal to $\inf\{j \geq 0: D_j \in {\mathsf B}\}$, so the hitting probability 
			\begin{align*}
				h(x):=P_x[\tau_{{\mathsf B}} < \infty],
			\end{align*}
			is a harmonic extension from $\Ind{\mathsf B}$, which goes to $0$ when $\vert x \vert \to \infty$. By the maximum principle, we thus have 
			\begin{align*}
				\max_{x \in \mathsf B^c} h(x) < 1.
			\end{align*}
			For every $x \in {\mathsf B}$, because 
			\begin{align*}
				P_x[D_2 \in {\mathsf B}^c] \geq \frac{1}{16(d-1)^4},
			\end{align*}
			the random walk uniformly has some probability to escape to infinity. We thus have
			\begin{align}\label{eq.def.q0}
				q_0:=\max_{y\in\mathsf B}h(y)<1.
			\end{align}
			Then by the strong Markov property, every additional visit of $\mathsf B$ has at most a probability $q_0$, so the number of visits of $\mathsf B$ has a geometric decay
			\begin{align*}
				P[\#\{0 \leq j \leq M-1:D_j \in \mathsf B\} \ge k] \le q_0^{k-1} P\Ll[\exists j\le M-1: D_j\in\mathsf B\Rr].
			\end{align*}
			Recall that $X_0$ uniformly in $\{x \in \mathbb{Z}^{d-1}: \vert x \vert \leq \frac{M}{10}\}$, so  we further have
			\begin{align}\label{eq.D.hitting}
				\begin{split}
					P[\#\{0 \leq j \leq M-1:D_j \in \mathsf B\} \ge k] &\leq q_0^{k-1} \sum_{j=0}^{M-1}P\Ll[D_j\in\mathsf B\Rr] \\
					&\lesssim q_0^{k-1} \cdot M \cdot M^{1-d} \\
					&\lesssim q_0^{k-1} \cdot M^{2-d}.
				\end{split}
			\end{align}
			Actually, the distribution of $D_j$ is a convolution of uniform distribution in $\Ind{ \vert \cdot \vert \leq \frac{M}{10}}$ with a heat kernel of typical radius $\sqrt{j} \ll M$, when $M$ is large. That is why we can insert an upper bound $M^{1-d}$ above.
			
			In summary, \eqref{eq.NXX.1}, \eqref{eq.NXX.2}, \eqref{eq.D.hitting} together proves
			\begin{align}\label{eq.nw.EIT.pre}
				P[\vert L(X) \cap L(X')\vert \ge k]  \lesssim M^{2-d} \cdot q_0^{k}.
			\end{align}
			
			\smallskip
			
			\textit{Step~2: change of measure.} Notice that, the oriented route $L(X)$ constructed above can exit $B_M$, which is not desired. We thus use the method of reject, and restrict the measure of route on a probability $Q$
			\begin{align*}
				\frac{\d Q}{\d P} = \frac{\Ind{X \subset B_M}}{P[X \subset B_M]}. 
			\end{align*} 
			The trajectory of $X$ has a slight drift of slope $\frac{1}{2(d-1)}$ in all the direction $e_2,e_3, \cdots, e_d$, with a typical fluctuation $\sqrt{M}$. Since the starting point $X_0$ has a distance $\frac{9 M}{20}$ to the boundary, so it is very rare to happen $\{X \nsubseteq B_M\}$ for large $M$. The full argument can be done by the standard reflex principle or Chernoff bound. Here we omit the details and just claim
			\begin{align*}
				\liminf_{M \to \infty} P[X \subset B_M] \geq \frac{1}{2}.
			\end{align*}
			Under the restricted measure, we have the following bound for two independent samples
			\begin{align*}
				Q[\vert L(X) \cap L(X')\vert \ge k] &= \frac{1}{P[X \subset B_M]^2} P[\vert L(X) \cap L(X')\vert \ge k, X \subset B_M, X' \subset B_M]\\
				&\leq 4 P[\vert L(X) \cap L(X')\vert \ge k]\\
				&\lesssim 4 M^{2-d} \cdot q_0^{k}.
			\end{align*}
			That is, the statement is still valid under $Q$. 
			
			Finally, concerning the statement for the measure of oriented route on $\mathscr{G}_{M, \ell_{0}}$, we define a mapping of shift and scaling $\Phi : B_M \to Z_{M, \ell_{0}}$
			\begin{align*}
				\Phi (x) = \ell_{0} \cdot \Ll( x - \Ll\lfloor \frac{M}{2} \Rr\rfloor e_1 \Rr).
			\end{align*} 
			Recall that we choose $M$ to be odd for convenience, so this is a bijection. The measure $\mu$ in the statement is then induced by the mapping
			\begin{align*}
				\forall \Gamma \in \mathscr{G}_{M, \ell_{0}}, \qquad \mu(\Gamma) := Q[X = \Phi^{-1} (\Gamma) ].
			\end{align*}
			One can verify that
			\begin{align*}
				\sum_{\substack{\Gamma, \Gamma' \in \mathscr{G}_{M, \ell_{0}} \\ \vert \Gamma \cap \Gamma' \vert \ge k}}\mu\otimes\mu\big(\Gamma, \Gamma'\big) = 	Q[\vert L(X) \cap L(X')\vert \ge k],
			\end{align*}
			so \eqref{eq.nw.EIT} is verified.
		\end{proof}
		\begin{remark}
			We choose to express $\mu$ in convention of a general measure rather than a probability, although the latter is also proved. This aims to distinguish that the oriented route $\Gamma$ does not rely on $\PP$. 
		\end{remark}

		\subsubsection{Allocation of flow}
		Now we construct the random flow. The construction is inspired by \cite[Section~2]{BPP98}, which samples the oriented route first, then picks the one with crossing path by method of reject. The weight in our construction is slightly more complicated as the percolation is revealed by the exploration process in Definition~\ref{def.expl}.

		We keep $\mu$ in Proposition~\ref{prop.EIT} for the measure of oriented route, and use $\PP$ for other random elements including the percolation. We then denote by 
		\begin{align*}
			\PP_{\sigma, U}[(\,\cdot \,)] := \PP[(\sigma, U) \in (\,\cdot \,)], 
		\end{align*}
		the induced probability for $(\sigma, U)$. As stated in Definition~\ref{def.expl}, this probability is a uniform product measure over $(\{+,-\} \otimes [0,1])^{\otimes \Gamma}$ and is independent of the percolation $\a$. As we use the coupling space of the percolation, we also abuse the notation by writing 
		\begin{align*}
			\Eb[(\,\cdot \,) \, \vert \a] \equiv \Eb[(\,\cdot \,) \, \vert \p_c\text{-percolation}, \p\text{-percolation}].
		\end{align*}
		
		We first introduce two constants in our construction.
		\begin{lemma}
			Given a route $\Gamma$, define that 
			\begin{align}\label{eq.def.st}
				q:=\P_{\p_c}[Q(z_t)\text{ is pre-good}], \qquad s_t:=\PP\big[\skl^{\sigma_t}(\CC_t,v_t)\text{ is typical}\,\big|\,\mathcal G,\CC_t,\sigma_t\big].
			\end{align}
			$q$ does not depend on the index $t$, and they satisfy 
			\begin{align}\label{eq.st.bound}
				\forall z_t \in \Gamma, \qquad q, s_t \geq 1- \frac{\epsilon_0}{4}.
			\end{align}
		\end{lemma}
		\begin{proof}
			The pre-goodness only depends on the cube $Q(z)$ introduced in \eqref{eq.def.SQH.ell0.QSH}, thus is i.i.d. and does not dependent on the index. The estimate in \eqref{eq.st.bound} is the consequence of \eqref{eq.proba.pregood}, \eqref{eq.proba.k-nest}, and Remark~\ref{rmk.two_good_cluster}.
		\end{proof}
		
		\begin{definition}[Weight of random flow]\label{def.weight}
			Given the percolation $\a$, sample the couple 
			\begin{align*}
				\Big(\Gamma, \expl(\Gamma, \a; \sigma, U)\Big),
			\end{align*}
			following measure
			\begin{align}\label{eq.weight.couople}
				\frac{1}{W}  \omega_\Gamma\, \zeta_\Gamma\,  \d (\mu \otimes \PP_{\sigma, U}),
			\end{align}
			where the weight is defined by 
			\begin{align}\label{eq.weight}
				\omega_\Gamma:=\prod_{z_t \in \Gamma}\frac{\Ind{Q(z_t)\text{ pre-good}}}{q},\qquad
				\zeta_\Gamma:=\prod_{\substack{z_t \in \Gamma \\ t < \vert \Gamma \vert - 1}}\frac{\Ind{\skl^{\sigma_t}(\CC_t,v_t)\text{ typical}}}{s_t},
			\end{align}
			and the normalization factor $W$ is defined as
			\begin{align}\label{eq.W}
				W := \sum_{\Gamma \in \mathscr{G}_{M, \ell_{0}}} \mu(\Gamma) \int \omega_\Gamma\, \zeta_\Gamma\, \d \PP_{\sigma, U}.
			\end{align}
			When $W = 0$, \eqref{eq.weight.couople} is set $0$ by default. 	The couple $(\Gamma, \expl(\Gamma, \a; \sigma, U))$ then defines a path $\gamma(\Gamma, \a; \sigma, U)$ as Definition~\ref{def.path}, and further entails a unit flow
			\begin{align}\label{eq.def.flow_EIT}
				\theta(e) := \frac{1}{W}\sum_{\Gamma \in \mathscr{G}_{M, \ell_{0}}} \mu(\Gamma) \Eb\Ll[\omega_\Gamma\, \zeta_\Gamma\, \Ind{e \in \gamma(\Gamma, \a; \sigma, U)} \vert \a\Rr]. 
			\end{align}
		\end{definition}

		We then justify the well-definedness.
		\begin{lemma}\label{lem.weight.mean}
			$W$ only depends on $\a$ and satisfies $\Eb[W]=1$. Moreover, when ${W > 0}$, \eqref{eq.weight.couople} defines a probability measure, and $\expl(\Gamma, \a; \sigma, U)$ survives under this measure. Then, $\gamma(\Gamma, \a; \sigma, U)$ is a self-avoiding path connecting $\gamma$ connects $[\partial \Lambda_{L-\ell_{0}/2}]^{-}$ and $[\partial  \Lambda_{L-\ell_{0}/2}]^{+}$, and unit flow given by \eqref{eq.def.flow_EIT} is well-defined.
		\end{lemma}
		\begin{proof}
			By the definition \eqref{eq.W}, it is clear that $\a$ is the only remaining random element, as the other variables have been integrated. The expectation of $W$ can be calculated using Fubini's lemma and conditional expectation  
			\begin{align*}
				\Eb[W] = \sum_{\Gamma \in \mathscr{G}_{M, \ell_{0}}} \mu(\Gamma)\Eb[ \omega_\Gamma\, \Eb[\zeta_\Gamma\,  \, \vert \G]]. 
			\end{align*}

			Define the filtration
			\begin{align}\label{eq.def.filtration}
				\mathcal F_{t-1} := \sigma\Ll(\mathcal G,  \Ll(T(z_u,z_{u+1}), \sigma_u, U_u \Rr)_{0 \leq u \leq t-1}\Rr).
			\end{align}
			The tube $T(z_t,z_{t+1})$ is independent of $\mathcal F_{t-1}$ given $\mathcal G$, so 
			\begin{align*}
				\Eb\Ll[\frac{\Ind{\skl^{\sigma_t}(\CC_t,v_t)\text{ typical}}}{s_t} \, \Big\vert \mathcal F_{t-1}\Rr]=1.
			\end{align*}
			Hence the product in $\zeta_\Gamma$ has a $(\mathcal F_t)$-martingale structure, and $\Eb[\zeta_\Gamma\mid\mathcal G]=1$.
			The pre-good events of distinct cubes are independent, since the cubes are disjoint. This implies that $\Eb[\omega_\Gamma]=1$. 
			Summing against $\mu$ gives the claim $\Eb[W]=1$.
			
			When $W > 0$, we notice that the sum over weight \eqref{eq.weight} is 
			\begin{align*}
				\sum_{\Gamma \in \mathscr{G}_{M, \ell_{0}}} \mu(\Gamma) \int \frac{1}{W} \omega_\Gamma\, \zeta_\Gamma\, \d \PP_{\sigma, U} = 1.
			\end{align*}
			So it defines a probability measure. The weight \eqref{eq.weight} requires the pre-goodness along $\Gamma$ and all the sprinkling clusters are typical in the exploration, thus it survives and can define  $\gamma(\Gamma, \a; \sigma, U)$ following Definition~\ref{def.path}.
			
			The statement about $\gamma(\Gamma, \a; \sigma, U)$  follows Lemma~\ref{lem.flow} and Proposition~\ref{prop.path}. Moreover, we notice the following identity thanks to the independence between $\PP_{\sigma, U}$ and $\a$.
			\begin{align}\label{eq.id.PP_sigmaU}
				\Eb[(\,\cdot \,) \vert \a] = \int (\,\cdot \,) \d \PP_{\sigma, U}.
			\end{align}
			Therefore, the unit flow is obtained just a average over all the possible $\sigma$ and $U$.
		\end{proof}

		\subsection{Intersection between flows}\label{subsec.intersection}
		
		Once the random flow is constructed, we can derive the expression of the associated Dirichlet energy.
		\begin{lemma}\label{lem.flow_intersection.EIT}
			For the unit flow $\theta$ given in \eqref{eq.def.flow_EIT}, its Dirichlet energy satisfies the following expression  
			\begin{align}\label{eq.flow_intersection.EIT}
				\mathscr{E}_*(\Lambda_{L-\ell_{0}/2}, \a,\theta) = \frac{1}{W^2} \sum_{\Gamma, \Gamma' \in \mathscr{G}_{M, \ell_{0}}} \mu\otimes\mu\big(\Gamma, \Gamma'\big)   \Eb\Ll[ \omega_\Gamma \omega_{\Gamma'} \zeta_\Gamma \zeta'_{\Gamma'} \vert \gamma \cap \gamma'\vert\, \vert \a\Rr].   
			\end{align}
			Here $(\sigma, U), (\sigma', U')$ are two independent copies living in $\PP$, and $\gamma, \gamma'$ are the random paths determined respectively by $\Big(\Gamma, \expl(\Gamma, \a; \sigma, U)\Big)$ and $\Big(\Gamma', \expl(\Gamma, \a; \sigma', U')\Big)$. $\zeta'_{\Gamma'}$ is a version of \eqref{eq.weight} associated to $ (\sigma', U')$.
		\end{lemma}
		\begin{proof}
			This is a direct application of Lemma~\ref{lem.flow}. Notice that $\gamma$ is sampled following the weight in \eqref{eq.weight.couople}
			\begin{align*}
				\frac{1}{W}  \omega_\Gamma\, \zeta_\Gamma\,  \d (\mu \otimes \PP_{\sigma, U}),
			\end{align*}
			so we double all the random elements except $\a$.
		\end{proof}
		
		As we have seen, the EIT property in Proposition~\ref{prop.EIT} has a factor $M^{2-d}$ as the mean-field mixing property of $(\Gamma, \Gamma')$. We still need a similar mixing property for the scale below  $\ell_{0}$; see Remark~\ref{rmk.heurictic.mixing}. This is the difference between the present work and \cite{BPP98}. A main result in this subsection is about the term $\Eb[\zeta_\Gamma \zeta'_{\Gamma'} \vert \gamma \cap \gamma'\vert \, \vert \G]$.
		
		\begin{prop}\label{prop.weight_intersection}
			Under the same setting as Lemma~\ref{lem.flow_intersection.EIT}, there exist two finite positive constants $C_{\ref{eq.bound.weight}}, K_{\ref{eq.bound.weight_intersection}}$ depending on $\epsilon_0, d$, such that the following two estimates hold
			\begin{align}
				\Eb\Ll[\zeta_\Gamma \zeta'_{\Gamma'}  \, \vert \G\Rr] &\leq C_{\ref{eq.bound.weight}} \Ll(1 + \epsilon_0\Rr)^{\vert \Gamma \cap \Gamma' \vert}, \label{eq.bound.weight}\\
				\Eb\Ll[\zeta_\Gamma \zeta'_{\Gamma'} \vert \gamma \cap \gamma'\vert \, \vert \G\Rr] &\leq K_{\ref{eq.bound.weight_intersection}} \ell_{0}^{8-d} \cdot \vert \Gamma \cap \Gamma' \vert \cdot \Ll(1 + \epsilon_0\Rr)^{\vert \Gamma \cap \Gamma' \vert}. \label{eq.bound.weight_intersection}
			\end{align}
		\end{prop}
		
		\subsubsection{Analysis of coalescence}
		Let us simplify at first the problem in Proposition~\ref{prop.weight_intersection}, and consider two independent exploration along a same route 
		\begin{align*}
			\Gamma = \Gamma' = I = \{z_0,z_1, \cdots \}.
		\end{align*}
		For the random variables $X$, we also use the superscript $X'$ to represent its version in $\expl(I, \a; \sigma', U')$. We observe that, the intersection can be reduced to the coalescence.
		
		\begin{lemma}[Coalescence between two exploration processes]\label{lem.coal.I}
			Given a route $I \subset Z_{M, \ell_{0}}$ and assume the pre-goodness along it, then define the coalescence between $\expl(I, \a; \sigma, U)$ and  $\expl(I, \a; \sigma', U')$ at the step $t$
			\begin{align}\label{eq.coal.I}
				\mathsf{E}_t := \{\CC_t \cap \CC'_t \neq \emptyset\} \cup \{\skl^{\sigma_t}(\CC_t,v_t) \cap \skl^{\sigma'_t}(\CC'_t,v_t) \neq \emptyset\}.
			\end{align}
			It satisfies that 
			\begin{align*}
				\vert \gamma \cap \gamma' \vert  \leq (K_{\mathrm{spn}}  + K_{\mathrm{skl}}) \ell_{0}^2 \cdot \#\{t: \mathsf{E}_t\}. 
			\end{align*}
		\end{lemma}
		\begin{proof}
			In the proof of \eqref{eq.path.length}, we have seen that the sprinkling paths and spanning paths are disjoint segments. Therefore, we have 
			\begin{align*}
				\vert \gamma \cap \gamma' \vert &\leq \sum_{t} \vert \gamma^{\spn}_t \cap \gamma^{',\spn}_t \vert +  \vert \gamma^{\skl}_t \cap \gamma^{',\skl}_t \vert\\
				&\leq (K_{\mathrm{spn}}  + K_{\mathrm{skl}}) \ell_{0}^2 \cdot  \#\{t: \gamma^{\spn}_t \cap \gamma^{',\spn}_t \neq \emptyset, \text{ or } \gamma^{\skl}_t \cap \gamma^{',\skl}_t \neq \emptyset \}.
			\end{align*}
			The second line is the result of \eqref{eq.path.length.spn_skl}. The coalescence is the necessary condition for the intersection of path (see \eqref{eq.gamma.spn} and \eqref{eq.gamma.skl}), 
			\begin{align*}
				\gamma^{\spn}_t &\subset \cC^*_t \subset \CC_t, \\
				\gamma^{\skl}_t &\subset \skl^{\sigma_t}(\cC^*_t,v_t) \subset \skl^{\sigma_t}(\CC_t,v_t),
			\end{align*}
			so we can further relax the inequality above and obtain the desired result. 
		\end{proof}

		We then study the coalescence. The heuristics to keep in mind says that, \emph{``in two independent exploration processes, it pays probability $\ell_0^{6-d}$ to coalesce, and every coalescence only lasts a length of exponential tail''}. All our effort is to turn this message rigorous in math. Several choices of parameters in Definition~\ref{def.spanning.z} and Definition~\ref{def.goodskl} will play the role in the proof. For the application in the future paragraphs, the statement also considers a variant case of exploration process (Case $\mathbf{II}$ below).
		
		\begin{prop}\label{prop.length.path.conditional}
			There exists a finite positive constant $C_{\ref{eq.length.path.conditional}}$ depending on $d$ and $\epsilon_0$, such that the following estimate 
			\begin{align}\label{eq.length.path.conditional}
				\Eb\Ll[\zeta_I \zeta'_{I} \#\{t: \mathsf{E}_t\} \, \vert \G\Rr] &\leq 	C_{\ref{eq.length.path.conditional}} \ell_{0}^{6-d} \cdot \vert I \vert \cdot \Ll(1 + \epsilon_0\Rr)^{\vert I \vert},
			\end{align}
			holds under the same setting as Lemma~\ref{lem.coal.I} and for $(\CC_0$, $\CC'_0)$ sampled under two cases:
			\begin{itemize}
				\item Case $\mathbf{I}$:  $\CC_0$, $\CC'_0$ are independent copies of $\kk$-nests uniformly sampled in $\mathscr{S}(z_0)$.
				\item Case $\mathbf{II}$:  $z_0$ has two different neighbor sites $z_{-1}, z'_{-1} \in Z_{M, \ell_{0}} \setminus I$, and there live two $\kk$-nests $\CC_{-1}, \CC'_{-1}$. The $\kk$-nests $\CC_0$, $\CC'_0$ are respectively obtained by the sprinkling from $\skl^{\sigma_{-1}}(\CC_{-1},v_{-1})$ and  $\skl^{\sigma'_{-1}}(\CC'_{-1},v'_{-1})$.
			\end{itemize}
		\end{prop}
		\begin{proof}
			Throughout the proof, we exceptionally use $A \lesssim B$ for $A \leq C B$ with a constant $C$ depending on $\epsilon_0, d$. We use the shorthand notations 
			\begin{align*}
				\skl^{\sigma_t}(\CC_t) \equiv \skl^{\sigma_t}(\CC_t,v_t), \qquad \qquad		\skl^{\sigma'_t}(\CC'_t) \equiv \skl^{\sigma'_t}(\CC'_t,v_t), 
			\end{align*}
			and 
			\begin{align}\label{eq.short.skl_typical}
				L_t := \{\skl^{\sigma_t}(\CC_t) \text{ is typical}\}, \qquad \qquad L'_t := \{\skl^{\sigma'_t}(\CC'_t) \text{ is typical}\},
			\end{align}
			when the direction is clear. 
			
			Using the bound $s_t \geq 1 -\epsilon_0/4$ in \eqref{eq.st.bound}, we just drop the denominators in  $\zeta_I \zeta'_{I}$, and obtain a point-wise bound 
			\begin{align*}
				\zeta_I \zeta'_{I} \#\{t: \mathsf{E}_t\} &\leq \Ll(\prod_{t} \Ind{L_t \cap L'_t} \Rr) \Ll(\sum_{t} \Ind{\mathsf{E}_t} \Rr)  (1 -\epsilon_0/4)^{-2 \vert I \vert} \\
				&\leq \sum_{t}\Ll(\Ind{\mathsf{E}_t} \prod_{u \leq t} \Ind{L_u \cap L'_u} \Rr) (1 + \epsilon_0)^{\vert I \vert}. 
			\end{align*}
			We further use a random variable $X_t$ to indicate the state of coalescence and non-coalescence
			\begin{align}\label{eq.def.X}
				\begin{split}
					\{X_t = 1\} &:= \mathsf{E}_t \cap \bigcap_{u \leq t} (L_u \cap L'_u), \\
					\{X_t = 0\} &:= \mathsf{E}^c_t \cap \bigcap_{u \leq t} (L_u \cap L'_u), \\
					\{X_t = -1\} &:= \bigcup_{u \leq t} L^c_u \cup (L'_u)^c. \\
				\end{split}
			\end{align}
			Then it suffices to justify that 
			\begin{align}\label{eq.length.path.conditional.reduction}
				\Eb\Ll[\sum_{t} \Ind{X_t = 1} \, \vert \G\Rr] \leq 	C_{\ref{eq.length.path.conditional}} \ell_{0}^{6-d} \cdot \vert I \vert.
			\end{align}
			This requires a detailed study of the probability of coalescence, which can be dominated by a Markov chain. We define the probability 
			\begin{align*}
				\mathbf{Q}[(\, \cdot \, )] := \PP\Ll[ (\, \cdot \, )  \Big\vert \G \Rr],
			\end{align*}
			and the filtration 
			\begin{align*}
				\G_t := \sigma \Ll( (X_u)_{0 \leq u \leq t}\Rr).
			\end{align*}
			We claim that, 
			there exists a constant $\delta(\epsilon_0) \in (0,1)$, such that the following transition probability holds:
			\begin{align}
				\mathbf{Q}[X_0 = 1] &\lesssim \ell_0^{6-d}. \label{eq.Markov.0}\\
				\mathbf{Q}[X_{t+1} = 1 \vert \G_t] &\leq (1-\delta) \Ind{X_t \geq 0}. 	 \label{eq.Markov.escape}\\
				\mathbf{Q}[X_{t+1} = 1 \vert X_t = 0, \G_{t-1}] &\lesssim \ell_0^{6-d}. \label{eq.Markov.merge}
			\end{align}
			Here the equation \eqref{eq.Markov.merge} is understood as 
			\begin{align}\label{eq.Markov.merge.def}
				\mathbf{Q}[X_{t+1} = 1 \vert X_t = 0, \G_{t-1}] := \frac{\mathbf{Q}[X_{t+1} = 1 , X_t = 0 \, \vert \G_{t-1}] }{\mathbf{Q}[X_t = 0 \, \vert \G_{t-1}] },
			\end{align}
			and $\{X_t = 0\}$ implies that $\{X_u \geq 0\}$ for all $u \leq t-1$.

			The transition probability is the heart of the proof, because it concludes \eqref{eq.length.path.conditional.reduction}. Heuristically, \eqref{eq.Markov.0}, \eqref{eq.Markov.escape}, \eqref{eq.Markov.merge} suggest a stochastic domination by a Poisson point process of intervals, whose density is of order $\ell_0^{6-d}$, and the length independently has a geometric decay of the parameter $\delta$. That is, we consider the random variable $N_t$ as the number of rush of $1$ in the sequence of $X$ starting from $t$
			\begin{align*}
				N_0 &:= \max\{k \in \N: X_0 = 1, X_{1} = 1, \cdots, X_{k-1} = 1\},\\
				\forall t \in \N_+, \quad N_t &:= \max\{k \in \N: X_{t-1}=0, X_t = 1, X_{t+1} = 1, \cdots, X_{t+k-1} = 1\}.
			\end{align*}
			Clearly, the following identity holds
			\begin{align*}
				\sum_{t} \Ind{X_t = 1} = \sum_{t} N_t.
			\end{align*}
			We have a tail estimate for $N_t, t \geq 1$
			\begin{align}\label{eq.N.tail.1}
				\begin{split}
					\mathbf{Q}[N_t \geq k] &= \mathbf{Q}[X_{t-1}=0, X_t = 1, X_{t+1} = 1, \cdots, X_{t+k-1} = 1]\\
					&=\mathbf{Q}[X_{t-1}=0] \, \cdot \mathbf{Q}[X_t = 1 \, \vert X_{t-1}=0] \\
					&\qquad \cdots \cdot \mathbf{Q}[ X_{t+k-1} = 1 \, \vert X_{t+k-2} = 1, \cdots , X_t = 1, X_{t-1}=0]\\
					&\lesssim \ell_{0}^{6-d} (1-\delta)^{k-1}.
				\end{split}
			\end{align}
			Here we apply \eqref{eq.Markov.merge} to the term $\mathbf{Q}[X_t = 1 \, \vert X_{t-1}=0]$, and \eqref{eq.Markov.escape} to the other conditional probabilities. There is also a similar estimate for $N_0$
			\begin{align}\label{eq.N.tail.0}
				\begin{split}
					\mathbf{Q}[N_0 \geq k] &= \mathbf{Q}[X_0 = 1, X_{1} = 1, \cdots, X_{k-1} = 1]\\
					&=\mathbf{Q}[X_0 = 1] \, \cdot \mathbf{Q}[X_1 = 1 \, \vert X_0 = 1] \\
					&\qquad \cdots \cdot \mathbf{Q}[ X_{k-1} = 1 \, \vert X_{k-2} = 1, \cdots , X_1 = 1, X_0=1]\\
					&\lesssim \ell_{0}^{6-d} (1-\delta)^{k-1}.
				\end{split}
			\end{align}
			The only difference is an application of \eqref{eq.Markov.0} to the term $\mathbf{Q}[X_0 = 1]$. Therefore, the geometric decay in \eqref{eq.N.tail.0} and \eqref{eq.N.tail.1} entails that
			\begin{align*}
				\forall t \in \N, \qquad  \Eb\Ll[N_t \, \vert \G\Rr] = \sum_{k = 1}^\infty \mathbf{Q}[N_0 \geq k] \lesssim \sum_{k = 1}^\infty \ell_{0}^{6-d} (1-\delta)^{k-1}\lesssim \delta^{-1} \cdot \ell_{0}^{6-d},
			\end{align*}
			and 
			\begin{align*}
				\Eb\Ll[\sum_{t} \Ind{X_t = 1} \, \vert \G\Rr] = \Eb\Ll[\sum_{t} N_t \, \vert \G\Rr] \lesssim \delta^{-1} \cdot \ell_0^{6-d} \cdot \vert I \vert.
			\end{align*}
			This verifies \eqref{eq.length.path.conditional.reduction}, so it finishes the proof of \eqref{eq.length.path.conditional}.

			\medskip

			We then prove \eqref{eq.Markov.0}, \eqref{eq.Markov.escape}, and \eqref{eq.Markov.merge} in three steps. Recall the definition \eqref{eq.def.S_child} and \eqref{eq.def.V_coal}, we define 
			\begin{align}\label{eq.def.S_coal}
				\mathscr{S}_{\mathrm{coal}}(z_t) := \{\cC'' \in \mathscr{S}(z_t): \CC_{t} \xleftrightarrow{\skl^{\sigma_t}(\CC_t)} \cC''\},
			\end{align}
			and write the shorthand notation of offspring spanning clusters 
			\begin{align}\label{eq.def.S_child_short}
				\mathscr{S}_{\mathrm{off}}(z_{t+1}) &\equiv \mathscr{S}_{\mathrm{off}}(\CC_{t}, \overrightarrow{z_t z_{t+1}}).
			\end{align}
			Notice that $\mathscr{S}_{\mathrm{off}}(z_{t+1})$ is actually the set of offspring clusters of $\CC_{t}$, but we choose the index $t+1$ to make the notation consistent with $\mathscr{S}_{\mathrm{off}}(z_{t+1}) \subset \mathscr{S}(z_{t+1})$.

			\medskip
			\textit{Step~1: estimate of $\mathbf{Q}[X_0 = 1]$.} As the statement contains two cases, we study them separately. Actually, both cases share the same spirit, as they pose a mixing property to the initial condition.
			
			\smallskip
			\textit{Case~$\mathbf{I}$.} The event $\{X_0 = 1\}$ can be dominated by 
			\begin{align}\label{eq.A0.I.decom}
				\mathbf{Q}[X_0 = 1] \leq \mathbf{Q}[\CC_0 \cap \CC'_0 \neq \emptyset] +  \mathbf{Q}[\CC_0 \cap \CC'_0 = \emptyset, \skl(\CC_0) \cap \skl(\CC'_0) \neq \emptyset].
			\end{align}
			
			For the first term, because $\CC_0$ and $\CC'_0$ are chosen independently among the $\kk$-nest, we have
			\begin{align}\label{eq.A0.I.decom.1}
				\mathbf{Q}[\CC_0 \cap \CC'_0 \neq \emptyset] = \frac{\sum_{j=1}^{\kk} {\# \mathscr{S}(z_{0}) \choose \kk } {\kk \choose j} {{\# \mathscr{S}(z_{0}) - \kk} \choose {\kk-j} }}{{\# \mathscr{S}(z_{0}) \choose \kk }^2} \lesssim \frac{1}{\# \mathscr{S}(z_{0})} \lesssim \ell_{0}^{6-d}.
			\end{align}
			Actually, the denominator has an order $(\# \mathscr{S}(z_{0}))^{2\kk}$, while the nominator only has an order $(\# \mathscr{S}(z_{0}))^{2\kk - 1}$. The last inequality uses the property of pre-good cube in \eqref{eq.pregood}.
			
			For the second term, we notice that 
			\begin{align*}
				\Ind{\skl(\CC_0) \cap \skl(\CC'_0) \neq \emptyset} \leq \sum_{\cC' \in \mathscr{S}(z_0)} \Ind{\skl(\CC_0) \cap \skl(\cC') \neq \emptyset, \cC' \subset \CC'_0}.
			\end{align*}
			Using the first moment estimate, we make the decomposition, and obtain that 
			\begin{multline}\label{eq.A0.I.decom.2.pre}
				\mathbf{Q}[\CC_0 \cap \CC'_0 = \emptyset, \skl(\CC_0) \cap \skl(\CC'_0) \neq \emptyset]\\
				\leq \sum_{\substack{S \text{ is } \kk\text{-nests}}} \mathbf{Q}[\CC_0 = S] \sum_{\substack{\cC' \in \mathscr{S}(z_0) \\ \cC' \cap S = \emptyset}}  \mathbf{Q}[\cC' \subset \CC'_0] \mathbf{Q}[\skl(S) \cap \skl(\cC') \neq \emptyset].
			\end{multline}
			The formula also relies on the independence choice between $\CC_0$ and $\CC'_0$. Since the choice of $ \CC'_0$ is uniform, we have
			\begin{align*}
				\mathbf{Q}[\cC' \subset \CC'_0] \leq \frac{\kk}{\# \mathscr{S}(z_{0})} \lesssim \ell_{0}^{6-d}.
			\end{align*}
			The uniform bound allows us to do the summation
			\begin{align*}
				&\sum_{\substack{\cC' \in \mathscr{S}(z_0) \\ \cC' \cap S = \emptyset}} \mathbf{Q}[\cC' \subset \CC'_0] \mathbf{Q}[\skl(S) \cap \skl(\cC') \neq \emptyset]\\
				&\lesssim \ell_{0}^{6-d}\sum_{\substack{\cC' \in \mathscr{S}(z_0) \\ \cC' \cap S = \emptyset}} \mathbf{Q}[\skl(S) \cap \skl(\cC') \neq \emptyset] \\
				&\lesssim \ell_{0}^{6-d} \Ll(\ell_{0}^{-d+2} \min\{\kk C_{\mathrm{b}} \ell^{2}_{0}, C_{\mathrm{b}} C_{\ref{eq.pregood}} \ell^{d-4}_{0}\} \cdot C_{\mathrm{b}} C_{\ref{eq.pregood}} \ell^{d-4}_{0} +   \ell_{0}^{-2d+6}  \kk C_{\mathrm{b}} \ell^{2}_{0} \cdot \Ll(C_{\mathrm{b}} C_{\ref{eq.pregood}} \ell^{d-4}_{0}\Rr)^2\Rr)\\
				&\lesssim \ell_{0}^{6-d}.
			\end{align*}
			In the last line, we use the estimate of coalescence \eqref{eq.AAB} and the regular boundary of typical spanning clusters in Definition~\ref{def.spanning.z}. We put this estimate back to \eqref{eq.A0.I.decom.2.pre}, and conclude that 
			\begin{align}\label{eq.A0.I.decom.2}
				\mathbf{Q}[\CC_0 \cap \CC'_0 = \emptyset, \skl(\CC_0) \cap \skl(\CC'_0) \neq \emptyset] \lesssim \ell_{0}^{6-d}.
			\end{align}
			
			Finally, viewing \eqref{eq.A0.I.decom}, then the bounds \eqref{eq.A0.I.decom.1}, \eqref{eq.A0.I.decom.2} together justify \eqref{eq.Markov.0} for Case~$\mathbf{I}$.
			
			\smallskip
			\textit{Case~$\mathbf{II}$.}
			Let us reveal at first $(\CC_0, \skl^{\sigma_0}(\CC_0))$. The event $\mathsf{E}_0$ has an inclusion that 
			\begin{align}\label{eq.A.inclusion1}
				\mathsf{E}_0 \subset \{\mathscr{S}'_{\mathrm{off}}(z_0) \cap \mathscr{S}_{\mathrm{off}}(z_0) \neq \emptyset\} \cup \{\mathscr{S}'_{\mathrm{off}}(z_0) \cap \mathscr{S}_{\mathrm{coal}}(z_0) \neq \emptyset\}.
			\end{align}
			See \eqref{eq.def.S_coal} and \eqref{eq.def.S_child_short} for the notations. Let us explain this identity. As the $\kk$-nest $\CC_0,\CC'_{0}$ are chosen respectively among $\mathscr{S}_{\mathrm{off}}(z_0)$ and $\mathscr{S}'_{\mathrm{off}}(z_0)$, the first subset in \eqref{eq.A.inclusion1} contains the case that $\CC_0$ intersects $\CC'_{0}$. Moreover, viewing the definition \eqref{eq.def.V_coal} and \eqref{eq.def.S_coal}, only the sprinkling clusters issued from $\mathscr{S}_{\mathrm{coal}}(z_0)$ can intersect $\skl^{\sigma_0}(\CC_0)$, which is contained by the second subset in \eqref{eq.A.inclusion1}. Therefore, the inclusion in \eqref{eq.A.inclusion1} is valid.

			The inclusion \eqref{eq.A.inclusion1} can go further as 
			\begin{align}\label{eq.A.inclusion2}
				\mathsf{E}_0 \subset \Ll\{\pivo\Ll(\CC'_{-1}, \bigcup_{\cC'' \in \mathscr{S}_{\mathrm{off}}(z_0)} \cC'' \cup \bigcup_{\cC''' \in \mathscr{S}_{\mathrm{coal}}(z_0)} \cC'''\Rr) \neq \emptyset  \Rr\}.
			\end{align}
			Clearly, if the {\rhs} of \eqref{eq.A.inclusion1} holds, there should be open pivotal edges between the two sets, which is \eqref{eq.A.inclusion2}.

			We now start to calculate the upper bound of the probability. To simplify the notation, we note 
			\begin{align}\label{eq.def.S_uion}
				\forall S \subset \mathscr{S}(z_0), \qquad \widehat{S} := \bigcup_{\cC'' \in S} \cC''.
			\end{align} 
			Using \eqref{eq.A.inclusion2}, the probability of $\{X_0 = 1\}$ can thus be written as 
			\begin{align*}
				\{X_0 = 1\} &= L_{-1} \cap L_0 \cap L'_{-1} \cap L'_0 \cap \mathsf{E}_0 \\
				&\subset L_{-1} \cap L_0 \cap  \Ll\{\pivo\Ll(\CC'_{-1}, \widehat{\mathscr{S}_{\mathrm{off}}(z_0)} \cup \widehat{\mathscr{S}_{\mathrm{coal}}(z_0)}\Rr) \neq \emptyset  \Rr\}.
			\end{align*} 
			We drop the terms $L'_0, L'_{-1}$ here. As there is still randomness in $\widehat{\mathscr{S}_{\mathrm{off}}(z_0)} \cup \widehat{\mathscr{S}_{\mathrm{coal}}(z_0)}$, we make a further decomposition in function of the configuration of offspring clusters and coalescent clusters
			\begin{align}\label{eq.Markov.0.decom}
				\begin{split}
					\mathbf{Q}[X_0 = 1]  &\leq \sum_{\substack{S_1 \subset \mathscr{S}(z_0) \\ m_{\mathrm{off}}  \leq  \# S_1 \leq M_{\mathrm{off}} }} \sum_{\substack{S_2 \subset \mathscr{S}(z_0) \\ \# S_2 \leq M_{\mathrm{coal}} }} \\
					& \qquad \mathbf{Q}\Ll[ \Ll\{\pivo\Ll(\CC'_{-1}, \widehat{S_1} \cup \widehat{S_2} \Rr) \neq \emptyset  \Rr\} \cap L_{-1} \cap L_0 \cap \{\mathscr{S}_{\mathrm{off}}(z_0) = S_1, \mathscr{S}_{\mathrm{coal}}(z_0) = S_2\}\Rr].
				\end{split}
			\end{align}
			The condition of typical sprinkling clusters then enters; see Definition~\ref{def.goodskl} and Proposition~\ref{prop.proba.k-nest}. Notice that the event
			\begin{align*}
				\Ll\{\pivo\Ll(\CC'_{-1}, \widehat{S_1} \cup \widehat{S_2} \Rr) \neq \emptyset  \Rr\},
			\end{align*}
			depends on the percolation  in the tube $T(x'_{-1},z_0)$, while the event
			\begin{align*}
				L_{-1} \cap L_0 \cap \{\mathscr{S}_{\mathrm{off}}(z_0) = S_1, \mathscr{S}_{\mathrm{coal}}(z_0) = S_2\},
			\end{align*}
			depends on the percolation in $T(z_{-1},z_0) \cup T(z_0,z_1)$. By the assumption $z_{-1}, z'_{-1} \in Z_{M, \ell_{0}} \setminus I$ and Lemma~\ref{lem.skl.independent}, the two events are independent, so we obtain
			\begin{multline}\label{eq.Markov.0.decom.plus}
				\mathbf{Q}[X_0 = 1] \leq  \sum_{\substack{S_1 \subset \mathscr{S}(z_0) \\ m_{\mathrm{off}}  \leq  \# S_1 \leq M_{\mathrm{off}} }} \sum_{\substack{S_2 \subset \mathscr{S}(z_0) \\ \# S_2 \leq M_{\mathrm{coal}} }}
				\mathbf{Q}\Ll[\Ll\{\pivo\Ll(\CC'_{-1}, \widehat{S_1} \cup \widehat{S_2} \Rr) \neq \emptyset  \Rr\}\Rr] \\ 
				\cdot \mathbf{Q}\Ll[ L_{-1} \cap L_0 \cap \{\mathscr{S}_{\mathrm{off}}(z_0) = S_1, \mathscr{S}_{\mathrm{coal}}(z_0) = S_2\} \Rr].
			\end{multline} 
			Using the first moment estimate in \eqref{eq.sprinkling.AB}
			\begin{align*}
				&\mathbf{Q}\Ll[ \Ll\{\pivo\Ll(\CC'_{-1}, \widehat{S_1} \cup \widehat{S_2} \Rr) \neq \emptyset  \Rr\} \Rr]  \\
				&\lesssim \ell_{0}^{-2} \cdot  \ell_{0}^{-d+4} (\kk C_{\mathrm{b}} \ell_{0}^2)  \cdot ((M_{\mathrm{off}} + M_{\mathrm{coal}})C_{\mathrm{b}} \ell_{0}^2)\\
				&\lesssim \ell_{0}^{6-d},
			\end{align*}
			we obtain the upper bound that 
			\begin{equation*}
				\begin{split}
					\mathbf{Q}[X_0 = 1]	& \lesssim \sum_{\substack{S_1 \subset \mathscr{S}(z_0) \\ m_{\mathrm{off}}  \leq  \# S_1 \leq M_{\mathrm{off}} }} \sum_{\substack{S_2 \subset \mathscr{S}(z_0) \\ \# S_2 \leq M_{\mathrm{coal}} }} \ell_{0}^{6-d} \mathbf{Q}\Ll[ L_{-1} \cap L_0 \cap \{\mathscr{S}_{\mathrm{off}}(z_0) = S_1, \mathscr{S}_{\mathrm{coal}}(z_0) = S_2\} \Rr]\\
					&\lesssim  \sum_{\substack{S_1 \subset \mathscr{S}(z_0)}} \sum_{\substack{S_2 \subset \mathscr{S}(z_0)}} \ell_{0}^{6-d} \mathbf{Q}\Ll[  \{\mathscr{S}_{\mathrm{off}}(z_0) = S_1, \mathscr{S}_{\mathrm{coal}}(z_0) = S_2\} \Rr] \\
					&\lesssim \ell_{0}^{6-d}.
				\end{split}
			\end{equation*}
			Then \eqref{eq.Markov.0} for Case~$\mathbf{II}$ is justified.

			\begin{figure}[h!]
				\centering
				\includegraphics[width=0.6\textwidth]{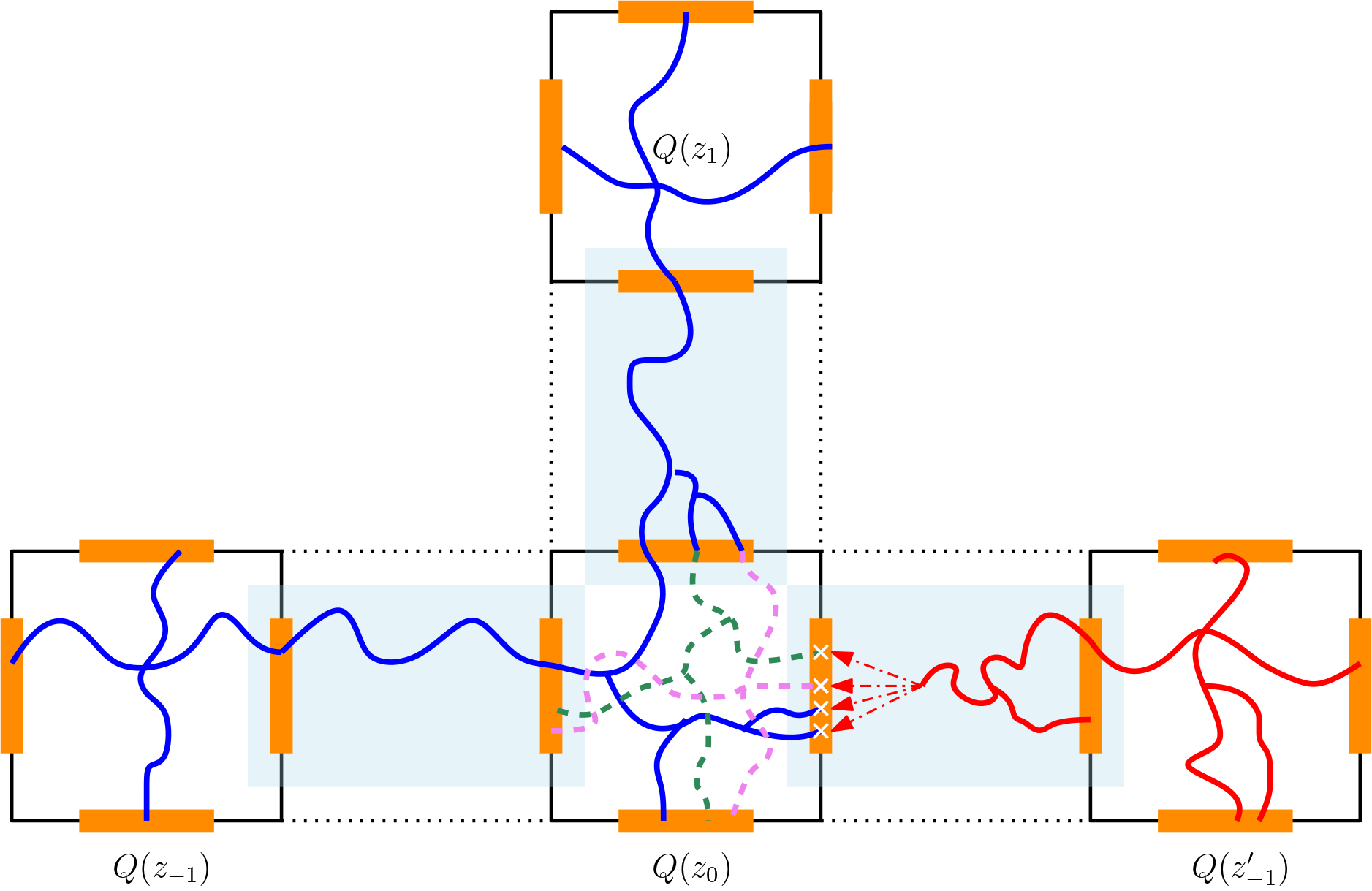}
				\caption{An illustration of the proof for Case~$\mathbf{II}$. The blue cluster is the one revealed in $\expl(I, \a; \sigma, U)$, while the red is the one in $\expl(I, \a; \sigma', U')$. We run the first exploration process, then the intersection occurs only when the second exploration touches certain part of the boundary, which is marked in white crosses.}\label{fig.A0_case2}
			\end{figure}

			\smallskip
			\textit{Step~2: estimate of $\mathbf{Q}[X_{t+1} = 1 \vert \G_t]$.} Actually, the following bound holds for some $\delta(\epsilon_0) \in (0,1)$ 
			\begin{align}\label{eq.Markov.escape.lower}
				\mathbf{Q}[X_{t+1} = 0 \vert \G_t] \geq \delta \Ind{X_t \geq 0}.
			\end{align}
			This means, whatever the state, it always has a probability $\delta$ to jump to the disjoint clusters state at next step. This probability can be written as product
			\begin{align}\label{eq.Markov.escape.decom}
				\begin{split}
					\mathbf{Q}[X_{t+1} = 0 \vert \G_t] &= \mathbf{Q}[G_{t+1} \vert \G_t] \cdot  \mathbf{Q}[ F_{t+1} \vert G_{t+1}, \G_t], \\
					G_{t+1} &:= \{\CC'_{t+1} \cap \CC_{t+1} = \emptyset \},\\
					F_{t+1} &:= L_{t+1} \cap L'_{t+1} \cap \{\skl^{\sigma'_{t+1}}(\CC'_{t+1}) \cap \skl^{\sigma_{t+1}}(\CC_{t+1}) = \emptyset\}.
				\end{split} 
			\end{align}
			The events indicate that the two spanning or sprinkling clusters exist, but they do not intersect. We then estimate the two terms, and the key is to use the independence between $\expl(I,\a; \sigma, U)$ and $\expl(I,\a; \sigma', U')$. For convenience, we just focus on the meaningful case $\{X_t \geq 0\}$.
			
			A first key estimate is the following conditional probability
			\begin{align}\label{eq.Markov.escape.key1}
				\mathbf{Q}[G_{t+1}   \vert \mathscr{S}_{\mathrm{off}}(z_{t+1}), \mathscr{S}'_{\mathrm{off}}(z_{t+1}), \CC_{t+1}, \G_t]
				\geq \frac{1}{{\# \mathscr{S}'_{\mathrm{off}}(z_{t+1}) \choose \kk}} \geq \Ll(\frac{1}{M_{\mathrm{off}}}\Rr)^{\kk}.
			\end{align}
			This comes from the additional randomness of the exploration process in Definition~\ref{def.expl}. In Step~2 there, we pick $\kk$ typical spanning clusters among all the offspring of $\cC^{*,'}_{t}$ to form the $\kk$-nest $\CC'_{t+1}$. Recall \eqref{eq.pigeonhole.k-nest}, the good cluster contains at least $2\kk$ offspring. Therefore, after excluding $ \CC_{t+1}$, there still remains at least $1$ possible configuration to form a $\kk$-nest disjoint from $\CC_{t+1}$, even when the worst case like $\CC'_{t} = \CC_{t}$ happens. We also use the fact that $U_{t+1}, U'_{t+1}$ are independent in two explorations, and ${\# \mathscr{S}'_{\mathrm{off}}(z_{t+1}) \leq M_{\mathrm{off}}}$ by Definition~\ref{def.goodskl}. The estimate \eqref{eq.Markov.escape.key1} then entails that 
			\begin{align}\label{eq.Markov.escape.decom1}
				\mathbf{Q}[G_{t+1} \vert \G_t] \geq \Ll(\frac{1}{M_{\mathrm{off}}}\Rr)^{\kk}.
			\end{align}
			
			Once we obtain disjoint $\CC'_{t+1}$ and $\CC_{t+1}$, we turn to $F_{t+1}$. A second key observation relies on Proposition~\ref{prop.coal}. If one inserts the typical order of clusters boundary in \eqref{eq.AAB}, as the calculation in \eqref{eq.Vcoal} shows, one can only show the probability of coalescence is $O(1)$, but cannot turn to a number strictly smaller than $1$. However,  \eqref{eq.AAB_plus} saves it. The price is to set the boundary of clusters separated with a distance $\ell_{0}$. That is the reason why we introduce (4) of Definition~\ref{def.spanning.z}, and include the randomness $\sigma$ in Definition~\ref{def.expl}.  Inspired by the discussion above, we then consider the following inclusion 
			\begin{align*}
				\begin{split}
					F_{t+1} &\supset F^{(0)}_{t+1} \cap F^{(1)}_{t+1} \cap F^{(2)}_{t+1},\\
					F^{(0)}_{t+1} &:= L_{t+1} \cap L'_{t+1} \cap \Ll\{ \skl^+(\CC_{t+1}), \skl^{-}(\CC'_{t+1}) \text{ are disjoint}\Rr\},\\
					F^{(1)}_{t+1} &:= \{ \sigma_{t+1} = +\},\\
					F^{(2)}_{t+1} &:= \{\sigma'_{t+1} = -\}.\\
				\end{split}
			\end{align*}
			Recall the notation of boundary in \eqref{eq.def.regular_bd_2} and \eqref{eq.def.SQH.ell0.bd.2}.
			Roughly, the {\rhs} above means $ \skl^+(\CC_{t+1})$ and $\skl^{-}(\CC'_{t+1})$ are the sprinkling clusters respectively generated from $\partial^\circ_{i, +} Q(z_{t+1})$ and $\partial^\circ_{i, -} Q(z_{t+1})$. The probability is then estimated as 
			\begin{align}\label{eq.Markov.escape.key2}
				\begin{split}
					\mathbf{Q}[ F_{t+1} \vert G_{t+1}, \G_t] &\geq \mathbf{Q}[ F^{(0)}_{t+1} \cap F^{(1)}_{t+1} \cap F^{(2)}_{t+1} \vert G_{t+1}, \G_t] \\
					&\geq \mathbf{Q}[ F^{(0)}_{t+1}   \vert G_{t+1}, \G_t] \cdot \mathbf{Q}[ F^{(1)}_{t+1}] \cdot \mathbf{Q}[ F^{(2)}_{t+1}].
				\end{split}
			\end{align}
			Here we also use the fact that the choices of $\sigma_{t+1}, \sigma'_{t+1}$ are independent, and are independent with other elements in the exploration. It remains to estimate the terms above one by one. 
			
			Concerning $\mathbf{Q}[ F^{(0)}_{t+1}   \vert G_{t+1}, \G_t]$, \eqref{eq.proba.k-nest} implies that
			\begin{align*}
				\mathbf{Q}[ F^{(0)}_{t+1}   \vert G_{t+1}, \G_t] &= \mathbf{Q}\Ll[ \Ll\{ \skl^+(\CC_{t+1}), \skl^{-}(\CC'_{t+1}) \text{ are typical}\Rr\}  \vert G_{t+1}, \G_t\Rr] \\
				&\qquad - \mathbf{Q}\Ll[ \Ll\{ \skl^+(\CC_{t+1}), \skl^{-}(\CC'_{t+1}) \text{ intersects}\Rr\}  \vert G_{t+1}, \G_t\Rr]\\
				&\geq 1 - \frac{\epsilon_0}{4} - \frac{\epsilon_0}{4} - \mathbf{Q}\Ll[ \Ll\{ \skl^+(\CC_{t+1}), \skl^{-}(\CC'_{t+1}) \text{ intersects}\Rr\}  \vert G_{t+1}, \G_t\Rr].
			\end{align*}
			For the event of intersection, \eqref{eq.def.SQH.ell0.bd.2} suggests that 
			\begin{align*}
				\dex\Ll(\PTR_{i, +} \CC_{t+1}, \PTR_{i, -} \CC'_{t+1}\Rr) \geq \frac{\ell_{0}}{8},  
			\end{align*}
			so \eqref{eq.AAB_plus} applies and yields
			\begin{align*}
				\mathbf{Q}\Ll[ \Ll\{ \skl^+(\CC_{t+1}), \skl^{-}(\CC'_{t+1}) \text{ intersects}\Rr\}  \vert G_{t+1}, \G_t\Rr] \lesssim \ell_{0}^{-2d + 6} \cdot (\ell_{0}^2)^2 \cdot \ell_{0}^{d-4} \lesssim \ell_{0}^{6-d}.
			\end{align*}
			This proves that 
			\begin{align}\label{eq.Markov.escape.key2.term0}
				\mathbf{Q}[ F^{(0)}_{t+1}   \vert G_{t+1}, \G_t] \gtrsim 1 - \frac{\epsilon_0}{2} - \ell_{0}^{6-d}.
			\end{align}
			
			The probability in $\mathbf{Q}[ F^{(1)}_{t+1} ], \mathbf{Q}[ F^{(2)}_{t+1} ]$ follows the setting in Definition~\ref{def.expl}, that $\sigma_{t+1}, \sigma'_{t+1}$ are independently and uniformly sampled in $\{+,-\}$ 
			\begin{align}\label{eq.Markov.escape.key2.term12}
				\mathbf{Q}[ F^{(1)}_{t+1} ] = \mathbf{Q}[ F^{(2)}_{t+1} ] =  \frac{1}{2}.
			\end{align}
			
			In conclusion, we have 
			\begin{align}\label{eq.Markov.escape.decom2}
				\mathbf{Q}[F_{t+1} \vert G_{t+1}] \geq \frac{1}{4}\Ll(1 - \frac{\epsilon_0}{2} - \ell_{0}^{6-d}\Rr). 
			\end{align}
			This also concludes \eqref{eq.Markov.escape} as 
			\begin{align*}
				\mathbf{Q}[X_{t+1} = 0 \vert \G_t] \geq \underbrace{\frac{1}{4}\Ll(1 - \frac{\epsilon_0}{2} - \ell_{0}^{6-d}\Rr)\Ll(\frac{1}{M_{\mathrm{off}}}\Rr)^{\kk}}_{=:\delta}.
			\end{align*}
			We set the {\rhs} as the $\delta$ in \eqref{eq.Markov.escape}.

			\smallskip
			\textit{Step~3: estimate of $\mathbf{Q}[X_{t+1} = 1 \vert X_t = 0, \G_{t-1}]$.} Recall the definition \eqref{eq.Markov.merge.def}
			\begin{align*}
				\mathbf{Q}[X_{t+1} = 1 \vert X_t = 0, \G_{t-1}] = \frac{\mathbf{Q}[X_{t+1} = 1, X_t = 0 \vert  \G_{t-1}]}{\mathbf{Q}[X_t = 0 \vert \G_{t-1}]}.
			\end{align*}
			The result \eqref{eq.Markov.escape.lower} implies that 
			\begin{align*}
				\mathbf{Q}[X_t = 0 \vert \G_{t-1}] \geq \delta,
			\end{align*}
			so it suffices to get an upper bound $\mathbf{Q}[X_{t+1} = 1, X_t = 0 \vert  \G_{t-1}]$.
			
			The estimate here is nearly as Step~1. We follow the notations and steps there until \eqref{eq.Markov.0.decom}, and obtain
			\begin{align*}
				&\mathbf{Q}[X_{t+1} = 1, X_t = 0 \vert  \G_{t-1}] \leq \sum_{\substack{S_1 \subset \mathscr{S}(z_{t+1}) \\ m_{\mathrm{off}}  \leq  \# S_1 \leq M_{\mathrm{off}} }} \sum_{\substack{S_2 \subset \mathscr{S}(z_{t+1}) \\ \# S_2 \leq M_{\mathrm{coal}} }} \\
				& \qquad \mathbf{Q}\Ll[  \Ll\{\pivo\Ll(\CC'_{t}, \widehat{S_1} \cup \widehat{S_2} \Rr) \neq \emptyset  \Rr\} \cap L_{t} \cap L_{t+1} \cap \{\mathscr{S}_{\mathrm{off}}(z_{t+1}) = S_1, \mathscr{S}_{\mathrm{coal}}(z_{t+1}) = S_2\} \cap \{X_t = 0\} \vert \G_{t-1}\Rr].
			\end{align*}
			
			Here $\{X_t = 0\}$, or $\mathsf{E}^c_t$ more precisely, plays a role. It implies that the events
			\begin{align*}
				\Ll\{\pivo\Ll(\CC'_{t}, \widehat{S_1} \cup \widehat{S_2} \Rr) \neq \emptyset  \Rr\} \text{ and } L_{t}, 
			\end{align*}
			are witnessed using disjoint percolation bonds in $T(z_t, z_{t+1})$. Therefore, we have
			\begin{align*}
				&\mathbf{Q}\Ll[ \Ll\{\pivo\Ll(\CC'_{t}, \widehat{S_1} \cup \widehat{S_2} \Rr) \neq \emptyset  \Rr\} \cap L_{t} \cap L_{t+1} \cap \{\mathscr{S}_{\mathrm{off}}(z_{t+1}) = S_1, \mathscr{S}_{\mathrm{coal}}(z_{t+1}) = S_2\} \cap \{X_t = 0\} \vert \G_{t-1}\Rr]\\
				&\leq \mathbf{Q}\Ll[ \Ll\{\pivo\Ll(\CC'_{t}, \widehat{S_1} \cup \widehat{S_2} \Rr) \neq \emptyset  \Rr\} \circ \Ll(L_{t} \cap L_{t+1} \cap \{\mathscr{S}_{\mathrm{off}}(z_{t+1}) = S_1, \mathscr{S}_{\mathrm{coal}}(z_{t+1}) = S_2\}\Rr) \vert \G_{t-1}\Rr]\\
				&\leq \mathbf{Q}\Ll[ \Ll\{\pivo\Ll(\CC'_{t}, \widehat{S_1} \cup \widehat{S_2} \Rr) \neq \emptyset  \Rr\} \Rr] \mathbf{Q}\Ll[L_{t} \cap L_{t+1} \cap \{\mathscr{S}_{\mathrm{off}}(z_{t+1}) = S_1, \mathscr{S}_{\mathrm{coal}}(z_{t+1}) = S_2\} \vert \G_{t-1}\Rr].
			\end{align*}
			Here the BK--Reimer inequality can be applied in the third line, because no events here rely on the bonds in $\G$ or $\G_{t-1}$. Then the structure is exactly like \eqref{eq.Markov.0.decom.plus}, so \eqref{eq.Markov.merge} can be justified.
			
			\begin{figure}[h!]
				\centering
				\includegraphics[width=0.5\textwidth]{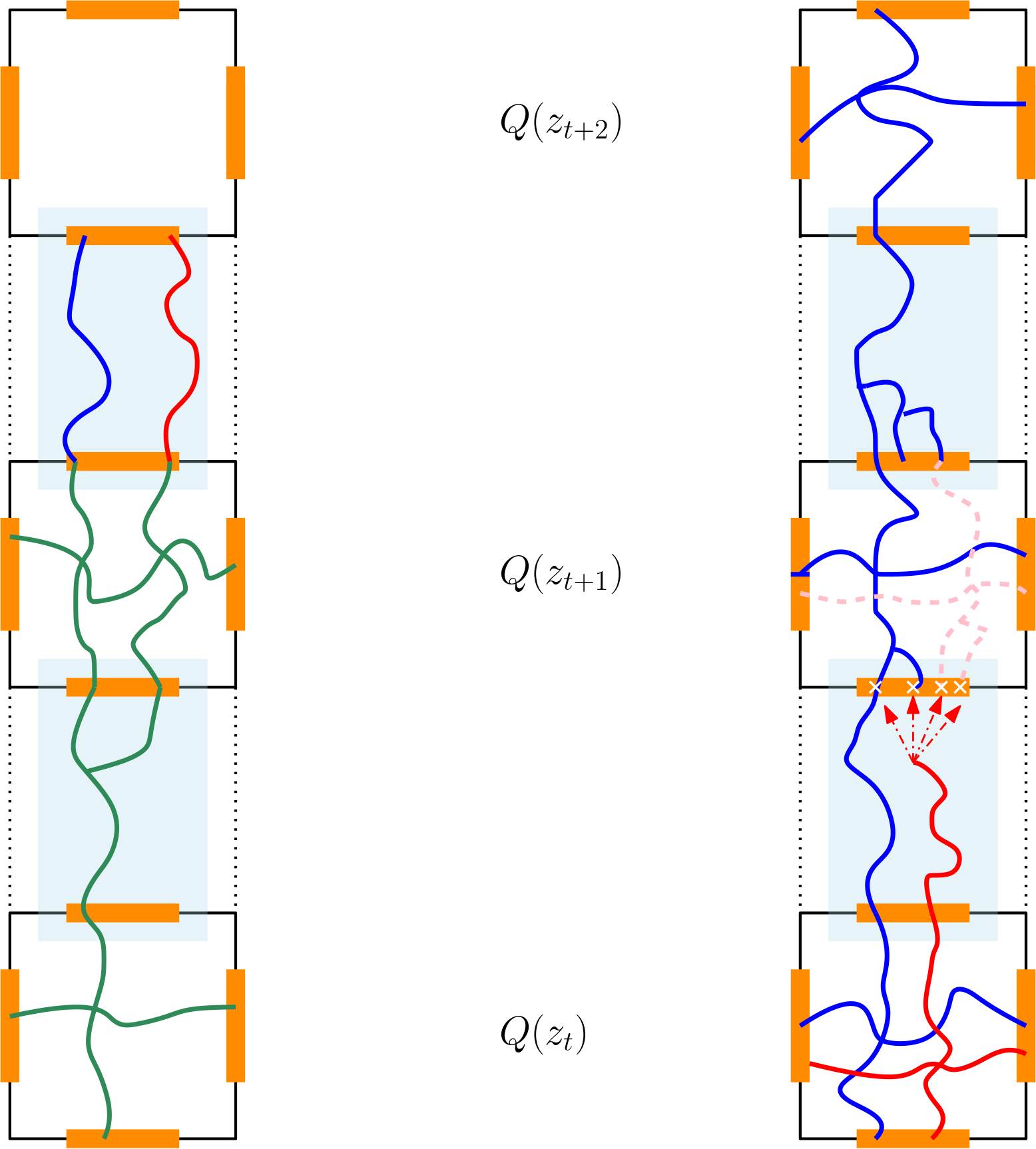}
				\caption{An illustration of the analysis in Step~2 and Step~3. The picture on the left is for $\mathbf{Q}[X_{t+1} = 1 \vert \G_t]$, which shows the repelling property between two independent explorations. Even the coalescence happens, it has a positive probability to end in next cube. The picture on the right is for ${\mathbf{Q}[X_{t+1} = 1 \vert X_t = 0, \G_{t-1}]}$, and the argument is similar as Case~$\mathbf{II}$ of  $\mathbf{Q}[X_{0} = 1]$}\label{fig.two_expl}
			\end{figure}

		\end{proof}

		\subsubsection{Weighted intersection between flows}
		
		We now extend the result to the general cases.
		\begin{proof}[Proof of Proposition~\ref{prop.weight_intersection}]
			The proof is divided into $3$ parts. We keep using the shorthand notations in the proof of Proposition~\ref{prop.length.path.conditional}.
			
			\textit{Step~1: decomposition of $\Gamma \cap \Gamma'$.} The key is to decompose $\Gamma \cap \Gamma'$ into disjoint segments $(I_j)_{j \in \N}$, and every segment $I_j$ is a common oriented route of the two explorations. Recall that the two oriented routes $\Gamma, \Gamma'$ are indexed as 
			\begin{align*}
				\Gamma  &= \{z_0, z_1, z_2, z_3, \cdots\},\\
				\Gamma'  &= \{z'_0, z'_1, z'_2, z'_3, \cdots\}.
			\end{align*}
			We also denote by $\operatorname{succ}(\Gamma,x)$ and $\operatorname{succ}(\Gamma',x)$ respectively the successor of the $x$ in $\Gamma$ and $\Gamma'$. Viewing the definition above, we naturally have 
			\begin{align*}
				\operatorname{succ}(\Gamma,z_t) = z_{t+1}.
			\end{align*}
			If $z_t \in \Gamma \cap \Gamma'$, we denote by $\phi(t)$ its associated index in $\Gamma'$, which implies the identity
			\begin{align*}
				z'_{\phi(t)} = z_t,
			\end{align*}
			and 
			\begin{align*}
				\operatorname{succ}(\Gamma',z_t) = \operatorname{succ}(\Gamma',z'_{\phi(t)}) = z'_{\phi(t)+1}.
			\end{align*}
			Then whether the two explorations follow a common route depends on the state of $z_{t+1}$ and $z'_{\phi(t)+1}$. We thus define the stopping times 
			\begin{align}\label{eq.def.tau}
				\begin{split}
					\tau_{-1}^* &:= -1, \\
					\forall j \in \N, \qquad \tau_j &:= \inf\{t > \tau_{j-1}^*: z_t \in \Gamma \cap \Gamma'\},\\
					\tau_j^* &:= \inf \{t \geq \tau_j: z_{t+1} \neq \operatorname{succ}(\Gamma', z_{t})\}.
				\end{split}
			\end{align}
			Then the intersection $\Gamma \cap \Gamma'$ into disjoint union:
			\begin{align}\label{eq.def.I}
				\Gamma \cap \Gamma' := \bigsqcup_{j = 0}^n I_j, \qquad I_j = [\tau_{j}, \tau_j^*] \cap \N.
			\end{align}
			Here $n$ stands for the total number of common routes.  Moreover, $I_j$ is not only a consecutive intersection, but common oriented route. This is the consequence of Definition~\ref{def.route.o}, which ensures that every tube contains a unique admissible direction, and the two processes explore every $I_j$ in the same order. Precisely, in the route $\Gamma'$, the segment $I_j$ is explored in the order $z'_{\phi(\tau_{j})}, z'_{\phi(\tau_{j}+1)}, \cdots, z'_{\phi(\tau_j^*)}$.
			
			\smallskip
			\textit{Step~2: proof of \eqref{eq.bound.weight}.} The proof relies on an induction. We use the shorthand notations in \eqref{eq.short.skl_typical}, then 
			\begin{align*}
				\zeta_\Gamma = \prod_{\substack{z_t \in \Gamma \\ t < \vert \Gamma \vert - 1}} \frac{\Ind{L_t}}{s_t}, \qquad \zeta'_{\Gamma'} = \prod_{\substack{z'_t \in \Gamma' \\ t < \vert \Gamma' \vert - 1}} \frac{\Ind{L'_t}}{s'_t}.
			\end{align*}
			We further define their truncated version as 
			\begin{align*}
				\zeta_t = \prod_{0 \leq u < t} \frac{\Ind{L_u}}{s_u}, \qquad \zeta'_{t} = \prod_{0 \leq u < t} \frac{\Ind{L'_u}}{s'_u},
			\end{align*}
			and the sigma-algebra along the exploration
			\begin{align*}
				\mathcal{F}_{t} := \sigma \Ll( \G,  (\CC_u, \sigma_u)_{0 \leq u \leq t}, \a(e)_{e \in \Ed(T(z_{u-1}, z_u)), 1 \leq u \leq t} \Rr), \\
				\mathcal{F}'_{t} := \sigma \Ll( \G,  (\CC'_u, \sigma'_u)_{0 \leq u \leq t}, \a(e)_{e \in \Ed(T(z'_{u-1}, z'_u)), 1 \leq u \leq t} \Rr).
			\end{align*}
			In another word, $\mathcal{F}_{t}$ contains the information exactly before the sprinkling cluster $\skl^{\sigma_t}(\CC_t)$. Notice that, $\zeta_t, \zeta'_{t}$ are respectively $\mathcal{F}_{t}, \mathcal{F}'_{t}$-measurable, so the notations above are consistent.
			
			Using these notations, we calculate \eqref{eq.bound.weight} at first conditioned on $\mathcal{F}_{\tau^*_{n}}, \mathcal{F'}_{\phi(\tau^*_{n})}$
			\begin{align*}
				\Eb\Ll[\zeta_\Gamma \zeta'_{\Gamma'}  \, \vert \G\Rr]
				&= \Eb\Ll[\zeta_{\tau^*_{n}} \zeta'_{\phi(\tau^*_{n})} \Eb\Ll[ \prod_{u \geq \tau^*_{n}} \frac{\Ind{L_u}}{s_u}  \prod_{u \geq \phi(\tau^*_{n})} \frac{\Ind{L'_u}}{s'_u}  \, \Big\vert \mathcal{F}_{\tau^*_{n}}, \mathcal{F'}_{\phi(\tau^*_{n})}\Rr] \, \Big\vert \G\Rr]\\
				&= \Eb\Ll[\zeta_{\tau^*_{n}} \zeta'_{\phi(\tau^*_{n})} \Eb\Ll[ \prod_{u \geq \tau^*_{n}} \frac{\Ind{L_u}}{s_u}  \, \Big\vert \mathcal{F}_{\tau^*_{n}}, \mathcal{F'}_{\phi(\tau^*_{n})}\Rr] \Eb\Ll[ \prod_{u \geq \phi(\tau^*_{n})} \frac{\Ind{L'_u}}{s'_u}  \, \Big\vert \mathcal{F}_{\tau^*_{n}}, \mathcal{F'}_{\phi(\tau^*_{n})}\Rr] \, \Big\vert \G\Rr]\\
				&= \Eb\Ll[\zeta_{\tau^*_{n}} \zeta'_{\phi(\tau^*_{n})} \, \Big\vert \G\Rr].
			\end{align*}
			Notice that, the explorations starting from site $z_{\tau^*_{n}} = z'_{\phi(\tau^*_{n})}$ use disjoint tubes.  Especially, at this site, the two tubes will be $T(z_{\tau^*_{n}}, z_{\tau^*_{n}+1})$ and $T(z'_{\phi(\tau^*_{n})}, z'_{\phi(\tau^*_{n})+1})$, and they are different due to \eqref{eq.def.tau}. We thus can apply at first the independence in the second line; see Lemma~\ref{lem.skl.independent}. Then the martingale structure can be applied respectively to the two explorations. We then further calculate the conditional expectation under $\mathcal{F}_{\tau_{n}}, \mathcal{F'}_{\phi(\tau_{n})}$
			\begin{align}\label{eq.prod.induction}
				\begin{split}
					\Eb\Ll[\zeta_\Gamma \zeta'_{\Gamma'}  \, \vert \G\Rr]
					&= \Eb\Ll[\zeta_{\tau_{n}} \zeta'_{\phi(\tau_{n})} \Eb\Ll[\prod_{\tau_{n} \leq u \leq \tau^*_{n}} \frac{\Ind{L_u}}{s_u}  \prod_{\phi(\tau_{n}) \leq u \leq \phi(\tau^*_{n})} \frac{\Ind{L'_u}}{s'_u} \, \Big\vert \mathcal{F}_{\tau_{n}}, \mathcal{F}'_{\phi(\tau_{n})}\Rr] \, \Big\vert \G\Rr]\\
					&\leq (1-\epsilon_0/4)^{-2 \vert I_n\vert} \Eb\Ll[\zeta_{\tau_{n}} \zeta'_{\phi(\tau_{n})}  \, \Big\vert \G\Rr].
				\end{split}
			\end{align}
			In the exploration of common route $I_n$, we just use the bound \eqref{eq.st.bound} to get a factor. We then iterate the bound above and obtain \eqref{eq.bound.weight}.
			
			\smallskip
			\textit{Step~3: proof of \eqref{eq.bound.weight_intersection}.} The proof is once again an induction. We use a variant definition of $\mathsf{E}_t$ in \eqref{eq.coal.I} for $z_t \in \Gamma \cap \Gamma'$
			\begin{align*}
				\tilde{\mathsf{E}}_t := \{\CC_t \cap \CC'_{\phi(t)} \neq \emptyset\} \cup \{\skl^{\sigma_t}(\CC_t,v_t) \cap \skl^{\sigma'_{\phi(t)}}(\CC'_{\phi(t)},v_{\phi(t)}) \neq \emptyset\}.
			\end{align*}
			The Lemma~\ref{lem.coal.I} reduces the object to proof
			\begin{align}\label{eq.bound.weight_intersection.reduction}
				\Eb\Ll[\zeta_\Gamma \zeta'_{\Gamma'}  \Ll(\sum_{z_t \in \cup_{j=1}^{n} I_j} \Ind{\tilde{\mathsf{E}}_t}\Rr) \, \Big\vert \G\Rr] \lesssim \ell_{0}^{6-d} \cdot \vert \Gamma \cap \Gamma' \vert \cdot \Ll(1 + \epsilon_0\Rr)^{\vert \Gamma \cap \Gamma' \vert}. 
			\end{align}
			The main term can be decomposed into two parts
			\begin{align*}
				\zeta_\Gamma \zeta'_{\Gamma'}  \Ll(\sum_{z_t \in \cup_{j=1}^{n} I_j} \Ind{\tilde{\mathsf{E}}_t}\Rr) = \mathbf{J.1} + \mathbf{J.2},
			\end{align*}
			where 
			\begin{align*}
				\mathbf{J.1} &= \zeta_\Gamma \zeta'_{\Gamma'} \Ll(\sum_{z_t \in \cup_{j=1}^{n-1} I_j} \Ind{\tilde{\mathsf{E}}_t}\Rr),\\
				\mathbf{J.2} &= \zeta_\Gamma \zeta'_{\Gamma'} \Ll(\sum_{z_t \in I_n} \Ind{\tilde{\mathsf{E}}_t}\Rr).
			\end{align*}
			
			For the term $\mathbf{I}$, we make a conditional expectation under $\mathcal{F}_{\tau_{n}}, \mathcal{F'}_{\phi(\tau_{n})}$, and obtain that 
			\begin{multline}\label{eq.weight.coal.I}
				\Eb\Ll[\mathbf{J.1} \, \Big\vert \mathcal{F}_{\tau_{n}}, \mathcal{F'}_{\phi(\tau_{n})}\Rr] =  \Ll(\sum_{z_t \in \cup_{j=1}^{n-1} I_j} \Ind{\tilde{\mathsf{E}}_t}\Rr) \Eb\Ll[\zeta_\Gamma \zeta'_{\Gamma'} \, \Big\vert \mathcal{F}_{\tau_{n}}, \mathcal{F'}_{\phi(\tau_{n})}\Rr] \\
				\leq (1+\epsilon_0)^{\vert I_n\vert} \zeta_{\tau_{n}} \zeta'_{\phi(\tau_{n})}  \Ll(\sum_{z_t \in \cup_{j=1}^{n-1} I_j} \Ind{\tilde{\mathsf{E}}_t}\Rr).
			\end{multline}
			Here the second line relies on the induction in \eqref{eq.prod.induction}.
			
			For the term $\mathbf{II}$, we apply the conditional expectation under $\mathcal{F}_{\tau_{n}-1}, \mathcal{F'}_{\phi(\tau_{n})-1}$, which gives that 
			\begin{align}\label{eq.weight.coal.II}
				\begin{split}
					&\Eb\Ll[\mathbf{J.2} \, \Big\vert \mathcal{F}_{\tau_{n}-1}, \mathcal{F'}_{\phi(\tau_{n})-1}\Rr] \\
					&=  \zeta_{\tau_{n}} \zeta'_{\phi(\tau_{n})} \Eb\Ll[\prod_{\tau_{n} \leq u \leq \tau^*_{n}} \frac{\Ind{L_u}}{s_u}  \prod_{\phi(\tau_{n}) \leq u \leq \phi(\tau^*_{n})} \frac{\Ind{L'_u}}{s'_u}   \Ll(\sum_{z_t \in I_n} \Ind{\tilde{\mathsf{E}}_t}\Rr) \, \Big\vert \mathcal{F}_{\tau_{n}-1}, \mathcal{F'}_{\phi(\tau_{n})-1}\Rr] \\
					&\leq C_{\ref{eq.length.path.conditional}} \ell_{0}^{6-d}  \cdot \Ll(1 + \epsilon_0\Rr)^{\vert I_n \vert} \cdot \vert I_n \vert \cdot \zeta_{\tau_{n}} \zeta'_{\phi(\tau_{n})}.
				\end{split}
			\end{align}
			Here we apply \eqref{eq.length.path.conditional} to the common oriented route $I_n$. In particular, the bound above follows Case~$\mathbf{II}$ there. We make the conditional expectation under $\mathcal{F}_{\tau_{n}-1}, \mathcal{F'}_{\phi(\tau_{n})-1}$ in order to make use of the sprinkling cluster from the step before, which then come to the common site $z_{\tau_n}$.
			
			We combine \eqref{eq.weight.coal.I} and \eqref{eq.weight.coal.II}, and apply \eqref{eq.weight} to conclude that 
			\begin{multline*}
				\Eb\Ll[\zeta_\Gamma \zeta'_{\Gamma'}  \Ll(\sum_{z_t \in \cup_{j=1}^{n} I_j} \Ind{\tilde{\mathsf{E}}_t}\Rr) \, \Big\vert \G\Rr]\\ 
				\leq \Ll(1 + \epsilon_0\Rr)^{\vert I_n \vert} \Eb\Ll[\zeta_{\tau_{n}} \zeta'_{\phi(\tau_{n})}  \Ll(\sum_{z_t \in \cup_{j=1}^{n-1} I_j} \Ind{\tilde{\mathsf{E}}_t}\Rr) \, \Big\vert \G\Rr] + C_{\ref{eq.length.path.conditional}} \ell_{0}^{6-d}  \cdot \Ll(1 + \epsilon_0\Rr)^{\vert I\vert} \cdot \vert I_n \vert.
			\end{multline*}
			An iteration of the estimate above yields \eqref{eq.bound.weight_intersection.reduction},  thus \eqref{eq.bound.weight_intersection}. We highlight that, in the last iteration involving $\sum_{z_t \in I_0} \Ind{\tilde{\mathsf{E}}_t}$, one needs \eqref{eq.length.path.conditional} with Case~$\mathbf{I}$ there.

		\end{proof}

		\subsection{Estimate of resistance}\label{subsec.R.estimate}
		With all the preparation above, we conclude Proposition~\ref{prop.lower.reduction}. Recall that, this proposition allows to conclude the lower bound $\ab(\p) \gtrsim (\p-\p_c)^3$; see the proof after its statement in Section~\ref{subsec.flow}. The proof of Proposition~\ref{prop.lower.reduction} is very similar as \cite[Section~2]{BPP98}, and the new input is the analysis below the length $\ell_{0}$ in Proposition~\ref{prop.weight_intersection}.
		
		\begin{proof}[Proof of Proposition~\ref{prop.lower.reduction}]
			The proof is divided into three steps. We shift the object from $\Lambda_{L}$ to $\Lambda_{L-\ell_{0}/2}$, but it will not change the result as $L \to \infty$.	
			
			\textit{Step~1: choice of $\epsilon_0$.} The parameters fixed in Propositions~\ref{prop.proba.pregood} and ~\ref{prop.proba.k-nest} all depend on $\epsilon_0$. Now we give its choice. The constrain comes from the competition between the geometric decay in Proposition~\ref{prop.EIT} and the geometric growth in \eqref{eq.bound.weight_intersection}. Recall that $q_0$ there only depends on $d$, and we need the EIT property to eat all the geometric growth terms. Therefore, we choose $\epsilon_0  \in(0,\frac14]$ small enough, depending on $d$, so that
			\begin{equation}\label{eq.nw.eps0}
				\Theta_0:=(1+ 100\epsilon_0)^4 q_0 < 1.
			\end{equation}
			Afterwards, the constants in Propositions~\ref{prop.proba.pregood}, ~\ref{prop.proba.k-nest}, ~\ref{prop.weight_intersection} only depend on $d$.
			
			\textit{Step~2: event $A_L$ and proof of \eqref{eq.lower.Alower}.} We propose the good event $A_L$ as
			\begin{align}\label{eq.def.AL}
				A_L:=\Ll\{W\ge\frac12\Rr\}.
			\end{align}
			The proof relies on  the Paley--Zygmund inequality
			\begin{align}\label{eq.W2.PZ}
				\PP\Ll[W\ge\frac12\Rr] = \PP\Ll[W \ge \frac12 \Eb[W]\Rr] \ge \frac14 \frac{\Eb[W]^2}{\Eb[W^2]} \ge  \frac{1}{4\Eb[W^2]}. 
			\end{align}
			Here we use the fact $\Eb[W]=1$ from Lemma~\ref{lem.weight.mean}. It suffices to obtain a second moment estimate uniformly with respect to $L$. Then the calculation is very similar as \eqref{lem.flow_intersection.EIT}, and we adapt the notations of independent copies there
			\begin{align*}
				\Eb[W^2] &=\sum_{\Gamma, \Gamma' \in \mathscr{G}_{M, \ell_{0}}} \mu\otimes\mu\big(\Gamma, \Gamma'\big)   \Eb\Ll[ \omega_\Gamma \omega_{\Gamma'} \Eb[\zeta_\Gamma \zeta'_{\Gamma'} \, \vert \G] \Rr]\\
				&\leq \sum_{\Gamma, \Gamma' \in \mathscr{G}_{M, \ell_{0}}} \mu\otimes\mu\big(\Gamma, \Gamma'\big)   \Eb\Ll[ \omega_\Gamma \omega_{\Gamma'} \Rr] \Ll(1 + \frac{\epsilon_0}{2}\Rr)^{\vert \Gamma \cap \Gamma' \vert}.
			\end{align*}
			We use \eqref{eq.bound.weight} in the second line. Notice the independence between different cubes in the indicator of $\omega_\Gamma$ and $q \geq 1-\epsilon_0/4$ from \eqref{eq.st.bound}, we then have
			\begin{align}\label{eq.bound.omega_Gamma}
				\Eb[\omega_\Gamma\omega_{\Gamma'}]=q^{-\vert \Gamma \cap \Gamma' \vert}\le(1-\epsilon_0/4)^{-\vert \Gamma \cap \Gamma' \vert}.
			\end{align}
			Therefore, by Proposition~\ref{prop.EIT} and \eqref{eq.nw.eps0},
			\begin{align*}
				\Eb[W^2]&\leq \sum_{\Gamma, \Gamma' \in \mathscr{G}_{M, \ell_{0}}} \mu\otimes\mu\big(\Gamma, \Gamma'\big)    \Ll(1 + \frac{\epsilon_0}{2}\Rr)^{2 \vert \Gamma \cap \Gamma' \vert}  \\
				&\leq 1+\sum_{k\ge1} C q_0^{k}M^{2-d} \Ll(1 + \frac{\epsilon_0}{2}\Rr)^{2k} \\
				&\leq 1+\sum_{k\ge1}\Theta_0^{k} M^{2-d}\le 1+CM^{2-d}\le2 .
			\end{align*}
			The last line relies on the choice of parameter in \eqref{eq.nw.eps0}. We put this bound back to \eqref{eq.W2.PZ}, and concludes a uniform lower bound $\PP[A_L] \geq \frac{1}{8} =: c_A$.

			\textit{Step~3: proof of \eqref{eq.lower.Rupper}.} 	We combine the expression in \eqref{eq.flow_intersection.EIT} and the lower bound of $W$ on $A_L$, then estimate it under $\Eb[(\, \cdot \, ) \vert \G]$
			\begin{align*}
				&\Eb\Ll[\mathscr E_*(\Lambda_{L-\ell_{0}/2},\a,\theta)\Ind{A_L} \, \vert \G\Rr]\\
				&\leq 4 \sum_{\Gamma, \Gamma' \in \mathscr{G}_{M, \ell_{0}}} \mu\otimes\mu\big(\Gamma, \Gamma'\big)   \omega_\Gamma \omega_{\Gamma'}  \Eb\Ll[\zeta_\Gamma \zeta'_{\Gamma'} \vert \gamma \cap \gamma'\vert\, \vert \G\Rr]\\
				&\leq 4 K_{\ref{eq.bound.weight_intersection}} \ell_{0}^{8-d} \sum_{\Gamma, \Gamma' \in \mathscr{G}_{M, \ell_{0}}} \mu\otimes\mu\big(\Gamma, \Gamma'\big)  \omega_\Gamma \omega_{\Gamma'}  \cdot \vert \Gamma \cap \Gamma' \vert \cdot \Ll(1 + \frac{\epsilon_0}{2}\Rr)^{\vert \Gamma \cap \Gamma' \vert}.
			\end{align*}
			In the second line above, we also use the estimate \eqref{eq.bound.weight_intersection}, and the fact that  $\omega_\Gamma \omega_{\Gamma'}$ are both $\G$-measurable. 
			We then further take the expectation under $\PP$ and use the bound \eqref{eq.bound.omega_Gamma} again
			\begin{align*}
				&\Eb\big[\mathscr E_*(\Lambda_{L-\ell_{0}/2},\a,\theta)\Ind{A_L}\big] \\
				&\lesssim \ell_{0}^{8-d} \sum_{k\ge 1} k q_0^{k}M^{2-d} \Ll(1 + \frac{\epsilon_0}{2}\Rr)^{2k}\\
				&\lesssim \ell_{0}^{8-d} M^{2-d} \sum_{k\ge 1} k \Theta_0^{k}\\
				&\lesssim \ell_{0}^6 L^{2-d}.
			\end{align*}
			The last line relies once again the choice of parameter in \eqref{eq.nw.eps0}, which ensures the sum is finite. The relation $M = \frac{L}{\ell_{0}}$ is also inserted; see \eqref{eq.def.M}. Since $\ell_{0} = n_0 \ell_{*}$ and now $n_0$ is a constant only depending on $d$, this bound justifies \eqref{eq.lower.Rupper}.
		\end{proof}

		We finally prove the main theorem.
		\begin{proof}[Proof of Theorem~\ref{thm.main}]
			Combining the upper bound in Proposition~\ref{prop.upper} and the lower bound in Proposition~\ref{prop.lower}, we have 
			\begin{align*}
				(\p-\p_c)^3 \lesssim \ab(\p) \lesssim(\p-\p_c)^3.
			\end{align*}
			Then we recall the Einstein relation in \eqref{eq.def_sigma} and $\theta(\p) \approx (\p-\p_c), \p \searrow \p_c$ in \cite{hara1990mean}
			\begin{align*}
				(\p-\p_c)^2 \lesssim \sigma^2(\p) = \frac{2\ab(\p)}{\theta(\p)} \lesssim (\p-\p_c)^2.
			\end{align*}
			This concludes Theorem~\ref{thm.main}.
		\end{proof}

		\section{Intrinsic distances and slab thresholds}\label{sec.chem}
		
		In this section we prove Theorem~\ref{thm.time_constant} (in Section~\ref{sec.chem.dist}) and Theorem~\ref{thm.slab} (in Section~\ref{sec.chem.slab}).
		Throughout, $d\ge9$ and \eqref{eq:two_point_assumption} holds, as in these theorems.
		Recall that $|\cdot|$ denotes the $\ell^\infty$ norm, and we write $|\cdot|_1$ for the $\ell^1$ norm.
		Recall also that $\ell_*(\p)=\lfloor(\p-\pc)^{-1/2}\rfloor$.
		
		\subsection{Intrinsic distances}\label{sec.chem.dist}
		
		For $x,y\in\Zd$, we denote by $\din(x,y)$ the length of the shortest open path from $x$ to $y$, with the convention $\din(x,y)=+\infty$ if $x$ and $y$ are not connected by an open path.
		For $\p>\pc$, we write $\theta(\p):=\Pp[0\in\cCf]$.
		
		We first recall the results of Garet and Marchand \cite{garet2004asymptotic,garet2007large}, which define the time constant $\mu(\p)$ and the norm $\|\cdot\|_\p$ appearing in the introduction.
		Let $\p>\pc$.
		For every $u\in\Zd\setminus\{0\}$ there exists $\mu_\p(u)\in[0,\infty)$ such that, $\Pp$-almost surely,
		\begin{equation}\label{eq.chem.timeconstant}
			\lim_{n\to\infty,\ 0\leftrightarrow nu}\frac{\din(0,nu)}{n}=\mu_\p(u)
		\end{equation}
		(\cite[Corollary~3.3]{garet2004asymptotic}).
		The function $\mu_\p$ extends to a norm on $\R^d$ (\cite[Theorem~3.2, Corollary~4.2]{garet2004asymptotic}), and it is invariant under all permutations and reflections of the coordinates (\cite[Theorem~3.2\,(1)]{garet2004asymptotic}).
		We write
		\[
		\mu(\p):=\mu_\p(e_1),\qquad \|x\|_\p:=\frac{\mu_\p(x)}{\mu(\p)}\quad(x\in\R^d),
		\]
		where $\mu(\p)>0$ because $\mu_\p$ is a norm.
		Then $\|\cdot\|_\p$ is a norm with $\|e_i\|_\p=1$ for every $i\in\{1,\dots,d\}$ by the symmetry of $\mu_\p$, and $\mu_\p(x)=\mu(\p)\|x\|_\p$.
		The constant $\mu(\p)$ is the \emph{time constant}.
		Garet and Marchand also proved in \cite{garet2007large} that the convergence \eqref{eq.chem.timeconstant} holds with exponentially small error probabilities: for every $\varepsilon>0$,
		\begin{equation}\label{eq.chem.GM07}
			\limsup_{|x|_1\to\infty}\frac{1}{|x|_1}\log \Pp\Big[0\leftrightarrow x,\ \frac{\din(0,x)}{\mu(\p)\|x\|_\p}\notin(1-\varepsilon,1+\varepsilon)\Big]<0
		\end{equation}
		(\cite[Theorem~1.2]{garet2007large}).
		Together, \eqref{eq.chem.timeconstant} and \eqref{eq.chem.GM07} give the asymptotic 
		\begin{align*}
			\din(0,x)=(1+o(1))\mu(\p)\|x\|_\p, \text{ on } \{0\leftrightarrow x\},
		\end{align*}
		recalled in the introduction. These results are qualitative in $\p$, and Theorem~\ref{thm.time_constant} determines the order of $\mu(\p)$ as $\p\downarrow\pc$.
		
		The next proposition says that $\|\cdot\|_\p$ is comparable to $|\cdot|$, uniformly in $\p$.
		
		\begin{lemma}\label{lem.chem.alldir}
			For every $\p\in(\pc,1]$ and $x\in\R^d$,
			\[
			|x|\ \le\ \|x\|_\p\ \le\ |x|_1 .
			\]
			In particular, $\|x\|_\p\asymp|x|$, uniformly in $x\in\R^d$ and $\p\in(\pc,1]$.
		\end{lemma}
		
		\begin{proof}
			For the upper bound, the triangle inequality and $\|e_i\|_\p=1$ give 
			\begin{align*}
				\|x\|_\p\le\sum_{i=1}^d|x\cdot e_i|\,\|e_i\|_\p=|x|_1.
			\end{align*}
			For the lower bound, choose $i$ such that $|x\cdot e_i|=|x|$, and let $x'$ be the point obtained from $x$ by reflecting all the coordinates except the $i$-th.
			Then $\|x'\|_\p=\|x\|_\p$ by the invariance of $\|\cdot\|_\p$ under reflections, and $\frac12(x+x')=(x\cdot e_i)e_i$.
			Hence
			\[
			|x|=\|(x\cdot e_i)e_i\|_\p\le\tfrac12\big(\|x\|_\p+\|x'\|_\p\big)=\|x\|_\p. \qedhere
			\]
		\end{proof}

		By Lemma~\ref{lem.chem.alldir} and \eqref{eq.chem.GM07}, the intrinsic distance between two far apart points of the infinite cluster is of order $\mu(\p)$ times their distance, so it remains to determine the order of $\mu(\p)$; this is the content of Theorem~\ref{thm.time_constant}.
		The lower bound comes from the exponential decay of the slightly subcritical two-point function. The upper bound comes from the crossing paths used in the lower bound for the conductivity (Proposition~\ref{prop.lower.reduction}).
		
		\subsubsection{The lower bound}\label{sec.chem.lower}
		
		The first lemma compares the probability of a short connection at two different parameters.
		
		\begin{lemma}\label{lem.chem.monotone}
			For every $0<\p'<\p\le1$, $x,y\in\Zd$ and integer $m\ge0$,
			\[
			\Pp\big[\din(x,y)\le m\big]\ \le\ \Big(\frac{\p}{\p'}\Big)^{m}\ \P_{\p'}\big[\din(x,y)\le m\big].
			\]
		\end{lemma}
		
		\begin{proof}
			Let $\omega \equiv (\omega_e)_{e \in \Ed} \sim\Pp$, and let $\xi \equiv (\xi_e)_{e \in \Ed}$ be i.i.d.\ Bernoulli random variables of parameter $\p'/\p$, independent of $\omega$. We denote by $\P$ the probability space containing $\omega$ and $\xi$. Then the random variables
			\begin{align*}
				\omega'_e:=\omega_e\xi_e,
			\end{align*}
			define a configuration, which follows the law $\omega' \equiv (\omega'_e)_{e \in \Ed} \sim\P_{\p'}$ and satisfies $\omega'\le\omega$. 
			
			We study the event 
			\begin{align*}
				A:=\{\din(x,y)\le m\}.
			\end{align*}
			For $\omega\in A$, let $\gamma(\omega)$ be a shortest open path from $x$ to $y$ with $|\gamma(\omega)|\le m$. Then in the coupled space $(\omega, \xi)$, we observe an inclusion that 
			\begin{align*}
				\{ \omega' \in A\} \supset \Big(\{ \omega \in A \} \cap \{ \forall e \in \gamma(\omega), \xi_e = 1\}\Big).
			\end{align*}
			We take the expectation under $\P$ and obtain that 
			\begin{align*}
				\P_{\p'}[A] &\ge \P[\{ \omega \in A \} \cap \{ \forall e \in \gamma(\omega), \xi_e = 1\}] \\
				&=  \E\Ll[\Ind{\omega\in A} \P[\{ \forall e \in \gamma(\omega), \xi_e = 1\} \vert \omega] \Rr] \\ 
				&= \E\big[\Ind{\omega\in A}\,(\p'/\p)^{|\gamma(\omega)|}\big] \\
				&\ge(\p'/\p)^m\,\Pp[A]. 
			\end{align*}
			The first line is the consequence of inclusion above. The second line takes the conditional expectation under a given $\omega$. Then the independence between $\omega$ and $\xi$ yields the third line. For the forth line,  $(\p'/\p)^{|\gamma(\omega)|}\ge(\p'/\p)^m$ comes from the condition $|\gamma(\omega)|\le m$ on $A$. This justifies the statement in the lemma.
		\end{proof}
		
		The second lemma is the input from the high-dimensional literature.
		
		\begin{lemma}\label{lem.chem.HMS}
			There exists $c_{\ref{lem.chem.HMS}}>0$ such that for every $\p\in(0,\pc]$ and $n\ge1$,
			\[
			\P_{\p}\big[0\leftrightarrow ne_1\big]\lesssim\exp\Big(-c_{\ref{lem.chem.HMS}}\,(\pc-\p)^{1/2}\,n\Big).
			\]
		\end{lemma}
		
		\begin{proof}
			This is from  \cite[Theorem~1.1]{hutchcroft2023high} applied to $x=ne_1$.
		\end{proof}
		
		The third lemma is the converse of the lower large deviation estimate in \cite{garet2007large}; it only uses the almost sure convergence \eqref{eq.chem.timeconstant}.
		
		\begin{lemma}\label{lem.chem.converse}
			Let $\p\in(\pc,1]$ and $a>\mu(\p)$.
			Then
			\[
			\liminf_{n\to\infty}\Pp\big[\din(0,ne_1)\le an\big]\ \ge\ \theta(\p)^2>0 .
			\]
		\end{lemma}
		
		\begin{proof}
			Define the event 
			\begin{align*}
				E_n:=\{0\in\cCf,\ ne_1\in\cCf\}.
			\end{align*}
			Because both events $\{0\in\cCf\}$ and $\{ne_1\in\cCf\}$ are increasing, we have
			\begin{align*}
				\Pp[E_n]\ge\theta(\p)^2,
			\end{align*} 
			by the FKG inequality. Since $\cCf$ is unique, $E_n\subset\{0\leftrightarrow ne_1\}$.
			By \eqref{eq.chem.timeconstant} with $u=e_1$, almost surely on $\{0\in\cCf\}$ we have $\din(0,ne_1)\le an$ for all large enough $n$ such that $ne_1\in\cCf$, because $a>\mu(\p)$.
			Therefore $\Ind{E_n\cap\{\din(0,ne_1)>an\}}\to0$ almost surely, and, thus,  
			\begin{align*}
				\Pp[E_n\cap\{\din(0,ne_1)>an\}]\to 0.
			\end{align*}
			Combining the estimates above, we conclude that
			\[
			\Pp\big[\din(0,ne_1)\le an\big]\ge\Pp[E_n]-\Pp\big[E_n\cap\{\din(0,ne_1)>an\}\big]\ge\theta(\p)^2-o(1). \qedhere
			\]
		\end{proof}
		
		With the lemmas above, we prove the lower bound of the time constant.
		\begin{prop}\label{prop.chem.lower}
			For every $\p\in(\pc,1]$, $\mu(\p)\gtrsim(\p-\pc)^{-1/2}$.
		\end{prop}
		
		\begin{proof}
			Throughout this proof, we write 
			\begin{align*}
				\delta:=\p-\pc.
			\end{align*}
			Since $\din(0,ne_1)\ge n$ whenever $0\leftrightarrow ne_1$, \eqref{eq.chem.timeconstant} gives $\mu(\p)\ge1$.
			This proves the claim when $\delta\ge\pc/2$.
			
			Assume now that $\delta<\pc/2$, and we write 
			\begin{align*}
				\p_-:=\pc-\delta\in(\pc/2,\pc).
			\end{align*}
			Let $a>0$ and $n\ge1$.
			Since $\{\din(0,n e_1)\le an\}$ is contained in $\{0\leftrightarrow n e_1\}$, Lemmas~\ref{lem.chem.monotone} and~\ref{lem.chem.HMS} (with $\p'=\p_-$ and $m=\lfloor an\rfloor$) give
			\[
			\Pp\big[\din(0,n e_1)\le an\big]\le\Big(\frac{\p}{\p_-}\Big)^{an}\P_{\p_-}\big[0\leftrightarrow n e_1\big]\lesssim\Big(\frac{\p}{\p_-}\Big)^{an}\exp\Big(-c_{\ref{lem.chem.HMS}}\,\delta^{1/2}n\Big).
			\]
			Since $\delta<\pc/2$, we have the following elementary calculation
			\[
			\frac{\p}{\p_-}=1+\frac{2\delta}{\pc-\delta}\le1+\frac{4\delta}{\pc}\le e^{4\delta/\pc}.
			\]
			By choosing $a:=\frac{c_{\ref{lem.chem.HMS}}\pc}{8}\,\delta^{-1/2}$, it satisfies $4a\delta/\pc=\frac{c_{\ref{lem.chem.HMS}}}{2}\delta^{1/2}$, and we obtain
			\[
			\Pp\big[\din(0,n e_1)\le an\big]\lesssim\exp\Big(-\frac{c_{\ref{lem.chem.HMS}}}{2}\,\delta^{1/2}n\Big)\xrightarrow{n\to\infty}0 .
			\]
			If $a>\mu(\p)$, Lemma~\ref{lem.chem.converse} would give a positive lower bound on the $\liminf$ of the left-hand side, which is a contradiction.
			Hence $\mu(\p)\ge a\gtrsim\delta^{-1/2}$.
		\end{proof}
		
		\subsubsection{The upper bound}\label{sec.chem.upper}
		
		The proof of the upper bound relies on the existence of geodesic path within desired scaling. 
		
		\begin{lemma}\label{lem.chem.crossing}
			There exists a constant $C_{\ref{lem.chem.crossing}}(d) < +\infty$, such that the event $E_{L}$
			\begin{align*}
				E_L := \{\text{ there exists an open path connecting } [\partial\Lambda_L]^- \text{ to } [\partial\Lambda_L]^+ \text{ within } C_{\ref{lem.chem.crossing}}\,\ell_*(\p)\,L \text{ bonds }\}
			\end{align*}
			holds with a positive probability for every $\p\in(\pc,1]$
			\[
			\liminf_{L\to\infty}\Pp[E_L]\ge \frac{1}{8} .
			\]
		\end{lemma}
		
		\begin{proof}
			The event $A_L = \{ W \geq \frac{1}{2}\}$ defined in \eqref{eq.def.AL} is a specific subset of $E_L$-type event. In Step~2 of the proof of Proposition~\ref{prop.lower.reduction}, we have showed that 
			\begin{align*}
				\liminf_{L\to\infty}\Pp[A_L]\ge \frac{1}{8}.
			\end{align*}
			We then verify the existence of the crossing geodesic path. Recall that $W$ defined in \eqref{eq.W} is the total weight of the crossing random flow, sampled over the couple of oriented route and exploration process $(\Gamma, \expl(\Gamma, \a; \sigma, U))$; see the details in Definition~\ref{def.weight}. Then, $W>\frac{1}{2}$ suggests the existence of at least one such crossing path $\gamma \equiv \gamma(\Gamma, \a; \sigma, U)$ defined in \eqref{eq.def.path}. The estimate in \eqref{eq.path.length} suggests that 
			\begin{align*}
				\vert \gamma \vert \leq \Ll(K_{\mathrm{spn}} + K_{\mathrm{skl}}\Rr)\ell_{0}^2 \vert \Gamma \vert \leq \underbrace{2\Ll(K_{\mathrm{spn}} + K_{\mathrm{skl}}\Rr) n_0}_{=:C_{\ref{lem.chem.crossing}}}\ell_{*}(\p) L.
			\end{align*} 
			The second inequality follows the facts
			\begin{itemize}
				\item $\vert \Gamma \vert \leq 2 M$ for the oriented route (Definition~\ref{def.route.o});
				\item $M = \frac{L}{\ell_{0}}$ (\eqref{eq.def.M});
				\item $\ell_{0} = n_0 \ell_{*}(\p)$ (\eqref{eq.ell0}).
			\end{itemize}
			Besides, the constants $K_{\mathrm{spn}}, K_{\mathrm{skl}}, n_0$ only depend on $d$; see Step~1 in the proof of Proposition~\ref{prop.lower.reduction}. Therefore, by choosing $C_{\ref{lem.chem.crossing}}$ as above, we have $A_L \subset E_L$ and justify the claim.
		\end{proof}
		
		The second lemma says that, with high probability, no pair of points in $\Lambda_L$ at distance of order $L$ in the direction $e_1$ is connected by an unusually short path.
		
		\begin{lemma}\label{lem.chem.pairs}
			For every $\p\in(\pc,1]$ and $\varepsilon\in(0,1)$, the event 
			\begin{align*}
				F_L := \{\exists\,x,y\in\Lambda_L:\ (y-x)\cdot e_1\ge L/2,\ x\leftrightarrow y,\ \din(x,y)\le(1-\varepsilon)\mu(\p)\|y-x\|_\p\}
			\end{align*}
			satisfies
			\[
			\lim_{L\to\infty}\Pp[F_L]=0 .
			\]
		\end{lemma}
		
		\begin{proof}
			The event $\{x\leftrightarrow y,\ \din(x,y)\le(1-\varepsilon)\mu(\p)\|y-x\|_\p\}$ is contained in the event inside the probability in \eqref{eq.chem.GM07} (with $y-x$ in place of $x$).
			Hence \eqref{eq.chem.GM07} and translation invariance give $B>0$ and $R_0<\infty$ such that, for all $x,y\in\Zd$ with $|y-x|_1\ge R_0$,
			\[
			\Pp\big[x\leftrightarrow y,\ \din(x,y)\le(1-\varepsilon)\mu(\p)\|y-x\|_\p\big]\le e^{-B|y-x|_1}.
			\]
			If $(y-x)\cdot e_1\ge L/2$, then $|y-x|_1\ge L/2$. A union bound for $L\ge2R_0$ then yields that 
			\begin{align*}
				\Pp[F_L] &\leq \sum_{x,y\in\Lambda_L:\ (y-x)\cdot e_1\ge L/2} \Pp\Ll[\{x\leftrightarrow y,\ \din(x,y)\le(1-\varepsilon)\mu(\p)\|y-x\|_\p\}\Rr] \\
				&\leq L^{2d}e^{-BL/2} \xrightarrow{L \to \infty} 0.
			\end{align*}
		\end{proof}
		
		\begin{prop}\label{prop.chem.upper}
			For every $\p\in(\pc,1]$, $\mu(\p)\lesssim(\p-\pc)^{-1/2}$.
		\end{prop}
		
		\begin{proof}
			Fix $\varepsilon\in(0,1)$. Let $E_L$ be the event in Lemma~\ref{lem.chem.crossing}, and let $F_L$ be the event in Lemma~\ref{lem.chem.pairs}. The estimates in two lemmas imply that
			\begin{align*}
				\liminf_{L \to \infty}\Pp[E_L\cap F^c_L] \geq \liminf_{L \to \infty}\Pp[E_L] - \lim_{L \to \infty}\Pp[F_L] \geq \frac{1}{8}.
			\end{align*}
			Fix a configuration in $E_L\cap F^c_L$, we consider one $\gamma$ as an open path in $\Lambda_L$ from $x\in[\partial\Lambda_L]^-$ to $y\in[\partial\Lambda_L]^+$ with $|\gamma|\le C_{\ref{lem.chem.crossing}}\,\ell_*(\p)\,L$.
			
			On the one hand, such pair $(x,y)$ satisfies the condition in $F^c_L$, which gives 
			\begin{align*}
				\din(x,y)>(1-\varepsilon)\mu(\p)\|y-x\|_\p \ge (1-\varepsilon)\mu(\p) (y-x)\cdot e_1\ge\frac{(1-\varepsilon)\mu(\p)L}{2}.
			\end{align*}
			The second inequality comes from Lemma~\ref{lem.chem.alldir}.
			
			On the other hand, the condition $E_L$ implies that 
			\begin{align*}
				\din(x,y)\le|\gamma|\le C_{\ref{lem.chem.crossing}}\,\ell_*(\p)\,L.
			\end{align*}
			Therefore, we conclude that 
			\[
			(1-\varepsilon)\,\mu(\p)\,\frac L2<\din(x,y)\le C_{\ref{lem.chem.crossing}}\,\ell_*(\p)\,L,
			\]
			which implies that
			\begin{align*}
				\mu(\p)<2C_{\ref{lem.chem.crossing}}\,\ell_*(\p)/(1-\varepsilon).
			\end{align*} 
			Using $\ell_*(\p)\le(\p-\pc)^{-1/2}$ completes the proof.
		\end{proof}
		
		\begin{proof}[Proof of Theorem~\ref{thm.time_constant}]
			The lower bound is Proposition~\ref{prop.chem.lower}, and the upper bound is Proposition~\ref{prop.chem.upper}.
		\end{proof}
		
		\subsection{Slab thresholds}\label{sec.chem.slab}
		We prove Theorem~\ref{thm.slab} in this section.
		
		\subsubsection{The lower bound}\label{sec.slab.lower}
		
		\begin{prop}\label{prop.slab.lower}
			For every $n\ge1$,
			\[
			\pc\big(\ZZ^{d-1}\times\{0,\dots,n\}\big)-\pc\ \gtrsim\ n^{-2}.
			\]
		\end{prop}
		
		\begin{proof}
			Up to a permutation of the coordinates, $\ZZ^{d-1}\times\{0,\dots,n\}$ is the slab $S_{0,n}$ of \eqref{eq.def.Smn_Lj}.
			Let $c_{**}$ be the constant in Lemma~\ref{lem.A.upperbound.super} and let $\p=\pc+c_{**}n^{-2}$. We denote by 
			\begin{align*}
				\cC_n(k e_1) := \{x \in S_{0,n}: x \xleftrightarrow{S_{0,n}} k e_1\},
			\end{align*}
			the cluster issued from $k e_1$ and restricted in $S_{0,n}$. Recall $A_{\p}(0,n;k,j)$ defined in \eqref{eq.def.Amnkj}, then Lemma~\ref{lem.A.upperbound.super} yields that 
			\begin{align*}
				\E_{\p}[\vert  \cC_n(k e_1) \vert] = \sum_{j=0}^nA_{\p}(0,n;k,j) \lesssim \sum_{j=0}^n (n+1) \lesssim (n+1)^2
			\end{align*}
			because each term satisfies $A_{\p}(0,n;k,j) \lesssim n+1$, for every $0\le j \le n$. This expected size is at most of order $n^2$, so $S_{0,n}$ has no infinite cluster at $\p$, that is, $\pc(S_{0,n})\ge\p$.
		\end{proof}
		
		\subsubsection{The upper bound}\label{sec.slab.upper}
		
		\begin{prop}\label{prop.slab.upper}
			For every $n\ge1$,
			\[
			\pc\big(\ZZ^2\times\{0,\dots,n\}^{d-2}\big)-\pc\ \lesssim\ n^{-2}.
			\]
		\end{prop}
		
		\begin{proof}
			The proof relies on the exploration process in Definition~\ref{def.expl}. Given a parameter $\epsilon_0 \in (0,1)$, we fix the parameters in Remarks~\ref{rmk.proba.pregood} and ~\ref{rmk.proba.k-nest}. We also set
			\begin{align*}
				n = \ell_{0} = n_0 \ell_{*},
			\end{align*}
			which is equivalent to 
			\begin{align*}
				\p = \p_c + (n/n_0)^{-2}.
			\end{align*}
			We partition $\ZZ^2\times\{0,\dots,n\}^{d-2}$ as cubes of side length $n$, and then adapt the  $(Q,S,H,T)$-structure in \eqref{eq.def.SQH.ell0.QSH} with a shift $- \lfloor \frac{n e_1}{2} \rfloor$. A variant of the exploration process in Definition~\ref{def.expl} can be implemented in the coarsened sites $n\ZZ^2\times\{0\}^{d-2}$: 
			\begin{enumerate}
				\item It starts from $\mathbf{C} = \{0\}$ and a $\kk$-nest $\mathcal{C}(0)$ in $Q(0)$ if it is pre-good, and there is no route $\Gamma$ fixed a prior.
				\item Exploration: Pick one site $z \in \mathbf{C}$. It has an associated $\kk$-nest denoted by $\mathcal{C}(z)$. For all neighbor sites $z' \stackrel{n}{\sim} z$, if 
				\begin{itemize}
					\item $z' \notin \mathbf{C}$,
					\item $Q(z')$ is pre-good,
					\item the sprinkling cluster $\skl(\mathcal{C}(z), \overrightarrow{zz'})$ in $T(z,z')$ is typical,
				\end{itemize}
				add $z'$ in $\mathbf{C}$ and choose uniformly its associated $\kk$-nest $\mathcal{C}(z')$.
				\item Every site in $\mathbf{C}$ makes the exploration above at most once.
			\end{enumerate}
			The outcome will be a connected component $\mathbf{C} \subset n\ZZ^2\times\{0\}^{d-2}$. Using the order of exploration process above, for every site $z \in \mathbf{C}$, there is a unique route $\Gamma(0,z)$ in the exploration to connect $0$ and $z$, because every site is visited at most once. As a consequence, $\mathbf{C}$ has a structure of tree. Following the same argument in Definition~\ref{def.path}, one can extract a cluster 
			\begin{align*}
				\cC \subset \ZZ^2\times\{0,\dots,n\}^{d-2},
			\end{align*}
			to connect $Q(z)_{z \in \mathbf{C} }$. 
			
			Using \eqref{eq.st.bound}, we observe that $\mathbf{C}$ stochastically dominates a site percolation of parameter $(1 - \epsilon_0)^{2}$. By choosing $\epsilon_0$ such that
			\begin{align*}
				(1 - \epsilon_0)^{2} > \p_c(\ZZ^2 \text{ site-percolation}),
			\end{align*}
			then $\mathbf{C}$ is more connected than a supercritical $\ZZ^2$ site percolation, so we have 
			\begin{align*}
				\P_{\p}[ \# \mathbf{C} = + \infty] > 0.
			\end{align*}
			Therefore, $\cC$ also has a positive probability to be an infinite cluster
			\begin{align*}
				\P_{\p}[0 \xleftrightarrow{\cC } \infty] > 0.
			\end{align*}
			This justifies that
			\begin{align*}
				\pc(\ZZ^2\times\{0,\dots,n\}^{d-2})\le\p= \p_c + (n/n_0)^{-2},
			\end{align*}
			which is desired.

		\end{proof}
		
		\begin{proof}[Proof of Theorem~\ref{thm.slab}]
			The second inequality in \eqref{eq.main_slab} holds since $\ZZ^2\times\{0,\dots,n\}^{d-2}$ is a subgraph of $\ZZ^{d-1}\times\{0,\dots,n\}$, the third one is Proposition~\ref{prop.slab.upper}, and the first one is proved in Section~\ref{sec.slab.lower}.
		\end{proof}

		\appendix

		\section{Two-point function in box and annulus}
		In this appendix, we estimate some two-point functions in boxes and annulus. Recall that $\vert\cdot\vert$ denotes the $\ell^\infty$ norm on $\mathbb R^d$.
		
		Throughout this appendix, we write
		\begin{align*}
			\H := \{x \in \Zd : x \cdot  e_1 \geq 0\}, \qquad \H_k :=\{x \in \Zd : x \cdot  e_1 \geq -k\} \quad (k \in \ZZ),
		\end{align*}
		\begin{equation*}
			\mathbb S_n := \H \setminus \H_{-(n+1)} = \{x \in \H : x\cdot e_1 \leq n\}
		\end{equation*}
		for the upper-half space, its shifts, and the slab of height $n$, so that $\partial \H_{-n} = \{x \in \Zd : x \cdot  e_1 = n\}$. We also write $\Lambda_n := [-n,n]^d \cap \Zd$ for the $\ell^\infty$-box of radius $n$ centered at the origin.
		
		We first list some results from other papers that we will use:
		\begin{lemma}[Lemmas~4.1 and 4.2 of \cite{van2025one}]\label{lem.Nn_moments}
			Let $\mathcal N_n := \#\{x \in \partial \H_{-n} : 0 \xleftrightarrow{\mathbb S_n} x\}$. There exist $C_{\ref{lem.Nn_moments}}>0$ and $c_{\ref{lem.Nn_moments}}>0$ such that, for every $n \geq 1$,
			\begin{align*}
				\E_{\p_c}[\mathcal N_n^2] \leq C_{\ref{lem.Nn_moments}}\, n, \qquad \E_{\p_c}[\mathcal N_n] \geq c_{\ref{lem.Nn_moments}}\, \frac{1}{n}.
			\end{align*}
		\end{lemma}
		
		For $S \subset \Zd$ containing $0$, let $\partial S := \{x\in S: \exists y\notin S \text{ with } |x-y|_1=1\}$ denote its inner vertex boundary (so that, in particular, $\partial \H_{-n}$ agrees with the notation already introduced), and recall
		\begin{align*}
			\psi_p(S) := \sum_{x\in\partial S} \P_p\Ll[0\xleftrightarrow{S}x\Rr].
		\end{align*}
		Its bound is in \eqref{eq.expansion_2}.
		
		We record the following identity, which follows from reversibility of connection probabilities together with the invariance of $\H$ under translations in the directions $ e_2,\ldots, e_d$: for every $m\geq0$ and every $p\in\H$ with $p\cdot e_1=m$,
		\begin{align}\label{eq.psi_identity}
			\psi_{\p_c}(\H_m) = \sum_{w\in\H:\, w\cdot e_1=m} \P_{\p_c}\Ll[0\xleftrightarrow{\H}w\Rr] = \sum_{z\in\partial\H} \P_{\p_c}\Ll[p\xleftrightarrow{\H}z\Rr].
		\end{align}
		Indeed, the translation $x\mapsto x+m e_1$ maps $\H_m$ onto $\H$, sending $0$ to a point at height $m$ and $\partial\H_m$ onto $\partial\H$; combined with reversibility and the transverse translation invariance just mentioned, this shows that all three quantities above coincide.

		\begin{lemma}[Theorem~1.4 of \cite{panis2026}]\label{lem.halfspace_twopoint}
			Let $d>6$ and assume that the critical two-point function satisfies $\P_{\p_c} \big[ x\stackrel{\Zd}{\longleftrightarrow} y \big] \asymp \vert x-y\vert^{2-d}$. For $x,y \in \H$, let $r_{x,y} := \min(x \cdot  e_1, |x-y|)$. There exist $c_{\ref{lem.halfspace_twopoint}}>0$ and $C_{\ref{lem.halfspace_twopoint}}>0$ such that, for every $x,y \in \H$,
			\begin{align*}
				c_{\ref{lem.halfspace_twopoint}}\, \frac{(1+r_{x,y})(1+r_{y,x})}{1+|x-y|^d} \leq \P_{\p_c} \Ll[ x\stackrel{\H}{\longleftrightarrow} y \Rr] \leq C_{\ref{lem.halfspace_twopoint}}\, \frac{(1+r_{x,y})(1+r_{y,x})}{1+|x-y|^d}.
			\end{align*}
		\end{lemma}
		
		\begin{lemma}[Theorem~12 of \cite{CH20}]\label{lem.box_twopoint}
			Fix $\alpha_i, \beta_i > 0$ for $1 \leq i \leq d$, and fix $M>1$. For each $n$, let
			\begin{align*}
				R_n := [-\alpha_1 n, \beta_1 n] \times \cdots \times [-\alpha_d n, \beta_d n].
			\end{align*}
			There exists $c_{\ref{lem.box_twopoint}}(M,(\alpha_i),(\beta_i)) > 0$ such that, uniformly in $n$ and in $x \in R_n$,
			\begin{align*}
				\tau_{MR_n}(0,x) \geq c_{\ref{lem.box_twopoint}}\,\|x\|^{2-d}.
			\end{align*}
		\end{lemma}
		
		\begin{lemma}[Theorem~1.2 of \cite{van2025one}, Proposition~3.1 of \cite{CH20}]\label{lem.onearm_half}
			There exist $C_{\ref{lem.onearm_half}}>0$ and $c_{\ref{lem.onearm_half}}>0$ such that, for every $n \geq 1$,
			\begin{align*}
				\frac{c_{\ref{lem.onearm_half}}}{n^3} \leq \P_{\p_c}\Ll[0 \xleftrightarrow{\H} \H \setminus \Lambda_n \Rr] \leq \frac{C_{\ref{lem.onearm_half}}}{n^3}.
			\end{align*}
		\end{lemma}

		\subsection{Boundary of clusters in upper-half space}
		Define the box 
		\begin{align*}
			Q_{n,rn} := [0,n] \times (-rn,rn) \times \cdots \times (-rn,rn) \cap \Zd \qquad (r\geq1),
		\end{align*}
		and, for $i \in \{1,\ldots,d\} \cup \{-2,\ldots,-d\}$, its faces
		\begin{align*}
			\partial_i Q_{n,rn} := \begin{cases}
				\{y \in Q_{n,rn} : y \cdot  e_1 = n\} & \text{if } i=1, \\
				\{y \in Q_{n,rn} : y \cdot  e_i = rn\} & \text{if } 2 \leq i \leq d, \\
				\{y \in Q_{n,rn} : y \cdot  e_{-i} = -rn\} & \text{if } -d \leq i \leq -2.
			\end{cases}
		\end{align*}
		The union of all these faces is denoted by
		\begin{align*}
			\widehat{\partial Q}_{n,rn} := \bigcup_{\substack{1 \leq \vert i \vert \leq d \\ i \neq -1}} \partial_i Q_{n,rn}.
		\end{align*}
		
		\begin{figure}[h!]
			\centering
			\includegraphics[width=0.8\textwidth]{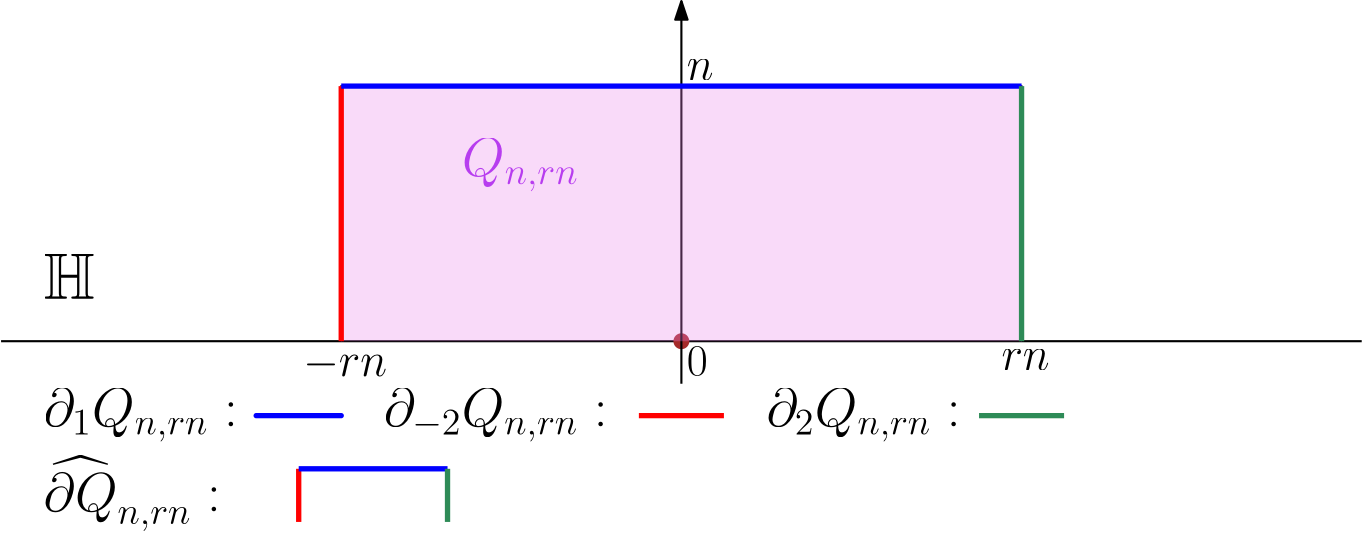}
			\caption{An illustration of the box $Q_{n,rn}$ and its faces.}\label{fig.box_faces}
		\end{figure}
		As illustrated in Figure~\ref{fig.box_faces}, $\widehat{\partial Q}_{n,rn}$ consists of every face of $Q_{n,rn}$ except the one lying on $\partial\H$.
		
		In this part, we prove a regularity estimate (Lemma~\ref{lem.boundary_regularity}) comparing the total mass of the two-point function on different faces of $Q_{n,rn}$, which will be used later.
		
		\begin{lemma}\label{lem.boundary_firstmoment}
			Let $r\geq1$. There exists $c_{\ref{lem.boundary_firstmoment}}>0$ such that, for every $n \in \N_+$,
			\begin{align}\label{eq.boundary_firstmoment}
				\sum_{y \in \widehat{\partial Q}_{n,rn}} \P_{\p_c}\Ll[0 \xleftrightarrow{Q_{n,rn}} y\Rr] \geq \frac{c_{\ref{lem.boundary_firstmoment}}}{n}.
			\end{align}
		\end{lemma}
		\begin{proof}
			Note that $Q_{n,rn} \subset \mathbb S_n$. Decomposing an open path from $0$ to $\partial\H_{-n}$ within $\mathbb S_n$ at its first intersection with $\widehat{\partial Q}_{n,rn}$, and applying the BK inequality,
			\begin{align}\label{eq.firstmoment_step1}
				\sum_{z \in \partial \H_{-n}} \P_{\p_c}\Ll[0 \xleftrightarrow{\mathbb S_n} z\Rr] \leq \Ll(\sum_{y\in\widehat{\partial Q}_{n,rn}} \P_{\p_c}\Ll[0 \xleftrightarrow{Q_{n,rn}} y\Rr]\Rr) \cdot \sup_{y \in Q_{n,rn}} \sum_{z \in \partial \H_{-n}} \P_{\p_c}\Ll[y \xleftrightarrow{\mathbb S_n} z\Rr].
			\end{align}
			
			We bound the supremum. Fix $y \in Q_{n,rn}$, and let $h:=y\cdot e_1\in[0,n]$. Let $\breve\H_n := \{x\in\Zd : x\cdot e_1 \leq n\} \supset \mathbb S_n$; monotonicity in the domain gives $\P_{\p_c}[y \xleftrightarrow{\mathbb S_n} z] \leq \P_{\p_c}[y \xleftrightarrow{\breve\H_n} z]$ for every $z \in \partial\H_{-n}=\partial\breve\H_n$. The map $R(x):=n e_1-x$ is a lattice isometry sending $\breve\H_n$ onto $\H$, $\partial\breve\H_n$ onto $\partial\H$, and $y$ to a point $R(y)$ at height $n-h$; it therefore preserves connection probabilities, so
			\begin{align*}
				\sum_{z \in \partial\breve\H_n} \P_{\p_c}\Ll[y \xleftrightarrow{\breve\H_n} z\Rr] = \sum_{z'\in\partial\H} \P_{\p_c}\Ll[R(y) \xleftrightarrow{\H} z'\Rr] \overset{\eqref{eq.psi_identity}}{=} \psi_{\p_c}(\H_{n-h}).
			\end{align*}
			Since $h\in[0,n]$, we have $n-h\in[0,n]$, so \eqref{eq.expansion_2} gives $\psi_{\p_c}(\H_{n-h}) \lesssim 1$, uniformly in $y\in Q_{n,rn}$. This shows that the supremum in \eqref{eq.firstmoment_step1} is at most $1$ up to a constant.
			
			Combining \eqref{eq.firstmoment_step1} with $\sum_{z\in\partial\H_{-n}} \P_{\p_c}[0 \xleftrightarrow{\mathbb S_n} z] \geq c_{\ref{lem.Nn_moments}}/n$ from Lemma~\ref{lem.Nn_moments}, we obtain \eqref{eq.boundary_firstmoment}.
		\end{proof}
		
		\begin{lemma}\label{lem.boundary_regularity}
			Let $r>1$. There exists a finite positive constant $C_{\ref{lem.boundary_regularity}}(r)$ such that the following estimate holds for all $n \in \N_+$
			\begin{align}\label{eq.boundary_regularity}
				\forall i \text{ with } 2 \leq \vert i \vert \leq d, \qquad \sum_{y \in \partial_i Q_{n,rn}} \P_{\p_c}\Ll[0 \xleftrightarrow{Q_{n,rn}} y\Rr] \leq C_{\ref{lem.boundary_regularity}} \sum_{y \in \partial_1 Q_{n,rn}} \P_{\p_c}\Ll[0 \xleftrightarrow{Q_{n,rn}} y\Rr].
			\end{align}
			As a consequence, there exists $c_{\ref{lem.boundary_regularity}}(r)>0$ satisfying
			\begin{align}\label{eq.boundary_0_mean_segment}
				\sum_{y \in \partial_1 Q_{n,rn}} \P_{\p_c}\Ll[0 \xleftrightarrow{Q_{n,rn}} y\Rr] \geq \frac{c_{\ref{lem.boundary_regularity}}}{n}.
			\end{align}
		\end{lemma}
		\begin{proof}
			Throughout this proof, write $Q\equiv Q_{n,rn}$ and $\partial_i Q\equiv\partial_i Q_{n,rn}$. By the symmetry of $Q$ among the directions $\pm 2,\ldots,\pm d$, it suffices to treat the case $i=2$. Let
			\begin{align*}
				S := \{x \in Q : x \cdot  e_2 = \lfloor (r-1)n \rfloor \}, \qquad Q' := \{x \in Q : x \cdot  e_2 \geq \lfloor (r-1)n \rfloor\},
			\end{align*}
			so that $S$ is the hyperplane in $Q$ at distance $n$ from $\partial_2 Q$, and $Q'$ is the part of $Q$ beyond $S$. Decomposing an open path realizing $0 \xleftrightarrow{Q} y$ (for $y \in \partial_2 Q$) at its last intersection $z$ with $S$, and applying the BK inequality,
			\begin{align*}
				\sum_{y \in \partial_2 Q} \P_{\p_c}\Ll[0 \xleftrightarrow{Q} y\Rr] \leq \sum_{z \in S} \P_{\p_c}\Ll[0 \xleftrightarrow{Q} z\Rr] \sum_{y \in \partial_2 Q} \P_{\p_c}\Ll[z \xleftrightarrow{Q'} y\Rr].
			\end{align*}
			Since $r>1$, $Q'$ is contained in the translate $Q_z$ of $Q$ obtained by placing $Q$'s near face at $S$ and pointing $Q$'s height-direction along $ e_2$ instead of $ e_1$. Monotonicity in the domain and translation/rotation invariance of $\P_{\p_c}$ then give, for every $z \in S$,
			\begin{align}\label{eq.boundaryCrossingRotation}
				\sum_{y \in \partial_2 Q} \P_{\p_c}\Ll[z \xleftrightarrow{Q'} y\Rr] \leq \sum_{y \in \partial_2 Q} \P_{\p_c}\Ll[z \xleftrightarrow{Q_z} y\Rr] = \sum_{y \in \partial_1 Q} \P_{\p_c}\Ll[0 \xleftrightarrow{Q} y\Rr].
			\end{align}
			\begin{figure}[h!]
				\centering
				\includegraphics[width=0.9\textwidth]{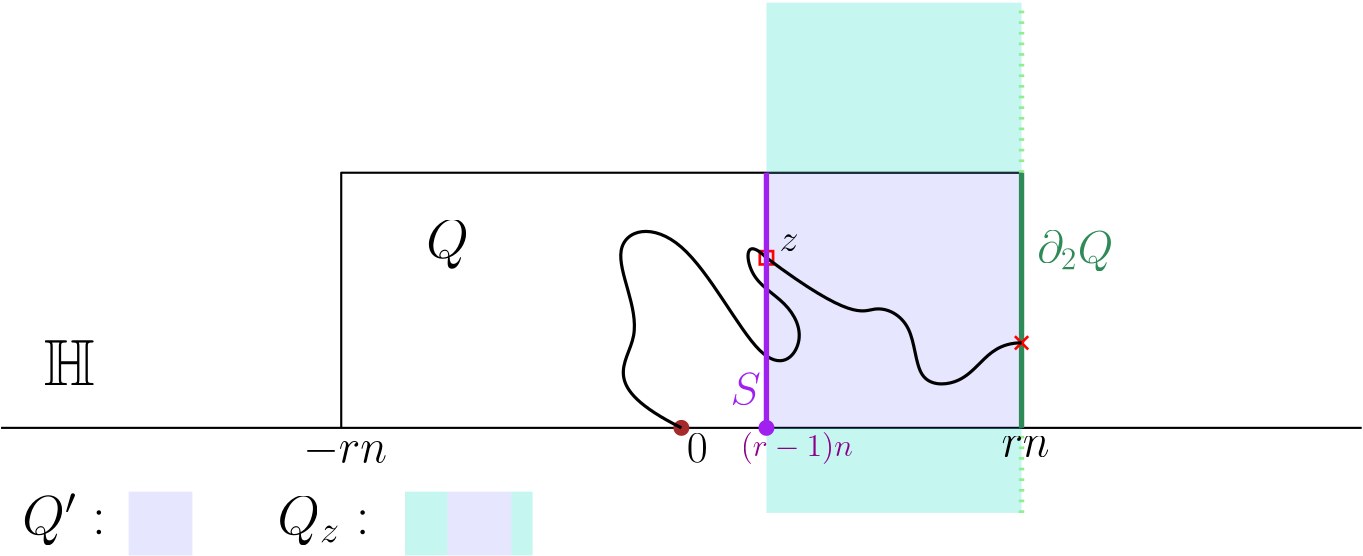}
				\caption{An illustration of the argument in \eqref{eq.boundaryCrossingRotation}. The path starts at $0$ and the small red square marks the point $z$, the last intersection point of the path with $S$.}\label{fig.boundaryCrossingRotation}
			\end{figure}
			Combining the last two displays,
			\begin{align}\label{eq.boundary_regularity_middle}
				\sum_{y \in \partial_2 Q} \P_{\p_c}\Ll[0 \xleftrightarrow{Q} y\Rr] \leq \Ll(\sum_{y \in \partial_1 Q} \P_{\p_c}\Ll[0 \xleftrightarrow{Q} y\Rr]\Rr) \cdot \Ll(\sum_{z \in S} \P_{\p_c}\Ll[0 \xleftrightarrow{Q} z\Rr]\Rr).
			\end{align}
			It remains to bound $\sum_{z\in S} \P_{\p_c}[0 \xleftrightarrow{Q} z]$ by a constant. For every $z \in S$, we have $z \cdot  e_1 \in [0,n]$ and $\vert z \vert \geq z \cdot  e_2 = \lfloor (r-1)n \rfloor$. Then Lemma~\ref{lem.halfspace_twopoint} gives, uniformly for $z \in S$,
			\begin{align*}
				\P_{\p_c}\Ll[0 \xleftrightarrow{\H} z\Rr] \leq C_{\ref{lem.halfspace_twopoint}}\, \frac{1+n}{1+\vert z\vert^d} \leq C_1(r)\, \frac{n^{1-d}}{(r-1)^d}
			\end{align*}
			for some $C_1(r)>0$. Since $S$ has one coordinate fixed and the remaining $d-1$ coordinates ranging over $[0,n]$ or $(-rn,rn)$, there is $C_2(r)>0$ such that $\vert S \vert \leq C_2\,n^{d-1}$. Therefore, relaxing $\P_{\p_c}[0 \xleftrightarrow{Q} z] \leq \P_{\p_c}[0 \xleftrightarrow{\H} z]$,
			\begin{align*}
				\sum_{z \in S} \P_{\p_c}\Ll[0 \xleftrightarrow{Q} z\Rr] \leq \sum_{z \in S} \P_{\p_c}\Ll[0 \xleftrightarrow{\H} z\Rr] \leq C_2\,n^{d-1} \cdot C_1\, \frac{n^{1-d}}{(r-1)^d} = \frac{C_1C_2}{(r-1)^d} =: C_3(r),
			\end{align*}
			which is a finite constant depending only on $r,d$. Plugging this into \eqref{eq.boundary_regularity_middle} and setting $C_{\ref{lem.boundary_regularity}}(r):=C_3$ yields \eqref{eq.boundary_regularity}.
			
			Finally, summing \eqref{eq.boundary_regularity} over $i$ with $2\leq\vert i\vert\leq d$ and combining with Lemma~\ref{lem.boundary_firstmoment} (with constant $c_{\ref{lem.boundary_firstmoment}}(r)$),
			\begin{align*}
				\frac{c_{\ref{lem.boundary_firstmoment}}}{n} \leq \sum_{y \in \widehat{\partial Q}} \P_{\p_c}\Ll[0 \xleftrightarrow{Q} y\Rr] \leq \Ll(2(d-1)C_{\ref{lem.boundary_regularity}}+1\Rr) \sum_{y \in \partial_1 Q} \P_{\p_c}\Ll[0 \xleftrightarrow{Q} y\Rr],
			\end{align*}
			which gives \eqref{eq.boundary_0_mean_segment} with $c_{\ref{lem.boundary_regularity}}(r) := c_{\ref{lem.boundary_firstmoment}}/\big(2(d-1)C_{\ref{lem.boundary_regularity}}+1\big)$.
		\end{proof}
		
		We now turn to a variant estimate, in which $Q_{n,rn}$ is replaced by a cube $x+\Lambda_n$ centered at a point $x\in\H$ that need not lie on $\partial\H$. For $x\in\H$ and $i\in\{1,\ldots,d\}\cup\{-2,\ldots,-d\}$, define the faces
		\begin{align*}
			\partial_i(x+\Lambda_n) := \begin{cases}
				\{y \in (x+\Lambda_n)\cap\H : y \cdot  e_1 = x\cdot e_1+n\} & \text{if } i=1, \\
				\{y \in (x+\Lambda_n)\cap\H : y \cdot  e_i = x\cdot e_i+n\} & \text{if } 2 \leq i \leq d, \\
				\{y \in (x+\Lambda_n)\cap\H : y \cdot  e_{-i} = x\cdot e_{-i}-n\} & \text{if } -d \leq i \leq -2,
			\end{cases}
		\end{align*}
		and $\widehat{\partial(x+\Lambda_n)} := \bigcup_{1\leq\vert i\vert\leq d,\, i\neq-1} \partial_i(x+\Lambda_n)$.
		
		\begin{lemma}\label{lem.Nn_mean_upper}
			Let $\mathcal N_n := \#\{x \in \partial \H_{-n} : 0 \xleftrightarrow{\mathbb S_n} x\}$, as in Lemma~\ref{lem.Nn_moments}. There exists $C_{\ref{lem.Nn_mean_upper}}>0$ such that, for every $n\geq1$,
			\begin{align*}
				\E_{\p_c}[\mathcal N_n] \leq \frac{C_{\ref{lem.Nn_mean_upper}}}{n}.
			\end{align*}
		\end{lemma}
		\begin{proof}
			Since $\{\mathcal N_n>0\} = \{0\xleftrightarrow{\mathbb S_n}\partial\H_{-n}\}$ and $\mathbb S_n\subset\H$, monotonicity in the domain gives $\{\mathcal N_n>0\}\subset\{0\xleftrightarrow{\H}\partial\H_{-n}\}$. Every $x\in\partial\H_{-n}$ satisfies $x\cdot e_1=n$, so $\vert x\vert\geq n>n-1$, i.e.\ $x\in\Lambda_{n-1}^c$; hence $\{0\xleftrightarrow{\H}\partial\H_{-n}\}\subset\{0\xleftrightarrow{\H}\Lambda_{n-1}^c\}$. For $n\geq2$, Lemma~\ref{lem.onearm_half} applied at $n-1$ gives, with $C_{\ref{lem.onearm_half}}>0$ as in that lemma,
			\begin{align*}
				\P_{\p_c}[\mathcal N_n>0] \leq \P_{\p_c}\Ll[0\xleftrightarrow{\H}\Lambda_{n-1}^c\Rr] \leq \frac{C_{\ref{lem.onearm_half}}}{(n-1)^3} \leq \frac{8C_{\ref{lem.onearm_half}}}{n^3},
			\end{align*}
			using $n-1\geq n/2$ for the last inequality; the case $n=1$ is trivial since $\P_{\p_c}[\mathcal N_1>0]\leq1$. Combining this with the Cauchy--Schwarz inequality,
			\begin{align*}
				\E_{\p_c}[\mathcal N_n] = \E_{\p_c}\Ll[\mathcal N_n\Ind{\mathcal N_n>0}\Rr] \leq \E_{\p_c}[\mathcal N_n^2]^{1/2}\,\P_{\p_c}[\mathcal N_n>0]^{1/2},
			\end{align*}
			and using $\E_{\p_c}[\mathcal N_n^2]\leq C_{\ref{lem.Nn_moments}}\,n$ from Lemma~\ref{lem.Nn_moments}, we get
			\begin{align*}
				\E_{\p_c}[\mathcal N_n] \leq \big(C_{\ref{lem.Nn_moments}}\,n\big)^{1/2}\Ll(\frac{8C_{\ref{lem.onearm_half}}}{n^3}\Rr)^{1/2} = \big(8C_{\ref{lem.Nn_moments}}C_{\ref{lem.onearm_half}}\big)^{1/2}\, n^{-1},
			\end{align*}
			which proves the claim with $C_{\ref{lem.Nn_mean_upper}}:=\big(8C_{\ref{lem.Nn_moments}}C_{\ref{lem.onearm_half}}\big)^{1/2}$.
		\end{proof}
		
		\begin{lemma}\label{lem.cube_face_upper_bound}
			There exists $C_{\ref{lem.cube_face_upper_bound}}>0$ such that, for every $n\in\N_+$ and every $x\in\H$ with $x\cdot e_1 =: k \in [0,n]$,
			\begin{align*}
				\sum_{y\in\widehat{\partial(x+\Lambda_n)}} \P_{\p_c}\Ll[x \xleftrightarrow{(x+\Lambda_n)\cap\H} y\Rr] \leq C_{\ref{lem.cube_face_upper_bound}}\, \frac{1+k}{n}.
			\end{align*}
		\end{lemma}
		
		\begin{proof}
			Fix $i$ with $1\leq\vert i\vert\leq d$, $i\neq-1$; write $\mathbf u :=  e_i$ if $i>0$ and $\mathbf u:=- e_{-i}$ if $i<0$, so that $F:=\partial_i(x+\Lambda_n) = \{y\in(x+\Lambda_n)\cap\H : y\cdot\mathbf u = x\cdot\mathbf u+n\}$. Let $m:=\lfloor n/2\rfloor$ and
			\begin{align*}
				M := \{z\in(x+\Lambda_n)\cap\H : z\cdot\mathbf u = x\cdot\mathbf u+m\}.
			\end{align*}
			Decomposing an open path realizing $x\xleftrightarrow{(x+\Lambda_n)\cap\H} y$ (for $y\in F$) at its last intersection $z$ with $M$, and applying the BK inequality,
			\begin{align}\label{eq.cube_step1}
				\sum_{y\in F} \P_{\p_c}\Ll[x\xleftrightarrow{(x+\Lambda_n)\cap\H} y\Rr] \leq \Ll(\sum_{z\in M} \P_{\p_c}\Ll[x\xleftrightarrow{\H} z\Rr]\Rr) \cdot \sup_{z\in M}\sum_{y\in F} \P_{\p_c}\Ll[z\xleftrightarrow{\H_z^{\mathbf u}} y\Rr],
			\end{align}
			where $\H_z^{\mathbf u} := \{w\in\Zd : (w-z)\cdot\mathbf u\geq0\}$; the domain has been relaxed from $(x+\Lambda_n)\cap\H$ to $\H$ in the first factor, and from $(x+\Lambda_n)\cap\H$ to $\H_z^{\mathbf u}$ in the second factor, both by monotonicity in the domain (since $F\subset\H_z^{\mathbf u}$ for every $z\in M$, as $F$ lies at $\mathbf u$-distance $n-m$ from $M$).
			
			\begin{figure}[h!]
				\centering
				\includegraphics[width=0.6\textwidth]{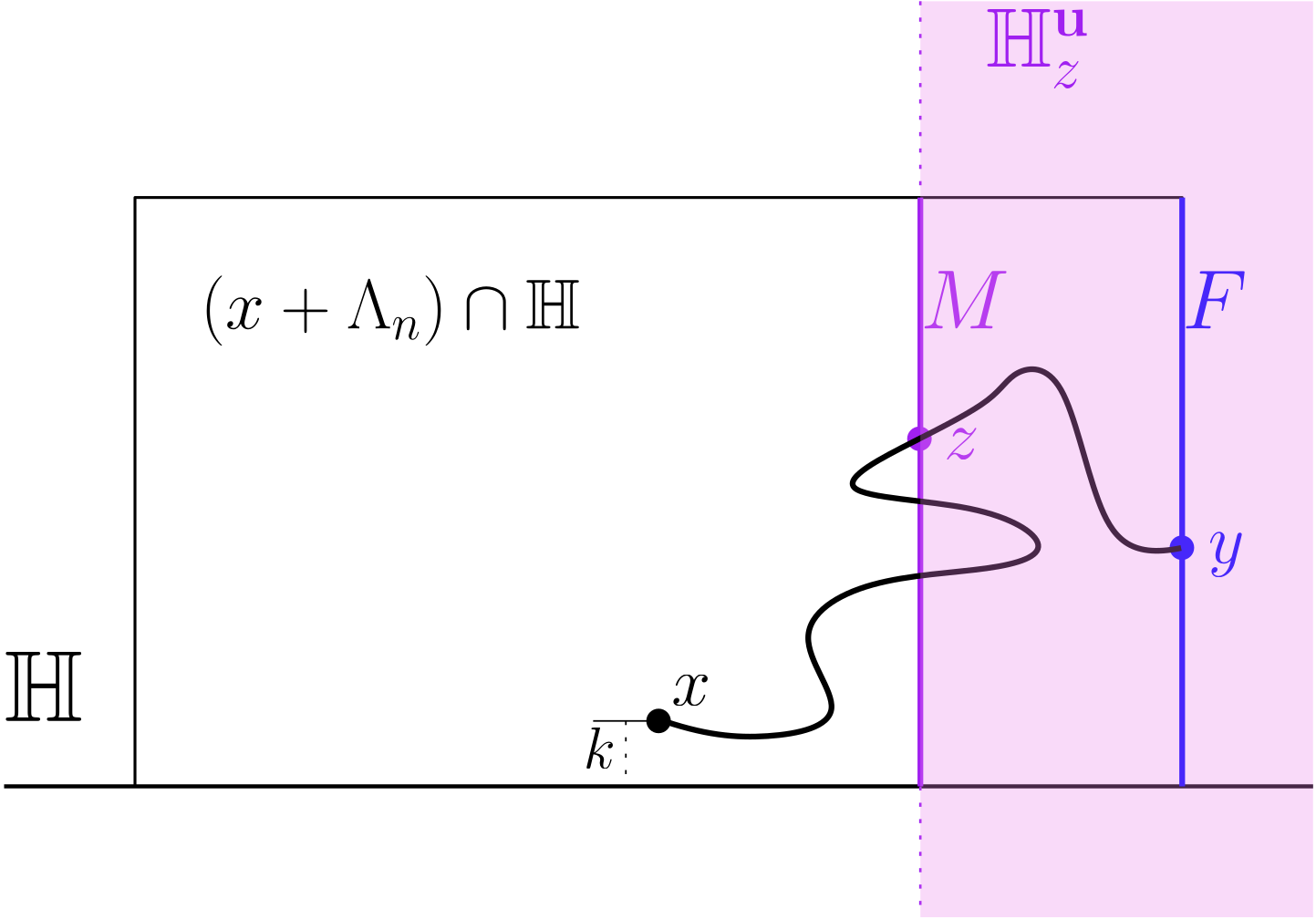}
				\caption{An illustration of the decomposition in \eqref{eq.cube_step1}.}\label{fig.cube_face_decomposition}
			\end{figure}
			See Figure~\ref{fig.cube_face_decomposition} for an illustration of this decomposition.
			
			\emph{Bounding the first factor.} For $z\in M$, since $(z-x)\cdot\mathbf u = m$ and $z\in x+\Lambda_n$, we have $m\leq\vert x-z\vert\leq n$, so $\vert x-z\vert\geq m\geq(n-1)/2$. Writing $j:=z\cdot e_1$, Lemma~\ref{lem.halfspace_twopoint} gives, uniformly for $z\in M$,
			\begin{align*}
				\P_{\p_c}\Ll[x\xleftrightarrow{\H} z\Rr] \leq C_{\ref{lem.halfspace_twopoint}}\, \frac{(1+\min(k,\vert x-z\vert))(1+\min(j,\vert x-z\vert))}{1+\vert x-z\vert^d} \leq C_{\ref{lem.halfspace_twopoint}}\, \frac{(1+k)(1+n)}{1+((n-1)/2)^d}.
			\end{align*}
			Since $M$ has one coordinate fixed and the remaining $d-1$ coordinates ranging over intervals of length $O(n)$, there is $C_1>0$ such that $\vert M\vert \leq C_1\, n^{d-1}$. Combining the two previous displays, there is $C_2>0$ such that
			\begin{align*}
				\sum_{z\in M} \P_{\p_c}\Ll[x\xleftrightarrow{\H} z\Rr] \leq C_2\,(1+k).
			\end{align*}
			
			\emph{Bounding the second factor.} By the invariance of nearest-neighbor percolation on $\Zd$ under permutations and sign changes of the coordinate axes, $\H_z^{\mathbf u}$ is isomorphic, via a lattice isometry sending $z$ to $0$, to the standard half-space $\H$; under this isomorphism $F$ corresponds to the hyperplane at height $n-m$ from $0$. Hence
			\begin{align*}
				\sum_{y\in F} \P_{\p_c}\Ll[z\xleftrightarrow{\H_z^{\mathbf u}} y\Rr] \leq \sum_{w\in\H:\, w\cdot e_1=n-m} \P_{\p_c}\Ll[0\xleftrightarrow{\H} w\Rr] = \E_{\p_c}[\mathcal N_{n-m}],
			\end{align*}
			where the inequality enlarges $F$ to the full hyperplane at that height. Since $m=\lfloor n/2\rfloor$, we have $n-m\geq n/2$, so Lemma~\ref{lem.Nn_mean_upper} gives, uniformly for $z\in M$,
			\begin{align*}
				\sum_{y\in F} \P_{\p_c}\Ll[z\xleftrightarrow{\H_z^{\mathbf u}} y\Rr] \leq \E_{\p_c}[\mathcal N_{n-m}] \leq \frac{C_{\ref{lem.Nn_mean_upper}}}{n-m} \leq \frac{2C_{\ref{lem.Nn_mean_upper}}}{n}.
			\end{align*}
			
			\emph{Conclusion.} Plugging the two bounds into \eqref{eq.cube_step1},
			\begin{align*}
				\sum_{y\in F} \P_{\p_c}\Ll[x\xleftrightarrow{(x+\Lambda_n)\cap\H} y\Rr] \leq C_2\,(1+k) \cdot \frac{2C_{\ref{lem.Nn_mean_upper}}}{n} = \frac{2C_2C_{\ref{lem.Nn_mean_upper}}(1+k)}{n}.
			\end{align*}
			Since this bound does not depend on the choice of $i$, summing over the $2(d-1)+1$ faces $i\in\{1,\ldots,d\}\cup\{-2,\ldots,-d\}$ gives
			\begin{align*}
				\sum_{y\in\widehat{\partial(x+\Lambda_n)}} \P_{\p_c}\Ll[x\xleftrightarrow{(x+\Lambda_n)\cap\H} y\Rr] \leq (2d-1)\cdot\frac{2C_2C_{\ref{lem.Nn_mean_upper}}(1+k)}{n},
			\end{align*}
			which proves the lemma with $C_{\ref{lem.cube_face_upper_bound}} := 2(2d-1)C_2 C_{\ref{lem.Nn_mean_upper}}$.
		\end{proof}

		Now we compute the 3-point function, which will be useful.
		Write $\tau(x,y):=\P_{\p_c}[x\overset{\Zd}{\longleftrightarrow}y]$ for the two-point function on $\Zd$.
		
		\begin{lemma}\label{lem.threepoint_upper}
			Assume that $\tau(x,y)\lesssim \vert x-y\vert^{2-d}$ for every $x\neq y\in\Zd$. There exists $C_{\ref{lem.threepoint_upper}}>0$ such that, for pairwise distinct $x,y,z\in\Zd$,
			\begin{equation*}
				\P_{\p_c}[x\leftrightarrow y,\ x\leftrightarrow z] \leq C_{\ref{lem.threepoint_upper}}\,\min\{\vert x-y\vert,\vert y-z\vert,\vert x-z\vert\}^{4-d}\,\max\{\vert x-y\vert,\vert y-z\vert,\vert x-z\vert\}^{2-d}.
			\end{equation*}
		\end{lemma}
		\begin{proof}
			Let $a\leq b\leq c$ denote $\vert x-y\vert,\vert x-z\vert,\vert y-z\vert$ sorted increasingly; the triangle inequality gives $b\leq c\leq2b$. By the BK inequality, decomposing at the branch point of a spanning tree of the cluster of $x$,
			\begin{equation*}
				\P_{\p_c}[x\leftrightarrow y,\ x\leftrightarrow z] \leq \sum_{u\in\Zd}\tau(x,u)\tau(u,y)\tau(u,z).
			\end{equation*}
			Partition $\Zd$ by which of $x,y,z$ is nearest to $u$. On the parts nearest $x$ or nearest $y$, the triangle inequality gives $\vert u-z\vert\geq\vert x-z\vert/2=b/2$, resp.\ $\vert u-z\vert\geq\vert y-z\vert/2=c/2$, so the hypothesis gives $\tau(u,z)\lesssim  b^{2-d}$ throughout either part, whence these two parts contribute
			\begin{equation*}
				\lesssim  b^{2-d}\sum_{u\in\Zd}\tau(x,u)\tau(u,y) \lesssim  b^{2-d}\,a^{4-d},
			\end{equation*}
			using the standard bubble bound $\sum_u\tau(x,u)\tau(u,y)\lesssim (1+\vert x-y\vert)^{4-d}$ ($d>4$; this follows from the hypothesis by splitting at $\vert x-u\vert\lessgtr a/2$ and using $\#\{u:\vert x-u\vert=k\}\lesssim k^{d-1}$). On the part nearest $z$, the triangle inequality gives $\vert x-u\vert\geq\vert x-z\vert/2=b/2$, so the hypothesis gives $\tau(x,u)\lesssim  b^{2-d}$ throughout, and this part contributes $\lesssim  b^{2-d}\sum_u\tau(u,y)\tau(u,z)\lesssim  b^{2-d}c^{4-d}\leq b^{2-d}a^{4-d}$, using $c\geq a$. Summing the three parts, and recalling $b\asymp\max\{a,b,c\}$, proves the lemma.
		\end{proof}

		\subsection{Upper bound for the 2-point function in cubes and annulus}\label{subsct:Upper bound for the 2-point function in cubes and annulus}
		In this subsection, we give upper bounds for some 2-point functions in cubes and annulus, which are foundations of our computation.
		
		Fix $\alpha_i,\beta_i>0$ for $1\leq i\leq d$, and for each $n$ let
		\begin{align*}
			R_n := [-\alpha_1n,\beta_1n]\times\cdots\times[-\alpha_dn,\beta_dn]
		\end{align*}
		as in Lemma~\ref{lem.box_twopoint}. For $x\in R_n$, let
		\begin{align*}
			l_x := \min_{1\leq i\leq d}\min\{\beta_in-x\cdot e_i,\ x\cdot e_i+\alpha_in\}
		\end{align*}
		denote the minimal distance between $x$ and $\partial R_n$ (when $\alpha_i=\beta_i=1$ for all $i$, $R_n=\Lambda_n$ and this recovers $l_x=n-\vert x\vert$). We first record two general upper bounds for the 2-point function inside a single box $R_n$, depending on the boundary distance of one or both endpoints; the constants below depend only on $d$, uniformly over all choices of $(\alpha_i),(\beta_i)$.
		
		\begin{lemma}\label{lem.box_twopoint_onesided}
			There exists $C_{\ref{lem.box_twopoint_onesided}}>0$ such that, for every $n\in\N_+$, every $\alpha_i,\beta_i>0$ ($1\leq i\leq d$), and every $x,y\in R_n$ with $x\neq y$,
			\begin{align*}
				\P_{\p_c}\Ll[x\overset{R_n}\leftrightarrow y\Rr] \leq C_{\ref{lem.box_twopoint_onesided}}\,(1+l_x)\,\vert x-y\vert^{1-d}.
			\end{align*}
		\end{lemma}
		\begin{proof}
			Without loss of generality, we may assume that the face of $R_n$ nearest to $x$ is $\{z\cdot e_1=\beta_1n\}$. Then $R_n\subset\H_x:=\{z\in\Zd:z\cdot e_1\leq\beta_1n\}$, and monotonicity in the domain gives $\P_{\p_c}[x\overset{R_n}\leftrightarrow y]\leq\P_{\p_c}[x\overset{\H_x}\leftrightarrow y]$. Lemma~\ref{lem.halfspace_twopoint} then gives, using $\min(l_x,\vert x-y\vert)\leq l_x$ and $\min(\cdot,\vert x-y\vert)\leq\vert x-y\vert$ for the other factor,
			\begin{align*}
				\P_{\p_c}\Ll[x\overset{\H_x}\leftrightarrow y\Rr] \leq C_{\ref{lem.halfspace_twopoint}}\,\frac{(1+l_x)(1+\vert x-y\vert)}{1+\vert x-y\vert^d} \leq 2C_{\ref{lem.halfspace_twopoint}}\,(1+l_x)\,\vert x-y\vert^{1-d},
			\end{align*}
			using $1+\vert x-y\vert\leq2\vert x-y\vert$ (as $\vert x-y\vert\geq1$). This proves the lemma with $C_{\ref{lem.box_twopoint_onesided}}:=2C_{\ref{lem.halfspace_twopoint}}$.
		\end{proof}
		
		\begin{lemma}\label{lem.box_twopoint_twosided}
			There exists $C_{\ref{lem.box_twopoint_twosided}}>0$ such that, for every $n\in\N_+$, every $\alpha_i,\beta_i>0$ ($1\leq i\leq d$), and every $x,y\in R_n$ with $x\neq y$,
			\begin{align*}
				\P_{\p_c}\Ll[x\overset{R_n}\leftrightarrow y\Rr] \leq C_{\ref{lem.box_twopoint_twosided}}\,(1+l_x)(1+l_y)\,\vert x-y\vert^{-d}.
			\end{align*}
		\end{lemma}
		\begin{proof}
			By symmetry between $x$ and $y$, we may assume $l_x\leq l_y$.
			
			\emph{Case 1: $\vert x-y\vert\leq2l_y$.} Lemma~\ref{lem.box_twopoint_onesided} gives
			\begin{align*}
				\P_{\p_c}\Ll[x\overset{R_n}\leftrightarrow y\Rr] \leq C_{\ref{lem.box_twopoint_onesided}}\,(1+l_x)\,\vert x-y\vert^{1-d} \leq 2C_{\ref{lem.box_twopoint_onesided}}\,(1+l_x)(1+l_y)\,\vert x-y\vert^{-d}.
			\end{align*}
			
			\emph{Case 2: $\vert x-y\vert>2l_y\geq2l_x$.} Let $\rho:=\lfloor\vert x-y\vert/2\rfloor\geq l_y\geq l_x$. Since $\vert x-y\vert\geq2\rho>\rho$, we have $y\notin x+\Lambda_\rho$, so every open path inside $R_n$ from $x$ to $y$ must exit $(x+\Lambda_\rho)\cap R_n$; decomposing at its first such exit point $z$ and applying the BK inequality,
			\begin{align*}
				\P_{\p_c}\Ll[x\overset{R_n}\leftrightarrow y\Rr] \leq \Ll(\sum_{z\in\partial(x+\Lambda_\rho)\cap R_n} \P_{\p_c}\Ll[x\overset{(x+\Lambda_\rho) \cap R_n}\longleftrightarrow z\Rr]\Rr) \cdot \sup_{z\in\partial(x+\Lambda_\rho)\cap R_n} \P_{\p_c}\Ll[z\overset{R_n}\leftrightarrow y\Rr].
			\end{align*}
			
			\begin{figure}[h!]
				\centering
				\includegraphics[width=0.8\textwidth]{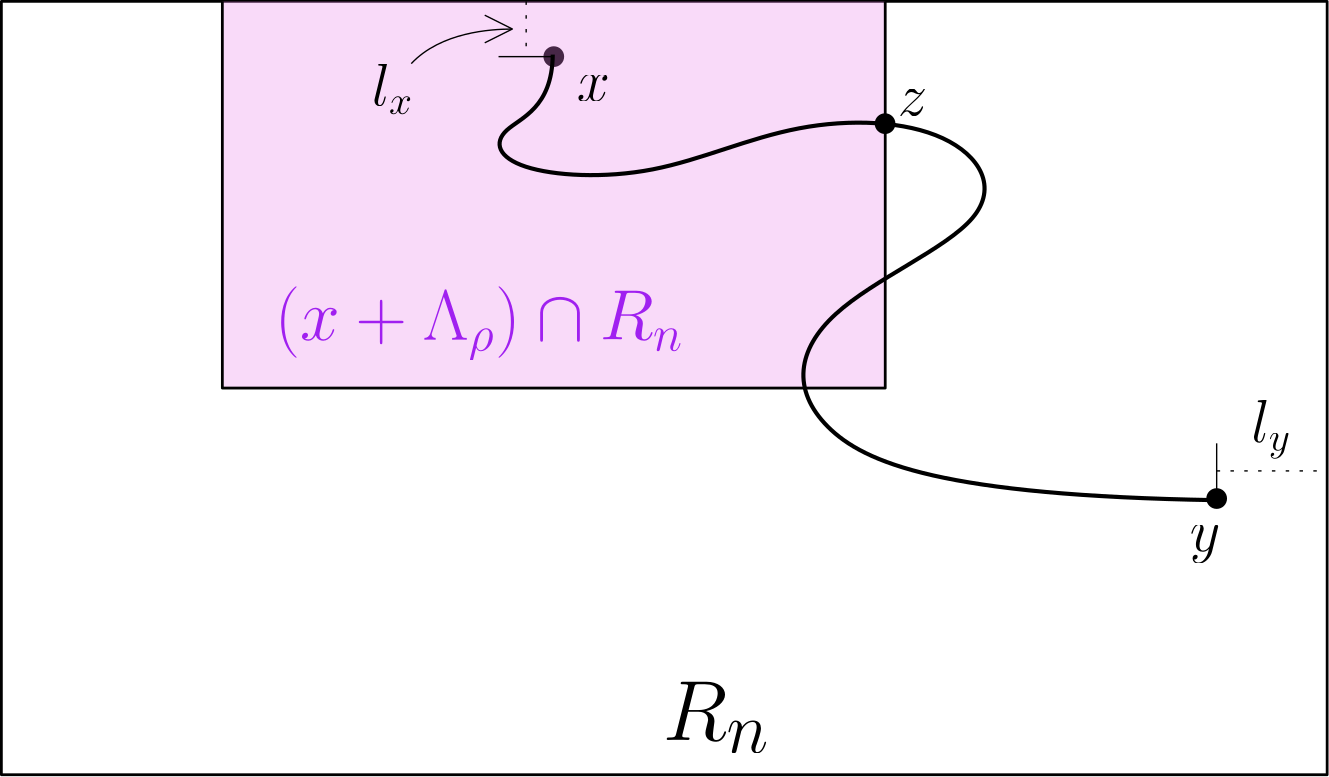}
				\caption{A illustration of Case~2 of the proof of Lemma~\ref{lem.box_twopoint_twosided}.}\label{fig.box_twopoint_twosided}
			\end{figure}
			See Figure~\ref{fig.box_twopoint_twosided} for an illustration of this decomposition.
			
			For the sum, $R_n\subset\H_x$ together with Lemma~\ref{lem.cube_face_upper_bound}  gives
			\begin{align*}
				\sum_{z\in\partial(x+\Lambda_\rho)\cap R_n} \P_{\p_c}\Ll[x\overset{(x+\Lambda_\rho)\cap R_n}\longleftrightarrow z\Rr] \leq \sum_{z\in\widehat{\partial(x+\Lambda_\rho)}} \P_{\p_c}\Ll[x\overset{(x+\Lambda_\rho)\cap\H_x}\longleftrightarrow z\Rr] \leq C_{\ref{lem.cube_face_upper_bound}}\,\frac{1+l_x}{\rho}.
			\end{align*}
			For the supremum, every $z\in\partial(x+\Lambda_\rho)$ satisfies $\vert x-z\vert=\rho$, so
			\begin{align*}
				\vert y-z\vert \geq \vert x-y\vert-\vert x-z\vert \geq 2\rho-\rho = \rho \geq \vert x-y\vert/2-1;
			\end{align*}
			Lemma~\ref{lem.box_twopoint_onesided} then gives, for some $C_1>0$,
			\begin{align*}
				\P_{\p_c}\Ll[z\overset{R_n}\leftrightarrow y\Rr] \leq C_{\ref{lem.box_twopoint_onesided}}\,(1+l_y)\,\vert y-z\vert^{1-d} \leq C_1\,(1+l_y)\,\vert x-y\vert^{1-d}.
			\end{align*}
			Combining the two bounds and using $\rho\geq\vert x-y\vert/2-1\geq\vert x-y\vert/4$ (as $\vert x-y\vert\geq2$ in this case),
			\begin{align*}
				\P_{\p_c}\Ll[x\overset{R_n}\leftrightarrow y\Rr] \leq C_{\ref{lem.cube_face_upper_bound}}\,\frac{1+l_x}{\rho}\cdot C_1\,(1+l_y)\,\vert x-y\vert^{1-d} \leq C_2\,(1+l_x)(1+l_y)\,\vert x-y\vert^{-d},
			\end{align*}
			for some $C_2>0$. Combining both cases proves the lemma.
		\end{proof}
		
\begin{figure}[h]
			\centering
			\includegraphics[width=0.8\textwidth]{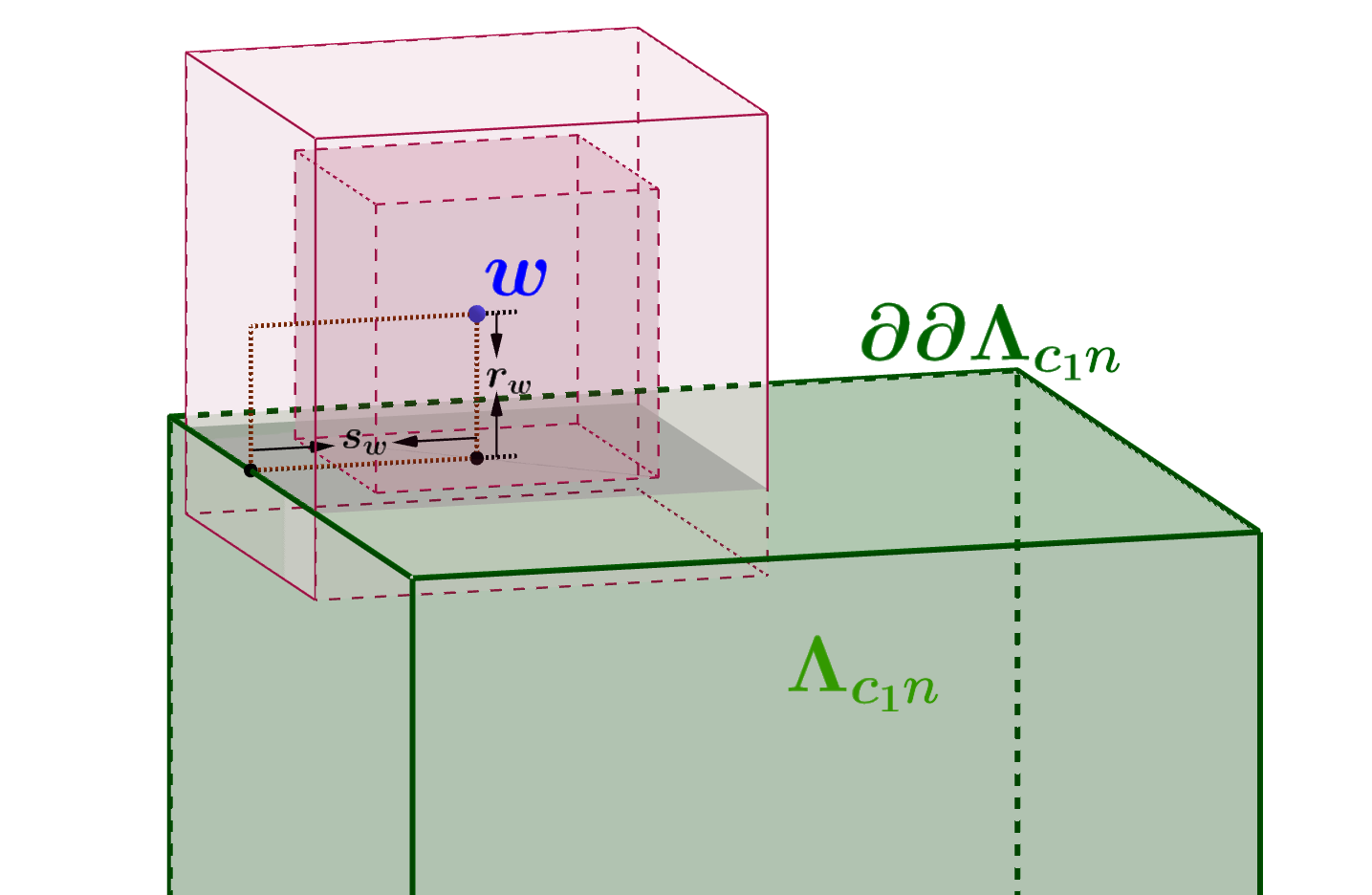} 
			\caption{Illustration of the condition $r_w < s_w$. The distances $r_w$ and $s_w$ correspond to the radii of the smallest $L^\infty$ balls (indicated by the pink cubes) centered at $w$ that intersect $\partial \Lambda_{c_1 n}$ and $\partial \partial \Lambda_{c_1 n}$, respectively. TODO: relabel this figure, replacing every occurrence of $L$ by $n$ (i.e.\ $\Lambda_{c_1L}\to\Lambda_{c_1n}$ and $\partial\partial\Lambda_{c_1L}\to\partial\partial\Lambda_{c_1n}$); $w$, $r_w$, $s_w$ are unchanged.}
			\label{figebr}
		\end{figure}

\begin{figure}[h]
			\centering
			\includegraphics[width=0.4\textwidth]{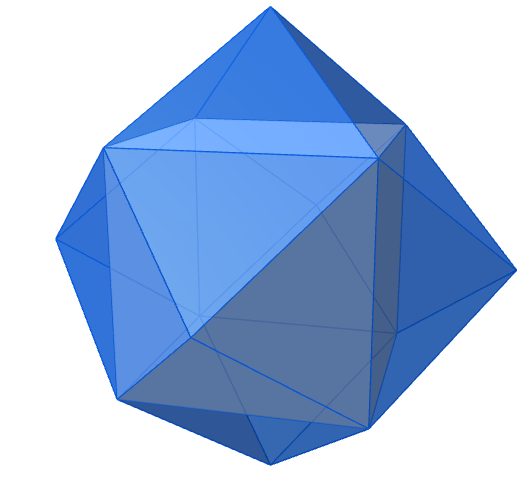} 
			\caption{Geometric interpretation of $\RD_{c_1 n}$. In general dimension $d$, $\RD_{c_1 n}$ is the union of $2d$ hyper-pyramids, each based on a face of the inner cube $\Lambda_{c_1 n}$ and reaching a height of $c_1 n$. This figure illustrates the $d=3$ case, where $\RD_{c_1 n}$ forms a solid rhombic dodecahedron with the interior of the central cube $\Lambda_{c_1 n}$ excluded.}
			\label{ktfu}
		\end{figure}

		For $c_1\in(0,1/5)$ and $n\in\N_+$, define
		\begin{equation*}
			\operatorname{int} \Lambda_{c_1 n} := \Lambda_{c_1n} \setminus \partial \Lambda_{c_1n}
		\end{equation*}
		and define the annulus
		\begin{align*}
			\mathrm{Ann}_{c_1n} := \Lambda_n\setminus \operatorname{int}\Lambda_{c_1n}
		\end{align*}
		For each $w\in\mathrm{Ann}_{c_1n}$, let $r_w$ be the minimal $\ell^\infty$ distance between $w$ and $\partial\Lambda_{c_1n}$, and let $s_w$ be the minimal $\ell^\infty$ distance between $w$ and $\partial\partial\Lambda_{c_1n}$, the set of points of $\partial\Lambda_{c_1n}$ at which at least two coordinates simultaneously attain absolute value $c_1n$ (i.e., the codimension-$2$ skeleton of $\partial\Lambda_{c_1n}$). Define
		\begin{equation}
			\begin{aligned}
				&\RD_{c_1n} := \{ w \in \mathrm{Ann}_{c_1n}: r_w < s_w \}\\
				&= \{ (x_1, \dots, x_d) \in \ZZ^d: \forall i,j \in \{ 1, \dots, d \}, (i \neq j \Rightarrow |x_i|+|x_j| <2c_1 n) \} \setminus \operatorname{int} \Lambda_{c_1 n}.
			\end{aligned}
		\end{equation}
		See Figures~\ref{figebr} and~\ref{ktfu} for the geometric visualizations of $r_w$, $s_w$, and the set $\RD_{c_1 n}$.
		
		
		For $w \in \mathrm{Ann}_{c_1n}$, write $a_{(1)}(w) \geq a_{(2)}(w) \geq \cdots \geq a_{(d)}(w)$ for the decreasing rearrangement of $\vert w\cdot e_1\vert,\ldots,\vert w\cdot e_d\vert$; since $w \notin \operatorname{int}\Lambda_{c_1n}$, we have $a_{(1)}(w) = \vert w\vert \geq c_1 n$, and
		\begin{align}\label{eq.rw_sw_formula}
			r_w = a_{(1)}(w) - c_1n, \qquad s_w = c_1n - a_{(2)}(w).
		\end{align}
		In particular, for $w \in \RD_{c_1n}$, the condition $r_w<s_w$ reads $a_{(1)}(w)+a_{(2)}(w) < 2c_1n$, so
		\begin{align}\label{eq.w_bound}
			\vert w \vert = a_{(1)}(w) < 2c_1n, \qquad \forall w \in \RD_{c_1n}.
		\end{align}
		
		For $x\in\Lambda_n$, let $l_x := n-\vert x\vert$ denote the mini
		mal $\ell^\infty$ distance between $x$ and $\partial\Lambda_n$.
		
		\begin{lemma}\label{lem.RD_bound}
			Let $c_1 \in (0,1/5)$. There exists $C_{\ref{lem.RD_bound}} >0$ such that, for every $n\in\N_+$, every $w \in \RD_{c_1 n}$, and every $x \in \Lambda_n\setminus\Lambda_{3n/4}$,
			\begin{equation}
				\P_{\p_c} \Ll[ x \stackrel{\mathrm{Ann}_{c_1n}}{\longleftrightarrow} w \Rr] \leq C_{\ref{lem.RD_bound}}\, n^{-d+1} (1+l_x)(1+r_w)\, s_w^{-1}.
			\end{equation}
		\end{lemma}
		\begin{proof}
			Without loss of generality, we may assume that $a_{(1)}(w)$ is attained at coordinate $1$ with a positive sign, i.e.\ using \eqref{eq.rw_sw_formula},
			\begin{align*}
				w\cdot e_1 = c_1n+r_w, \qquad \max_{2\leq j\leq d} \vert w\cdot e_j\vert = a_{(2)}(w) = c_1n-s_w.
			\end{align*}
			Write $\H^{(1)} := \{y\in\Zd : y\cdot e_1\geq c_1n\}$, a translate of $\H$. Since $c_1<1/5$, $\vert w\vert<2c_1n<2n/5$ and $s_w\leq c_1n<n/5$, so every $z\in w+\Lambda_{s_w}$ satisfies $\vert z\vert<3n/5$; since $\vert x\vert\geq3n/4$, this shows in particular $x\notin w+\Lambda_{s_w}$.
			
			From the following figure we can see
			\begin{align}\label{eq.local_halfspace_match}
				(w+\Lambda_{s_w})\cap\mathrm{Ann}_{c_1n} = (w+\Lambda_{s_w})\cap\H^{(1)}:
			\end{align}
			the box $w+\Lambda_{s_w}$ is too small to reach $\partial\partial\Lambda_{c_1n}$ (which sits at distance $s_w$ from $w$).
			
			\begin{figure}[h!]
				\centering
				\includegraphics[width=0.7\textwidth]{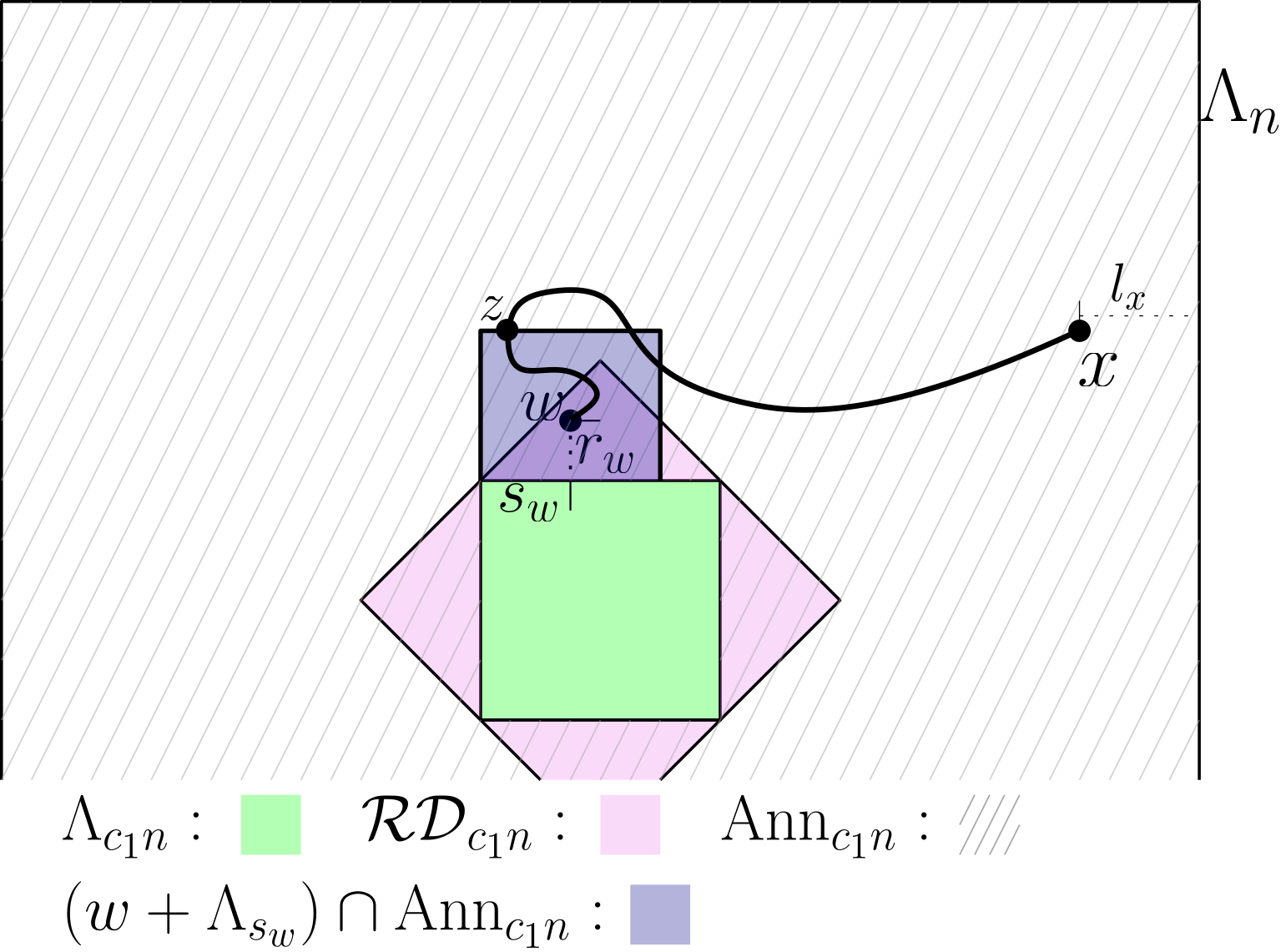}
				\caption{A illustration of the proof of Lemma~\ref{lem.RD_bound}.}\label{fig.RD_bound_decomposition}
			\end{figure}
			
			Since $\mathrm{Ann}_{c_1n}\subset\Lambda_n$, monotonicity in the domain together with Lemma~\ref{lem.box_twopoint_onesided} gives, for every $z\neq x$,
			\begin{align*}
				\P_{\p_c}\Ll[x\overset{\mathrm{Ann}_{c_1n}}\leftrightarrow z\Rr] \leq \P_{\p_c}\Ll[x\overset{\Lambda_n}\leftrightarrow z\Rr] \leq C_{\ref{lem.box_twopoint_onesided}}\,(1+l_x)\,\vert x-z\vert^{1-d}.
			\end{align*}
			Combined with $\vert x\vert\geq3n/4$ and $\vert z\vert<3n/5$, which give $\vert x-z\vert\geq\vert x\vert-\vert z\vert\geq3n/20$, this yields, for some $C_1>0$,
			\begin{align}\label{eq.reflection_bound_general}
				\P_{\p_c}\Ll[x\overset{\mathrm{Ann}_{c_1n}}\leftrightarrow z\Rr] \leq C_1\,(1+l_x)\, n^{1-d}, \qquad \forall z \text{ with } \vert z\vert<3n/5.
			\end{align}
			
			Decomposing an open path from $w$ to $x$ realizing $x\overset{\mathrm{Ann}_{c_1n}}\leftrightarrow w$ at its last exit $z$ from $w+\Lambda_{s_w}$, and applying the BK inequality together with \eqref{eq.reflection_bound_general},
			\begin{align}\label{eq.main_decomp_RD}
				\P_{\p_c}\Ll[x\overset{\mathrm{Ann}_{c_1n}}\leftrightarrow w\Rr] \leq C_1\,(1+l_x)\,n^{1-d} \cdot \sum_{z\in\partial(w+\Lambda_{s_w})\cap\mathrm{Ann}_{c_1n}} \P_{\p_c}\Ll[z\overset{(w+\Lambda_{s_w})\cap\mathrm{Ann}_{c_1n}}\leftrightarrow w\Rr].
			\end{align}
			By \eqref{eq.local_halfspace_match}, the sum ranges over $\partial(w+\Lambda_{s_w})\cap\H^{(1)}$. Then Lemma~\ref{lem.cube_face_upper_bound} gives
			\begin{align*}
				\sum_{z\in\partial(w+\Lambda_{s_w})\cap\mathrm{Ann}_{c_1n}} \P_{\p_c}\Ll[z\overset{(w+\Lambda_{s_w})\cap\mathrm{Ann}_{c_1n}}\leftrightarrow w\Rr] \leq C_2\,\frac{1+r_w}{s_w},
			\end{align*}
			for some $C_2>0$. Plugging this into \eqref{eq.main_decomp_RD} proves the lemma with $C_{\ref{lem.RD_bound}}:=C_1C_2$.
		\end{proof}

		\subsection{Lower bound for a 2-point function in cubes}
		
		In this subsection, we aim to prove Proposition~\ref{lem.box_twopoint_lower}. Recall that, for $x\in\Lambda_n$, $l_x := n-\vert x\vert$ denotes the minimal distance between $x$ and $\partial\Lambda_n$, and that $\partial\partial\Lambda_n$ denotes the set of points of $\partial\Lambda_n$ at which at least two coordinates simultaneously attain absolute value $n$.
		
		\begin{prop}\label{lem.box_twopoint_lower}
			For every $\alpha>0$ there exists $c_{\ref{lem.box_twopoint_lower}}(\alpha)>0$ such that the following holds for every $n\in\N_+$. Let $x\in\partial\Lambda_n$ be such that the distance from $x$ to $\partial\partial\Lambda_n$ is at least $\alpha n$, and let $y\in\Lambda_n$ be such that $l_y\geq\alpha n$. Then
			\begin{align*}
				\P_{\p_c}\Ll[x\overset{\Lambda_n}\leftrightarrow y\Rr] \geq \frac{c_{\ref{lem.box_twopoint_lower}}}{n^{d-1}}.
			\end{align*}
		\end{prop}
		
		Fix $\alpha>0$ and let $n,x,y$ be as in Proposition~\ref{lem.box_twopoint_lower}. By the invariance of nearest-neighbor percolation under permutations and sign changes of the coordinate axes, we may assume $x\cdot e_1=n$. Throughout the remainder of this section, fix $h:=\lfloor\alpha n/4\rfloor$ and write
		\begin{align*}
			Q(x) &:= \big\{z\in\Zd : z\cdot e_1\in[n-h,n],\ \vert z\cdot e_j-x\cdot e_j\vert<2h\ \text{for}\ 2\leq j\leq d\big\}, \\
			F &:= \{z\in Q(x) : z\cdot e_1=n-h\}, \\
			\fC &:= \{u\in Q(x) : x\overset{Q(x)}\longleftrightarrow u\}.
		\end{align*}
		Since $x$ is at distance at least $\alpha n$ from $\partial\partial\Lambda_n$, we have $Q(x)\subset\Lambda_n$; thus $F$ is the unique face of $Q(x)$ lying in the interior of $\Lambda_n$, and $\fC$ is the cluster of $x$ within $Q(x)$. For $r\in\N_+$ and $w\in F$, write $w^{+r}:=w-r e_1$.
		
		\begin{figure}[h!]
			\centering
			\includegraphics[width=0.95\textwidth]{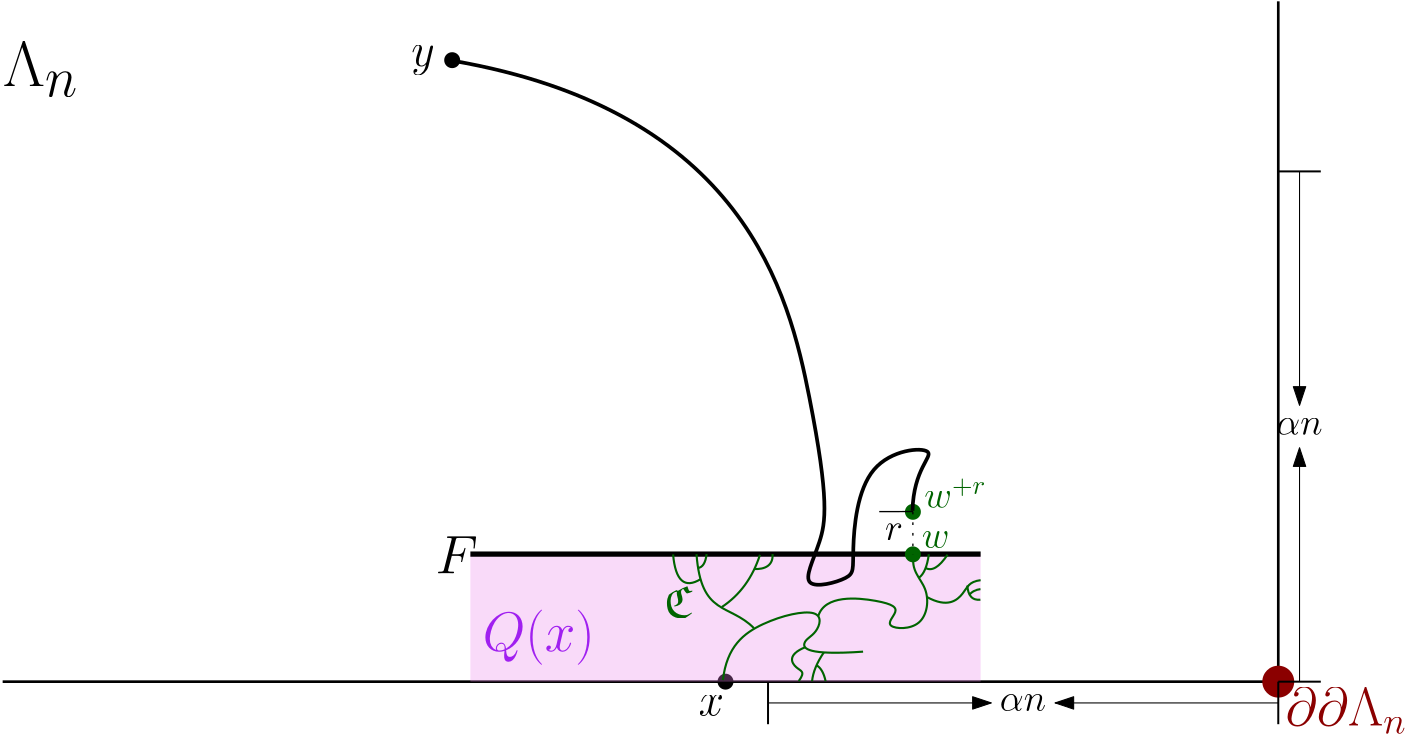}
			\caption{An illustration of the proof of Proposition~\ref{lem.box_twopoint_lower}.}\label{fig.box_twopoint_lower_setup}
		\end{figure}

		We first show how Proposition~\ref{lem.box_twopoint_lower} follows from two Lemma, whose proofs are postponed to the end of this section.
		
		For $r\in\N_+$, define
		\begin{align*}
			\EE_{\ref{lem.err_boundary}}(\fC,r) := \sum_{\substack{w,w'\in\fC\cap F \\ w\neq w'}} \max\{r,\vert w-w'\vert\}^{4-d}.
		\end{align*}
		
		\begin{lemma}\label{lem.err_boundary}
			There exists $C_{\ref{lem.err_boundary}}(\alpha)>0$ such that, for every $n\in\N_+$ and every $r\in\N_+$,
			\begin{align*}
				\E_{\p_c}\Ll[\EE_{\ref{lem.err_boundary}}(\fC,r)\Rr] \leq C_{\ref{lem.err_boundary}}\, r^{6-d}\, n^{-1}.
			\end{align*}
		\end{lemma}
		
		For $r\in\N_+$, define
		\begin{align*}
			\EE_{\ref{lem.err_cluster}}(\fC,r) := \sum_{w\in\fC\cap F}\ \sum_{v\in\fC} \max\{r,\vert w-v\vert\}^{2-d}.
		\end{align*}
		
		\begin{lemma}\label{lem.err_cluster}
			There exists $C_{\ref{lem.err_cluster}}(\alpha)>0$ such that, for every $n\in\N_+$ and every $r\in\N_+$,
			\begin{align*}
				\E_{\p_c}\Ll[\EE_{\ref{lem.err_cluster}}(\fC,r)\Rr] \leq C_{\ref{lem.err_cluster}}\, r^{6-d}\, n^{-1}.
			\end{align*}
		\end{lemma}

		\begin{proof}[Proof of Proposition~\ref{lem.box_twopoint_lower}, given Lemmas~\ref{lem.err_boundary} and~\ref{lem.err_cluster}]
			Let $r\in\N_+$ be a parameter, to be fixed at the end of the proof depending only on $\alpha,d$. We work throughout conditionally on $\fC$, viewed as a fixed admissible realization of the cluster of $x$ within $Q(x)$. Recall that, given $\fC$, the edges of $\Zd\setminus\EDGE(\fC)$ (writing $\EDGE(\fC)$ for the edges of $\Zd$ with at least one endpoint in $\fC$) remain i.i.d.\ $\mathrm{Bernoulli}(\p_c)$, independent of $\fC$.
			
			For $w\in\fC\cap F$, define
			\begin{equation*}
				A(w) := \Big\{w\overset{\Lambda_n\setminus\EDGE(\fC)}\longleftrightarrow y\Big\}\setminus\bigcup_{w'\in(\fC\cap F)\setminus\{w\}}\Big\{w\overset{\Lambda_n\setminus\EDGE(\fC)}\longleftrightarrow w'\Big\}.
			\end{equation*}
			Since $w\in\fC$, the event $A(w)$ implies $x\overset{\Lambda_n}\longleftrightarrow y$. Moreover, the events $A(w)$, $w\in\fC\cap F$, are pairwise disjoint, because $A(w)\cap A(w')$ for $w\neq w'$ would connect $w$ to $w'$ through $y$ off $\EDGE(\fC)$, contradicting the exclusion in $A(w)$. Hence
			\begin{equation}\label{eq.lower_step1}
				\P_{\p_c}\big[x\overset{\Lambda_n}\longleftrightarrow y\ \big|\ \fC\big] \geq \sum_{w\in\fC\cap F} \P_{\p_c}\big[A(w)\mid\fC\big].
			\end{equation}
			
			Fix $w\in\fC\cap F$ and set $w^{+r}:=w-r e_1$. Let $I:=\{w,w- e_1,\ldots,w-r e_1\}$ and let $U$ be the edges of $\Lambda_n\setminus\EDGE(\fC)$ with at least one endpoint in $I$, so $\vert U\vert\leq(r+1)\cdot2d$. Define
			\begin{equation*}
				A(w^{+r}) := \Big\{w^{+r}\overset{\Lambda_n\setminus\EDGE(\fC)}\longleftrightarrow y\Big\}\setminus\bigcup_{w'\in\fC\cap F}\Big\{w^{+r}\overset{\Lambda_n\setminus\EDGE(\fC)}\longleftrightarrow w'\Big\}.
			\end{equation*}
			Exactly as in the proof of Proposition~\ref{prop.extension} (with $I,U$ as above in place of those defined there), a witnessing configuration for $A(w^{+r})$ can be surgered, by opening the segment $w\to w^{+r}$ and the last edge of $U$ met by the witnessing path while closing the rest of $U$, into a configuration witnessing $A(w)$. this yields
			\begin{equation}\label{eq.lower_surgery}
				\P_{\p_c}\big[A(w)\mid\fC\big] \geq c_2(r)\, \P_{\p_c}\big[A(w^{+r})\mid\fC\big], \qquad c_2(r):=\min\{\p_c,1-\p_c\}^{(r+1)\cdot2d}.
			\end{equation}
			
			Decomposing, by the BK inequality, an open path from $w^{+r}$ to $y$ at its first intersection with $\fC$ (if it fails to avoid $\EDGE(\fC)$), and bounding the exclusion in $A(w^{+r})$ by a union bound,
			\begin{equation}\label{eq.lower_step3}
				\begin{aligned}
					\P_{\p_c}\big[A(w^{+r})\mid\fC\big] \geq \P_{\p_c}\big[w^{+r}\overset{\Lambda_n}\longleftrightarrow y\big]&- \sum_{v\in\fC}\P_{\p_c}\big[w^{+r}\overset{\Lambda_n}\longleftrightarrow v\big]\P_{\p_c}\big[v\overset{\Lambda_n}\longleftrightarrow y\big]\\
					& - \sum_{w'\in\fC\cap F}\P_{\p_c}\big[w^{+r}\leftrightarrow y,\ w^{+r}\leftrightarrow w'\big].
				\end{aligned}
			\end{equation}
			
			Since $v\in\fC\subset Q(x)$ and $w'\in\fC\cap F\subset Q(x)$ satisfy $v\cdot e_1,w'\cdot e_1\geq n-h$ while $w^{+r}\cdot e_1=n-h-r$, the triangle inequality gives $\vert w^{+r}-v\vert \geq \max\{r,\vert w-v\vert\}$ and $\vert w^{+r}-w'\vert \geq \max\{r,\vert w-w'\vert\}$. Since $x\in\partial\Lambda_n$ and $l_y\geq\alpha n$, $\vert x-y\vert\geq l_y\geq\alpha n$, and since $v,w,w'$ lie within a $l_\infty$ distance $\alpha n/2$ of $x$, we get $\vert w^{+r}-y\vert,\vert v-y\vert,\vert w'-y\vert \asymp_\alpha n$. Thus, uniformly over $v\in\fC$,
			\begin{equation}\label{esbvkur}
				\P_{\p_c}\big[w^{+r}\overset{\Lambda_n}\longleftrightarrow v\big]\P_{\p_c}\big[v\overset{\Lambda_n}\longleftrightarrow y\big] \lesssim_{\alpha} n^{2-d}\,\max\{r,\vert w-v\vert\}^{2-d}.
			\end{equation}
			Lemma~\ref{lem.threepoint_upper} applied to $w^{+r},y,w'$ gives, uniformly over $w'\in\fC\cap F$,
			\begin{equation}\label{esbvkur2}
				\P_{\p_c}\big[w^{+r}\leftrightarrow y,\ w^{+r}\leftrightarrow w'\big] \lesssim_{\alpha} n^{2-d}\,\max\{r,\vert w-w'\vert\}^{4-d}.
			\end{equation}
			Moreover, since $w^{+r}$ and $y$ both have a distance at least $\alpha n/4$ from the boundary $\partial \Lambda_n$, Lemma~\ref{lem.box_twopoint} gives
			\begin{equation}\label{esbvkur3}
				\P_{\p_c}[w^{+r}\overset{\Lambda_n}\longleftrightarrow y] \gtrsim_{\alpha} n^{2-d}
			\end{equation}

			Summing over $w \in \fC \cap F$ for \eqref{eq.lower_step3}, applying the inequalities \eqref{esbvkur}, \eqref{esbvkur2}, and \eqref{esbvkur3}, and recalling the definitions of $\EE_{\ref{lem.err_cluster}}(\fC,r)$ and $\EE_{\ref{lem.err_boundary}}(\fC,r)$, we obtain
			\begin{equation}
				\sum_{w \in \fC \cap F} \P_{\p_c}\big[A(w^{+r})\mid\fC\big] \geq c_4 |\fC \cap F| n^{2-d} - C_3 n^{2-d}\,\EE_{\ref{lem.err_cluster}}(\fC,r) - C_3 n^{2-d}\,\EE_{\ref{lem.err_boundary}}(\fC,r)
			\end{equation}
			for some $C_3 := C_3(\alpha) >0$ and $c_4 := c_4 (\alpha) >0$.
			Then applying \eqref{eq.lower_step1} and \eqref{eq.lower_surgery} gives
			\begin{equation}\label{eq.lower_conditional}
				\P_{\p_c}\big[x\overset{\Lambda_n}\longleftrightarrow y\ \big|\ \fC\big] \geq c_2 n^{2-d} \Big(c_4 |\fC \cap F|- C_3 \EE_{\ref{lem.err_cluster}}(\fC,r) - C_3 \EE_{\ref{lem.err_boundary}}(\fC,r)\Big).
			\end{equation}
			
			Taking expectation on both sides of \eqref{eq.lower_conditional} and applying Lemmas~\ref{lem.err_cluster} and~\ref{lem.err_boundary} gives
			\begin{equation}
				P_{\p_c}\big[x\overset{\Lambda_n}\longleftrightarrow y\big] \geq c_2 n^{2-d} \Big(c_4 \E_{\p_c} \big[ |\fC \cap F| \big]- C_3 r^{6-d} n^{-1} (C_{\ref{lem.err_cluster}} +C_{\ref{lem.err_boundary}})\Big).
			\end{equation}
			Lemma~\ref{lem.boundary_regularity} (applied with $2h,h$ in place of $rn,n$) gives
			\begin{equation*}
				\E_{\p_c}[\vert\fC\cap F\vert]=\sum_{w\in F}\P_{\p_c}[x\overset{Q(x)}\longleftrightarrow w]\geq c_{\ref{lem.boundary_regularity}}(2)/h\gtrsim_\alpha n^{-1}
			\end{equation*}
			Hence
			\begin{equation*}
				\P_{\p_c}\big[x\overset{\Lambda_n}\longleftrightarrow y\big] \geq c_2(r)\Big(c_5(\alpha) - C_6(\alpha)\,r^{6-d}\Big)\,n^{1-d}
			\end{equation*}
			for some $c_5 (\alpha), C_6(\alpha)>0$.
			
			Fix $r=r(\alpha)\in\N_+$ large enough that $C_6 (\alpha)\,r(\alpha)^{6-d}\leq c_5(\alpha)/2$. Then
			\begin{equation*}
				\P_{\p_c}\big[x\overset{\Lambda_n}\longleftrightarrow y\big] \geq \frac{c_5(\alpha)}{2}\,c_2(r(\alpha))\cdot n^{1-d} =: \frac{c_{\ref{lem.box_twopoint_lower}}(\alpha)}{n^{d-1}},
			\end{equation*}
			which proves the proposition.
		\end{proof}
		
		It remains to prove Lemmas~\ref{lem.err_boundary} and~\ref{lem.err_cluster}, which we do in the remainder of this subsection.

		\begin{proof}[Proof of Lemma~\ref{lem.err_boundary}]
			Throughout, $l_u$ denotes the distance from $u$ to $\partial Q(x)$. After translating $x$ to the origin and reflecting the height direction, $Q(x)$ is of the form $R_h$ (with $(\alpha_i),(\beta_i)$ depending only on the fixed ratio $1{:}4$ between the height and the transverse half-width of $Q(x)$), so Lemmas~\ref{lem.box_twopoint_onesided} and~\ref{lem.box_twopoint_twosided} apply to $Q(x)$ in place of $R_n$, with constants depending only on $d$. Taking expectation and writing $k:=\vert w-w'\vert$,
			\begin{equation*}
				\E_{\p_c}\big[\EE_{\ref{lem.err_boundary}}(\fC,r)\big] = \sum_{w\neq w'\in F} \P_{\p_c}\big[x\overset{Q(x)}\leftrightarrow w,\ x\overset{Q(x)}\leftrightarrow w'\big]\, \max\{r,k\}^{4-d}.
			\end{equation*}
			Since $l_x=l_w=l_{w'}=0$, decomposing at the branch point $u\in Q(x)$ of the two connections (BK) and applying Lemma~\ref{lem.box_twopoint_twosided} to each of the three resulting factors,
			\begin{equation*}
				\begin{aligned}
					\P_{\p_c}\big[x\leftrightarrow w,\ x\leftrightarrow w'\big] &\leq \sum_{u\in Q(x)} \P_{\p_c}[x\leftrightarrow u]\P_{\p_c}[u\leftrightarrow w]\P_{\p_c}[u\leftrightarrow w']\\
					&\lesssim  \sum_u (1+l_u)^3\,\vert x-u\vert^{-d}\vert u-w\vert^{-d}\vert u-w'\vert^{-d}.
				\end{aligned}
			\end{equation*}
			Split the sum by which of $x,w,w'$ is nearest to $u$:
			
			\emph{$u$ nearest $w$} ($\vert x-u\vert\sim h$, $\vert u-w'\vert\sim k$, $l_u\lesssim1+\vert u-w\vert$):
			\begin{equation*}
				\sum_{\vert u-w\vert\lesssim k} (1+\vert u-w\vert)^3\, h^{-d}\,\vert u-w\vert^{-d}\,k^{-d} \lesssim  h^{-d}k^{-d}\sum_{m=1}^k m^{d-1+3-d} \sim h^{-d}k^{3-d}.
			\end{equation*}
			
			\emph{$u$ nearest $w'$}: symmetric, same bound $h^{-d}k^{3-d}$.
			
			\emph{$u$ nearest $x$} ($\vert u-w\vert\sim\vert u-w'\vert\sim h$, $l_u\lesssim1+\vert x-u\vert$):
			\begin{equation*}
				\sum_{\vert x-u\vert\lesssim h}(1+\vert x-u\vert)^3\,\vert x-u\vert^{-d}\,h^{-2d} \lesssim  h^{-2d}\sum_{m=1}^h m^{d-1+3-d} \sim h^{3-2d} \leq h^{-d}k^{3-d} \quad(k\leq h).
			\end{equation*}
			Hence $\P_{\p_c}[x\leftrightarrow w,x\leftrightarrow w']\lesssim  h^{-d}k^{3-d}$. Since $\vert F\vert\sim h^{d-1}$ and $F$ has $\sim k^{d-2}$ points at distance $k$,
			\begin{equation*}
				\E_{\p_c}\big[\EE_{\ref{lem.err_boundary}}(\fC,r)\big] \lesssim  h^{d-1}\cdot h^{-d}\sum_{k=1}^{h} k^{d-2}\,k^{3-d}\max\{r,k\}^{4-d} = h^{-1}\sum_k k\max\{r,k\}^{4-d}.
			\end{equation*}
			Splitting at $k\lessgtr r$: $\sum_{k\leq r}k\cdot r^{4-d}\sim r^{6-d}$, and $\sum_{k>r}k^{5-d}\lesssim  r^{6-d}$ (converges since $d>6$). Thus $\E_{\p_c}[\EE_{\ref{lem.err_boundary}}(\fC,r)]\lesssim_{\alpha} h^{-1}r^{6-d}\sim n^{-1}r^{6-d}$.
		\end{proof}
		
		\begin{proof}[Proof of Lemma~\ref{lem.err_cluster}]
			As above, $l_u$ denotes distance to $\partial Q(x)$, and Lemmas~\ref{lem.box_twopoint_onesided} and~\ref{lem.box_twopoint_twosided} apply to $Q(x)$ (in the form $R_h$ described in the proof of Lemma~\ref{lem.err_boundary}), with constants depending only on $d$. Taking expectation,
			\begin{equation*}
				\E_{\p_c}\big[\EE_{\ref{lem.err_cluster}}(\fC,r)\big] = \sum_{w\in F}\sum_{v\in Q(x)} \P_{\p_c}\big[x\overset{Q(x)}\leftrightarrow w,\ x\overset{Q(x)}\leftrightarrow v\big]\,\max\{r,\vert w-v\vert\}^{2-d}.
			\end{equation*}
			Decomposing at the branch point $u$ (BK) as before,
			\begin{equation*}
				\P_{\p_c}[x\leftrightarrow w,x\leftrightarrow v] \leq \sum_u \P_{\p_c}[x\leftrightarrow u]\P_{\p_c}[u\leftrightarrow w]\P_{\p_c}[u\leftrightarrow v].
			\end{equation*}
			Fix $u$ and sum over $v\in Q(x)$ first, using the one-sided bound $\P_{\p_c}[u\leftrightarrow v]\lesssim (1+l_u)\vert u-v\vert^{1-d}$ (Lemma~\ref{lem.box_twopoint_onesided}) and the convolution estimate $\sum_v\vert u-v\vert^{1-d}\max\{r,\vert w-v\vert\}^{2-d}\lesssim \max\{r,\vert u-w\vert\}^{3-d}$:
			\begin{equation*}
				\sum_{v\in Q(x)}\P_{\p_c}[u\leftrightarrow v]\max\{r,\vert w-v\vert\}^{2-d} \lesssim  (1+l_u)\max\{r,\vert u-w\vert\}^{3-d}.
			\end{equation*}
			Using Lemma~\ref{lem.box_twopoint_twosided} on the remaining two factors ($l_x=l_w=0$), and splitting by which of $x,w$ is nearest $u$:
			
			\emph{$u$ nearest $w$} ($\vert x-u\vert\sim h$, write $m:=\vert u-w\vert$, then $l_u \lesssim m$):
			\begin{equation*}
				h^{-d}\sum_{m=1}^h m^{d-1}(1+m)^3\,m^{-d}\max\{r,m\}^{3-d} = h^{-d}\sum_m m^2\max\{r,m\}^{3-d} \lesssim  h^{-d}r^{6-d}
			\end{equation*}
			(same split at $m\lessgtr r$ as in Lemma~\ref{lem.err_boundary}, using $d>6$).
			
			\emph{$u$ nearest $x$} ($\vert u-w\vert\sim h$, write $m:=\vert x-u\vert$; since $r\leq h$, $\max\{r,h\}\sim h$):
			\begin{equation*}
				h^{-d}\cdot h^{3-d}\sum_{m=1}^h m^{d-1}(1+m)^3\,m^{-d} \sim h^{3-2d}\cdot h^3 = h^{6-2d} \lesssim h^{-d}r^{6-d} \quad(\text{since } d>6,\ r\leq h).
			\end{equation*}
			Hence $\sum_v\P_{\p_c}[x\leftrightarrow w,x\leftrightarrow v]\max\{r,\vert w-v\vert\}^{2-d}\lesssim  h^{-d}r^{6-d}$, uniformly in $w\in F$. Summing over $w\in F$ ($\vert F\vert\sim h^{d-1}$),
			\begin{equation*}
				\E_{\p_c}\big[\EE_{\ref{lem.err_cluster}}(\fC,r)\big] \lesssim  h^{d-1}\cdot h^{-d}r^{6-d} = h^{-1}r^{6-d} \sim n^{-1}r^{6-d},
			\end{equation*}
			which proves the lemma.
		\end{proof}

		\section{Estimates of variance in Section~\ref{subsec.regular.var}}\label{sec.appendix_variance}
		
		In this appendix, we complete the proofs of Propositions~\ref{prop:VARf1} and~\ref{prop:VARf2}. In the proofs of these propositions, the conditional variances of $\ERE_1(x,\fC,r)$ and $\ERE_2(x,\fC,r)$ are bounded, up to a constant depending on $d$ and $c_1$, by the sum of four terms (\eqref{sunzi1}--\eqref{sunzi4}, \eqref{sunzi5}--\eqref{sunzi8}). The main text only bounds the part of the third term with $u\in\Lambda_n\setminus\Lambda_{n/2}$. Here we bound all the remaining contributions: this is the content of Propositions~\ref{prop.var_term_1}--\ref{prop.var_term_4} below, which, together with the estimates in the main text, prove Propositions~\ref{prop:VARf1} and~\ref{prop:VARf2}.
		
		Throughout this appendix, $d>6$, $c_1\in(0,1/5)$, $r\in\N_+$, and we fix $x\in\partial_{1}^\circ \Lambda_n$. The integer $n$ is assumed to be sufficiently large depending on $d$, $c_1$, and $r$. We write $\lesssim$ (resp.\ $\lesssim_{c_1}$) for an inequality up to a multiplicative constant depending only on $d$ (resp.\ on $d$ and $c_1$), and never on $n$, $r$ or $x$.
		
\begin{figure}[h!]
			\centering
			\includegraphics[width=0.75\textwidth]{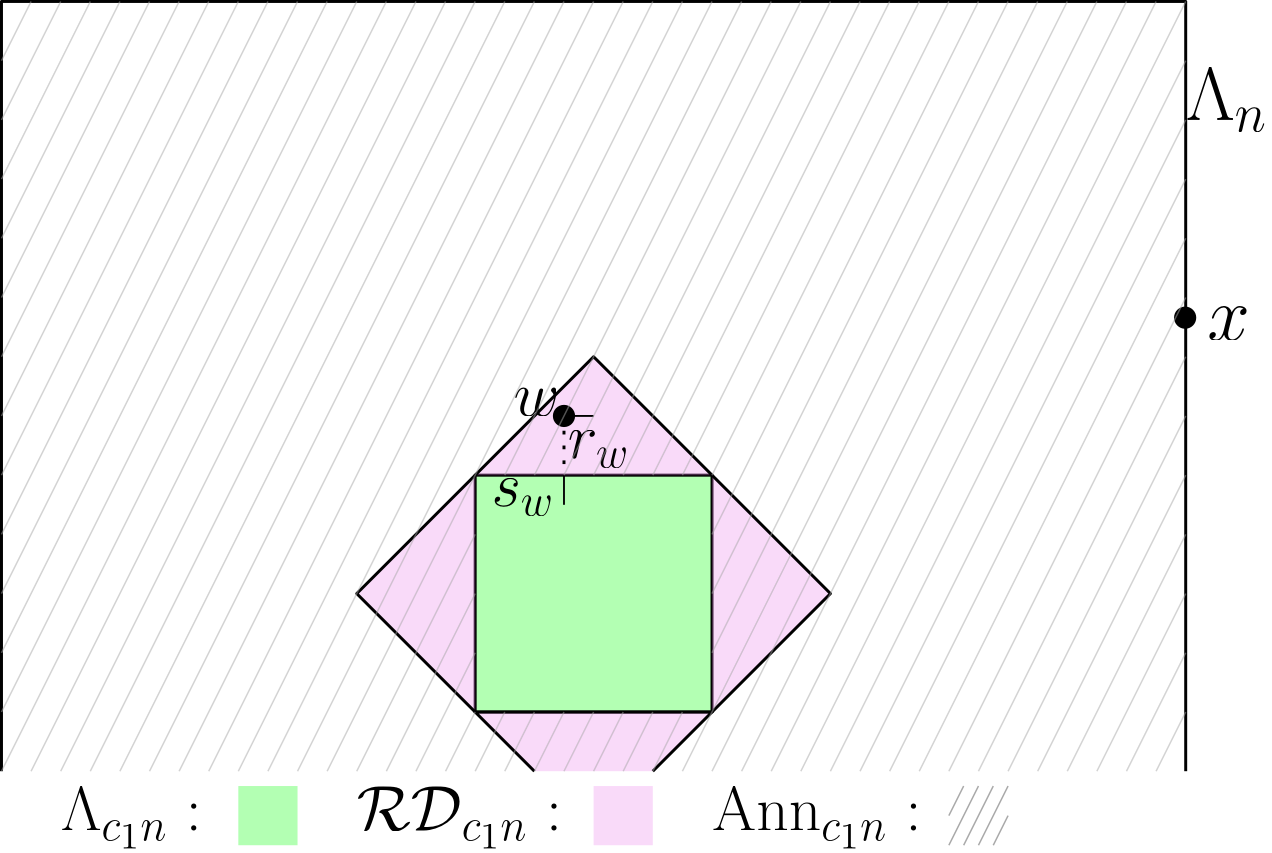}
			\caption{The geometric setting of this appendix: the cubes $\Lambda_n$ and $\Lambda_{c_1n}$, the set $\RD_{c_1n}$, and the fixed point $x\in\partial_{1}^\circ \Lambda_n$.}\label{fig.variance_geometry}
		\end{figure}

		Recall that $\vert\cdot\vert$ is the $\ell^\infty$ norm, that $\Lambda_n=[-n,n]^d\cap\Zd$, and that $\partial\Lambda_n$ is the inner vertex boundary of $\Lambda_n$. For $u\in\Lambda_n$, $l_u:=n-\vert u\vert$ is the distance from $u$ to $\partial\Lambda_n$. Recall also that $\operatorname{int}\Lambda_{c_1n}=\Lambda_{c_1n}\setminus\partial\Lambda_{c_1n}$ and $\mathrm{Ann}_{c_1n}=\Lambda_n\setminus\operatorname{int}\Lambda_{c_1n}$. The fixed point $x\in\partial_{1}^\circ \Lambda_n$ satisfies $\vert x\vert=n$, i.e.\ $l_x=0$. The set $\RD_{c_1n}\subset\mathrm{Ann}_{c_1n}$ and the distances $r_w,s_w$ (for $w\in\mathrm{Ann}_{c_1n}$) are those defined in Appendix~\ref{subsct:Upper bound for the 2-point function in cubes and annulus}. See Figure~\ref{fig.variance_geometry}.
		
		We also recall that, for $D\subset\Zd$ and $u,v\in D$, $\tau_D(u,v):=\P_{\p_c}\big[u\stackrel{D}{\longleftrightarrow}v\big]$, and that, for $u\in\mathrm{Ann}_{c_1n}$, $g(u):=\P_{\p_c}\big[u\stackrel{\mathrm{Ann}_{c_1n}}{\longleftrightarrow}\partial\Lambda_{c_1n}\big]$. Finally, we recall
		\begin{align*}
			\rho_1 &:= \sum_{y\in\Lambda_n\setminus\Lambda_{3n/4}}\max\{r,\vert x-y\vert\}^{-d+2}\,g(y),\\
			\rho_2 &:= \sum_{y\in\partial\Lambda_n}\max\{r,\vert x-y\vert\}^{-d+4}\,g(y),
		\end{align*}
		and, for $u\in\mathrm{Ann}_{c_1n}$,
		\begin{align*}
			F_1(u) &:= \sum_{y\in\Lambda_n\setminus\Lambda_{3n/4}}\max\{r,\vert x-y\vert\}^{-d+2}\,\tau_{\mathrm{Ann}_{c_1n}}(u,y),\\
			F_2(u) &:= \sum_{y\in\partial\Lambda_n}\max\{r,\vert x-y\vert\}^{-d+4}\,\tau_{\mathrm{Ann}_{c_1n}}(u,y).
		\end{align*}
		These are the quantities $\rho_i$ and $F_i$ appearing in the proofs of Propositions~\ref{prop:VARf1} and~\ref{prop:VARf2}.

		The following four propositions bound the four terms mentioned above.
		
		\begin{prop}\label{prop.var_term_1}
			For $i\in\{1,2\}$,
			\begin{align*}
				n^3\sum_{u,w\in\mathrm{Ann}_{c_1n}}g(w)\,\tau_{\mathrm{Ann}_{c_1n}}(w,u)\,\tau_{\mathrm{Ann}_{c_1n}}(u,x)\,F_i(w)\,F_i(u)\lesssim_{c_1}r^{6-d}.
			\end{align*}
		\end{prop}
		
		\begin{prop}\label{prop.var_term_2}
			For $i\in\{1,2\}$,
			\begin{align*}
				n^3\sum_{u,w\in\mathrm{Ann}_{c_1n}}g(u)\,\tau_{\mathrm{Ann}_{c_1n}}(u,w)\,\tau_{\mathrm{Ann}_{c_1n}}(u,x)\,F_i(w)^2\lesssim_{c_1}r^{6-d}.
			\end{align*}
		\end{prop}
		
		\begin{prop}\label{prop.var_term_3}
			For $i\in\{1,2\}$,
			\begin{align*}
				n^3\rho_i\sum_{u\in\mathrm{Ann}_{c_1n}\cap\Lambda_{n/2}}g(u)\,\tau_{\mathrm{Ann}_{c_1n}}(u,x)\,F_i(u)\lesssim_{c_1}r^{6-d}.
			\end{align*}
		\end{prop}
		
		\begin{prop}\label{prop.var_term_4}
			For $i\in\{1,2\}$,
			\begin{align*}
				\sum_{w\in\mathrm{Ann}_{c_1n}}g(w)\,F_i(w)^2\lesssim_{c_1}r^{6-d}.
			\end{align*}
		\end{prop}

		We first record the bounds on $\rho_i$ and $F_i$ that were proved in the main text.
		
		\begin{lemma}\label{lem.rho_F_basic}
			For $i\in\{1,2\}$,
			\begin{align*}
				\rho_i\lesssim_{c_1}1, \qquad F_i(u)\lesssim\max\{r,\vert x-u\vert\}^{4-d}\quad\text{for every }u\in\mathrm{Ann}_{c_1n}.
			\end{align*}
		\end{lemma}
		
		The bounds on $\rho_1$ and $F_1$ (resp.\ $\rho_2$ and $F_2$) are proved in the proof of Proposition~\ref{prop:VARf1} (resp.\ Proposition~\ref{prop:VARf2}); see in particular \eqref{eq.F1_bound} and \eqref{eq.F2_bound}.
		
		The following is a refined bound on $F_i$ in $\RD_{c_1n}$:
		
		\begin{lemma}\label{lem.F_RD}
			For $i\in\{1,2\}$ and every $w\in\RD_{c_1n}$,
			\begin{align*}
				F_i(w)\lesssim n^{4-d}\,(1+r_w)\,s_w^{-1}.
			\end{align*}
		\end{lemma}
		
		\begin{proof}
			Fix $w\in\RD_{c_1n}$. Every $y$ in the sums defining $F_1(w)$ and $F_2(w)$ lies in $\Lambda_n\setminus\Lambda_{3n/4}$, so Lemma~\ref{lem.RD_bound} gives
			\begin{align*}
				\tau_{\mathrm{Ann}_{c_1n}}(w,y)\lesssim n^{1-d}\,(1+l_y)\,(1+r_w)\,s_w^{-1}.
			\end{align*}
			Since $l_y=0$ for $y\in\partial\Lambda_n$, we obtain
			\begin{align*}
				F_1(w)&\lesssim n^{1-d}\,(1+r_w)\,s_w^{-1}\sum_{y\in\Lambda_n\setminus\Lambda_{3n/4}}\max\{r,\vert x-y\vert\}^{2-d}\,(1+l_y),\\
				F_2(w)&\lesssim n^{1-d}\,(1+r_w)\,s_w^{-1}\sum_{y\in\partial\Lambda_n}\max\{r,\vert x-y\vert\}^{4-d}.
			\end{align*}
			It remains to show that both sums over $y$ are $\lesssim n^3$. For the first one, since $\vert x\vert=n$ we have $l_y\leq\vert x-y\vert$. Using $\#\{y\in\Lambda_n:\vert x-y\vert=m\}\lesssim m^{d-1}$ and splitting at $m\lessgtr r$,
			\begin{align*}
				\sum_{y\in\Lambda_n\setminus\Lambda_{3n/4}}\max\{r,\vert x-y\vert\}^{2-d}\,(1+l_y)\lesssim \sum_{m=1}^{2n}m\cdot m^{d-1}\max\{r,m\}^{2-d}\lesssim r^3+n^3\lesssim n^3,
			\end{align*}
			For the second one, using $\#\{y\in\partial\Lambda_n:\vert x-y\vert=m\}\lesssim m^{d-2}$ and splitting at $m\lessgtr r$,
			\begin{align*}
				\sum_{y\in\partial\Lambda_n}\max\{r,\vert x-y\vert\}^{4-d}\lesssim \sum_{m=1}^{2n}m^{d-2}\max\{r,m\}^{4-d}\lesssim r^3+n^3\lesssim n^3.
			\end{align*}
			This proves the lemma.
		\end{proof}

		The proofs of the Propositions~\ref{prop.var_term_1}, \ref{prop.var_term_2}, \ref{prop.var_term_3}, and~\ref{prop.var_term_4} will use the following bounds. 
		
		For $g$, we will use the following two upper bounds. Since $g(u)=\P_{\p_c}[u\in\fC]$ for $u\in\mathrm{Ann}_{c_1n}$, Proposition~\ref{PxinfC} gives
		\begin{equation}\label{eq.g_outer}
			g(u)\lesssim_{c_1} (1+l_u)\,n^{-3},\qquad\forall\,u\in\Lambda_n\setminus\Lambda_{n/2}.
		\end{equation}
		Moreover, the $\ell^\infty$ distance from $u\in\mathrm{Ann}_{c_1n}$ to $\partial\Lambda_{c_1n}$ is $r_u$, so on the event $\{u\stackrel{\mathrm{Ann}_{c_1n}}{\longleftrightarrow}\partial\Lambda_{c_1n}\}$, $u$ is connected to $\Zd\setminus(u+\Lambda_{r_u-1})$, and \eqref{eq.one_arm} gives
		\begin{equation}\label{eq.g_inner}
			g(u)\lesssim(1+r_u)^{-2},\qquad\forall\,u\in\mathrm{Ann}_{c_1n}.
		\end{equation}

		Using the bounds \eqref{eq.g_outer} and \eqref{eq.g_inner} on $g$, the bounds on $\tau_{\mathrm{Ann}_{c_1n}}$ from Lemmas~\ref{lem.box_twopoint_onesided}, \ref{lem.box_twopoint_twosided}, and~\ref{lem.RD_bound}, and the bounds on $\rho_i$ and $F_i$ in Lemmas~\ref{lem.rho_F_basic} and~\ref{lem.F_RD}, we now prove Propositions~\ref{prop.var_term_1}, \ref{prop.var_term_2}, \ref{prop.var_term_3} and~\ref{prop.var_term_4}.
		
		In the remaining part of this appendix, we will ignore the dependency on $c_1$ and write $\lesssim$ for an inequality up to a multiplicative constant depending only on $d$ and $c_1$.

		\begin{proof}[Proof of Proposition~\ref{prop.var_term_1}]
			Throughout this proof, fix $i\in\{1,2\}$, write $F:=F_i$, and for $v\in\Lambda_n$ let $M(v):=\max\{r,\vert x-v\vert\}$. By Lemma~\ref{lem.rho_F_basic},
			\begin{equation}\label{eq.F_pointwise}
				F(v)\lesssim M(v)^{4-d},\qquad\forall\,v\in\mathrm{Ann}_{c_1n}.
			\end{equation}
			This is the only property of $F$ that is used below, so the proof is the same for $i=1$ and $i=2$. We split the sum in the proposition according to whether $w\in\Lambda_n\setminus\Lambda_{n/2}$ or $w\in\mathrm{Ann}_{c_1n}\cap\Lambda_{n/2}$ , and show that each of the two contributions is $\lesssim_{c_1}r^{6-d}$.
			
			\medskip
			\noindent\emph{Case 1: $w\in\Lambda_n\setminus\Lambda_{n/2}$.} Since $\vert x\vert=n$, we have $l_w\leq\vert x-w\vert$, so \eqref{eq.g_outer} gives $n^3g(w)\lesssim1+\vert x-w\vert$. Together with \eqref{eq.F_pointwise}, the contribution of this case is
			\begin{equation*}
				\lesssim\sum_{u\in\mathrm{Ann}_{c_1n}}\tau_{\mathrm{Ann}_{c_1n}}(u,x)\,F(u)\,\Phi(u),\qquad\Phi(u):=\sum_{w\in\Lambda_n}(1+\vert x-w\vert)\,\tau_{\Lambda_n}(w,u)\,M(w)^{4-d}.
			\end{equation*}
			Fix $u\in\mathrm{Ann}_{c_1n}$, and write $k:=\vert x-u\vert$ and $M_k:=\max\{r,k\}$. We claim that
			\begin{equation}\label{eq.Phi_bound}
				\Phi(u)\lesssim(1+l_u)\,M_k^{6-d}.
			\end{equation}
			By Lemma~\ref{lem.box_twopoint_onesided}, $\tau_{\Lambda_n}(w,u)\lesssim(1+l_u)(1+\vert w-u\vert)^{1-d}$ for every $w\in\Lambda_n$. Write $m:=\vert x-w\vert$ and $p:=\vert w-u\vert$. We split the sum defining $\Phi(u)$ into three parts.
			\begin{itemize}
				\item $m<k/2$. Then $p\geq k-m>k/2$, so $\tau_{\Lambda_n}(w,u)\lesssim(1+l_u)(1+k)^{1-d}$. Using $\#\{w\in\Lambda_n:\vert x-w\vert=m\}\lesssim(1+m)^{d-1}$, this part contributes
				\begin{equation*}
					\lesssim(1+l_u)(1+k)^{1-d}\sum_{m<k/2}(1+m)^{d}\max\{r,m\}^{4-d}.
				\end{equation*}
				If $k\leq2r$, the sum is at most $r^{4-d}\sum_{m<k/2}(1+m)^d\lesssim r^{4-d}(1+k)^{d+1}$, so the contribution is $\lesssim(1+l_u)\,r^{4-d}(1+k)^2\lesssim(1+l_u)\,r^{6-d}$. If $k>2r$, splitting at $m\lessgtr r$, the sum is $\lesssim r^{4-d}\sum_{m\leq r}(1+m)^d+\sum_{r<m<k/2}m^4\lesssim r^5+k^5\lesssim k^5$, so the contribution is $\lesssim(1+l_u)\,k^{6-d}$. In both cases it is $\lesssim(1+l_u)M_k^{6-d}$.
				\item $k/2\leq m<2k$. Here $1+m\lesssim1+k$ and $M(w)\geq\max\{r,k/2\}\gtrsim M_k$, so $M(w)^{4-d}\lesssim M_k^{4-d}$. Moreover $p\leq m+k<3k$, and $\#\{w\in\Lambda_n:\vert w-u\vert=p\}\lesssim(1+p)^{d-1}$ gives $\sum_{w:\,p<3k}(1+p)^{1-d}\lesssim\sum_{p<3k}1\lesssim1+k$. Hence this part contributes $\lesssim(1+l_u)(1+k)^2M_k^{4-d}\lesssim(1+l_u)M_k^{6-d}$, where we used $1+k\lesssim M_k$.
				\item $m\geq2k$. Then $p\geq m-k\geq m/2$, so $\tau_{\Lambda_n}(w,u)\lesssim(1+l_u)(1+m)^{1-d}$. Using $\#\{w\in\Lambda_n:\vert x-w\vert=m\}\lesssim(1+m)^{d-1}$, this part contributes
				\begin{equation*}
					\lesssim(1+l_u)\sum_{m=2k}^{2n}(1+m)\max\{r,m\}^{4-d}.
				\end{equation*}
				If $2k\leq r$, splitting at $m\lessgtr r$ and using $d>6$, the sum is $\lesssim r^{4-d}\sum_{m\leq r}(1+m)+\sum_{m>r}m^{5-d}\lesssim r^{6-d}$. If $2k>r$, the sum is $\lesssim\sum_{m\geq2k}m^{5-d}\lesssim k^{6-d}$. In both cases it is $\lesssim M_k^{6-d}$.
			\end{itemize}
			Summing the three parts proves \eqref{eq.Phi_bound}.
			
			Now, since $l_x=0$, Lemma~\ref{lem.box_twopoint_twosided} gives $\tau_{\mathrm{Ann}_{c_1n}}(u,x)\leq\tau_{\Lambda_n}(x,u)\lesssim(1+l_u)(1+k)^{-d}$. Using this, \eqref{eq.F_pointwise}, \eqref{eq.Phi_bound} and $l_u\leq k$, the contribution of Case~1 is
			\begin{equation*}
				\lesssim\sum_{u\in\mathrm{Ann}_{c_1n}}(1+k)^{2-d}\,M_k^{10-2d}\lesssim\sum_{k=0}^{2n}(1+k)\max\{r,k\}^{10-2d}\lesssim r^{12-2d}\leq r^{6-d},
			\end{equation*}
			where we used $\#\{u\in\Lambda_n:\vert x-u\vert=k\}\lesssim(1+k)^{d-1}$ and used $d>6$.
			
			\medskip
			\noindent\emph{Case 2: $w\in\mathrm{Ann}_{c_1n}\cap\Lambda_{n/2}$.} For such $w$ we have $\vert x-w\vert\geq n/2$, so $F(w)\lesssim n^{4-d}$ by \eqref{eq.F_pointwise}, and $g(w)\lesssim(1+r_w)^{-2}$ by \eqref{eq.g_inner}. Hence the contribution of this case is
			\begin{equation*}
				\lesssim n^{7-d}\sum_{w\in\mathrm{Ann}_{c_1n}\cap\Lambda_{n/2}}(1+r_w)^{-2}\,\Psi(w),\qquad\Psi(w):=\sum_{u\in\mathrm{Ann}_{c_1n}}\tau_{\mathrm{Ann}_{c_1n}}(w,u)\,\tau_{\mathrm{Ann}_{c_1n}}(u,x)\,F(u).
			\end{equation*}
			We first bound the sum over $w$ of $(1+r_w)^{-2}$. For $w\in\mathrm{Ann}_{c_1n}$ with $\vert w\vert=\rho$ we have $r_w=\rho-c_1n$, and $\#\{w:\vert w\vert=\rho\}\lesssim\rho^{d-1}$. Hence
			\begin{equation}\label{eq.inner_sum}
				\sum_{w\in\mathrm{Ann}_{c_1n}\cap\Lambda_{n/2}}(1+r_w)^{-2}\lesssim\sum_{\rho=c_1n}^{n/2}\rho^{d-1}(1+\rho-c_1n)^{-2}\lesssim n^{d-1}.
			\end{equation}

			Next, we bound $\Psi(w)$ uniformly in $w\in\mathrm{Ann}_{c_1n}\cap\Lambda_{n/2}$, splitting the sum over $u$ according to $\vert x-u\vert<n/4$ or $\vert x-u\vert\geq n/4$. Write $k:=\vert x-u\vert$ and $M_k:=\max\{r,k\}$.
			\begin{itemize}
				\item $k<n/4$. Then $\vert w-u\vert\geq\vert x-w\vert-k>n/4$, so Lemma~\ref{lem.box_twopoint_onesided} gives $\tau_{\mathrm{Ann}_{c_1n}}(w,u)\leq\tau_{\Lambda_n}(u,w)\lesssim(1+l_u)n^{1-d}$. As in Case~1, $\tau_{\mathrm{Ann}_{c_1n}}(u,x)\lesssim(1+l_u)(1+k)^{-d}$, and $F(u)\lesssim M_k^{4-d}$. Using $l_u\leq k$ and $\#\{u\in\Lambda_n:\vert x-u\vert=k\}\lesssim(1+k)^{d-1}$, this part contributes
				\begin{equation*}
					\lesssim n^{1-d}\sum_{k}(1+k)^{d-1}(1+k)^{2-d}M_k^{4-d}=n^{1-d}\sum_k(1+k)\max\{r,k\}^{4-d}\lesssim n^{1-d}\,r^{6-d},
				\end{equation*}
				where, for the last inequality, we split at $k\lessgtr r$ and used $d>6$.
				\item $k\geq n/4$. Here $\tau_{\mathrm{Ann}_{c_1n}}(u,x)\leq\tau_{\Lambda_n}(x,u)\lesssim k^{1-d}\lesssim n^{1-d}$ by Lemma~\ref{lem.box_twopoint_onesided}, and $F(u)\lesssim n^{4-d}$. Moreover, $\sum_{u\in\Lambda_n}\tau_{\mathrm{Ann}_{c_1n}}(w,u)\leq\sum_{u\in\Lambda_n}\tau_{\Zd}(w,u)\lesssim\sum_{p\leq2n}(1+p)^{d-1}(1+p)^{2-d}\lesssim n^2$. Hence this part contributes $\lesssim n^{7-2d}$.
			\end{itemize}
			Therefore $\Psi(w)\lesssim n^{1-d}r^{6-d}+n^{7-2d}$. Combining this with \eqref{eq.inner_sum}, the contribution of Case~2 is
			\begin{equation*}
				\lesssim n^{7-d}\cdot n^{d-1}\big(n^{1-d}r^{6-d}+n^{7-2d}\big)=n^{7-d}\,r^{6-d}+n^{13-2d} \lesssim r^{6-d},
			\end{equation*}
			where, in the last inequality, we used  $d\geq7$ and that $n$ is sufficiently large.
			
			Adding the contributions of the two cases proves the proposition.
		\end{proof}

		\begin{proof}[Proof of Proposition~\ref{prop.var_term_2}]
			Fix $i\in\{1,2\}$, write $F:=F_i$, and let $M(v):=\max\{r,\vert x-v\vert\}$ for $v\in\Lambda_n$; recall \eqref{eq.F_pointwise}. We split the sum in the proposition according to whether $u\in\Lambda_n\setminus\Lambda_{n/2}$ or $u\in\mathrm{Ann}_{c_1n}\cap\Lambda_{n/2}$.
			
			\medskip
			\noindent\emph{Case 1: $u\in\Lambda_n\setminus\Lambda_{n/2}$.} Since $l_u\leq\vert x-u\vert$, \eqref{eq.g_outer} gives $n^3g(u)\lesssim1+\vert x-u\vert$. By Lemma~\ref{lem.box_twopoint_onesided}, $\tau_{\mathrm{Ann}_{c_1n}}(u,x)\lesssim(1+\vert x-u\vert)^{1-d}$, and $\tau_{\mathrm{Ann}_{c_1n}}(u,w)\leq\tau_{\Zd}(u,w)\lesssim(1+\vert u-w\vert)^{2-d}$. Hence, by the standard bubble estimate, for every $w\in\mathrm{Ann}_{c_1n}$,
			\begin{equation*}
				n^3\sum_{u\in\Lambda_n\setminus\Lambda_{n/2}}g(u)\,\tau_{\mathrm{Ann}_{c_1n}}(u,x)\,\tau_{\mathrm{Ann}_{c_1n}}(u,w)\lesssim\sum_{u\in\Zd}(1+\vert x-u\vert)^{2-d}(1+\vert u-w\vert)^{2-d}\lesssim(1+\vert x-w\vert)^{4-d}.
			\end{equation*}
			Therefore, by \eqref{eq.F_pointwise} and $\#\{w\in\Lambda_n:\vert x-w\vert=m\}\lesssim(1+m)^{d-1}$, the contribution of Case~1 is
			\begin{equation*}
				\lesssim\sum_{w\in\mathrm{Ann}_{c_1n}}(1+\vert x-w\vert)^{4-d}M(w)^{8-2d}\lesssim\sum_{m=0}^{2n}(1+m)^3\max\{r,m\}^{8-2d}\lesssim r^{12-2d}\leq r^{6-d},
			\end{equation*}
			where we split at $m\lessgtr r$ and used $d>6$.
			
			\medskip
			\noindent\emph{Case 2: $u\in\mathrm{Ann}_{c_1n}\cap\Lambda_{n/2}$.} For such $u$ we have $\vert x-u\vert\geq n/2$, so $\tau_{\mathrm{Ann}_{c_1n}}(u,x)\lesssim n^{1-d}$ by Lemma~\ref{lem.box_twopoint_onesided}, and $g(u)\lesssim(1+r_u)^{-2}$ by \eqref{eq.g_inner}. We claim that
			\begin{equation*}
				\Theta(u):=\sum_{w\in\mathrm{Ann}_{c_1n}}\tau_{\mathrm{Ann}_{c_1n}}(u,w)\,F(w)^2\lesssim n^{3-d}.
			\end{equation*}
			Write $m:=\vert x-w\vert$ and $p:=\vert u-w\vert$. If $p\leq n/4$, then $m\geq n/4$, so $F(w)^2\lesssim n^{8-2d}$, and, since $\tau_{\mathrm{Ann}_{c_1n}}(u,w)\lesssim(1+p)^{2-d}$, these $w$ contribute $\lesssim n^{8-2d}\sum_{p\leq n/4}(1+p)\lesssim n^{10-2d}$. If $p>n/4$, then $\tau_{\mathrm{Ann}_{c_1n}}(u,w)\lesssim n^{2-d}$, and since $(1+m)^{d-1}\max\{r,m\}^{8-2d}\lesssim\max\{r,m\}^{7-d}\leq1$, these $w$ contribute $\lesssim n^{2-d}\sum_{m\leq2n}(1+m)^{d-1}\max\{r,m\}^{8-2d}\lesssim n^{3-d}$. Since $d\geq7$, $n^{10-2d}\leq n^{3-d}$, which proves the claim. Hence, by \eqref{eq.inner_sum}, the contribution of Case~2 is
			\begin{equation*}
				\lesssim n^3\sum_{u\in\mathrm{Ann}_{c_1n}\cap\Lambda_{n/2}}(1+r_u)^{-2}\,n^{1-d}\,n^{3-d}\lesssim n^{7-2d}\cdot n^{d-1}=n^{6-d}\leq r^{6-d}.
			\end{equation*}
			Adding the two cases proves the proposition.
		\end{proof}
		
		\begin{proof}[Proof of Proposition~\ref{prop.var_term_3}]
			Fix $i\in\{1,2\}$ and write $F:=F_i$; recall \eqref{eq.F_pointwise}. By Lemma~\ref{lem.rho_F_basic}, $\rho_i\lesssim1$, so it suffices to show that
			\begin{equation*}
				n^3\sum_{u\in\mathrm{Ann}_{c_1n}\cap\Lambda_{n/2}}g(u)\,\tau_{\mathrm{Ann}_{c_1n}}(u,x)\,F(u)\lesssim r^{6-d}.
			\end{equation*}
			For such $u$ we have $\vert x-u\vert\geq n/2$, so $\tau_{\mathrm{Ann}_{c_1n}}(u,x)\lesssim n^{1-d}$ by Lemma~\ref{lem.box_twopoint_onesided}, $F(u)\lesssim n^{4-d}$ by \eqref{eq.F_pointwise}, and $g(u)\lesssim(1+r_u)^{-2}$ by \eqref{eq.g_inner}. We will split the sum according to whether $u\in\RD_{c_1n}$ or not (note that $\RD_{c_1n}\subset\Lambda_{2c_1n}\subset\Lambda_{n/2}$ by \eqref{eq.w_bound}), and we first prove two counting estimates:
			\begin{equation}\label{eq.RD_count}
				\sum_{w\in\RD_{c_1n}}s_w^{-2}\lesssim n^{d-2}\log n,\qquad\sum_{w\in(\mathrm{Ann}_{c_1n}\cap\Lambda_{n/2})\setminus\RD_{c_1n}}(1+r_w)^{-2}\lesssim n^{d-2}\log n.
			\end{equation}
			For the first one, recall that $s_w=c_1n-a_{(2)}(w)$ and $r_w<s_w$ for $w\in\RD_{c_1n}$. Fix $s\geq1$ and $r<s$. If $r_w=r$ and $s_w=s$, then one coordinate of $w$ is $\pm(c_1n+r)$, and the other $d-1$ coordinates lie in $[-(c_1n-s),c_1n-s]$ with maximal absolute value $c_1n-s$; there are $\lesssim n^{d-2}$ such $w$. Hence
			\begin{equation*}
				\sum_{w\in\RD_{c_1n}}s_w^{-2}\lesssim\sum_{s=1}^{c_1n}s\cdot n^{d-2}s^{-2}\lesssim n^{d-2}\log n.
			\end{equation*}
			For the second one, if $r_w=r\leq c_1n$ and $w\notin\RD_{c_1n}$, then $a_{(2)}(w)\geq c_1n-r$, so one coordinate of $w$ is $\pm(c_1n+r)$, another one has absolute value in $[c_1n-r,c_1n+r]$, and the others are arbitrary in $[-(c_1n+r),c_1n+r]$; there are $\lesssim r\,n^{d-2}$ such $w$. If $r_w=r>c_1n$, there are $\lesssim r^{d-1}$ such $w$. Hence
			\begin{equation*}
				\sum_{w\in(\mathrm{Ann}_{c_1n}\cap\Lambda_{n/2})\setminus\RD_{c_1n}}(1+r_w)^{-2}\lesssim\sum_{r\leq c_1n}r\,n^{d-2}(1+r)^{-2}+\sum_{c_1n<r\leq n/2}r^{d-3}\lesssim n^{d-2}\log n.
			\end{equation*}
			
			\medskip
			\noindent\emph{Case 1: $u\in\RD_{c_1n}$.} By Lemma~\ref{lem.RD_bound}, $\tau_{\mathrm{Ann}_{c_1n}}(u,x)\lesssim n^{1-d}(1+r_u)s_u^{-1}$, and by Lemma~\ref{lem.F_RD}, $F(u)\lesssim n^{4-d}(1+r_u)s_u^{-1}$. Together with $g(u)\lesssim(1+r_u)^{-2}$, we get $g(u)\,\tau_{\mathrm{Ann}_{c_1n}}(u,x)\,F(u)\lesssim n^{5-2d}s_u^{-2}$, so by \eqref{eq.RD_count} the contribution of Case~1 is
			\begin{equation*}
				\lesssim n^3\cdot n^{5-2d}\cdot n^{d-2}\log n=n^{6-d}\log n.
			\end{equation*}
			
			\medskip
			\noindent\emph{Case 2: $u\in(\mathrm{Ann}_{c_1n}\cap\Lambda_{n/2})\setminus\RD_{c_1n}$.} By \eqref{eq.RD_count}, the contribution of Case~2 is
			\begin{equation*}
				\lesssim n^3\cdot n^{1-d}\cdot n^{4-d}\sum_{u}(1+r_u)^{-2}\lesssim n^{8-2d}\cdot n^{d-2}\log n=n^{6-d}\log n.
			\end{equation*}
			Since $n$ is sufficiently large, $n^{6-d}\log n\lesssim r^{6-d}$, which proves the proposition.
		\end{proof}
		
		\begin{proof}[Proof of Proposition~\ref{prop.var_term_4}]
			Fix $i\in\{1,2\}$ and write $F:=F_i$; recall \eqref{eq.F_pointwise}. We split the sum in the proposition into three parts: $w\in\Lambda_n\setminus\Lambda_{n/2}$, $w\in\RD_{c_1n}$, and $w\in(\mathrm{Ann}_{c_1n}\cap\Lambda_{n/2})\setminus\RD_{c_1n}$.
			
			\medskip
			\noindent\emph{Case 1: $w\in\Lambda_n\setminus\Lambda_{n/2}$.} Since $l_w\leq\vert x-w\vert$, \eqref{eq.g_outer} gives $g(w)\lesssim(1+\vert x-w\vert)\,n^{-3}$. Together with \eqref{eq.F_pointwise} and $\#\{w\in\Lambda_n:\vert x-w\vert=m\}\lesssim(1+m)^{d-1}$, the contribution of this case is
			\begin{equation*}
				\lesssim n^{-3}\sum_{m=0}^{2n}(1+m)^d\max\{r,m\}^{8-2d}\lesssim n^{-3}\big(r^{9-d}+n^2\big)\lesssim r^{6-d},
			\end{equation*}
			where we split at $m\lessgtr r$, used $m^{8-d}\leq m$ for $m>r$, and used that $n\geq r$ and $n$ is sufficiently large.
			
			\medskip
			\noindent\emph{Case 2: $w\in\RD_{c_1n}$.} By \eqref{eq.g_inner} and Lemma~\ref{lem.F_RD}, $g(w)F(w)^2\lesssim(1+r_w)^{-2}n^{8-2d}(1+r_w)^2s_w^{-2}=n^{8-2d}s_w^{-2}$. Hence, by \eqref{eq.RD_count}, the contribution of this case is $\lesssim n^{8-2d}\cdot n^{d-2}\log n=n^{6-d}\log n$.
			
			\medskip
			\noindent\emph{Case 3: $w\in(\mathrm{Ann}_{c_1n}\cap\Lambda_{n/2})\setminus\RD_{c_1n}$.} Here $\vert x-w\vert\geq n/2$, so $F(w)\lesssim n^{4-d}$, and $g(w)\lesssim(1+r_w)^{-2}$. By \eqref{eq.RD_count}, the contribution of this case is $\lesssim n^{8-2d}\cdot n^{d-2}\log n=n^{6-d}\log n$.
			
			Since $n$ is sufficiently large, $n^{6-d}\log n\lesssim r^{6-d}$, which proves the proposition.
		\end{proof}

		\bigskip
		\textbf{Acknowledgments.} Gu was supported by the National Natural Science Foundation of China (Nos. 12622111, 12595284, and 12595280).
		TH was supported by a Packard Fellowship in Science and Engineering and NSF grant DMS-1928930.
		This work was initiated at 12th ICSAA (September 8th - 12th, 2025) in Bucharest, and was developed during 54th Saint-Flour Summer School (June 29th - July 11th, 2026). Gu also thanks Wenhao Zhao (EPFL) and Yuchen Liu (PKU) for inspiring discussions at the early stage of this project.

		\bigskip
		\textbf{Statement of AI use.} This project was initiated September 2025. The first draft of the manuscript, including fully detailed proofs, was prepared without AI assistance in writing or in research. Once this was done, we used commercially available AI assistants to review the manuscript. This led to the discovery of many errors (mostly very small, but with some more significant) that were fixed partly with AI assistance; the writing remains predominantly human.
		The most significant AI use was in Section~\ref{sec.lower.large}: The first version contained a mistake that the AI assistant claimed it could not see how to fix, specifically in the analysis of the flow we had constructed. We suggested that it could probably be fixed using a modified definition of the flow based on the EIT approach of \cite{BPP98}, and from this suggestion the assistant was indeed immediately able to fix the argument. We then rewrote the section by hand with this output as a rough guide (and with most of the section preserved from the earlier version). After finishing the paper, we used AI to formalize our results in lean, taking many results from the literature as given; this formalization is available at \url{https://github.com/tomhutchcroft/near-critical-diffusivity}.
		AI assistants were also used to produce the numerical simulations reported in Figure~\ref{fig:numerics}.

		\bibliographystyle{abbrv}
		\bibliography{Ref}
		
	\end{document}